\input graphicx



\input pstricks
\input pst-node
\input pst-plot
\newcount\chapnomb \chapnomb=1
\newcount\parnomb \parnomb=1
\newcount\nomb \nomb=1

\ifx\optionkeymacros\undefined\else\endinput\fi

\catcode`\Œ=\active\defŒ{{\aa}}       
\catcode`\º=\active\defº{\int}        
\catcode`\=\active\def{\c c}        
\catcode`\¶=\active\def¶{\partial}    
\catcode`\Ä=\active\defÄ{\oint}       
\catcode`\Æ=\active\defÆ{\triangle}   
\catcode`\Â=\active\defÂ{\neg}        
\catcode`\µ=\active\defµ{\mu}         
\catcode`\¿=\active\def¿{{\o}}        
\catcode`\¹=\active\def¹{\pi}         
\catcode`\Ï=\active\defÏ{{\oe}}       
\catcode`\§=\active\def§{{\ss}}       
\catcode`\ =\active\def {\dagger}     
\catcode`\Ã=\active\defÃ{\sqrt}       
\catcode`\·=\active\def·{\Sigma}      
\catcode`\Å=\active\defÅ{\approx}     
\catcode`\½=\active\def½{\Omega}      
\catcode`\£=\active\def£{{\it\$}}     
\catcode`\°=\active\def°{\infty}      
\catcode`\¤=\active\def¤{{\S}}        
\catcode`\¦=\active\def¦{{\P}}        
\catcode`\¥=\active\def¥{\bullet}     
\catcode`\»=\active\def»{\leavevmode\raise.585ex\hbox{\b a}}      
\catcode`\¼=\active\def¼{\leavevmode\raise.6ex\hbox{\b o}}        
\catcode`\­=\active\def­{\not=}       
\catcode`\²=\active\def²{\leq}        
\catcode`\³=\active\def³{\geq}        
\catcode`\Ö=\active\defÖ{\div}        
\catcode`\É=\active\defÉ{{\dots}}     
\catcode`\¾=\active\def¾{{\ae}}       
\catcode`\Ç=\active\defÇ{\ll}         
\catcode`\Ò=\active\defÒ{``}          
\catcode`\Á=\active\defÁ{!`}          
\catcode`\¢=\active\def¢{\rlap/c}     
\catcode`\Ô=\active\defÔ{`}           
\catcode`\Õ=\active\defÕ{'}           


\catcode`\=\active\def{{\AA}}       
\catcode`\'=\active\def'{\c C}        
\catcode`\¯=\active\def¯{{\O}}        
\catcode`\¸=\active\def¸{\Pi}         
\catcode`\Î=\active\defÎ{{\OE}}       
\catcode`\®=\active\def®{{\AE}}       
\catcode`\×=\active\def×{\diamond}    
\catcode`\¡=\active\def¡{\accent'27}  
\catcode`\Ó=\active\defÓ{''}          
\catcode`\±=\active\def±{\pm}         
\catcode`\È=\active\defÈ{\gg}         
\catcode`\À=\active\defÀ{?`}          
\catcode`\Ð=\active\defÐ{--}          
\catcode`\Ñ=\active\defÑ{---}         


\catcode`\Š=\active\defŠ{\"a}        
\catcode`\'=\active\def'{\"e}        
\catcode`\•=\active\def•{\"{\i}}     
\catcode`\š=\active\defš{\"o}        
\catcode`\Ÿ=\active\defŸ{\"u}        
\catcode`\Ø=\active\defØ{\"y}        
\catcode`\€=\active\def€{\"A}        
\catcode`\…=\active\def…{\"O}        
\catcode`\†=\active\def†{\"U}        
\catcode`\‡=\active\def‡{\'a}        
\catcode`\Ž=\active\defŽ{\'e}        
\catcode`\'=\active\def'{\'{\i}}     
\catcode`\—=\active\def—{\'o}        
\catcode`\œ=\active\defœ{\'u}        
\catcode`\ƒ=\active\defƒ{\'E}        
\catcode`\ˆ=\active\defˆ{\`a}        
\catcode`\=\active\def{\`e}        
\catcode`\"=\active\def"{\`{\i}}     
\catcode`\˜=\active\def˜{\`o}        
\catcode`\=\active\def{\`u}        
\catcode`\Ë=\active\defË{\`A}        
\catcode`\‹=\active\def‹{\~a}        
\catcode`\–=\active\def–{\~n}        
\catcode`\›=\active\def›{\~o}        
\catcode`\Ì=\active\defÌ{\~A}        
\catcode`\"=\active\def"{\~N}        
\catcode`\Í=\active\defÍ{\~O}        
\catcode`\‰=\active\def‰{\^a}        
\catcode`\=\active\def{\^e}        
\catcode`\"=\active\def"{\^{\i}}     
\catcode`\™=\active\def™{\^o}        
\catcode`\ž=\active\defž{\^u}        

\let\optionkeymacros\null
\catcode`\@=11
\def\W@{\immediate\write\sixt@@n}
\def\spaces@{\space\space\space\space\space}
\def\spaces@@{\spaces@\spaces@\spaces@}

\voffset=-8mm
\hoffset=-1mm
\hsize=16cm
\vsize=23cm

\newdimen\ex@
\newif\iftop@
\newif\ifbot@

\def\thfrac#1#2#3#4#5{\ifmmode{\vbox{\hbox{$#4$}\kern#2\p@}
            \above#1\ex@\vbox{\kern#3\p@\hbox{$#5$}}}\else\relax\fi}
\def\thunderline#1#2#3{\ifmmode\vtop{\ialign{##\crcr
     $\hfil\displaystyle{#3}\hfil$\crcr
     \noalign{\kern#2\p@\nointerlineskip}
     \leaders\hrule height#1\ex@\hfill\crcr}}
     \else$\setbox\z@\hbox{#3} \dp\z@=0pt \thunderline{#1}{#2}{\box\z@}$\fi}
\def\thoverline#1#2#3{\ifmmode\vbox{\ialign{##\crcr
     \leaders\hrule height
     #1\ex@\hfill\crcr\noalign{\kern#2\p@\nointerlineskip}
     $\hfil\displaystyle{#3}\hfil$\crcr}}
     \else$\thoverline{#1}{#2}{\hbox{#3}}$\fi}
\def\thsemiboxed#1#2#3{\ifmmode\setbox\z@\hbox{$\displaystyle{#3}$}
     \hbox{\lower#1\ex@\hbox{\lower#2\p@\hbox{\lower\dp\z@
     \hbox{\vbox{\hbox{\vbox{\box\z@\kern#2\p@}%
     \kern#2\p@\vrule width#1\ex@}\hrule height#1\ex@}}}}}
     \else$\setbox\z@\hbox{#3}\dp\z@=0pt\thsemiboxed{#1}{#2}{\box\z@}$\fi}
\def\thboxed#1#2#3{\ifmmode\setbox\z@\hbox{$\displaystyle{#3}$}
     \hbox{\lower#1\ex@ \hbox{\lower#2\p@\hbox{\lower\dp\z@
     \hbox{\vbox{\hrule height#1\ex@\hbox{\vrule width#1\ex@
     \kern#2\p@\vbox{\kern#2\p@\box\z@\kern#2\p@}%
     \kern#2\p@\vrule width#1\ex@}\hrule height#1\ex@}}}}}
     \else$\setbox\z@\hbox{#3}\dp\z@=0pt\thboxed{#1}{#2}{\box\z@}$\fi}
\def\thsquared#1#2#3{\ifmmode\setbox\z@\hbox{$\displaystyle{#3}$}
       \vcenter{\hrule height#1\ex@\hbox{\vrule width#1\ex@
        \kern#2\p@\vbox{\kern#2\p@\box\z@\kern#2\p@}%
        \kern#2\p@\vrule width#1\ex@}\hrule height#1\ex@}
      \else$\setbox\z@\hbox{#3}\dp\z@=0pt\thsquared{#1}{#2}{\box\z@}$\fi}

\def\NoBlackBoxes{\global\overfullrule\z@}
\def\BlackBoxes{\global\overfullrule5\p@}
\def\,{\relax\ifmmode\mskip\thinmuskip\relax\else\kern.16667em\fi}
\def\!{\relax\ifmmode\mskip-\thinmuskip\relax\else\kern-.16667em\fi}
\def\>{\relax\ifmmode\mskip\medmuskip\relax\else\kern.222222em\fi}
\def\;{\relax\ifmmode\mskip\thickmuskip\relax\else\kern.27777em\fi}
\def\negmedspace{\relax\ifmmode\mskip-\medmuskip\relax\else\kern-.222222em\fi}
\def\negthickspace{\relax\ifmmode\mskip-\thickmuskip\relax\else
 \kern-.27777em\fi}
\def\snug{\unskip\kern-\mathsurround}
\let\thinspace\,
\let\negthinspace\!
\let\medspace\>
\let\thickspace\;

\def\topsmash{\top@true\bot@false\smash@}
\def\botsmash{\top@false\bot@true\smash@}
\def\smash{\top@true\bot@true\smash@}
\def\smash@{\relax\ifmmode\def\next{\mathpalette\mathsm@sh}\else
 \let\next\makesm@sh\fi\next}
\def\finsm@sh{\iftop@\ht\z@\z@\fi\ifbot@\dp\z@\z@\fi\box\z@}

\newfam\msafam  
\newfam\msbfam
\newfam\eufmfam  
\newfam\sffam  
\newfam\bifam

\newdimen\t@ille

\def\taille#1{\t@ille=#1pt
\p@=.1\t@ille
\ex@=.02\t@ille
\font\tenrm=cmr10 at \t@ille                  
\font\sevenrm=cmr7 at .7\t@ille
\font\fiverm=cmr5 at .5\t@ille
\font\teni=cmmi10 at \t@ille                  
\font\seveni=cmmi7 at .7\t@ille
\font\fivei=cmmi5 at .5\t@ille
\font\tensy=cmsy10 at \t@ille                 
\font\sevensy=cmsy7 at .7\t@ille
\font\fivesy=cmsy5 at .5\t@ille
\font\tenex=cmex10 at \t@ille                 
\font\tenbf=cmbx10 at \t@ille  								   			 
\font\sevenbf=cmbx7 at .7\t@ille
\font\fivebf=cmbx5 at .5\t@ille
\font\tentt=cmtt10 at \t@ille                 
\font\tenbi=cmbxti10 at \t@ille
\font\tensl=cmsl10 at \t@ille                 
\font\tenit=cmti10 at \t@ille                 
\font\tensf=cmss10 at \t@ille                 
\font\tenmsa=msam10 at \t@ille
\font\sevenmsa=msam7 at .7\t@ille
\font\fivemsa=msam5 at .5\t@ille
\font\tenmsb=msbm10 at \t@ille                
\font\sevenmsb=msbm7 at .7\t@ille
\font\fivemsb=msbm5 at .5\t@ille
\font\teneufm=eufm10 at \t@ille
\font\seveneufm=eufm7 at .7\t@ille
\font\fiveeufm=eufm5 at .5\t@ille
\textfont0=\tenrm \scriptfont0=\sevenrm \scriptscriptfont0=\fiverm
\textfont1=\teni  \scriptfont1=\seveni  \scriptscriptfont1=\fivei
\textfont2=\tensy \scriptfont2=\sevensy \scriptscriptfont2=\fivesy
\textfont3=\tenex \scriptfont3=\tenex   \scriptscriptfont3=\tenex
\textfont\itfam=\tenit
\textfont\slfam=\tensl
\textfont\bffam=\tenbf   \scriptfont\bffam=\sevenbf
\scriptscriptfont\bffam=\fivebf
\textfont\ttfam=\tentt
\textfont\msafam=\tenmsa \scriptfont\msafam=\sevenmsa
\scriptscriptfont\msafam=\fivemsa
\textfont\msbfam=\tenmsb \scriptfont\msbfam=\sevenmsb
\scriptscriptfont\msbfam=\fivemsb
\textfont\eufmfam=\teneufm \scriptfont\eufmfam=\seveneufm
\scriptscriptfont\eufmfam=\fiveeufm
\textfont\sffam=\tensf
\textfont\bifam=\tenbi
\abovedisplayskip=1.2\t@ille plus .3\t@ille minus .9\t@ille
\abovedisplayshortskip=0pt plus .3\t@ille
\belowdisplayskip=1.2\t@ille plus .3\t@ille minus .9\t@ille
\belowdisplayshortskip=.7\t@ille plus .3\t@ille minus .4\t@ille
\jot=.3\t@ille
\smallskipamount=.3\t@ille plus .1\t@ille minus .1\t@ille
\medskipamount=.6\t@ille plus .2\t@ille minus .2\t@ille
\bigskipamount=1.2\t@ille plus .4\t@ille minus .4\t@ille
\normalbaselineskip=1.2\t@ille
\normallineskip=.1\t@ille
\normallineskiplimit=0pt
\normalbaselines
\tenrm}
\taille{10}

\let\teneuf=\teneufm

\let\euffam=\eufmfam

\def\euf{\fam\euffam\teneuf}

\def\hexnumber@#1{\ifcase#1 0\or1\or2\or3\or4\or5\or6\or7\or8\or9\or
 A\or B\or C\or D\or E\or F\fi}
\edef\bffam@{\hexnumber@\bffam}
\edef\msa@{\hexnumber@\msafam}
\edef\msb@{\hexnumber@\msbfam}

\mathchardef\arr@bas="0040

\mathchardef\Gamma="0000
\mathchardef\Delta="0001
\mathchardef\Theta="0002
\mathchardef\Lambda="0003
\mathchardef\Xi="0004
\mathchardef\Pi="0005
\mathchardef\Sigma="0006
\mathchardef\Upsilon="0007
\mathchardef\Phi="0008
\mathchardef\Psi="0009
\mathchardef\Omega="000A
\mathchardef\varGamma="0100
\mathchardef\varDelta="0101
\mathchardef\varTheta="0102
\mathchardef\varLambda="0103
\mathchardef\varXi="0104
\mathchardef\varPi="0105
\mathchardef\varSigma="0106
\mathchardef\varUpsilon="0107
\mathchardef\varPhi="0108
\mathchardef\varPsi="0109
\mathchardef\varOmega="010A
\mathchardef\boldGamma="0\bffam@00
\mathchardef\boldDelta="0\bffam@01
\mathchardef\boldTheta="0\bffam@02
\mathchardef\boldLambda="0\bffam@03
\mathchardef\boldXi="0\bffam@04
\mathchardef\boldPi="0\bffam@05
\mathchardef\boldSigma="0\bffam@06
\mathchardef\boldUpsilon="0\bffam@07
\mathchardef\boldPhi="0\bffam@08
\mathchardef\boldPsi="0\bffam@09
\mathchardef\boldOmega="0\bffam@0A

\mathchardef\boxdot="2\msa@00
\mathchardef\boxplus="2\msa@01
\mathchardef\boxtimes="2\msa@02
\mathchardef\square="0\msa@03
\mathchardef\blacksquare="0\msa@04
\mathchardef\centerdot="2\msa@05
\mathchardef\lozenge="0\msa@06
\mathchardef\blacklozenge="0\msa@07
\mathchardef\circlearrowright="3\msa@08
\mathchardef\circlearrowleft="3\msa@09
\mathchardef\rightleftharpoons="3\msa@0A
\mathchardef\leftrightharpoons="3\msa@0B
\mathchardef\boxminus="2\msa@0C
\mathchardef\Vdash="3\msa@0D
\mathchardef\Vvdash="3\msa@0E
\mathchardef\vDash="3\msa@0F
\mathchardef\twoheadrightarrow="3\msa@10
\mathchardef\twoheadleftarrow="3\msa@11
\mathchardef\leftleftarrows="3\msa@12
\mathchardef\rightrightarrows="3\msa@13
\mathchardef\upuparrows="3\msa@14
\mathchardef\downdownarrows="3\msa@15
\mathchardef\upharpoonright="3\msa@16

\mathchardef\downharpoonright="3\msa@17
\mathchardef\upharpoonleft="3\msa@18
\mathchardef\downharpoonleft="3\msa@19
\mathchardef\rightarrowtail="3\msa@1A
\mathchardef\leftarrowtail="3\msa@1B
\mathchardef\leftrightarrows="3\msa@1C
\mathchardef\rightleftarrows="3\msa@1D
\mathchardef\Lsh="3\msa@1E
\mathchardef\Rsh="3\msa@1F
\mathchardef\rightsquigarrow="3\msa@20
\mathchardef\leftrightsquigarrow="3\msa@21
\mathchardef\looparrowleft="3\msa@22
\mathchardef\looparrowright="3\msa@23
\mathchardef\circeq="3\msa@24
\mathchardef\succsim="3\msa@25
\mathchardef\gtrsim="3\msa@26
\mathchardef\gtrapprox="3\msa@27
\mathchardef\multimap="3\msa@28
\mathchardef\therefore="3\msa@29
\mathchardef\because="3\msa@2A
\mathchardef\doteqdot="3\msa@2B

\mathchardef\triangleq="3\msa@2C
\mathchardef\precsim="3\msa@2D
\mathchardef\lesssim="3\msa@2E
\mathchardef\lessapprox="3\msa@2F
\mathchardef\eqslantless="3\msa@30
\mathchardef\eqslantgtr="3\msa@31
\mathchardef\curlyeqprec="3\msa@32
\mathchardef\curlyeqsucc="3\msa@33
\mathchardef\preccurlyeq="3\msa@34
\mathchardef\leqq="3\msa@35
\mathchardef\leqslant="3\msa@36
\mathchardef\lessgtr="3\msa@37
\mathchardef\backprime="0\msa@38
\mathchardef\risingdotseq="3\msa@3A
\mathchardef\fallingdotseq="3\msa@3B
\mathchardef\succcurlyeq="3\msa@3C
\mathchardef\geqq="3\msa@3D
\mathchardef\geqslant="3\msa@3E
\mathchardef\gtrless="3\msa@3F
\mathchardef\sqsubset="3\msa@40
\mathchardef\sqsupset="3\msa@41
\mathchardef\vartriangleright="3\msa@42
\mathchardef\vartriangleleft ="3\msa@43
\mathchardef\trianglerighteq="3\msa@44
\mathchardef\trianglelefteq="3\msa@45
\mathchardef\bigstar="0\msa@46
\mathchardef\between="3\msa@47
\mathchardef\blacktriangledown="0\msa@48
\mathchardef\blacktriangleright="3\msa@49
\mathchardef\blacktriangleleft="3\msa@4A
\mathchardef\vartriangle="0\msa@4D
\mathchardef\blacktriangle="0\msa@4E
\mathchardef\triangledown="0\msa@4F
\mathchardef\eqcirc="3\msa@50
\mathchardef\lesseqgtr="3\msa@51
\mathchardef\gtreqless="3\msa@52
\mathchardef\lesseqqgtr="3\msa@53
\mathchardef\gtreqqless="3\msa@54
\mathchardef\Rrightarrow="3\msa@56
\mathchardef\Lleftarrow="3\msa@57
\mathchardef\veebar="2\msa@59
\mathchardef\barwedge="2\msa@5A
\mathchardef\doublebarwedge="2\msa@5B
\mathchardef\angle="0\msa@5C
\mathchardef\measuredangle="0\msa@5D
\mathchardef\sphericalangle="0\msa@5E
\mathchardef\varpropto="3\msa@5F
\mathchardef\smallsmile="3\msa@60
\mathchardef\smallfrown="3\msa@61
\mathchardef\Subset="3\msa@62
\mathchardef\Supset="3\msa@63
\mathchardef\Cup="2\msa@64

\mathchardef\Cap="2\msa@65

\mathchardef\curlywedge="2\msa@66
\mathchardef\curlyvee="2\msa@67
\mathchardef\leftthreetimes="2\msa@68
\mathchardef\rightthreetimes="2\msa@69
\mathchardef\subseteqq="3\msa@6A
\mathchardef\supseteqq="3\msa@6B
\mathchardef\bumpeq="3\msa@6C
\mathchardef\Bumpeq="3\msa@6D
\mathchardef\lll="3\msa@6E

\mathchardef\ggg="3\msa@6F

\mathchardef\circledS="0\msa@73
\mathchardef\pitchfork="3\msa@74
\mathchardef\dotplus="2\msa@75
\mathchardef\backsim="3\msa@76
\mathchardef\backsimeq="3\msa@77
\mathchardef\complement="0\msa@7B
\mathchardef\intercal="2\msa@7C
\mathchardef\circledcirc="2\msa@7D
\mathchardef\circledast="2\msa@7E
\mathchardef\circleddash="2\msa@7F
\def\ulcorner{\delimiter"4\msa@70\msa@70 }
\def\urcorner{\delimiter"5\msa@71\msa@71 }
\def\llcorner{\delimiter"4\msa@78\msa@78 }
\def\lrcorner{\delimiter"5\msa@79\msa@79 }

\mathchardef\y@n="0\msa@55
\mathchardef\ch@ckmark="0\msa@58
\mathchardef\c@rcledR="0\msa@72
\mathchardef\m@ltese="0\msa@7A
\def\yen{\relax\ifmmode \y@n\else $\m@th\y@n$ \fi}
\def\checkmark{\relax\ifmmode \ch@ckmark\else $\m@th\ch@ckmark$ \fi}
\def\circledR{\relax\ifmmode \c@rcledR\else $\m@th\c@rcledR$ \fi}
\def\maltese{\relax\ifmmode \m@ltese\else $\m@th\m@ltese$ \fi}

\mathchardef\lvertneqq="3\msb@00
\mathchardef\gvertneqq="3\msb@01
\mathchardef\nleq="3\msb@02
\mathchardef\ngeq="3\msb@03
\mathchardef\nless="3\msb@04
\mathchardef\ngtr="3\msb@05
\mathchardef\nprec="3\msb@06
\mathchardef\nsucc="3\msb@07
\mathchardef\lneqq="3\msb@08
\mathchardef\gneqq="3\msb@09
\mathchardef\nleqslant="3\msb@0A
\mathchardef\ngeqslant="3\msb@0B
\mathchardef\lneq="3\msb@0C
\mathchardef\gneq="3\msb@0D
\mathchardef\npreceq="3\msb@0E
\mathchardef\nsucceq="3\msb@0F
\mathchardef\precnsim="3\msb@10
\mathchardef\succnsim="3\msb@11
\mathchardef\lnsim="3\msb@12
\mathchardef\gnsim="3\msb@13
\mathchardef\nleqq="3\msb@14
\mathchardef\ngeqq="3\msb@15
\mathchardef\precneqq="3\msb@16
\mathchardef\succneqq="3\msb@17
\mathchardef\precnapprox="3\msb@18
\mathchardef\succnapprox="3\msb@19
\mathchardef\lnapprox="3\msb@1A
\mathchardef\gnapprox="3\msb@1B
\mathchardef\nsim="3\msb@1C
\mathchardef\napprox="3\msb@1D
\mathchardef\ncong="3\msb@1D
\def\napprox{\not\approx}
\mathchardef\varsubsetneq="3\msb@20
\mathchardef\varsupsetneq="3\msb@21
\mathchardef\nsubseteqq="3\msb@22
\mathchardef\nsupseteqq="3\msb@23
\mathchardef\subsetneqq="3\msb@24
\mathchardef\supsetneqq="3\msb@25
\mathchardef\varsubsetneqq="3\msb@26
\mathchardef\varsupsetneqq="3\msb@27
\mathchardef\subsetneq="3\msb@28
\mathchardef\supsetneq="3\msb@29
\mathchardef\nsubseteq="3\msb@2A
\mathchardef\nsupseteq="3\msb@2B
\mathchardef\nparallel="3\msb@2C
\mathchardef\nmid="3\msb@2D
\mathchardef\nshortmid="3\msb@2E
\mathchardef\nshortparallel="3\msb@2F
\mathchardef\nvdash="3\msb@30
\mathchardef\nVdash="3\msb@31
\mathchardef\nvDash="3\msb@32
\mathchardef\nVDash="3\msb@33
\mathchardef\ntrianglerighteq="3\msb@34
\mathchardef\ntrianglelefteq="3\msb@35
\mathchardef\ntriangleleft="3\msb@36
\mathchardef\ntriangleright="3\msb@37
\mathchardef\nleftarrow="3\msb@38
\mathchardef\nrightarrow="3\msb@39
\mathchardef\nLeftarrow="3\msb@3A
\mathchardef\nRightarrow="3\msb@3B
\mathchardef\nLeftrightarrow="3\msb@3C
\mathchardef\nleftrightarrow="3\msb@3D
\mathchardef\divideontimes="2\msb@3E
\mathchardef\varnothing="0\msb@3F
\mathchardef\nexists="0\msb@40
\mathchardef\mho="0\msb@66
\mathchardef\eth="0\msb@67
\mathchardef\eqsim="3\msb@68
\mathchardef\beth="0\msb@69
\mathchardef\gimel="0\msb@6A
\mathchardef\daleth="0\msb@6B
\mathchardef\lessdot="3\msb@6C
\mathchardef\gtrdot="3\msb@6D
\mathchardef\ltimes="2\msb@6E
\mathchardef\rtimes="2\msb@6F
\mathchardef\shortmid="3\msb@70
\mathchardef\shortparallel="3\msb@71
\mathchardef\smallsetminus="2\msb@72
\mathchardef\thicksim="3\msb@73
\mathchardef\thickapprox="3\msb@74
\mathchardef\approxeq="3\msb@75
\mathchardef\succapprox="3\msb@76
\mathchardef\precapprox="3\msb@77
\mathchardef\curvearrowleft="3\msb@78
\mathchardef\curvearrowright="3\msb@79
\mathchardef\digamma="0\msb@7A
\mathchardef\varkappa="0\msb@7B
\mathchardef\hslash="0\msb@7D
\mathchardef\hbar="0\msb@7E
\mathchardef\backepsilon="3\msb@7F

\mathchardef\bbk="0\msb@7C
\mathchardef\bbA="0\msb@41
\mathchardef\bbB="0\msb@42
\mathchardef\bbC="0\msb@43
\mathchardef\bbD="0\msb@44
\mathchardef\bbE="0\msb@45
\mathchardef\bbF="0\msb@46
\mathchardef\bbG="0\msb@47
\mathchardef\bbH="0\msb@48
\mathchardef\bbI="0\msb@49
\mathchardef\bbJ="0\msb@4A
\mathchardef\bbK="0\msb@4B
\mathchardef\bbL="0\msb@4C
\mathchardef\bbM="0\msb@4D
\mathchardef\bbN="0\msb@4E
\mathchardef\bbO="0\msb@4F
\mathchardef\bbP="0\msb@50
\mathchardef\bbQ="0\msb@51
\mathchardef\bbR="0\msb@52
\mathchardef\bbS="0\msb@53
\mathchardef\bbT="0\msb@54
\mathchardef\bbU="0\msb@55
\mathchardef\bbV="0\msb@56
\mathchardef\bbW="0\msb@57
\mathchardef\bbX="0\msb@58
\mathchardef\bbY="0\msb@59
\mathchardef\bbZ="0\msb@5A

\def\k{\relax\ifmmode\bbk\else $\bbk$\fi}
\def\F{\relax\ifmmode\bbF\else $\bbF$\fi}
\def\N{\relax\ifmmode\bbN\else $\bbN$\fi}
\def\Z{\relax\ifmmode\bbZ\else $\bbZ$\fi}
\def\Q{\relax\ifmmode\bbQ\else $\bbQ$\fi}
\def\R{\relax\ifmmode\bbR\else $\bbR$\fi}
\def\C{\relax\ifmmode\bbC\else $\bbC$\fi}
\def\T{\relax\ifmmode\bbT\else $\bbT$\fi}
\def\arrabas{\relax\ifmmode\arr@bas\else $\m@th\arr@bas$\fi}

\let\B\=
\let\D\.

\edef\msafam@{\hexnumber@\msafam}
\mathchardef\dabar@"0\msafam@39
\def\dashrightarrow{\mathrel{\dabar@\dabar@\mathchar"0\msafam@4B}}
\def\dashleftarrow{\mathrel{\mathchar"0\msafam@4C\dabar@\dabar@}}
\edef\msbfam@{\hexnumber@\msbfam}

\def\widehat#1{\setbox\z@\hbox{$\m@th#1$}%
 \ifdim\wd\z@>\tw@ em\mathaccent"0\msbfam@5B{#1}%
 \else\mathaccent"0362{#1}\fi}
\def\widetilde#1{\setbox\z@\hbox{$\m@th#1$}%
 \ifdim\wd\z@>\tw@ em\mathaccent"0\msbfam@5D{#1}%
 \else\mathaccent"0365{#1}\fi}

\def\newcodes@{\catcode`\\=12 \catcode`\{=12 \catcode`\}=12 \catcode`\#=12
 \catcode`\%=12\relax}
\def\oldcodes@{\catcode`\\=0 \catcode`\{=1 \catcode`\}=2 \catcode`\#=6
 \catcode`\%=14\relax}
\def\comment{\newcodes@\endlinechar=10 \comment@}
{\lccode`\!=`\\
\lowercase{\gdef\comment@#1^^J{\comment@@#1!endcomment\comment@@@}%
\gdef\comment@@#1!endcomment{\futurelet\next\comment@@@}%
\gdef\comment@@@#1\comment@@@{\ifx\next\comment@@@\let
\next=\comment@\else\def\next{\oldcodes@\endlinechar=`\^^M\relax}%
 \fi\next}}}

\font\truetenrm=cmr10
\footline={\hss\truetenrm\folio\hss}

\catcode`\@=12

\font\para=cmbx12 at 18pt

\font\soustitre=cmbx12 at 13pt

\font\toto=msam10 at 20pt
\font\tot=msbm10 at 10pt

\NoBlackBoxes

\def\OO#1{{\reg O}\ifmmode _{#1}\else $_{#1}$\fi}
 
\let \lra=\longrightarrow

\def\hra{\lhook\joinrel\lra}
\def\strot{\mathchoice{\vrule width-1pt height 8.5pt depth 3.5pt}%
   {\vrule width-1pt height 8.5pt depth 3.5pt}{\vrule width-0.7pt height 5.95pt depth 2.45pt}%
   {\vrule width-0.5pt height 4.25pt depth 1.75pt}}

\def\retr#1,{\leavevmode\kern-5mm{\bf #1}}
\def\prop {\retr Proposition {\the\chapnomb}.{\the\parnomb}.{\the\nomb},\global\advance\nomb by 1\par\nobreak\smallskip}
\def\th {\retr ThŽorme {\the\chapnomb}.{\the\parnomb}.{\the\nomb},\global\advance\nomb by 1\par\nobreak\smallskip}
\def\defi {\retr DŽfinition {\the\chapnomb}.{\the\parnomb}.{\the\nomb},\global\advance\nomb by 1\par\nobreak\smallskip}
\def\defis {\retr DŽfinitions {\the\chapnomb}.{\the\parnomb}.{\the\nomb},\global\advance\nomb by 1\par\nobreak\smallskip}
\def\coro {\retr Corollaire {\the\chapnomb}.{\the\parnomb}.{\the\nomb},\global\advance\nomb by 1\par\nobreak\smallskip}
\def\lem {\retr Lemme {\the\chapnomb}.{\the\parnomb}.{\the\nomb},\global\advance\nomb by 1\par\nobreak\smallskip}

\def\propa {\retr Proposition A{\the\chapnomb}.{\the\nomb},\global\advance\nomb by 1\par\nobreak\smallskip}
\def\tha {\retr ThŽorme A{\the\chapnomb}.{\the\nomb},\global\advance\nomb by 1\par\nobreak\smallskip}
\def\defia {\retr DŽfinition A{\the\chapnomb}.{\the\nomb},\global\advance\nomb by 1\par\nobreak\smallskip}
\def\defisa {\retr DŽfinitions A{\the\chapnomb}.{\the\nomb},\global\advance\nomb by 1\par\nobreak\smallskip}
\def\coroa {\retr Corollaire A{\the\chapnomb}.{\the\nomb},\global\advance\nomb by 1\par\nobreak\smallskip}
\def\lema {\retr Lemme A{\the\chapnomb}.{\the\nomb},\global\advance\nomb by 1\par\nobreak\smallskip}

\def\rems{\retr Remarques :,\par\nobreak\smallskip}
\def\rem{\retr Remarque :,\par\nobreak\smallskip}

\def\theo {\retr ThŽorme {\the\nomb},\global\advance\nomb by 1\par\nobreak\smallskip}
\def\defin {\retr DŽfinition {\the\nomb},\global\advance\nomb by 1\par\nobreak\smallskip}
\def\lemm {\retr Lemme {\the\nomb},\global\advance\nomb by 1\par\nobreak\smallskip}
\def\corol {\retr Corollaire {\the\nomb},\global\advance\nomb by 1\par\nobreak\smallskip}

\def\dst{\displaystyle}

\def\cucuplus{\mathrel{\lower .9pt\rlap{$\square$}+}}

\def\<{\mathopen<}
\def\>{\mathclose>}
\def\P{\hbox{\doub P}}
\def\rom#1{\uppercase\expandafter{\romannumeral #1}} 
\long\def\art#1{{\parindent0pt\item{#1}}\hangindent=7mm\hangafter=-20}
\long\def\artart#1{{\parindent0pt\item{#1}}\hangindent=12mm\hangafter=-20}

\def\momo#1\over#2{\mathrel{\mathop{#2}\limits^{\lower5pt\hbox{$\scriptstyle #1$}}}}

\smallskipamount4pt plus1pt minus1pt
\medskipamount8pt plus2pt minus2pt
\bigskipamount16pt plus4pt minus4pt
\baselineskip16pt

\catcode`\Ê=\active
\def\dfleche#1#2{\smash{\mathop{\hbox to 9mm{\rightarrowfill}}\limits^{#1}_{#2}}}
\def\bfleche#1#2{\llap{$\vcenter{\hbox{$\scriptstyle#1$}}$}\left\downarrow\vbox to 4.5mm{}\right.
  \rlap{$\vcenter{\hbox{$\!\scriptstyle#2$}}$}}
\def\gfleche#1#2{\smash{\mathop{\hbox to 9mm{\leftarrowfill}}\limits^{\ #1}_{\ #2}}}
\def\hfleche#1#2{\llap{$\vcenter{\hbox{$\scriptstyle#1$}}$}\left\uparrow\vbox to 4.5mm{}\right.
  \rlap{$\vcenter{\hbox{$\!\scriptstyle#2$}}$}}

\def\diagram#1{{\def\normalbaselines{\baselineskip18pt \lineskip3pt \lineskiplimit3pt} \matrix{#1}}}

\def\carre#1,#2,#3,#4;#5,#6,#7,#8.{\diagram{#1&\dfleche{#5}{}&#2\cr
  \bfleche{#6}{}&&\bfleche{}{#7}\cr
  #3&\dfleche{#8}{}&#4\cr}}
\def\cartesien#1,#2,#3;#4,#5,#6,#7.{\carre{#2\times_#3 #1},#1,#2,#3;
  #4,#5,#6,#7.}
\def\invcarre#1,#2,#3,#4;#5,#6,#7,#8.{\diagram{#1&\gfleche{#5}{}&#2\cr
  \hfleche{#6}{}&&\hfleche{}{#7}\cr
  #3&\gfleche{#8}{}&#4\cr}}
\def\rest#1{\lower 2pt \hbox{$|$}_{#1}}

\def\aa{{\euf a}}
\def\qq{{{\euf P}\strot}}
\def\rr{{\euf r}}
\def\bb{{\euf b}}

\def\pmod#1{\allowbreak\mkern5mu({\rm mod}\,\,#1)}

\def\fmin#1{\buildrel #1\over {\hbox{\raise 3pt\hbox{\tot p}\hskip-2pt\toto\char 032}}}
\def\zeta{{\mathchar"0110\strot}}

\font\smcap=cmcsc10 at 10pt
\font\para=cmbx12 at 18pt
\def\O{\hbox{$\cal O$}}

\def\dst{\displaystyle}

\def\gfP{\relax\ifmmode\bbP\else $\bbP$\fi}
\def\gP{{\euf P}}
\def\gQ{{\euf Q}}
\def\P{{\euf p}}

\def\Log{{\rm Log}}
\def\llog{{\bf log}}
\def\aa{{\euf a}}
\def\bb{{\euf b}}
\def\cc{{\euf c}}
\def\qq{{\euf q}}
\def\rr{{\euf r}}

\def\notdivi{\,/\!\!\!\!\;|\,}

\def\ff{{\euf f}}
\def\okp{{O_{\P}}}
\def\okpp{{O_{(\P)}}}

\def\notdiv{\,\nmid\,}
\def\notdivi{\,\nmid\,}

\def\pmodast#1{\allowbreak\mkern18mu({\rm mod^*}\,\,#1)}

\def\lra{\longrightarrow}

\def\ggP{{\bf P}}
\def\m{{\euf m}}
\def\n{{\euf n}}

\def\qed{\hfill \hbox{\psset{xunit=0.15cm,yunit=0.15cm}
\psline(-3.293,0.707)(-0.293,0.707)
\psline(-2.5858,1.4142)(0.4142,1.4142)
\psline(-3,0)(-0.8787,2.1213)
\psline(-2,0)(0.1213,2.1213)}}

\def\prop {{\soustitre Proposition ({\the\chapnomb}.{\the\nomb})}\global\advance\nomb by 1}
\def\th {{\soustitre ThŽorme ({\the\chapnomb}.{\the\nomb})}\global\advance\nomb by 1}
\def\defi {{\soustitre DŽfinition ({\the\chapnomb}.{\the\nomb})}\global\advance\nomb by 1}
\def\defis {{\soustitre DŽfinitions ({\the\chapnomb}.{\the\nomb})}\global\advance\nomb by 1}
\def\coro {{\soustitre Corollaire ({\the\chapnomb}.{\the\nomb})}\global\advance\nomb by 1}
\def\lem {{\soustitre Lemme ({\the\chapnomb}.{\the\nomb})}\global\advance\nomb by 1}
\def\rem {{\soustitre Remarque}}
\def\rems {{\soustitre Remarques}}

\newcount\chapnomb \chapnomb=0
\newcount\parnomb \nomb=1

\def\gaga{\the\pageno}


\def\argb{(0.2)}
\def\argc{(0.3)}
\def\argd{(0.4)}
\def\argf{(0.6)}
\def\argg{(0.7)}
\def\argh{(0.8)}
\def\argi{(0.9)}
\def\argj{(0.10)}
\def\argk{(0.11)}
\def\argl{(0.12)}
\def\argm{(0.14)}
\def\argo{(0.13)}
\def\argp{(0.16)}
\def\aargp{(0.17)}

\def\argq{(1.1)}
\def\argr{(1.2)}
\def\argt{(1.4)}
\def\argu{(1.5)}
\def\argw{(1.7)}
\def\argx{(1.8)}
\def\argz{(1.10)}
\def\argaa{(1.11)}

\def\argac{(2.2)}
\def\argad{(2.3)}
\def\argaf{(2.5)}
\def\argag{(2.6)}
\def\argah{(2.7)}
\def\argaj{(2.9)}
\def\argak{(2.10)}
\def\argal{(2.11)}
\def\argam{(2.12)}
\def\argan{(2.13)}
\def\argao{(2.14)}
\def\argap{(2.15)}
\def\argaq{(2.16)}
\def\argar{(2.17)}
\def\argat{(2.19)}
\def\argau{(2.20)}

\def\argav{(3.1)}
\def\argaw{(3.2)}
\def\argay{(3.4)}
\def\argba{(3.6)}
\def\argbc{(3.8)}
\def\aaargbg{(3.12)}

\def\argbf{(4.1)}
\def\argbg{(4.2)}
\def\argbi{(4.4)}
\def\argbj{(4.5)}
\def\argbk{(4.6)}
\def\argbl{(4.7)}
\def\argbm{(4.8)}
\def\argbn{(4.9)}
\def\argbo{(4.10)}
\def\argbq{(4.12)}
\def\argbr{(4.13)}
\def\argbs{(4.14)}

\def\aargbu{(5.2)}
\def\argbv{(5.3)}
\def\argbw{(5.4)}
\def\argbx{(5.5)}
\def\argby{(5.6)}
\def\aargby{(5.7)}
\def\argbz{(5.8)}
\def\argcb{(5.10)}
\def\argcc{(5.11)}
\def\argcd{(5.12)}
\def\argce{(5.13)}
\def\argcf{(5.14)}
\def\argcg{(5.15)}
\def\argch{(5.16)}
\def\argci{(5.17)}

\def\argcj{(6.1)}
\def\argck{(6.2)}
\def\argcl{(6.3)}

\def\argcn{(7.1)}
\def\argco{(7.2)}
\def\argcp{(7.3)}
\def\argcq{(7.4)}
\def\argcr{(7.5)}
\def\argcs{(7.6)}
\def\argct{(7.7)}
\def\argcw{(7.10)}
\def\argcy{(7.12)}
\def\argda{(7.14)}
\def\argdc{(7.16)}

\def\argde{(8.1)}
\def\argdf{(8.2)}
\def\argdg{(8.3)}
\def\argdh{(8.4)}
\def\argdi{(8.5)}
\def\argdj{(8.6)}
\def\argdk{(8.7)}
\def\argdl{(8.8)}
\def\1argdm{(8.9)}
\def\argdm{(8.10)}
\def\argdn{(8.11)}
\def\argdo{(8.12)}
\def\argdr{(8.15)}
\def\argds{(8.16)}

\def\argdt{(9.1)}
\def\argdu{(9.2)}
\def\argdx{(9.5)}
\def\argdz{(9.7)}
\def\argea{(9.8)}
\def\argeb{(9.9)}
\def\argec{(9.10)}

\def\argef{(10.1)}
\def\argeg{(10.2)}
\def\argeh{(10.3)}

\def\argei{(11.1)}
\def\argek{(11.3)}
\def\aargek{(11.4)}
\def\3argek{(11.5)}
\def\argel{(11.6)}
\def\aargel{(11.7)}
\def\3argel{(11.8)}
\def\4argel{(11.9)}
\def\5argel{(11.10)}
\def\6argel{(11.11)}
\def\7argel{(11.12)}
\def\8argel{(11.13)}
\def\argem{(11.14)}
\def\argen{(11.15)}
\def\argeo{(11.16)}
\def\argeq{(11.18)}

\def\argfb{(12.8)}
\def\argfd{(12.10)}
\def\argfe{(12.11)}

\def\argfp{(13.2)}
\def\argfr{(13.4)}
\def\argfu{(13.7)}
\def\aargfu{(13.8)}
\def\argfv{(13.9)}
\def\argfw{(13.10)}
\def\argfx{(13.11)}
\def\argfz{(13.13)}
\def\arggd{(13.17)}
\def\argge{(13.18)}
\def\arggf{(13.19)}
\def\arggi{(13.22)}
\def\arggj{(13.23)}
\def\arggl{(13.25)}
\def\arggm{(13.26)}
\def\arggo{(13.28)}
\def\arggp{(13.29)}
\def\arggq{(13.30)}
\def\arggr{(13.31)}

\def\arggu{(14.2)}
\def\aaarggu{(14.3)}
\def\aaarggv{(14.4)}
\def\aaarggx{(14.6)}
\def\arggv{(14.7)}
\def\arggw{(14.8)}
\def\arggx{(14.9)}
\def\arggy{(14.10)}
\def\arggz{(14.11)}
\def\argha{(14.12)}

\def\argfg{(\uppercase\expandafter{\romannumeral 1}.{\romannumeral 1})}
\def\argfh{(\uppercase\expandafter{\romannumeral 1}.{\romannumeral 2})}
\def\argfi{(\uppercase\expandafter{\romannumeral 1}.{\romannumeral 3})}
\def\argfj{(\uppercase\expandafter{\romannumeral 1}.{\romannumeral 4})}
\def\argfk{(\uppercase\expandafter{\romannumeral 1}.{\romannumeral 5})}
\def\argfl{(\uppercase\expandafter{\romannumeral 1}.{\romannumeral 6})}
\def\argfm{(\uppercase\expandafter{\romannumeral 1}.{\romannumeral 7})}
\def\argfn{(\uppercase\expandafter{\romannumeral 1}.{\romannumeral 8})}

\vskip 1cm
\pageno=0
\centerline {\para Les quartiers de la lune de Troie}
\bigskip

\centerline {\para ou}
\bigskip
\centerline {\para Le bouclier d'Achille}

\bigskip

\centerline{\soustitre Essai sur le corps de classe, par}
\medskip

\centerline{Jacques BoŽchat et Maurice Mischler}

\vfill\eject
\pageno=-1
\centerline {\para Introduction}
\bigskip
La thŽorie du corps de classe occupe une place privilŽgiŽe dans les mathŽmatiques. En effet, cet ensemble de thŽormes est ˆ la base de plusieurs pans de ce qui se fait de plus pointu actuellement en thŽorie des nombres. Mais notre objectif Žtait plus prŽcis : nous avons, il y a deux ans, rŽdigŽ une preuve du thŽorme de Catalan-Mih$\breve{ \rm a}$ilescu. Et dans cette preuve, nous avions besoin de l'existence du corps de Hilbert d'un corps de nombres et du ThŽorme de $\check{\rm C}$ebotarev. On voit rapidement  qu'il faut pour cela une bonne partie des thŽormes principaux de la thŽorie du corps de classe global. Nous avons donc fait un sŽminaire sur le sujet pour comprendre cette thŽorie. 

Le but est donc de prouver de la manire la plus directe ces deux rŽsultats. Plusieurs approches sont possibles : l'approche adŽlique, l'approche cohomologique et l'approche classique. L'approche classique consisterait ˆ relire dans le texte les oeuvres de  $\check{\rm C}$ebotarev et de Takagi. Mais il serait dommage d'oublier ce qui s'est fait par la suite visant ˆ amŽliorer la comprŽhension profonde, notamment l'application d'Artin et le quotient de Herbrand. En revanche, les approches adŽliques et cohomologiques nous ont paru un peu ŽloignŽes du problme initial. Restait donc une ligne un peu mŽdiane utilisant la cohomologie, mais uniquement cyclique et en ``snobant'' les adles. NŽanmoins, rongŽ par le remords, et voyant que la thŽorie vu avec les adles pouvait se dŽduire assez facilement de ce que nous avions dŽjˆ fait, nous l'avons mis tout de mme au Chapitre 13.

Pour pouvoir lire cet exposŽ avec aisance, il serait prŽfŽrable d'avoir suivi au moins un cours de thŽorie ŽlŽmentaire des nombres. C'est-ˆ-dire conna"tre la notion de corps de nombres; la thŽorie de Galois sur ceux-ci, sur les corps finis et sur les corps $\P$ -adiques; le thŽorme de Dirichlet sur les unitŽs; le thŽorme ``$n=efr$'' sur les extensions d'idŽaux premiers dans une extension galoisienne de corps de nombres; les notions de groupes de dŽcomposition et d'inertie; les normes absolues et relatives d'ŽlŽments et d'idŽaux ainsi que quelques rŽsultats basiques de l'analyse rŽelle et complexe. Ces rŽsultats seront tout de mme rappelŽs, simplement pour parler le mme langage avec le lecteur.

Nous espŽrons que le lecteur prendra plaisir ˆ parcourir (ou ˆ Žtudier) ce texte et que cette thŽorie cessera d'effrayer les gens, car il est vrai qu'elle est un peu dure pour des novices et trop standard pour des mathŽmaticiens actifs, donc peu de gens prennent la peine de regarder en dŽtail tout cela et c'est bien dommage !

Bien sžr, le point de vue que nous prŽsentons ici est largement inspirŽ de divers ouvrages ou articles. Notamment [Jan], [La2] et [Neu]  mais l'approche est un peu diffŽrente, un peu plus directe, plusieurs petites erreurs ont ŽtŽ corrigŽes (on espre ne pas en avoir ajoutŽes) et surtout quelques dŽveloppements du genre ``left to the reader'' Žclaircis et un peu ŽtoffŽs.

Voyons un peu la structure de notre texte~:

Il faut considŽrer le Chapitre 0 comme une bo"te ˆ malice dans laquelle se trouvent les rŽsultats importants sur lesquelles la thŽorie que nous prŽsenterons sera construite. On peut aussi le considŽrer comme ce qu'on met dans notre sac ˆ dos avant de partir faire une excursion en montagne~: il y a un peu de tout et c'est un peu comprimŽ car le sac est toujours trop petit !

Le Chapitre 1 est un chapitre d'Žchauffement sur un rŽsultat technique qui n'est utilisŽ qu'une fois dans le Chapitre 2. On a hŽsitŽ de mettre tout cela en appendice, mais personne ne lit les appendices et en plus, il y a tout de mme certains raisonnements qui seront revus par la suite.

Dans le Chapitre 2, on dŽmontre de manire analytique ce qu'on appelle la premire inŽgalitŽ du corps de classe. On utilisera un peu d'analyse complexe, mais ˆ un niveau assez basique, il faut essentiellement conna"tre la notion de fonctions holomorphes et mŽromorphes. On montre au passage que l'application d'Artin (vue au Chapitre 0) est surjective

Dans le Chapitre 3, on prouve le ThŽorme de $\check{\rm C}$ebotarev faible et quelques rŽsultats comme le thŽorme de Dirichlet sur les progressions arithmŽtiques ainsi que des rŽponses sur le comportement modulo $p$ de polyn™mes irrŽductibles dans $\Z[X]$. La fin du chapitre donne de jolis rŽsultats sur les diffŽrentes manires de dŽcomposer pour un idŽal premier (par exemple, le ThŽorme de Bauer).

Le Chapitre 4 parle de cohomologie cyclique (depuis le dŽbut), du quotient de Herbrand et on donne quelques calculs dans le cas d'extensions cycliques de corps de nombres.

Le Chapitre 5 est difficilement dŽfinissable~: on calcule essentiellement l'indice $[K^*:N(L^*)K^*_{\m}]$, mais \c ca ne vous dira pas grand-chose. En revanche, nous faisons une digression permettant de dŽfinir l'exponentielle et le logarithmes sur les corps $\P$-adiques et une autre digression pour savoir quand un ŽlŽment est une norme d'une extension finie dans les $\P$-adiques.

Pour le Chapitre 6, les calculs faits aux chapitres 4 et 5 permettent, avec l'Žtude approfondie d'un diagramme du tonnerre, de prouver l'ŽgalitŽ fondamentale du corps de classes pour les extensions cycliques. Cela implique un thŽorme connu sous le nom de ``ThŽorme de la Norme de Hasse''. 

Avec le Chapitre 7, nous entrons dans le monde des extensions abŽliennes avec la notion de ``$K$-modules admissibles''. Cela nous permet de dŽmontrer le grand ``ThŽorme de rŽciprocitŽ d'Artin''. En corollaire, on prouve le ThŽorme de $\check{\rm C}$ebotarev fort et le fameux ``ThŽorme de Kronecker-Weber''.

C'est dans le Chapitre 8, que nous apercevons la lumire~: on y dŽfinit la notion de sous-groupe de congruence, ce qui nous permet d'Žnoncer le ``ThŽorme d'existence du corps de classe''. Ce rŽsultat, fondamental, nous occupera jusqu'au Chapitre 10. Nous faisons quelques rŽductions pour pouvoir attaquer le problme dans des cas plus faciles.

Dans le Chapitre 9, nous nous concentrons sur un cas ``plus simple''~: les extensions de Kummer. Nous calculons un nouvel indice. Et comme interlude, nous prouvons le cŽlbre thŽorme de rŽciprocitŽ quadratique, tout en Žtant conscient que cette preuve est probablement la plus compliquŽe de ce rŽsultat classique.

Au Chapitre 10, nous prouvons le thŽorme d'existence, nous dŽfinissons aussi l'application d'Artin dans le cas des extensions galoisiennes non abŽliennes. Nous utiliserons cette application au Chapitre 12.

Le Chapitre 11 Žtait initialement consacrŽ ˆ la construction du corps de Hilbert. Pour cela, il faut montrer que Òle conducteur est admissible". Pour parvenir ˆ ce rŽsultat, nous dŽfinissons le "symbole de norme rŽsiduelle", notŽ $\theta_\P$. Ce symbole nous permettra en outre de prouver au Chapitre 14 les thŽormes du corps de classe local. Nous prouvons donc en plus un certain nombre de rŽsultats pas directement utiles pour l'admissibilitŽ du conducteur. Enfin, nous pouvons construire le corps de Hilbert.

Tout idŽal d'un corps de nombres devient principal dans son corps de Hilbert. Cette propriŽtŽ remarquable est dŽmontrŽe dans le Chapitre 12. On sort un peu des chemins arithmŽtiques pour faire une petite incursion dans la thŽorie des groupes, pour bien sžr revenir au rŽsultat qui nous intŽresse.

Les notions d'idle et d'adle ont ŽtŽ introduit pour formuler la thŽorie du corps de classe pour les extensions infinies. Nous avons donc dŽcidŽ au Chapitre 13 de Òrecoller les morceaux" pour des lecteurs voulant peut-tre se lancer dans cette voie.

Enfin, voyant qu'il ne fallait plus trop travailler pour obtenir les rŽsultats du corps de classe local (en caractŽristique 0), nous avons introduit ce dernier chapitre pour cueillir encore ces rŽsultats, non sans s'tre assurŽ que tout corps local est bien le localisŽ d'un corps de nombres.

Le premier appendice parle des nombres premiers de la forme $x^2+n y^2$  ou $x^2+xy+m y^2$ et on y prouve que sous certaines conditions, il y a une infinitŽ de tels nombres premiers. On utilise pour la preuve l'existence du corps de la classe d'un groupe qu'on nomme $H_{\O}$ dans des corps quadratiques.

Le sujet du second appendice est le symbole de Hilbert. Nous en donnons les premires dŽfinitions et les premires propriŽtŽs. 
\medskip
Enfin, on pourrait nous ajouter qu'il serait judicieux de parler de la version cohomologique des thŽormes du corps de classe introduite par J. Tate. Nous y avons pensŽ, mais il nous semble que cela alourdirait notre propos. Mais peut-tre qu'un jour, nous apporterons un appendice ˆ ce texte quand nous aurons trouvŽ une manire ŽlŽgante de prŽsenter cela.

On espre bien sžr que chaque lecteur apprendra quelque chose dans ce texte et qu'il n'hŽsitera pas ˆ nous signaler des erreurs, des maladresses  ou des suggestions d'amŽliorations. Toute remarque est ˆ envoyer ˆ l'adresse :

\centerline{maurice.mischler@romandie.com ou maurice.mischler@vd.educanet2.ch}

Par souci de ne pas trop charger le ficher informatique, nous avons enlevŽ un certain nombre d'illustrations non mathŽmatiques. Vous pouvez les voir sur le site

\centerline{http://mathmontmus.romandie.com/resource/12252/262704}

\vskip1cm






\vfill\eject

\parindent 0pt

\NoBlackBoxes
\centerline{\para Table des matires}

\vskip 2cm

{\bf Chapitre 0 :  Rappels et premiers exemples} \dotfill 1
\bigskip
{\bf Chapitre 1 : Un rŽsultat sur $j(x,{\euf K})$}\dotfill 19
\bigskip
{\bf Chapitre 2 : SŽries de Dirichlet et premire inŽgalitŽ du corps de classe}\dotfill 26
\bigskip
{\bf Chapitre 3 : ThŽorme de $\check{\rm C}$ebotarev} \dotfill 39
\bigskip
{\bf Chapitre 4 :   Cohomologie des groupes cycliques et quotient de Herbrand}\dotfill 50
\bigskip
{\bf Chapitre 5 :  Un calcul d'indice}\dotfill 60
\bigskip
{\bf Chapitre 6 :   L'ŽgalitŽ fondamentale du corps de  classe pour les extensions cycliques et 
 thŽorme de la norme de Hasse }\dotfill 70
\bigskip
{\bf Chapitre 7 : La loi de rŽciprocitŽ d'Artin }\dotfill 77
\bigskip
{\bf Chapitre 8 : PrŽparation ˆ la formation des classes }\dotfill 85
\bigskip
{\bf Chapitre 9 :  Quelques rŽsultats sur la thŽorie des $n$-extensions de Kummer,  et calcul d'un nouvel indice }\dotfill 94
\bigskip
{\bf Chapitre 10 : Le thŽorme principal du corps de classe }\dotfill 103
\bigskip
{\bf Chapitre 11 :  Symbole de restes normiques, conducteur  et corps de classe de Hilbert }\dotfill 111
\bigskip
{\bf Chapitre 12 :  Capitulation des idŽaux d'un corps nombres dans son corps de Hilbert }\dotfill 125
\bigskip
{\bf Chapitre 13 :  InterprŽtation idŽlique}\dotfill 135
\bigskip
{\bf Chapitre 14 :  Corps de classe local}\dotfill 152
\bigskip
{\bf Appendice 1 : Deux mots sur les corps quadratiques  et sur des reprŽsentations de nombres premiers}\dotfill 161
\bigskip
{\bf Appendice 2 : Deux mots sur le symbole de Hilbert }\dotfill 166
\bigskip
{\bf Glossaire et symboles }\dotfill 169
\bigskip
{\bf Index }\dotfill 172
\bigskip
{\bf Bibliographie }\dotfill 175

\vfill\eject

\vfill\eject

\bigskip
\centerline{\para Chapitre 0 :}
\medskip
\centerline{\para Rappels et premiers exemples }
\bigskip
\parindent 0.5cm
\pageno=1
Le dŽbut de ce premier chapitre est une espce de mise en commun des outils et rŽsultats mathŽmatiques qui sont dits Òbien connus". Le plus simple serait de dire~: Òon considre connus les rŽsultats de [La1], [Sam], [Mar], [Fr-Tay] et [Nar]Ó. Mais c'est un peu court et pas trs gentil pour le lecteur. La volontŽ est aussi de se mettre d'accord sur les notations. Le lecteur connaisseur en thŽorie algŽbrique des nombres pourra sauter les premires considŽrations jusqu'au paragraphe intitulŽ ``$K$-modules, application d'Artin et compagnie''. Ce n'est qu'ˆ partir de lˆ que nous donnerons toutes les preuves (ou au moins les rŽfŽrences).
\bigskip
\centerline{\soustitre GŽnŽralitŽ sur les corps de nombres et l'algbre ŽlŽmentaire}
\bigskip
Les ensembles $\N=\{0,1,2,3,\ldots,\}$, $\Z$, $\Q$, $\R$ et $\C$ seront supposŽe connus. Les notions de groupes, anneaux, idŽaux, corps, modules, espaces vectoriels, algbres,... aussi. Si $A$ est un anneau. On note $A^*$ ou $U(A)$ l'ensemble des inversibles de $A$.

On rappelle la donnŽe des thŽormes d'isomorphismes. Chaque fois que nous dirons ``par les thŽormes d'isomorphismes...'' nous nous rŽfŽrerons ˆ ceci~:
\medskip
{\soustitre ThŽormes d'isomorphismes :}
\medskip
{\sl \art{a)}Soit $G$, $G'$ des groupes et $f\, :\, G\to G'$ un homomorphisme de groupe de noyau $H$. Alors $f$ induit un isomorphisme $f'\,:\,  G/H\to {\rm im} (f)$ qui factorise $f$ en la suite d'homomorphisme
$$G\buildrel p\over\lra G/H\buildrel f'\over\lra {\rm im}(f)\buildrel i\over \lra G',$$
o $p$ et $i$ sont les projections, respectivement les injections canoniques. Plus gŽnŽralement, si $H'$ est un sous-groupe normal de $G'$ et $H=f^{-1}(H')$, on en dŽduit un homomorphisme injectif~:
$$\overline{f}\, :\, G/H\lra G'/H'$$
qui est un isomorphisme si $f$ est surjectif.

\art{b)}Soient $G$ un groupe, et $H_{1}\supset  H_{2}$ des sous-groupes normaux de $G$ ($H_{2}$ est alors automatiquement normal dans $H_{1}$). Alors on a un isomorphisme
$$(G/H_{2})/(H_{1}/H_{2})\simeq G/H_{1}.$$

\art{c)}Soient $H_{1}, H_{2}$ des sous-groupes d'un groupe $G$. Supposons que $H_{1}\subset \{x\in G\mid xH_{2}x^{-1}=H_{2}\}$, $H_{1}\cap H_{2}$ est alors un sous-groupe normal de $H_{1}$ et $H_{1}H_{2}=H_{2}H_{1}$ est un sous-groupe de $G$ dans lequel $H_{2}$ est normal. Alors on a~:
$$H_{1}/(H_{1}\cap H_{2})\simeq (H_{1}H_{2})/H_{2}.$$

\qed

}

\bigskip

Un {\it corps de nombres} est un corps de dimension finie vu comme $\Q$-espace vectoriel. Il sera souvent notŽ $K,L,E,M$ ou $H$ et cette dimension se note $[K:\Q]$, pour un corps de nombres $K$ (attention, si $G\supset H$ sont des groupes, $[G:H]$ sera $|G/H|$, le contexte permettra de distinguer). Pour simplifier, nous considŽrerons que tous ces corps sont inclus dans le mme corps algŽbriquement clos, disons $\C$. On peut montrer que pour tout corps de nombres $K$, il existe $\theta\in\C$  tel que $K=\Q(\theta)$. Si $L\subset K\subset \Q$ sont des corps de nombres, alors on a $[L:\Q]=[L:K]\cdot [K:\Q]$ (le mme rŽsultat est aussi vrai pour des corps quelconques par exemple des corps finis). 

Si $K$ est un corps de nombres, on note souvent $O_K$ l'anneau des ŽlŽments de $K$ entiers sur $\Q$. On suppose connu que $O_K$ est un anneau de {\it Dedekind}, c'est-ˆ-dire noethŽrien, intŽgralement clos (donc intgre) et tout idŽal premier non nul de $O_K$ est maximal. On peut montrer que dans ce cas-lˆ, l'ensemble $I_{K}$ des idŽaux fractionnaires de $K$ est un groupe abŽlien librement engendrŽ par les idŽaux premiers de $O_{K}$. Remarquons que souvent on dira idŽal de $K$ plut™t que de $O_{K}$. La graphie est souvent $\aa$ ou $\bb$ pour des idŽaux  et $\P$ et $\gP$ pour les idŽaux premiers ou les places infinies. Nous reviendrons plus tard sur la notions de places. Les unitŽs de $O_{K}$ devraient se noter $O_{K}^*$ ou $U(O_{K})$... mais nous les noterons $U_{K}$.

\newcount\gaaga \gaaga=\gaga

Soit $m\in \N$, $m>1$. On note $\zeta_{m}$ une racine $m$-me primitive de l'unitŽ. Le corps $K=\Q(\zeta_{m})$ est appelŽ le {\it $m$-ime corps cyclotomique}. On sait que $[\Q(\zeta_{m}):\Q]=\varphi(m)$, o $\varphi$ est l'indicateur d'Euler. Il est bien connu, mais toujours assez dŽlicat ˆ prouver que dans ce cas $O_{K}=\Z[\zeta_{m}]$.

Soit $L/K$ une extension de corps de nombres de degrŽ $n$. Soit $\P$ un idŽal premier de $K$ et $\gP$ un idŽal premier de $L$. On dit que $\gP$ est {\it au-dessus} de $\P$ si $\P=\gP\cap O_{K}$, ou ce qui est Žquivalent, $\gP$ appara"t dans la dŽcomposition de l'idŽal $\P O_{L}$ et on Žcrit dans ce cas $\gP |\P$. On peut alors identifier $O_{K}/\P$ ˆ un sous-corps de $O_{L}/\gP$ (ces corps sont d'ailleurs finis). On notera $f(\gP/\P)=[O_{L}/\gP:O_{K}/\P]$, qu'on appelle le {\it degrŽ rŽsiduel de $\gP/\P$}.  Si on Žcrit $\P O_{L}=\gP_{1}^{e_{1}}\cdots \gP_{r}^{e_{r}}$ et notant $f_{i}$ pour $f(\gP_{i}/\P)$, on a $\sum_{i=1}^re_{i}f_{i}=n$.  Le nombre entier $e_{i}$, notŽ $e(\gP_{i}/\P)$ s'appelle {\it l'indice de ramification de $\gP_{i}/\P$}. On dit que $\gP$ n'est pas ramifiŽ dans $L$ (ou ne se ramifie pas dans $L$) si $e_{i}=1$ pour tout $i$. On peut montrer que le nombre de  $\P$ qui sont ramifiŽ est fini. Par exemple, si $K=\Q$ et $L=\Q(\zeta_{m})$, et si $p$ est un nombre premier, alors l'idŽal $(p)$ ramifie dans $L$ si et seulement si $p$ divise $m$.
Enfin, si $K\subset L\subset E$ sont des corps de nombres, et si $\P\subset \gP\subset\ggP$ sont des idŽaux premiers des $K,L,E$ respectivement, on a $e(\ggP/\P)=e(\ggP/\gP)\cdot e(\gP/\P)$; et il en est de mme pour les $f$.

\headline={\hfill \smcap Rappels et premiers exemples \hfill}

\bigskip

\centerline{\soustitre Extensions galoisiennes}
\bigskip
\newcount\gaagb \gaagb=\gaga

Soit $L/K$ une extension algŽbrique de corps d'indice fini $n$. On dit qu'elle est {\it galoisienne} si $|{\rm Aut}_{K}(L)|=[L:K]=n$, o ${\rm Aut}_{K}(L)$ est l'ensemble des $K$-automorphismes de $L$. Dans ce cas, ${\rm Aut}_{K}(L)$ se note ${\rm Gal}(L/K)$ (le groupe de Galois de $L/K$). Si $L$ et $K$ sont des corps finis, alors $L/K$ est une extension galoisienne, mieux, elle est cyclique (le groupe de Galois est un groupe cyclique) engendrŽ par {\it l'automorphisme de Frobenius} $x\mapsto x^q$, o $q=|K|$. Supposons maintenant que $L$ et $K$ soient des corps de nombres. On peut voir que tout $\sigma\in {\rm Gal}(L/K)$ agit transitivement sur les idŽaux premiers de $L$ qui sont au-dessus d'un idŽal $\P$ de $K$ fixŽ. Cela implique que si $\P O_{L}=\gP_{1}^{e_{1}}\cdots \gP_{r}^{e_{r}}$, on a $e_{1}=e_{2}=\cdots =e_{r}=:e$ et $f_{1}=f_{2}=\cdots =f_{r}=:f$, et ainsi, $e\cdot f\cdot r=n$. Quand nous dirons ``la thŽorie de Galois implique que...'', nous ferons rŽfŽrence ˆ un des rŽsultats suivants~:
\medskip\goodbreak
{\soustitre ThŽormes de Galois}
\medskip
Si $K$ et $L$ sont des corps, on note $KL$ le plus petit corps contenant
$K$ et $L$. 

{\sl 
\art{a)}Soit $K\subset L\subset E$ des corps. Supposons que $E/K$ soit une extension
galoisienne de groupe $G$. Alors $E/L$ est galoisienne et ${\rm Gal}(E/L)=\{g\in G\mid g|_L={\rm
Id}_L\}:=H$. De plus $L/K$ est galoisienne si et seulement si $H$ est un sous-groupe normal de
$G$ et dans ce cas, ${\rm Gal}(L/K)\simeq G/H$. Inversement, si $H_{1}$ est un sous-groupe de $G$, le corps $L_{1}:={\rm Fix}(H_{1})=\{x\in E\mid h(x)=x\ \forall h\in H_{1}\}$, qu'on appelle {\it le corps fixe par $H_{1}$} est tel que ${\rm Gal}(E/L_{1})=H_{1}$.

\art{b)}Soit $K\subset L$ et $K\subset E$ deux extensions de corps. On suppose que
$L/K$ est galoisienne. Alors $EL/L$ est aussi galoisienne et ${\rm Gal}(EL/E)\simeq {\rm
Gal}(L/L\cap E)\subset{\rm Gal}(L/K)$. Cet isomorphisme est donnŽ par la restriction ˆ $L$. L'application 
$R\, : {\rm Gal}(EL/E)\to {\rm Gal}(L/K)$ ainsi dŽfinie est un homomorphisme injectif.

\art{c)}Soit $K\subset L$ et $K\subset E$ deux extensions galoisiennes de corps
telles que $L\cap E=K$. Alors $EL/K$ est une extension galoisienne et ${\rm Gal}(EL/K)\simeq
{\rm Gal}(L/K)\times {\rm Gal}(E/K) $.\qed

} 

\bigskip

\centerline{\soustitre Normes}

\bigskip

\art{a)}{\bf La norme absolue}

\newcount\gaagc \gaagc=\gaga

\art{}Soit $K$ un corps de nombres et $\aa$ un idŽal de $K$. Alors l'anneau $O_{K}/\aa$ est fini et son cardinal se note $\N(\aa)$, {\it la norme absolue de $\aa$}. On peut voir que $\N(\aa\cdot\bb)=\N(\aa)\cdot\N(\bb)$, pour tout idŽal $\aa$ et $\bb$. On prolonge la dŽfinition pour tout idŽal fractionnaire de $K$.

\art{b)}{\bf La norme relative d'un ŽlŽment}

\art{}Soit $L/K$ une extension de corps de nombres et $\alpha\in L$. L'application $\mu_{\alpha}\, :\, L\to L$ dŽfinie par $\mu_{\alpha}(\beta)=\alpha\cdot \beta$ est un endomorphisme $K$-linŽaire de $L$. On pose $N_{L/K}(\alpha)=\det(\mu_{\alpha})$. Il est clair que si $\alpha\in K$, $N_{L/K}(\alpha)=\alpha^n$. On peut voir aussi que $N_{L/K}(\alpha\cdot\beta)=N_{L/K}(\alpha)\cdot N_{L/K}(\beta)$, pour tout $\alpha,\beta\in L$. Si $K\subset L\subset E$ sont des corps de nombres et $\alpha\in E$, on a $N_{E/K}(\alpha)=N_{L/K}(N_{E/L}(\alpha))$. De plus, si $\alpha\in K^*$, on a $\N(\alpha \cdot O_{K})=\left | N_{K/\Q}(\alpha)\right |$. Si $\sigma_{1},\ldots ,\sigma_{n}$ sont les $K$-morphismes de $L$ dans $\C$, on a $N_{L/K}(\alpha)=\prod_{i=1}^n\sigma_{i}(\alpha)$. En particulier si $L/K$ est galoisienne de groupe $G$, on a $N_{L/K}(\alpha)=\prod_{\sigma\in G}\sigma(\alpha)$. 

\art{c)}{\bf La norme relative d'un idŽal}

\art{}Soit $L/K$ une extension de corps de nombres d'indice $n$. Si $\P$ et $\gP$ sont des idŽaux premiers de $K$ et $L$ respectivement tels que $\gP|\P$, on pose $N_{L/K}(\gP)=\P^{f(\gP/\P)}$, et on prolonge multiplicativement cette norme ˆ tous les idŽaux fractionnaires de $L$ (puisque les idŽaux premiers engendrent $I_K$). On voit immŽdiatement que si $\aa$ est un idŽal fractionnaire de $K$, on a $N_{L/K}(\aa)=\aa^n$. De mme, pour tout $a\in L$, on $N_{L/K}(a\cdot O_{L})=N_{L/K}(a)\cdot O_{K}$, o $N_{L/K}$ est la norme relative dŽfinie prŽcŽdemment. Enfin, si l'extension $L/K$ est galoisienne de groupe $G$, et $\aa$ est un idŽal fractionnaire de $L$, alors $N_{L/K}(\aa)=\prod_{\sigma\in G}\sigma(\aa)$.
\bigskip
\centerline{\soustitre Ramification et dŽcomposition et automorphisme de Frobenius}
\bigskip
Soit $L/K$ une extension de corps de nombres et $\P$ un idŽal premier de $K$. On dit que $\P$ {\it se dŽcompose totalement dans $L$} si $e(\gP/\P)=f(\gP/\P)=1$ pour tout idŽal premier $\gP$ de $L$ tel que $\gP|\P$. On a le rŽsultat suivant~:
\medskip

{\soustitre Lemme ``dŽcomposition-ramification''}
\medskip
{\sl 

\art{a)}Soit $L_1/K$ et $L_2/K$ deux extensions de corps de nombres. Alors l'ensemble des idŽaux premiers de
$K$ qui se dŽcomposent compltement (resp. qui ne ramifient pas) dans $L_1L_2$ est l'ensemble de idŽaux premiers de
$K$ qui se dŽcomposent compltement (resp. qui ne ramifient pas) dans $L_1$ et dans $L_2$. 

\art{b)}Soit $L/K$ une extension de corps de nombres et $E/K$ la plus petite extension galoisienne contenant $L$.
Alors l'ensemble des idŽaux premiers de $K$ qui se dŽcomposent compltement (resp. qui ne ramifient pas) dans $L$ est
l'ensemble de idŽaux premiers de $K$ qui se dŽcomposent compltement (resp. qui ne ramifient pas) dans $E$.

}
\medskip
Cf. [Mar, Thm. 31 + Corollary, pp. 107-108]
\bigskip
Soit $L/K$ une extension galoisienne de corps de nombres de groupe $G$, $\P$ un idŽal premier de $K$ et $\gP$ un idŽal premier de $L$ tel que $\gP |\P$. On dŽfinit le {\it groupe de dŽcomposition de $\gP$ sur $\P$ (ou de $\gP$ sur $K$)},
 $$Z(\gP/\P)=Z(\gP/K):=\{\sigma\in G\mid \sigma(\gP)=\gP\}.$$
\newcount\gaagd \gaagd=\gaga

Si $\gP_{1},\gP_{2} |\P$, alors il existe $\sigma\in G$ tel que $Z(\gP_{1}/\P)=\sigma^{-1}Z(\gP_{2}/\P)\sigma$. Ainsi, si l'extension est abŽlienne ($G$ est abŽlien), alors tous les $Z(\gP/\P)$ sont Žgaux si $\P$ est fixŽ, et on note alors ce groupe $Z(L/\P)$, ou mme $Z(\P)$ s'il n'y a pas d'ambig¬\"uitŽ, et on l'appelle le {\it groupe de dŽcomposition de $\P$ sur $L$}. Revenons au cas gŽnŽral ($G$ non nŽcessairement abŽlien); nous avons vu que $[L:K]=n=e\cdot f\cdot r$. On peut aussi voir que $[G: Z(\gP/\P)]=r$ et $|Z(\gP/\P)|=e\cdot f$. Si $\sigma\in Z(\gP/\P)$, alors il dŽtermine un $\overline{\sigma}\in \overline{G}:={\rm Gal}((O_{L}/\gP)/(O_{K}/\P))$, et l'application $\sigma\mapsto \overline{\sigma}$ est un homomorphisme surjectif de $G$ sur $\overline{G}$. Le noyau de cette application se note $T(\gP/\P)$ ou $T(\gP/K)$ et s'appelle le {\it groupe d'inertie} de $\gP/\P$ ou de $\gP/K$. On a aussi, comme pour $Z$, $T(\sigma(\gP)/\P))=\sigma T(\gP/\P)\sigma^{-1}$ pour tout $\sigma\in G$. On a donc $|T(\gP/\P)|=e$ et
$Z(\gP/\P)/T(\gP/\P)\simeq \overline{G}$, de cardinal $f$. On a donc, pour tout $\sigma\in
T(\gP/\P)$, $\sigma(x)\equiv x\pmod \gP$ pour tout $x\in O_L$. Supposons que $\P$ ne ramifie
pas, c'est-ˆ-dire $e=1$ et donc le groupe d'inertie est trivial et donc l'application $\sigma\mapsto\overline{\sigma}$
est un isomorphisme de $Z(\gP/\P)$ sur $\overline{G}$. Nous avons vu que le groupe de Galois $\overline{G}$ est un groupe cyclique avec un gŽnŽrateur privilŽgiŽ qui est l'application $\nu\mapsto \nu^{\N(\P)}$ pour tout $\nu\in O_L/\gP$ et l'unique ŽlŽment de $Z(\gP/\P)$ qui correspond ˆ cet automorphisme s'appelle aussi {\it l'automorphisme de Frobenius de $\gP/\P$}. On le note ${\rm Frob}(\gP/\P)$. Il est caractŽrisŽ
comme l'ŽlŽment de $G$ qui satisfait~:

\newcount\gaage \gaage=\gaga

$${\rm Frob}(\gP/\P)(x)\equiv x^{\N(\P)}\pmod \gP\quad \hbox{ pour tout }x\in O_L.$$

On a aussi ${\rm Frob}(\sigma(\gP)/\P)=\sigma\ {\rm Frob}(\gP/\P)\ \sigma^{-1}$ et donc l'ensemble $\{
{\rm Frob}(\gP/\P)\mid
\gP |\P\}$, qu'on note ${\rm Fr}_{L/K}(\P)$ est une classe de conjugaison dans $G$. Si $G$ est
abŽlien, alors ${\rm Frob}(\gP/\P)$ ne dŽpend que de $\P$, on le notera ${\rm Frob}_{L/K}(\P)$, et
on a

$${\rm Frob}_{L/K}(\P)(x)\equiv x^{\N(\P)}\pmod {\P O_L}\quad \hbox{ pour tout }x\in O_L.$$

De plus, si $K\subset M\subset L$ sont des corps de nombres et $\P\subset\gP\subset{\bf P}$ sont
des idŽaux premiers de $K$, $M$ et $L$ respectivement. Alors ${\rm Frob}(\gP/\P)={\rm
Frob}({\bf P}/\P)|_{M}$, s'ils sont dŽfinis. De mme, $Z({\bf P}/\gP)\subset Z({\bf P}/\P)$ et $T({\bf P}/\gP)\subset T({\bf P}/\P)$.

Enfin, soient $K\subset L$ et $K\subset E$  deux extensions de corps de nombres telles que $L/K$ soit galoisienne. On sait par la thŽorie de Galois que ${\rm Gal}(LE/E)$ peut tre vu, via la restriction ˆ $L$, comme un sous-groupe de ${\rm Gal}(L/K)$. Mais, il y a mieux~: si $\gP$ est un idŽal premier de $LE$, alors $Z(\gP/E)$ (resp. $T(\gP/E)$) est isomorphe (via la mme restriction) ˆ $Z(\gP\cap L/L\cap E)$ (resp.  ˆ $T(\gP\cap L/L\cap E)$) et peut tre vu comme un sous-groupe de $Z(\gP\cap L/K)$ (resp. de $T(\gP\cap L/K)$).
\bigskip

\centerline{\soustitre Places et complŽtions}

Soit $K$ un corps de nombres. Une valeur absolue est une application $|\ |\, :\, K\to\R$, satisfaisant les conditions, pour tout $x,y\in K$~:

\art{(1)}$|x|\geq 0$ et $|x|=0\Leftrightarrow x=0$,

\art{(2)}$|xy|=|x||y|$,

\art{(3)}$|x+y|\leq |x|+|y|$.

Si on remplace la condition (3) par la condition plus forte

\art{(3)'}$|x+y|\leq \max (|x|,|y|)$,

on dit que la valeur absolue est {\it non archimŽdienne}, et {\it archimŽdienne} sinon.

Deux valeurs absolue $|\  |_{1}$ et $|\ |_{2}$ sont dites Žquivalentes s'il existe $c,d\in\R$ tels que pour tout $x\in K$ on ait $c|x|_{1}\leq |x|_{2}\leq d|x|_{1}$. Si deux valeurs absolues sont Žquivalentes, elles induisent sur $K$ la mme topologie. L'ensemble des classes d'Žquivalences des valeurs absolues de $K$ s'appellent les {\it places} de $K$. Si $K$ est un corps de nombres, nous allons donner l'ensemble de ses places.
\bigskip\goodbreak
\art{a)}{\bf Les places finies} (non archimŽdiennes).

A chaque idŽal premier $\P$ de $K$, on associe une valuation $v_{\P}$ dŽfinie de la manire suivante~: si $\aa$ est un idŽal fractionnaire de $K$, on peut Žcrire de manire unique $\aa=\P^r\cdot\aa'$ o $\P$ ne divise pas $\aa'$ (on note $\P\notdiv\aa'$). Alors on dŽfinit $v_{\P}(\aa)=r$. Si $x\in K$, on pose $v_{\P}(x)=v_{\P}(x\cdot O_{K})$. La valeur absolue associŽe ˆ cette valuation est dŽfinie ainsi~:
$$|x|_{\P}=\N(\P)^{-v_{\P}(x)},\quad x\in K.$$
On peut montrer que cette valeur absolue est non archimŽdienne, que si $\qq\ne\P$, les valeurs absolues $|\  |_{\qq}$ et $|\  |_{\P}$ sont non-Žquivalentes et que toute valeurs absolue non-archimŽdienne sur $K$ est Žquivalente ˆ une de celles-ci. On note $\gfP_{0}(K)$ l'ensemble des places finies.

\art{b)}{\bf Les places infinies} (archimŽdiennes).

Supposons que $[K:\Q]=n$. On sait qu'il existe $n$ plongements $\varphi\, :\, K\to\C$ (c'est-ˆ-dire des $\Q$-homomorphismes (injectifs)).  A chacun de ces plongements on associe une valeur absolue~: 
$$|x|_{\varphi}=|\varphi(x) |$$
o $x\in K$ et $|\  |$ est la norme complexe. Si $\varphi(K)\subset \R$, on dit que $\varphi$ est {\it rŽelle}. Si $\varphi(K)\not\subset \R$, on dit que $\varphi$ est {\it complexe}. Si $\varphi$ est complexe, le conjuguŽ complexe $\overline{\varphi}$ de $\varphi$ et $\varphi$ dŽfinissent la mme place. Autrement, ces valeurs absolues sont non-Žquivalentes. Ainsi, si $n=r+2s$, o $r$ est le nombre de plongements rŽels de $K$ et $2s$ est le nombre de plongements complexes, on a en tout $r+s$ places infinies.

Si $K$ est un corps de nombres, on peut montrer qu'il n'y a pas d'autres places. Donc, en rŽsumŽ, il y a une infinitŽ de places finies (autant que d'idŽaux premiers) et un nombre fini de places infinies...

 On note $\gfP_{\infty}(K)$ (resp $\gfP_{\R}(K)$, $\gfP_{\C}(K)$) l'ensemble des places infinies (resp. rŽelles, complexes) de $K$. Souvent, une place infinie se notera $\P$, comme pour les places finies. 

L'ensemble de toutes les places se note bien sžr $\gfP(K)$.
\newcount\gaagf \gaagf=\gaga

\bigskip

Revenons un instant sur les places finies. Si $\P\in \gfP_{0}(K)$, on note $\bbK_{\P}$ le complŽtŽ topologique (on attribue une limite ˆ chaque suite de Cauchy) relativement ˆ la valeur absolue dŽfinie par $\P$. C'est un corps dit ``local'' sur lequel $v_{\P}$ et $|\ |_{\P}$ se prolongent. On dŽfinit 
$$\okp=\{ x\in \bbK_{\P}\mid v_{\P}(x)\geq 0\}\quad \hbox{ et }\quad \widehat{\P}= \{ x\in \bbK_{\P}\mid v_{\P}(x)> 0\}.$$
$\okp$ est un anneau local d'idŽal maximal $\widehat{\P}$. Cet idŽal est principal, on note souvent $\pi$ un gŽnŽrateur de $\widehat{\P}$ qu'on appelle ``uniformisante''; et tout idŽal de $\okp$ est du type $\pi^k\cdot\okp$. On considre que $K\subset \bbK_{\P}$ et on note $\okpp:=\okp\cap K=\{{\alpha\over\beta}\in K\mid \alpha,\beta\in O_K\hbox{ et }\beta\in O_K\setminus\P\}$, le localisŽ de $O_{K}$ en $\P$. C'est aussi anneau local et son idŽal maximal se note $\widetilde{\P}$, il est aussi principal et chaque idŽal est du type $\pi^k\cdot\okpp$, pour une uniformisante qu'on notera aussi parfois $\pi$ (lorsque nous ne devrons pas utiliser les deux). On a $\P\cdot\okpp=\widetilde{\P}$ et $\P\cdot \okp=\widetilde{\P}\cdot\okp=\widehat{\P}$. On a aussi, pour tout $k\in\N$, $k>0$~:
$$O_{K}/\P^k\simeq\okpp/\widetilde{\P}^k\simeq\okp/\widehat{\P}.$$

Selon l'humeur et le besoin du moment, il aurait aussi ŽtŽ possible de dŽfinir $\okp$ comme la limite du systme projectif $\left \{O_{K}/\P^k,\delta_{k+1,k}\right\}$, o $\delta_{k+1,k}\, :\, O_{K}/\P^{k+1}\to O_{K}/\P^k$ est l'homomorphisme canonique $x\pmod {\P^{k+1}}\mapsto x\pmod{\P^{k}}$. 

Dans  le cas o $K=\Q$, on retrouve bien sžr les nombres $p$-adiques habituels $\Q_{p}$.

Si $L/K$ est une extension galoisienne de corps de nombres, $\P$ et $\gP$ des idŽaux premiers de $K$ et $L$ respectivement tels que $\gP|\P$. Alors $\bbL_{\gP}/\bbK_{\P}$ est aussi une extension galoisienne de groupe de Galois canoniquement isomorphe ˆ $Z(\gP/\P)$. Et donc la norme vaut $N_{\bbL_{\gP}/\bbK_{\P}}(x)=\prod_{\sigma\in Z(\gP/\P)}\sigma(x)$. Les autres normes se dŽfinissent de manire identique.

Si la place $\P$ est infinie, la situation est plus simple~: $\bbK_{\P}=\R$ si la place est rŽelle et $\bbK_{\P}=\C$ sinon.  Nous aurons besoin de considŽrations plus fines sur les complŽtions, mais nous regarderons ces choses au moment o nous en aurons besoin !

\bigskip
\goodbreak

\centerline{\soustitre Quelques thŽormes}
\bigskip

Tout d'abord un thŽorme facile, trs souvent utilisŽ~:
\bigskip
{\soustitre ThŽorme chinois}

{\sl Soit $A$ un anneau commutatif, $\aa$ et $\bb$ des idŽaux de $A$ copremiers (i.e. $\aa+\bb=A$). Alors on a l'isomorphisme~:

$$A/(\aa\bb)\simeq A/\aa\times A/\bb.$$

Cf. [Sam, Lemme 1, \S 1.3, p. 22]\qed 

}
\bigskip
Un autre thŽorme plus compliquŽ dont il est toujours utile de relire (ou de se souvenir de) la preuve~:
\bigskip
{\soustitre ThŽorme des unitŽs de Dirichlet}
\medskip
{\sl Soit $K$ un corps de nombres tel que $[K:\Q]=r+2s$, o $r$ est le nombre de plongements rŽels et 2s, le nombre de plongements complexes. Soit $U_{K}$ le groupe des ŽlŽments inversibles de l'anneau $O_{K}$. Alors on a l'isomorphisme~:
$$U_{K}\simeq W\times\Z^{r+s-1},$$
o $W$ est l'ensemble des racines de l'unitŽ que $K$ contient.

Cf. [Sam, ThŽorme 1, \S 4.4, p. 72].\qed

}

\bigskip
Un autre grand classique~:

\medskip
{\soustitre ThŽorme 90 de Hilbert}
\medskip
{\sl Soit $L/K$ une extension cyclique de corps d'indice fini. Mettons que ${\rm Gal}(L/K)=<\sigma>$. Soit $x\in L$.  Alors $N_{L/K}(x)=1$ si et seulement s'il existe $y\in L^*$ tel que $x={y\over \sigma (y)}$.

Cf. [La1, Thm. 6.6.1, p. 298].\qed

}
\bigskip

A partir de maintenant les choses sŽrieuses commencent~:
\bigskip
\centerline{\soustitre $\taille {14}K$-modules, application d'Artin et compagnie}
\bigskip

Soit $K$ un corps de nombres. Nous allons dŽfinir un objet important. NŽanmoins la nomenclature n'est pas vraiment uniforme dans la
littŽrature~: Janusz les nomme {\it Modulus}, Lang les nomme {\it Cycle}, Neukirch, {\it Modul}, Koch, {\it ErklŠrungsmodul},
Washington, {\it Divisors} Lorenz {\it Modul}... il a bien fallu choisir. La notion de module existe dŽjˆ, mais le nom nous a paru
assez bon tout de mme, aprs d'‰pres discussions, nous nous sommes mis d'accord avec {\it $K$-module} (vous n'allez tout de mme pas confondre avec la notion de $K$-espace vectoriel...).
\bigskip\goodbreak

\defis

Soit $K$ un corps de nombres. Un {\it $K$-module} est une application $$\m\, :\, \gfP(K)\to\N$$ avec les propriŽtŽs suivantes~: 
\art{a)}$\m(\P)=0$ sauf pour un nombre fini de $\P$,

\art{b)}$m(\P)=0$ si $\P\in\gfP_{\C}(K)$,

\art{c)}$m(\P)=0$ ou 1 si $\P\in\gfP_{\R}(K)$.
\newcount\gaagg\gaagg=\gaga

L'usage est d'Žcrire $\m$ comme le produit formel
$$\m=\prod_{\P\in\gfP(K)}\P^{\m(\P)}=\m_{0}\cdot\m_{\infty},$$
o $\m_{0}$ est identifiŽ ˆ un idŽal (l'idŽal $\prod_{\P\in\gfP_{0}(K)}\P^{\m(\P)}$), et $\m_{\infty}$ un sous-ensemble de $\gfP_{\R}(K)$ ( c'est l'ensemble $\{\P\in\gfP_{\R}(K)\mid m(\P)=1\}$) .

Si $\m$ et $\m'$ sont des $K$-modules, on dŽfinit ${\rm pgcd}(\m,\m')$ comme suit~:
$${\rm pgcd}(\m,\m')=\prod_{\P\in\gfP(K)}\P^{\min (\m(\P),\m'(\P))}={\rm pgcd}(\m_{0},\m'_{0})\cdot(\m_{\infty}\cap\m'_{\infty}).$$
On dŽfinit  ${\rm ppcm}(\m,\m')$ de la mme manire en rempla\c cant $\min$ par $\max$ et $\cap$ par $\cup$. Si ${\rm pgcd}(\m,\m')=O_{K}\cdot\emptyset=:{\euf 1}$, on dit que $\m$ et $\m'$ sont {\it premiers entre eux}. Enfin, on dŽfinit $\m\cdot\m':=(\m_{0}\cdot\m'_{0})\cdot (\m_{\infty}\cup\m'_{\infty})$. On vŽrifie facilement que
$$\m\cdot \m'={\rm pgcd}(\m,\m')\cdot{\rm ppcm}(\m,\m').$$

Posons maintenant $S_{0}(\m)=\{\P\in\gfP_{0}(K)\mid \m(\P)>0\},\ S_{\infty}(\m)=\{\P\in\gfP_{\infty}(K)\mid \m(\P)>0\}$ et $S(\m)=S_{0}(\m)\cup S_{\infty}(\m)$. Si $\P\in S(\m)$, on Žcrit $\P|\m$, et on Žcrit $\P\notdiv\m$ dans le cas contraire. 

Si $\m$ et $\m'$ sont des $K$-modules, on dit que $\m|\m'$ si $\m(\P)\leq \m'(\P)$, pour tout $\P\in \gfP(K)$. Dans ce cas-lˆ, il existe un $K$-module $\n$ tel que $\m\cdot\n=\m'$; ce $\n$ n'est pas unique (ˆ cause des places infinies), mais ce n'est pas grave ! on peut prendre par exemple celui dont les places infinies sont disjointes avec celles de $\m$.

Soit $\P\in\gfP_{0}(K)$ et $n\in\N$, $n>0$. On pose 
$$K^*_{\P^n}=\{\alpha\in K^*\mid v_{\P}(\alpha-1)\geq n\}.$$
Soit $\P\in\gfP_{\R}(K)$. Supposons que le plongement associŽ ˆ $\P$ soit $\sigma\, :\, K\to\R$;  on pose 
$$K^*_{\P}=\{\alpha\in K^*\mid\sigma(\alpha)>0\}.$$
Et enfin, si $\m$ est un $K$-module, on posera
$$K^*_{\m}=\bigcap_{\P\in S(\m)}K^*_{\P^{\m(\P)}}.$$

Soit maintenant $x,y\in K^*$, on Žcrira
$$x\equiv y\pmodast \m\iff x\cdot y^{-1}\in K^*_{\m}\iff x\equiv y\pmodast {\P^{\m(\P)}}\ \forall \P\in S(\m).$$
Cela permet d'Žcrire
$$K^*_\m=\{x\in K^*\mid x\equiv 1 \pmodast\m\}.$$
Cette relation d'Žquivalence est cruciale dans tout ce qui va suivre. C'est une relations d'Žquivalence ``multiplicative'', c'est-ˆ-dire $x\equiv y\pmodast\m$ implique que $x\cdot z\equiv y\cdot z\pmodast\m$, pour tout $z\in K$, mais pas que $x+z \equiv\ y+z \pmodast\m$ (par exemple, si $K=\Q$ et $\m=5\cdot \emptyset$, on a ${1\over 3}\equiv 7\pmodast \m$, mais ${1\over 3}+3\not\equiv 7+3\pmodast \m$). Si $\m$ est rŽduit ˆ une seule place, nous noterons $\P^n$ pour $\m$, en oubliant les Žventuels $\emptyset$ ou $O_K$ qui ne feraient qu'alourdir la notation.

Remarquons que $v_\P(x-1)\geq n$ veut dire que $x\in 1+\P^{n}\cdot O_{(\P)}=1+\widetilde{\P}^n$ o, rappelons-le,  $O_{(\P)}$ est le
localisŽ de $O_{K}$ en $\P$. Cela veut aussi dire que $x={a\over b}$ avec $(a,\P)=(b,\P)=1$ (i.e.
$v_\P(a)=v_\P(b)=0$) et $a\equiv b\pmod{\P^{m(\P)}}$.

Remarquons encore que si $x,y\in K^*$, $n\in\N$, $n>0$, et $\P\in\gfP_{0}(K)$, alors 
$$\eqalign{x \equiv y\pmodast{\P^n}&\Leftrightarrow{x-y\over y}\in\widetilde{\P}^n\Leftrightarrow x-y \in y\cdot \widetilde{\P}^n = \widetilde{\P}^{n+v_{\P}(y)}\cr &\Leftrightarrow v_\P(x- y)\geq n+v_{\P}(y)\Leftrightarrow x\equiv y\pmod{\P^{n+v_\P(y)}},\cr }\eqno{(*)}$$
La dernire Žquivalence n'Žtant bien sžr vraie que si $x,y\in O_K$. De plus, si $x \equiv y\pmodast{\P^n}$, alors il est Žvident que $v_{\P}(x)=v_{\P}(y)$.

Si $\P\in \gfP_{\R}(K)$, dire que $x\equiv y \pmodast \P$ veut simplement dire que $\sigma(x)$ et $\sigma(y)$ ont le mme signe, si $\sigma$ est le plongement attachŽ ˆ $\P$.

Soit $\m$ un $K$-module et $x\in K$. Nous dirons que $x$ est un {\it $\m$-entier} si $v_\P(x)\geq 0$ pour tout $\P\in S_0(\m)$. L'ensemble des $\m$-entiers se note $O_{(\m)}$ et on vŽrifie que
$$O_{(\m)}=\bigcap_{\P|\m_0}O_{(\P)},$$
est un anneau commutatif.

\bigskip

\lem

{\sl Soit $K$ un corps de nombres, $\P\in\gfP_{0}(K)$, $u\in O_{(\P)}$ et $n\in\N\setminus\{0\}$. Alors il existe $a\in O_{K}$ tel que $a\equiv u\pmodast{\P^n}$.

}
{\bf Preuve}

Puisque $u\in O_{(\P)}$, on a $u={\alpha\over\beta}$, avec $\alpha\in O_{K}$ et $\beta\in O_{K}\setminus\P$. Puisque $\P$ est un idŽal maximal de $O_{K}$, on a $\beta O_{K}+\P=O_{K}$. On montre alors facilement par rŽcurrence que $\beta O_{K}+\P^n=O_{K}$. Il existe donc $\beta'\in O_{K}$ tel que $\beta\cdot\beta'\equiv 1\pmod{\P^n}$. En posant $a=\alpha\cdot\beta'$, on vŽrifie que ${a-u\over u}=\beta\cdot\beta'-1\in \P^n\subset \widetilde{\P}^n$.\qed

\bigskip

Vient maintenant un thŽorme qui va souvent tre utilisŽ. Il est un peu moins fort que le thŽorme d'approximation faible, voilˆ pourquoi nous l'avons appelŽ le
\bigskip\goodbreak
\th\ {\bf (ThŽorme d'approximation dŽbile)}

{\sl Soit $K$ un corps de nombres, $\m$ et $\m'$ des $K$-modules, et $y$, $z\in K$. Alors 
$$\hbox{il existe $x\in K$ satisfaisant }\left\{\eqalign{x&\equiv y\pmodast\m\cr\ x&\equiv z\pmodast{\m'}\cr}\right .\iff y\equiv z\pmodast{{\rm pgcd}(\m,\m')}.$$

}

{\bf Preuve}

$Ò\Rightarrow"$ : par hypothse, $xy^{-1}\in K^*_\m$ et $zx^{-1}\in K^*_{\m'}$. Ainsi, $y^{-1}z\in K^*_\m\cdot K^*_{\m'}\subset K_{{\rm pgcd}(\m,\m')}$. Cette dernire inclusion est une vŽrification facile~: soient $\alpha\in K^*_{\m}$, $\beta\in K^*_{\m'}$ et $\P\in \gfP_{0}$. Si  $\min(\m(\P),\m'(\P))=0$, il n'y a rien ˆ vŽrifier, supposons donc que $\min(\m(\P),\m'(\P))>0$. On a, par hypothse $\alpha-1\in \widetilde{\P}^{\m(\P)}\subset \widetilde{\P}^{\min(\m(\P),\m'(\P))}$. On a la mme chose pour $\beta-1$. On a donc
$$\widetilde{\P}^{\min(\m(\P),\m'(\P))}\ni (\alpha-1)(\beta-1)=\alpha\beta-1+\underbrace{(1-\alpha)+(1-\beta)}_{\in \widetilde{\P}^{\min(\m(\P),\m'(\P))}}.$$
On a donc $\alpha\beta\in K_{ \widetilde{\P}^{\min(\m(\P),\m'(\P))}}$. Et pour les places infinies, c'est encore plus facile~: si $\P\in \m_{\infty}\cap\m'_{\infty}$, on se souvient que dire que $x\equiv y\pmodast {\P}$ veut simplement dire que $\sigma(x)$ et $\sigma(y)$ ont le mme signe, o $\sigma$ est le plongement associŽ ˆ $\P$; donc, on a $y\equiv z\pmodast\P$.

$Ò\Leftarrow"$ :  pour les places finies, c'est le thŽorme chinois ! voilˆ comment on pourrait rŽsumer cette partie. Le lecteur peut trs bien faire cela en exercice. Mais ce n'est pas le le genre de la maison que de laisser les parties un peu ``laborieuses''. 

Supposons donc que $\m_{0}=\prod_{i=1}^r\P_{i}^{\m(\P_{i})}\cdot\prod_{j=1}^s\qq_{j}^{\m(\qq_{j})}$ et que $\m'_{0}=\prod_{i=1}^r\P_{i}^{\m'(\P_{i})}\cdot\prod_{k=1}^{s'}\rr_{k}^{\m'(\rr_{k})}$. Sachant que $y\equiv z\pmodast {\P_{i}^{t_{i}}}$, o $t_{i}=\min(\m(\P_{i}),\m'(\P_{i}))$, il faut rŽsoudre le systme
$$\left\{\eqalign{\alpha&\equiv y\pmodast {\P_{i}^{\m(\P_{i})}}\quad i=1,\ldots , r_{1}\cr
\alpha&\equiv z\pmodast{\P_{i}^{\m'(\P_{i})}}\quad i=r_{1}+1,\ldots , r\cr
\alpha&\equiv y\pmodast{\qq_{j}^{\m(\qq_{j})}}\quad j=1,\ldots , s\cr
\alpha&\equiv z\pmodast{\rr_{k}^{\m'(\rr_{k})}}\quad k=1,\ldots , s'.\cr}\right .\eqno{(+)}$$
En ayant supposŽ que $r=r_{1}+r_{2}$ et que $\m(\P_{i})\geq\m'(\P_{i})$ si $i=1,\ldots ,r_{1}$ et que $\m'(\P_{i})>\m(\P_{i})$ si $i=r_{1}+1,\ldots , r$. 

On se souvient que chaque anneau $O_{(\P_{i})}$, $O_{(\qq_{j})}$ ou $O_{(\rr_{k})}$ possde une uniformisante qu'on notera $\pi_{\P_{i}},\pi_{\qq_{j}}$ ou $\pi_{\rr_{k}}$. ConsidŽrons l'ŽlŽment 
$$\lambda=\prod_{i=1}^{r_{1}}\pi_{\P_{i}}^{v_{\P_{i}}(y)}\cdot \prod_{i=r_{1}+1}^{r}\pi_{\P_{i}}^{v_{\P_{i}}(z)}\cdot\prod_{j=1}^{s}\pi_{\qq_{j}}^{v_{\qq_{j}}(y)}\cdot\prod_{k=1}^{s'}\pi_{\rr_{k}}^{v_{\rr_{k}}(z)}.$$

Divisons chaque terme de chaque Žquivalence du systme $(+)$ par $\lambda$. On a maintenant que ${y\over \lambda}\in O_{(\P_{i})}^*, O_{(\qq_{j})}^*$ pour $i=1,\ldots , r_{1}$ et $j=1\ldots , s$; de mme,   ${z\over \lambda}\in O_{(\P_{i})}^*, O_{(\rr_{k})}^*$ pour $i=r_{1}+1,\ldots , r$ et $k=1\ldots , s'$.
En vertu du lemme prŽcŽdent, de la remarque $(*)$ prŽcŽdent le Lemme \argb\  et du thŽorme chinois, il existe $a_{i},b_{j},c_{k}\in O_{K}$ et surtout $v\in O_{K}$ satisfaisant le systme

$$\left\{\eqalign{v&\equiv {y\over\lambda}\equiv a_{i}\pmod {\P_{i}^{\m(\P_{i})}}\quad i=1,\ldots , r_{1}\cr
v&\equiv {z\over \lambda}\equiv a_{i}\pmod{\P_{i}^{\m'(\P_{i})}}\quad i=r_{1}+1,\ldots , r\cr
v&\equiv {y\over\lambda}\equiv b_{j}\pmod{\qq_{j}^{\m(\qq_{j})}}\quad j=1,\ldots , s\cr
v&\equiv {z\over\lambda}\equiv c_{k}\pmod{\rr_{k}^{\m'(\rr_{k})}}\quad k=1,\ldots , s'\cr}\right .$$

Ainsi, en posant $\alpha=v\cdot\lambda$, on rŽsout le systme $(+)$. Donc, la partie Òplaces finie'' est rŽsolue.

Pour les places infinies, il suffit de voir la chose suivante~: soient $\sigma_{1},\ldots ,\sigma_{r}$ les plongements infinis rŽels de $K$ et $i\in \{1,\ldots , r\}$. Alors il existe $\alpha_{i}$ tel que $\sigma_{j}(\alpha_{i})>0$ pour $j\ne i$ et $\sigma_{i}(\alpha_{i})<0$. En effet, si on doit trouver un $\beta$ tel que $\sigma_i(\beta)$ doit avoir un signe prescrit, on pose $\varepsilon_i=0$ si $\sigma_i(\beta)$ doit tre positif et $\varepsilon_i=1$ si $\sigma_i(\beta)$ doit tre nŽgatif; le nombre
$$\beta=\alpha_1^{\varepsilon_1}\cdots\alpha_r^{\varepsilon_r}$$
rŽpond ˆ la question. Trouvons donc ces $\alpha_i$. Nous savons que $K=\Q(\theta)$, pour $\theta\in\C$. On suppose sans limiter la gŽnŽralitŽ que $\sigma_1(\theta)<\sigma_2(\theta)<\cdots <\sigma_r(\theta)$. Choisissons $a_0,\ldots ,a_r\in\Q$ tels que
$$a_0<\sigma_1(\theta)<a_1<\sigma_2(\theta)<\cdots <a_{r-1}<\sigma_r(\theta)<a_r.$$
On vŽrifie facilement que
$$\alpha_i={\theta-a_{i-1}\over \theta-a_i}$$
possde les propriŽtŽs requises.

Maintenant, il s'agit de recoller les morceaux~: supposons qu'$\alpha$ est une solution pour la partie finie et $\beta$ est une solution pour la partie infinie, posons en outre $M=\N(\P_1\cdots \P_r\cdot\qq_1\cdots\qq_s\cdot\rr_1\cdots \rr_{s'})$. Alors le nombre
$$x=\alpha+M^N\cdot\beta$$
avec $N$ suffisamment grand pour que ${\rm sgn}(\sigma_i(x))={\rm sgn}(\sigma_i(\beta))$ pour tous les $i$ nŽcessaires et pour que les $v_{\P_i}({M^N\beta\over \alpha})$ , $v_{\qq_j}({M^N\beta\over \alpha})$,  $v_{\rr_k}({M^N\beta\over \alpha})$ soient suffisamment grand pour que $x\equiv \alpha\pmodast{ \P^{\max(\m(\P),\m'(\P))}}$ pour tout $\P\in S({\rm ppcm}(\m,\m'))=\{\P_1,\ldots ,\P_r,\qq_1,\ldots,\qq_s,\rr_1,\ldots ,\rr_{s'}\}$.\qed
\bigskip
Voici maintenant un corollaire que nous utiliserons souvent~:
\bigskip

\goodbreak
\coro

{\sl Soient $\m$ et $\m'$ des $K$-modules premiers entre eux. Alors on a l'isomorphisme
$$K^*/K^*_{\m\m'}\simeq K^*/K^*_\m\times K^*/K^*_{\m'}.$$

}

{\bf Preuve}

L'homomorphisme~:
$$\eqalign{ K^*&\longrightarrow K^*/K^*_\m\times K^*/K^*_\n\cr \alpha&\longmapsto (\alpha K^*_\m,\alpha K^*_\n)\cr}$$
est surjectif. En effet, soit $(\beta,\gamma)\in
K^*\times K^*$. Par le thŽorme d'approximation dŽbile, il existe $\alpha\in K^*$ tel que $\alpha\equiv\beta \pmodast \m$ et $\alpha\equiv \gamma\pmodast{ \m'}$ et donc $\alpha
K^*_\m=\beta K^*_\m$ et $\alpha K^*_\n=\gamma K^*_\n$ ce qui prouve la surjectivitŽ. Le noyau de cet
homomorphisme est $K_\m^*\cap K_\n^*=K^*_{\m\n}$ (cela suit directement des dŽfinitions). On a donc un isomorphisme $K^*/K^*_{\m\n}\simeq K^*/K^*_\m\times K^*/K^*_\n$. \qed
\bigskip\goodbreak

\newcount\gaagh\gaagh=\gaga

\defis

Soit $K$ un corps de nombres. On rappelle que $I_K$ est le groupe des idŽaux fractionnaires. Soit $S$ un ensemble fini de places finies de $K$. On note 
$$I_K^S=\{{\euf a}\in I_K\mid v_\P({\euf a})=0\hbox{ pour tout }\P\in S\}.$$
Si $\m$ est un $K$-module, on note $\thboxed 25{I_K^{S_0(\m)}=I_K(\m)}$, et mme parfois, on note $I_K(\m)=I(\m)$, quand il n'y a pas d'ambigu•tŽ. 

Soit $L/K$ une extension abŽlienne de corps de nombres et $S$ un ensemble fini de places finies de $K$ contenant toutes celles qui ramifient dans $L$. On dŽfinit {\it l'application d'Artin}
$$\thboxed 25{\Phi_{L/K}=\Phi_{L/K}^S\; : I_K^S \rightarrow {\rm Gal}(L/K)}$$
comme suit~: si $\P\not\in S$ est un idŽal premier, $\Phi_{L/K}(\P)= {\rm Frob}_{L/K}(\P)$; o   ${\rm Frob}_{L/K}(\P)$ est l'homomorphisme de Frobenius dŽfinit plus haut; et on prolonge $\Phi_{L/K}$ en un homomorphisme de groupe sur $I_K^S$ tout entier, car les idŽaux premiers qui ne sont
pas dans $S$ engendrent $I_K^S$.  Si $\m$ est un $K$-module, l'application d'Artin $\Phi_{L/K}^{S_0(\m)}$ se note $\Phi_{L/K}^\m$ (ou mme $\Phi_{L/K}$ s'il n'y a pas d'ambigu•tŽ). Disons sans ambages  que nous montrerons que cette application est surjective (ThŽorme \argaq) et  il sera beaucoup question du noyau de cette application dans le reste de ce texte. Remarquons un petit rŽsultat trivial sur ce noyau~: si $\P$ est un idŽal {\bf premier} dans l'ensemble de dŽfinition de $\Phi_{L/K}$, alors on a
$$\eqalign{\P\in\ker(\Phi_{L/K})&\iff Z(\gP/\P)=\{{\rm Id}_L\}\ \forall\ \gP|\P\cr &\iff f(\gP/\P)=1\ \forall\ \gP|\P\cr &\iff\P\hbox{ se dŽcompose totalement dans $L$}\cr &\iff N_{L/K}(\gP)=\P\ \forall\ \gP|\P.\cr}$$
Mais n'allez surtout pas croire que ce noyau n'est fait que d'idŽaux qui sont des normes d'idŽaux de $L$, mais nous allons voir tout bient™t qu'on peut rapidement trouver une inclusion (cf. Corollaire \argg).

Supposons maintenant que $L/K$ est une extension quelconque de corps de nombres. On note $\widetilde{S}=\widetilde{S}_L=\{\gP\in \gfP_0(L)\mid \exists \P\in S ,\gP|\P\}$. Si $\m$ est un $K$-module, on note 
$$\widetilde{\m}=\widetilde{\m}_L=\underbrace{(\m_0\cdot
O_L)}_{=\tilde\m_0}\cdot\underbrace{\prod_{\gP\in \gfP_\R(L), \gP|_K\in \m_\infty}\gP}_{=\tilde\m_\infty}.$$
On dit que $\widetilde{\m}$ est le {\it $L$-module engendrŽ par $\m$}. Et on a bien sžr
$$I(\widetilde{\m})=I_L(\widetilde{\m})=I_L^{\widetilde{S_0(\m)}}.$$

\newcount\gaagi\gaagi=\gaga

\bigskip
\centerline{\soustitre Quelques rŽsultats rapidement prouvables}
\bigskip

\bigskip

\th

{\sl Soit $L/K$ et $E/K$ deux extensions finies d'un corps de nombres $K$. Supposons de plus que $L/K$ soit une
extension abŽlienne de groupe de Galois $G$. La thŽorie de Galois nous apprend que l'extension $EL/E$ est aussi
abŽlienne et que $H:={\rm Gal}(EL/E)$ s'identifie ˆ un sous-groupe de $G$ via la restriction ˆ $L$. Notons $R\, :\,
H\rightarrow G$ cette restriction. Soit $S$ un ensemble finie de places finies de $K$ contenant toutes celles qui ramifient dans $L$, on a~:
$$\Phi_{L/K}^S\circ N_{E/K}=R\circ \Phi_{EL/E}^{\widetilde{S}_E}.$$

}

{\bf Preuve}

On est dans la situation suivante~:

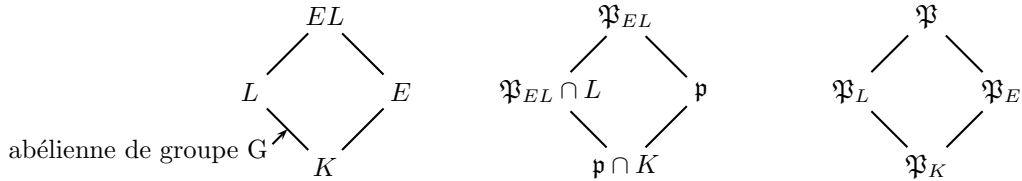
\vglue 2.5cm

\rput(5,0){\rnode{K}{$ K  $}} 
\rput(4,1){\rnode{L}{$ L  $}} 
\rput(6,1){\rnode{E}{$ E  $}} 
\rput(5,2){\rnode{EL}{$ EL  $}} 
\ncline[nodesep=3pt]{K}{L}
\Aput{abŽlienne de groupe G } 
\psline{->}(4.3,0.3)(4.5,0.5)
\ncline[nodesep=3pt]{K}{E}
\ncline[nodesep=3pt]{L}{EL}
\ncline[nodesep=3pt]{E}{EL}
\rput(9,0){\rnode{pcapk}{$ \P\cap K  $}} 
\rput(8,1){\rnode{pelcapl}{$ \gP_{EL}\cap L  $}} 
\rput(10,1){\rnode{p}{$ \P  $}} 
\rput(9,2){\rnode{pel}{$ \gP_{EL}  $}} 
\ncline[nodesep=3pt]{pcapk}{pelcapl}
\ncline[nodesep=3pt]{pcapk}{p}
\ncline[nodesep=3pt]{pelcapl}{pel}
\ncline[nodesep=3pt]{p}{pel}
\rput(13,0){\rnode{pk}{$ \gP_K  $}} 
\rput(12,1){\rnode{pl}{$ \gP_{L} $}} 
\rput(14,1){\rnode{pe}{$ \gP_E  $}} 
\rput(13,2){\rnode{pp}{$ \gP  $}} 
\ncline[nodesep=3pt]{pk}{pl}
\ncline[nodesep=3pt]{pk}{pe}
\ncline[nodesep=3pt]{pl}{pp}
\ncline[nodesep=3pt]{pe}{pp}
\vskip 1cm
Remarquons tout d'abord que $\Phi_{EL/E}$ est bien dŽfini sur $I_E^{\widetilde{S}_E}$.  En effet, soit $\P\not\in \widetilde{S}_E$. Il faut montrer que $\P$ ne ramifie pas dans $EL$. Soit donc $\gP_{EL}$ un idŽal de $EL$ au-dessus de $\P$. On a vu dans la partie ÒRamification..." que $T(\gP_{EL}/\P)$ pouvait tre vu comme un sous-groupe de $T((\gP_{EL}\cap L)/(\gP_{EL}\cap K))$. Or, par dŽfinition de $\widetilde{S}_E$ et de $S$, $\gP_{EL}\cap K=\P\cap K$ ne ramifie pas dans $L$, ce qui montre que $T((\gP_{EL}\cap L)/(\gP_{EL}\cap K))=\{ {\rm Id}_L\}$, et donc que $T(\gP_{EL}/\P)=\{{\rm Id}_{EL}\}$ ce qui montre que $\P$ ne ramifie pas dans $EL$.

Soit, comme dans l'illustration ci-dessus $\gP$ un idŽal premier de $O_{EL}$. Posons $\gP_L=\gP\cap O_L$,
$\gP_E=\gP\cap O_E$ et $\gP_K=\gP\cap O_K$ les idŽaux au-dessous de $\gP$. On suppose que $\gP_K\not\in S$.
Posons $\sigma=\Phi_{EL/E}(\gP_E)$,
$\tau=\Phi_{L/K}(\gP_K)$, $q=\N(\gP_K)$ et $N_{E/K}(\gP_E)=\gP_K^f$ o $f=f(\gP_E/\gP_K)$. Alors on a $\N(\gP_E)=q^f$. Par
dŽfinition de $\sigma$, on a $\sigma(x)\equiv x^{\N(\gP_E)}\pmod{\gP}$ pour tout $x\in O_{EL}$, c'est-ˆ-dire que $\sigma(x)\equiv
x^{q^f}\pmod{\gP}$ pour tout $x\in O_{EL}$; cela implique que $\sigma(x)\equiv x^{q^f}\pmod{\gP_L}$ pour tout $x\in O_{L}$, car $\sigma|_L$
est un automorphisme de $L$. De mme, $\tau(x)\equiv x^q\pmod{\gP_L}$, donc, en composant $f$ fois, on trouve $\tau^f(x)\equiv
x^{q^f}\pmod{\gP_L}$. On en dŽduit que $R\circ \sigma=\tau^f$, car les congruences qu'on vient de montrer donnent que $R\circ
\sigma$ et $\tau^f$ induisent le mme ŽlŽment dans le groupe de Galois de l'extension $(O_L/\gP_L)/(O_K/\gP_K)$. Puisque $\gP_K$
ne ramifie pas dans $O_L$, cette opŽration est injective sur les ŽlŽment de $Z(L/\gP_K)$ dont $R\circ
\sigma$ et $\tau^f$ font partie (cf. partie ÒRamification..."). On en dŽduit que

$$R\circ\Phi_{EL/E}(\gP_E)=\Phi_{L/K}(\gP_K)^f=\Phi_{L/K}(\gP_K^f)=\Phi_{L/K}\circ N_{E/K}(\gP_E).$$
On conclut par multiplicativitŽ.\qed

\bigskip

\coro

{\sl Soit $L/K$ une extension abŽlienne de corps de nombres, et $\Phi_{L/K}\ : \ I_K^S\rightarrow {\rm Gal}(L/K)$, l'application
d'Artin, alors on a $N_{L/K}(I_L^{\widetilde{S}})\subset \ker \Phi_{L/K}$.

}

{\bf Preuve}

On applique le thŽorme prŽcŽdent ˆ $L=E$. Cela donne $\Phi_{L/K}\circ N_{L/K}=\Phi_{L/L}={\rm Id}_L$, ce qui prouve le corollaire.\qed

\bigskip\goodbreak

\defi

Soit $K$ un corps de nombres. On dŽfinit $P=P_{K}\subset I_{K}$ le sous-groupe des idŽaux fractionnaires principaux de $K$. Il est bien connu que le groupe quotient $I_{K}/P_{K}$ est un groupe fini (cf. [Sam, Thm. 2, chap IV, \S 3, p.71]) et Žvidemment rŽduit au groupe trivial si $O_{K}$ est un anneau principal. On appelle ce groupe fini {\it le groupe des classes d'idŽaux}; son cardinal se note $h=h_{K}$.

Soit  $\m$ un $K$-module, on note $P(\m)=P\cap I(\m)$. On note encore $K^*(\m)=\{x\in K^*\mid
O\cdot x\in P(\m)\}=\{x\in K^*\mid v_\P(x)=0$ pour tout $\P \in S_0(\m)\}=\cap_{\P\in S_0(\m)}O^*_{(\P)}$.

On note $P_\m=\{ O_K\cdot x\mid x\in K^*_\m\}$. Remarquons que $O_K\cdot x=O_K\cdot y$ avec $x\in K^*_\m$ n'implique pas forcŽment que
$y\in  K^*_\m$. Nous allons aussi beaucoup Žtudier  $I(\m)/P_\m$ qu'on appellera {\it groupe de classes radiales modulo
$\m$}. En allemand Hasse l'appelle {\it Strahlklassengruppe modulo $\m$}, et en anglais, on nomme cela {\it
Ray class group modulo
$\m$}. Remarquons que si $\m=O_K\cdot\emptyset={\euf 1}$, le groupe des classes radiales est le groupe des classes usuel dont on sait qu'il est fini. On note aussi $h_\m$ le cardinal $\left |I(\m)/P_\m\right |$. Nous allons montrer que ce cardinal est fini et mme en donner une formule. Enfin, on notera $U_\m$ pour $U_K\cap K^*_\m$.
\newcount\gaagj\gaagj=\gaga

\bigskip\goodbreak
\lem

{\sl Soient $K$ un corps de nombres et  $\euf a,\bb$ deux idŽaux fractionnaires de $K$. ConsidŽrons
${\cal C}$ la classe de
$\euf a$ modulo $P_K$. Alors il existe ${\euf c}\in{\cal C}$ tel que $\cc$ soit premier ˆ $\euf b$, 
c'est-ˆ-dire $v_\P({\euf b})\cdot v_\P( \cc)=0$ pour tout idŽal premier $\P$.

}

{\bf Preuve}

C'est aussi un corollaire du thŽorme chinois. Supposons que $\aa$ soit un Òvrai" idŽal de $K$. Posons ${\euf a}=\prod_{\P\in\bbP_0}\P^{v_\P({\euf a})}$, $v_\P(\aa)\in \N$, et $ \bb=\prod_{\P\in\bbP_0}\P^{v_\P({\bb})}$. Soit $V$ l'ensemble des idŽaux premiers qui divisent $\euf a$ ou
$\bb$. Pour tout $\P\in V$, on choisit $x_\P\in
\P^{v_\P({\euf a})}\setminus\P^{v_\P({\euf a})+1}$. Par le thŽorme chinois, il existe $a\in O_K$ tel que
$a\equiv x_\P\pmod{\P^{v_\P({\euf a})+1}}$. L'idŽal fractionnaire ${1\over a}{\euf a}$ rŽpond ˆ la
question, car $v_\P( {1\over a}{\euf a})=v_\P({1\over a})+v_\P({\euf a})=0$ pour tout $\P\in V$; et si
$\P\not\in V$, alors $v_\P(\bb)=0$.

Si $\aa$ est un idŽal fractionnaire quelconque, alors $\aa=\aa'\cdot\aa''^{-1}$, o $\aa$ et $\aa'$ sont des idŽaux de $K$. En construisant $a'$ et $a''\in O_K$ comme ci-dessus pour $\aa'$ et $\aa''$ respectivement, on voit que l'idŽal $\cc={a''\over a'}\cdot\aa$ possde les propriŽtŽs requises. \qed

\bigskip
\lem

{\sl Soient $K$ un corps de nombres et $\m$ un $K$-module. Alors on a l'isomorphisme

$$K^*(\m)/K^*_\m\simeq \prod_{\P\in S_0(\m)}\left (O_{(\P)}/\widetilde{\P}^{\m(\P)}\right )^*\times \{\pm 1\}^s,$$
o $s=|S_\infty(\m)|$ est le nombre de place infinie divisant $\m$.

}
\medskip
{\bf Preuve}

Puisque $K^*(\m)=\cap_{\P\in S_0(\m)}O^*_{(\P)}$, tout $x\in K^*(\m)$ est tel que $( x\;{\rm mod} \;\widetilde{\P}^{\m(\P)})\in\left (O_{(\P)}/\widetilde{\P}^{\m(\P)}\right )^*$. En effet, $x\in O_{(\P)}\setminus \widetilde{\P}$, donc $x\cdot O_{(\P)}+\widetilde{\P}=O_{(\P)}$ ($\widetilde{\P}$ est un idŽal maximal), par suite et par rŽcurrence, on montre que  $x\cdot O_{(\P)}+\widetilde{\P}^{\m(\P)}=O_{(\P)}$, donc il existe $y\in O_{(\P)}$ et $\alpha\in \widetilde{\P}^{\m(\P)}$ tels que $yx+\alpha=1$, ce qui veut dire que $( x\;{\rm mod} \;\widetilde{\P}^{\m(\P)})\in\left (O_{(\P)}/\widetilde{\P}^{\m(\P)}\right )^*$. Ainsi, l'homomorphisme
$$x\mapsto( x\;{\rm mod} \;\widetilde{\P}^{\m(\P)})_{\P\in S_{0} (\m)}\times ({\rm sgn}(\sigma_{\P}(x)))_{\P\in S_{\infty}(\m)}$$
est bien dŽfini de $K^*(\m)$ dans le groupe de droite de l'isomorphisme cherchŽ. Puisqu'on travaille dans $\left (O_{(\P)}/\widetilde{\P}^{\m(\P)}\right )^*$,  dire que $x\equiv y\pmod{\widetilde{\P}^{\m(\P)}}$ revient ˆ dire que $x\equiv y\pmodast{{\P}^{\m(\P)}}$. Finalement,  le thŽorme d'approximation dŽbile montre que notre homomorphisme est surjectif et son noyau est clairement $K^*_{\m}$.\qed

\bigskip

\lem

{\sl Soient $K$ un corps de nombres et $\aa$ un idŽal de $K$. Alors 
$$\left |\left (O_{K}/\aa\right )^*\right |=\N(\aa)\cdot\prod_{\P|\aa}\left (1-{1\over\N(\P)}\right )$$

}

{\bf Preuve}

C'est une gŽnŽralisation de la formule sur l'indicateur d'Euler. Par le thŽorme chinois et le fait que pour tout anneau produit, on a  $(A\times B)^*=A^*\times B^*$, on a $\left (O_{K}/\aa\right )^*=\prod_{\P|\aa}\left (O_{K}/\P^{v_{\P}(\aa)}\right )^*$. Commenons donc par prouver que si  $\aa=\P^k$, $k\in \N$, alors          
$\left |\left (O_{K}/\P^{k-1}\right )^*\right |=\N(\P^k)-\N(\P^{k-1})$. On procde par rŽcurrence sur $k$. Si $k=1$, c'est Žvident, car $O_{K}/\P$ est un corps. Rappelons le fait ŽlŽmentaire suivant~: tout homomorphisme d'anneau $f\, :\, A\to B$ dŽfinit un homomorphisme de groupe $f^*\, :\, A^*\to B^*$ tel que $\ker(f^*)=(1+\ker(f))\cap A^*$; et si $f$ est surjectif et que $(1+\ker(f))\subset A^*$, alors $f^*$ est aussi surjective. Ainsi, supposons le thŽorme vrai pour $k-1$. L'homomorphisme surjectif $f\, : \,O_{K}/\P^k\to O_{K}/\P^{k-1}$ induit un homomorphisme $f^*\, :\, \left (O_{K}/\P^k\right )^*\to \left (O_{K}/\P^{k-1}\right )^*$. Son noyau est $1+\P^{k-1}/\P^k$ (car on montre que $1+\P^{k-1}/\P^k\subset \left (O_{K}/\P^k\right )^*$ par un argument similaire ˆ celui du dŽbut du lemme prŽcŽdent. Ainsi, 
$$\left |\left (O_{K}/\P^k\right )^*\right |=\left |\left (O_{K}/\P^{k-1}\right )^*\right |\cdot \left | 1+\P^{k-1}/\P^k\right |\buildrel (*)\over =\left |\left (O_{K}/\P^{k-1}\right )^*\right |\cdot\N(\P)\buildrel\rm hyp.\ de\ rec.\over =\N(\P^k)-\N(\P^{k-1}).$$
 L'ŽgalitŽ $(*)$ vient du fait que $ \left | 1+\P^{k-1}/\P^k\right |= \left | \P^{k-1}/\P^k\right |=\left |O_K/\P\right |=\N(\P)$. Et on conclut, si $\aa$ est quelconque~:
$$\left |\left (O_{K}/\aa\right )^*\right |=\prod_{\P |\aa}\N(\P^k)-\N(\P^{k-1})=\N(\aa)\cdot\prod_{\P|\aa}\left (1-{1\over\N(\P)}\right ).$$\qed

\goodbreak
\th

{\sl Soient $K$ un corps de nombres et $\m$ un $K$-module. Alors  $h_\m$ est fini. Plus prŽcisŽment,

$$h_\m={2^s\cdot \N(\m_0)\cdot \prod_{\P|\m_0}(1-{1\over \N(\P)})\over [U_K:U_\m]}\cdot h,$$

o $s$ est le nombre de places infinies divisant $\m$ et $h$ est cardinal du groupe des classes usuel.

}

{\bf preuve}

On considre l'application composŽe $I_K(\m)\buildrel {\rm incl.}\over\hookrightarrow I\buildrel{\rm
nat.}\over \rightarrow I_K/P_K$. Elle est clairement surjective, car dans n'importe quelle classe d'idŽaux
usuelle, il est possible de trouver un reprŽsentant premier ˆ un idŽal fractionnaire fixŽ (ici, il
s'agit de $\m_0$ (cf. lemme \argi). Le noyau est clairement $P_K(\m)$. Ainsi, $I_K(\m)/P_K(\m)\simeq I_K/P_K$. Or, $P_\m\subset P_K(\m)$. Donc, on a une suite exacte de groupes abŽliens~:

$$1\lra P_K(\m)/P_\m\lra I_K(\m)/P_\m\lra I_K(\m)/P_K(\m)\lra 1.$$

Il suffit ainsi de montrer que $P_K(\m)/P_{\m}$ est fini, de calculer son cardinal et conclure car $h_\m=|P_K(\m)/P_\m|\cdot h$, par ce qui prŽcde. L'application $\eqalign{K^*(\m)&\rightarrow P(\m)\cr
x&\mapsto Ox\cr}$ est par dŽfinition un homomorphisme surjectif. Donc en composant avec la projection canonique, on a un homomorphisme
surjectif $K^*(\m)\rightarrow P(\m)/P_m$. Le noyau de cet homomorphisme est clairement $U_K\cdot K_\m^*$.
Les Lemmes \argj\ et \argk\ nous montrent que $\left | K^*(\m)/K^*_\m\right | =2^s\cdot \N(\m_0)\cdot \prod_{\P|\m_0}(1-{1\over \N(\P)})$. 

D'autre part,  $K^*(\m)\supset U_K\cdot K^*_\m\supset K^*_\m$. Par les thŽormes d'isomorphismes (partie c)), on  a alors que $U_K\cdot K^*_\m/ K^*_\m\simeq U_K/U_K\cap K^*_\m$; et on a par dŽfinition que $U_\m=U_K\cap K^*_\m$. Cela montre
le thŽorme.\qed

\medskip
Voici un lemme trs important quand on parle de cyclotomie, et c'est ce que nous allons faire les prochains thŽormes~:
\medskip
\lem

{\sl Soit $m\in \N$ et $K$ un corps de nombre contenant une racine $m$-ime de l'unitŽ $\zeta_m$. Soit $p\in\gfP_{0}(\Q)$ tel que $p\notdiv   m$. Soit $\P$ un idŽal de $K$ au-dessus de $p$, c'est-ˆ-dire que $\P\cap\Z=p\Z$.
L'application $\phi \, : \,  O_K\lra O_K/\P:=\F$ envoie Žvidemment $\Z$ sur $\F_p=\Z/p\Z$ et donc $|\F|=\N(\P)=p^f=q$ o $f=[O_K/\P:\Z/p\Z]$. Alors on a $m|q-1$. Plus prŽcisŽment, le groupe $\{ \zeta_m^i\mid 1\leq i\leq m\}\subset O_K$ est envoyŽ injectivement par $\phi$. Son image est donc un sous-groupe cyclique d'ordre de $m$ de $\F^*$.

}

{\bf Preuve}

Il suffit de vŽrifier que $\zeta_m\not\equiv \zeta_m^i\pmod\P$ pour tout $2\leq i\leq m$. Soit
$f=\prod_{i=1}^m(x-\alpha_i)$ un polyn™me. Alors le polyn™me dŽrivŽ ŽvaluŽ en $\alpha_1$ vaut
$f'(\alpha_1)=\prod_{i=2}^m(\alpha_1-\alpha_i)$. On applique cela ˆ
$f=X^m-1=\prod_{i=1}^m(X-\zeta_m^i)$. On trouve alors
$\prod_{i=2}^m(\zeta_m-\zeta_m^i)=m\zeta_m^{m-1}$. Puisque $\P$ est un idŽal premier et que ni $m$,
ni $\zeta_m$ ne sont dans $\P$, on en dŽduit que $\prod_{i=2}^m(\zeta_m-\zeta_m^i)\not\in\P$, donc
aucun $\zeta_m-\zeta_m^i$ n'est dans $\P$.\qed

\bigskip
Nous allons maintenant montrer un exemple important qu'on pourrait qualifier d'exemple
gŽnŽrique. En effet,  c'est le premier rŽsultat qui donne un lien entre le groupe de Galois d'une
extension et un groupe de classes radiales. Tous les autres grands rŽsultats de la thŽorie du corps de classe
seront ÒinspirŽ" de cet exemple.
\bigskip

\th

{\sl Soit $m\in\Z$, $m>1$, $m\not\equiv 2\pmod 4$. Posons $K=\Q$, $L=\Q(\zeta_m)$ et $S=\{p\in \bbP_0(\Q)\mid
p|m\}$. Alors Žvidemment, $I_K^S=I_K(\m)$, avec $\m=m\Z\cdot\infty$ o $\infty$ est l'unique place infinie de $\Q$. Alors on a 

$$I_K(\m)/P_\m\simeq {\rm Gal}(L/K).$$

} 

{\bf Preuve}

L'ensemble $\{{a\over b}\in \Q^*\mid a,b\in\Z$, positifs, premiers ˆ $m\}$ est identifiable ˆ $I_K^S$, car
tout idŽal fractionnaire de $\Q$ est engendrŽ par deux ŽlŽments, et on choisit le positif. Soit $p\not\in
S$. Par dŽfinition de l'application d'Artin, on a dans ce cas, $\Phi_{L/K}(p)=\sigma_p\in {\rm Gal}(L/K)$
caractŽrisŽ par $\sigma_p(\zeta_m)\equiv \zeta_m^p\pmod\P$ pour tout idŽal premier $\P$ de $L$ au-dessus de $p$. Le Lemme $\argo$ nous montre alors que $\sigma_p(\zeta_m)= \zeta_m^p$.  On compose $\Phi_{L/K}$ avec l'isomorphisme ${\rm Gal}(L/K)\buildrel\approx\over \lra (\Z/m\Z)^*$; on notera cette composition $\Psi_{L/K}$. Il est clair que $\Psi_{L/K}({a\over b})={\overline{a}\over \overline{b}}\in (\Z/m\Z)^*$, o $\overline{a}$ et $\overline{b}$ sont les classes de $a$ et $b$ modulo $m$. Cette application est Žvidemment surjective. Et on a $\ker \Phi_{L/K}=\ker \Psi_{L/K}=\{{a\over b}\in I_K^S\mid a\equiv b\pmod m\}\simeq P_\m$ o $\m=m\Z\cdot\infty$. On a donc prouvŽ que $I_K(\m)/P_\m\simeq {\rm Gal}(L/K)$.\qed

\bigskip\goodbreak

\defi

Soit $L/K$ une extension de corps de nombres. Cette extension est appelŽe {\it cyclotomique} s'il existe
$m\in\N,$ $m\geq 1$ tel que $K\subset L\subset K(\zeta_m)$. Puisque l'extension $K(\zeta_m)/K$ est
abŽlienne, l'extension $L/K$ est en particulier abŽlienne aussi.

\bigskip\goodbreak
\th

{\sl Soit $K\subset L\subset K(\zeta_m)$ une extension cyclotomique de corps de nombres. Soit
$\m=\m_0\cdot\m_\infty$, un $K$-module tel que $m|\m_0$ (c'est-ˆ-dire que $m\cdot O_K\supset \m_0$) et
tel que $\m_\infty$ est l'ensemble de toutes les places rŽelles de $K$. Alors l'application d'Artin
$\Phi^\m_{L/K}\; :\; I_K(\m)\rightarrow G:={\rm Gal(L/K)}$ est bien dŽfinie et $P_\m\subset \ker
\Phi^\m_{L/K}$.

}

{\bf Preuve}

Les idŽaux de premiers de $K$ qui ramifient dans $L$ ramifient aussi dans $K(\zeta_m)$, donc divisent $m$. Donc, $\Phi^\m_{L/K}$ est bien dŽfinie. Puisque $\Phi^\m_{L/K}=R\circ
\Phi^\m_{K(\zeta_m)/K}$, o $R$ est la restriction ${\rm Gal}(K(\zeta_m)/K)\rightarrow {\rm Gal}(L/K)$
(ThŽorme \argf), on peut supposer que $L=K(\zeta_m)$. On prend donc $L=K(\zeta_m)$. Appelons $i\; :\; {\rm
Gal}(K(\zeta_m)/K)\rightarrow (\Z/m\Z)^*$ l'homomorphisme injectif usuel. Si $\P$ est un idŽal premier
avec $\P\in I_K(\m)$, on sait que $\Phi^\m_{L/K}(\P)=\sigma$ tel que $\sigma(\zeta_m)\equiv
\zeta_m^{\N(\P)}\pmod \gP$ o $\gP$ est n'importe quel idŽal premier de $K(\zeta_m)$ au-dessus de $\P$.
Par le Lemme \argo, on en dŽduit que $\sigma(\zeta_m)=\zeta_m^{\N(\P)}$. Ainsi,
$i(\Phi^\m_{L/K}(\P))=\overline{\N(\P)}$, o $\overline{\N(\P)}$ est la classe de $\N(\P)$ modulo $m$. Par
multiplicitŽ, on en dŽduit que $i(\Phi^\m_{L/K}({{\euf a}}))={\overline{\N (\euf
a})}$ pour tout ${{\euf a}}\in I_K(\m)$. Soit $xO_K\in P_\m$. On peut
supposer que $x\in K_\m^*$. Par dŽfinition de $K_\m^*$, $x$ est totalement positif, donc
$N_{K/\Q}(x)=|N_{K/\Q}(x)|\in\Q_+^*$. Donc, $\N(xO_K)=N_{K/\Q}(x)$, et ainsi,

$$i(\Phi^\m_{L/K}(xO_K))=\overline{N_{K/\Q}(x)}.\eqno{(i)}$$  

D'autre part, $x\equiv 1\pmodast \m$, cela implique puisque $m\cdot O_K\supset \m_0$, que $x-1\in
m\cdot\bigcap_{\P|m}O_{(\P)}$. Soit $E/K$ une extension telle que $E/\Q$ soit
galoisienne. Notons $\widetilde{O}=O_E$. Pour tout ŽlŽment $\tau\in {\rm Gal}(E/\Q)$, on a $\tau(x)-1\in
m\cdot \bigcap_{\gP|m}\widetilde{O}_{(\gP)}$ ($\gP$ est bien sžr un idŽal premier de $\widetilde{O}$).
Par consŽquent, puisque $N_{K/\Q}(x)$ est un produit de certains de ces $\tau(x)$, on aura $N_{K/\Q}(x)-1\in
m\cdot \bigcap_{\gP|m}\widetilde{O}_{(\gP)}\cap\Q=m\cdot \bigcap_{p|m}\Z_{(p)}:=m\cdot\Z_{(m)}$.
Cela implique que
$N_{K/\Q}(x)\equiv 1\pmodast m$ ce qui implique par $(i)$ et puisque $\Z_{(m)}/m\Z_{(m)}\simeq \Z/m\Z$ que
$i(\Phi^\m_{L/K}(xO_K))=\overline{1}$, et donc (puisque $i$ est injectif) que $\Phi^\m_{L/K}(xO_K)$ est
l'unitŽ du groupe de Galois ${\rm Gal}(L/K)$ et donc que $P_\m\subset \ker
\Phi^\m_{L/K}$.\qed
\bigskip
Remarquons que ce thŽorme sera redŽmontrŽ Òen passant" au Chapitre 7 (Proposition \argcs). Mais pour une raison technique, nous aurons besoin de ce rŽsultat pour montrer le thŽorme de $\check{\rm C}$ebotarev pour les extensions cyclotomiques (Proposition \argaw), c'est pourquoi, nous l'avons dŽjˆ mis ici.
\bigskip\goodbreak
Voici encore un rŽsultat intŽressant en lui-mme qui relve plut™t de la thŽorie de Galois et que nous utiliserons 2 fois~: la premire au Chapitre 3 (ThŽorme \aaargbg) et la seconde au Chapitre 12 (ThŽorme~\argfe)~:
\goodbreak\bigskip
\th

{\sl Soit $K\subset E\subset L$ des corps de nombres. On suppose $L/K$ galoisienne et on pose $G={\rm Gal}(L/K)$ et $H={\rm Gal}(L/E)\subset G$. On note $X$ l'ensemble des classes ˆ droite de $G$ modulo $H$. Il est clair que $|X|=[E:K]$. Soit $\P$ un idŽal premier de $K$ et $\gP$ un idŽal premier de $L$ au dessus de $\P$. On suppose $\gP$ non ramifiŽ sur $\P$. Soit encore $\sigma={\rm Frob}(\gP/\P)$ qui est le gŽnŽrateur de $Z(\gP/\P)$, le groupe de dŽcomposition de $\gP$ sur $\P$. On fait agir $Z(\gP/\P)$ sur $X$ qui se dŽcompose en $r$ orbites $C_1,\ldots ,C_r$ de longueur respectivement $f_1,\ldots ,f_r$. Pour fixer le esprits, il existe $\tau_1,\ldots ,\tau_r\in G$ tels que 
$$C_i=\{H\cdot\tau_i,H\cdot\tau_i\cdot\sigma,\ldots ,H\cdot\tau_i\cdot\sigma^{f_i-1}\},$$
pour $i=1,\ldots ,r$. Alors, en posant $\P_i=\tau_i(\gP)\cap E$, pour  $i=1,\ldots ,r$, on a la liste exacte des idŽaux premiers de $E$ au-dessus de $\P$, c'est-ˆ-dire $\P\cdot O_E=\P_1\cdots\P_r$. Et de plus, on a $f_i=f(\P_i/\P)$.

}

{\bf Preuve}

Il est clair que les $\P_i$ sont des idŽaux premiers au-dessus de $\P$. Montrons d'abord que $\P_i$ ne dŽpend que de $C_i$. Fixons donc une de ces classes $C_i$. Tout d'abord, si $H\cdot\tau_i=H\cdot\tau_i'$, alors $\tau_i'=h\cdot \tau_i$ pour un $h\in H$. Alors on a, puisque $h|_E$ est l'identitŽ~:
$$\tau_i'(\gP)\cap E=h(\tau_i(\gP))\cap E=h(\tau_i(\gP))\cap h(E)\buildrel \rm h\ inj.\over =h(\tau_i(\gP)\cap E)=\tau_i(\gP)\cap E.$$
D'autre part, si on prend un ŽlŽment de $C_i$, disons $H\cdot \tau_i\cdot\sigma^j$, alors on a~:
$$(\tau_i\cdot\sigma^j)(\gP)\cap E=\tau_i(\sigma^j(\gP))\cap E\buildrel\sigma\in Z(\gP/\P)\over =\tau_i(\gP)\cap E.$$
Donc, $\P_i$ ne dŽpend bien que de $C_i$. Fixons donc un $H\cdot\tau\in C_i$. Alors $f_i=|C_i|=$ l'indice du stabilisateur de $H\cdot\tau$ dans  $Z(\gP/\P)$. Ce stabilisateur est 
$$\{\eta\in Z(\gP/\P)\mid H\cdot\tau\cdot\eta=H\cdot\tau\}=\{\eta\in Z(\gP/\P)\mid\eta\in \tau^{-1}\cdot H\cdot \tau\}=Z(\gP/\P)\cap \tau^{-1}\cdot H\cdot \tau.$$
Ainsi,
$$\eqalign{|C_i|&=[Z(\gP/\P):Z(\gP/\P)\cap \tau^{-1}\cdot H\cdot \tau]=[\tau \cdot Z(\gP/\P)\cdot\tau^{-1}:\tau \cdot Z(\gP/\P)\cdot\tau^{-1}\cap H]\cr
&=[Z(\tau(\gP)/\P):Z(\tau(\gP)/\P)\cap H]=[Z(\tau(\gP)/\P):Z(\tau(\gP)/\underbrace{\tau(\gP)\cap E}_{=\P_i})]\cr
&={f(\tau(\gP)/\P)\over f(\tau(\gP)/\P_i)}\buildrel f\ \rm mult.\over =f(\P_i/\P).\cr}$$

Donc les $\P_i$ on la propriŽtŽ cherchŽe. En effet, si $\P_0$ est un idŽal premier de $E$ au-dessus de $\P$, alors il existe $\tau\in G$ tel que $\P_0=\tau (\P)\cap E$ (propriŽtŽ galoisienne). Donc $\P_0=\P_i$, o $i$ est tel que $C_i$ est l'orbite qui contient $H\cdot\tau$. Enfin, 
$$\sum_{i=1}^r f(\P_i/\P)=\sum_{i=1}^r|C_i|=|X|=[E:K].$$
De sorte que les $\P_i$ sont deux ˆ deux disjoints.\qed


\vfill\eject

\global\advance\chapnomb by 1
\nomb=1

\centerline{\para Chapitre 1 :}
\medskip
\centerline{\para Un rŽsultat sur $\taille{18} j(x,{\euf K})$ }
\bigskip
Ce chapitre est en fait consacrŽ ˆ la preuve d'un unique thŽorme qui sera utilisŽ au chapitre suivant. Voici l'ŽnoncŽ de ce rŽsultat~:
\bigskip
\th

{\sl
Soit $K$ un corps de nombres et $\m$ un $K$-module. Alors il existe une constante $\rho_{\m}>0$ telle que
si $\euf K$ est une classe de $I_K(\m)$ modulo $P_\m$, et si $j(x,{\euf K})$ est le nombre d'idŽaux (entiers) ${\euf
a}\in{\euf K}$ tels que $\N({\euf a})\leq x$. Alors 
$$j(x,{\euf K})=\rho_{\m}\cdot x+O(x^{1-{1\over n}})$$
o $n=[K:\Q]$, et si $f(x)$ et $g(x)$ sont des fonctions rŽelles, on Žcrit $f(x)=O(g(x))$ lorsqu'il existe une constante $B>0$ tel que $|f(x)|<B\cdot |g(x)|$.

}
\newcount\gaagk\gaagk=\gaga

\bigskip

\headline={\hfill \phantom{ouuh}\hfill}
{\bf Notations}

Pour tout ce chapitre, fixons $K$ un corps de nombres, $[K:\Q]=n=r+2s$, o $r$ est le nombres de
places infinies rŽelles et $s$ le nombres de places infinies complexes. Soit $\m=\m_0\m_\infty$ un
$K$-module. Posons $r_0$ le nombre de places rŽelles qui divisent $\m$. Posons $$\eqalign{v=v_\m\ :
K&\rightarrow
\R^n\simeq \R^r\times\C^s\cr \alpha&\mapsto (\sigma_1(\alpha),\ldots
,\sigma_{r_0}(\alpha),\sigma_{r_0+1}(\alpha),\ldots ,\sigma_{r}(\alpha),\sigma_{r+1}(\alpha),\ldots
,\sigma_{r+s}(\alpha))\cr }$$ avec
$\sigma_1,\ldots ,\sigma_{r_0}$ les plongements correspondant aux places (rŽelles) qui divisent $\m$,
$\sigma_1,\ldots ,\sigma_{r}$, les plongements correspondant aux places rŽelles et $\sigma_{r+1},\ldots
,\sigma_{r+s}$ les plongements correspondant aux places complexes; on pose pour $j=1,\ldots , s$,
$\sigma_{r+s+j}=\overline{\sigma_{r+j}}$. On dŽfinit 
encore 
$$\eqalign{l\ : K^*&\rightarrow
\R^{s+r}\cr \alpha &\mapsto (\log|\sigma_1(\alpha)|,\ldots ,
\log|\sigma_r(\alpha)|,2\log|\sigma_{r+1}(\alpha)|,\ldots ,2 \log|\sigma_{r+s}(\alpha)|)\cr }$$ 
et 

$$\eqalign{l_0 \ : {\R^*}^n\simeq{\R^*}^r\times{\C^*}^s&\rightarrow \R^{r+s}\cr (x_1,\ldots ,x_r,y_1,\ldots
,y_s)&\mapsto (\log|x_1|,\ldots ,\log|x_r|,2\log|y_1|,\ldots ,2\log|y_s|)\cr}$$
Il est Žvident que si $\alpha\ne 0$, on a $l(\alpha)=l_0(v(\alpha))$. Enfin on pose 

$$\eqalign{N_0\ : {\R^*}^n\simeq{\R^*}^r\times{\C^*}^s&\rightarrow \R\cr
(x_1,\ldots ,x_r,y_1,\ldots ,y_s)&\mapsto |x_1|\cdot\cdots\cdot |x_r|\cdot |y_1|^2\cdot\cdots\cdot
|y_s|^2.\cr}$$
Si on pose $N (\alpha)=|N_{K/\Q}(\alpha)|$, on a $N(\alpha)=N_0(v(\alpha))$.
\bigskip

Fixons encore $\euf K$ une classe de $I(\m)/P_\m$. On remarque dŽjˆ qu'il existe $\euf b$ un idŽal entier
(i.e. inclus dans $O_K$) dans la classe ${\euf K}^{-1}$. En effet, si ${\euf c}\over {\euf d}$ est un idŽal
fractionnaire d'une classe quelconque de $I(\m)/P_\m$ avec $\euf c$ et $\euf d$ des idŽaux entiers, alors
il existe $t\in\N$ tel que ${\euf d}^t\in P_\m$, puisque $I(\m)/P_\m$ est fini (cf. ThŽorme~\argl). Alors ${{\euf c}\over {\euf d}}\cdot{\euf d}^t={\euf c}\cdot{\euf d}^{t-1}$ est un idŽal entier dans
la mme classe que ${\euf c}\over{\euf d}$. Fixons donc un de ces $\euf b$ dans ${\euf K}^{-1}$. Soit ${\euf a}\in {\euf K}$, un idŽal. On a alors ${\euf a}\cdot{\euf b}=(\alpha)\in P_\m$, i.e. $\alpha\in K^*_\m\cap O_K$.
Fixons enfin $\alpha_0\in O_K$ tel que $\alpha_0\equiv
0\pmod {\euf b}$ et $\alpha_0\equiv 1\pmod {{\euf m}_0}$ (c'est possible gr‰ce au thŽorme
chinois, et puisque $\bb$ et $\m_0$ sont premiers entre eux).
\bigskip\goodbreak

\lem

{\sl Sous les mmes notations que prŽcŽdemment, on a $j(x,{\euf K})$ est Žgal au nombre
d'idŽaux principaux $(\alpha)$ avec 

\art{a)} $\alpha\equiv\alpha_0\pmod {\m_0\cdot {\euf b}}$

\art{b)} $\sigma_i(\alpha)>0$\quad $i=1,\ldots , r_0$

\art{c)} $0<N(\alpha)\leq x\cdot\N({\euf b})$.

}

{\bf Preuve}

Si ${\euf a}\in {\euf K}$ est un idŽal tel que $\N({\euf a})\leq x$, alors ${\euf a}\cdot {\euf b}=(\alpha)\in
P_\m$, donc  $\alpha$ est tel que  $\sigma_i(\alpha)>0$ pour $i=1,\ldots , r_0$, $\alpha\equiv 1\pmod
{\m_0}$; d'autre part $\alpha\in {\euf b}$ et $N(\alpha)=\N({\euf a}\cdot{\euf b})=\N({\euf
a})\cdot\N({\euf b})\leq x\cdot\N({\euf b})$. Cela prouve que $\alpha$ satisfait les conditions a), b)
et c). RŽciproquement, si $\alpha$ vŽrifie ces trois conditions, on a ${\euf a}={\euf b}^{-1}\cdot
(\alpha)\in{\euf K}$ est un idŽal entier, car ${\euf a}={\euf b}^{-1}\cdot
(\alpha)\subset \bb^{-1}\cdot\bb=O_K$; et $\N(\aa)=\N(\bb)^{-1}\cdot N(\alpha)\leq x$. Ce qui prouve notre lemme.\qed 

\bigskip

\headline={\hfill \smcap Un rŽsultat sur $\taille{10} j(x,{\euf K})$\hfill}

Notons temporairement $U$ les unitŽs de $O_K$. Si $\alpha$ satisfait les conditions a), b) et c), alors
l'ensemble des
$\beta$ tels que $(\alpha)=(\beta)$ et qui vŽrifient aussi a), b) et c) est exactement l'ensemble des
$\alpha\cdot u$, avec $u\in U_\m$ ($=U\cap K_\m^*$).  La preuve du thŽorme des unitŽs de Dirichlet
montre que $l(U)$ est un sous-$\Z$-module libre de rang $r+s-1$ de $\R^{r+s}$ et le sous-espace qu'il
engendre est l'hyperplan $H$ des $(x_i)_{1\leq i\leq r+s}$ tels que $\sum_{i=1}^{r+s} x_i=0$; et
alors $U\simeq W\times l(U)$ o $W$ est le sous-groupe de $U$ formŽ des racines de l'unitŽ (il
est fini et cyclique). Puisque $U_\m$ est d'indice fini dans $U$ (cf. ThŽorme \argl, il a les mmes propriŽtŽs
que $U$, c'est-ˆ-dire que l'on a 
$U_\m\simeq W_\m\times l(U_\m)$ o $W_\m=W\cap U_\m$, et $l(U_\m)$ est un $\Z$ module de rang $r+s-1$ qui
engendre le $\R$-sous-espace $H$. Notons $w_\m$ le cardinal de $W_\m$ et on choisit $e_1,\ldots
,e_{r+s-1}$ une $\Z$-base de $l(U_\m)$ (donc aussi une $\R$-base de $H$) que l'on complte en une
$\R$-base de $\R^{r+s}$ en ajoutant $e_0=(\underbrace{1,\ldots , 1}_{r\ {\rm fois}},2,\ldots ,2)$.
En vertu de la premire remarque de ce paragraphe, on peut, pour chaque idŽal principal comme
dans le Lemme \argr, choisir un gŽnŽrateur $\alpha$, qui outre a), b), c) vŽrifie la condition 
$$l(\alpha)=c_0 e_0+\sum_{i=1}^{r+s-1} c_i e_i\quad\hbox{ avec } 0\leq c_i<1\hbox{ pour }i=1,\ldots
,r+s-1.$$
Dans ce cas $\alpha$ est entirement dŽterminŽ ˆ un facteur $w\in W_\m$ prs. On a ainsi montrŽ le
\bigskip

\lem

{\sl Sous les mmes notations que dans le paragraphe prŽcŽdent, on a $w_\m\cdot j(x,{\euf K})$
est le nombre d'ŽlŽments $\alpha\in\alpha_0+{\euf b}\m_0$ tels que 

$$\displaylines{\sigma_i(\alpha)>0\hbox{ si } 1\leq i\leq r_0,\quad 0<N_0(v(\alpha))\leq x\cdot\N({\euf
b})\quad \hbox{ et }\cr l(\alpha)=c_0 e_0+\sum_{i=1}^{r+s-1} c_i e_i\quad\hbox{ avec } 0\leq c_i<1\hbox{
pour }i=1,\ldots ,r+s-1.\cr}$$

\qed

}
\bigskip

Interrompons-nous un instant dans notre discours pour Žnoncer un petit lemme sur le plongement $v$ vu au
dŽbut de ce chapitre.
\bigskip
\lem

{\sl Soit $K/\Q$ un corps de nombres. Alors $v(O_K)$ est un {\it $\Z$-rŽseau plein} de $\R^n$, c'est-ˆ-dire
un sous-$\Z$-module libre de rang $n$ de $\R^n$ contenant une $\R$-base de $\R^n$.

}
{\bf Preuve} 

L'homomorphisme $v$ est clairement injectif, donc $v(O_K)$ est Žvidemment un $\Z$-module de
rang $n$. Reste ˆ voir qu'il est plein. Soit $\omega_1,\ldots ,\omega_n$ une $\Z$-base de $O_K$. Il suffit
de voir que la matrice 
$$A:=\pmatrix{\sigma_1(\omega_1)&\cdots&\cdots &\sigma_1(\omega_n)\cr \vdots&&&\vdots\cr
\sigma_r(\omega_1)&\cdots&\cdots &\sigma_r(\omega_n) \cr \Re\sigma_{r+1}(\omega_1)&\cdots&\cdots
&\Re\sigma_{r+1}(\omega_n)\cr \Im\sigma_{r+1}(\omega_1)&\cdots&\cdots&
\Im\sigma_{r+1}(\omega_n)\cr
\vdots&&&\vdots \cr 
\Re\sigma_{r+s}(\omega_1)&\cdots&\cdots&
\Re\sigma_{r+s}(\omega_n)\cr \Im\sigma_{r+s}(\omega_1)&\cdots&\cdots&
\Im\sigma_{r+s}(\omega_n)\cr}$$ est de dŽterminant non nul. On vŽrifie sans peine (en utilisant la relation
$2\Im(z)+i\cdot z=i\cdot \overline{z}$)   que ce dŽterminant vaut
$2^{-s}\cdot i^s\cdot\det(\sigma_i(\omega_j)_{1\leq i,j\leq n}$ avec une Òbonne" numŽrotation des $\sigma_i$.
donc $|\det(A)|=2^{-s}\cdot |d(K)|^{1\over 2}$, o $d(K)$ est le discriminant de $K$ sur $\Q$ qui est 
non nul (cf.  [Sam, \S 2.7, Proposition 3, p. 47]).\qed

\bigskip

{\bf Notation}

On notera $\Gamma$ l'ensemble des $x\in \R^s\times\C^s\simeq \R^n$ tels que 

\art{i)}les $r_0$ premires coordonnŽes de $x$ sont $>0$.

\art{ii)}$0<N_0(x)\leq 1$.

\art{iii)}$l_0(x)=c_0 e_0+\sum_{i=1}^{r+s-1} c_i e_i$, avec $0\leq c_i<1$ pour $1\leq i\leq r+s-1$.

\bigskip

On voit facilement que si $t>0$, $l_0(tx)=l_0(x)+\log(t)\cdot e_0$ (ceci gr‰ce ˆ la dŽfinition judicieuse
de $e_0$ !!). De sorte que $x$ remplit la condition iii) si et seulement si $tx$ la remplit. Donc
l'ensemble $t\cdot\Gamma$ est l'ensemble des $x\in \R^s\times\C^s $ qui satisfont les conditions i) et
iii) prŽcŽdentes et la condition $0<N_0(x)\leq t^n$. Posons $\Lambda=v({\euf b}\m_0)$. C'est un
$\Z$-rŽseau plein de $\R^n$ (c'est-ˆ-dire un sous-$\Z$-module libre de rang $n$ de $\R^n$ contenant une
$\R$-base de $\R^n$), car
$v(O_K)$ en est un gr‰ce au lemme prŽcŽdent et
${\euf b}\m_0$ est d'indice fini dans
$O_K$. Posons enfin $\Lambda_0=v(\alpha_0+{\euf b}\m_0)=v(\alpha_0)+\Lambda$, le translatŽ de $\Lambda$
par $v(\alpha_0)$. Ce qui prŽcde se traduit alors ainsi~:

\bigskip

\lem

{\sl Sous les mmes notations, soit pour tout $t>0$, $M(t)$ le nombre d'ŽlŽments de $\Lambda_0$ contenu
dans $t\cdot\Gamma$, alors on a
$$w_\m\cdot j(x,{\euf K})=M(t)\quad\hbox{ o } t^n=x\cdot\N({\euf b}).$$

}\qed

\bigskip

On va faire maintenant un Òrappel" sur la mesure de Jordan (voir [Apo] pour les dŽtails).
\newcount\gaagl\gaagl=\gaga

\bigskip

{\bf Rappel}

Soit $A\subset\R^n$ bornŽ. Un {\it pavŽ} dans $\R^n$ est une produit d'intervalles bornŽs de $\R$. Le
volume d'un tel pavŽ est le produit des longueurs de ces intervalles (les intervalles peuvent tre
ouverts, fermŽs, semi-ouverts ou un point (dans ce cas, le volume est nul)). Le {\it volume extŽrieur} de
$A$, $\overline{v}(A)$, est l'infimum des $\sum_{i=1}^N{\rm vol}(P_i)$ pris sur tout recouvrement
$\{P_1,\ldots ,P_N\}$ de $A$ par un nombre fini de pavŽs. Le {\it volume intŽrieur } $\underline{v}(A)$
est le supremum des $\sum_{i=1}^N{\rm vol}(P_i)$ pris sur toute famille  de pavŽs $\{P_1,\ldots ,P_N\}$
telles que
$\buildrel\hskip 3pt \circ\over P_i\cap \buildrel\hskip 3pt \circ\over P_j=\emptyset$ pour tout $i\ne j$ ($\buildrel\hskip 3pt \circ\over  P_i$ veut dire l'intŽrieur de $P_i$) et $\cup_i P_i\subset\buildrel\hskip 3pt \circ\over A$. On a toujours $ \underline{v}(A)\leq \overline{v}(A)$. On dit alors que $A$ est {\it $J$-mesurable} ($J$ pour Jordan) si
$\underline{v}(A)= \overline{v}(A)$. On pose alors ${\rm vol(A)}=\underline{v}(A)= \overline{v}(A)$. Remarquons que si $M$ est une matrice $n\times n$, et $P$ un pavŽ, alors ${\rm vol}(M\cdot P)=|\det(M)|\cdot {\rm vol(P)}$, ainsi le volume de tout parallŽlotope  est connu.

\bigskip\goodbreak

 \th
 
{\sl Si $A\subset \R^n$ est bornŽ, il est $J$-mesurable si et seulement si $\overline{v}(\partial A)=0$ (o $\partial A$, le bord de $A$ veut dire $\overline{A}\setminus \buildrel\hskip 3pt \circ\over A$, et $\overline{A}$ est l'adhŽrence de $A$).

}

{\bf Preuve}

Cf. [Apo, Thm. 14.9, p. 397].\qed

\bigskip

\defi

Si $(X, d_X)$ et $(Y, d_Y)$ sont des espaces mŽtriques $f\, :\, X\rightarrow Y$ est dite {\it lipschitzienne}
s'il existe $M>0$ tel que $d_Y(f(x),f(y))\leq  M\cdot d_X(x,y)$ pour tout $x,y\in X$. On appellera $M$
{\it constante de Lipschitz}. Toute fonction lipschitzienne est en particulier continue. Par exemple, si $f\, :\,
\R^n\rightarrow\R^m$ est telle que les dŽrivŽes partielles existent et sont continues, alors $f$ est
lipschitzienne sur tout intervalle compact.

\bigskip

\th

{\sl Si $f\, :\, [0;1]^k\rightarrow\R^n$ est lipschitzienne et $k<n$. Alors $\overline{v}(f([0,1]^k))=0$.

}

{\bf Preuve}

On a donc par hypothse il existe $M>0$ tel que $|f(x)-f(y)|\leq M|x-y|$ pour tout $x,y\in [0;1]^k$.
DŽcoupons $[0;1]^k$ en $N^k$ sous-cubes de c™tŽs ${1\over N}$. Si $C$ est l'un de ces cubes, alors $f(C)$
sera de diamtre $\leq M\cdot\sqrt k\cdot {1\over N}$. Si on pose $B=(M\cdot \sqrt{k}+2)^n$, alors $f(C)$
rencontre au plus $B$ sous-cubes de $\R^n$ de la forme $[{i_1\over N};{i_1+1\over N}]\times\cdots\times
[{i_n\over N};{i_n+1\over N}]$, avec $i_1,\ldots ,i_n\in\Z$. Ainsi $f([0;1]^k)$ est inclus dans la rŽunion
d'au plus $B\cdot N^k$ cubes de c™tŽ ${1\over N}$, donnant un volume infŽrieur ou Žgal ˆ $B\cdot{N^k\over
N^n}={B\over N^{n-k}}<\varepsilon$ si $N$ est assez grand. Ce qui prouve que le volume extŽrieur est aussi
petit que l'on veut.\qed

\bigskip

\defi

Une partie $A\subset \R^n$ est dite {\it $(n-1)$-Lipschitz paramŽtrisable} si $A$ est contenu dans la
rŽunion d'un nombre fini d'images d'applications lipschitzienne $f_i\, :\,[0;1]^{n-1}\rightarrow\R^n$.

\bigskip
\goodbreak
\th

{\sl Soit $A\subset \R^n$ bornŽe. On suppose que $\partial A$ est $(n-1)$-Lipschitz paramŽtrisable. Soit
$\Lambda$ un rŽseau (ou translatŽ de rŽseau) dans $\R^n$. Posons $N(t)$ le nombre de points de $\Lambda$
contenu dans $t\cdot A$. Soit ${\rm vol}(\Lambda)$ le volume d'un parallŽlotope
 fondamental de $\Lambda$. Alors $A$ est $J$-mesurable et

$$N(t)={{\rm vol}(A)\over {\rm vol}(\Lambda)}\cdot t^n+O(t^{n-1}).$$
On rappelle que $f(x)=O(g(x))$ veut dire qu'il existe une constante $B$ telle que $|f(x)|\leq  B \cdot
|g(x)|$.

}

{\bf Preuve}

Le ThŽorme \argx\ et un petit raisonnement sur les recouvrements d'une union finie d'ensemble montrent que
$\overline{v}(\partial A)=0$ donc par le ThŽorme \argw\ que $A$ est  $J$-mesurable. Soit $P$ le parallŽlotope fondamental
de $\Lambda$. Posons $n(t)$ le nombre de $x\in \Lambda$ tels que $x+P\subset t\cdot \buildrel\hskip 3pt
\circ\over A$ et $s(t)$ le nombre de $x\in \Lambda$ tels que $(x+P)\cap t\cdot\overline{A}\ne\emptyset$. Il
est clair que $n(t)\leq N(t)\leq s(t)$, et donc que $n(t)\cdot{\rm vol}(P)\leq {\rm vol}(t\cdot
A)=t^n\cdot{\rm vol}(A)\leq s(t)\cdot{\rm vol}(P)$, ou encore $n(t)\leq {{\rm vol}(A)\over {\rm
vol}(\Lambda)}\cdot t^n\leq s(t)$. Cela implique que $|N(t)-{{\rm vol}(A)\over {\rm vol}(\Lambda)}\cdot
t^n|\leq s(t)-n(t)$. Il suffit donc de montrer que $s(t)-n(t)=O(t^{n-1})$. 

Nous avons que $s(t)-n(t)$ est le nombre de $x\in \Lambda$ tels que $(x+P)\cap( \partial tA=t\partial A)\ne\emptyset$.
Par hypothse $\partial A$ est contenu dans la
rŽunion d'un nombre fini d'images d'applications lipschitzienne $f_i\, :\,[0;1]^{n-1}\rightarrow\R^n$. Il suffit
donc de montrer que si $f\, :\,[0;1]^{n-1}\rightarrow\R^n$ est lipschitzienne, alors le nombre $R(t)$ des $x\in
\Lambda$ tels que $(x+P)\cap t\cdot f([0;1]^{n-1})\ne\emptyset$ est un $O(t^{n-1})$. 

Divisons $[0;1]^{n-1}$ en $[t]^{n-1}$ sous-cubes de c™tŽ ${1\over [t]}$ (ne pas confondre l'intervalle
$[0;1]$ avec la partie entire de $t$ qu'on note $[t]$). Soit $C$ un de ces sous-cubes. On a que $f(C)$
est de diamtre $\leq M\cdot\sqrt{n-1}\cdot {1\over [t]}\leq {2M\sqrt{n-1}\over t}$ (si $t\geq 3$). Donc,
$t\cdot{\rm diam}(f(C))={\rm diam}(t\cdot f(C))\leq 2M\sqrt{n-1}$. Cela implique qu'il existe une
constante $B>0$, indŽpendante de $C$ et de $t$ telle que $t\cdot f(C)$ rencontre au plus $B$ rŽgions de
la forme $x+P$ ($x\in \Lambda$). Donc, $t\cdot f([0;1]^{n-1})$ rencontre au plus $B\cdot [t]^{n-1}\leq
B\cdot t^{n-1}$ rŽgions de la forme $x+P$ ($x\in P$), ce qui veut dire que $R(t)\leq B\cdot t^{n-1}$. Cela
prouve notre thŽorme.\qed
\bigskip

\th

{Soit le $\Gamma$ du Lemme \argu. Alors $\Gamma$ est bornŽ et $\partial
\Gamma$ est $(n-1)$-Lipschitz paramŽtrisable.

}

{\bf Preuve}

Montrons tout d'abord que $\Gamma$ est bornŽ. Soit $x=(x_1,\ldots , x_r, y_1,\ldots , y_s)\in
 \R^r\times\C^s\simeq \R^n$. La condition ii) de la dŽfinition de $\Gamma$
se traduit en $x\in {\R^*}^r\times {\C^*}^s $ et $\sum\log|x_i|+2\sum\log|y_i|\leq 0$. Cela implique que
$c_0\cdot n\leq 0$ (car la somme des coefficients d'un vecteur est linŽaire, les sommes des coefficients
des $e_i$ est nulle pour $i=1,\ldots ,r+s-1$ et que celle des coefficients de $e_0$ vaut $n$); ce qui
est Žquivalent ˆ $c_0\leq 0$. Ainsi, $\Gamma$ est l'ensemble des $x\in {\R^*}^r\times {\C^*}^s$ tels que
les $r_0$ premires composantes de $x$ sont positives, c'est la condition (i) et
$l_0(x)=c_0e_0+\sum_{i=1}^{r+s-1} c_i e_i$ avec
$0\leq c_i< 1$ et $c_0\leq 0$, c'est la condition qu'on appellera (ii)'. On en dŽduit que les composantes (dans la base
canonique de
$\R^{r+s}$) de $l_0(x)$ sont bornŽs supŽrieurement et par consŽquent aussi les $|x_i|$ et les $|y_i|$ aussi
(ˆ cause des logarithmes); c'est-ˆ-dire que $\Gamma$ est bornŽ.

Montrons maintenant que $\partial \Gamma$ est $(n-1)$-Lipschitz paramŽtrisable. Pour cela, il suffit de
le montrer pour $\Gamma_0$ qui est l'ensemble des $x\in {\R^*}^r\times {\C^*}^s$ tels que
les $r$ premires composantes de $x$ sont positives et qui satisfait (ii)'; en effet, si
$x=(x_1,\ldots, x_i,\ldots , x_r, y_1,\ldots , y_s)$ satisfait la condition (ii)', alors $(x_1,\ldots,
-x_i,\ldots , x_r, y_1,\ldots , y_s)$ le satisfait aussi, donc $\Gamma$ est symŽtrique par rapport aux
hyperplans
$x_i=0$, $i=r_0+1,\ldots , r$ et si d'un c™tŽ, le bord est $(n-1)$-Lipschitz paramŽtrisable, l'autre c™tŽ le
sera aussi. On a en outre que ${\rm vol}(\Gamma)=2^{r-r_0}\cdot{\rm vol}(\Gamma_0)$. Ecrivons
$e_i=(e_i^{(1)},\ldots ,e_i^{(r+s)})$, $i=1,\ldots , r+s-1$. Soit $x=(x_1,\ldots , x_r, y_1,\ldots ,
y_s)\in\Gamma_0$. Les conditions pour $x$ et $l_0(x)$ sont alors que $x_1,\ldots , x_r>0$, et
$\log(x_i)=\sum_{k=1}^{r+s-1}c_k\cdot e_k^{(i)}+c_0$, $1\leq i\leq r$ et
$2\, \log|y_j|=\sum_{k=1}^{r+s-1}c_k\cdot e_k^{(r+j)}+2c_0$, $1\leq j\leq s$ avec $c_k\in [0;1[$, $k=1,\ldots
, r+s-1$ et $c_0\in ]-\infty;0]$. Remplaons $e^{c_{0}}$ par $c_{r+s}$, on a alors $c_{r+s}\in ]0;1]$. Posons aussi,
pour $1\leq j\leq s$, $y_j=\rho_j e^{i\theta_j}$ ( $\rho_j>0$ et $0\leq\theta_j<2\pi$). On peut alors
dŽcrire $\Gamma_0$ comme Žtant l'ensemble des $(x_1,\ldots , x_r,\rho_1 e^{i\theta_1},\ldots , \rho_s
e^{i\theta_s})\in \R^r\times\C^s$ tels que 

$$\left .\eqalign{ x_i&=c_{r+s}\cdot e^{\sum_{k=1}^{r+s-1} c_k e_k^{(i)}}\ 1\leq i\leq r\cr
\rho_j&= c_{r+s}\cdot e^{{1\over 2}\sum_{k=1}^{r+s-1} c_k e_k^{(r+j)}}\ 1\leq j\leq s\cr 
\theta_j&=2\pi c_{r+s+j}\hskip 1.5cm 1\leq j\leq s\cr}\right\} (*)$$
avec $c_{r+s}\in]0;1]$ et $c_k\in [0;1[$ pour toutes les valeurs de $k\ne r+s$ ($1\leq k\leq n$). Soit

$$\eqalign{f:\ [0,1]^n&\longrightarrow \R^r\times \C^s\cr (c_1,\ldots ,c_n)&\longmapsto (x_1,\ldots ,
x_r,\rho_1 e^{i\theta_1},\ldots ,\rho_s e^{i\theta_s})\cr}$$
donnŽe par les relations $(*)$. Il est clair que $f$ est continue et $f([0;1[^{r+s-1}\times ]0;1]\times
[0;1[^{s-1})=\Gamma_0$ et donc $f([0;1]^n)=\overline{\Gamma}_0$ cela vient du fait que image de tout
compact est un compact et l'image de l'adhŽrence est inclu dans l'adhŽrence de l'image. On a aussi
$f(]0;1[^n)\subset \Gamma_0$. La diffŽrence entre le cube $[0;1]^n$ et son intŽrieur est $2n$ cubes fermŽ
de dimension $n-1$; de plus $f$ est lipschitzienne car partout continžment dŽrivable. Il suffit donc pour
conclure de montrer que $f(]0;1[^n)\subset \buildrel\hskip 3pt \circ\over \Gamma_0$, ou encore que $f\
:]0;1[^n \rightarrow \R^n$ est une application ouverte (i.e. l'image d'un ouvert est un ouvert). Et pour
cela, on observe que $f$ est la composŽe des quatre applications suivantes, manifestement ouvertes~: 
$]0;1[^n\buildrel f_1\over\rightarrow \R^n\buildrel f_2\over\rightarrow \R^n \buildrel
f_3\over\rightarrow \R^r\times ]0;\infty[^s\times\R^s\buildrel f_4\over\rightarrow\R^r\times\C^s$.
dŽfinies de la faon suivante~:

\art{a)}$f_1((t_1,\ldots ,t_n))=(t_1,\ldots , t_{r+s-1},\log (t_{r+s}),t_{r+s+1},\ldots , t_n)$.

\art{b)}$f_2$ est l'application linŽaire $(u_1,\ldots ,u_n)\mapsto (u_1,\ldots ,u_n)\cdot M$, o
\bigskip\goodbreak

\def\strat{\vrule height5pt depth5pt width0pt}
\def\strot{\vrule height9pt depth2pt width0pt}

$\hskip5.6cm \raise-14pt\hbox{$\overbrace{\phantom{\hbox to 73pt
{\hrulefill}}}^{r+s
}$}$

$$M=\pmatrix{\hbox to 14pt{\hrulefill}&e_1&\hbox to 14pt{\hrulefill}&0&\cdots &0\cr
&\vdots&\phantom{iuhoi}\smash{\vrule height 20pt depth40pt}&\vdots& &\vdots \cr
\hbox to 14pt{\hrulefill}&e_{r+s-1}&\hbox to 14pt{\hrulefill}& \vdots& &\vdots\cr
\strat \hbox to 14pt{\hrulefill}&e_0&\hbox to 14pt{\hrulefill}& 0&\cdots &0\cr
\noalign{\hrule}
\strot & 0_{s\times r+s}&&&I_s&\cr}$$
L'application $f_2$ est ouverte, car $M$ est inversible.

\art{c)}Pour $f_3$, on applique $x\mapsto e^x$ aux $r$ premires coordonnŽes, $x\mapsto e^{{1\over 2}x}$
aux $s$ suivantes et on multiplie les $s$ dernires par $2\pi$.

\art{d)}$f_4((x_1,\ldots , x_r,\rho_1,\ldots ,\rho_s,\theta_1,\ldots ,\theta_s))=(x_1,\ldots , x_r,
\rho_1 e^{i\theta_1},\ldots ,\rho_s e^{i\theta_s})$, l'application $(\rho,\theta)\mapsto \rho
e^{i\theta}$ Žtant une application ouverte de $]0;\infty[\times\R$ dans $\C$, car l'image d'un cube ouvert
$(\Delta \rho ,\Delta\theta)$ donne le domaine suivant~:
\vglue 1cm

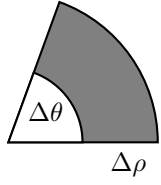
\pspicture(-5,-3)(10,2)
 \pswedge[fillstyle=solid,fillcolor=gray]{2}{0}{70} 
 \pswedge[fillstyle=solid,fillcolor=white]{1}{0}{70}
 \rput(1.6,-0.3){$\Delta \rho$}
 \rput(0.5,0.4){$\Delta \theta$}
\endpspicture
\bigskip

{\bf Preuve du ThŽorme \argq} 

On a~:

$$\eqalign{j(x,{\euf K})&\buildrel \rm Lemme\ \argu\over = {1\over w_{\m}}\cdot M(x^{1\over n}\cdot \N(\bb)^{1\over n})\cr
&\buildrel\rm Thm.\ \argz\ et\ \argaa\over ={1\over w_{\m}}\cdot\left ({{\rm Vol}(\Gamma)\over {\rm Vol(\Lambda_{0})}}\cdot\N(\bb)\cdot x+O((x^{1\over n}\cdot\N(\bb)^{1\over n})^{n-1})\right )\cr
&=\rho_{\m}\cdot x+O(x^{1-{1\over n}})\cr}$$
o $\rho_{\m}={{\rm Vol}(\Gamma)\over w_{\m}}\cdot {\N(\bb)\over {\rm Vol}(\Lambda_{0})}$ ne dŽpend que de $\m$. En effet, $\Gamma$ et $w_{\m}$ ne dŽpendent que de $\m$ de manire Žvidente; l'idŽal $\bb$ est dans la classe ${\euf K}^{-1}$ et  dŽpend bien sžr de $\euf K$, mais ${\rm Vol}(\Lambda_{0})={\rm Vol}(\Lambda)={\rm Vol}(v(\bb\m_{0}))=\N(\bb)\cdot\N(\m_{0})\cdot{\rm Vol}(v(O_{K}))$. Ainsi,  ${\N(\bb)\over {\rm Vol}(\Lambda_{0})}$ ne dŽpend aussi que de $\m$.\qed


\vfill\eject
\global\advance\chapnomb by 1
\nomb=1

\centerline{\para Chapitre 2 :  SŽries de Dirichlet et}
\medskip
\centerline{\para premire inŽgalitŽ du corps de classe }

\headline{\hfill\phantom{ouuh}\hfill}

\bigskip
Dans ce chapitre, nous ferons un peu d'analyse complexe pour prŽparer le chapitre suivant o nous dŽmontrerons le thŽorme de  $\check{\rm C}$ebotarev. Le but est aussi de prouver la premire inŽgalitŽ du corps de classe (ThŽorme \argau)
\bigskip
\defi

Tout le monde conna"t l'exponentielle complexe $e^z=\sum_{n=0}^\infty {z^n\over n!}$ qui converge dans
tout le plan complexe. Si $t>0$ est un nombre rŽel strictement positif, on pose $t^s=e^{s\log(t)}$ o 
$\log (t)$ est le logarithme usuel des nombres rŽels. Une {\it sŽrie de Dirichlet} est une fonction
complexe du type 

$$\sum_{n=1}^\infty {a(n)\over n^s},$$ 
avec $a(n),s\in\C$.
\bigskip

Soit  $b\geq 0$, $\delta>0$ et $0<\varepsilon<{\pi\over 2}$. On dŽfinit aussi 

$$D(b,\delta,\varepsilon)=\{s\in\C\mid \Re(s)\geq b+\delta\hbox{ et }|\arg(s-b)|\leq {\pi\over
2}-\varepsilon\}.$$
\bigskip

Graphiquement, cela donne ceci~:

\vskip4cm
\psset{xunit=0.75cm,yunit=0.75cm}
\pspicture(-5,-4)(15,0)
\psline{->}(0,-2)(0,5)

\psline(2,5)(2,-5)\rput(1.7,0.4){$b$}
\psline(3.5,4)(3.5,-4)
\psline[linestyle=dashed](3.5,1)(2,0)(3.5,-1)
\psline[fillstyle=solid, fillcolor=lightgray](8,4)(3.5,1)(3.5,-1)(8,-4)
\psline{->}(-2,0)(10,0)
\rput(7,1.4){$D(b,\delta,\varepsilon)$}
\psline{->}(2.6,-3)(2,-3)\psline{->}(3,-3)(3.5,-3)\rput(2.8,-3){$\delta$}
\psarc[arcsepB=2pt]{->}(3.5,1){1}{33}{90} 
\rput(4.4,2.3){$\varepsilon$}
\endpspicture
 
\vskip1cm
\th
\medskip
Soit $\sum_{n=1}^\infty {a(n)\over n^s}$ une sŽrie de Dirichlet. Posons $s(n)=\sum_{i=1}^n a(i)$. On
suppose qu'il existe $a>0$ et $b\geq 0$ tels que  $|s(n)|\leq a\cdot n^b$ pour tout $n\geq 1$. Alors la
sŽrie converge uniformŽment dans $D(b,\delta,\varepsilon)$ pour tout $\delta,\varepsilon$ tel que
$\delta>0$ et $0<\varepsilon<{\pi\over 2}$. Cela implique que cette sŽrie dŽfinit une fonction holomorphe
sur le demi-plan $\Re(s)>b$ (cf. [Con, Th. 2.1, p.147]).

{\bf Preuve}

On posera $\sigma=\Re(s)$. Donc $|t^s|=t^\sigma$ si $t>0$. On va utiliser le critre de convergence de
Cauchy. Soit donc $u$ et $v$ des entiers tels que $v\geq u+1>2$; et $s\in D(b,\delta,\varepsilon)$. On a 
$$\eqalign{\left |\sum_{n=u}^v {a(n)\over n^s}\right |&=\left |\sum_{n=u}^v{s(n)-s(n-1)\over n^s}\right
|=\left|\sum_{n=u}^v{s(n)\over n^s}-\sum_{n=u-1}^{v-1}{s(n)\over (n+1)^s}\right |\cr
&=\left |{s(v)\over v^s}-{s(u-1)\over u^s}+\sum_{n=u}^{v-1}s(n)\left ({1\over n^s}-{1\over (n+1)^s}\right
)\right |\cr
&\leq \left | {s(v)\over v^s}\right |+\left | {s(u-1)\over u^s}\right |+\sum_{n=u}^{v-1}|s(n)|\cdot\left
|{1\over n^s}-{1\over (n+1)^s}\right |\cr}$$
Or, ${1\over n^s}-{1\over (n+1)^s}=s\cdot\int_{n}^{n+1}{dt\over t^{s+1}}$, (on se souvient que $\int
(g+ih)\;dt =\int g\; dt+i\int h\; dt$ si le domaine d'intŽgration est rŽel). Ainsi, 

$$\eqalign{\left |\sum_{n=u}^v {a(n)\over n^s}\right |&\leq {a\over v^{\sigma-b}}+{a\over
u^{\sigma-b}}+\sum_{n=u}^{v-1}|s|\cdot a\cdot n^b\cdot\left |\int_{n}^{n+1}{dt\over t^{s+1}}\right |\cr
&\leq {2a\over u^{\sigma-b}}+|s|\cdot a\cdot \sum_{n=u}^{v-1}\left |\int_{n}^{n+1}{t^b\cdot dt\over
t^{s+1}}\right |\cr
&\leq {2a\over u^{\sigma-b}}+|s|\cdot a\cdot\int_{u}^{v}{ dt\over
t^{\sigma+1-b}}\cr
&\leq {2a\over u^{\sigma-b}}+|s|\cdot a\cdot\int_{u}^{\infty}{ dt\over
t^{\sigma+1-b}}={2a\over u^{\sigma-b}}+{|s|\cdot a\over (\sigma-b) u^{\sigma-b}}.\cr}$$
Or, puisque $s$ est dans $D$, on a ${|s|\over \sigma-b}\leq {|s-b|+b\over \sigma-b}\leq {1\over\cos
(\theta)}+{b\over\delta}$, o $\theta=\arg(s-b)$. Puisque $|\theta|\leq {\pi\over 2}-\varepsilon$, on voit facilement que 
${1\over\cos(\theta)}\leq M=:{1\over\cos({\pi\over 2}-\varepsilon)}$. Ainsi, 

$$\left |\sum_{n=u}^v {a(n)\over n^s}\right |\leq {2a+aM+{ab\over\delta}\over
u^{\sigma-b}}<{2a+aM+{ab\over\delta}\over u^\delta}\buildrel {\rm
uniform\acute ement}\over\lra 0\hbox{ lorsque } u\rightarrow \infty.$$
\qed

\bigskip
\newcount\gaagm\gaagm=\gaga

\headline{\hfill\smcap  SŽries de Dirichlet et premire inŽgalitŽ du corps de classe\hfill}

\th
\medskip
{\sl La fonction $\zeta(s)=\sum_{n=1}^\infty {1\over n^s}$ pour $\Re(s)>1$ se prolonge en une fonction
mŽromorphe, encore notŽe $\zeta(s)$ sur le demi-plan $\Re(s)>0$. Il y a un seul p™le en $s=1$ qui est
simple et de rŽsidu 1.

}

{\bf preuve}

En vertu du thŽorme prŽcŽdent (pour $s(n)=n$, et $a=b=1$), on trouve que $\zeta(s)$ converge dans le
demi-plan $\Re(s)>1$ et y dŽfinit une fonction holomorphe. ConsidŽrons le sŽrie
$\zeta_2(s)=\sum_{n=1}^\infty {(-1)^{n+1}\over n^s}$. On voit que $|s(n)|\leq 1$. Par le thŽorme
prŽcŽdent, avec
$a=1$ et $b=0$, on a que $\zeta_2(s)$ est holomorphe sur le demi-plan $\Re(s)>0$. Lorsque $\Re(s)>1$, la
convergence est absolue. Dans ce domaine, on trouve 

$$\eqalign{\zeta_2(s)&={1\over 1^s}-{1\over 2^s}+{1\over 3^s}-{1\over 4^s}+\cdots={1\over 1^s}+{1\over
2^s}-{2\over 2^s}+{1\over 3^s}+{1\over 4^s}-{2\over 4^s}+{1\over 5^s}+{1\over 6^s}-{2\over
6^s}+\cdots\cr
&({1\over 1^s}+{1\over 2^s}+{1\over 3^s}+{1\over 4^s}+\cdots )-2({1\over 2^s}+{1\over 4^s}+{1\over
6^s}+\cdots )=\zeta(s)-{2\over 2^{s}}\cdot \zeta(s)=(1-{1\over 2^{s-1}})\cdot\zeta(s),\cr}$$
ainsi $\zeta(s)={\zeta_2(s)\over (1-{1\over 2^{s-1}})}$ lorsque $\Re(s)>1$. Mais par ce qui prŽcde, on a
que la fonction ${\zeta_2(s)\over (1-{1\over 2^{s-1}})}$ est mŽromorphe sur le demi-plan $\Re(s)>0$ dont
les p™les (s'il y en a) sont parmi les zŽros de $(1-{1\over 2^{s-1}})$, donc de la forme $1+{2k\pi\over
\log(2)}\cdot i$, $ k\in\Z$.  On considre ensuite $\zeta_3(s)={1\over 1^s}+{1\over 2^s}-{2\over
3^s}+{1\over 4^s}+{1\over 5^s}-{2\over 6^s}+\cdots$. Si cette dŽfinition ne vous convient pas, vous
pouvez poser $a(n)=1-\sum_{j=0}^2 \zeta^{n\cdot j}$, avec $\zeta=\zeta_3=e^{2i\pi\over 3}$ et
$\zeta_3(s)=\sum_{n=1}^\infty {a(n)\over n^s}$. Evidemment, il y aura toujours collusion de notation
entre les $\zeta_n$ qui sont les racines de l'unitŽ et la fonction $\zeta$ de Riemann, mais le contexte
permettra toujours de comprendre de quoi on parle. Le thŽorme prŽcŽdent s'applique ˆ nouveau et donc,
la fonction $\zeta_3(s)$ est une fonction holomorphe sur $\Re(s)>0$; et lorsque $\Re(s)>1$, on a 
$$\zeta_3(s)={1\over 1^s}+{1\over 2^s}+{1\over 3^s}-{3\over 3^s}+{1\over 4^s}+{1\over 5^s}+{1\over
6^s}-{3\over 6^s}+\cdots=\zeta(s)-{3\over 3^s}\zeta(s)=(1-{1\over 3^{s-1}})\zeta(s).$$
Donc ${\zeta_3(s)\over (1-{1\over 3^{s-1}})}$ est un prolongement mŽromorphe de $\zeta$ sur $\Re(s)>0$.
Donc les p™les sont parmi les nombres $1+{2l\pi\over \log (3)}\cdot i$, $l\in\Z$. Par unicitŽ du
prolongement analytique (cf.  [Con, Cor. 3.8, p. 70]) on a ${\zeta_2(s)\over (1-{1\over
3^{s-1}})}={\zeta_3(s)\over (1-{1\over 3^{s-1}})}$ en dehors de leurs p™les. De plus, puisque ces fonctions
sont dŽfinies sur des disques entourants ces p™les, si $s$ est un p™le de l'une de ces fonctions, il doit
tre aussi un p™le de l'autre. Donc, $s=1+{2k\pi\over \log (2)}\cdot i=1+{2l\pi\over \log (3)}\cdot i$ 
pour certains $k$ et $l\in \Z$; donc ${k\over\log(2)}= {l\over\log(3)}$, ou encore $2^l=3^k$; donc $k=l=0$.
On en dŽduit que le seul p™le Žventuel est $s=1$; et c'est bien un p™le car $\sum_{n=1}^\infty {1\over
n^\sigma}\rightarrow \infty$ si $\sigma\rightarrow  1$ par valeurs $>1$. Ce p™le est simple, car les
zŽros de $(1-{1\over 2^{s-1}})$ sont simples (la dŽrivŽe qui vaut ${\log(2)\over 2^{s-1}}$ ne s'annule
pas en $s=1$).
\goodbreak

Calculons le rŽsidu en $s=1$. On veut donc calculer $\lim_{s\rightarrow 1}(s-1)\zeta(s)$. Pour ce calcul,
il suffit de prendre des $s$ rŽels $>1$. On considre le graphe de $y=x^{-s}$.

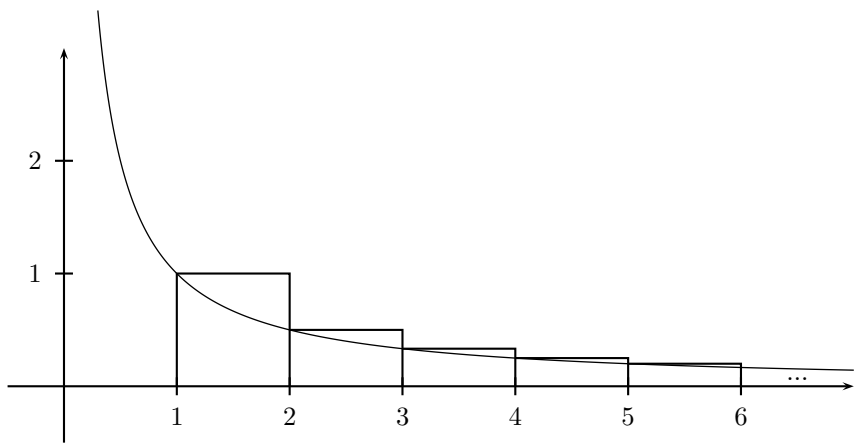
\vglue 5cm 

\psset{xunit=1.5cm,yunit=1.5cm}
\psaxes[Dx=1,Dy=1]{->}(0,0)(-0.5,-0.5)(7,3)
\psplot[plotpoints=1000,linecolor=black,linewidth=0.5pt]{0.3}{7}{ 1  x  div }
\psline(1,0)(1,1)(2,1)(2,0)
\psline(2,0.5)(3,0.5)(3,0)
\psline(3,0.3333)(4,0.3333)(4,0)
\psline(4,0.25)(5,0.25)(5,0)
\psline(5,0.2)(6,0.2)(6,0)
\rput(6.5,0.07){...}

\vskip 1cm
La somme de aires des rectangles donne Žvidemment $\zeta(s)$ et on voit que ${1\over s-1}=\int_{1}^\infty
{dx\over x^s}\leq \zeta(s)$. Cet autre dessin 

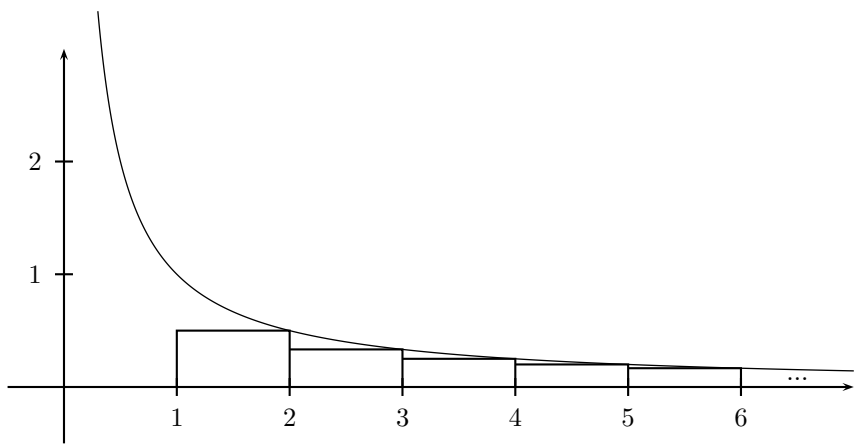
\vglue 6cm 

\psset{xunit=1.5cm,yunit=1.5cm}
\psaxes[Dx=1,Dy=1]{->}(0,0)(-0.5,-0.5)(7,3)
\psplot[plotpoints=1000,linecolor=black,linewidth=0.5pt]{0.3}{7}{ 1  x  div }
\psline(1,0)(1,0.5)(2,0.5)(2,0)
\psline(2,0.3333)(3,0.3333)(3,0)
\psline(3,0.25)(4,0.25)(4,0)
\psline(4,0.2)(5,0.2)(5,0)
\psline(5,0.1666)(6,0.1666)(6,0)
\rput(6.5,0.07){...}

\vskip 1cm
 montre que $\zeta(s)-1\leq \int_{1}^\infty {dx\over x^s}={1\over s-1}$. On en dŽduit que 

$${1\over s-1}\leq \zeta(s)\leq {s\over s-1},$$
ainsi, $1\leq (s-1)\zeta(s)\leq s$ et donc $\lim_{s\rightarrow 1}(s-1)\zeta(s)=1$.\qed

\bigskip
{\soustitre Rappel}

{\sl
Soit $K$ un corps de nombres et $\m$ un $K$-module. Alors il existe une constante $\rho_{\m}>0$ telle que
si $\euf K$ est une classe de $I(\m)$ modulo $P_\m$, et si $j(x,{\euf K})$ est le nombre d'idŽaux ${\euf
a}\in{\euf K}$ tel que $\N({\euf a})\geq x$. Alors 
$$j(x,{\euf K})=\rho_{\m}\cdot x+O(x^{1-{1\over d}})$$
o $d=[K:\Q]$.

}

\medskip
C'est le ThŽorme \argq\ qui Žtait l'objet du Chapitre 1.
\bigskip
\newcount\gaagn\gaagn=\gaga

\defi

Soit $K$ un corps de nombres, $\m$ un $K$-module et $\euf K$ une classe de $I(\m)$ modulo $P_\m$.
On dŽfinit 
$$\zeta_\m(s,{\euf K})=\sum_{{\euf a}\in{\euf K}}{1\over\N({\euf a})^s}.$$
Evidemment, cette dŽfinition n'a de sens que si le convergence de cette sŽrie est absolue. C'est ce que
nous allons voir incessamment.
\bigskip\goodbreak
\th
{\sl 

La sŽrie $\zeta_\m(s,{\euf K})$ converge absolument sur le demi-plan $\Re(s)>1$ et se prolonge en une
fonction mŽromorphe sur $\Re(s)>1-{1\over d}$ ($d=[K:\Q]$) avec un seul p™le en $s=1$; ce p™le est
simple de rŽsidu $\rho_\m$. 

}

{\bf Preuve}

ConsidŽrons la sŽrie de Dirichlet $\sum_{n=1}^\infty {a(n)\over n^s}$, avec $a(n)$, le nombre d'idŽaux
${\euf a}$ dans $\euf K$ tels que $\N({\euf a})=n$. Dans ce cas, $s(n)=\sum_{i=1}^n a(i)=j(n,{\euf K})$.
Par le ThŽorme \argq, ${s(n)\over n}$ est bornŽ, pour $n\geq 1$. Donc il existe $A>0$ tel que
$|s(n)|=s(n)\leq A\cdot n$.  Le ThŽorme \argac\  s'applique et dit que la sŽrie $\sum_{n=1}^\infty {a(n)\over n^s}$ dŽfinit une fonction holomorphe sur le demi-plan $\Re (s)>1$. Clairement, la convergence dans ce demi-plan est absolue, donc l'ordre des termes n'importe pas et cela justifie qu'on Žcrive simplement
$\sum_{{\euf a}\in{\euf K}}{1\over\N({\euf a})^s}$. 

ConsidŽrons maintenant la sŽrie $\sum_{n=1}^\infty
{b(n)\over n^s}$, avec $b(n)=a(n)-\rho_\m$ (le $a(n)$ Žtant celui de tout ˆ l'heure). Pour cette nouvelle
sŽrie, on a $s(n)=\sum_{i=1}^n b(i)=j(n,{\euf K})-\rho_\m\cdot n$. Par le ThŽorme \argq, il existe
$B>0$ tel que $s(n)\leq B\cdot n^{1-{1\over d}}$. A nouveau par le ThŽorme \argac, la sŽrie dŽfinit une
fonction holomorphe (notons-la $f$) sur le demi-plan $\Re(s)>1-{1\over d}$. Lorsque $\Re(s)>1$, on voit
que $f(s)+\rho_\m\cdot\zeta(s)=\sum_{n=1}^\infty {a(n)\over n^s}=\zeta_\m (s,{\euf K})$. Donc, en vertu
du ThŽorme \argad, $\zeta_\m (s,{\euf K})$ se prolonge en une fonction mŽromorphe sur $\Re(s)>1-{1\over d}$,
avec un seul p™le en $s=1$ qui est simple et de rŽsidu $\rho_\m$.\qed

\bigskip

\defis

\art{a)} Soit $G=(G,*,1)$ un groupe abŽlien fini. On appelle {\it caractre} de $G$ tout homomorphisme
$\chi \,  : G\, \lra \C^*$. On note ${\bf 1}$ le caractre de $G$ qui envoie tout ŽlŽment de $G$ sur $1$.
L'ensemble $\widehat{G}$ des caractres est lui-mme un groupe isomorphe ˆ $G$~:
$\chi_1\chi_2(g):=\chi_1(g)\cdot\chi_2(g)$. Remarquons que l'on a pour tout $\chi\in \widehat{G}$~:

$$\sum_{g\in G}\chi(g)=\cases{0& si $\chi\ne {\bf 1}$\cr |G|& si $\chi= {\bf 1}.$\cr}$$

\art{}En effet, c'est clair si $\chi={\bf 1}$. Supposons que $\chi\ne {\bf 1}$, donc il existe $h\in G$ tel que
$\chi(h)\ne 1$. L'application $g\mapsto hg$ est clairement une bijection de $G$ dans lui-mme. Donc, 

$$\sum_{g\in G}\chi(g)=\sum_{g\in G}\chi(hg)=\chi(h)\cdot \sum_{g\in G}\chi(g).$$

Si $\sum_{g\in G}\chi(g)\ne 0$, on en dŽduirait que $\chi(h)=1$, ce qui est contradiction. 

De plus, puisque $G$ et  $\widehat{G}$ sont isomorphes, on a aussi pour tout $g\in G$

$$\sum_{\chi\in \widehat{G}}\chi(g)=\cases{0& si $g\ne 1$\cr |G|& si $g= 1.$\cr}$$

\newcount\gaago\gaago=\gaga

\art{b)}Soit $K$ un corps de nombres, $\m$ un $K$-module et $\chi$ un caractre de $I(\m)/P_\m$. On dŽfinit 

$$L_\m(s,\chi)=\sum_{{\euf a}\in I(\m)}{\chi({\euf a})\over \N({\euf a})^s},$$

\art{}o $\chi({\euf a}):=\chi({\euf a}\bmod {P_\m})$. Comme tout ˆ l'heure, cette sŽrie n'est dŽfinie que si
la convergence est absolue. C'est Žvidemment ce qu'on va voir !
\bigskip
\th

{\sl La sŽrie $L_\m(s,\chi)$ converge absolument sur le demi-plan $\Re(s)>1$ et admet un
prolongement mŽromorphe sur $\Re(s)>1-{1\over d}$. Si $\chi\ne {\bf 1}$ (le caractre unitŽ), il
n'y a pas de p™le (donc le prolongement est holomorphe). Si $\chi= {\bf 1}$, il y a un p™le unique, simple en $s=1$
et de rŽsidu $h_\m\rho_\m$, (souvenons nous de $h_\m$, c'est le cardinal du groupe $I(\m)/P_\m$).

}

{\bf Preuve}

La convergence est absolue sur $\Re(s)>1$, car $\left |{\chi({\euf a})\over \N({\euf a})^s}\right |={1\over
\N({\euf a})^\sigma}$, o $\sigma=\Re(s)$ et parce que $\sum_{{\euf a}\in I(\m)}{1\over \N({\euf
a})^\sigma}$ $=\sum_{{\euf K}\in I(\m)/P_\m}\zeta_\m(\sigma,{\euf K})$; cette somme est finie, car
$I(\m)/P_\m$ est fini de cardinal $h_\m$ (cf. ThŽorme \argl). Toujours si $\Re(s)>1$, en rŽarrangeant les termes, on obtient

$$L_\m(s,\chi)=\sum_{{\euf K}\in I(\m)/P_\m}\chi({\euf K})\cdot\zeta_\m(s,{\euf K}).$$
On conclut en vertu du ThŽorme \argaf\ et en se souvenant que 

$$\sum_{{\euf K}\in I(\m)/P_\m}\chi({\euf K})=\cases{h_\m&si $\chi= {\bf 1}$\cr 0&sinon.\cr} $$
\qed
\bigskip

\defis

\art{a)}On notera $\zeta_\m(s)$ au lieu de $L_\m(s,{\bf 1})=\sum_{{\euf K}\in I(\m)/P_\m}\zeta_\m(s,{\euf K})$. On l'appellera la {\it fonction zta de $\m$}. Et si $\m={\euf 1}=O_K\cdot \emptyset$, on note cette fonction $\zeta_K(s)$ qui est {\it la fonction zta} de $K$.

\art{b)}Soit $z\in\C$. Posons $\Log(z)=\log|z|+i\arg (z)$ o $-\pi<\arg(z)\leq\pi$, o $\log(\cdot )$ est la fonction logarithme usuelle sur les nombres rŽels. On appelle $\Log(\cdot)$ la {\it
branche principale du logarithme}. On peut voir que $e^{\Log(z)}=z$, pour tout $z\in\C^*$. Si $w\in\C$ est tel que $e^{w}=z$, alors il existe $k\in\Z$ tel que $w=\Log(z)+2k\pi i$. 
De plus, si $|z|<1$, alors $-\Log(1-z)=\Log((1-z)^{-1})=\sum_{n=1}^\infty {z^n\over n}$.

\art{c)}Si $(a_n)_{n=1}^\infty$ est une suite de nombres complexes non nuls, nous dirons que $\prod_{n=1}^\infty a_n$ {\it converge} si $\lim_{N\rightarrow\infty}\prod_{n=1}^N a_n$ existe et est {\bf non nulle}. Nous utiliserons le rŽsultat suivant~:  si la somme $\sum_{n=1}^\infty \Log (a_n)$ converge absolument, alors $\prod_{n=1}^\infty a_n$ converge. Nous dirons
dans ce cas que le produit $\prod_{n=1}^\infty a_n$ {\it converge absolument}.

\bigskip\goodbreak

\th

{\sl Soit $K$, $\m$ et $\chi$ comme avant. Alors si $\Re(s)>1$, on a la
reprŽsentation 
$$L_\m(s,\chi)=\prod_{\P\notdivi\m}\left (1-{\chi(\P)\over \N(\P)^s}\right )^{-1},$$
le produit convergeant absolument. Cela prouve que $L_\m (s,\chi)\ne 0$.

Posons $\llog L_\m (s,\chi)=\sum_{\P\notdivi
\m}\Log\left(\left(1-{\chi(\P)\over\N(\P)^s}\right)^{-1}\right)$. Alors $\llog L_\m (s,\chi)$ est holomorphe sur
$\Re(s)>1$ et il existe $g_\chi(s)$ une fonction holomorphe sur le demi-plan $\Re(s)>{1}$, prolongeable
holomorphiquement au voisinage de $1$ telle que 
$$\llog L_\m (s,\chi)=\sum_{\P\notdivi\m}{\chi(\P)\over \N(\P)^s}+g_\chi(s),$$
pour tout $s\in\C$ tel que $\Re(s)>1$. 

Voir le lemme suivant pour une explication de l'Žcriture $\llog L_\m (s,\chi)$. Remarquons encore qu'ici (et pour ce chapitre), l'Žcriture $\P\notdiv\m$ ne concerne que les places finies.

}

{\bf Preuve}

Tout d'abord, la sŽrie $\sum_{\P\not\hskip 1.9pt |\hskip 0.9 pt\m}{\chi(\P)\over \N(\P)^s}$ converge
absolument sur $\Re(s)>1$, car c'est une partie de $\sum_{{\euf a}\in I(\m)} {\chi(\aa)\over \N(\aa)^s}$,
qui converge dans ce domaine en vertu du ThŽorme \argah. D'autre part, $\sum_{\P\notdivi \m}\sum_{j=2}^\infty
{\chi(\P^j)\over j\cdot \N(\P)^{sj}}$ converge absolument si $\sigma:=\Re(s)>{1\over 2}$, car 

$$\eqalign{\sum_{\P\notdivi \m}\sum_{j=2}^\infty{1\over j\cdot \N(\P)^{\sigma j}}&\leq  \sum_{\P\notdivi
\m}{1\over 2}\cdot {1\over \N(\P)^{2\sigma}}\cdot \sum_{j=0}^\infty {1\over \N(\P)^{\sigma
j}}={1\over 2}\cdot \sum_{\P\notdivi \m}{1\over \N(\P)^{2\sigma}}\cdot {1\over 1-{1\over\N(\P)^\sigma}}\cr
&={1\over 2}\cdot \sum_{\P\notdivi \m}{1\over \N(\P)^{\sigma}}\cdot{1\over \N(\P)^\sigma-1}\leq {1\over 2}\cdot
\sum_{\P\notdivi \m}{1\over \N(\P)^{\sigma}}\cdot{4\over\N(\P)^\sigma}.\cr}$$
La dernire inŽgalitŽ vient de $4(\N(\P)^\sigma-1)=\N(\P)^\sigma+3\N(\P)^\sigma-4\geq \N(\P)^\sigma
+3\sqrt{2}-4\geq \N(\P)^\sigma$. On a donc prouvŽ que 

$$\sum_{\P\notdivi \m}\sum_{j=2}^\infty{1\over j\cdot \N(\P)^{\sigma j}}\leq 2\sum_{\P\notdivi
\m}{1\over \N(\P)^{2\sigma}}\leq 2\zeta_\m(2\sigma).$$
Et donc $g_\chi(s):=\sum_{\P\notdivi
\m}\sum_{j=2}^\infty{\chi(\P^j)\over j\cdot \N(\P)^{s j}}$ est une fonction holomorphe sur $\Re(s)>{1\over 2}$, car elle converge absolument et uniformŽment sur tout compact de ce domaine. Maintenant, en appliquant la sŽrie de la branche principale du logarithme, avec $z={\chi(\P)\over\N(\P)^{s}}$ (qui est de module $<1$ si $\Re(s)>1$), on trouve que

$$\sum_{\P\not\hskip 1.9pt |\hskip 0.9 pt\m}{\chi(\P)\over \N(\P)^s}+g_\chi(s)=\sum_{\P\notdivi
\m}\sum_{j=1}^\infty{\chi(\P^j)\over j\cdot \N(\P)^{s j}}=\sum_{\P\notdivi
\m}\Log\left(\left(1-{\chi(\P)\over\N(\P)^s}\right)^{-1}\right)=\llog  L_\m (s,\chi)$$
qui converge absolument et est holomorphe si $\Re(s)>1$. En choisissant une numŽrotation quelconque des $\P\notdiv
\m$ (c'est possible de le faire maintenant qu'on a la convergence absolue), on a

$$e^{\llog  L_\m (s,\chi)}=\lim_{n\rightarrow\infty}e^{\sum_{i=1}^n\Log\left (\left(
1-{\chi(\P_i)\over\N(\P_i)^s}\right )^{-1}\right )}=\lim_{n\rightarrow\infty}\prod_{i=1}^n\left (1-{\chi(\P_i)\over
\N(\P_i)^s}\right)^{-1}=\prod_{\P\notdivi m}\left (1-{\chi(\P)\over
\N(\P)^s}\right)^{-1}.$$
Il reste ˆ vŽrifier que le dernier produit est Žgal ˆ $L_\m(s,\chi)$. Il suffit de montrer que 
$$L_\m(s,\chi)=\lim_{N\rightarrow\infty}\prod_{\eqalign{\noalign{\vskip-4pt}\scriptstyle
\P\notdivi\m\phantom{88}\cr
\noalign{\vskip-5pt} \scriptstyle\N(\P)\leq N}} \left (1-{\chi(\P)\over
\N(\P)^s}\right)^{-1}.$$
Or, 

$$\prod_{\eqalign{\noalign{\vskip-4pt}\scriptstyle
\P\notdivi\m\phantom{88}\cr \noalign{\vskip-5pt} \scriptstyle\N(\P)\leq N}} \left (1-{\chi(\P)\over
\N(\P)^s}\right)^{-1}=\prod_{\eqalign{\noalign{\vskip-4pt}\scriptstyle
\P\notdivi\m\phantom{88}\cr \noalign{\vskip-5pt} \scriptstyle\N(\P)\leq N}} \left (\sum_{j=0}^\infty
{\chi(\P)^j\over \N(\P)^{sj}}\right)=\sum_{{\euf a}\in S_N}{\chi({\euf a})\over\N({\euf a})^s},$$
o $S_N$ est l'ensemble des idŽaux entiers, premiers ˆ $\m$ dont les diviseurs premiers sont de norme infŽrieure ou
Žgale ˆ $N$; la dernire ŽgalitŽ venant de la convergence absolue de la somme $\sum_{j=0}^\infty
{\chi(\P)^j\over \N(\P)^{sj}}$ et de la multiplicativitŽ de $\chi$ et de la norme. Ainsi,

$$\left | \prod_{\eqalign{\noalign{\vskip-4pt}\scriptstyle
\P\notdivi\m\phantom{88}\cr \noalign{\vskip-5pt} \scriptstyle\N(\P)\leq N}} \left (1-{\chi(\P)\over
\N(\P)^s}\right)^{-1}-L_\m(s,\chi)\right |\leq \sum_{\eqalign{\noalign{\vskip-2pt}\scriptstyle
\N({\euf a})>N\phantom{8}\cr \noalign{\vskip-8pt} \scriptstyle\aa\in I{(\m)}}}{1\over{\N({\euf a})^{\Re(s)}}}\longrightarrow
0\hbox{ si $N\rightarrow \infty,$} $$
en vertu du ThŽorme \argah.\qed

\bigskip
\goodbreak
\newcount\gaagp\gaagp=\gaga

\lem

{\sl Si $X\subset \C$ est un ouvert connexe et simplement connexe et si $f\ :\ X\longrightarrow \C^*$ est une fonction holomorphe et non nulle, alors il existe, sur $X$, un logarithme de cette fonction, c'est-ˆ-dire une fonction holomorphe $\widetilde{f}$ telle que $e^{\widetilde{f}(x)}=f(x)$. De plus, si on fixe une valeur de $\widetilde{f}$ en un point de $X$, alors cette fonction est unique. Remarquons que si $f(x)=f(y)$, cela n'implique
pas forcŽment que $\widetilde{f}(x)=\widetilde{f}(y)$. Lors du thŽorme prŽcŽdent, la fonction $\llog  L_\m (s,\chi)$ Žtait un logarithme de la fonction $L_\m (s,\chi)$, voilˆ pourquoi, nous avons Žcrit $\llog  L_\m (s,\chi)$ plut™t que $\llog (L_\m (s,\chi))$. 

}

{\bf Preuve}

Cf. [Ru, Thm 13.18, p. 263]

\bigskip 
\coro

{\sl Soit $K$, $\m$ comme avant. Il existe $h(s)$ une fonction holomorphe dans un voisinage $V$ de $s=1$ et
$g(s)$ une fonction holomorphe sur $\Re(s)>{1\over 2}$ telles que si $s$ est dans $V$ et que
$\Re(s)>1$, alors on a~:

$$\llog\zeta_\m(s)=-\Log(s-1)+h(s)=\sum_{\P\notdivi\m}{1\over\N(\P)^s}+g(s)\eqno{(1)}$$ 

}

{\bf Preuve}

Soit $\chi$ un caractre de $I(\m)/P_\m$. Si $\Re(s)>1$, alors il est clair que $\Log(s-1)+\llog L_\m(s,\chi)$
est un logarithme de $(s-1)\cdot L_\m(s,\chi)$. Si on prend $\chi={\bf 1}$, alors, vu les dŽfinitions qu'on a
faites et en vertu du ThŽorme \argah, la fonction $(s-1)\cdot\zeta_\m(s)$ est dŽfinie, holomorphe et non nulle sur
un voisinage de $s=1$. Donc, cette fonction possde un logarithme sur ce voisinage, en vertu du lemme
prŽcŽdent. Notons $h(s)$ le logarithme qui co•ncide avec $\Log(s-1)+\llog\zeta_\m(s)$ sur $\Re(s)>1$.
On en dŽduit la premire ŽgalitŽ. La seconde dŽcoule du ThŽorme \argaj.\qed 

\bigskip

\defis

Soient $f_1$ et $f_2$ des fonctions complexes dŽfinies sur $\Re(s)>1$. On Žcrira $f_1\sim f_2$ si
$g(s):=f_1(s)-f_2(s)$ est une fonction holomorphe sur $\Re(s)>1$ et prolongeable en une fonction toujours
holomorphe au voisinage de $s=1$.

Soit $K$ un corps de nombres et $S$ un ensemble d'idŽaux premiers de $O_K$. On dira que $S$ est de {\it densitŽ
de Dirichlet $\delta$} si l'une des conditions suivantes est satisfaite~:

$$\leqalignno{ \lim_{s\rightarrow 1^+}{\sum_{\P\in S}{1\over \N(\P)^s}\over -\Log(s-1)} &=\delta& 1) \cr
\sum_{\P\in S}{1\over
\N(\P)^s}\sim -\delta\cdot\Log(s-&1)&2)\cr}$$

Il est clair que la condition 2) implique la condition 1). Si $S$ satisfait la
condition 2), on dit de plus que $S$ est {\it rŽgulier}. Le Corollaire \argal\ implique que l'ensemble des
idŽaux premiers de $I(\m)$ est  rŽgulier de densitŽ 1. Il est clair que si $S$ est fini, sa densitŽ est nulle et que si $S$ et $S'$ sont disjoints de densitŽ $\delta$ et $\delta'$, alors la densitŽ de $S\cup S'$ vaut $\delta+\delta'$. Cela implique, puisque l'ensemble de idŽaux premiers $I(\m)$ contient tous les idŽaux premiers de $K$ sauf un nombre fini, que si $S$ est de densitŽ $\delta$, on a $0\leq\delta\leq 1$ et que l'ensemble des idŽaux premiers de $K$ est rŽgulier de densitŽ 1. On montre aussi facilement que si $S\subset S'$ sont de densitŽ $\delta$ et $\delta'$, alors $\delta\leq\delta'$. 

\bigskip
\newcount\gaagq\gaagq=\gaga

\goodbreak

\lem

{\sl  Soit $K$ un corps de nombres. Alors l'ensemble des idŽaux premiers $\P$ tels que $f(\P/\Q)>1$ est
un ensemble rŽgulier de densitŽ de Dirichlet nulle.

}

{\bf Preuve}

Notons $A$ l'ensemble de ces idŽaux. On va montrer que $\sum_{\P\in A}{1\over \N(\P)^s}$ est holomorphe
au voisinage de $s=1$. Pour tout $\P\in A$, la fonction $s\mapsto {1\over \N(\P)^s}$ est
clairement holomorphe au voisinage de $s=1$. D'autre part, puisque $\P\in A$, il existe $p\in \gfP_{0}(\Q)$ tel que $\N(\P)=p^{f(\P,\Q)}\geq p^2$. De plus, il y a au plus $[K:\Q]$ nombres premiers au-dessous de $\P$. Ainsi 

$$\sum_{\P\in A}{1\over \N(\P)^\sigma}\leq [K:\Q]\cdot \sum_{p\in \gfP_0(\Q)}{1\over
p^{2 \sigma}}\leq [K:\Q]\cdot\zeta(2 \sigma),\hbox{ o }\sigma=\Re(s)>{1\over 2}.$$

Ce qui montre , en vertu du ThŽorme \argad, que $\sum_{\P\in A}{1\over \N(\P)^s}$ converge absolument et uniformŽment sur tout compact du demi-plan  $\Re(s)>{1\over 2}$, donc est une fonction holomorphe au voisinage de $1$, cela prouve le lemme.\qed

\bigskip\goodbreak

\lem

{\sl Soit $L/K$ une extension de corps de nombres. Alors l'ensemble $S=\{\P\in \gfP_{0}(K)\mid \P$ se
dŽcompose compltement dans $L\}$ est rŽgulier de densitŽ $[E:K]^{-1}$ o $E$ est la cl™ture galoisienne de
$L/K$, c'est-ˆ-dire la plus petite extension $E/K$, galoisienne, contenant $L$.

}

{\bf Preuve}

Si $\P\in \gfP(K)$ se dŽcompose compltement dans $L$, alors, en vertu du lemme ``dŽcomposition-ramification'' du Chapitre 0, page \the\gaagc, il se dŽcompose
compltement dans $E$. On a donc $f(\gP/\P)=1$ pour tout idŽal premier $\gP$ de $E$ au-dessus de $\P$ (il y en a donc $[E:K]$ et $\N(\P)=\N(\gP)$). On trouve alors~:
$$[E:K]\cdot \sum_{\P\in S}{1\over \N(\P)^s}=\sum_{\gP\in S_1}{1\over \N(\gP)^s},$$
o $S_1=\{\gP\in \gfP_{0}(E)\mid  f(\gP/\gP\cap K)=1$ et $\gP$ est non ramifiŽ sur
$K\}$. L'ensemble $S_1$ est rŽgulier de densitŽ de Dirichlet 1, car pour obtenir $S_1$, on enlve ˆ tous les idŽaux premiers de $E$ un nombre fini d'idŽaux (ceux qui se ramifient sur $K$) et un sous-ensemble de ceux qui sont de degrŽ sur $\Q$ supŽrieur ˆ 1 qui est de densitŽ 0 (cf. Lemme \argan\ et utilisant la relation $f(\gP/\gP\cap \Q)=f(\gP/\gP\cap K)\cdot f(\gP\cap K/\gP\cap \Q)$). Cela prouve que $S$ est de densitŽ ${1\over [E:K]}$.\qed
\bigskip\goodbreak

\coro

{\sl Soit $L/K$ une extension galoisienne de corps de nombres de groupe de Galois $G$. 

\art{a)} Soit $H\subset G$ un
sous-groupe normal. Alors l'ensemble 
$$S=\{\P\in\gfP_{0}(K)\mid \P\hbox{ est non ramifiŽ dans $L$ et }{\rm Frob}(\gP/\P)\in H\hbox{ pour tout
}\gP\in\gfP_{0}(L)\hbox{ au-dessus de }\P\}$$
est rŽgulier de densitŽ ${1\over [G:H]}={|H|\over |G|}$.

\art{b)} De plus, si $H_1$ et $H_2$ sont des sous-groupes normaux de $G$, alors l'ensemble
$$S'=\{\P\in\gfP_{0}(K)\mid \P\hbox{ est non ramifiŽ dans $L$ et }{\rm Frob}(\gP/\P)\in H_1\cup H_1\hbox{ pour tout
}\gP\in\gfP_{0}(L),\ \gP|\P\}$$
est rŽgulier de densitŽ ${|H_1\cup H_2|\over |G|}$. Ce rŽsultat se gŽnŽralise ˆ un nombre fini de sous-groupes.

}

{\bf Preuve}

Pour a), posons $F$ le corps fixe par $H$. Puisque $H$ est un sous-groupe normal de $G$, la thŽorie de Galois nous dit que l'extension $F/K$ est galoisienne de groupe de Galois $G/H$. Si $\P_0\in \gfP_{0}(F)$ est un idŽal premier au-dessus de $\P$, par hypothse, on
a ${\rm Frob}(\P_0/\P)={\rm Frob}(\gP/\P)|_{F}={\rm Id_F}$ dans $G/H$, pour tout idŽal $\gP\in \gfP(L)$ au-dessus de
$\P_0$. Or, on sait que ${\rm Frob}(\P_0/\P)$ engendre $Z(\P_{0}/\P)$ qui est de cardinal $f(\P_{0}/\P)$. Ainsi,  $f(\P_0,\P)=1$, ce qui veut dire que $\P$ se dŽcompose totalement dans $F$. En
rŽsumŽ, $S$ est l'ensemble des idŽaux premiers de $O_K$ qui se dŽcomposent totalement dans $O_F$. Le Lemme \argao\ nous apprend alors que $S$ est rŽgulier de densitŽ ${1\over [F:K]}={1\over [G:H]}$.

Prouvons b). Les sous-groupes $H_1$, $H_2$ et $H_1\cap H_2$ correspondent ˆ des ensembles $S_1$,
$S_2$ et $S_1\cap S_2$  de densitŽ de Dirichlet respectivement de ${|H_1|\over |G|}$, ${|H_2|\over |G|}$ et
${|H_1\cap H_2|\over |G|}$. Il est clair que $S'=S_1\cup S_2$ et que 

$$\sum_{\P\in S'}{1\over \N(\P)^s}=\sum_{\P\in S_1}{1\over \N(\P)^s}+\sum_{\P\in S_2}{1\over \N(\P)^s}-\sum_{\P\in
S_1\cap S_2}{1\over \N(\P)^s}\sim -{|H_1|+|H_2|-|H_1\cap H_2|\over |G|}\cdot\Log(s-1)$$
Ainsi, $S'$ est rŽgulier de densitŽ de Dirichlet ${|H_1|+|H_2|-|H_1\cap H_2|\over |G|}={|H_1\cup H_2|\over |G|}$.
\qed

\medskip\goodbreak
Voici maintenant un rŽsultat annoncŽ depuis longtemps~:
\medskip
\th \ {\bf (SurjectivitŽ de l'application d'Artin)}

{\sl Soit $L/K$ une extension abŽlienne de groupe $G$ et $\m$ un $K$-module divisible par toutes les places de $K$ qui ramifient dans $L$. Alors l'application d'Artin
$$\Phi_{L/K}\ :I_K(\m)\longrightarrow G$$
est surjective.

}

\bigskip

{\bf Preuve}

Soit $\sigma\in G$. Notons $H=<\sigma>$ le sous-groupe engendrŽ par $\sigma$. Par le Corollaire \argap, partie a),
l'ensemble $S=\{\P\in\gfP_{0}(K)\mid \P\hbox{ est non ramifiŽ dans $L$ et }{\rm Frob}(\gP/\P)\in H\hbox{ pour tout
}\gP\in\gfP_{0 }(L)\hbox{ au-dessus de}\ \P\}$ est de densitŽ ${1\over [G:H]}$. Supposons par l'absurde que $\sigma$ ne
soit pas atteint par $\Phi_{L/K}$. Notons $\cup H_i$ la rŽunion de tous les sous-groupes propres de $H$. Il est
Žvident que $\cup H_i\subset H$, mais que $\cup H_i\ne H$, car $\sigma\not\in \cup H_i$. Donc, par hypothse absurde, $S$ est
inclu dans $S'=\{\P\in\gfP_0(K)\mid \P\hbox{ est non ramifiŽ dans $L$ et }{\rm Frob}(\gP/\P)\in\cup H_i\hbox{ pour
tout }\gP\in\gfP_0(L)\hbox{ au-dessus de}\ \P\}$. Ainsi, on trouve, gr‰ce au Corollaire \argap, partie b)~:

$$\delta (S)\leq \delta(S')={|\cup H_i|\over |G|}<{|H|\over |G|}=\delta(S),$$
ce qui est une contradiction.\qed
\bigskip
Voici un autre rŽsultat important qui montre en substance que la connaissance des idŽaux premiers qui se dŽcomposent
totalement dŽtermine une extension galoisienne. 
\bigskip
\th

{\sl Soient $L_1/K$ et $L_2/K$ deux extensions galoisiennes de corps de nombres (plongŽes dans $\C$).
Posons, pour $i=1,2$, $S_i=\{\P\in\gfP_0(K)\mid \P$ se dŽcompose compltement dans $L_i\}$. Alors les
conditions suivantes sont Žquivalentes~:

\art{a)} $L_1\subset L_2$

\art{b)} $S_2\subset S_1$

\art{c)} $S_2\setminus S_1$ est de densitŽ de Dirichlet nulle.

}

{\bf Preuve}

a)$\Longrightarrow$ b) $\Longrightarrow$ c) est trivial. Prouvons c)$\Longrightarrow$ a). Posons
$S=\{\P\in\gfP_0(K)\mid \P$ se dŽcompose compltement dans $L_1L_2\}$. Le Lemme \argao\ nous dit
que $\delta (S)={1\over [L_1L_2:K]}$ et $\delta (S_i)={1\over [L_i:K]}$, pour $i=1,2$. Par le Lemme ``dŽcomposition-ramification'', page \the\gaagc, $S=S_{1}\cap S_{2}$, et ainsi $S_2$ est la rŽunion disjointe de $S$ et de $S_2\setminus S_1$. Donc, par hypothse, $\delta(S_2)=\delta(S)$. Ce qui montre que $[L_1L_2:K]=[L_2:K]$, donc $L_1L_2=L_{2}$, ce qui veut dire que
$L_1\subset L_2$.\qed

\bigskip

\coro

{\sl Avec les mmes notations, On a Žquivalence entre 

\art{a)}$L_1=L_2$

\art{b)}$S_1=S_2$

\art{c)}$S_1\triangle S_2$ est de densitŽ de Dirichlet nulle (ici $\triangle$ dŽsigne la diffŽrence
symŽtrique).

}

{\bf Preuve}

C'est immŽdiat.\qed

\bigskip  

\th

{\sl Soit $K$ un corps de nombres, $\m$ un $K$-module et $I(\m)\supset H\supset P_\m$, $H$ Žtant un
sous-groupe. Soit $S\subset H\cap \gfP_0 (K)$ ayant une densitŽ de Dirichlet $\delta$. Posons $h=[I(\m):H]$.
Alors on a $$\delta\leq {1\over h}.$$

Supposons de plus que $\delta >0$. Soit $\chi$ un caractre de $I(\m)/H$. On note $\widehat{I(\m)/H}$
l'ensemble de ces caractres. Si $\chi\ne{\bf 1}$, alors on a

$$L_\m(1,\chi)\ne 0,$$
et pour chaque classe $\euf K$ de $I(\m)/H$, l'ensemble $\{\P\in \gfP_0(K)\mid \P\in {\euf K}\}$
est rŽgulier, de densitŽ ${1\over h}$.

Remarquons que $L_{\m}(s,\chi)$ n'a ŽtŽ dŽfinie que lorsque $\chi$ est un caractre de $I(\m)/P_\m$. 
On associera tout caractre $\chi$ de  $I(\m)/H$ ˆ $\widetilde{\chi}=\chi\circ\varphi$ o $\varphi$ est l'homomorphisme
canonique $I(\m)/P_\m\rightarrow I(\m)/H$.

}

{\bf Preuve}

Soit $\chi\in\widehat{I(\m)/H}$. En identifiant $\chi$ et $\widetilde{\chi}$, et gr‰ce
au ThŽorme \argaj, on sait que pour tout $s$ tel que $\Re(s)>1$, on a~:

$$\llog L_\m (s,\chi)=\sum_{\P\not\hskip 1.9pt |\hskip 0.9 pt\m}{\chi(\P)\over
\N(\P)^s}+g_\chi(s),\eqno{(i)}$$
avec $g_\chi(s)$ une fonction holomorphe au voisinage de $s=1$. Si $\chi\ne{\bf 1}$, on note $n(\chi)$ l'ordre
du zŽro de $L_\m(s,\chi)$. Gr‰ce au ThŽorme \argah, on sait que $n(\chi)\geq 0$. On a $n(\chi)=0$ si et seulement
si $L_\m(1,\chi)\ne 0$. Ainsi, $L_\m(s,\chi)=(s-1)^{n(\chi)}\cdot f_\chi(s)$ o $f_\chi(s)$ est une fonction
holomorphe sur un voisinage de 1 et $f_\chi(1)\ne 0$. En vertu du Lemme \argak, il existe une fonction 
logarithme (qu'on notera $\log f_\chi$) tel que pour tout $s$ tel que $\Re(s)>1$, on ait~: 

$$\llog L_\m(s,\chi)=n(\chi)\cdot \Log(s-1)+\log f_\chi(s),$$
avec $\log f_\chi(s)$ holomorphe sur un voisinage de 1. 

Si $\chi={\bf 1}$, le Corollaire \argal\ nous dit qu'il existe une fonction $f(s)$ holomorphe au voisinage de 1,
telle que pour tout $s$ tel que $\Re(s)>1$, on ait

$$\llog L_\m(s,{\bf 1})=-\Log(s-1)+f(s).\eqno{(ii)}$$
Ainsi, en sommant sur tous les $\chi$ et en remarquant que $\sum_{\chi\in
\widehat{I(\m)/H}}\chi(\P)=\cases{0& si
$\P\not\in H$\cr h&si $\P\in H$\cr}$  \goodbreak(voir DŽfinition \argag), on a pour tout $s$ tel $\Re(s)>1$~:

$$h\cdot\sum_{\P\in H\cap \gfP_{0} (K)}{1\over \N(\P)^s}=\sum_{\chi\in
\widehat{I(\m)/H}}\llog L_\m(s,\chi)+g_0(s)=-(1-\sum_{\chi\ne {\bf 1}}n_\chi)\Log(s-1)+\widetilde g(s),$$
o $g_0(s)$ et $\widetilde g(s)$ sont des fonctions holomorphes au voisinage de $s=1$  (vous aurez
remarquŽ que $g_0$ est la somme des $-g_\chi$ et  $\widetilde g$ est la somme de $g_0$, des $\log(f_\chi)$
et de $f$). Par hypothse, si $s>1$ est rŽel, il existe une fonction $h(s)\rightarrow 0$ lorsque
$s\rightarrow 1$ telle que~:
$$ \sum_{\P\in S}{1\over \N(\P)^s}=-(\delta+h(s))\cdot\Log(s-1).$$
Si $s>1$, on a 

$$0\leq \sum_{\P\in (H\cap \gfP_{0}(K))\setminus S}{1\over \N(\P)^s}=-(\underbrace{{1\over h}(1-\sum_{\chi\ne
{\bf 1}}n(\chi))-\delta}_{(*)}-h(s))\cdot\Log(s-1)+\widetilde g(s).$$
Il est clair que  $(*)$ doit tre positif, car sinon l'inŽgalitŽ $\leq$ devient fausse pour des
$s\rightarrow 1$. On en dŽduit que ${1\over h}-\delta\buildrel (**)\over \geq {1\over h}\sum_{\chi\ne {\bf
1}} n(\chi)\geq 0$, ce qui prouve que $\delta\leq {1\over h}$.

Supposons ˆ prŽsent que $\delta >0$. On trouve, en rŽutilisant $(**)$ que ${1\over h}>{1\over
h}-\delta\buildrel (**)\over\geq  {1\over h}\cdot \sum_{\chi\ne {\bf 1}} n(\chi)$; et ainsi $\sum_{\chi\ne
{\bf 1}} n(\chi)<1$, ce qui montre que $n(\chi)=0$, pour tout $\chi\ne{\bf 1}$. Par dŽfinition des
$n(\chi)$, cela implique que pour tout  $\chi\ne{\bf 1}$,
$$L_\m(1,\chi)\ne 0\quad \hbox{ et donc que } \llog L_\m(s,\chi)\hbox{ est holomorphe au voisinage de
$s=1$}.$$

Soit ${\euf K}$ une classe de $I(\m)$ modulo $H$, et soit ${\euf a}\in{\euf K}$. Il est Žvident que si $\P\in \gfP_{0}(K)\cap I(\m)$, on a $\P\in {\euf K}\iff\aa^{-1}\P\in H$. Ainsi, $\sum_{\chi\in
\widehat{I(\m)/H}}\chi^{-1}(\aa)\cdot\chi(\P)=\cases{h& si
$\P\in {\euf K}$\cr 0&sinon\cr}$. Et donc,  en rŽutilisant $(i)$
et $(ii)$, on trouve~:

$$\eqalign{h\cdot\sum_{\P\in {\euf K}\cap \gfP_0(K)}{1\over
\N(\P)^s}&=\sum_{\chi\in\widehat{I(\m)/H}}\chi^{-1}({\euf a})\cdot\sum_{\P\in I(\m)\cap
\gfP_0(K)}{\chi(\P)\over \N(\P)^s}\cr
&=\sum_{\chi\in\widehat{I(\m)/H}}\chi^{-1}({\euf
a})\underbrace{\llog L_\m(s,\chi)}_{\hbox{holomorphe au vg de 1 si $\chi\ne {\bf 1}$}}+\hat g(s)\cr
&=-\Log(s-1)+\hat{\hat g}(s),\cr}$$
o $\hat g(s)$ et $\hat{\hat g}(s)$ sont des fonctions holomorphes au voisinage de $1$. Cela prouve que
${\euf K}\cap \gfP_0(K)$ est rŽgulier de densitŽ ${1\over h}$ et donc le thŽorme.\qed

\bigskip\goodbreak

\bigskip
Maintenant, nous pouvons enfin Žnoncer le premier rŽsultat important de la thŽorie du corps de classe~:
\bigskip

\th \ {\soustitre (premire inŽgalitŽ du corps de classe)}

{\sl  Soit $L/K$ une extension galoisienne, $\m$ un $K$-module, alors

$$h:=[I_K(\m):P_\m\cdot N_{L/K}(I_L(\widetilde{\m}))]\leq [L:K].$$

De plus, si on pose $H=P_\m\cdot N_{L/K}(I_L(\widetilde{\m}))$ et soit $\chi\in\widehat{I_K(\m)/H}$, $\chi\ne {\bf 1}$,
alors $L_\m(1,\chi)\ne 0$ et si ${\euf K}\in I_K(\m)/H$ alors ${\euf K}\cap\gfP_{0}(K)$ est rŽgulier de densitŽ
${1\over h}$.

 }

{\bf preuve }

Posons $S=\{\P\in\gfP_0(K)\mid \P$ se dŽcompose compltement dans $L$ et $\P\notdiv \m\}$. Le Lemme \argao\ montre que
$S$ est de densitŽ $\delta:={1\over [L:K]}$ (les premiers divisant $\m$, ne changent rien ˆ l'affaire,
puisqu'il n'y en a qu'un nombre fini). On applique alors le ThŽorme \argat\ dans ce cas (on peut, car souvenons-nous que tout idŽal premier totalement dŽcomposŽ est une norme), cela  montre notre
thŽorme.\qed

\vfill\eject


\global\advance\chapnomb by 1
\nomb=1

\centerline{\para Chapitre 3 : }
\centerline{\para ThŽorme de $\taille{18}\check{\bf C}$ebotarev }
\bigskip
Nous allons (comme le titre peut le faire imaginer) dŽmontrer le thŽorme de $\check{\rm C}$ebotarev.
{\it A priori}, ce thŽorme ne nous sera pas utile pour prouver les principaux rŽsultats de la thŽorie
du corps de classe. Mais plusieurs rŽsultats obtenus aux chapitres prŽcŽdents nous permettrons de prouver
un version affaiblie de ce thŽorme. Nous aurons nŽanmoins recourt ˆ un rŽsultat que nous montrerons
ultŽrieurement pour passer de la version faible ˆ la version forte de ce thŽorme. Mais foin de
considŽration un peu oiseuse, voici de quoi il s'agit~:
\headline={\hfill \phantom{ouuh}\hfill}
\bigskip

\th \ {\bf (ThŽorme de $\check{\bf C}$ebotarev)}

{\sl Soit $L/K$ une extension galoisienne de corps de nombres de groupe de Galois $G$. Soit $C$ une
classe de conjugaison de $G$. Alors l'ensemble 

$$A=\{\P\mid \P\hbox{ est un idŽal premier de $K$ non ramifiŽ dans $L$
 avec }{\rm Fr}_{L/K}(\P)=C\}$$
est rŽgulier et sa densitŽ de Dirichlet vaut ${|C|\over |G|}$.

}
\bigskip
On rappelle que si $\P$ est un idŽal de $O_K$, ${\rm Fr}_{L/K}(\P)$ est l'ensemble $\{ {\rm
Frob}(\gP/\P)\mid\gP|\P\}$, qui est une classe de conjugaison, car tous les ${\rm Frob}(\gP/\P)$ sont
conjuguŽs entre eux si $\P$ est fixŽ; et  enfin ${\rm Frob}(\gP/\P)$ est l'unique l'ŽlŽment de
$G$ qui satisfait~:

$${\rm Frob}(\gP/\P)(x)\equiv x^{\N(\P)}\pmod \gP\quad \hbox{ pour tout }x\in O_L.$$
Tout cela a ŽtŽ vu au Chapitre 0 en plus grand dŽtail.

\medskip

Nous allons montrer ce rŽsultat en trois Žtapes~: la premire lorsque $L/K$ est une extension cyclotomique
(c'est ˆ dire que $L\subset K(\zeta)$ pour une certaine racine de l'unitŽ $\zeta$, la deuxime lorsque
$L/K$ est abŽlienne et la troisime dans le cas quelconque. Nous n'aurons besoin de la thŽorie du corps de
classe que dans le troisime cas. Le premier cas n'est qu'une compilation de rŽsultats dŽjˆ obtenus
\bigskip

\prop

{\sl Le thŽorme de $\check{\sl C}$ebotarev est vrai dans le cas o $L/K$ est une extension cyclotomique.}

{\bf Preuve}

Soit $C$ une classe de conjugaison de $G$. Puisque $G$ est abŽlien, $C=\{\sigma\}$, $\sigma\in G$. Soit $\m$ un $K$-module satisfaisant les hypothses du ThŽorme \argp. Nous savons que, dans notre cas, l'application d'Artin $\Phi^\m_{L/K}\ :I_K(\m)\longrightarrow G$ est surjective (cf. Proposition \argaq), ainsi il existe ${\euf a}\in I_K(\m)$ tel que $\Phi^\m_{L/K}({\euf a})=\sigma$. Posons $H=\ker \Phi^\m_{L/K}$. Nous avons vu au ThŽorme \argp\ que dans notre cas, nous avons $P_\m\subset H\subset I_K(\m)$. Posons  $S=\{\P\in\gfP_0(K)\mid \P$ est non ramifiŽ dans $L$ et ${\rm Fr}_{L/K}(\P)=\{Id_G\}\}$. Le Corollaire \argap, partie a) nous montre que $S$ est rŽgulier de densitŽ de Dirichlet ${1\over |G|}$, car dans notre cas, ${\rm Fr}_{L/K}(\P)$ est rŽduit ˆ un seul ŽlŽment qui est ${\rm Frob}(\gP/\P)$ o $\gP$ est
n'importe quel idŽal de $O_L$ au-dessus de $\P$. 
Il est enfin clair que $S\subset H\cap \gfP_0(K)$. Soit ${\euf K}$ la classe de $\euf a$ dans $I_K(\m)/H$. Le
ThŽorme \argat\ s'applique alors et on a donc que $S'=\{\P\in \gfP_0(K)\mid \P\in {\euf
K}\}=\{\P\in \gfP_0(K)\mid \Phi^\m_{L/K}(\P)={\sigma}\}=\{ \P\in \gfP_0(K)\mid {\rm Fr}_{L/K}(\P)=C\}$ est
rŽgulier de densitŽ ${1\over [I_K(\m):H]}={1\over |G|}$. Cela prouve notre proposition.\qed
\bigskip\goodbreak

On peut dŽjˆ montrer le thŽorme de Dirichlet sur les progressions arithmŽtiques~:
\bigskip

\coro

{\sl Si $m$ et $n$ sont des nombres entiers premiers entre eux, $n>1$, alors il existe une infinitŽ de nombres
premiers $p$ tels que $p\equiv m\pmod n$. }

{\bf Preuve}

Posons $K=\Q$ et $L=\Q(\zeta_n)$. Posons encore $\sigma_m$ l'ŽlŽment de ${\rm Gal}(L/K)$ qui envoie
$\zeta_n$ sur $\zeta_n^m$. Par le thŽorme de $\check{\sl C}$ebotarev appliquŽ ˆ $L/K$ (qui est une
extension cyclique), l'ensemble des premiers $p$ de $\gfP_{0}(\Q)$ tels que ${\rm Frob}(\P/p)=\sigma_m$ ($\P$ est un idŽal premier de $O_L$ au-dessus de $p$) est de densitŽ ${1\over [L:K]}\ne 0$, donc est infini. Or, de tels $p$ sont congrus ˆ $m$ modulo $n$, car  ${\rm Frob}(\P/p)=\sigma_m$ implique en particulier que
$\zeta_n^p\equiv \zeta_n^m\pmod \P$ ce qui veut dire (Lemme \argo) que $\zeta_n^p=\zeta_n^m$ ou encore que
$p\equiv m\pmod n$.\qed
\bigskip
{\soustitre PrŽparations au cas abŽlien}
\bigskip
\headline={\hfill \smcap ThŽorme de $\check {\rm C}$ebotarev\hfill}
\lem

{\sl Soit $L/K$ une extension de corps de nombres et $m>1$ un entier naturel. Alors il existe une
extension $M/K$ cyclotomique et cyclique de degrŽ $m$ telle que $M\cap L=K$ (les extension $L/K$, $M/K$
telle que $M\cap L=K$ sont dites {\it linŽairement disjointes}).

}

{\bf Preuve}

Voyons tout d'abord qu'il suffit de prouver le lemme pour $K=\Q$~: supposons donc que ce soit vrai dans ce cas-lˆ et montrons le cas gŽnŽral. Supposons donc $L/K$ comme dans l'hypothse de ce lemme et qu'il existe $M/\Q$ une extension cyclique cyclotomique  de degrŽ
$m$ telle que $M\cap L=\Q$. On a donc le diagramme suivant~:
\vglue 5cm

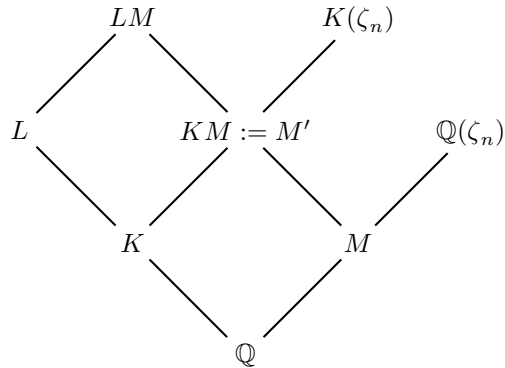
\rput(5,0){\rnode{Q}{$ \Q  $}} 
\rput(4,1){\rnode{K}{$K$}}
\rput(3,2){\rnode{L}{$L$}}
\rput(4,3){\rnode{LM}{$LM$}}
\rput(5,2){\rnode{KM}{$KM:=M'$}}
\rput(6,3){\rnode{Kzet}{$K(\zeta_{n})$}}
\rput(7,2){\rnode{Qzet}{$\Q(\zeta_{n})$}}
\rput(6,1){\rnode{M}{$M$}}
\ncline[nodesep=3pt]{Q}{K} 
\ncline[nodesep=3pt]{K}{L} 
\ncline[nodesep=3pt]{L}{LM} 
\ncline[nodesep=3pt]{LM}{KM} 
\ncline[nodesep=3pt]{K}{KM} 
\ncline[nodesep=3pt]{KM}{Kzet} 
\ncline[nodesep=3pt]{M}{Qzet} 
\ncline[nodesep=3pt]{Q}{M}
\ncline[nodesep=3pt]{M}{KM}


L'extension $(KM=M')/K$ est Žvidemment cyclotomique, car $M'\subset K(\zeta)$. De plus, puisque $L\cap M=\Q$, on a bien sžr $M\cap K=\Q$. La thŽorie de Galois dit que l'extension $M'/K$ est galoisienne de groupe isomorphe ˆ ${\rm Gal}(M/(M\cap K))={\rm Gal}(M/\Q)$ qui est, par hypothse, un groupe cyclique de degrŽ $m$. En rŽsumŽ, $M'/K$ est cyclotomique cyclique de degrŽ $m$. Reste ˆ voir que $M'\cap L=K$. De la mme manire qu'avant et puisque $L\cap M=\Q$, on a que ${\rm Gal}(LM/L)$ est isomorphe ˆ ${\rm Gal}(M/(M\cap L))={\rm Gal}(M/\Q)$ qui est cyclique d'ordre $m$. En appliquant ˆ nouveau la thŽorie de Galois au diagramme
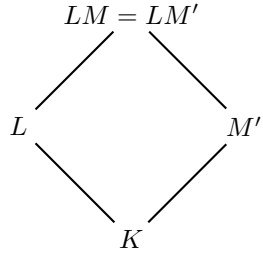
\vglue 3cm
\rput(5,0){\rnode{K}{$K $}} 
\rput(4,1){\rnode{L}{$L$}}
\rput(5,2){\rnode{LM}{$LM=LM'$}}
\rput(6,1){\rnode{Mpr}{$M'$}}
\ncline[nodesep=3pt]{K}{Mpr} 
\ncline[nodesep=3pt]{K}{L} 
\ncline[nodesep=3pt]{L}{LM} 
\ncline[nodesep=3pt]{Mpr}{LM} 

on trouve que ${\rm Gal}(LM'/L)\simeq{\rm Gal}(M'/(M'\cap L))$. Donc $\left |{\rm Gal}(M'/(M'\cap L))\right
|=\left | {\rm Gal}(M'/K)\right |=m$. Cela prouve que $M'\cap L=K$.
\medskip
Prouvons alors le cas $K=\Q$. Soit donc $L/\Q$ une extension de degrŽ finie. Si $L_1/\Q$ et $L_2/\Q$ 
sont des sous-extension cyclotomique de $L/\Q$, alors $L_1L_2/\Q$ en est aussi une. Il y a donc, parmi
les sous-extensions de $L/\Q$, une extension cyclotomique maximale $L_0/\Q$. Soit $n$ tel que $L_0\subset
\Q(\zeta_n)$. Soit $p\in \gfP_{0}(\Q)$ tel que $p\notdiv n$ et $p\equiv 1\pmod m$ (l'existence d'un tel $p$ est un cas particulier du ThŽorme de Dirichlet sur les progressions arithmŽtiques, ou alors de la remarque ci-aprs). Par maximalitŽ de $L_0$, on a $\Q(\zeta_p)\cap L\subset L_0\subset \Q(\zeta_n)$. Donc $\Q(\zeta_p)\cap L\subset\Q(\zeta_p)\cap\Q(\zeta_n)=\Q$. Comme $\Q(\zeta_p)/\Q$ est cyclique de degrŽ $p-1$, il y a par la thŽorie de Galois et puisque tout sous-groupe d'un groupe cyclique est normal, il existe une sous-extension $M/\Q$ cyclique d'ordre $m$, c'est le sous-corps fixe par l'unique sous-groupe d'ordre ${p-1\over m}$ de ${\rm Gal}(\Q(\zeta_p)/\Q)$.\qed

\bigskip

\rem

Soit $m$ et $n$ des nombres entiers, le fait de l'existence d'un nombre premier tel que $p\notdiv n$ et tel que $p\equiv 1\pmod m$ est effectivement un cas particulier du ThŽorme de Dirichlet sur les progressions arithmŽtiques, mais on peut le dŽmontrer ҈ la main" de manire semblable ˆ la preuve historique d'Euclide sur l'infinitŽ des nombres premiers. Dire que $p\equiv 1\pmod m$ est Žquivalent au fait que
$\F_p^*$ contient un ŽlŽment d'ordre $m$ (car $\F_p^*$ est un groupe cyclique ( [Jac1, Theorem 2.18, p.132])); cela veut dire que $\F_p^*$ contient une racine du polyn™me cyclotomique $\Phi_m$, c'est-ˆ-dire que l'Žquation $\Phi_m(X)\equiv 0\pmod p$ a une solution. Montrons qu'il y a une infinitŽ de tels $p$ (parmi ceux-lˆ, il y en a qui ne divisent pas $n$ et le tour est jouŽ). Si $p_1,\ldots, p_r$ sont des nombres premiers tels que $\Phi_m(X)\equiv 0$ a une solution modulo $p_1,\ldots, p_r$ (c'est possible, car l'Žquation
$\Phi_m(X)=0,\pm 1$ n'a qu'un nombre fini de solutions); alors il existe $k\in\Z$ tel que $|\Phi_m(m\cdot
p_1\cdot\cdots\cdot p_r\cdot k)|>1$. Si $p|\Phi_m(m\cdot p_1\cdot\cdots\cdot p_r\cdot k)$, alors $p\ne
p_1,\ldots ,p_r$ et $p/\!\! | m$, car $\Phi_m(m\cdot p_1\cdot\cdots\cdot p_r\cdot k)\equiv 1$ modulo
$p_1,\ldots , p_r$ et $m$ (le coefficient constant de tout polyn™me cyclotomique est 1). Cela prouve
notre rŽsultat. 
\bigskip 
\newcount\gaagr\gaagr=\gaga

\lem

{\sl Soit $m$ et $n$ des entiers tels que $n|m$, on pose $T(m,n)$ le nombre des ŽlŽments du groupe cyclique
d'ordre $m$ qui ont pour ordre un multiple de $n$. Si $n=p_1^{a_1}\cdot\cdots\cdot p_r^{a_r}$ et
$m=p_1^{b_1}\cdot\cdots \cdot p_r^{b_r}$ avec $a_i\leq b_i$ pour tout $i=1,\ldots , r$. Alors 

$$T(m,n)=m\cdot\prod_{i=1}^r(1-p_i^{a_i-1-b_i}).$$

}
\goodbreak
{\bf Preuve}

Gr‰ce au thŽorme chinois, on voit que $T(m,n)=\prod_{i=1}^rT(p_i^{b_i},p_i^{a_i})$.
De plus, $T(p^b,p^a)=$ nbre d'Žl. d'ordre $p^a$ + nbre d'Žl. d'ordre $p^{a+1}+\cdots+$ nbre d'Žl. d'ordre
$p^{b}=\varphi(p^a)+\varphi(p^{a+1})+\cdots+\varphi(p^b)=(p^{a}-p^{a-1})+(p^{a+1}-p^a)+\cdots
+(p^b-p^{b-1})=p^b-p^{a-1}=p^b\cdot (1-p^{a-1-b})$, ce qui montre alors notre lemme.\qed

\bigskip\goodbreak

\prop 

{\sl Le thŽorme de $\check {\rm C}$ebotarev affaibli est vrai pour les extensions abŽlienne. Par
thŽorme de $\check {\rm C}$ebotarev affaibli, on entend qu'on a la densitŽ de Dirichlet, mais pas la
rŽgularitŽ.

}
{\bf Preuve}

Soit $L/K$ une extension abŽlienne de corps de nombres, $G={\rm Gal}(L/K)$ et $\sigma\in G$ d'ordre $n$.
Soit $m$ un multiple de $n$. Soit $M/K$ une extension cyclotomique, cyclique de degrŽ $m$ telle que $M\cap L=K$; le Lemme \argay\ nous en assure l'existence. Soit $\tau\in{\rm Gal}(M/K)$ tel que $n$ divise l'ordre de $\tau$. Posons $F=LM$. La thŽorie de Galois nous dit que ${\rm Gal}(F/K)\simeq{\rm Gal}(L/K)\times{\rm Gal}(M/K)$. Soit $\rho\in{\rm Gal}(F/K)$ l'unique ŽlŽment tel que $\rho|_L=\sigma$ et $\rho|_M=\tau$. Posons $E=F^\rho$ (le sous-corps de $F$ fixe par $\rho$). Il est clair que $M\cap E=M^\tau$ (le
sous-corps de $M$ fixe par $\tau$). A priori, on a $F\supset EM\supset E$. On va voir, qu'en fait
$F=EM$. On a les ŽgalitŽs~: $[F:E]=$ ordre de $\rho$ = ordre de $\tau=[M:M^\tau]=[M:M\cap E]=[EM:E]$ (la
dernire ŽgalitŽ venant de la thŽorie de Galois). Donc, $F=EM$. On en dŽduit que $F/E$ est une
extension cyclotomique; en effet, puisque $M/K$ est cyclotomique, on a $M\subset K(\zeta)$ pour une
racine de l'unitŽ $\zeta$.  Donc, $F=EM\subset E\, K(\zeta)=E(\zeta)$, elle est aussi cyclique (engendrŽe par $\rho$), par la thŽorie de Galois. Le thŽorme de $\check {\rm C}$ebotarev pour les extensions cyclotomiques s'applique ˆ $F/E$, et donc, l'ensemble
$$A'_\tau=\{{\euf q}\in\gfP_{0}(E)\mid {\euf q}\hbox{ ne ramifie pas dans $F$ et } \Phi_{F/E}({\euf
q})=\rho\}$$
a une densitŽ de Dirichlet $\delta(A'_\tau)={1\over [F:E]}$. Si on enlve ˆ $A'_\tau$ les idŽaux qui sont
ramifiŽs sur $K$, on trouve la mme densitŽ. De plus, si on ne prend que les idŽaux $\euf q$ de $A'_\tau$
tels que $f({\euf q}/{\euf q}\cap K)=1\ (\Leftrightarrow \N({\euf q}\cap K)=\N({\euf q}))$, cela ne change
rien non plus ˆ la densitŽ (cf. Lemme \argan). Ainsi
$$A''_\tau=\{{\euf q}\in\gfP_{0}(E)\mid {\euf q}\hbox{ ne ramifie pas dans $F$, n'est pas ramifiŽ sur $K$,
}\N({\euf q}\cap K)=\N({\euf q})\hbox{ et } \Phi_{F/E}({\euf
q})=\rho\}$$
est tel que $\delta(A'_\tau)=\delta(A''_\tau)$. Soit encore
$$A_\tau=\{\P\in\gfP_{0}(K)\mid {\euf p}\hbox{ ne ramifie pas dans $F$ et } \Phi_{F/K}(\P)=\rho\}.$$

Soit ${\euf q}\in A''_\tau$ et $\P={\euf q}\cap K$. Alors $\P\in A_\tau$. En effet, par propriŽtŽ de $A''_\tau$, $\P$ ne ramifie pas dans $F$. Soit $x\in O_{F}$. On a $\Phi_{F/K}(\P)(x)\equiv x^{\N(\P)}\pmod{\gP}$, o $\gP$ est n'importe quel idŽal premier de $F$ au-dessus de $\P$, en particulier pour ceux au-dessus de $\qq$. Or, par hypothse, $\N(\P)=\N(\qq)$. Donc, pour tout $\gP$ au-dessus de $\qq$ et tout $x\in O_{F}$, on a
$$\Phi_{F/K}(\P)(x)\equiv x^{\N(\P)}\equiv x^{\N(\qq)}\equiv\Phi_{F/E}(\qq)(x)\pmod\gP.$$
Ainsi, $\Phi_{F/K}(\P)=\Phi_{F/E}(\qq)=\rho$ par hypothse sur $\qq$ et puisque l'automorphisme de Frobenius est caractŽrisŽ par ces congruences.

RŽciproquement, soit
$\P\in A_\tau$ et ${\euf q}\in\gfP_{0}(E)$ au-dessus de $\P$. En particulier, on a que ${\euf q}$
ne ramifie pas sur $K$ et ne ramifie pas dans $F$. Comme
$\Phi_{E/K}(\P)=\Phi_{F/K}(\P)|_{E}=\rho|_E={\rm Id}_E$, on en dŽduit que $f({\euf q}/\P)=1$, car
$\Phi_{E/K}(\P)={\rm Frob}({\euf q}/\P)$ qui est d'ordre $f({\euf q}/\P)$. Ainsi, $\N({\euf
q})=\N(\P)=\N({\euf q}\cap K)$. En particulier, cette ŽgalitŽ nous permet de montrer comme avant que
$\Phi_{F/E}({\euf q})=\Phi_{F/K}(\P)=\rho$, donc ${\euf q}\in A''_\tau$. On vient de voir aussi que
chaque $\P\in A_\tau$ est totalement dŽcomposŽ dans $E$, c'est ˆ dire que $\P$ a exactement $[E:K]$
idŽaux premiers $\euf q$ au-dessus de lui. Et donc, si $s>1$, on a $[E:K]\sum_{\P\in
A_\tau}\N(\P)^{-s}=\sum_{{\euf q}\in A''_\tau}\N({\euf q})^{-s}$. Cela prouve que la densitŽ de Dirichlet
de $A_\tau$ vaut 
$$\delta(A_\tau)={1\over [E:K]}\cdot \delta(A''_\tau)={1\over [E:K]}\cdot {1\over [F:E]}={1\over [F:K]}.$$
Il Žvident que pour tout $\tau \in {\rm Gal}(M/K)$, on a $$A_\tau\subset C_\sigma=\{\P\in \gfP(K)\mid
\P\hbox{ non ramifiŽ dans $L$ et }\Phi_{L/K}(\P)=\sigma\}.$$

D'autre part, les $A_\tau$ sont tous disjoints, lorsque $\tau$ varie. Posons $\euf T$ l'ensemble des $\tau\in {\rm Gal}(M/K)$
tels que l'ordre de~$\sigma$ $(=n)$ divise l'ordre de $\tau$. On a $\bigsqcup_{\tau\in{\euf T}} A_\tau\subset
C_\sigma$. Ainsi, pour $s>1$, on a~:

$$1\geq {\sum_{\P\in C_\sigma}\N(\P)^{-s}\over \sum_{\P\in \gfP_{0}(K)}\N(\P)^{-s}}\geq {\sum_{\tau\in{\euf
T}}\sum_{\P\in A_{\tau}}\N(\P)^{-s}\over \sum_{\P\in \gfP_{0}(K)}\N(\P)^{-s}},$$
et donc, en vertu de cette inŽgalitŽ et du calcul de $\delta(A_\tau)$, on trouve~:

$$\liminf_{s\rightarrow 1^+}{\sum_{\P\in C_\sigma}\N(\P)^{-s}\over \sum_{\P\in \gfP_0(K)}\N(\P)^{-s}}\geq
{T(m,n)\over [\underbrace{F}_{=LM}:K]}={1\over [L:K]}\cdot {T(m,n)\over [M:K]}={1\over [L:K]}\cdot
{T(m,n)\over m}.$$

Si $n=p_1^{a_1}\cdot\cdots\cdot p_r^{a_r}$, on choisit $m=p_1^{b_1}\cdot\cdots\cdot p_r^{b_r}$ avec des
$b_i$ aussi grands qu'on veut. Ainsi, gr‰ce au lemme prŽcŽdent, on a ${T(m,n)\over
m}=\prod_{i=1}^r(1-p_i^{a_i-1-b_i})\rightarrow 1$. Par consŽquent, pour tout $\varepsilon >1$, il existe
$s_1>1$ tel que si $1<s<s_1$, on ait ${\sum_{\P\in C_\sigma}\N(\P)^{-s}\over \sum_{\P\in
\gfP_{0}(K)}\N(\P)^{-s}}>{1\over [L:K]}-\varepsilon$. Ce qui veut dire, en posant $f_\sigma(s)={\sum_{\P\in
C_\sigma}\N(\P)^{-s}\over \sum_{\P\in
\gfP_{0}(K)}\N(\P)^{-s}}$ que l'on a

$$\liminf_{s\rightarrow 1} f_\sigma(s)\geq {1\over [L:K]}.$$

On va montrer que $\limsup_{s\rightarrow 1}f_\sigma(s)\leq {1\over [L:K]}$. Tout d'abord, remarquons que
pour tout $s>1$, on a $\sum_{\sigma\in G}f_\sigma(s)=1$ (c'est Žvident, car $G$ est abŽlien et l'application d'Artin est surjective). C'est ce
raisonnement qui nous permettra de conclure~: ce qu'on vient de voir nous dit que pour tout
$\varepsilon>0$, il existe $s_\sigma(\varepsilon)$ tel que $f_\sigma(s)>{1\over |G|}-{\varepsilon\over
|G|-1}$, si $1<s<s_\sigma(\varepsilon)$. Fixons
$\sigma_0$ et soit $s_0(\varepsilon)=\min_{\sigma\ne\sigma_0}s_\sigma(\varepsilon)$. Pour tout
$1<s<s_0(\varepsilon)$, on a ${|G|-1\over
|G|}-(|G|-1)\cdot{\varepsilon\over |G|-1}+f_{\sigma_0}(s)<\sum_{\sigma\in G}f_\sigma=1$, ce qui montre que 
$$f_{\sigma_0}(s)<{1\over |G|}+\varepsilon.$$
On a donc prouvŽ que pour tout $\sigma\in G$, on a
$$\limsup_{s\rightarrow 1} f_\sigma(s)\leq {1\over |G|}={1\over [L:K]},$$
ce qui prouve que $\lim_{s\to 1}{\sum_{\P\in C_\sigma}\N(\P)^{-s}\over \sum_{\P\in \gfP_{0}(K)}\N(\P)^{-s}}={1\over |G|}$, et donc le thŽorme de $\check {\rm C}$ebotarev affaibli, dans le cas abŽlien, car l'ensemble des idŽaux premiers de $K$ est de densitŽ 1 (voir DŽfinitions \argam). \qed
\bigskip
Attention ! ici, nous n'avons prouvŽe que la partie Òfaible" du thŽorme , c'est-ˆ-dire que nous avons la
densitŽ de Dirichlet, mais pas la rŽgularitŽ.

\bigskip

\th

{\sl Le $\check {\rm C}$ebotarev affaibli est vrai pour les extension galoisienne quelconque est vrai.
Mieux~: si le ThŽorme de $\check {\rm C}$ebotarev (fort) est vrai pour les extensions abŽliennes, il est
vrai pour les extensions galoisiennes quelconques.

}

{\bf Preuve}

Soit donc $L/K$ une extension galoisienne de corps de nombres de groupe $G$ et soit $C$ une classe de
conjugaison de $G$. Posons

$$A=\{\P\in \gfP_{0}(K)\mid \P\hbox{ est non ramifiŽ dans $L$
 avec }{\rm Fr}_{L/K}(\P)=C\}.$$
Pour la preuve, fixons, $\tau\in C$ et $K'=L^\tau$ (le corps fixe par $\tau$). La thŽorie de Galois nous
dit que $L/K'$ est galoisienne et ${\rm Gal}(L/K')=<\tau>$ cyclique. Posons encore

$$D'=\{{\euf q}\in \gfP_0(K')\mid {\euf q}\hbox{ ne ramifie pas dans $L$, non ramifiŽ sur $K$, $f({\euf
q}/{\euf q}\cap K)=1$ et $\Phi_{L/K'}({\euf q})=\tau$}\}.$$
Soit ${\euf q}\in D'$ et $\P={\euf q}\cap K$. Soit encore $\gP\in\gfP_0(L)$ au-dessus de $\euf q$. Puisque
$f({\euf q}/{\euf q}\cap K)=1$, on a $\N({\euf q})=\N(\P)$, et donc, ${\rm Frob}(\gP/\P)\equiv
x^{\N(\P)}=x^{\N({\euf q})}\equiv {\rm Frob}(\gP/{\euf q})\pmod\gP$. Ainsi, puisque les automorphismes de
Frobenius sont caractŽrisŽs par ces congruences, on a ${\rm Frob}(\gP/\P)={\rm Frob}(\gP/{\euf
q})=\Phi_{L/K'}({\euf q})=\tau$. On a aussi $e(\gP/\P)=e(\gP/{\euf q})\cdot e({\euf q}/\P)=1$. Donc,
$\P\in A$. Remarquons encore au passage que l'ordre de
${\rm Frob}(\gP/{\euf q})=f(\gP/{\euf q})$. D'autre part, l'ordre de $\tau= {\rm Frob}(\gP/{\euf q})$ vaut
$[L:K']$. Donc $\gP$ est le seul idŽal de $L$ au-dessus de $\euf q$ $(*)$.

RŽciproquement, soit $\P\in A$. Alors il existe $\gP\in\gfP_0(L)$ au-dessus de $\P$ tel que ${\rm
Frob}(\gP/\P)=\tau$. Alors $\qq:=\gP\cap K'$ n'est pas ramifiŽ sur $K$, ni dans $L$; et ${\rm
Frob}({\euf q}/\P)={\rm Frob}(\gP/\P)|_{K'}=\tau|_{K'}={\rm Id}_{K'}$. Cela montre que $f({\euf q}/\P)=1$
(c'est l'ordre du Frobenius). Donc $\N(\P)=\N({\euf q})$ et donc, par le mme argument que tout ˆ
l'heure, on a ${\rm Frob}(\gP/\P)={\rm Frob}(\gP/{\euf q})=\Phi_{L/K'}({\euf q})=\tau$. Ainsi, ${\euf
q}\in D'$.

Ainsi, on a prouvŽ que ${\euf q}\mapsto {\euf q}\cap K$ est une application surjective de $D'\to A$.
L'affirmation $(*)$ prouve de plus que pour chaque $\P\in A$, le nombre de ${\euf q}$ dans $D'$ au-dessus
de $\P$ est Žgal au nombre de $\gP\in\gfP_0(L)$ au-dessus de $\P$ tel que ${\rm Frob}(\gP/\P)=\tau$. Soit
$\gP_0$ l'un de ces $\gP$. Soit $\sigma\in G$. Soit $C_G(\tau)$ le centralisateur de $\tau$ dans $G$
(c'est-ˆ-dire les ŽlŽments de $G$ qui commutent avec $\tau$). De la relation
${\rm Frob}(\sigma(\gP_0))/\P)=\sigma\tau\sigma^{-1}$, on aura ${\rm
Frob}(\sigma(\gP_0))/\P)=\tau$ si et seulement si $\sigma\in C_G(\tau)$. Il est donc clair que
$Z(\gP_0/\P)=\{\nu\in G\mid \nu(\gP_0)=\gP_0\}$ est un sous-groupe de $C_G(\tau)$
 et chaque $\sigma(\gP_0)$ est comptŽ $|Z(\sigma(\gP_0)/\P)|=|Z(\gP_0/\P)|$ fois. En dŽfinitive, le
nombre de ${\euf q}\in D'$ au-dessus de $\P\in A$ est ${|C_G(\tau)|\over |Z(\gP_0/\P)|}$.

D'autre part, $|C|=[G:C_G(\tau)]$ et $|Z(\gP_0/\P)|=f(\gP_0/\P)=f(\gP_0/{\euf q})=[L:K']$. Ainsi,
$${|C_G(\tau)|\over |Z(\gP_0/\P)|}={|C_G(\tau)|\cdot |C|\over [L:K']\cdot |C|}={|G|\over [L:K']\cdot
|C|}={[L:K]\over [L:K']\cdot |C|}.$$
Ce qui montre que si $\Re(s)>1$, on a $\sum_{\P\in A}\N(\P)^{-s}={|C|\cdot [L:K']\over
[L:K]}\cdot\sum_{{\euf q}\in D'}\N({\euf q})^{-s}$. D'autre part, si 
$$D=\{{\euf q}\in \gfP(K')\mid {\euf q}\hbox{ ne ramifie pas dans $L$ et }\Phi_{L/K'}({\euf q})=\tau\},$$
on a $\sum_{{\euf q}\in D}\N({\euf q})^{-s}=\sum_{{\euf q}\in D'}\N({\euf q})^{-s}+g(s)$ o $g$ est une
fonction holomorphe sur un voisinage de 1; en effet, pour obtenir $D$, on a ajoutŽ ˆ $D'$ un nombre fini
d'idŽaux (ceux qui ramifient sur $K$) et ceux dont le $f>1$, qui est un ensemble de densitŽ nulle (cf.
Lemme \argan). Par le thŽorme de $\check{\rm C}$ebotarev appliquŽ ˆ l'extension $L/K'$ et ˆ
$D$ (cette extension est cyclique, donc abŽlienne; mais attention, pas forcŽment cyclique cyclotomique,
donc, on ne peut donc pas faire l'Žconomie de la Proposition \argba), l'ensemble $D$ a une densitŽ de
Dirichlet de ${1\over [L:K']}$. Donc $D'$ a aussi une densitŽ de ${1\over [L:K']}$ et finalement,  $A$ a
une densitŽ de Dirichlet $|C|\over [L:K]$.\qed
\bigskip

\rem

Dans la dŽmonstration prŽcŽdente, on voit que si $D$ est rŽgulier, alors $A$ aussi. Il suffit donc de
prouver que le thŽorme de $\check {\rm C}$ebotarev fort pour les extensions abŽlienne. La seule preuve que nous connaissons utilise un thŽorme fondamental de la thŽorie du corps de classe qui s'Žnonce comme suit~:
\bigskip
\lem

{\sl Soit $L/K$ une extension abŽlienne de corps de nombres. Alors il existe un $K$-module $\m$,
divisible par tous les idŽaux premiers de $O_K$ qui ramifient dans $L$, tel que 

$$P_\m\subset \ker (\Phi_{L/K}\, :\, I_K(\m)\rightarrow {\rm Gal}(L/K)).$$

}

{\bf Preuve}

C'est le ThŽorme \argda\ (qu'on appelle ``thŽorme de rŽciprocitŽ d'Artin''). On espre que le lecteur aura vŽrifiŽ qu'on n'a pas utilisŽ le
thŽorme de $\check {\rm C}$ebotarev fort pour les extensions abŽliennes pour prouver ce rŽsultat.\qed
\bigskip
\goodbreak
\coro

{\sl Le thŽorme $\check {\rm C}$ebotarev fort est vrai pour les extensions abŽliennes

}

{\bf Preuve}

Soit $C$ une classe de conjugaison de $G={\rm Gal}(L/K)$ abŽlien. Puisque par le lemme prŽcŽdent il existe $\m$ tel que $P_\m\subset \ker (\Phi_{L/K})$, on fait le mme (exactement le mme) raisonnement que pour la preuve du thŽorme de $\check {\rm C}$ebotarev fort pour les extensions cyclotomiques (Proposition \argaw) et on trouve que
l'ensemble $\{\P\in \gfP_{0}(K)\mid \P$ ne divise pas $\m$ et ${\rm Fr}_{L/K}=C\}$ est rŽgulier et de densitŽ ${1\over |G|}$. Comme on a vu au lemme prŽcŽdent que $\m$ peut tre divisible par tout les premiers qui ramifient dans $L$, l'ensemble $\{\P\in \gfP_0(K)\mid \P$ ne ramifie pas dans $L$ et ${\rm Fr}_{L/K}=C\}$ est aussi rŽgulier et de densitŽ ${1\over |G|}$ (la diffŽrence entre les deux ensembles est un ensemble fini). Cela prouve le corollaire.\qed
\bigskip
\centerline{\soustitre INTERLUDE}
\bigskip
Maintenant vient  quelques rŽsultats ``annexes''  qui ne sont pas directement nŽcessaires pour la preuve des grands rŽsultats que nous nous proposons de dŽmontrer, mais qui nous ont paru digne d'intŽrt. Voici un premier thŽorme dont l'ŽnoncŽ a probablement effleurŽ  chacun de nous au moins une fois dans sa vie et qui a peut-tre ŽtŽ ˆ la base du thŽorme de $\check {\rm C}$ebotarev.
\bigskip

\th

{\sl Soit $f\in \Z[X]$ irrŽductible (sur $\Z[X]$, ou sur $\Q[X]$, c'est Žgal gr‰ce au lemme de Gauss (cf. [La1, ThŽorme 4.2.3, p.191])) de degrŽ $n>1$. Soit $p\in\gfP_{0}(\Q)$ un nombre premier. Notons $\overline{f}$ la rŽduction modulo $p$ de $f$.
Alors les affirmations suivantes sont vraies~:

\art{a)}Il existe une infinitŽ de premiers $p$ pour lesquels $\overline{f}$ est totalement scindŽ dans
$\F_p[X]$ (i.e. est le produit de polyn™mes de degrŽ 1)

\art{b)}Il existe une infinitŽ de premiers $p$ tels que $\overline{f}$ n'a pas de racine dans $\F_p$.

\art{c)}Si $n$ est premier, il existe une infinitŽ de $p$ pour lesquels $\overline{f}$ est
irrŽductible dans $\F_p[X]$.

}

{\bf Preuve}

Posons $\theta=\theta_1,\theta_2,\ldots ,\theta_n$ les racines de $f$ dans $\C$. Soit $L/\Q$ le corps des
racines de $f$ (splitting field en anglais cf. [Jac1, p. 225 et suiv.] pour plus de dŽtails) et $E=\Q(\theta)$. Posons $G={\rm Gal}(L/\Q)$ et $H={\rm Gal}(L/E)$ ($=\{\sigma\in G\mid \sigma(\theta)=\theta\}$). L'application $\tau H\mapsto \tau (\theta)$ est une bijection bien dŽfinie de $\{\tau H\mid\tau\in G\}$ sur
$\{\theta_1,\ldots ,\theta_n\}$ (ils ont le mme nombre d'ŽlŽments et il y a surjection, $G$ agissant
transitivement sur les racines). C'est une bijection de $G$-ensembles (action ˆ gauche).

\art{a)} Soit $p$ un nombre premier ne divisant pas le discriminant de $f$ (c'est le discriminant de
$\Z[\theta]$ comme $\Z$-module qui vaut aussi $\prod_{1\leq i<j\leq n}(\theta_i-\theta_j)^2$); {\it a
fortiori}, $p$ ne divise pas le discriminant de
$E/\Q$, donc ne ramifie pas dans $E$, ni dans $L$ (par dŽfinition du corps des racines et gr‰ce au Lemme ``dŽcomposition-ramification'', page \the\gaagc). Soit $\P$ un idŽal premier de $L$ au-dessus de $p$. Posons $\overline{G}={\rm Gal}((O_L/\P)/\F_p)$, $\sigma={\rm Frob}(\P/p)$, $\overline{\theta_i}$ la rŽduction modulo $\P$ des
$\theta_i$ et $\overline{\sigma}$ le gŽnŽrateur de $\overline{G}$ (aussi appelŽ Frobenius). Puisque $p$ ne divise pas le discriminant de $f$ (ce qui veut dire que le discriminant de $\overline f$ est non nul) les
$\overline{\theta_i}$ sont tous disjoints. Il est clair que (par dŽfinition de $\sigma$)
$\overline{\sigma}$ permute les $\overline{\theta_i}$ de la mme manire que $\sigma$ permute les
$\theta_i$. Les racines de $\overline{f}$ qui sont dans $\F_p$ sont celles laissŽes fixes par
$\overline{\sigma}$. En particulier, $\overline{f}$ est scindŽ totalement dans $\F_p[X]$ si et seulement
si $\overline{\sigma}(\overline{\theta_i})=\overline{\theta_i}$ donc si et seulement si
$\sigma(\theta_i)=\theta_i$ pour $i=1,\ldots ,n$; ce qui veut dire que $\sigma=1$ dans $G$ qui est
Žquivalent au fait que $p$ se dŽcompose totalement dans $L$ (ou dans $E$, cf. Lemme ``dŽcomposition-ramification'', page \the\gaagc). Or, ces $p$ ont une densitŽ de $1\over [L:\Q]$ (cf. Lemme \argao). Ce qui prouve qu'ils sont une infinitŽ. Remarquons, qu'ici on n'a pas utilisŽ le thŽorme de $\check {\rm C}$ebotarev.

\art{b)}Sous les mmes notations qu'en a), dire que $\overline{f}$ n'a pas de racine  dans $\F_p$ revient
ˆ dire que $\overline{\sigma}(\overline{\theta_i})\ne \overline{\theta_i}$ ou encore
$\sigma(\theta_i)\ne\theta_i$ pour $i=1,\ldots ,n$. Cela veut dire que $\sigma(\tau H)\ne \tau H$, ou
encore $\sigma\not\in \tau H\tau^{-1}$ pour tout $\tau\in G$. On va montrer que de tels $\sigma$
existent toujours. Le nombre de classes des conjuguŽs de $H$ est l'indice dans $G$ du normalisateur de
$H$ (les $\tau$ tels que $\tau H\tau^{-1}=H$). Comme ce normalisateur contient $H$, cet indice est
$\leq {|G|\over |H|}$. Et alors (petite subtilitŽ avec l'identitŽ qui est dans toutes les classes)~: 

$$|\cup_{\tau}\tau H\tau^{-1}|\leq {|G|\over |H|}\cdot (|H|-1)+1=|G|-{|G|\over |H|} +1<|G|,$$

car ${|G|\over |H|}=n>1$. Donc de tels $\sigma$ existent. Appliquant le thŽorme de $\check {\rm
C}$ebotarev pour $C$, la classe de conjugaison d'un de ces sigmas, on voit qu'il existe une infinitŽ de $p$ tels
${\rm Fr}_{L/\Q}(p)=C$ ce qui montre la partie b).

\art{c)}Il est clair que $\overline{f}$ est irrŽductible $\Leftrightarrow$ $\overline{\sigma}$ agit
cycliquement sur les $\overline{\theta_i}$ (si $(\overline{\theta_1},\ldots ,\overline{\theta_r})$ est un
cycle de l'action sur les $\overline{\theta_i}$, le polyn™me $(X-\overline{\theta_1})\cdots
(X-\overline{\theta_r})$ est un polyn™me qui divise $\overline{f}$ et qui est dans
$\F_p[X]$, car laissŽ fixe par $\overline\sigma$). Donc $\overline{f}$ est irrŽductible $\Leftrightarrow$
$\sigma$ agit cycliquement sur les $\theta_i$ $\Rightarrow$ l'ordre de $\sigma$ est $n$. RŽciproquement,
si $\sigma$ est d'ordre $n$ et si de plus, $n$ est premier, alors $\sigma$ agit
nŽcessairement cycliquement sur les $\theta_i$ (en effet, l'ordre de $\sigma$ est Žgal au ppcm de
l'ordre des sous-cycles disjoints de la permutation des $\theta_i$ par l'action de $\sigma$, car les
$\theta_i$ engendrent $L$, et comme cet ordre est premier, cela veut dire qu'il n'y a qu'un cycle); donc,
$\overline{\sigma}$ agit cycliquement sur les $\overline{\theta_i}$, donc $\overline{f}$ est irrŽductible. Comme $n||G|$ et que $n$ est premier, le thŽorme de Cauchy (qui est un corollaire du premier thŽorme de Sylow (cf. [La1, Thm. 1.6.2, p. 35]), mais qu'on peut trs bien montrer pour lui-mme, mais je ne le ferai pas, ce n'est pas digne de vous si vous avez pu lire jusqu'ici) nous assure l'existence d'un ŽlŽment $\tau$ de $G$ d'ordre $n$. Chaque ŽlŽment de $C$, la classe de conjugaison de $\tau$, est aussi d'ordre~$n$. Le thŽorme de $\check {\rm C}$ebotarev appliquŽ ˆ $C$ nous assure donc une infinitŽ de $p\in\gfP_{0}(\Q)$ tels que le ${\rm Frob}(\P/p)$ agit cycliquement sur les $\theta_i$, ce qui veut dire que $\overline{f}$ est irrŽductible dans $\F_p[X]$.\qed

\bigskip   
\rem
\medskip
Dans la partie c) du thŽorme prŽcŽdent l'hypothse $n$ premier est cruciale. En effet, le polyn™me
$X^4+1$ est irrŽductible dans $\Z[X]$ (c'est le $8^e$ polyn™me cyclotomique) et il n'est irrŽductible
dans aucun~$\F_p$~: 

\art{a)} Dans $\F_2$, il vaut  $(X+1)^4$, 

\art{b)}Si $p\equiv 1\pmod 8$, il est totalement scindŽ. En effet, on a vu ˆ
la partie a) du thŽorme prŽcŽdent que $X^4+1$ se scinde totalement si et seulement si $p$ se
dŽcompose totalement dans $\Q[X]/(X^4+1)\simeq \Q(\zeta_8)$, si et seulement si $f(\P/p)=1$ pour tout
idŽal $\P$ au-dessus de $p$. Or, $f(\P/p)$ est l'ordre de $p$ modulo~8. En effet, ${\rm Frob}_{\Q(\zeta_{8})/\Q}(p)=:\sigma_{p}$, est tel que $\sigma_{p}(\zeta_{8})=\zeta_{8}^p$ (cf. preuve du ThŽorme \argm); l'ordre de $\sigma_{p}$ est  l'ordre de $p$ modulo 8 mais vaut aussi $f(\P/p)$ par dŽfinition. Dans notre cas, cet ordre vaut justement 1, ce qui montre l'affirmation.

\art{c)}Si $p\not \equiv 1\pmod 8$, on a tout de mme $p^2\equiv 1\pmod 8$. Cela veut dire qu'ici
$f(\P/p)$, qui est l'ordre de ${\rm Frob}_{\Q(\zeta_8)/\Q}(p)$ vaut $2$, il ne peut donc pas agir
cycliquement sur les racines de $X^4+1$, ce qui veut dire que $\overline{X^4+1}$ n'est pas irrŽductible.
\bigskip

Le second rŽsultat que nous allons prŽsenter est une gŽnŽralisation du ThŽorme \argar\ et quelques extras.
\bigskip\goodbreak
\defi
\medskip
Soit $L/K$ une extension de corps de nombres. Notons $S(L/K)=\{\P\in \gfP_0(K)\mid \P$ se dŽcompose compltement dans $L\}$. Rappelons qu'on a montrŽ au Lemme \argao\ que $\delta(S(L/K))={1\over [E:K]}$, o $E$ est la cl™ture galoisienne de $L/K$. On pose ensuite $\widetilde{S}(L/K)=\{P\in \gfP_0(K)\mid \P$ ne ramifie pas dans $L$ et il existe (au moins) un premier $\gP$ de $L$ au-dessus de $\P$ avec $f(\gP/\P)=1\}$. Il est clair que si $L/K$ est une extension galoisienne, alors $S(L/K)=\widetilde{S}(L/K)$

\newcount\gaags\gaags=\gaga
\bigskip

\th
\medskip
{\sl Soit $L/K$ une extension de corps de nombres. Alors $\widetilde{S}(L/K)$ a une densitŽ de Dirichlet

$$\delta(\widetilde{S}(L/K))\geq {1\over [L:K]},$$

avec ŽgalitŽ si et seulement si $L/K$ est galoisienne.

}

{\bf preuve}

Soit $E/K$ la cl™ture galoisienne de $L/K$. Notons $G={\rm Gal}(E/K)$ et $H={\rm Gal}(E/L)\subset G$. Soit $\P$ un idŽal premier de $K$ non ramifiŽ dans $L$. Par le Lemme "dŽcomposition-ramification" de la page \the\gaagc, nous savons que $\P$ ne ramifie pas non plus dans $E$. Soit encore $\gP$ un idŽal premier au-dessus de $\P$ et $\sigma={\rm Frob}_{E/K}(\gP/\P)$. Par le ThŽorme \aargp, $\P\in \widetilde{S}(L/K)$ si et seulement s'il existe $\tau\in G$ tel que $H\cdot \tau\cdot\sigma=H\cdot\tau$, i.e. $\sigma\in\tau^{-1}\cdot H\cdot\tau$. Donc $\P\in \widetilde{S}(L/K)$ si et seulement si $ \sigma \in \bigcup_{\tau\in G}\tau^{-1}\cdot H\cdot\tau$. En observant que $ \bigcup_{\tau\in G}\tau^{-1}\cdot H\cdot\tau=\bigcup_{h\in H} C_h$, o $C_h$ est la classe de conjugaison de $h$, en s'inspirant de la preuve du Corollaire \argap, et gr‰ce au thŽorme de $\check {\rm C}$ebotarev, on en dŽduit que $\widetilde{S}(L/K)$ a une densitŽ de Dirichlet qui vaut 
$${\left |\bigcup_{\tau\in G}\tau^{-1}\cdot H\cdot\tau\right |\over |G| }\geq {|H|\over |G|}={1\over [L:K]},$$
avec ŽgalitŽ si et seulement si $H$ est normal dans $G$, si et seulement si $L/K$ est galoisienne.\qed

\bigskip

\defi
\medskip

Soit $L/K$ une extension de corps de nombres. On dira qu'un idŽal premier $\P$ de $K$ non ramifiŽ dans $L$ a une dŽcomposition du type $(f_1,\ldots ,f_r)$, avec $f_1\leq f_2\leq \cdots\leq f_r$, si $\P\cdot O_L=\gP_1\cdots\gP_r$, avec $f(\gP_i/\P)=f_i$, pour $i=1,\ldots ,r$.

\bigskip
\th
\medskip
{\sl Soit $L/K$ une extension de corps de nombres. Posons $A=\{\P\in\gfP_0(K)\mid \P$ a une dŽcomposition du type $(f_1,\ldots ,f_r)\}$. Alors 
$$A\ne\emptyset\iff A\hbox{ a une densitŽ de Dirichlet strictement positive, et donc $|A|=\infty$.}$$

}
{\bf Preuve}

Prouvons la partie $\Rightarrow$ (l'autre Žtant bien sžr triviale). Soit $E/K$ la cl™ture galoisienne de $L/K$. Notons $G={\rm Gal}(E/K)$ et $H={\rm Gal}(E/L)$. Soit aussi $\P\in A$ (supposŽ non vide). Soit encore $\gP$ un idŽal de $E$ au-dessus de $\P$ et $\sigma={\rm Frob}_{E/K}(\gP/\P)$. Gr‰ce au ThŽorme \aargp, on sait qu'il existe, pour $i=1,\ldots , r$, $\tau_i\in {\rm Gal}(E/K)$ tel que les $C_i=\{H\cdot\tau_i,H\cdot\tau_i\cdot\sigma,\ldots ,H\cdot\tau_i\cdot\sigma^{f_i-1}\}$ forment les orbites de l'action de $<\sigma>$ sur les classes ˆ droite de $G$ modulo $H$ et les $\P_i:=\tau_i(\gP)\cap L$ sont exactement les idŽaux premiers de $E$ au-dessus de $\P$. Et finalement $f(\P_i/\P)=f_i$, pour $i=1,\ldots,r$. ConsidŽrons $C=\{\rho^{-1}\cdot\sigma\cdot\rho\mid \rho\in G\}$ la classe de conjugaison de $\sigma$. Par le thŽorme de $\check {\rm C}$ebotarev (ThŽorme \argav), l'ensemble $B:=\{ \qq\in\gfP_0(K)\mid {\rm Fr}_{E/K}(\qq)=C\}$ est de densitŽ de Dirichlet $\delta(B)>0$. Soit $\qq\in B$ et $\gQ\in \gfP_0(E)$ au dessus de $\qq$. Par dŽfinition de $B$, il existe $\rho\in G$ tel que ${\rm Frob}_{E/K}(\gQ/\qq)=\rho^{-1}\cdot\sigma\cdot\rho$, i.e. ${\rm Frob}_{E/K}(\rho(\gQ)/\qq)=\sigma={\rm Frob}_{E/K}(\gP/\P)$. Cela veut dire que les orbites de l'action de $<\sigma>$ sont Žvidemment les mmes (= les $C_i$), et donc gr‰ce au ThŽorme \aargp\ que $\qq$ a une dŽcomposition du mme type que $\P$ dans $E$. Ainsi, $\qq\in A$.  On a montrŽ que $B\subset A$ et donc $0<\delta(B)\leq\delta(A)$.\qed
\bigskip
Maintenant nous allons donner une gŽnŽralisation du ThŽorme \argar\ :
\bigskip
\th

{\sl Soit $K$ un corps de nombres et $L/K$ et $M/K$ des extensions finie de $K$. 

\art{a)}Supposons que $M/K$ soit une extension galoisienne. Alors on a~:
$$L\subset M\iff S(M/K)\subset S(L/K)\iff \delta( S(M/K)\setminus S(L/K))=0.$$

\art{b)}Supposons que $L/K$ soit une extension galoisienne. Alors on a~:
$$L\subset M\iff \widetilde{S}(M/K)\subset S(L/K)\ (= \widetilde{S}(L/K))\iff \delta(  \widetilde{S}(M/K)\setminus S(L/K))=0.$$

}

{\bf Preuve}

Toute les implications $\Rightarrow$ sont triviales.

Montrons la partie b). Il suffit de prouver que $\delta(  \widetilde{S}(M/K)\setminus S(L/K))=0$ et  $L/K$ abŽlienne $\Rightarrow L\subset M$. Soit $E/K$ une extension galoisienne qui contient $L$ et $M$. En vertu de la correspondance de Galois, il suffit donc de montrer que ${\rm Gal}(E/M)\subset {\rm Gal}(E/L)$, ou encore que $\sigma|_L={\rm Id}_L$ pour tout $\sigma\in {\rm Gal}(E/M)$. Soit donc  $\sigma\in {\rm Gal}(E/M)\subset {\rm Gal}(E/K)$. Soit $C=\{\tau\sigma\tau^{-1}\mid \tau\in {\rm Gal}(E/K)\}$, la classe de conjugaison de $\sigma$. Par le thŽorme de $\check {\rm C}$ebotarev, l'ensemble des $\P\in \gfP_{0}(K)$ tels que ${\rm Fr}_{E/K}(\P)=C$ est de densitŽ strictement positive. Soit donc un tel $\P$ et $\gP\in\gfP_{0}(E)$ tel que ${\rm Frob}_{E/K}(\gP/\P)=\sigma$ et notons $\gP'=\gP\cap O_{M}$. On sait que $\sigma |_{M}={\rm Frob}_{M/K}(\gP'/\P)$. Or, puisque $\sigma\in{\rm Gal}(E/M)$, on a donc ${\rm Frob}_{M/K}(\gP'/\P)={\rm Id}_{M}$, i.e. $f(\gP'/\P)=1$ et donc $\P\in \widetilde{S}(M/K)$. Par hypothse, $\delta(  \widetilde{S}(M/K)\setminus S(L/K))=0$, donc puisque l'ensemble dans lequel nous avons ŽtŽ puiser $\P$ est de densitŽ strictement positive, on peut supposer que $\P\in S(L/K)$. Donc ${\rm Frob}_{L/K}(\gP''/\P)={\rm Id}_{L}$, o $\gP''=\gP\cap O_{L}$. Mais ${\rm Frob}_{L/K}(\gP''/\P)=\sigma|_{L}$, donc $\sigma|_{L}={\rm Id}_{L}$, ce qui montre la partie b).

Prouvons la partie a). Par hypothse, on a $ \delta( S(M/K)\setminus S(L/K))=0$ et $M/K$ galoisienne. On doit montrer que $L\subset M$. Soit $L'/K$ la cl™ture galoisienne de $L/K$. En vertu du Lemme ÒdŽcomposition-ramification'', page \the \gaagc, on a que $S(L/K)=S(L'/K)$. De plus, puisque $M/K$ est galoisienne, on a $S(M/K)=\widetilde{S}(M/K)$. L'hypothse s'Žcrit alors $\delta(  \widetilde{S}(M/K)\setminus S(L'/K))=0$. La partie b) nous montre alors que $L'\subset M$ et donc $L\subset L'\subset M$.\qed

\bigskip
\rem

\medskip
La partie a) de ce thŽorme est connue sous le nom de ÒThŽorme de Bauer".

\vfill\eject

\global\advance\chapnomb by 1
\nomb=1

\centerline{\para Chapitre 4 : }
\medskip
\centerline{\para Cohomologie des groupes cycliques et quotient de Herbrand}
\bigskip

Ici, nous allons donner quelques rudiments de cohomologie cyclique en dŽtail pour introduire le quotient de Herbrand. Le lecteur expŽrimentŽ (ou pressŽ) voudra bien ne retenir que les Lemmes \argbi\ ˆ \argbm, et et passer ˆ la section ``{\bf Calculs explicites dans le cas d'extensions cycliques}''
\medskip

Soit $G$ un groupe notŽ multiplicativement et $A$ un $G$-module, c'est-ˆ-dire un groupe abŽlien $(A,*)$
muni d'une action $G\times A\rightarrow A$, $(\sigma,a)\mapsto \sigma(a)$ avec la propriŽtŽ que
$Id_G(a)=1(a)=a$, $(\sigma\tau)(a)=\sigma(\tau(a))$ et $\sigma(a*b)=\sigma(a)*\sigma(b)$. Chaque ŽlŽment
$\sigma$ de $G$ est associŽ ˆ un unique ŽlŽment de ${\rm End}_\Z(A)$ qu'on notera encore $\sigma$.  Puisque
$A$ est un groupe abŽlien, on le munit aussi via cette action d'une structure de $\Z[G]$-module~:
$(\sigma_1+\sigma_2)(a):=\sigma_1(a)*\sigma_2(a)$. 

Supposons ˆ partir de maintenant que $G=<\sigma>$ est cyclique d'ordre $n$. Posons

$$\Delta=1-\sigma\hskip 0.5cm\hbox{et}\hskip0.5cm N=1+\sigma+\sigma^2+\cdots +\sigma^{n-1}\in\Z[G].$$ 

Posons aussi~:

$$\eqalign{\Delta |A\, :\, A&\longrightarrow A\cr a&\longmapsto a-\sigma(a)\ \hbox{ en notation
additive}\cr a&\longmapsto {a\over \sigma(a)}\ \hbox{ en notation
multiplicative,}\cr}$$

et 

$$\eqalign{N |A\, :\, A&\longrightarrow A\cr a&\longmapsto \sum_{i=0}^{n-1}\sigma^i(a)\ \hbox{ en
notation additive}\cr a&\longmapsto \prod_{i=0}^{n-1}\sigma^i(a)\ \hbox{ en notation
multiplicative.}\cr}$$
\headline={\hfill \phantom{ouuh}\hfill}
\newcount\gaagt\gaagt=\gaga

Par exemple, si $L/K$ est une extension cyclique de groupe $G$, si $A=L^*$, alors $N|A=N_{L/K}$ est la
norme usuelle et si $A=L$, alors $N|A={\rm Tr}_{L/K}$ est la trace usuelle.

Puisque $G$ est cyclique d'ordre $n$, il est Žvident que $\Delta N=N\Delta=0$; donc ${\rm
Im}(\Delta)\subset {\rm ker}(N)$ et ${\rm Im}(N)\subset\ker \Delta$. Remarquons aussi que $\Delta$ dŽpend
de $\sigma$, mais si on prend un autre gŽnŽrateur de $G$, il a la mme image et le mme noyau (tout cela
vient de la relation $(1-\sigma)(1+\sigma+\cdots +\sigma^{i-1})=1-\sigma^i$). On dŽfinit alors 

$$\thboxed 15{H^0(A)=\ker(\Delta|A)/N(A)\hskip0.5cm \hbox{ et }\hskip0.5cm H^1(A)=\ker(N|A)/\Delta(A)}$$

Si $f\ :A\rightarrow B$ est un homomorphisme de $G$-module (ce qui veut dire $f(\sigma(a))=\sigma(f(a))$ pour tout $a\in A$), alors on a $f\Delta=\Delta f$ et $fN=Nf$. On
en dŽduit que $f(\ker (\Delta |A))\subset \ker (\Delta |B)$ et $f(\Delta(A))\subset \Delta(B)$; les mmes
relations sont vŽrifiŽes si on remplace $\Delta$ par $N$. Cela implique que $f$ induit des homomorphismes
$f_i\ : H^i(A)\rightarrow H^i(B)$.

Sauf mention express, les $G$-modules seront notŽs additivement.
\bigskip\goodbreak
\lem \ {\bf (Lemme de l'hexagone)}
\medskip
{\sl Soit $0\rightarrow A\buildrel f\over \rightarrow B \buildrel g\over \rightarrow C\rightarrow 0$ une
suite exacte de $G$-module. Alors il existe des homomorphismes $\delta_0\ : H^0(C)\rightarrow H^1(A)$ et
$\delta_1\ : H^1(C)\rightarrow H^0(A)$ tels que l'hexagone suivant soit exact partout 
\goodbreak
\vglue 2cm

\hbox{\hskip 5cm
\rput(0,0){\rnode{h1c}{$H^1(C) $}} 
\rput(2,1){\rnode{h0a}{$H^0(A)$}}
\rput(4,1){\rnode{h0b}{$H^0(B)$}}
\rput(6,0){\rnode{h0c}{$H^0(C)$}}
\rput(4,-1){\rnode{h1a}{$H^1(A)$}}
\rput(2,-1){\rnode{h1b}{$H^1(B)$}}
\ncline[nodesep=3pt]{->}{h1c}{h0a}
\Aput{$\delta_{1}$} 
\ncline[nodesep=3pt]{->}{h0a}{h0b}
\Aput{$f_{0}$} 
\ncline[nodesep=3pt]{->}{h0b}{h0c}
\Aput{$g_{0}$} 
\ncline[nodesep=3pt]{->}{h0c}{h1a}
\Aput{$\delta_{0}$} 
\ncline[nodesep=3pt]{->}{h1a}{h1b}
\Aput{$f_{1}$} 
\ncline[nodesep=3pt]{->}{h1b}{h1c}
\Aput{$g_{1}$} }


\vskip 1.5cm

}

\goodbreak

{\bf Preuve}

DŽfinissons $\delta_0$~: le diagramme suivant est commutatif~:

 \vglue 2cm

\psset{xunit=0.8cm,yunit=0.8cm}

\hbox{\hskip 4cm
\rput(0,0){\rnode{zer1}{$0 $}} 
\rput(0,2){\rnode{zer2}{$0$}}
\rput(2,0){\rnode{a1}{$A$}}
\rput(2,2){\rnode{a2}{$A$}}
\rput(4,0){\rnode{b1}{$B$}}
\rput(4,2){\rnode{b2}{$B$}}
\rput(6,0){\rnode{c1}{$C$}}
\rput(6,2){\rnode{c2}{$C$}}
\rput(8,0){\rnode{zer3}{$0 $}} 
\rput(8,2){\rnode{zer4}{$0$}}
\ncline[nodesep=3pt]{->}{zer1}{a1}
\ncline[nodesep=3pt]{->}{zer2}{a2}
\ncline[nodesep=3pt]{->}{a2}{a1}
\Aput{$\Delta |A $} 
\ncline[nodesep=3pt]{->}{b2}{b1}
\Aput{$\Delta |B$} 
\ncline[nodesep=3pt]{->}{c2}{c1}
\Aput{$\Delta |C$} 
\ncline[nodesep=3pt]{->}{a1}{b1}
\Aput{$f $} 
\ncline[nodesep=3pt]{->}{a2}{b2}
\Aput{$f $} 
\ncline[nodesep=3pt]{->}{b1}{c1}
\Aput{$g$} 
\ncline[nodesep=3pt]{->}{b2}{c2}
\Aput{$g$} 
\ncline[nodesep=3pt]{->}{c1}{zer3}
\ncline[nodesep=3pt]{->}{c2}{zer4}}


\vskip 2cm

\headline={\hfill \smcap Cohomologie des groupes cycliques et quotient de Herbrand\hfill}

Soit $c+N(C)\in H^0(C)$. Alors $\Delta(c)=0$ (puisque par hypothse $c\in\ker(\Delta|C)$). Nous allons
Òchasser dans le diagramme"~: puisque $g$ est surjective, il existe $b\in B$ tel que $g(b)=c$. Par
commutativitŽ du diagramme et puisque $\Delta(c)=0$, on a $g(\Delta(b))=0$; il existe donc $a\in A$
tel que $f(a)=\Delta(b)$. D'autre part, $0=N(\Delta(b))=N(f(a))=f(N(a))$; cela implique que $N(a)=0$, par
injectivitŽ de $f$. On posera donc 
$$\delta_0(c+N(C))=a+\Delta(A)\in H^1(A).$$
 Reste ˆ voir que si $c'+N(C)=c+N(C)$, alors $a-a'\in\Delta(A)$, avec $a,a'$ et $b,b'$ dŽfinis par $g(b)=c$,
$\Delta(b)=f(a)$ et $g(b')=c'$, $\Delta(b')=f(a')$. Par hypothse, $c-c'=N(d)$ pour un $d\in C$; soit $b''$
tel que $d=g(b'')$ ($g$ est surjective). On a $g(b-b'-N(b''))=0$, donc $b-b'-N(b'')=f(a'')$ pour un
$a''\in A$. Appliquons $\Delta$ ˆ cette dernire ŽgalitŽ. On trouve $f(a-a')=f(a)-f(a')=f(\Delta(a''))$,
ce qui implique, par injectivitŽ de $f$ que $a-a'=\Delta(a'')$. Cela montre que $\delta_0$ est bien
dŽfinie et on vŽrifie que c'est un homomorphisme. 

En Žchangeant les r™les de $\Delta$ et de $N$, on a la mme preuve pour dŽfinir $\delta_1$. 

Regardons l'exactitude de l'hexagone pour les fonctions ayant un indice $0$, les autres se dŽduisent en
permutant $N$ et $\Delta$. Cette vŽrification se fait en 6 Žtapes~:

\art{i)}Puisque $g\circ f=0$, on a $g_0\circ f_0=0$ ce qui implique que ${\rm Im}(f_0)\subset\ker(g_0)$.

\art{ii)}Montrons l'inclusion inverse. Soit $b+N(B)\in\ker(g_0)$. Alors $b\in\ker(\Delta |B)$ et $g(b)\in
N(C)$, posons donc $g(b)=N(c)$, $c\in C$. Il faut montrer que $b\equiv f(a)\pmod{N(B)}$, avec $a\in\ker(\Delta|A)$. Il
faut donc trouver $b_0\in B$ tel que $b=f(a)+N(b_0)$. Soit $b_{0}\in B$ tel que $c=g(b_0)$ ($g$ est surjective). On a donc
$g(b)=N(c)=N(g(b_0))=g(N(b_0))$ ce qui veut dire que $b-N(b_0)\in \ker(g)={\rm Im}(f)$. Il existe donc
$a\in A$ tel que $b=f(a)+N(b_0)$; appliquant $\Delta$ ˆ cette dernire ŽgalitŽ, on trouve
$0=\Delta(f(a))+0=f(\Delta(a))$, ce qui implique ($f$ est injective) que $a\in \ker(\Delta | A)$. On a donc
bien $b+N(B)=f(a)+N(B)$, avec $a\in \ker(\Delta |A)$. Cela prouve que $\ker(g_0)\subset {\rm Im}(f_0)$.

\art{iii)} Montrons que ${\rm Im}(g_0)\subset\ker(\delta_0)$. Il faut appliquer $\delta_0$ ˆ un ŽlŽment
de la forme $g(b)+N(C)$ avec $b\in \ker(\Delta|B)$. Par dŽfinition, $\delta_0(g(b)+N(C))=a+\Delta(A)$,
avec $f(a)=\Delta(b)=0$, donc, $a=0$, par injectivitŽ de $f$. Cela prouve que
$\delta_0(g(b)+N(C))=0+\Delta(A)\in \ker(\delta_{0})$.

\art{iv)} Montrons que $\ker(\delta_0)\subset {\rm Im}(g_0)$. Soit $c\in \ker(\Delta |C)$ tel que
$c+N(C)\in \ker(\delta_0)$. Cela veut dire qu'il existe $b\in B$ et $a\in A$ tels que
$g(b)=c$, $\Delta(b)=f(a)$ et $a\in \Delta(A)$; disons $a=\Delta(a')$, avec $a'\in A$. Il suffit de
montrer que $c=g(b')$, avec $\Delta(b')=0$. On a $\Delta(b)=f(a)=f(\Delta(a'))=\Delta(f(a'))$; posons
alors $b'=b-f(a) $. On a $\Delta(b')=0$ et $g(b')=g(b)=c$, ce qui montre que $c+N(C)=g_0(b'+N(B))$. 

\art{v)} Montrons que ${\rm Im}(\delta_0)\subset \ker(f_1)$. Soit $c\in \ker(\Delta|C)$. On a
$f_1(\delta_{0}(c+N(C)))=f(a)+\Delta(B)$ o $a$ est tel que $f(a)=\Delta(b)$ et $g(b)=c$; on en dŽduit donc
que $f(a)\in \Delta(B)$ et donc $f_1(\delta_{0}(c+N(C)))=0+\Delta(B)$.

\art{vi)}Il reste ˆ voir que $\ker(f_1)\subset {\rm Im}(\delta_0)$. Soit $a\in\ker(N |A)$ tel que
$f_1(a+\Delta(a))=0$, i.e. tel que $f(a)\in \Delta(B)$. Alors $f(a)=\Delta(b)$ pour un $b\in B$. Posons
$c=g(b)$. On a $\Delta(c)=\Delta(g(b))=g(\Delta(b))=g(f(a))=0$. Donc $c+N(C)\in H^0(C)$ et
$\delta_0(c+N(C))=a+\Delta(A)$.

Cela achve la preuve de ce lemme.\qed

\bigskip

\lem
\medskip
{\sl Si $A\buildrel f\over \rightarrow B \buildrel g\over \rightarrow C$ est une composition de deux
homomorphismes de $G$-modules, alors on a $(g\circ f)_{i}=g_{i}\circ f_{i}$ ($i=0,1$) et si $f$ est
un isomorphisme alors $f_i$ aussi ($i=0,1$) (il en est de mme pour $g$). D'autre part, si $A\buildrel
f,g\over \rightarrow B$ sont des homomorphismes de $G$-modules, alors $(f+g)_i=f_i+g_i$  ($i=0,1$). Enfin,
si on a le diagramme commutatif suivant (les deux lignes Žtant exactes)~:

\vskip1.9cm


\psset{xunit=0.8cm,yunit=0.8cm}

\hbox{\hskip 5cm\rput(0,0){\rnode{zer1}{$0 $}} 
\rput(0,2){\rnode{zer2}{$0$}}
\rput(2,0){\rnode{a1}{$A'$}}
\rput(2,2){\rnode{a2}{$A$}}
\rput(4,0){\rnode{b1}{$B'$}}
\rput(4,2){\rnode{b2}{$B$}}
\rput(6,0){\rnode{c1}{$C'$}}
\rput(6,2){\rnode{c2}{$C$}}
\rput(8,0){\rnode{zer3}{$0 $}} 
\rput(8,2){\rnode{zer4}{$0$}}
\ncline[nodesep=3pt]{->}{zer1}{a1}
\ncline[nodesep=3pt]{->}{zer2}{a2}
\ncline[nodesep=3pt]{->}{a2}{a1}
\Aput{$h $} 
\ncline[nodesep=3pt]{->}{b2}{b1}
\Aput{$i $} 
\ncline[nodesep=3pt]{->}{c2}{c1}
\Aput{$j $} 
\ncline[nodesep=3pt]{->}{a1}{b1}
\Aput{$f' $} 
\ncline[nodesep=3pt]{->}{a2}{b2}
\Aput{$f $} 
\ncline[nodesep=3pt]{->}{b1}{c1}
\Aput{$g'$} 
\ncline[nodesep=3pt]{->}{b2}{c2}
\Aput{$g$} 
\ncline[nodesep=3pt]{->}{c1}{zer3}
\ncline[nodesep=3pt]{->}{c2}{zer4}}
\medskip
Alors le prisme ˆ base hexagonale suivant est aussi exact et commutatif~:

\vglue 1.2cm
\psset{xunit=1cm,yunit=1cm}

\hbox{\hskip5cm
\rput(0,-1.2){\rnode{h1cpr}{$H^1(C') $}} 
\rput(2,-0.2){\rnode{h0apr}{$H^0(A')$}}
\rput(4,-0.2){\rnode{h0bpr}{$H^0(B')$}}
\rput(6,-1.2){\rnode{h0cpr}{$H^0(C')$}}
\rput(4,-2.2){\rnode{h1apr}{$H^1(A')$}}
\rput(2,-2.2){\rnode{h1bpr}{$H^1(B')$}}
\ncline[nodesep=2pt]{->}{h1cpr}{h0apr}
\Bput{$\delta_{1}'$} 
\ncline[nodesep=3pt]{->}{h0apr}{h0bpr}
\Aput{$f_{0}'$} 
\ncline[nodesep=3pt]{->}{h0bpr}{h0cpr}
\Bput{$g_{0}'$} 
\ncline[nodesep=3pt]{->}{h0cpr}{h1apr}
\Aput{$\delta_{0}'$} 
\ncline[nodesep=3pt]{->}{h1apr}{h1bpr}
\Aput{$f_{1}'$} 
\ncline[nodesep=3pt]{->}{h1bpr}{h1cpr}
\Aput{$g_{1}'$} 
\rput(0,0){\rnode{h1c}{$H^1(C) $}} 
\rput(2,1){\rnode{h0a}{$H^0(A)$}}
\rput(4,1){\rnode{h0b}{$H^0(B)$}}
\rput(6,0){\rnode{h0c}{$H^0(C)$}}
\rput(4,-1){\rnode{h1a}{$H^1(A)$}}
\rput(2,-1){\rnode{h1b}{$H^1(B)$}}
\ncline[nodesep=3pt]{->}{h1c}{h0a}
\Aput{$\delta_{1}$} 
\ncline[nodesep=3pt]{->}{h0a}{h0b}
\Aput{$f_{0}$} 
\ncline[nodesep=3pt]{->}{h0b}{h0c}
\Aput{$g_{0}$} 
\ncline[border=2pt, nodesep=3pt]{->}{h0c}{h1a}
\Bput{$\delta_{0}$} 
\ncline[nodesep=3pt]{->}{h1a}{h1b}
\Aput{$f_{1}$} 
\ncline[border=2pt , nodesep=3pt]{->}{h1b}{h1c}
\Bput{$g_{1}$} 
\ncline[nodesep=3pt]{->}{h1c}{h1cpr}
\Bput{$j_{1}$} 
\ncline[nodesep=3pt]{->}{h0a}{h0apr}
\Bput{$h_{0}$} 
\ncline[nodesep=3pt]{->}{h0b}{h0bpr}
\Bput{$i_{0}$} 
\ncline[nodesep=3pt]{->}{h0c}{h0cpr}
\Aput{$j_{0}$} 
\ncline[nodesep=3pt]{->}{h1a}{h1apr}
\Bput{$h_{1}$} 
\ncline[nodesep=3pt]{->}{h1b}{h1bpr}
\Bput{$i_{1}$} }
\vglue 3cm

}

{\bf Preuve}

C'est une consŽquence directe des dŽfinitions et du Lemme de l'hexagone.\qed
\bigskip\goodbreak

\defi
\medskip
Soit $A$ un $G$-module tel que $H^0(A)$ et $H^1(A)$ sont finis. On appelle {\it quotient de
Herbrand} le rapport 
$$q(A)={|H^1(A)|\over |H^0(A)|}$$

\bigskip
\newcount\gaagu\gaagu=\gaga

\lem
\medskip
{\sl Soit $0\rightarrow A\buildrel f\over \rightarrow B \buildrel g\over \rightarrow C\rightarrow 0$ une
suite exacte de $G$-module. Si deux des quotients $q(A), q(B), q(C)$ sont dŽfinis, alors le troisime
aussi, et dans ce cas, on a

$$q(B)=q(A)\cdot q(C).$$

}
\goodbreak
{\bf preuve}

Regardons encore une fois le diagramme du Lemme de l'hexagone~: 

\vskip 1cm

\hbox{\hskip5cm
\rput(0,0){\rnode{h1c}{$H^1(C) $}} 
\rput(2,1){\rnode{h0a}{$H^0(A)$}}
\rput(4,1){\rnode{h0b}{$H^0(B)$}}
\rput(6,0){\rnode{h0c}{$H^0(C)$}}
\rput(4,-1){\rnode{h1a}{$H^1(A)$}}
\rput(2,-1){\rnode{h1b}{$H^1(B)$}}
\ncline[nodesep=3pt]{->}{h1c}{h0a}
\Aput{$\delta_{1}$} 
\ncline[nodesep=3pt]{->}{h0a}{h0b}
\Aput{$f_{0}$} 
\ncline[nodesep=3pt]{->}{h0b}{h0c}
\Aput{$g_{0}$} 
\ncline[nodesep=3pt]{->}{h0c}{h1a}
\Aput{$\delta_{0}$} 
\ncline[nodesep=3pt]{->}{h1a}{h1b}
\Aput{$f_{1}$} 
\ncline[nodesep=3pt]{->}{h1b}{h1c}
\Aput{$g_{1}$} }

\vskip 2cm

Il est clair que par exemple si $q(A)$ et $q(B)$ sont dŽfinis, alors $q(C)$ l'est aussi, car
$|H^1(C)|=|\ker(\delta_1)||{\rm Im}(\delta_1)|\buildrel \rm exactitude\over = |{\rm Im}(g_1)||{\rm
Im}(\delta_1)|\leq |H^1(B)||H^0(A)|<\infty$. On montre de la mme manire que $|H^1(C)|<\infty$; donc
$q(C)$ est dŽfini. Dans la mme veine de raisonnement, on a 
$$|H^0(A)||H^0(C)||H^1(B)|=|\ker(f_0)||{\rm Im}(f_0)||\ker(\delta_0)||{\rm Im}(\delta_0)||\ker(g_1)||{\rm
Im}(g_1)|$$
et
$$\eqalign{|H^1(A)||H^1(C)||H^0(B)|&=|\ker(f_1)||{\rm Im}(f_1)||\ker(\delta_1)||{\rm
Im}(\delta_1)||\ker(g_0)||{\rm Im}(g_0)|\cr
&\buildrel \rm exactitude\over = |{\rm Im}(\delta_0)||\ker(g_1)||{\rm Im}(g_1)||\ker(f_0)||{\rm
Im}(f_0)||\ker(\delta_0)|.\cr}$$
En quotientant la seconde ŽgalitŽ par la premire, on trouve $q(A)\cdot q(C)\cdot q(B)^{-1}$ ˆ gauche et 1 ˆ droite,
ce qui montre notre lemme.\qed
\bigskip

\coro

Si $A=A_1\oplus\cdots \oplus A_n$ est une somme directe de $G$-modules, alors pour $i=0,1$, on a
$H^i(A)=H^i(A_1)\oplus\cdots\oplus H^i(A_n)$, ainsi 
$$q(A)=\prod_{i=1}^n q(A_i).$$
$\phantom{vive les maths}$\qed
\bigskip

\lem
\medskip

{\sl Si $A$ est un $G$-module fini, alors $q(A)=1$

}

{\bf Preuve}

On a

$$q(A)={[\ker(N|A):{\rm Im}(\Delta|A)]\over [\ker(\Delta|A):{\rm Im}(N|A)]}={|\ker(N|A)|\over |{\rm Im}
(\Delta |A)|}\cdot{|{\rm Im}(N|A)|\over|\ker(\Delta|A)|}={|A|\over |A|}=1.$$
$\phantom{vive les maths}$\qed
\bigskip

\coro
\medskip

{\sl  Soit $A$ un $G$-module et $B$ un sous-$G$-module de $A$ tels que $[A:B]<\infty$. Alors $q(A)$ est dŽfini si et seulement si $q(B)$ est dŽfini et alors, $q(A)=q(B)$.

}

{\bf Preuve}

On a la suite exacte

$$0\longrightarrow B\buildrel\rm incl.\over \longrightarrow A\buildrel\rm proj.\over\longrightarrow A/B \longrightarrow 0$$
et on conclut par le Lemme \argbi\ et le Lemme \argbk.\qed

\goodbreak\bigskip

\prop
\medskip

{\sl Supposons que $G=<\sigma>$, avec $n=|G|=m\cdot d$. Soit $R$ un anneau intgre de caractŽristique
0. Soit $A=\bigoplus_{i=1}^d R u_i$, un $R$-module de base $u_1,\ldots , u_d$. On suppose que  $G$ agit
par permutation des $u_i$ de la manire suivante~: $u_{i+1}=\sigma (u_i)$, pour $i=1,\ldots ,d-1$ et
$u_1=\sigma(u_d)$ (on dira alors que $A$ est un {\it module de permutation}). Le $R$-module $A$ devient
ainsi un
$G$-module. On suppose de plus que
$R/mR$ est fini. Alors $q(A)$ est dŽfini et 
$$q(A)=[R:mR]^{-1}.$$
En particulier, si $R=\Z$, $q(A)={1\over m}$.

}

{\bf Preuve}

Remarquons que $N(\sum_{i=1}^d\lambda_i u_i)=\sum_{i=1}^d\lambda_iN(u_i)$ et que, pour tout $i=1,\ldots , d$, on a $N(u_i)=\sum_{j=0}^{n-1}\sigma^j(u_i)=m\cdot\sum_{j=1}^d u_j$. Ainsi, Žtant en caractŽristique 0, on a~:

$$\ker (N|A)=\left \{\sum_{i=1}^d\lambda_i u_i\mid \sum_{i=1}^d\lambda_i=0\right \}\eqno{(i)}$$
et

$${\rm Im}(N|A)= m\cdot R\cdot\sum_{i=1}^d u_i.\eqno {(ii)}$$

D'autre part,  avec la convention que $u_i=u_j$ si $i\equiv j\pmod d$, on a $\Delta(\sum_{i=1}^d\lambda_i
u_i)=\sum_{i=1}^d
\lambda_i(u_i-u_{i+1})=(\lambda_1-\lambda_d)u_1+(\lambda_2-\lambda_1) u_2+\cdots
+(\lambda_d-\lambda_{d-1}) u_d$, alors on a~:

$$\ker(\Delta |A)=\left\{\sum_{i=1}^d\lambda_i u_i\mid \lambda_1=\cdots =\lambda_d\right\}=R\cdot
\sum_{i=1}^d u_i\eqno{(iii)}$$

Enfin, tout ŽlŽment de ${\rm Im}(\Delta|A)$ est aussi un ŽlŽment de $\ker(N |A)$
($(\lambda_1-\lambda_d)+(\lambda_2-\lambda_1)+\cdots (\lambda_d-\lambda_{d-1})=0$). Et rŽciproquement, si
$\sum_{i=1}^d\mu_i u_i\in \ker (N |A)$, on a vu en $(i)$ que $\sum_{i=1}^d\mu_i=0$, donc en posant
$\lambda_i=\sum_{j=1}^i\mu_{j}$, pour $i=1,\ldots ,d$, on a $\mu_i=\lambda_i-\lambda_{i-1}$, avec la convention
que $\lambda_0=0$ et on a $\sum_{i=1}^d\mu_i u_i=\Delta(\sum_{i=1}^d\lambda_i u_i)\in {\rm
Im}(\Delta|A)$. Ce qui montre que 

$${\rm Im}(\Delta|A)=\ker(N|A)\eqno{(iv}).$$
Par les ŽgalitŽs $(ii)$ et $(iii)$ on trouve que $H^0(A)={\ker(\Delta|A)\over {\rm Im}(N|A)}=R/mR$ et par
l'ŽgalitŽ  $(iv)$ on a  $H^1(A)={\ker (N|A)\over {\rm Im}(\Delta|A)}=\{0\}$. Ainsi $q(A)={1\over
[R:mR]}$.\qed

\bigskip

\centerline{\soustitre Calculs explicites dans le cas d'extensions cycliques}
\bigskip

Fixons pour ce paragraphe $L/K$ une extension cyclique de corps de nombres avec $G={\rm
Gal}(L/K) =\  <\sigma>$, $|G|=[L:K]=n$. On suppose encore que $[L:\Q]=r+2s$ et $[K:\Q]=r'+2s'$, o $r$
(resp.
$r'$) est le nombre de plongements rŽels de $L$ (resp. de $K$) dans $\C$,  et $2s$ et $2s'$ le nombre de
plongements complexes de $L$ (resp. de $K$) dans $\C$.
\bigskip
\defis
\medskip
Soit $\P$ une place infinie de $K$ (il y en a $r'+s'$).  On dit que $\P$ {\it ramifie dans $L$} si
elle est rŽelle et s'il existe une place complexe de $L$ qui prolonge $\P$. Comme nous sommes dans un
contexte galoisien, c'est donc le cas pour toutes les places de $L$ qui prolongent $\P$. Dans ce cas, nous
noterons $e_\P=2$ et $f_\P=1$. Si $\P$ est une place infinie non ramifiŽe, nous noterons $e_\P=f_\P=1$.
Cette dŽfinition implique que dans tous les cas (fini ou infini), nous avons

$$n=e_\P\cdot f_\P\cdot r_\P,$$
o $r_\P$ est le nombre de places qui prolongent $\P$. Le lecteur attentif aura remarquŽ que que si $\P$ est une place finie et l'extension galoisienne, on note $f_\P$ pour $f(\gP/\P)$, o $\gP$ est n'importe quel idŽal premier au-dessus de $\P$, de mme pour $e_\P$.
\bigskip
\newcount\gaagv\gaagv=\gaga

Soit $\m$ un $K$-module. Il est clair que les groupes suivants sont de $G$-modules~:

$$L^*,I_L,O_L,U_{L}, I_L(\widetilde{\m}).$$
(Voir le chapitre 0 pour les dŽfinitions). On supposera de plus que $\m$ est divisible par tous les
idŽaux premiers qui ramifient dans $L$. Puisque l'application $\P\mapsto \P\cdot O_L$ est un
homomorphisme injectif de $I_K(\m)$ dans $I_L(\widetilde{\m})$, par abus, on identifiera $I_K(\m)$
avec son image dans $I_L(\widetilde{\m})$ qui est l'ensemble des idŽaux fractionnaires ${\euf a}$ de $I_L(\widetilde{\m})$ tels que $\sigma({\euf a})={\euf a}$. En effet, si ${\euf a}$ est dans l'image de $I_K(\m)$, alors $\sigma({\euf a})={\euf a}$, car $\sigma|_K$ est l'identitŽ. RŽciproquement, si ${\euf a}=\prod \gP_i^{a_i}$ est tel que $\sigma ({\euf a})={\euf a}$. On peut regrouper les $\gP_i$ qui sont au-dessus d'un mme $\P$ de $\gfP(K)$. On a donc ${\euf a}=\prod_\P\prod_{\gP_i|\P}\gP_i^{a_i}=:\prod_\P{\euf a}_\P$. Puisque $\sigma$ permute les $\gP_i$ au-dessus de $\P$, on a $\sigma({\euf a}_\P)={\euf a}_\P$. De plus, puisque $G$ agit transitivement sur les $\gP_i$ au-dessus de $\P$, il existe $a_\P\in\Z$ tel que ${\euf a}_\P=(\prod_{\gP_i|\P}\gP_i)^{a_\P}=\P^{a_\P}\cdot O_L$, car $\P$ ne ramifie pas dans $O_L$. Cela montre que $\euf a$ est dans l'image de $I_K(\m)$. 

\bigskip

\prop
\medskip

{\sl Sous les mmes notations que prŽcŽdemment,  On a~:

\art{a)}$H^0(I_L(\widetilde{\m}))=I_K(\m)/N_{L/K}(I_L(\widetilde{\m}))$

\art{b)}$H^1(I_L(\widetilde{\m}))=1$ 

\art{c)}$H^0(L^*)=K^*/N_{L/K}(L^*)$

\art{d)}$H^1(L^*)=1$.

}

{\bf Preuve}

Prouvons a). Dire que ${\euf a}\in \ker (\Delta|I_L(\widetilde{\m}))$ revient ˆ dire que $\sigma({\euf a})={\euf a}$
et donc, en vertu de la remarque qui prŽcde la proposition que ${\euf a}\in I_K(\m)$.

Prouvons c). On a $\ker(\Delta |L^*)=\{x\in L^*\mid \Delta(x)=1\}=\{x\in L^*\mid\sigma(x)=x\}=\{x\in
L^*\mid \tau(x)=x$ pour tout $\tau\in G\}={\rm Fix}(G)^*=K^*$.

Prouvons d). Dire que $H^1(L^*)=1$ revient ˆ dire que pour tout $x\in L^*$, si $N_{L/K}(x)=1$, alors il
existe $y\in L^*$ tel que $x={y\over \sigma(y)}$ et c'est le thŽorme 90 de Hilbert  (cf. Chapitre 0).

Prouvons b). Soit ${\euf a}\in \ker(N| I_L(\widetilde{\m}))$. Supposons, comme lors de la remarque prŽcŽdente que ${\euf a}=\prod_{\P\in\gfP_{0}(K)}{\euf a}_\P$ o ${\euf a}_\P$ est le produit de toutes les puissances des
idŽaux premiers de $L$ au-dessus de $\P$. Notons ${\euf a}_\P=\prod_{i=0}^{r-1}\gP_i^{a_i}$ de telle
manire que $a_0\ne 0$ et $\gP_i=\sigma^i(\gP_0)$, pour $i=0,\ldots ,r-1$ (si $a_0=0$ est inŽvitable,
cela veut dire que ${\euf a}_{\P}=O_L$ et on posera pour la suite ${\euf b_p}=O_L$, c'est notamment le cas quand $\P|\m$). Remarquons d'abord que $r>1$. En effet, si $r=1$, $\gP_0$ est le seul idŽal de $L$ au-dessus de $\P$, cela impliquerait que dans la factorisation de $N_{L/K}({\euf a})$, $\P$ appara"trait avec l'exposant $f(\gP_0|\P)\cdot a_0\ne 0$. Cela est impossible puisque ${\euf a}\in \ker(N| I_L(\widetilde{\m}))$ veut dire que $N_{L/K}({\euf a})=O_K$. D'autre part, et pour la mme raison sur l'exposant de $\P$, le fait que $N_{L/K}({\euf a})=O_K$ implique que $N_{L/K}({\euf a}_\P)=O_K$ pour tout $\P$. Calculons~: $O_K=N_{L/K}({\euf a}_\P)=\P^{f_\P\cdot\sum_{i=0}^{r-1}a_i}$, o $f_\P=f(\gP_i|\P)$ pour tout $i$.  Cela prouve que $\sum_{i=0}^{r-1}a_i=0$. Posons alors, pour tout $i=0,\ldots , r-1$, $c_i=\sum_{j=0}^{i}a_i$ et ${\euf b}_\P=\prod_{i=0}^{r-2}\gP_i^{c_i}$. On a
$$\Delta({\euf b}_\P)={\prod_{i=0}^{r-2}\gP_i^{c_i}\over
\prod_{i=1}^{r-1}\gP_i^{c_{i-1}}}=\gP_0^{a_0}\gP_1^{a_1}\cdots\gP_{r-2}^{a_{r-2}}\gP_{r-1}^{-c_{r-2}}={\euf a}_\P,$$
car $-c_{r-2}=a_{r-1}$. On a montrŽ donc que $\Delta(\prod_\P{\euf b}_\P)=\prod_\P{\euf a}_\P={\euf a}$,
ce qui veut dire que ${\euf a}\in \Delta( I_L(\widetilde{\m}))$ et donc que $\ker(N| I_L(\widetilde{\m}))=\Delta( I_L(\widetilde{\m}))$ et $H^1(I_L(\widetilde{\m}))=1$.\qed

\bigskip

\defis
\medskip

On a toujours  $L/K$ une extension cyclique de groupe $G$ et $\m$ un $K$-module (juste ici, ce n'est pas
nŽcessaire qu'il contienne les premiers qui ramifient). On dŽfinit les applications $\iota\ :\
L^*\rightarrow I_L$, $\alpha\mapsto \iota(\alpha)=\alpha \cdot O_L$, $j_\m \ :\ I_L\rightarrow I_L(\widetilde{\m})$ par
$j_\m(\gP)=\gP$ si $\gP\notdiv \m$ et $j_\m(\gP)=O_L$ si $\gP |\m$ et enfin, $f_\m=j_\m\circ\iota$.
Puisque $\m$ est un $K$-module, $\iota,j_\m$ et $f_\m$ sont des homomorphisme de $G$-modules. Posons $S$ l'ensemble des places de $L$ qui divisent $\m$ et ${L^*}^S=\ker(f_\m)=\{\alpha\in L^*\mid \iota(\alpha)$ n'est divisible que par des idŽaux premiers de $S\}$ (on les nomme parfois les $S$-unitŽs).
\bigskip
\newcount\gaagw\gaagw=\gaga

\lem
\medskip
{\sl Sous les mme notations et hypothses, si $q(\ker (j_\m))$ est dŽfini, alors $q(\iota({L^*}^S))$ aussi et
ils sont Žgaux. Si de plus $q(U_L)$ existe ($U_L$ est l'ensemble des unitŽs de $O_L$) alors $q({L^*}^S)$ aussi
et on a

$$q({L^*}^S)=q(U_L)\cdot q(\ker(j_\m)).$$

}

{\bf Preuve}

On a la suite exacte~: $1\rightarrow \iota({L^*}^S)\rightarrow \ker(j_\m)\rightarrow C\rightarrow 1$ o
$$C=\ker(j_\m)/\iota({L^*}^S)\simeq \ker(j_\m)/(\iota(L^*)\cap \ker(j_\m))\buildrel \rm thm.\ d'isom\over
\simeq (\ker(j_\m)\cdot \iota(L^*))/\iota(L^*),$$ le dernier terme est un sous-groupe de $I_L/\iota(L^*)$
qui est fini (cf. [Sam, Thm. 2, chap IV, \S 3, p.71]), ainsi $q(C)=1$, en vertu du Lemme \argbk\ et si $q(\ker(j_\m))$ existe alors
$q(\iota({L^*}^S))$ aussi et ils sont Žgaux, en vertu du Corollaire \argbl. On a aussi la suite exacte
$1\rightarrow U_L\buildrel \rm incl.\over\rightarrow {L^*}^S\buildrel\iota\over\rightarrow \iota({L^*}^S)\rightarrow 1$. Le
Lemme \argbi\ appliquŽ ˆ cette suite exacte nous permet de conclure.\qed

\bigskip\goodbreak

{\soustitre ComplŽments (bien utile et un peu redondant) sur les places infinies}
\smallskip
Soit $L/K$ une extension galoisienne (quelconque) de corps de nombre et $\gP$ une place infinie de
$L$ correspondant  ˆ un plongement $\varphi\ :\ L\rightarrow\C$ (rappelons que $\overline{\varphi}$
est aussi un plongement correspondant ˆ $\gP$ et donc que deux plongements complexes conjuguŽs
correspondent ˆ la mme place). Si
$\sigma\in G$, alors
$\sigma(\gP)$ est la place qui correspond au plongement $\varphi\circ\sigma^{-1}\ :\ L\rightarrow\C$ (ceci
pour avoir
$\sigma_1(\sigma_2(\gP))=(\sigma_1\circ\sigma_2)(\gP)$). De cette manire, $G$ agit transitivement sur
les places infinie de $L$ qui prolongent une mme place infinie de $K$. Soit $|\cdot |_\gP$, la valeur
absolue qui correspond ˆ $\gP$ ($|x|_\gP=|\varphi(x)|$, pour $x\in L$ et $\varphi\ :\ L\rightarrow \C$
est le plongement correspondant ˆ~$\gP$). On a alors
$|x|_{\sigma(\gP)}=|\varphi(\sigma^{-1}(x))|=|\sigma^{-1}(x)|_\gP$, autrement dit,
$|\sigma(x)|_\gP=|x|_{\sigma^{-1}(\gP)}$.

D'autre part, le {\it groupe de dŽcomposition de $\gP$}, notŽ $Z(\gP)$ est le sous-groupe des ŽlŽments de
${\rm Gal}(L/K)$ tels que $\sigma(\gP)=\gP$, donc tels que $\varphi\circ\sigma^{-1}=\varphi$ ou
$\overline\varphi$, ainsi ce sous-groupe est d'ordre 1 ou 2. Si $\P=\gP\cap K$, alors $|Z(\gP)|=2$ si et
seulement si $\P$ ramifie, en effet, tout plongement de $K$ s'Žtend en $n=[L:K]$ plongements de $L$. Si $\P$
est une place complexe, on aura $2n$ plongements au-dessus, donc $n$ places et $\sigma(\gP)=\gP$ si et seulement si
$\sigma=Id_L$. Si $\P$ est une place rŽelle qui ne ramifie pas (donc qui reste rŽelle), il y aura aussi
$n$ places au dessus et on est dans la mme situation. En revanche, si $\P$ ramifie, il n'y aura que
${n\over 2}$ places au-dessus, donc $|Z(\gP)|=2$, car ${\rm Gal}(L/K)$ agit transitivement sur toutes les
places au-dessus de $\P$. Remarquons, qu'ainsi $|Z(\gP)|=e_{\P}\cdot f_{\P}$ pour tout $\gP|\P$, comme pour les places finies.
\smallskip
{\soustitre Fin des complŽments}
\bigskip
\goodbreak
{\th (\bf Minkowski, si $K=\Q$, 1900), (Herbrand, gŽnŽral, 1930)}

{\sl Soit $L/K$ une extension galoisienne (quelconque) de groupe $G$ et $\gP_1,\ldots ,\gP_{r+s}$ les places
infinies de $L$ (les $r$ premires Žtant rŽelles, les $s$ suivantes Žtant complexes). Alors il existe
$\omega_1,\ldots ,\omega_{r+s}\in U_L$ tels que~:

\art{a)}$G$ permute les $\omega_i$ de la mme manire qu'il permute les $\gP_i$ ($\omega_i\leftrightarrow
\gP_i$ est un homomorphisme de $G$-ensemble).

\art{b)}On a $\omega_1\cdots \omega_{r+s}=1$ et c'est la seule relation (sur $\Z$ entre les $\omega_i$),
cela veut dire que si on prend $r+s-1$ $\omega_i$, il sont linŽairement indŽpendants sur $\Z$, ou encore
$\omega_1^{a_1}\cdots \omega_{r+s}^{a_{r+s}}=1\Leftrightarrow a_1=\cdots =a_{r+s}$.

\art{c)} Si $W$ est le sous-$\Z$-module engendrŽ par les $\omega_i$, alors $W$ est un $G$-module
d'indice fini dans $U_L$.

}

{\bf Preuve} 

Il est clair que $c)$ dŽcoule de $b)$ par le thŽorme des unitŽs de Dirichlet (cf. Chapitre 0). 

Montrons a). Soit $\P$ une place infinie de $K$ et choisissons $\gP$ une place de $L$ qui prolonge $\P$.
En regardant la preuve classique du thŽorme des unitŽs de Dirichlet (cf. [Sam, Thm. 1, chap. IV, ¤4, p.72]), on peut trouver  un ŽlŽment $\omega'_\gP$ dans $U_L$ tel que $|\omega'_\gP|_{\euf Q}<1$ pour toute place infinie ${\euf Q}\ne \gP$ de $L$. Posons encore $\omega''_\gP=\prod_{\tau\in Z(\gP)}\tau(\omega'_\gP)$. Si ${\euf Q}$ est une place infinie de $L$ diffŽrente de $\gP$, alors on a
$$|\omega''_\gP|_{\euf Q}=\prod_{\tau\in Z(\gP)}|\tau(\omega'_\gP)|_{\euf Q}=\prod_{\tau\in
Z(\gP)}|\omega'_\gP|_{\tau^{-1}({\euf Q})}<1,$$
car $\tau^{-1}({\euf Q})\ne \gP$ si $\tau\in Z(\gP)$. Soit maintenant $\euf Q$ une place infinie de $L$
au-dessus de $\P$. On choisit $\rho\in G$ tel que $\rho(\gP)={\euf Q}$ et on pose $\omega''_{\euf
Q}=\rho(\omega''_\gP)$; c'est indŽpendant du choix de $\rho$, car $\omega''_\gP$ est invariant par
$Z(\gP)$. Si $\euf Q$ une place infinie de $L$
au-dessus de $\P$ diffŽrente de $\gP$, et si $\euf R$ est une place infinie de $L$ diffŽrente de $\euf
Q$, alors 
$$|\omega''_{\euf Q}|_{\euf R}=|\rho(\omega''_\gP)|_{\euf R}=|\omega_\gP|_{\rho^{-1}({\euf R})}<1. $$

En procŽdant ainsi pour chaque place infinie de $K$, on obtient des $\omega''_1,\ldots ,\omega''_{r+s}\in
U_L$ qui sont permutŽs par $G$ de la mme manire que les $\gP_1,\ldots ,\gP_{r+s}$. On a en outre que
$\rho_1(\rho_2(\omega''_\gP))=(\rho_1\rho_2)(\omega''_\gP)$, faisant de l'application $\gP_i\mapsto
\omega''_{\gP_i}:=\omega''_i$ un isomorphisme de $G$-ensemble. Puisque $|\omega''_i|_{\gP_j}<1$ pour
$i\ne j$, tout choix de $r+s-1$ quelconque des $\omega''_i$ sont toujours $\Z$-linŽairement indŽpendants
(cf. preuve du thŽorme de Dirichlet). Notons $\P_1,\ldots ,\P_{r'+s'}$ les places infinies de $K$ et pour
chaque $i=1,\ldots ,r'+s'$, soit $v_i:=\prod_{\gP_j|\P_i}\omega''_j$. Alors $v_i,\ldots ,v_{r'+s'}\in
U_K$, car ils sont invariants par $G$. Comme $U_K$ est de rang $r'+s'-1$ sur $\Z$, il existe $a_1,\ldots
, a_{r'+s'}\in \Z$ tels que $\prod_{i=1}^{r'+s'} v_i^{a_i}=1$. Ces $a_i$ sont tous non nuls, car
$r'+s'-1$ quelconques des $v_1,\ldots , v_{r'+s'}$ sont linŽairement indŽpendants. On remplace les
$\omega''_i$ par $\omega_i:={\omega''_i}^{a_j}$ o $\gP_i|\P_j$ pour $1\leq i\leq r+s$. On a maintenant
$\prod_{i=1}^{r+s}\omega_i=1$; et par choix des $a_i$, les $\omega_i$ sont permutŽs de le mme manire
que les $\gP_i$ (c'est pour a qu'on a du faire une incursion par les $v_i$). Et c'est la seule
relation car $r+s-1$ parmi les $\omega_i$ sont linŽairement indŽpendants.\qed

\bigskip

\lem
\medskip

{\sl Si $L/K$ est une extension cyclique, alors le quotient de Herbrand de
$U_L$ vaut
$$q(U_L)={[L:K]\over 2^{r_0}},$$
o $r_0$ est le nombre de places infinies de $K$ qui ramifient dans $L$.

}
{\bf Preuve}

Posons $W$ le $G$-module engendrŽ par les $\omega_i$ du thŽorme prŽcŽdent. Pour chaque place infinie
$\P$ de $K$, on forme (abstraitement) le $\Z$-module libre $A_\P=\oplus_{i=1}^{r_\P}\Z u_{\P,i}$ o
$r_\P$ est le nombre de places infinies de $L$ qui prolongent $\P$. Puis, $A=\oplus_{\P\in \gfP_{\infty}(K)}A_\P$. On fait
agir $G$ de telle manire que chaque $A_\P$ est un module de permutation (cf. Proposition \argbm), donc $G$ agit transitivement sur les $u_{\P,i}$, c'est possible, puisque $r_\P$ divise $[L:K]=|G|$. On considre le $G$-homomorphisme $A\rightarrow W$ dŽfini en envoyant les $u_{\P,i}$ sur les $\omega_i$ de manire cohŽrente avec l'action de $G$. C'est un homomorphisme surjectif de noyau  $\Z\cdot
(\sum_{\P\in\gfP_{\infty}(K)}\sum_{i=1}^{r_\P}u_{\P,i})=\Z$ avec l'action triviale de $G$. En bref, on a la suite exacte
$$0\rightarrow \Z\rightarrow A\rightarrow W\rightarrow 1.$$

Or, $\Z$ est un module de permutation (au sens de la Proposition \argbm, avec le $d=1$). Donc, $q(\Z)$ existe et vaut $q(\Z)={1\over [L:K]}$. D'autre part, $A=\oplus A_\P$, et chaque $A_\P$ est un module de
permutation (avec $d=r_\P$), donc, $q(A_\P)$ existe et en vertu du Corollaire \argbk,
$q(A)=\prod_\P q(A_\P)$. En outre, puisque $[L:K]=r_\P e_\P f_\P$, on a $q(A_\P)={r_\P\over
[L:K]}={1\over e_\P f_\P}={1\over |Z(\gP)|}$ pour
$\gP |\P$. Mais, on a vu lors du rappel bien utile que 
$|Z(\gP)|=\cases{1&si $\P$ ne ramifie pas dans $L$\cr 2& sinon\cr}$. Ainsi, 
$$q(A)={1\over 2^{r_0}}$$
o $r_0$ est le nombre de place infinies de $K$ qui ramifient dans $L$. Enfin, puisque $W$ est d'indice
fini dans $U_L$ (partie c) du ThŽorme \argbr), en vertu du Corollaire \argbl\ et du Lemme \argbi, on a

$$ q(U_L)=q(W)={q(A)\over q(\Z)}={[L:K]\over 2^{r_0}}.$$

\qed

\bigskip

\th
\medskip
 
{\sl Soit $L/K$ une extension cyclique de corps de nombres. Soit $\m=\m_0\cdot \m_\infty$ un $K$-module
tel que $\m_\infty$ contienne toutes les places infinies de $K$ qui ramifient dans $L$. Soit $S$
l'ensemble de places de $L$ qui divisent $\m$. Alors $q({L^*}^S)$ existe et vaut
$$q({L^*}^S)=[L:K]\cdot \prod_{\P|\m}{1\over e_\P f_\P}.$$

}
{\bf Preuve}

On a montrŽ au Lemme \argbq, que si $q(U_L)$ et $q(\ker (j_\m))$ existaient, alors $q({L^*}^S)$ aussi et que $q({L^*}^S)=q(U_L)\cdot q(\ker(j_\m))$. On a montrŽ au Lemme \argbs\ que $q(U_L)$ existait et valait
${[L:K]\over 2^{r_0}}$. Il suffit donc de calculer $q(\ker (j_\m))$; et c'est le plus facile ˆ voir~:

Par dŽfinition, $\ker(j_\m)$ est le groupe abŽlien libre engendrŽ par les idŽaux premiers de $L$ qui sont
dans $S$. Notons $A(\P)$ le groupe abŽlien libre engendrŽ par les idŽaux premiers de $L$ qui sont
au-dessus de $\P$. Alors $A(\P)$ est un sous-$G$-module de $\ker (j_\m)$ et $\ker (j_\m)=\oplus_{\P|\m_{0}} A(\P)$. Chaque $A(\P)$ est un module de permutation (au sens de la Proposition \argbm, avec $d=r_\P$) et $G$ agit transitivement sur la base formŽe des idŽaux premiers de $L$ au-dessus de $\P$. A nouveau, puisque $[L:K]=r_\P e_\P f_\P$, la Proposition \argbm\ nous montre  $q(A(\P))$ existe et vaut ${1\over e_\P f_\P}$. Donc, en vertu du Corollaire \argbj, on a $q(\ker (j_\m))$ existe et vaut 
$$q(\ker (j_\m))=\prod_{\P|\m_0} q(A(\P))=\prod_{\P|\m_0}{1\over e_\P f_\P}.$$
Enfin, en se souvenant que $2^{r_0}=\prod_{\P|\m_\infty} e_\P f_\P$, on trouve
$$q({L^*}^S)=q(U_L)\cdot q(\ker(j_\m))=[L:K]\cdot\prod_{\P|\m_\infty} {1\over  e_\P f_\P}\cdot
\prod_{\P|\m_0}{1\over e_\P f_\P}=[L:K]\cdot \prod_{\P|\m}{1\over e_\P f_\P}.$$

\qed

\vfill\eject

\global\advance\chapnomb by 1
\nomb=1

\centerline{\para Chapitre 5 :}
\medskip
\centerline{\para Un calcul d'indice }
\bigskip

Dans ce chapitre, nous allons calculer (comme son nom l'indique) un indice. Cet indice para"t sorti de nulle part, mais il sera crucial pour prouver l'ŽgalitŽ fondamentale du corps de classe au chapitre suivant; et cette dernire sera une des briques importantes pour dŽmontrer la rŽciprocitŽ d'Artin. Ici, le lecteur ferait bien de se souvenir des dŽfinitions faites au Chapitre 0 sur les $K$-modules aux pages  \the\gaagg\ et suivantes.
\bigskip
\defi
\medskip

Soit $L/K$ une extension cyclique de corps de nombres de groupe $G=<\sigma>$. Posons $N=N_{L/K}$ la norme de $L$ sur $K$. Soit $\m$ un $K$-module. On pose 
$$a(\m)=[K^*:N(L^*)K_\m^*]$$
o, on le rappelle, $K^*_\m=\{x\in K^*\mid x\equiv 1\pmodast \m \}$. 
\bigskip
\headline={\hfill \phantom{ouuh}\hfill}
\lem
\medskip
{\sl Soit $L/K$ une extension quelconque de corps de nombres, $\m=\m_{0}\cdot\m_{\infty}$ un $K$-module et $\widetilde{\m}=\widetilde{\m}_{0}\cdot\widetilde{\m}_{\infty}$, le $L$-module engendrŽ par $\m$. Alors on a

\art{a)}$L^*_{\tilde{\m}_{0}}\cap K=K^*_{\m_{0}}$.

\art{b)}$$N(L^*_{\tilde \m})\subset K^*_\m.$$

}

{\bf Preuve}

Prouvons a). L'inclusion $\supset$ est Žvidente. Prouvons l'autre inclusion. Soit $\alpha\in L^*_{\tilde{\m}_{0}}\cap K$ et soit $\P\in\gfP(K)$ tel que $\P |\m_{0}$. Posons $n=n_{\P}$ l'exposant de $\P$ dans $\m$. Soit $\gP\in\gfP(L)$ tel que $\gP |\P$. Alors l'exposant de $\gP$ dans $\widetilde{m}_{0}$ vaut $n\cdot e$, o $e=e(\gP/\P)$ est l'indice de ramification de $\gP/\P$. Dire que $\alpha\in L^*_{\tilde{\m}_{0}}\cap K$ implique que $\alpha=1+x$, avec $x\in K$ et $v_{\gP}(x)\geq n\cdot e$. Or, il est Žvident que $v_{\P}(x)=e\cdot v_{\gP}(x)$. Donc $v_{\P}(x)\geq n$, ce qui veut dire que $\alpha=1+x\in K^*_{\P^n}$, ceci pour tout $\P |\m_{0}$. Donc $\alpha\in  \cap_{\P|\m_{0}} K^*_{\P^{n_\P}}=K^*_{\m_{0}}$. Remarquons que si on remplace $\m_{0}$ par $\m$ avec d'Žventuelles places ˆ l'infini, cette ŽgalitŽ est fausse si une des places infinie divisant $\m$ devient complexe partout par exemple.

Prouvons b). Soit $x\in L^*_{\tilde \m}$. Si $\P|\m$ est une place infinie (donc rŽelle) correspondant ˆ un plongement $\sigma\, :\, K\to \R$. Soient $\sigma_1,\ldots ,\sigma_r, \sigma_{r+1},\overline{ \sigma}_{r+1}, \ldots ,\sigma_{r+s},\overline{ \sigma}_{r+s}$ les extensions de $\sigma$ en plongements de $L$ dans $\C$ tels que $\sigma_1,\ldots ,\sigma_r$ soient rŽelles et les autres complexes. Alors $\sigma_1,\ldots ,\sigma_r$ correspondent aux places infinies $\gP_{1},\ldots ,\gP_{r}$ qui divisent $\widetilde{\m}$, donc $\sigma_{i}(x)>0$ pour $i=1,\ldots, r$, et donc
$$\sigma(N(x))=\prod_{i=1}^r\sigma_{i}(x)\cdot \prod_{j=1}^s \sigma_{r+j}\cdot\overline{ \sigma}_{r+j}>0,$$
ce qui montre que $N(L^*_{\tilde \m})\subset K^*_{\m_{\infty}}$
Regardons maintenant le cas des places finies. Soit $E/L$ l'enveloppe galoisienne de $L/K$ (c'est-ˆ-dire la plus petite extension galoisienne $E/K$ qui contienne $L$). Si $G={\rm Gal}(E/K)$ et $H={\rm Gal}(E/L)$, on a $N(x)=\prod_{\sigma}\sigma(x)$, o $\sigma$ parcourt un systme de reprŽsentants de classe de $G$ modulo $H$ (ˆ gauche). Soit $\widetilde{\widetilde{\m}}_{0}$ le $E$-module engendrŽ par $\m_{0}$. Puisque $x\in L^*_{\tilde{\m}_{0}}\subset E^*_{\tilde{\tilde\m}_{0}}$ et que $\widetilde{\widetilde{\m}}_{0}$ est invariant par $G$, on a aussi $\sigma(x)\in E^*_{\tilde{\tilde\m}_{0}}$ pour tout  $\sigma\in G$. Ainsi, $N(x)=\prod_{\sigma}\sigma(x)\in E^*_{\tilde{\tilde\m}_{0}}\cap K\buildrel\rm a)\over =K^*_{\m_{0}}$. Et cela montre que  $N(L^*_{\tilde \m})\subset K^*_{\m_{\infty}}\cap  K^*_{\m_{0}}= K^*_{\m}$.\qed

\bigskip\goodbreak
 
\newcount\gaagx\gaagx=\gaga

\lem
\medskip

{\sl Si $\m$ et $\n$ sont des $K$-modules premiers entre eux, alors on a~:

$$a(\m\cdot \n)=a(\m)\cdot a(\n).$$

}

{\bf Preuve}

On se souvient (Corollaire \argd) de l'isomorphisme ~: $K^*/K^*_{\m\n}\to K^*/K^*_\m\times K^*/K^*_\n$. Le passage au quotient induit un homomorphisme surjectif $$f\, :\, K^*/K_{\m\n}\to K^*/(N(L^*)K^*_\m)\times K^*/(N(L)K^*_\n).$$
 Pour prouver le lemme, il suffit de montrer que $\ker(f)=N(L^*)K^*_{\m\n}/K^*_{\m\n}$. Puisque
$K_\m^*\cap K_\n^*=K^*_{\m\n}$, il est Žvident que $N(L^*)K^*_{\m\n}/K^*_{\m\n}\subset \ker(f)$.
RŽciproquement, soit $\alpha\cdot K^*_{\m\n}\in\ker(f)$. On a donc $\alpha\cdot K^*_\m\subset N(L^*)K_\m^*$
et  $\alpha\cdot K^*_\n\subset N(L^*)K_\n^*$. Cela veut dire qu'il existe $\beta_1$, $\beta_2\in L^*$ tels
que $\alpha\equiv N(\beta_1)\pmodast\m$ et $\alpha\equiv N(\beta_2)\pmodast\n$. Puisque $\m$ et $\n$ sont premiers entre eux, $\widetilde\m$ et $\widetilde\n$ le sont aussi. En vertu du ThŽorme d'approximation dŽbile (cf. ThŽorme \argc), il existe $\beta\in L^*$ tel que $\beta\equiv\beta_1\pmodast{\widetilde\m}$ et $\beta\equiv\beta_2\pmodast{\widetilde\n}$,
c'est-ˆ-dire $\beta\cdot\beta_1^{-1}\in L^*_{\tilde \m}$ et $\beta\cdot\beta_2^{-1}\in L^*_{\tilde \n}$. Puisque  $N(L^*_{\tilde \m})\subset K^*_\m$  (et idem pour $\n$) (cf. lemme prŽcŽdent),
on a alors $N(\beta)\equiv N(\beta_1)\pmodast\m$ et  $N(\beta)\equiv N(\beta_2)\pmodast\n$, ainsi, $\alpha\equiv
N(\beta)\pmodast\m$ et $\pmodast\n$, donc (ˆ nouveau gr‰ce au thŽorme d'approximation dŽbile) $\alpha\equiv
N(\beta)\pmodast{\m\cdot\n}$; ce qui montre que $\alpha\in N(L^*)\cdot K_{\m\n}^*$ et le thŽorme est prouvŽ. \qed

\bigskip

En vertu du lemme qu'on vient de voir, le calcul de $a(\m)$ se rŽduit au cas o $\m=\P^m$ avec $m\geq 1$ entier pour les places finies et $\m=\P$ pour les places infinies.

\bigskip

\lem
\medskip

{\sl Si $\P$ est une place infinie rŽelle et $\m=\P$, alors 
$$a(\m)=e_\P=e_\P f_\P.$$
(voir DŽfinitions \argbn).
}

\headline={\hfill \smcap Un calcul d'indice \hfill}
{\bf preuve}

Si $\varphi$ est le plongement correspondant ˆ $\P$, alors $K^*_\m$ est le noyau de l'application
surjective
$$K^*\buildrel \varphi\over\rightarrow \R \buildrel{\rm nat.}\over\rightarrow \R^*/\R^*_+=\{\pm 1\}.$$
Ainsi, $[K^*:K^*_\m]=2$. 

Supposons que $\P$ ramifie dans $L$. Soit $\gP_1,\ldots ,\gP_r$ les places de $L$ qui
prolongent $\P$. Puisque $\P$ ramifie, elles sont toutes complexes et $r={n\over 2}$, o $n=[L:K]$. Soit $\sigma_1,\ldots
,\sigma_r$ les plongements correspondants aux $\gP_i$. Alors, pour tout $x\in L^*$, on a
$\varphi(N(x))=\prod_{i=1}^r\sigma_i(x)\overline{\sigma_i(x)}>0$. Donc, $N(L^*)\subset K^*_\m$, et donc
$N(L^*)K^*_\m=K^*_\m$, donc $a(\m)=2=e_\P$.

Supposons que $\P$ ne ramifie pas dans $L$. A nouveau $\gP_1,\ldots ,\gP_r$ sont les places de $L$ qui
prolongent $\P$ et $\sigma_1,\ldots ,\sigma_r$ les plongements correspondants; dans ce cas, $r=n$ et 
$\varphi(N(x))=\prod_{i=1}^r\sigma_i(x)$. Par le thŽorme d'approximation dŽbile, il est possible de trouver $x\in L^*$ tel que $\sigma_1(x)<0$ et $\sigma_i(x)>0$, pour $i=2,\ldots ,r$. Cela prouve que
$N(L^*)\not\subset K^*_\m$, ce qui montre que $N(L^*)K^*_\m=K^*$ et donc $a(\m)=1=e_\P$.\qed

\bigskip

Passons aux places finies~:

Tout d'abord un petit lemme gentil~:
\medskip
\lem

{\sl Soit $A$ un anneau de Dedekind ne possŽdant qu'un nombre fini d'idŽaux premiers. Alors $A$ est principal

}

{\bf Preuve}

Soit $\aa$ un idŽal de $A$. Puisque $A$ est de Dedekind, il existe des idŽaux premiers $\P_1,\ldots ,\P_s$ et des $r_1,\ldots ,r_s$ uniques tels que $\aa=\P_1^{r_1}\cdots\P_s^{r_s}$. Pour $i=1,\ldots , s$, fixons $a_i\in \P_i^{r_i}\setminus\P_i^{r_i+1}$; c'est possible dans un anneau de Dedekind (si $\P_i^{r_i}=\P_i^{r_i+1}$, alors $O_K=\P_i$...). Par le thŽorme chinois, il existe $a\in A$ tel que $a\equiv a_i\pmod{\P_i^{r_i+1}}$, pour tout $i=1,\ldots , s$. Alors on a $aA=\P_1^{r_1}\cdots\P_s^{r_s}=\aa$, car pour tout idŽal premier de $A$, on a $v_\P(\aa)=v_\P(a)=v_\P(aA)$.\qed

\bigskip\goodbreak

\lem
\medskip
{\sl Soit $\P$ une place finie et $\m=\P^m$ ($m\geq 1$). Alors 

\art{a)}$[K^*:N(L^*)\cdot K^*(\m)]=f_\P$

\art{b)}$a(\m)=f_\P\cdot [K^*(\m):(K^*(\m)\cap N(L^*))\cdot K^*_\m]$.

}

{\bf Preuve}

Montrons que b) dŽcoule de a). On a
$a(\m)=[K^*:N(L^*)K^*_\m]=[K^*:N(L^*)K^*(\m)][N(L^*)K^*(\m):N(L^*)K^*_\m]\buildrel {\rm a)}\over
=f_\P\cdot [N(L^*)K^*(\m):N(L^*)K^*_\m]$. L'application (qui est l'inclusion) $K^*(\m)\rightarrow
N(L^*)K^*(\m)$ donne un isomorphisme
$$K^*(\m)/K^*(\m)\cap N(L^*)K^*_\m\longrightarrow N(L^*) K^*(\m)/N(L^*)K^*_\m.$$
Et, finalement, puisque $K^*_\m\subset K^*(\m)$, on a 
$$K^*(\m)\cap N(L^*)K^*_\m=K^*(\m)K^*_\m\cap N(L^*)K^*_\m=(K^*(\m)\cap N(L^*))\cdot K^*_\m$$
ce qui montre la partie b).

Montrons la partie a). Pour simplifier l'Žcriture, notons $R=O_K$ et $R_{(\P)}$ le localisŽ de $R$ en
$\P$. Il est bien connu que $R_{(\P)}$ est un anneau de valuation discrte (c'est une des propriŽtŽ fondamentale des anneaux de Dedekind), donc local et
principal. Notons $\pi$, {\it l'uniformisante de $R_{(\P)}$}, c'est-ˆ-dire l'ŽlŽment qui engendre
$\P R_{(\P)}$, l'unique idŽal maximal de $R_{(\P)}$. On remarque que dans notre cas,
$K^*(\m)=R^*_{(\P)}$. Ainsi, tout ŽlŽment $x$ de $K^*$ s'Žcrit de manire unique $x=u\cdot \pi^k$, avec
$u\in K^*(\m)$ et $k\in\Z$. Ce qui montre que $K^*$ est en bijection avec $\Z\times K^*(\m)$. Notons
maintenant $T=O_L$ et $T_{(\P)}=(R-\P)^{-1}\cdot T$. L'ensemble des idŽaux premiers de $T_{(\P)}$ est
en bijection avec les idŽaux de $T$ qui ne rencontrent pas $(R-\P)$. Ainsi, si $\gP_1,\ldots , \gP_r$
sont les idŽaux de $T$ au-dessus de $\P$, alors $\gP_1T_{(\P)},\ldots, \gP_rT_{(\P)}$ est l'ensemble des
idŽaux premiers de $T_{(\P)}$. D'autre part, puisque $T$ est un
anneau de Dedekind, alors $T_{(\P)}$ l'est aussi. Donc, il est principal (cf.  Lemme \argbx). Posons $\gP_iT_{(\P)}=(\pi_i)$, pour $i=1,\ldots ,r$. Ainsi, chaque
ŽlŽment de $L^*$ s'Žcrit $v\cdot\pi_1^{k_1}\cdots\pi_r^{k_r}$, avec $v\in T^*_{(\P)}$ et $k_1,\ldots ,
k_r\in \Z$. Si $f=f(\gP_i/\P)=f_\P$, alors on a $N(\gP_i)=\P^f$, $i=1,\ldots ,r$, et donc $N(\pi_i)=u_i\cdot
\pi^f$, avec $u_i\in R^*_{(\P)}=K^*(\m)$. On en dŽduit que 
$$N(L^*)\cdot K^*(\m)=\{\pi^{fk}\mid k\in\Z\}\cdot K^*(\m)\simeq f\Z\times K^*(\m).$$
D'o $[K^*:N(L^*)\cdot K^*(\m)]=[\Z\times K^*(\m):f\Z\times K^*(\m)]=f=f_\P$.\qed

\bigskip

\defi
\medskip

Mettons-nous sous les mmes hypothses que prŽcŽdemment, c'est-ˆ-dire $L/K$ est une extension cyclique de corps de nombres de groupe de Galois $G$, de cardinal $n$. On prend $\P$ une place finie et on considre le $K$-module $\m=\P^m$ avec $m\in\N$. Soit $\gP$ une place de $L$ au-dessus de $\P$. Nous notons $\bbK_\P$ et $\bbL_\gP$ les corps locaux associŽs aux places $\P$ et $\gP$ respectivement. L'extension $\bbL_\gP/\bbK_\P$ est aussi cyclique de groupe de Galois
canoniquement isomorphe ˆ $Z(\gP)=\{\sigma\in G\mid \sigma(\gP)=\gP\}$ de cardinal $e_\P f_\P:=n_\P$ [cf.
Fr-Tay Th. 21, p.118]. Notons encore $O_\P=O_{\bbK_{\P}}$ l'anneau de valuation de $\bbK_\P$ et $O_\gP$ celui de
$\bbL_\gP$. Soit
$\widehat{\P}$ l'unique idŽal maximal de $O_\P$ et $\widehat{\gP}$ celui de $O_\gP$. On note
$O_{(\P)}:=O_\P\cap K$ le localisŽ de $O_K$ en $\P$ (avant on l'avait notŽ $R$, mais c'Žtait parce qu'on
avait besoin de $T$...) et
$O_{(\gP)}:=O_\gP\cap L$ le localisŽ de $O_L$ en $\gP$. L'idŽal maximal de $O_{(\P)}$ se note
$\widetilde{\P}$ et celui de $O_{(\gP)}$ se note $\widetilde{\gP}$. Les unitŽs de
$O_{\P}$ devraient se noter $O^*_\P$, mais se notent $U_\P$ (attention de ne pas confondre avec le $U_\m$ de la DŽfinition \argh\ et du ThŽorme \argl) et celles de $O_\gP$ se notent $U_\gP$, les
unitŽs de $O_{(\P)}$ devraient se noter  $U_{(\P)}$ ou $O_{(\P)}^*$, mais dans notre cas, c'est $K^*(\m)$, o $\m=\P\cdot\m_\infty$. \newcount\gaagy \gaagy=\gaga
Enfin, pour $k\in\N$, on Žcrit $U_\P^{(k)}$ pour $1+\widehat{\P}^k\subset 1+\widehat{\P}\subset U_\P$, car $\widehat{\P}$ est le seul idŽal maximal de $O_\P$. On a aussi $\P O_{(\P)}=\widetilde{\P}$ et $\P O_\P=\widetilde{\P}O_\P=\widehat{\P}$ et $O_K/\P^k\simeq O_{(\P)}/\widetilde{\P}^k\simeq O_\P/\widehat{\P}^k$ [cf. Fr-Tay Th. 11 + Cor, p.77]. La norme $N_{\bbL_\gP/\bbK_\P}$ sera notŽe $N_\P$ ou mme $N$ s'il n'y a pas d'ambigu•tŽ.

Nous avons aussi besoin d'Žtendre la dŽfinition de $\pmodast$ sur $\bbK_{\P}^*$ : si $x,y\in\bbK_{\P}^*$ et $n>0$ est un entier, alors on dit que $x\equiv y\pmodast {\widehat{\P}^{n}}$ si $\dst {x-y\over y}\in \widehat{\P}^n$. Sur $\bbL_{\gP}$ on dŽfinit cette Žquivalence de la mme manire.
\bigskip

Avec tout ce petit monde, nous sommes prt ˆ Žnoncer le lemme suivant
\bigskip
\lem
\medskip

{\sl Soit $L/K$ une extension cyclique de corps de nombres, $\P\in\gfP_{0}(K)$, $m\in \N$, $m>0$ et $\m=\P^m$ un $K$-module. Alors on a
$$K^*(\m)/(K^*(\m)\cap N(L^*))K^*_\m\simeq U_\P/N(U_\gP) U_\P^{(m)}$$
o $\gP$ est un idŽal premier de $L$ au-dessus de $\P$.

}

{\bf Preuve}

Rappelons le fait ŽlŽmentaire suivant~: tout homomorphisme d'anneau $f\, :\, A\to B$ dŽfinit un
homomorphisme de groupe $f^*\ :\ A^*\to B^*$ tel que $\ker(f^*)=(1+\ker(f))\cap A^*$. L'homomorphisme
surjectif $O_{(\P)}\to O_{(\P)}/\widetilde{\P}^m$ induit donc un homomorphisme $O^*_{(\P)}=K^*(\m)\to
(O_{(\P)}/\widetilde{\P}^m)^*$. Puisque $O_{(\P)}$ est local, alors $1+\widetilde{\P}^m\subset K^*(\m)$,
donc cet homomorphisme est surjectif. On a aussi $(O_{(\P)}/\widetilde{\P}^m)^*\simeq
(O_\P/\widehat{\P}^m)^*=U_\P/(1+\widehat{\P}^m)=U_\P/U_\P^{(m)}$ (l'avant dernire ŽgalitŽ vient aussi du fait
que $O_\P$ est local). En particulier l'homomorphisme 
$$\eqalign{f\ :\ K^*(\m)&\longrightarrow U_\P/N(U_\gP) U_\P^{(m)}\cr x&\longmapsto x\;N(U_\gP)U_\P^{(m)}\cr}$$
est surjectif. Il nous reste ˆ faire la preuve que $\ker(f)=(K^*(\m)\cap N(L^*))K^*_\m$. 

En regardant les
dŽfinitions, on observe que $K^*_\m=1+\widetilde{\P}^m\subset 1+\widehat{\P}^m=U_\P^{(m)}$. Donc,
$K^*_\m\subset \ker(f)$. Soit $\alpha\in L^*$ est tel que $N_{L/K}(\alpha)\in K^*(\m)$, alors
$v_\P(N_{L/K}(\alpha))=0$, (voir le paragraphe suivant ou le Chapitre 0 pour se remŽmorer la dŽfinition de $v_\P$); cela implique que
$v_{\gP_i}(\alpha)=0$ pour tout idŽal $\gP_i$ de $L$ au-dessus de $\P$, donc que $\alpha\in L^*(\widetilde{\m})$ (se souvenir de la dŽfinition de $\widetilde{\m}$). Soit $\tau_1,\ldots ,\tau_r$ un systme de reprŽsentant de $G$ modulo $Z(\gP)$. Alors 
$$N_{L/K}(\alpha)=\prod_{\tau\in Z(\gP)}\prod_{i=1}^r\tau_i\tau(\alpha)=\prod_{\tau\in
Z(\gP)}\tau\left(\prod_{i=1}^r\tau_i(\alpha)\right)=N_\P\left(\prod_{i=1}^r\tau_i(\alpha)\right),$$
et $\prod_{i=1}^r\tau_i(\alpha)\in L^*(\widetilde{\m})\subset U_\gP$. On a donc montrŽ que $K^*(\m)\cap N(L^*)\subset
N_\P(U_\gP)\cap K^*\subset \ker(f)$. Et ainsi $(K^*(\m)\cap N(L^*))K^*_\m\subset \ker(f)$.

Montrons l'autre inclusion~: soit donc $\alpha\in \ker(f)$. Il existe donc $\beta\in U_\gP$ tel que
$\alpha N_\P(\beta)^{-1}\in U_\P^{(m)}$. Soient $\gP_1=\gP,\gP_2,\ldots ,\gP_r$ les idŽaux premiers de
$L$ au-dessus de $\P$. Comme $L^*$ est dense dans $\bbL_\gP$, il existe $\beta_0\in L^*$ tel que
$\beta_0\equiv\beta\pmodast{\widehat{\gP}^{me}}$, o $e=e_\P$ et $x\equiv y\pmodast{\widehat{\gP}^{me}}$
veut dire par extension que ${x\over y}\in 1+{\widehat{\gP}^{me}}$. Le thŽorme chinois pour $O_L$ nous assure
l'existence d'un $\gamma\in L^*$ tel que 
$$\gamma\equiv\beta_0\pmodast {\widehat{\gP}_1^{me}}\quad \hbox{et}\quad \gamma\equiv 1\pmodast {\widehat{\gP}_j^{me}}\ \hbox{lorsque $j>1$}.$$
Prenons, comme tout ˆ l'heure $\tau_1,\ldots ,\tau_r$ un systme de reprŽsentants de $G$ modulo $Z(\gP)$,
mais en plus, on impose que $\tau_1=Id_L$ et $\tau_j(\gP)=\gP_j$, pour tout $j>1$. Alors, si $j>1$ et
$\tau\in Z(\gP)$, on a $\tau_j^{-1}\tau(\gamma)\equiv 1\pmodast{\widehat{\gP}^{me}}$. Alors,
$$N_{L/K}(\gamma)=\prod_{j=1}^r\prod_{\tau\in Z(\gP)}\tau_j^{-1}\tau(\gamma)\equiv\prod_{\tau\in
Z(\gP)}\tau(\gamma)\equiv\prod_{\tau\in Z(\gP)}\tau(\beta)=N_\P(\beta)\pmodast{\widehat{\gP}^{me}}.$$
Puisque $N_{L/K}(\gamma)$ et $N_\P(\beta)$ sont dans  $\bbK_\P$ et que $\widehat{\gP}^e$ est le seul idŽal au-dessus de $\widehat{\P}$, cette dernire congruence est vraie modulo $\widehat{\P}^m$. On a donc prouvŽ que $\alpha N_\P(\beta)^{-1}\equiv \alpha N_{L/K}(\gamma)^{-1}\pmodast{\widehat{\P}^m}$ (avec la mme convention pour $\pmodast {\widehat{\P}^m}$). Et, puisque par hypothse $\alpha N_\P(\beta)^{-1}\in U_\P^{(m)}$, on a $\alpha N_{L/K}(\gamma)^{-1}  \in U_\P^{(m)}\cap K^*=K^*_\m$. Donc
$N_{L/K}(\gamma)\in \ker(f) K^*_\m\cap N(L^*)\subset K^*(\m)\cap N(L^*)$ et donc, $\alpha\in (K^*(\m)\cap
N(L^*))K^*_\m$.\qed

\bigskip
\centerline{\soustitre DIGRESSION}
\medskip

Interrompons un court instant notre propos pour un petit rŽsultat technique sur les sŽries logarithmes et
exponentielles sur $\bbK_\P$. Et rappelons que l'on dŽfinit formellement 
$$\exp(x)=\sum_{n=0}^\infty {x^n\over n!}\quad\hbox{et}\quad \log(1+x)=\sum_{n=1}^\infty
(-1)^{n+1}{x^n\over n}$$
\newcount\gaagz \gaagz=\gaga
De plus, si $p\Z$ est l'unique idŽal de $\Z$ au-dessous de $\P$, on note $e_0=v_\P(p)$.

Voici quelques rŽsultats ŽlŽmentaires sur $\bbK_\P$ (qu'on a dŽjˆ d'ailleurs vus au chapitre 0, pour la plupart)~:

Puisque $O_\P$ est un anneau de valuation discrte, tout ŽlŽment $x$ de $\bbK_\P$ s'Žcrit de manire unique $x=u\cdot \pi^t$ o $\pi$ est un gŽnŽrateur de $\widehat{\P}$, appelŽ {\it uniformisante} et
$t=:v_\P(x)\in\Z\cup\{\infty\}$ est la {\it valuation $\P$-adique de $x$} (qui Žtend celle sur $K$). On a les propriŽtŽs suivantes~:
\art{a)}$v_\P(x)=\infty\iff x=0$

\art{b)} $v_\P(xy)=v_\P(x)+v_\P(y)$

\art{c)} $v_\P(x+y)\geq\inf (v_\P(x),v_\P(y))\ \hbox{avec ŽgalitŽ si $v_\P(x)\ne v_\P(y).$}$

Sur $\bbK_\P$, on
dŽfinit une {\it valeur absolue} ou {\it norme} qui vaut 
$$|x|_\P=\N(\P)^{-v_\P(x)}.$$
Cette valeur absolue est {\it non-archimŽdienne}. De a), b) c), on trouve $|x|_\P=0\iff x=0$,
$|xy|_\P=|x|_\P\cdot |y|_\P$ et $|x+y|\leq \sup (|x|_\P,|y|_\P)$ avec ŽgalitŽ si $|x|_\P\ne |y|_\P$.
Cette dernire propriŽtŽ implique qu'une somme $\sum_{n=0}^\infty x_n$ converge dans $\bbK_\P$ pour cette valeur absolue si et seulement si $|x_n|_\P$ tend vers $0$.

\bigskip

\prop
\medskip

{\sl Dans $\bbK_\P$, les sŽries $\exp(x)$ et $\log(1+x)$ convergent si $v_\P(x)>{e_0\over p-1}$ (on peut voir en particulier que $\exp(1)=e$ n'existe pas dans $\bbK_\P$). De plus, si 
$v_\P(x)>{e_0\over p-1}$, on a 

\art{a)} $v_\P(\exp(x)-1)=v_\P(x)$

\art{b)}$v_\P(\log(1+x))=v_\P(x)$.

}

{\bf Preuve}

Soit $n\in\N$, $n\geq 1$. On vŽrifie facilement que $v_p(n!)=[{n\over p}]+[{n\over p^2}]+[{n\over
p^3}]+\cdots$. Ecrivons $n$ en base $p$~: $n=a_0+a_1p+\cdots +a_k p^k$, avec $a_k\ne 0$ et $0\leq
a_i\leq p-1$. Donc ${n\over p^l}=\underbrace{a_0+\cdots +a_{l-1}p^{l-1}\over p^l}_{<1}+a_l+\cdots + a_k
p^{k-l}$. Ainsi, $[{n\over p^l}]=a_l+\cdots +a_k p^{k-l}$ et donc
$$\eqalign{v_p(n!)&=a_1+a_2(1+p)+a_3(1+p+p^2)+\cdots +a_k(1+p+\cdots +p^{k-1})\cr
&=a_1{p-1\over p-1}+a_2{p^2-1\over p-1}+a_3{p^3-1\over p-1}+\cdots +a_k{p^k-1\over p-1}\cr
&={1\over p-1}(n-a_0-a_1-a_2-\cdots -a_k)={1\over p-1}(n-S),\cr}$$ 
o $S=\sum_{i=0}^k a_i$. On vŽrifie tout aussi facilement que $v_\P(n !)=v_p(n!)\cdot
v_\P(p)=v_p(n!)\cdot e_0$. Calculons donc 
$$ v_\P({x^n\over n!})=n\cdot v_\P(x)-e_0\cdot v_p(n!)=n\cdot\underbrace{(v_\P(x)-{e_0\over p-1})}_{>0\
\rm par\ hyp.}+{e_0\over p-1}\cdot S\longrightarrow\infty\quad \hbox{si }n\to\infty.$$ 
Cela prouve que $|{x^n\over n!}|_\P$ tend vers 0, donc, que la sŽrie $\exp(x)$ converge. Montrons la
partie a)~: puisque $\exp(x)-1=x+{x^2\over 2}+{x^3\over 6}+\cdots$, il suffit de montrer que si $n\geq
2$, on a $v_\P({x^n\over n!})>v_\P(x)$. Cela est vrai~:
$$v_\P({x^n\over n!})-v_\P(x)=\underbrace{(n-1)}_{>0}\underbrace{(v_\P(x)-{e_0\over p-1})}_{>0}+{e_0\over
p-1}\underbrace{(-1+\sum_{i=0}^k a_i)}_{\geq 0}>0,$$
toujours avec $n=a_0+a_1p+\cdots +a_k p^k$, le dŽveloppement de $n$ en base $p$. Donc a) est prouvŽ.

Regardons maintenant le sŽrie $\log(1+x)$. Puisqu'on a si bien ŽvaluŽ $v_\P({x^n\over n!})$,
profitons-en ! on voit que $v_\P({x^n\over n})=v_\P({x^n\over n!})+\underbrace{v_\P((n-1)!)}_{\geq
0}\to\infty$ si $n\to \infty$. Donc, $\log(1+x)$ converge si $v_\P(x)>{e_0\over p-1}$ (en fait on vient
d'observer que le domaine de convergence de $\log(1+x)$ est plus Žtendu que celui de $\exp(x)$, ce qui est
contraire ˆ la situation dans $\C$). 

Pour la partie b), il suffit d'observer comme pour la partie a) que 

$$v_\P({x^n\over n})-v_\P(x)=\underbrace{(n-1)}_{>0}\underbrace{(v_\P(x)-{e_0\over
p-1})}_{>0}+\underbrace{v_\P((n-1)!)}_{\geq 0}+{e_0\over p-1}\underbrace{(-1+\sum_{i=0}^k a_i)}_{\geq
0}>0.$$
\phantom{les maths sont belles mais le monde aussi, et au fond qui guide qui ?}\qed
\bigskip

\coro

\medskip

{\sl Si $m>{v_\P(p)\over p-1}$, alors $\exp(x)$ est un isomorphisme de groupe $\widehat{\P}^m\to
U_\P^{(m)}$. La rŽciproque Žtant $\log(x)$. 

}
{\bf Preuve}

Cela dŽcoule de la proposition prŽcŽdente et des identitŽs formelle $\exp(x+y)=\exp(x)\cdot\exp(y)$,
$\exp(\log(1+x))=1+x$ et $\log(\exp(x))=x$.
\qed
\bigskip

\centerline{\soustitre Fin de la digression}

\bigskip

\prop

\medskip

{\sl Soit $d>0$ un entier, et $m>v_\P(d)+{v_\P(p)\over p-1}$. Alors tout ŽlŽment de $U_\P^{(m)}$ est une
puissance $d$-ime d'un ŽlŽment de $U_\P$. En particulier, si $d=[\bbL_\gP,\bbK_\P]$ et
$m>v_\P(d)+{v_\P(p)\over p-1}$, alors on a
$$U_\P^{(m)}\subset N_\P(U_\gP).$$

}

{\bf Preuve}

Puisque $m>v_\P(d)+{v_\P(p)\over p-1}>{v_\P(p)\over p-1}$, en vertu du corollaire prŽcŽdent,
l'application $\log\, :\, U_\P^{(m)}\to\widehat{\P}^m$ est un isomorphisme de groupe. Soit $1+x\in 
U_\P^{(m)}$ et $y=\log(1+x)\in\widehat{\P}^m$ (donc, $v_\P(y)\geq m$). Ainsi,
$$v_\P({y\over d})=v_\P(y)-v_\P(d)\geq m-v_\P(d)>{v_\P(p)\over p-1}>0,$$ 
donc, d'une part, ${y\over d}\in \widehat{\P}$ (sa valuation est $\geq 1$) et on peut en prendre son
exponentielle. Posons donc $z=\exp({y\over d})\in 1+\widehat{\P}\subset U_\P$. On a alors~:
$z^d=\exp(d\cdot{y\over d})=\exp(y)=1+x$. Cela montre que tout ŽlŽment de $U_\P^{(m)}$ est une
puissance $d$-ime d'un ŽlŽment de $U_\P$. Montrons la seconde partie de la proposition. Si $1+x\in 
U_\P^{(m)}$, on vient de voir que $1+x=z^d$ avec $z\in U_\P$. Or $z^d=N_\P(z)$ si $d=[\bbL_\gP,\bbK_\P]$. Cela
montre la proposition.\qed
\bigskip\goodbreak
\coro
\medskip

{\sl Si $d=[\bbL_\gP:\bbK_\P]$ et si $m>v_\P(d)+{v_\P(p)\over p-1}$. Alors

$$a(\P^m)=f_\P\cdot [U_\P:N_\P(U_\gP)]$$

}

{\bf Preuve}

C'est maintenant du tout cuit~:

$$\eqalign{a(\P^m)&\buildrel \rm Lemme\ \argby\ b)\over =f_\P\cdot [K^*(\m):(K^*(\m)\cap N(L^*))\cdot K^*_\m]
\cr&\buildrel
\rm Lemme\ \argbz\over =f_\P\cdot [U_\P:N(U_\gP) U_\P^{(m)}]\buildrel \rm Proposition\ \argcc\over = f_\P\cdot
[U_\P:N_\P(U_\gP)].\cr}$$ 
\qed

\bigskip

Reprenons un peu de cohomologie cyclique (mais cette fois dans le cas o l'extension est $\bbL_\gP/\bbK_\P$,
le groupe de Galois est $Z(\gP)$ et le $Z(\gP)$-module est $U_\gP$). Suivant les dŽfinitions, on a
$\ker(\Delta | U_\gP)=\{ x\in U_\gP\mid \sigma(x)=x\ \forall x\in Z(\gP)\}=U_\gP\cap \bbK_\P=U_\P$. Ainsi,
$$[U_\P : N_\P(U_\gP)]=|H^0(U_\gP)|=q^{-1}(U_\gP)\cdot |H^1(U_\gP)|,\eqno{(*)}$$
o $q(U_\gP)$ est le quotient de Herbrand. Ainsi, il ne nous reste plus qu'ˆ calculer $q(U_\gP)$ et
$|H^1(U_\gP)|$.

Jusqu'ici l'extension $L/K$ Žtait supposŽe cyclique. Nous allons affaiblir cette hypothse pour obtenir des rŽsultats intŽressants (les ThŽormes \argch\ et \argci). Nous supposerons que $L/K$ est une extension galoisienne quelconque, mais nous imposerons seulement que $\bbL_\gP/\bbK_\P$ soit cyclique.
\bigskip
\goodbreak
\lem
\medskip

{\sl
Soit $L/K$ une extension galoisienne de corps de nombres, $\P\in\gfP_0(K)$, $\gP\in \gfP_0(L)$, $\gP |\P$ tels que  $\bbL_\gP/\bbK_\P$ soit cyclique. Alors on a~:
$$ |H^1(U_\gP)|=e_\P.$$

}
{\bf Preuve}

Par le ThŽorme 90 de Hilbert, on a $\ker(N |{U_\gP})=U_\gP\cap \Delta(\bbL^*_\gP)$. Mettons que
$Z(\gP)=<\tau >$. Soit $x\in \bbL^*_\gP$. Puisque $\tau(\gP)=\gP$, on a $v_\gP(\tau(x))=v_\gP(x)$; donc
$\Delta(x)={x\over\tau(x)}\in U_\gP$. Cela prouve que $\Delta(\bbL^*_\gP)\subset U_\gP$ et donc que $\ker(N
|{U_\gP})=\Delta(\bbL^*_\gP)$, et ainsi, $H^1(U_\gP)\simeq \Delta (\bbL^*_\gP)/ \Delta(U_\gP)$. D'autre
part, l'application
$$\bbL^*_\gP\buildrel\Delta\over \longrightarrow\Delta(\bbL^*_\gP)\to\Delta(\bbL^*_\gP)/\Delta(U_\gP)$$
est Žvidemment surjective. On prŽtend que le noyau est $\bbK^*_\P\cdot U_\gP$. Clairement, $\bbK^*_\P\cdot
U_\gP$ est inclu dans le noyau, car $\Delta (\bbK^*_\P)=\{1\}$. RŽciproquement, si $x$ est un ŽlŽment du
noyau, cela veut dire que $\Delta(x)=\Delta(u)$ pour un $u\in U_\gP$. Donc, $\Delta({x\over u})=1$, i.e.
${x\over u}=a\in \bbK^*_\P$, et donc $x=a\cdot u\in \bbK^*_\P\cdot U_\gP$. On a ainsi montrŽ que 
$$H^1(U_\gP)\simeq \bbL^*_\gP/(\bbK^*_\P\cdot U_\gP).$$
Si $\pi$ est une uniformisante (un gŽnŽrateur de $\widehat{\gP}$), l'application $\pi^k\cdot u\mapsto
(k,u)$ est un isomorphisme de $\bbL^*_\gP$ sur $\Z\times U_\gP$. Puisque $\pi^{e_\P}$ est un gŽnŽrateur de
$\widehat{\P}O_{\gP}$, tout gŽnŽrateur de $\widehat\P$ peut s'Žcrire $\pi^{e_\P}\cdot u$, o $u\in U_\gP$. Ainsi,  l'image de $\bbK^*_\P\cdot U_\gP$ via cet isomorphisme est $e_\P\Z\times U_\gP$.
D'o, $H^1(U_\gP)\simeq \Z/e_\P\Z$, ce qui montre le lemme.\qed
\bigskip\goodbreak
\lem
\medskip

{\sl
Soit $L/K$ une extension galoisienne de corps de nombres, $\P\in\gfP_0(K)$, $\gP\in \gfP_0(L)$, $\gP |\P$ tels que  $\bbL_\gP/\bbK_\P$ soit cyclique. Alors on a~:
$$q(U_\gP)=1.$$

}

{\bf Preuve}

Puisque $O_\gP$ est local, alors comme lors du Lemme \argbz, on voit que $U_\gP/U_\gP^{(m)}$ est le groupe des unitŽs de $O_\gP/\widehat{\gP}^{(m)}$ qui est fini, ainsi, $[U_\gP :U_{\gP}^{(m)}]<\infty$. Cela prouve, en vertu du Corollaire \argbl, que $q(U_\gP)=q(U_\gP^{(m)})$, pour tout $m>1$. Le Corollaire \argcb\ (appliquŽ ˆ $\bbL_\gP$) montre que si $m$ est assez grand, l'application $\log\, : \, U_\gP^{(m)}\to \widehat{\gP}^m$ est un isomorphisme de groupe. De plus, tout ŽlŽment du groupe de Galois $\Z(\gP)$ est continu pour la topologie $\gP$-adique (cela se vŽrifie aisŽment en utilisant le fait que $v_\gP(\sigma(x))=v_\gP(x)$ pour tout $x\in \bbL_\gP$ et $\sigma\in Z(\gP)$). Ainsi,
$\sigma(\log(x))=\log(\sigma(x))$ pour tout $\sigma\in Z(\gP)$, d'o l'application $\log$ est un
isomorphisme de $Z(\gP)$-module. On en dŽduit donc (gr‰ce au Lemme \argbi , dans le cas o $C$
est trivial) que $q(U_\gP)=q(\widehat{\gP}^m)$ si $m$ est assez grand. Ensuite, puisque $\widehat{\gP}^m$ est d'indice fini dans $O_\gP$, on en dŽduit, gr‰ce au Corollaire \argbl , que
$q(U_\gP)=q(O_\gP)$. Par le thŽorme de la base normale (cf. [La1, Thm. 6.13.1, p. 320]), il existe $\omega\in
\bbL_\gP$ tel que $\tau^1(\omega),\ldots , \tau^d(\omega)$ forme une base de $\bbL_\gP$ sur $\bbK_\P$, o
$d=[\bbL_\gP :\bbK_\P]$ et $\tau$ est un gŽnŽrateur de $Z(\gP)$. Quitte ˆ multiplier $\omega$ par un ŽlŽment de $\bbK_\P$ (par exemple une bonne puissance de $\pi^{e_\P}$), on peut supposer que $\omega\in O_\gP$. Soit $M$, le sous-$O_\gP$-module engendrŽ par les $\tau^i(\omega)$. Puisque $O_\gP$ et $M$ sont des $O_\P$-modules libres de mme rang. Donc $M$ est d'indice fini dans $O_\gP$. On en dŽduit, gr‰ce au Corollaire \argbl, que $q(U_\gP)=q(M)$. Finalement, $M$ est un $O_\P$-module de permutation (dans le sens de la Proposition \argbm ), avec dans notre cas, $m=1$, $G=Z(\gP)$ et $R=O_\P$. Et cela prouve que $q(U_\gP)=q(M)=1$.\qed

\bigskip

\th
\medskip
{\sl Soit $L/K$, une extension cyclique de corps de nombres, $\m$ un $K$-module tel que les exposants
des diviseurs premiers finis de $\m$ soient suffisamment grands, alors on a

$$a(\m)=[K^*:N_{L/K}(L^*) K^*_\m]=\prod_{\P|\m} e_\P f_\P.$$

}

{\bf Preuve}

Cela dŽcoule de ce qui prŽcde~: si $\m=\prod_{\P|\m_0}\P^{n_\P}\cdot
\prod_{\P|\m_\infty}\P^{n_\P}$, alors

$$\eqalign{a(\m)&\buildrel {\rm Lemme\ \argbv}\over =\prod_{\P|\m_0}a(\P^{n_\P})\cdot
\prod_{\P|\m_\infty}a(\P^{n_\P})\buildrel {\rm Lemme\ \argbw}\over = \prod_{\P|\m_0}a(\P^{n_\P})\cdot
\prod_{\P|\m_\infty} e_\P f_\P \cr
&\buildrel {\rm Corollaire\ \argcd}\over =\prod_{\P|\m_0}f_\P\cdot [U_\P:N_\P(U_\gP)]\cdot
\prod_{\P|\m_\infty} e_\P f_\P\buildrel(*)\over =\prod_{\P|\m_0}f_\P\cdot q^{-1}(U_\gP)\cdot |H^1(U_\gP)| \cdot
\prod_{\P|\m_\infty} e_\P f_\P \cr
&\buildrel {\rm Lemmes\ \argce\ et\ \argcf}\over =\prod_{\P|\m} e_\P f_\P.\cr}$$\qed

\bigskip\goodbreak

\centerline{\soustitre DIGRESSION}
\bigskip

\th
\medskip

{\sl Soit $L/K$, une extension galoisienne de corps de nombres et $\P$ un idŽal premier de $K$ non ramifiŽ dans $L$ et $\gP$ un idŽal de $L$ au-dessus de $\P$. Alors
$$U_\P=N_\P(U_\gP),$$
o $N_\P=N_{\bbL_\gP/\bbK_\P}$.

}
\bigskip\goodbreak
{\bf Preuve}

Puisque $\P$ ne ramifie pas, $Z(\P)$ est cyclique. Les Lemmes \argce\ et \argcf\ s'appliquent alors. On a donc $q(U_\gP)=1$ et  $|H^1(U_\gP)|=e_\P=1$. On en dŽduit que  $H^0(U_\gP)=1$, d'o la proposition, puisque par dŽfinition dans notre cas $H^0(U_\gP)=\ker(\Delta| U_\gP)/N_\P(U_\gP)=U_\P/N_\P(U_\gP)$ (voir le raisonnement aprs le Corollaire \argcd). Ce qui prouve la proposition.\qed

\bigskip
Un thŽorme semblable sera dŽmontrŽ au Chapitre 11 (Corollaire  \argeo)

\bigskip
\prop
\medskip

{\sl Soit $L/K$, une extension galoisienne de corps de nombres et $\P$ un idŽal premier de $K$ ramifiŽ ou non dans $L$ et $\gP$ un idŽal de $L$ au-dessus de $\P$ tel que $Z(\gP/\P)$ est cyclique. Alors

$$[\bbK^*_\P:N_{\bbL_\gP/\bbK_{\P}}(\bbL^*_\gP)]=[\bbL_\gP:\bbK_\P].$$

}

{\bf Preuve}

On a dŽjˆ vu au Lemme \argce\ que $\bbL^*_\gP\simeq \Z\times U_\gP$ (via une uniformisante). Ainsi,
$\bbL^*_\gP/U_\gP\simeq\Z$ comme $Z(\gP)$-module ($\Z$ Žtant un $Z(\gP)$-module trivial). Donc, en vertu de la Proposition \argbm, $q(\bbL^*_\gP/U_\gP)={1\over |Z(\gP)|}$. D'autre part, on sait (Lemme \argcf) que $q(U_\gP)=1$ et gr‰ce ˆ la suite exacte $1\to U_\gP\to \bbL^*_\gP\to \bbL^*_\gP/U_\gP\to 1$ et au Lemme \argbi, on tire que $q(\bbL^*_\gP)={1\over |Z(\gP)|}$. D'autre part, le ThŽorme 90 de Hilbert (ou la Proposition \argbo\ d)) nous dit que $|H^1(\bbL_\gP^*)|=1$. Finalement,
$${1\over |Z(\gP)|}=q(\bbL^*_\gP)={|H^1(\bbL_\gP^*)|\over |H^0(\bbL_\gP^*)|}={1\over |H^0(\bbL_\gP^*)|},$$
ce qui montre que $|\bbK^*_\P/N_{\bbL_\gP/\bbK_{\P}}(\bbL_\gP)|=|H^0(\bbL_\gP^*)|=|Z(\gP)|=[\bbL_\gP:\bbK_\P]$.
\qed
\bigskip
Un rŽsultat similaire sera aussi vu au Chapitre 11 (ThŽorme \argen)

\vfill\eject

\global\advance\chapnomb by 1
\nomb=1

\centerline{\para Chapitre 6 :}
\smallskip
\centerline{\para L'ŽgalitŽ fondamentale du corps de  classe }
\smallskip
\centerline{\para pour les extensions cycliques et }
\smallskip
\centerline{\para thŽorme de la norme de Hasse }
\bigskip

Dans ce chapitre, on va prouver que l'inŽgalitŽ montrŽe au ThŽorme \argau\ est en fait une
ŽgalitŽ dans le cas des extensions cycliques. Comme premier corollaire, on en dŽduira le thŽorme de la norme de Hasse.

Soit $L/K$ une extension cyclique de corps de nombres de groupe $G=<\sigma>$. Posons $N=N_{L/K}$ la
norme de $L$ sur $K$. Soit $\m$ un $K$-module. Rappelons que l'application $f=f_\m\ : L^*\to
I_L(\widetilde{\m})$ est la composition de l'application $\iota\ : L^*\to I_L$, $x\mapsto x O_L$ et de
l'application $j_\m\ : I_L\to I_L(\widetilde{\m})$, $j_\m(\gP)=\gP$ si $\gP\notdiv \m$ et $j_\m(\gP)=O_L$
si $\gP |\m$ et on prolonge par multiplicativitŽ. Clairement, $f$ est un homomorphisme de $G$-modules. On
cherche, dans un premier temps, des renseignements sur l'application $f_0\ :K^*/N(L^*)=H^0(L^*)\to
H^0(I_L(\widetilde{\m}))=I_K(\m)/N(I_L(\widetilde{\m}))$, associŽe ˆ $f$ (cf. dŽbut du Chapitre 4). 
\bigskip
{\bf Petite remarque :}
\medskip

Tout carrŽ commutatif

\vglue 1cm 
\psset{xunit=0.8cm,yunit=0.8cm}

\rput(8,0){\rnode{apr}{$A'$}}
\rput(8,2){\rnode{a}{$A$}}
\rput(10,0){\rnode{bpr}{$B'$}}
\rput(10,2){\rnode{b}{$B$}}
\ncline[nodesep=3pt]{->}{a}{apr}
\Bput{$\alpha $} 
\ncline[nodesep=3pt]{->}{a}{b}
\Aput{$f $} 
\ncline[nodesep=3pt]{->}{apr}{bpr}
\Aput{$f'$} 
\ncline[nodesep=3pt]{->}{b}{bpr}
\Aput{$\beta$} 

d'homomorphismes de $G$-modules se prolonge de manire unique en un diagramme commutatif aux lignes exactes

\vglue 2cm 
\psset{xunit=0.8cm,yunit=0.8cm}

\rput(4,0){\rnode{un1}{$1 $}} 
\rput(4,2){\rnode{un2}{$1$}}
\rput(6,0){\rnode{kerfpr}{$\ker(f')$}}
\rput(6,2){\rnode{kerf}{$\ker(f)$}}
\rput(8,0){\rnode{apr}{$A'$}}
\rput(8,2){\rnode{a}{$A$}}
\rput(10,0){\rnode{bpr}{$B'$}}
\rput(10,2){\rnode{b}{$B$}}
\rput(12,0){\rnode{cokerfpr}{${\rm coker}(f')$}} 
\rput(12,2){\rnode{cokerf}{${\rm coker}(f)$}}
\rput(14,0){\rnode{un3}{$1$}}
\rput(14,2){\rnode{un4}{$1$}}  
\ncline[nodesep=3pt]{->}{un2}{kerf}
\ncline[nodesep=3pt]{->}{un1}{kerfpr}
\ncline[nodesep=3pt]{->}{kerf}{kerfpr}
\Aput{$\delta $} 
\ncline[nodesep=3pt]{->}{kerf}{a}
\ncline[nodesep=3pt]{->}{kerfpr}{apr}
\ncline[nodesep=3pt]{->}{a}{apr}
\Aput{$\alpha $} 
\ncline[nodesep=3pt]{->}{a}{b}
\Aput{$f $} 
\ncline[nodesep=3pt]{->}{apr}{bpr}
\Aput{$f'$} 
\ncline[nodesep=3pt]{->}{b}{bpr}
\Aput{$\beta$} 
\ncline[nodesep=3pt]{->}{b}{cokerf}
\ncline[nodesep=3pt]{->}{bpr}{cokerfpr}
\ncline[nodesep=3pt]{->}{cokerf}{cokerfpr}
\Aput{$\gamma$} 
\ncline[nodesep=3pt]{->}{cokerf}{un4}
\ncline[nodesep=3pt]{->}{cokerfpr}{un3}

\vglue 1cm 
\headline{\hfill\phantom{ouuh}\hfill}
L'application $\delta=\alpha |_{\ker(f)}$ est bien dŽfinie, car $\alpha(\ker(f))\subset\ker(f')$. Puisque,
de plus, $\beta({\rm Im}(f))\subset {\rm Im}(f')$ et que par dŽfinition, ${\rm coker}(f)=B/{\rm Im}(f)$,
l'application $\gamma$ est obtenue par passages aux quotients de $\beta$. De plus, $\gamma$ est surjective
si $\beta$ l'est et $\delta$ est injective si $\alpha$ l'est.

{\bf fin de la petite remarque.}
\medskip
Puisque, clairement, on a $f(K^*)\subset I_K(\m)$, $f(K^*_\m)=P_\m$ et $f(N(L^*))\subset
N(I_L(\widetilde{\m}))$, on peut considŽrer (gr‰ce ˆ la petite remarque, utilisŽe 2 fois) le diagramme commutatif
suivant~: 

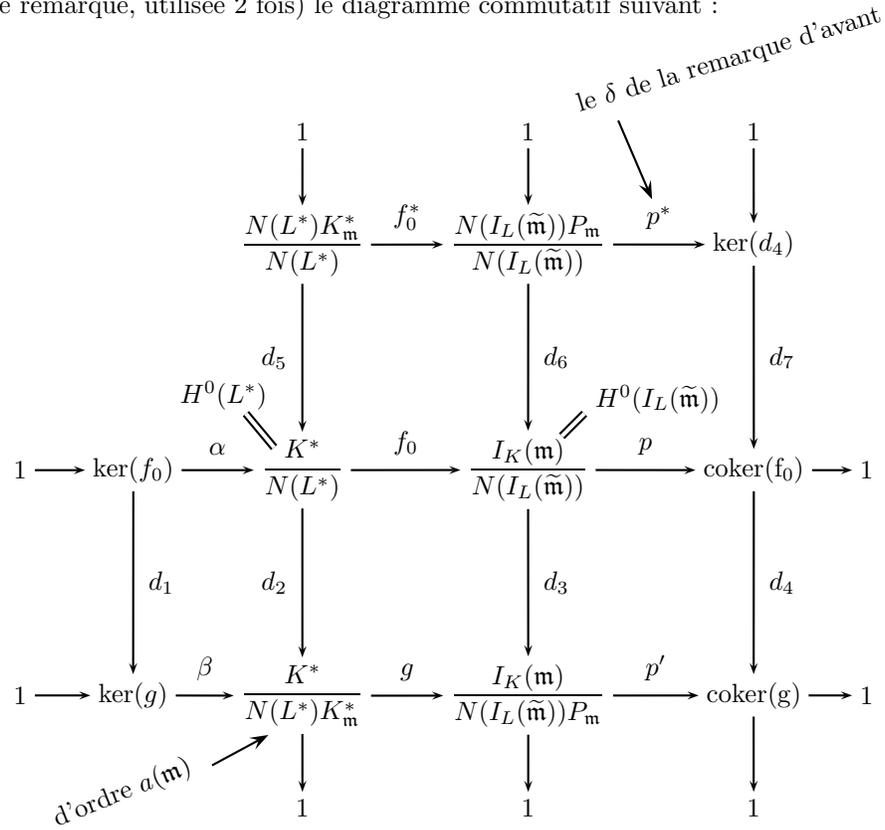
\vglue 11.5cm 

\psset{xunit=1.5cm,yunit=1.5cm}

\rput(6,1){\rnode{1b2}{$1  $}} 
\rput(4,1){\rnode{1b1}{$1  $}}
\rput(8,1){\rnode{1b3}{$1  $}} 
\rput(4,7){\rnode{1h1}{$1  $}} 
\rput(6,7){\rnode{1h2}{$1  $}} 
\rput(8,7){\rnode{1h3}{$1  $}} 
\rput(4,2){\rnode{c31}{$\dst {K^*\over N(L^*)K_{\m}^* } $}} 
\rput(4,4){\rnode{c21}{$\dst{K^*\over N(L^*)  }$}} 
\rput(4,6){\rnode{c11}{$\dst{N(L^*)K_{\m}^*\over N(L^*) } $}} 
\rput(6,2){\rnode{c32}{$\dst{I_{K}(\m)\over N(I_{L}(\widetilde{\m}))P_{\m} } $}} 
\rput(6,4){\rnode{c22}{$\dst{I_{K}(\m)\over N(I_{L}(\widetilde{\m})) } $}} 
\rput(6,6){\rnode{c12}{$\dst{N(I_{L}(\widetilde{\m}))P_{\m}\over N(I_{L}(\widetilde{\m})) } $}} 
\rput(8,2){\rnode{c33}{$\rm coker (g) $}} 
\rput(8,4){\rnode{c23}{$\rm coker  (f_{0}) $}}
\rput(8,6){\rnode{c13}{$\ker  (d_{4}) $}}
\rput(2.5,2){\rnode{c30}{$\ker  (g)$}}
\rput(2.5,4){\rnode{c20}{$\ker  (f_{0})$}}
\rput(1.5,2){\rnode{1g1}{$1$}}
\rput(1.5,4){\rnode{1g2}{$1$}}
\rput(9,2){\rnode{1d1}{$1$}}
\rput(9,4){\rnode{1d2}{$1$}}
\ncline[nodesep=3pt]{->}{c31}{1b1} 
\ncline[nodesep=3pt]{->}{c21}{c31}
\Bput{$d_{2}$} 
\ncline[nodesep=3pt]{->}{c11}{c21}
\Bput{$d_{5}$}  
\ncline[nodesep=3pt]{->}{1h1}{c11}
\ncline[nodesep=3pt]{->}{c32}{1b2}
\ncline[nodesep=3pt]{->}{c31}{c32}
\Aput{$g$}
\ncline[nodesep=3pt]{->}{c21}{c22}
 \Aput{$f_{0}$}
\ncline[nodesep=3pt]{->}{c22}{c32}
 \Aput{$d_{3}$}
 \ncline[nodesep=3pt]{->}{c11}{c12}
 \Aput{$f_{0}^*$}
 \ncline[nodesep=3pt]{->}{c12}{c22}
 \Aput{$d_{6}$}
 \ncline[nodesep=3pt]{->}{1h2}{c12}
\ncline[nodesep=3pt]{->}{c12}{c13}
\Aput{$p^*$}
\ncline[nodesep=3pt]{->}{c22}{c23}
\Aput{$p$}
\ncline[nodesep=3pt]{->}{c32}{c33}
\Aput{$p'$}
\ncline[nodesep=3pt]{->}{1h3}{c13}
\ncline[nodesep=3pt]{->}{c13}{c23}
\Aput{$d_{7}$}
\ncline[nodesep=3pt]{->}{c23}{c33}
\Aput{$d_{4}$}
\ncline[nodesep=3pt]{->}{c33}{1b3}
\ncline[nodesep=3pt]{->}{c30}{c31}
\Aput{$\beta $}
\ncline[nodesep=3pt]{->}{c20}{c21}
\Aput{$\alpha$}
\ncline[nodesep=3pt]{->}{c20}{c30}
\Aput{$d_{1}$}
\ncline[nodesep=3pt]{->}{1g1}{c30}
\ncline[nodesep=3pt]{->}{1g2}{c20}
\ncline[nodesep=3pt]{->}{c33}{1d1}
\ncline[nodesep=3pt]{->}{c23}{1d2}
\psline[doubleline=true,doublesep=1.5pt](6.5,4.5)(6.3,4.3) 
\rput(7.15,4.63){$H^0(I_{L}(\widetilde{\m}))$}
\psline[doubleline=true,doublesep=1.5pt](3.5,4.5)(3.75,4.2) 
\rput(3.3,4.66){$H^0(L^*)$}
\pcline[linecolor=white](7.7,7.36)(8,7.45) \aput{:U}{le $\delta$ de la remarque d'avant}
\psline[arrows=->,arrowsize=5pt](6.8,7.1)(7.1,6.4)
\pcline[linecolor=white](2,0.7)(3,1.1) \aput{:U}{d'ordre $a(\m)$}
\psline[arrows=->,arrowsize=5pt](3.2,1.38)(3.7,1.65)


o $f_0^*,f_0$ et $g$ sont induites par $f$. Les applications $d_5, d_6$, $d_7$, $\alpha$ et $\beta$ sont des inclusions, $d_2$ et $d_3$ sont les surjections naturelles. Les applications $d_1$, $d_4$ et $p^*$ sont les applications construites lors de la petite remarque. Puisque $d_3$ est surjectif, $d_4$ l'est aussi. Les trois colonnes sont exactes ainsi que les deux dernires lignes. Maintenant, nous allons faire une petite partie de chasse dans ce diagramme ainsi que quelques applications des thŽormes d'isomorphismes. Et en 6 affirmations, nous allons prouver le rŽsultat principal de ce chapitre.

 \newcount\gaagab\gaagab=\gaga
{\bf Premire affirmation~:} $p^*$ est surjective

En effet, soit $x\in \ker(d_4)$. Il existe $y\in I_K(\m)/N(I_L(\widetilde{\m}))$ tel que $p(y)=d_7(x)$.
On a donc, par commutativitŽ, $p'(d_3(y))=d_4(p(y))=d_4(d_7(x))=1$. Donc, $d_3(y)\in\ker(p')={\rm
Im}(g)$; il existe donc $z\in K^*/N(L^*)K^*_\m$ tel que $g(z)=d_3(y)$. Puisque $d_2$ est surjective, il
existe $u$ tel que $d_2(u)=z$. Alors, $d_3(f_0(u))=g(d_2(u))=g(z)=d_3(y)$. Cela implique que ${y\over
f_0(u)}\in\ker(d_3)={\rm Im}(d_6)$. Il existe donc $v$ tel que ${y\over f_0(u)}=d_6(v)$. Alors,
$d_7(p^*(v))=p(d_6(v))=p({y\over f_0(u)})=p(y)=d_7(x)$. Cela implique que $p^*(v)=x$, puisque $d_7$ est
injective. Cela montre notre premire affirmation.

\headline{\hfill \smcap L'ŽgalitŽ fond. du corps de  classe pour les ext. cycliques et 
 thm. de la norme de Hasse\hfill}

{\bf Deuxime affirmation~:}

$${\rm coker}(f_0)\simeq {\rm coker}(g)\eqno{(1)}.$$

En effet, de l'ŽgalitŽ $d_7\circ p^*\circ f_0^*=p\circ f\circ d_5=1$ et du fait que $d_7$ est injective,
on tire que $p^*\circ f_0^*=1$. D'autre part, $f(K^*_\m)=P_\m$, donc $f^*_0$ est surjective. Gr‰ce ˆ ces
deux faits et puisque $p^*$ est surjective, on en dŽduit que $\ker(d_4)=1$ et donc que $d_4$ est un
isomorphisme entre ${\rm coker}(f_0)$ et ${\rm coker}(g)$.
\medskip
On dŽfinit
\newcount\gaagaa \gaagaa=\gaga
$$\thboxed 25{n(\m)=[K^*_\m\cap\iota^{-1}(N(I_L(\widetilde{\m}))): K^*_\m\cap N(L^*)]}.$$

\goodbreak
{\bf Troisime affirmation~:}

$$|\ker(f_0) |=|\ker(g)|\cdot n(\m) \eqno{(2)}.$$

En effet, montrons dŽjˆ que $d_1$ est surjective. Soit donc $x\in \ker(g)$. Puisque $d_2$ est surjective,
il existe $y$ tel que $\beta(x)=d_2(y)$. On a alors $d_3(f_0(y))=g(d_2(y))=g(\beta(x))=1$. Cela prouve que
$f_0(y)\in\ker(d_3)={\rm Im}(d_6)$, il existe donc $z$ tel que $f_0(y)=d_6(z)$. Puisque $f_0^*$ est
surjective (deuxime affirmation), il existe $u$ tel que $f_0^*(u)=z$. Calculons~: $f_0(d_5(u))=d_6(f_0^*(u))=d_6(z)=f_0(y)$.
Donc ${y\over d_5(u)}\in \ker(f_0)$ et alors, $\beta(d_1({y\over d_5(u)}))=d_2(\alpha({y\over
d_5(u)}))\buildrel \alpha\ {\rm incl.}\over = d_2({y\over d_5(u)})=d_2(y)=\beta(x)$, ce qui prouve que
$x=d_1({y\over d_5(u)})$, puisque $\beta$ est injective. Ce qui montre que $d_1$ est surjective. On a
ainsi prouvŽ que 
$$|\ker(g)|={|\ker(f_0)|\over |\ker(d_1)|}.$$
Mais, puisque $\alpha$ est l'inclusion, on a $\ker(d_1)=\ker(f_0)\cap\ker(d_2)$. Soit $x\cdot
N(L^*)\in\ker(f_0)\cap\ker(d_2)$. Il existe donc $y\in K^*_\m$ tel que $x\cdot N(L^*)=y\cdot N(L^*)$ et
$f(y)\in N (I_L(\widetilde{\m}))$. C'est-ˆ-dire (voir la dŽfinition de $f$) que $\iota(y)\in 
N(I_L(\widetilde{\m}))$, ou encore, $y\in \iota^{-1}( N (I_L(\widetilde{\m})))\cap K^*_\m$. Enfin, si $y$
et $y_1\in K_\m^*$, alors $y\cdot N(L^*)=y_1\cdot N(L^*)\iff yy_1^{-1}\in N(L^*)$. Ainsi,
$$\ker(f_0)\cap\ker(d_2)={ K^*_\m\cap\iota^{-1}( N (I_L(\widetilde{\m})))\over K_{\m}^*\cap N(L^*)}$$
Ce qui montre que $|\ker(d_1)|=n(\m)$, ce qui montre la troisime affirmation.

\smallskip
Nous avons maintenant besoin d'un petit lemme facile, mais que nous rŽutiliserons par la suite~:
\medskip
\lem
\smallskip
{\sl Si $B$ est un sous-module d'indice fini d'un module $A$ et si $\beta$ est un homomorphisme de
module dŽfini sur $A$, alors on a
$$[A:B]=[\beta(A):\beta(B)]\cdot[\ker(\beta):B\cap \ker(\beta)].$$

}

{\bf Preuve}

$\beta$ induit un homomorphisme surjectif $A/B\to \beta(A)/\beta(B)$ dont le noyau est $\ker(\beta)\cdot
B/B\buildrel {\rm thm.\ d'isom.}\over \simeq \ker(\beta)/\ker(\beta)\cap B$. Et cela prouve notre lemme.\qed

\medskip
Posons $S$
l'ensemble des places de $L$ qui divisent $\m$ et ${L^*}^S$ l'ensemble des $S$-unitŽs qui est, rappelons-le,
le noyau $\ker(f_\m)=\{\alpha\in L^*\mid \iota(\alpha)$ n'est divisible que par des idŽaux premiers de
$S\}$. La suite exacte de $G$-modules et $G$-homomorphismes
$$1\to {L^*}^S\buildrel\gamma=\rm incl.\over\to L^*\buildrel f=f_\m\over\to I_L(\widetilde{\m})\buildrel\lambda\over\to
V:={\rm coker(f)}\to 1$$
se scinde en deux suites exactes (plus courtes)
$$1\to {L^*}^S\buildrel\gamma\over\to L^*\buildrel\alpha\over\to f(L^*)\to 1\quad\hbox{et}\quad
1\to f(L^*)\buildrel \beta\over\to I_L(\widetilde{\m})\buildrel\lambda\over\to
V\to 1,$$
o $\alpha$ est la surjection canonique et $\beta$ l'inclusion de $f(L^*)$ dans $I_L(\widetilde{\m})$. Ces
suites exactes correspondent aux hexagones exacts suivants (en se souvenant du Lemme \argbf\ et
au fait que $H^1(L^*)=H^1(I_L(\widetilde{\m}))=1$, (Proposition \argbo))~:
 \goodbreak

\vglue 2.5cm
\psset{xunit=0.9cm,yunit=0.9cm}
\rput(4,-1.2){\rnode{h1cpr}{$1 $}} 
\rput(7,-0.2){\rnode{h0apr}{$H^1(V)$}}
\rput(10,-0.2){\rnode{h0bpr}{$H^0(f(L^*))$}}
\rput(13,-1.2){\rnode{h0cpr}{$H^0(I_{L}(\widetilde{\m}))$}}
\rput(10,-2.2){\rnode{h1apr}{$H^0(V)$}}
\rput(7,-2.2){\rnode{h1bpr}{$H^1(f(L^*))$}}
\ncline[nodesep=2pt]{->}{h1cpr}{h0apr}
\ncline[nodesep=3pt]{->}{h0apr}{h0bpr}
\Aput{$\delta_{3}$} 
\ncline[nodesep=3pt]{->}{h0bpr}{h0cpr}
\Aput{$\beta_{0}$} 
\ncline[nodesep=3pt]{->}{h0cpr}{h1apr}
\Aput{$\lambda_{0}$} 
\ncline[nodesep=3pt]{->}{h1apr}{h1bpr}
\Aput{$\delta_{4}$} 
\ncline[nodesep=3pt]{->}{h1bpr}{h1cpr}
\rput(4,2){\rnode{h1c}{$1 $}} 
\rput(7,3){\rnode{h0a}{$H^1(f(L^*))$}}
\rput(10,3){\rnode{h0b}{$H^0({L^*}^S)$}}
\rput(13,2){\rnode{h0c}{$H^0(L^*)$}}
\rput(10,1){\rnode{h1a}{$H^0(f(L^*))$}}
\rput(7,1){\rnode{h1b}{$H^1({L^*}^S)$}}
\ncline[nodesep=3pt]{->}{h1c}{h0a}
\ncline[nodesep=3pt]{->}{h0a}{h0b}
\Aput{$\delta_{1}$} 
\ncline[nodesep=3pt]{->}{h0b}{h0c}
\Aput{$\gamma_{0}$} 
\ncline[border=2pt, nodesep=3pt]{->}{h0c}{h1a}
\Bput{$\alpha_{0}$} 
\ncline[nodesep=3pt]{->}{h1a}{h1b}
\Bput{$\delta_{2}$} 
\ncline[border=2pt , nodesep=3pt]{->}{h1b}{h1c}
\ncline[nodesep=3pt]{->}{h0c}{h0cpr}
\Aput{$f_{0}$} 
\ncline[doubleline=true,border=2pt , nodesep=3pt]{h1a}{h0bpr}
 
\vskip 3cm

On remarque que $f_0=\beta_0\circ\alpha_0$ (cf. Proposition \argbg).
\medskip

{\bf Quatrime affirmation~:}

On a

$$|{\rm coker(f_0)}|={|H^0(V)|\over |H^1(V)|}\cdot{|H^1({L^*}^S)|\over |H^1(f(L^*)|}\cdot |\ker(\beta_0)\cap {\rm
Im}(\alpha_0)|.\eqno{(3)}$$
\medskip
En effet, on a $|{\rm coker}(f_0)|=[H^0(I_L(\widetilde{\m})):{\rm
Im}(\beta_0\circ\alpha_0)]=[H^0(I_L(\widetilde{\m})):{\rm Im}(\beta_0)]\cdot [{\rm Im}(\beta_0):{\rm
Im}(\beta_0\circ\alpha_0)]$. D'autre part, on applique le Lemme \argcj, avec  $A=H^0(f(L^*))$ et $B={\rm
Im}(\alpha_0)$ et $\beta=\beta_0$. On trouve
$$[{\rm Im}(\beta_0): {\rm Im}(\beta_0\circ\alpha_0)]={[H^0(f(L^*)):{\rm Im}(\alpha_0)]\over
[\ker(\beta_0):\ker(\beta_0)\cap{\rm Im}(\alpha_0)]}.$$
D'autre part encore, $[H^0(f(L^*)):{\rm Im}(\alpha_0)]=|{\rm coker}(\alpha_0)|=|{\rm
Im}(\delta_2)|=|H^1({L^*}^S)|$. On a aussi, $[H^0(I_L(\widetilde{\m})):{\rm Im}(\beta_0)]=|{\rm
coker}(\beta_0)|=|{\rm Im}(\lambda_0)|$ et $|H^0(V)|=|{\rm Im}(\lambda_0)|\cdot |{\rm
Im}(\delta_4)|=|{\rm Im}(\lambda_0)|\cdot |H^1(f(L^*))|$. On obtient alors
$$|{\rm coker(f_0)}|={|H^0(V)|\over H^1(f(L^*))|}\cdot {|H^1({L^*}^S)|\over [\ker(\beta_0) : \ker(\beta_0)\cap{\rm
Im}(\alpha_0)]}.$$
Notons encore que $|\ker(\beta_0)|=|{\rm Im}(\delta_3)|=|H^1(V)|$. Cela prouve notre affirmation.

\medskip

{\bf Cinquime affirmation~:}
$$|\ker(f_0)|=|\ker(\beta_0)\cap {\rm Im}(\alpha_0)|\cdot{|H^0({L^*}^S)|\over |H^1(f(L^*))|}.\eqno{(4)}$$

En effet, $\ker(f_0)=\ker(\beta_0\circ\alpha_0)=\alpha_0^{-1}(\ker(\beta_0))$ est envoyŽ par $\alpha_0$ sur
$\ker(\beta_0)\cap {\rm Im}(\alpha_0)$ et a pour noyau $\ker(\alpha_0)\cap
\alpha_0^{-1}(\ker(\beta_0))=\ker(\alpha_0)$. Ainsi
$$|\ker(f_0)|=|\ker(\beta_0)\cap {\rm Im}(\alpha_0)|\cdot |\ker(\alpha_0)|=|\ker(\beta_0)\cap {\rm
Im}(\alpha_0)|\cdot |{\rm Im}(\gamma_0)|=|\ker(\beta_0)\cap {\rm Im}(\alpha_0)|\cdot{|H^0({L^*}^S)|\over
|H^1(f(L^*))|}.$$\qed
\medskip

{\bf Sixime affirmation~:}

$${|{\rm coker}(f_0)|\over |\ker(f_0)|}=q({L^*}^S).\eqno{(5)}$$

En effet, en combinant les Žquations (3) et (4), on trouve ${|{\rm coker}(f_0)|\over
|\ker(f_0)|}={q({L^*}^S)\over q(V)}$. Or, $V$ est fini, car $f(L^*)$ contient
$\iota(L^*_{\widetilde{\m}})=P_{\widetilde{\m}}$, donc $V\simeq I_L(\widetilde{\m})/f(L^*)$ est un
quotient de $I_L(\widetilde{\m})/P_{\widetilde{\m}}$ (groupe de classe gŽnŽralisŽ pour $L$) qui est fini
(cf. ThŽorme \argl) donc $q(V)=1$ (cf. Lemme \argbk). Cela prouve l'affirmation.
\bigskip

Nous voilˆ prt ˆ prouver l'ŽgalitŽ fondamentale. Mais, il faut encore mettre ensemble tous les
prŽparatifs et observer une ou deux choses. Du grand diagramme (p. \the\gaagab), on retient la ligne exacte~:
$$1\to\ker(g)\buildrel\beta\over\lra K^*/N(L^*)K_\m^*\buildrel g\over \lra I_K(\m)/N(I_L(\widetilde{\m}))
P_\m\buildrel p'\over \lra {\rm coker}(g)\to 1.$$
Alors, $[I_K(\m) :N(I_L(\widetilde{\m}))\cdot P_\m]=|{\rm coker}(g)|\cdot |{\rm Im(g)}|={|{\rm
coker}(g)|\over |\ker(g)|}\cdot a(\m)$ (voir au chapitre 5 pour la dŽfinition de $a(\m)$). En vertu des Žquation (1) et (2), c'est Žgal ˆ $a(\m)\cdot n(\m)\cdot {|{\rm coker}(f_0)|\over
|\ker(f_0)|}=a(\m)\cdot n(\m)\cdot q({L^*}^S)$ en vertu de l'Žquation (5). RŽsumons-nous, on a
$$[I_K(\m) :N(I_L(\widetilde{\m}))\cdot P_\m]=a(\m)\cdot n(\m)\cdot q({L^*}^S).\eqno{(6)}$$
Or, si $\m$ est divisible par toutes les places qui ramifient, on sait que 
$$q({L^*}^S)=[L:K]\cdot \prod_{\P|\m}{1\over e_\P f_\P} \quad \hbox{(cf. ThŽorme \argbs)}$$
et si les exposants des places finies de $\m$ sont assez grand, on a
$$a(\m)= \prod_{\P|\m}e_\P f_\P\quad \hbox{(cf. ThŽorme \argcg)}.$$
Supposons que ces deux conditions (sur $\m$) soient satisfaites. La premire ŽgalitŽ du corps de classe (cf. ThŽorme \argau) et ces deux derniers rŽsultats, combinŽs avec (6), nous donnent alors
$$n(\m)\cdot [L:K]=[I_K(\m) :N(I_L(\widetilde{\m}))\cdot P_\m]\leq [L:K].$$
On a donc montrŽ que, dans ce cas, $n(\m)=1$ et $[I_K(\m) :N(I_L(\widetilde{\m}))\cdot P_\m]= [L:K]$.
RŽsumons tout cela dans le 
\bigskip
\th {\bf (ŽgalitŽ fondamentale du corps de classe pour les ext. cycliques)}
\medskip
{\sl Si $L/K$ est une extension cyclique de corps de nombres et si $\m$ est un $K$-module divisible par
toutes les place ramifiŽes et dont les exposants des places finies qui le divisent sont suffisamment
grandes, alors on a~:

$$[I_K(\m) :N(I_L(\widetilde{\m}))\cdot P_\m]= [L:K]$$
et dans ces mmes conditions, on a $n(\m)=1$, ce qui peut s'Žnoncer comme suit~:

}
\medskip
\coro
\medskip

{\sl Sous les mmes hypothses, on a

$$K^*_\m\cap\iota^{-1}(N_{L/K}(I_L(\widetilde{\m})))=K^*_\m\cap N_{L/K}(L^*).$$

}
\qed
\bigskip
\coro {\soustitre (ThŽorme de la Norme de Hasse)}

{\sl Soit $L/K$ une extension cyclique de corps de nombres et $\alpha\in K^*$. Alors $\alpha$ est la norme
d'un ŽlŽment de $L$ si et seulement si $\alpha$ est un norme locale en toute idŽal premier $\P$ de $K$.
Autrement dit~:

$$\exists\beta\in L\hbox{ tel que } \alpha=N_{L/K}(\beta)\iff\forall\P\in \gfP(K),\exists\gP|\P\hbox{ et
}\beta_\P\in \bbL_{\gP}\hbox{ tel que }N_{\bbL_\gP/\bbK_\P}(\beta_\P)=\alpha,$$
o, rappelons-le, $\gfP(K)$ est l'ensemble des places (finie ou infinies) de $K$.

}

{\bf Preuve}

Si $\alpha=N_{L/K}(\beta)$, alors $\alpha=\prod_{\rho\in G}\rho(\beta)$, o $G={\rm Gal}(L/K)$. Soit $\P$
un idŽal premier de $K$ et $Z(\P)$, le groupe de dŽcomposition de $\P$, alors on sait que
$Z(\P)={\rm Gal}(\bbL_\gP/\bbK_\P)$ pour tout idŽal premier $\gP$ de $L$ au-dessus de $\P$. Si
$G=\sqcup_{i=1}^k\sigma_i Z(\P)$ est une dŽcomposition en classes de $G$ modulo $Z(\P)$, alors on a
$$\alpha=\prod_{\sigma\in Z(\P)}\sigma\left (\prod_{i=1}^k\sigma_i(\beta)\right
)=N_{\bbL_\gP/\bbK_\P}\left(\prod_{i=1}^k\sigma_i(\beta)\right ),$$  
et $\prod_{i=1}^k\sigma_i(\beta)\in L\subset
\bbL_\gP$.

RŽciproquement, supposons que $\alpha$ soit une norme locale partout. Montrons d'abord que l'idŽal
$(\alpha)=\alpha O_K$ est la norme d'un idŽal de $L$. Soit $\P$ un idŽal premier de $K$ tel que
$\P|(\alpha)$ et $\P^a$ est la puissance exacte de $\P$ dans la dŽcomposition de $(\alpha)$ en idŽaux
premiers. 
Si $\gP$ est un idŽal premier de $L$ au-dessus de $\P$ et $f=f(\gP/\P)$, alors, $N_{L/K}(\gP)=\P^f$.
Donc, si on arrive ˆ prouver que $f|a$, alors $\P^a=N_{L/K}(\gP^{a\over f})$ et on conclut par
multiplicativitŽ. On peut supposer par hypothse que
$$\alpha O_{\P}=N_{\bbL_\gP/\bbK_\P}(\beta_\P)\cdot O_{\P}=N_{\bbL_\gP/\bbK_\P}(\beta_\P\cdot O_{\gP}).$$

Notons $\widehat{\P}=\P O_{\P}$ et $\widehat{\gP}=\gP O_{\gP}$. Alors on a $\alpha
O_{\P}=\widehat{\P}^a$ et il existe $b\in\Z$ tel que $\beta_\P O_{\gP}=\widehat{\gP}^b$. Puisque
$N_{\bbL_\gP/\bbK_\P}(\widehat{\gP})=\widehat{\P}^f$, en Žlevant ˆ la puissance $b$, on trouve
$\widehat{\P}^a=\widehat{\P}^{bf}$. Cela montre que $a=bf$ donc que $f$ divise $a$. 

Soit $\m$ un $K$-module tel que $n(\m)=1$ (c'est possible, en vertu du ThŽorme \argck). Ecrivons
$\m=\prod_{\P_i\in \gfP(K)}\P_i^{b_i}$ et fixons, pour chaque $i$, $\gP_i$ une place au-dessus de $\P_i$,
$e_i=e(\gP_i/\P_i)$ et $\beta_i\in L_{\gP_i}$ tels que $\alpha=N_{\P_i}(\beta_i)$, o
$N_{\P_i}=N_{\bbL_{\gP_i}/\bbK_{\P_i}}$. Pour chaque $i$, on choisit un $\beta'_i\in L$ tel que
$\beta'_i\equiv\beta_i\pmodast{\widehat{\gP_i}^{b_i e_i}}$ (car $L$ est dense dans $L_{\gP_i}$, voir
encore la dŽfinition de $\widehat{\gP_i}^{b_i e_i}$ dans le cas infini, mais ici, a veut dire simplement
qu'ils ont le mme signe dans le plongement considŽrŽ). Par le ThŽorme d'approximation dŽbile (ThŽorme \argc), il existe
$\gamma\in L$ tel que pour chaque $i$, on ait
$$\left\{\eqalign{\gamma&\equiv\beta'_i\pmodast{\gP_i^{b_i e_i}}\cr \gamma&\equiv
1\pmodast{\gP^{b_ie_i}}\hbox{ pour les $\gP$ au-dessus de $\P_i$, $\gP\ne\gP_i$.}\cr}\right.\eqno{(*)}$$
Fixons un $i$ et considŽrons une dŽcomposition $G={\rm Gal}(L/K)=\bigsqcup_j\tau_jZ(\P_i)$ en choisissant
$\tau_i=1$. Ainsi, si $j\ne i$, $\tau_j\ne 1$ et donc, $\tau_j^{-1}(\gP_i^{b_i e_i})\ne \gP_i^{b_i e_i}$,
ce qui implique par $(*)$ que $\gamma\equiv 1\pmodast {\tau_j^{-1}(\gP_i^{b_i e_i})}$, ou encore
$\tau_j(\gamma)\equiv 1\pmodast {\gP_i^{b_i e_i}}$. On obtient alors~:
$$N_{L/K}(\gamma)=\prod_j\prod_{\sigma\in Z(\gP_i)}\tau_j\sigma(\gamma)\equiv\prod_{\sigma\in
Z(\gP_i)}\sigma(\gamma)=N_{\P_i}(\gamma)\equiv N_{\P_i}(\beta'_i)\pmodast{\gP_i^{b_i e_i}}.$$
On a aussi que $\alpha=N_{\P_i}(\beta_i)\equiv N_{\P_i}(\beta'_i)\pmodast{\widehat{\gP_i}^{b_i e_i}}$, mais
aussi modulo ${\gP_i^{b_i e_i}}$ (car elles sont dans $L^*$). On a donc montrŽ que
$\alpha\equiv N_{L/K}(\gamma)\pmodast{\gP_i^{b_i e_i}}$, donc aussi modulo $\P_i^{b_i}$, car ce sont des ŽlŽments de $K^*$. Ceci est vrai pour chaque $i$, ce qui implique que
$\alpha\equiv N_{L/K}(\gamma)\pmodast{\m}$, et alors, $\alpha N_{L/K}(\gamma)^{-1}\in K_\m^*$. Or on
a vu au dŽbut de la preuve que $(\alpha)$ Žtait la norme d'un idŽal de $L$, donc $(\alpha
N_{L/K}(\gamma)^{-1})$ aussi. Puisque $\alpha N_{L/K}(\gamma)^{-1}\in K_\m^*$ cet idŽal est en particulier
dans $I_L(\widetilde{\m})$. Cela montre que $\alpha N_{L/K}(\gamma)^{-1}\in
\iota^{-1}(N_{L/K}(I_L(\widetilde{\m})))$. Donc, puisqu'on a choisit $\m$ de sorte que $n(\m)=1$, le Corollaire \argcl, nous dit que $\alpha
N_{L/K}(\gamma)^{-1}\in N_{L/K}(L^*)$, et donc, $\alpha\in N_{L/K}(L^*)$.\qed

\bigskip

{\bf Remarque}

L'hypothse d'avoir une extension $L/K$ cyclique est obligatoire. En effet, Hasse ˆ montrŽ que dans
$\Q(\sqrt{-3},\sqrt{-39})$, le nombre $3$ n'est pas une norme, mais est une norme locale partout.
(Cf. [Has])

\vfill\eject
\global\advance\chapnomb by 1
\nomb=1

\centerline{\para Chapitre 7 :}
\medskip
\centerline{\para La loi de rŽciprocitŽ d'Artin }

\bigskip

La loi de rŽciprocitŽ d'Artin nous donne une description du noyau de l'application du mme nom. 

\bigskip

\defi
\medskip

Si $L/K$ est une extension abŽlienne de corps de nombres, on dit qu'un $K$-module est {\it admissible}
(pour $L/K$) s'il est divisible par tous les idŽaux premiers qui ramifient dans $L$ et si $P_\m\subset
\ker(\Phi_{L/K}|I_K(\m))$.
\newcount\gaagac\gaagac=\gaga
\bigskip
\headline={\hfill \phantom{ouuh}\hfill}
\lem
\medskip
{\sl Si $L/K$ est une extension abŽlienne de corps de nombres et si $\m$ est admissible, alors
$$\ker(\Phi_{L/K}|I_K(\m))=N_{L/K}(I_L(\widetilde{\m}))\cdot P_\m$$ et donc $\Phi_{L/K}$ induit un isomorphisme 

$$I_K(\m)/N_{L/K}(I_L(\widetilde{\m}))\cdot P_\m\simeq {\rm Gal}(L/K).$$

}

{\bf Preuve}

On sait que $N_{L/K}(I_L(\widetilde{\m}))\subset \ker(\Phi_{L/K}|I_K(\m))$ (cf. Corollaire \argg). Donc on a $N(I_L(\widetilde{\m}))\cdot P_\m\subset \ker(\Phi_{L/K}|I_K(\m))$. De plus,
$\Phi_{L/K}|I_K(\m)$ est surjective (cf. ThŽorme \argaq), de sorte que
$[I_K(\m):\ker(\Phi_{L/K}|I_K(\m))]=[L:K]$. Donc on a $[I_K(\m): N(I_L(\widetilde{\m}))\cdot P_\m]\geq
[L:K]$. Or, la premire inŽgalitŽ du corps de classe (cf. ThŽorme \argau) nous dit que
$[I_K(\m): N_{L/K}(I_L(\widetilde{\m}))\cdot P_\m]\leq [L:K]$. Cela prouve que
$\ker(\Phi_{L/K}|I_K(\m))=N(I_L(\widetilde{\m}))\cdot P_\m$ (puisqu'ils ont le mme indice dans
$I_K(\m)$ et que le premiers contient le second) et donc le lemme.\qed
\bigskip
\rem

Vous aurez remarquŽ que pour ce lemme nous n'avons pas eu besoin de l'ŽgalitŽ, mais seulement de la premire inŽgalitŽ.
\bigskip
{\soustitre Rappel}
\medskip
Soit $m$ un entier positif tel que $m\not\equiv 2\pmod 4$. Alors le $\Q$-module $(m)\cdot\infty$ est
admissible pour l'extension cyclotomique
$\Q(\zeta_m)/\Q$ ($\infty$ Žtant l'unique place infinie de $\Q$). Ce rŽsultat est le ThŽorme \argm. 

\bigskip

\lem
\medskip

{\sl Si $L/K$ est une extension abŽlienne de corps de nombres et si $\m$ est un $K$-module admissible, alors
tout $K$-module $\m'$ tel que $\m_0|\m_0'$ et $\m_\infty\subset\m'_\infty$ est aussi admissible.

}

{\bf Preuve}

C'est clair~: $\Phi_{L/K}|I_K(\m')$ est la restriction ˆ $I_K(\m')$ de $\Phi_{L/K}|I_K(\m)$ et
$P_{\m'}\subset P_\m\cap I_K(\m')$.\qed

\bigskip

\rem

Dans le rappel prŽcŽdent, on peut se passer de l'hypothse $m\not\equiv 2\pmod 4$. En effet, si
$m\equiv 2\pmod 4$, il est Žvident (et bien connu) que $\Q(\zeta_{m})=\Q(\zeta_{m\over 2})$. On a alors
$({m\over 2})\cdot\infty$ est admissible pour $\Q(\zeta_m)/\Q$. Et donc, par le lemme prŽcŽdent, on a
$(m)\cdot\infty$ est admissible pour $\Q(\zeta_m)/\Q$.
\bigskip
\lem
\medskip

{\sl Si $\m$ est un $K$-module admissible pour une extension abŽlienne $L/K$, alors $\m$ est aussi
admissible pour $E/K$ pour tout $L\supset E\supset K$.

}
{\bf Preuve}

C'est clair~: on sait que $\Phi_{E/K}=R\circ\Phi_{L/K}$ o $R:\ {\rm Gal}(L/K)\to{\rm Gal}(E/K)$ est
l'homomorphisme de restriction (cf. on a vu cela au Chapitre 0, p. \the\gaage).\qed

\bigskip
\headline={\hfill \smcap La loi de rŽciprocitŽ d'Artin \hfill}
\lem

{\sl Soit $L/K$ est une extension abŽlienne de corps de nombres et $\m$ un $K$-module admissible.
Soit $E/K$ une extension (quelconque) de corps de nombres. Si $\widetilde{\m}$ est l'extension de $\m$ ˆ
$E$ (voir Chapitre 0, p. \the\gaagi\ pour la prŽsentation de $\widetilde{\m}$). Alors $\widetilde{\m}$ est admissible pour $EL/E$ 

}
{\bf Preuve}

On sait que sur $I_E(\widetilde{\m})$, on a $R\circ\Phi_{EL}=\Phi_{L/K}\circ N_{E/K}$ o $R:\
{\rm Gal}(EL/E)\to {\rm Gal}(L/K)$ est l'homomorphisme de restriction (cf.
ThŽorme \argf). Pour conclure, il suffit de rappeler que $R$ est injective et de se
souvenir que $N_{E/K}(E^*_{\widetilde{\m}})\subset K^*_\m$ (cf. partie b) du Lemme \aargbu).
Donc $N_{E/K}(P_{\widetilde{\m}})\subset P_m$, ce qui prouve le lemme.\qed
\bigskip
\goodbreak
\prop
\medskip

{\sl Soit $m$ un nombre entier positif. Soit $K\subset E\subset K(\zeta_m)$, une extension de corps de
nombres. Si $\m'$ est l'extension ˆ $K$ du $\Q$-module $(m)\cdot\infty$ et si $\m$ est un $K$-module tel
que $\m_0'|\m_0$ et $\m'_\infty\subset\m_\infty$. Alors $\m$ est admissible pour $E/K$.

}

{\bf Preuve}

C'est un corollaire des 4 lemmes prŽcŽdents~: puisque $(m)\cdot\infty$ est admissible pour
$\Q(\zeta_m)/\Q$, par le Lemme \argcr, $\m'$ est admissible pour $K(\zeta_m)/K$. Par le Lemme \argcq, $\m'$ est admissible pour $E/K$ et par le Lemme \argcp, $\m$ est admissible pour $E/K$.\qed

\bigskip
On a dŽjˆ montrŽ cette proposition au Chapitre 0 (ThŽorme \argp), la preuve est ici plus simple et nous avions besoin de ce thŽorme pour montrer le thŽorme de $\rm\check C$ebotarev pour les extensions cyclotomiques. De plus, nous aurons parfois besoin des lemmes ayant permis la nouvelle preuve de ce rŽsultat.
\bigskip

\lem
\medskip

{\sl Soient $a$ et $r$ des entiers supŽrieurs ˆ 1 et $q$ un nombre premier. Alors il existe $p$ un nombre
premier tel que l'ordre de $a\bmod p$ soit $q^r$

}
{\bf Preuve}

On utilise la relation
$$\eqalignno{g(x)&={x^q-1\over x-1}={[(x-1)+1]^q-1\over x-1}&\cr
&=(x-1)^{q-1}+\cdots +\pmatrix{q\cr t\cr}(x-1)^{t-1}+\cdots +q&(*)\cr
&=x^{q-1}+x^{q-2}+\cdots +x+1.&\cr}$$
Si $n\geq 2$, $g(n)\geq 3$. En particulier, $g:=g(a^{q^{r-1}})\geq 3$. Soit $\ell$ un diviseur premier de
$g$. 

Si $\ell$ ne divise pas le dŽnominateur de $g$ qui est $a^{q^{r-1}}-1$, alors $p:=\ell \,|\, a^{q^r}-1$, ce
qui veut dire que
$a$ est d'ordre $q^r$ modulo~$p$, et le lemme est prouvŽ.

Si $\ell$ divise $a^{q^{r-1}}-1$, alors la relation $(*)$ prouve que $\ell=q$. On va montrer que dans ce
cas, $g$ n'est pas une puissance de $q$. On peut donc choisir un diviseur premier de $g$ diffŽrent de
$q$ et celui-ci conviendra par ce qui prŽcde.

\art{a)} Si $q$ est impair, dans la relation $(*)$ avec $x=a^{q^{r-1}}$, $(x-1)^{q-1}$ et les termes
$\pmatrix{q\cr t\cr}(x-1)^{t-1}$ pour $1<t<q$ sont divisibles par $q^2$, donc $q^2\notdiv g$. Donc si
$g$ est une puissance de $q$, on a $g=q$. Mais c'est impossible, car la relation $(*)$ montre que $g>q$.

\art{b)} Si $q=2$, $g=a^{2^{r-1}}+1$. Si $g$ est une puissance de 2, alors $a$ est impair, $a=2k+1$. Et
alors, $g=(2k+1)^{2^{r-1}}+1\equiv 2\pmod 4$. Donc $g=2$, mais on a vu que $g\geq 3$. 

Cela prouve que $g$ n'est pas une puissance de $q$ et donc le
lemme.\qed

\bigskip

\defi
\medskip

Dans un groupe abŽlien, $\sigma$ et $\tau$ sont dits {\it indŽpendants} si $<\sigma >\cap <\tau>=\{1\}$.
Deux entiers $a$ et $b$ sont dits {\it indŽpendants modulo $m$} si les classes de $a$ et de $b$ sont
indŽpendants dans $(\Z/m\Z)^*$.

\bigskip

\lem
\medskip

{\sl Soient $a,n$ des entiers supŽrieurs ˆ 1. Supposons que $n=q_1^{k_1}\cdots q_s^{k_s}$ (la
dŽcomposition de $n$ en puissance de nombres premiers). Alors il existe une infinitŽ d'entiers sans
facteurs carrŽs $m=p_1\cdots p_s\cdot p'_1\cdots p'_s$ avec $p_1<\cdots <p_s<p'_1< \cdots <p'_s$ premiers
et $p_1$ arbitrairement grand, tels que~:
\medskip
\art{1)}l'ordre de $a$ modulo $m$ est un multiple de $n$.

\art{2)} il existe $b$ tel que l'ordre de $b$ modulo $m$ est un multiple de $n$ et $a,b$ sont
indŽpendants modulo $m$.

}

{\bf Preuve}

On choisit successivement des entiers $r_1,\ldots , r_s, r'_1,\ldots ,r'_s$ tels que $r_i>k_i$ pour
$i=1,\ldots , s$ et des nombres premiers $p_1,\ldots ,p_s,p'_1,\ldots , p'_s$ tel que
$$q_1^{r_1}<p_1<q_2^{r_2}<p_2<\cdots <q^{r_s}<p_s<q_1^{r'_1}<p'_1<q_2^{r'_2}<\cdots <q_s^{r'_s}<p'_s$$
o, pour chaque $i$, $p_i$ (resp. $p'_i$) est choisit tel que l'ordre de $a$ modulo $p_i$ (
resp. $p'_i$) soit $q_i^{r_i}$ (resp. $q_i^{r'_i}$). C'est possible gr‰ce au Lemme \argct. Posons $m=p_1\cdots p_s\cdot p'_1\cdots p'_s$ ($p_1$ est arbitrairement grand, car on peut choisir $r_1$ arbitrairement grand). Clairement (thm. chinois), l'ordre de $a$ modulo $m$ est $q_1^{r'_1}\cdots q_s^{r'_s}$ qui est un multiple de $n$. Choisissons $b$ tel que 
$$\eqalign{b&\equiv a\pmod{p_1\cdots p_s}\cr
b&\equiv 1\pmod {p'_1\cdots p'_s}.\cr}$$
L'ordre de $b$ modulo $m$ est $q_1^{r_1}\cdots q_s^{r_s}$ qui est un multiple de $n$. VŽrifions encore que $a$ et $b$ sont indŽpendants~: supposons que $a^u\cdot b^v\equiv 1\pmod m$. Alors $1\equiv a^u\cdot b^v\equiv a^u\pmod{p'_1\cdots p'_s}$. Donc $q_1^{r'_1}\cdots q_s^{r'_s} | u$, car l'ordre de $a$ modulo $p'_1\cdots p'_s$ est $q_1^{r'_1}\cdots q_s^{r'_s}$. Donc, $a^u\equiv 1\pmod m$ et donc aussi $b^v\pmod 1\pmod m$.\qed
\bigskip
\prop
\medskip

{\sl Soit $L/K$ une extension abŽlienne de corps de nombres. Posons $n=[L:K]$ et $s$ un entier
supŽrieur ou Žgal ˆ 1. Soit $\P$ un idŽal premier de $K$ qui ne ramifie pas dans $L$. Alors il existe
$m\geq 1$ premier ˆ $s$ et ˆ $\P$ (i.e. $m\not\in \P$) tel que 
\medskip
\art{a)} si $E=K(\zeta_m)$, alors $\Phi_{E/K}(\P)$ a un ordre multiple de $n$.

\art{b)}$E\cap L=K$.

\art{c)}il existe $\tau\in {\rm Gal}(E/K)$ dont l'ordre est un multiple de $n$ et qui est indŽpendant  de
$\Phi_{E/K}(\P)$ dans ${\rm Gal}(E/K)$.

\art{d)}${\rm Gal}(E/K)\simeq {\rm Gal}(\Q(\zeta_{m})/\Q)$.

}

{\bf Preuve}

On applique le lemme prŽcŽdent ˆ $a=\N(\P)$ et $n=n$. On choisit  $M>1$ tel que $\Q(\zeta_M)\cap L$ soit
le plus grand sous-corps dans une extension cyclotomique de $\Q$. C'est possible, car si $E_1\subset
L\cap\Q(\zeta_{n_1})$ et $E_2\subset L\cap \Q(\zeta_{n_2})$, alors $L\cap\Q(\zeta_{{\rm ppcm}(n_1,n_2)})\supset E_1,E_2$ et on continue ainsi de suite, et a s'arrte forcŽment car le degrŽ de $L$ est fini. On choisit $m$ comme dans le lemme prŽcŽdent de sorte que $m$ soit premier ˆ $s$, $M$ et $\P$ (c'est possible car on a vu que le plus petit facteur premier de $m$ peut tre aussi grand qu'on veut). Alors, avec ce choix, on a $\Q\subset \Q(\zeta_m)\cap L\subset \Q(\zeta_m)\cap\Q(\zeta_M)=\Q$. Ainsi, $\Q(\zeta_m)\cap L=\Q$ et {\it a fortiori} $\Q(\zeta_m)\cap K=\Q$ aussi. La thŽorie de Galois nous dit alors que  ${\rm Gal}(L(\zeta_m)/L)\simeq{\rm Gal}(K(\zeta_m)/K)\simeq{\rm Gal}(\Q(\zeta_m)/\Q)\simeq(\Z/m\Z)^*$.  On a donc 
$$\eqalign{[L(\zeta_m):K(\zeta_m)]\cdot\varphi(m)&=[L(\zeta_m):K(\zeta_m)]\cdot [K(\zeta_m):K]=[L(\zeta_m):K]\cr &=[L(\zeta_m):L]\cdot [L:K]=\varphi(m)\cdot n.\cr}.$$
 Donc, $[L(\zeta_m):K(\zeta_m)]=n$. D'autre part, la thŽorie de Galois donne $[L:L\cap
K(\zeta_m)]=[L(\zeta_m):K(\zeta_m)]=n=[L:K]$, ce qui prouve que $L\cap K(\zeta_m)=K$, donc les parties $b)$ et $d)$ sont prouvŽes.

Par dŽfinition, $\Phi_{E/K}(\P)$ est l'automorphisme de $E/K$ caractŽrisŽ par $\zeta_m\mapsto
\zeta_m^{\N(\P)}=\zeta^a$ dont l'ordre est celui de $a$ modulo $m$, qui est un multiple de $n$ (par le
lemme prŽcŽdent). On a donc prouvŽ la partie~a).

ConsidŽrons le $b$ du lemme prŽcŽdent. On regarde l'automorphisme $\tau$ de $E/K$ dŽfinit par
$\tau(\zeta_m)=\zeta_m^b$. Il est clair que $\tau$ remplit les conditions de la partie c).\qed
\bigskip

\th{\soustitre (Lemme d'Artin)}
\medskip
{\sl Soit $L/K$ une extension cyclique de corps de nombres. $s>0$ un entier et $\P$ un idŽal premier de
$K$ non ramifiŽ dans $L$. Alors il existe $m>0$ un entier premier ˆ $s$ et ˆ $\P$ et une extension $F/K$
telle que 
\medskip
\art{a)}$L\cap F=K$

\art{b)}$L\cap K(\zeta_m)=K$

\art{c)}$L(\zeta_m)=F(\zeta_m)$

\art{d)}$\P$ est compltement dŽcomposŽ dans $F$.

}

{\bf Preuve}

Soit $m$ et $\tau$ comme dans la proposition prŽcŽdente. On choisit $\sigma$ un gŽnŽrateur de ${\rm
Gal}(L/K)=G$. Alors, par la thŽorie de Galois, $L(\zeta_m)/K$ est une extension abŽlienne de groupe
$G\times {\rm Gal}(K(\zeta_m)/K)$. Soit $F$ le corps fixe par le sous-groupe $H$ engendrŽ par les ŽlŽments
$$(\sigma,\tau)\hbox{ et }\Phi_{L(\zeta_m)/K}(\P)=(\Phi_{L/K}(\P),\Phi_{K(\zeta_m)/K}(\P)).$$

\art{a)} $L\cap F$ est le sous-corps de $L$ laissŽ fixe par les restrictions ˆ $L$ des ŽlŽments de $H$,
donc en particulier par $(\sigma,\tau)|_L=\sigma$. Donc $L\cap F$ est le sous-corps de $L$ laissŽ fixe
par $G$, et c'est Žvidemment $K$.

\art{b)} est donnŽ par la proposition prŽcŽdente.

\art{c)}On a $F(\zeta_m)=F\cdot K(\zeta_m)$ est le sous-corps de $L(\zeta_m)$ fixe par $H\cap (G\times
\{1\})$. Supposons que  $(\sigma,\tau)^u\cdot (\Phi_{L/K}(\P),\Phi_{K(\zeta_m)/K}(\P))^v\in G\times \{1\}$.
Comme (cf. proposition prŽcŽdente) $\tau$ et $\Phi_{K(\zeta_m)/K}(\P)$ sont indŽpendants, on a $\tau^u$ et $\Phi_{K(\zeta_m)/K}(\P))^v=1$. Donc, $u$ et $v$ sont des multiples de $n$ (toujours gr‰ce ˆ la proposition prŽcŽdente). Ainsi, $\sigma^u=1$ et $(\Phi_{L/K}(\P)^v=1$, ce qui implique que $H\cap (G\times \{1\})=\{1\}\times\{1\}$. Donc $F(\zeta_m)$ est le sous-corps de $L(\zeta_m)$ fixe par $\{1\}\times\{1\}$ qui est $L(\zeta_m)$.

\art{d)}Si on montre que le groupe de dŽcomposition de $\P$ dans $F$ est trivial, alors $\P$ se dŽcompose totalement dans $F$ et c'est (par dŽfinition) le sous-groupe engendrŽ par $\Phi_{F/K}(\P)$. Or, $\Phi_{F/K}(\P)=\Phi_{L(\zeta_m)/K}(\P)|_F=Id_F$ (par dŽfinition de $H$.) Cela montre partie d).\qed
\bigskip
En corollaire ˆ ce rŽsultat et ˆ sa preuve, on peut prouver le~:
\bigskip\goodbreak
\coro
\medskip

{\sl Soit $L/K$ une extension cyclique de corps de nombres $\P_1,\ldots ,\P_r$ des idŽaux premiers de
$K$, non ramifiŽs dans $L$, et pour chaque $i$, $F_i$ et $m_i$ comme dans le Lemme d'Artin associŽs ˆ $\P_i$. Alors, on peut choisir les $F_i$ et les $m_i$ de telle manire que si $F=F_1\cdots F_r$, alors on a $F\cap L=K$.

}

{\bf Preuve}

On construit $F_1, m_1,\tau_1$ comme pour le Lemme d'Artin et la
Proposition \argcw, avec $s=1$. Ensuite, supposons avoir construit
$F_{i-1},m_{i-1}$ et $\tau_{i-1}$, on construit $F_i, m_i$ et $\tau_i$ de telle manire que $s=m_1\cdots
m_{i-1}$. On rappelle que $\tau_i\in {\rm Gal}(K(\zeta_{m_i})/K)={\rm Gal}(\Q(\zeta_{m_1})/\Q):=G_i$ et
$\Phi_{K(\zeta_{m_i})/K}(\P_i)$ sont indŽpendants et d'ordre un multiple de $n=[L:K]$ et que $F_i$ est le
sous-corps de $L(\zeta_{m_1})$ laissŽ fixe par
$H_i$, le groupe engendrŽ par $(\sigma,\tau_i)$ et $\Phi_{L(\zeta_{m_i})/K}(\P_i)$, o $\sigma$ est un
gŽnŽrateur de $G={\rm Gal}(L/K)$. Mais, puisque les
$m_i$ sont premiers au $M$ de la proposition \argcw, on a vu que $L\cap K(\zeta_{m_i})=K$. On montre de la
mme manire que $L\cap K(\zeta_{m_1\cdots m_r})=K$ (car, puisque les $m_i$ sont premiers entre eux, on
a $\zeta_{m_1}\cdots\zeta_{m_s}=\zeta_{m_1\cdots m_r}$ et $m_1\cdots m_r$ est premier avec $M$). Ainsi, 
$L(\zeta_{m_1\cdots m_r})/K$ est une extension abŽlienne de groupe $G\times\prod_{j=1}^r G_j$, et dans
cette extension, pour tout $i$, $F_i$ est le corps fixe pour $\widetilde{H_i}:=H_i\times
\prod_{j\ne i}G_i$. Finalement, $F=F_1\cdots F_r$ est le corps fixe de $\cap_{i=1}^r\widetilde{H_i}$, et
donc $F\cap L$ est le sous-corps de $L$ fixŽ par les ŽlŽments de $\cap_{i=1}^r\widetilde{H_i}$ et ce
dernier groupe contient $(\sigma,\tau_1,\ldots ,\tau_s)$ dont la restriction ˆ $L$ est $\sigma$ qui
engendre $G$. Donc, $F\cap L=K$.\qed
\bigskip

Voici un gros thŽorme qui est le thŽorme de rŽciprocitŽ d'Artin pour les extensions cycliques.

\bigskip

\th
\medskip

{\sl Soit $L/K$ une extension cyclique de corps de nombres et $\m$ un $K$-module divisible par toutes les
places (finies ou infinies) qui ramifient dans $L$ et tel que 
$$[L:K]=[I_K(\m):P_\m\cdot N_{L/K}(I_L(\widetilde{\m})]$$
(ce qui est rŽalisŽ nous l'avons vu (ThŽorme \argck, ŽgalitŽ fondamentale) lorsque les exposants des
places finies qui divisent $\m$ sont suffisamment grands). Alors $\m$ est $K$-admissible.

}

{\bf Preuve}

Puisque $\Phi_{L/K}|I_K(\m)$ est surjective (cf. Proposition \argaq ), on a Žvidemment l'ŽgalitŽ
$[I_K(\m):\ker(\Phi_{L/K}|I_K(\m))]=[L:K]$, de sorte qu'il suffit de montrer que
$\ker(\Phi_{L/K}|I_K(\m))\subset P_\m\cdot N_{L/K}(I_L(\widetilde{\m})$, puisqu'ils ont le mme indice
dans $I_K(\m)$. Comme d'habitude, on pose $n=[L:K]$ et $\sigma$ un gŽnŽrateur de $G={\rm Gal}(L/K)$. Soit ${\euf a}\in I_K(\m)$ tel que $\Phi_{L/K}({\euf a})=1$. Ecrivons ${\euf
a}=\prod_{i=1}^r\P_i^{a_i}$, $a_i\in \Z$. On applique le Corollaire \argcy\ aux idŽaux $\P_1,\ldots ,\P_r$,
fournissant donc les $F_1,\ldots , F_r,m_1,\ldots , m_r$ et $F=F_1\cdots F_r$. Posons
$$\Phi_{L/K}(\P_i^{a_i})=\sigma^{d_i}\hbox{ avec }d_i\geq 0.$$ Le fait que $\Phi_{L/K}({\euf a})=1$ donne
que $\sigma^d=1$, o $d=d_1+\cdots + d_r$. Donc, $n|d$. Si $\m'$ est un $F$-module divisible par les idŽaux
premiers de $F$ qui ramifient dans $LF$. Alors, puisque $F\cap L=K$, l'application $\Phi_{LF/F}\ :
I_F(\m')\to {\rm Gal}(LF/F)\buildrel R\over\simeq G$ est une surjection, o $R$ est la restriction ˆ $L$.
On choisit un tel $\m'$ qui en plus est divisible par tous les $m_iO_F$ et par $\widetilde{\m}_F$ (l'extension
ˆ $F$ du $K$-module $\m$). Alors il existe ${\euf b}_0\in I_F(\m')\subset I_F(\widetilde{\m}_F)$ tel que
$R(\Phi_{LF/F}({\euf b}_0))=\sigma$. On a donc ${\euf b}:=N_{F/K}({\euf b}_0)\in I_K(\m)$ (mme dans
$I_K(\m\cdot (m_i))\subset I_K(\m)$). Par la relation $R\circ \Phi_{LF/F}=\Phi_{L/K}\circ N_{F/K}$ (cf.
ThŽorme \argf), on tire 
$$\Phi_{L/K}({\euf b})=\sigma.$$ 
L'idŽal ${\euf b}$ est la norme d'un idŽal dans $F_i$ (car $N_{F/K}=N_{F_i/K}\circ N_{F/F_i}$) et $\P_i$ aussi, car il est totalement dŽcomposŽ dans $F_i$. Ainsi, par multiplicativitŽ, il existe un idŽal fractionnaire ${\euf c}_i$ dans $F_i$ tel que $N_{F_i/K}({\euf c_i})=\P_i^{a_i}\cdot {\euf b}^{-d_i}$. En outre, notant $\widetilde{\m}_i$ l'extension de $\m$ ˆ $F_i$, on a que ${\euf c}_i\in I_{F_i}(\widetilde{\m}_i\cdot (m_i))$ (car ${\euf b}_0\in
I_F(\m')$ et $\P_i$ est premiers ˆ $m_i$ (cf. Proposition \argcw) et ˆ $\m$, car ${\euf a}\in I_K(\m)$). On a
donc, en vertu du ThŽorme \argf\ et de ce qui prŽcde~:
$$R_i\circ\Phi_{L F_i/F_i}({\euf c}_i)=\Phi_{L/K}(N_{F_i/K}({\euf c_i}))=\Phi_{L/K}(\P_i^{a_i})\cdot
\Phi_{L/K}({\euf b})^{-d_i}=\sigma^{d_i}\cdot\sigma^{-d_i}=1,$$
o $R_i$ est la restriction ${\rm Gal}(LF_i/F_i)\to{\rm Gal}(L/K)$. Comme cette restriction est un
application injective, on a $\Phi_{LF_i/F_i}({\euf c}_i)=1$. D'autre part, puisque
$L(\zeta_{m_i})=F_i(\zeta_{m_i})$ (cf. Lemme d'Artin, partie c)), on a $F_i\subset LF_i\subset
F_i(\zeta_{m_i})$. Donc, en vertu de la Proposition \argcs, tout $F_i$-module $\m_i''$ est admissible pourvu
que $\m_i''$ soit divisible par l'extension ˆ $F_i$ du $\Q$-module $(m_i)\cdot\infty$. Donc, en
choisissant $\m_i''$ le ppcm de $\widetilde{\m_i}$ et de l'extension ˆ $F_i$ de $(m_i)\cdot\infty$, on
a ${\euf c}_i\in I_{F_i}(\m_i'')$ et puisque $\m_i''$ est admissible et que ${\euf c}_i\in \ker(\Phi_{LF_i/F_i})$, alors ${\euf c}_i\in
P_{\m_i''}\cdot N_{LF_i/F_i}(I_{LF_i}(\widetilde{\m_i''}))$, o $\widetilde{\m_i''}$ est l'extension ˆ $LF_i$
de $\m_i''$. C'est-ˆ-dire
$${\euf c_i}=(\gamma_i)\cdot N_{LF_i/F_i}({\euf d}_i),\eqno{(*)}$$ 
o $\gamma_i\in {F_i}^*_{\m_i''}$ et ${\euf d}_i\in I_{LF_i}(\widetilde{\m_i''})$. Appliquant
$N_{F_i/K}$ de chaque c™tŽ de l'ŽgalitŽ $(*)$, on obtient
$$\P_i^{a_i}\cdot{\euf b}^{-d_i}=(N_{F_i/K}(\gamma_i))\cdot N_{LF_i/K}({\euf d}_i)=(\alpha_i)\cdot
N_{LF_i/K}({\euf d}_i)$$
o $\alpha_i=N_{F_i/K}(\gamma_i)\in N_{F_i/K}({F_i}^*_{\m_i''})\subset
N_{F_i/K}({F_i}^*_{\widetilde{\m_i}})\subset K^*_\m$ (la dernire inclusion vient de la partie b) du Lemme \aargbu). Ce qui veut dire que $(\alpha_i)\in P_\m$. En faisant le produit sur tous les $i$, on trouve~: 
$${\euf a}={\euf b}^d\cdot\prod_{i}(\alpha_i)\cdot\prod_{i}N_{LF_i/K}({\euf d}_i).$$
En posant ${\euf d}_i'=N_{LF_i/L}({\euf d}_i)\in I_L(\widetilde{\m})$, on trouve ${\euf
a}=\prod_{i}(\alpha_i)\cdot {\euf b}^d\cdot N_{L/K}(\prod_i{\euf d}_i')$. Finalement, puisque $n|d$ et que
$\euf b$ est un idŽal de $K$, premier ˆ $\m$, on a ${\euf b}^d=N_{L/K}({\euf b}^{d\over n}\cdot O_L)$ avec
${\euf b}\in I_K(\m)$ et ${\euf b}^{d\over n}\cdot O_L\in I_L(\widetilde{\m})$. En rŽsumŽ~:
$${\euf a}=\underbrace{\left (\prod_i\alpha_i\right )}_{\in P_\m}\cdot \underbrace{N_{L/K}(\prod_i{\euf
d}_i'\cdot {\euf b}^{d\over n}\cdot O_L)}_{\in N_{L/K}(I_L(\widetilde{\m}))}.$$\qed

\bigskip
Aprs ce gros morceau, on peut enfin prouver et Žnoncer le ThŽorme de rŽciprocitŽ d'Artin pour les
extensions abŽlienne. En rŽsumŽ, on s'y est pris d'abord avec les extensions cylotomiques (assez facile),
puis avec les extensions cycliques (trs difficile, on a besoin de l'ŽgalitŽ fondamentale et ˆ peu prs
tout ce qui prŽcde) et maintenant, on termine aisŽment avec les extensions abŽliennes. 
\bigskip

\th{\soustitre (ThŽorme de rŽciprocitŽ d'Artin)}
\medskip
{\sl Soit $L/K$ une extension abŽlienne de corps de nombres et $\m$ un $K$-module divisible par toutes les
places (finies et infinies) ramifiŽes dans $L$. Si les exposants des idŽaux premiers divisant $\m_0$
sont suffisamment grands, alors $\m$ est admissible pour $L/K$. 

Cela implique, en vertu du Lemme \argco, l'application $\Phi_{L/K}|_{I_K(\m)}$ est un
homomorphisme surjectif de noyau $P_\m\cdot N_{L/K}(I_L(\widetilde{\m}))$.

}

{\bf Preuve}

Posons $G={\rm Gal}(L/K)$ et $G=C_1\times\cdots\times C_s$, avec $C_i$ cycliques. Posons encore
$H_i=\prod_{j\ne i} C_j\times 1$ et $E_i$ le corps fixe par $H_i$. La thŽorie de Galois nous dit que
$E_i/K$ est une extension cyclique de groupe de Galois isomorphe ˆ $C_i$. Puisque $\m$ contient toutes les
places qui ramifient dans $L$, il contient toutes les places qui ramifient dans $E_i$, et, puisque les
exposants des idŽaux premiers divisant $\m_0$ sont suffisamment grands, alors en vertu du thŽorme
prŽcŽdent, $\m$ est admissible pour chaque extension $E_i/K$, ce qui veut dire que
$P_\m\subset\cap_{i=1}^s\ker(\Phi_{E_i/K}|I_K(\m))$. Mais, on montre facilement, puisque $L=E_1\cdots E_s$
que pour tout ŽlŽment ${\euf a}\in I_K(\m)$, on a $\Phi_{L/K}({\euf a})=(\Phi_{E_1/K}({\euf a}),\ldots
,\Phi_{E_s/K}({\euf a}))\in C_1\times\cdots\times C_s=G$. Donc,
$\ker(\Phi_{L/K}| I_K(\m))=\cap_{i=1}^s\ker(\Phi_{E_i/K}|I_K(\m))$. Cela prouve que $P_\m\subset
\ker(\Phi_{L/K}| I_K(\m))$ ce qui prouve le thŽorme.\qed

\bigskip

\coro
\medskip

{\sl Le ThŽorme de $\check{\rm C}$ebotarev fort est vrai

}

{\bf Preuve}

Le thŽorme prŽcŽdent est exactement le Lemme \argbc\ dont nous avions besoin pour le ThŽorme de
$\check{\rm C}$ebotarev fort.\qed

\bigskip
\coro
\medskip

{\sl Soit $L/K$ une extension abŽlienne et $E/K$ une extension galoisienne et $\m$ un $K$-module
admissible pour $L/K$. Supposons que
$$N_{E/K}(I_E(\widetilde{\m}'))\subset N_{L/K}(I_L(\widetilde{\m}))\cdot P_\m$$
o $\widetilde{\m}$, (resp. $\widetilde{\m}'$) est l'extension de $\m$ ˆ $L$ (resp. $E$). Alors $L\subset
E$.

}

{\bf Preuve}

En dehors d'un nombre fini d'idŽaux (ceux qui pourraient diviser $\m$), chaque idŽal premier de $K$
compltement dŽcomposŽ dans $E$ est contenu dans $N_{E/K}(I_E(\widetilde{\m}'))$ (c'est la norme de
n'importe quel idŽal au-dessus de lui). Donc, par hypothse, il est contenu dans
$N_{L/K}(I_L(\widetilde{\m}))\cdot P_\m=\ker(\Phi_{L/K}(I_K(\m)))$, donc il est totalement dŽcomposŽ dans
$L$. On a donc que tout idŽal premier de $K$ compltement dŽcomposŽ dans $E$ est totalement dŽcomposŽ dans
$L$ (sauf Žventuellement un nombre fini). Cela prouve que $L\subset E$, en vertu du ThŽorme \argar.\qed

\bigskip\goodbreak

{\bf Application}

Soit $L/\Q$ une extension abŽlienne de corps de nombres. Par la rŽciprocitŽ d'Artin, il existe
$\m=(m)\cdot\infty$ un $\Q$-module admissible pour $L/\Q$. On peut supposer $m>0$ entier. On sait, par
ailleurs qu'il est admissible pour $E/\Q$, avec $E=\Q(\zeta_m)$ et mme que $\ker(\Phi_{E/\Q}|I_\Q(\m))=P_\m$
(cf. ThŽorme~\argm). Or, la rŽciprocitŽ d'Artin nous montre que ce noyau est $P_\m\cdot
N_{E/\Q}(I_E(\widetilde{\m}'))$. Cela prouve que $N_{E/\Q}(I_E(\widetilde{\m}'))\subset P_\m\subset
P_\m\cdot N_{L/\Q}(I_L(\widetilde{m}))$, donc, en vertu du corollaire prŽcŽdent, $L\subset E$.
\bigskip
On a donc prouvŽ le

\bigskip

\th {\soustitre (ThŽorme de Kronecker-Weber)}
\medskip
{\sl Toute extension abŽlienne de $\Q$ est sous-extension d'une extension cyclotomique $\Q(\zeta_m)$.
}
\qed

\vfill\eject
\global\advance\chapnomb by 1
\nomb=1

\centerline{\para Chapitre 8 :}
\medskip
\centerline{\para PrŽparation ˆ la formation des classes }
\bigskip
Nous allons, dans ce chapitre, exposer les notions de sous-groupe de congruence et de conducteur, regarder quelques notions de bases sur eux. Puis nous allons Žnoncer thŽorme central de la thŽorie du corps de classe et faire quelques rŽductions.
\bigskip
\defi
\medskip
 Soit $K$ un corps de nombres. Un sous-groupe $H$ de $I_K$ est dit un {\it sous-groupe de congruence}
s'il existe un $K$-module $\m$ tel que $P_\m\subset H\subset I_K(\m)$. On dit alors que $H$ est un
{\it sous-groupe de congruence pour $\m$}.

\newcount\gaagad\gaagad=\gaga

\bigskip
Tout d'abord quelques remarques immŽdiates 
\bigskip
\headline={\hfill \phantom{ouuh}\hfill}
\rems

\medskip
\art{a)}Si $H$ est un sous-groupe de congruence pour $\m$, et $\m'$ est un $K$-module  multiple de $\m$
ayant les mme places finies, alors $H$ est un sous-groupe de congruence pour $\m'$. En effet,
dans ce cas, $P_{\m'}\subset P_\m$ et $I_K(\m)=I_K(\m)$.

\art{b)}Soit $H$ un sous-groupe de congruence pour $\m$ et $\m'$ est un autre $K$-module multiple de $\m$
(le produit de deux $K$-modules se fait par produit au niveau des idŽaux et par union au niveau des
places infinies, cf. Chapitre 0 pour plus de dŽtails). Alors $H\cap I_K(\m')$ est un sous-groupe de congruence pour $\m'$. En effet, il est
clair que $$\matrix{I_K(\m')&\subset &I_K(\m)\cr \cup&&\cup\cr P_{\m'}&\subset&P_\m\cr}.$$
Donc $P_{\m'}\subset P_\m\cap I_K(\m')\subset H\cap I_K(\m')\subset I_K(\m')$.

\bigskip

\defi
\medskip

On dira que deux sous-groupes de congruence $H$ et $H'$ sont {\it Žquivalents} s'il existe un $K$-module
$\m''$ tel que 
$$H\cap I_K(\m'')=H'\cap I_K(\m'').$$
Il est clair que si $\m'''$ est un multiple de $\m''$, alors on a aussi $H\cap I_K(\m''')=H'\cap
I_K(\m''')$ (puisque $I_K(\m''')\subset I_K(\m'')$). Ainsi, on a bien affaire ˆ une relation
d'Žquivalence. De plus, si $H$ (resp. $H'$) sont dŽfinis pour $\m$ (resp. $\m'$), alors on peut toujours choisir
pour $\m''$ un multiple commun de $\m$ et de $\m'$ (par exemple le ppcm de $\m,\m'$ et $\m''$), ainsi
$H\cap I_K(\m'')$ sera un sous-groupe de congruence pour $\m''$ .

\bigskip

\lem
\medskip

{\sl Si $\m_1$ et $\m_2$ sont des $K$-modules tels que $\m_1 |\m_2$ et pour $i=1,2$, $H_i$ est un
sous-groupe de congruence pour $\m_i$ tel que $H_2=H_1\cap I_K(\m_2)$, alors on a

\art{a)}$I_K(\m_2)/H_2\simeq I_K(\m_1)/H_1$

\art{b)}$H_1=H_2\cdot P_{\m_1}$.

En particulier la partie b) implique que si $H_2$ provient de $H_1$ et que $\m_1$ est fixŽ, alors $H_1$
est unique.

}
\goodbreak
{\bf Preuve}

On montre d'abord que $$I_K(\m_1)=I_K(\m_2)\cdot P_{\m_1}.\eqno{(*)}$$ L'inclusion
$\supset$ est triviale. RŽciproquement, soit $\aa_1\in I_K(\m_1)$. On peut l'Žcrire
$\aa_1=\aa\cdot \aa_2$ avec $\aa_2\in I_K(\m_2)$ et $\aa=\prod_i\P_i^{a_i}$, o les
$\P_i$ sont les idŽaux premiers qui divisent $\m_2$, mais pas $\m_1$. Par le thŽorme
d'approximation dŽbile (ThŽorme \argc), il existe, pour tout $i$, $\pi_i\in K$ tel que
$\pi_i\in\P_i\setminus\P_i^2$ et $\pi_i\equiv^* 1\pmod {\m_1}$. Alors $\alpha=\prod_i\pi_i^{a_i}\in
K^*_{\m_1}$ et $\aa (\alpha)^{-1}$ est premier ˆ $\m_2$ (par unicitŽ de la dŽcomposition d'un idŽal en
puissances d'idŽaux premiers). Ainsi,
$$\aa_1=\underbrace{\aa (\alpha)^{-1}}_{\in I_K(\m_2)}\cdot \underbrace{\aa_2}_{\in
I_K(\m_2)}\cdot\underbrace{(\alpha)}_{\in P_{\m_1}}\in I_K(\m_2)\cdot P_{\m_1}.$$
On a {\it a fortiori} $I_K(\m_1)=I_K(\m_2)\cdot H_1$ (car $P_{\m_1}\subset H_1\subset I_K(\m_1)$). Alors,
$$I_K(\m_2)/H_2=I_K(\m_2)/(H_1\cap I_K(\m_{2}))\buildrel{\rm thm.\ d'isom.}\over\simeq (I_K(\m_2)\cdot
H_1)/H_1=I_K(\m_1)/H_1,$$
prouvant a).

Pour la partie b), il est dŽjˆ clair que $H_2\cdot P_{\m_1}\subset H_1$. On va montrer que ces deux groupes
ont le mme indice dans $I_K(\m_1)$, ce qui permettra de conclure. Cette premire inclusion implique que 
$$H_2\subset (P_{\m_1}\cdot H_2)\cap I_K(\m_2)\subset H_1\cap I_K(\m_2)=H_2.$$
On en dŽduit que $H_2= (P_{\m_1}\cdot H_2)\cap I_K(\m_2)$. D'o,
$$\eqalign{I_K(\m_1)/H_1&\buildrel a)\over \simeq I_K(\m_2)/H_2=I_K(\m_2)/((P_{\m_1}\cdot H_2)\cap
I_K(\m_2))\buildrel{\rm thm.\ d'isom.}\over\simeq (\underbrace{I_K(\m_2)\cdot H_2}_{=\, I_K(\m_2)}\cdot P_{\m_1})/(H_2\cdot P_{\m_1})\cr&\buildrel
(*)\over = I_K(\m_1)/(H_2 P_{\m_1}).\cr}.$$ 
On a donc bien $[I_K(\m_1): H_2\cdot P_{\m_1}]=[I_K(\m_1): H_1]$, et c'est ce qu'il fallait pour montrer
la partie b).\qed
\bigskip
\headline={\hfill\smcap PrŽparation ˆ la formation des classes\hfill}
\lem
\medskip

{\sl Soit $H_1$ et $H_2$ des sous-groupes de congruence dŽfinis pour $\m_1$, respectivement $\m_2$,
Žquivalents. Alors si $\m={\rm pgcd}(\m_1,\m_2)$ (par convention, ${\rm pgcd}(\m,\m')={\rm
pgcd}(\m_0,\m'_0)\cdot (\m_\infty\cap\m'_\infty)$, pour tout $K$-modules $\m$ et $\m'$, cf. Chapitre 0 pour plus de dŽtails), il existe un sous-groupe de congruence $H$ dŽfinit pour
$\m$ tel que $H\cap I(\m_i)=H_i$. 

}

{\bf Preuve}

Soit $\m_3$ un $K$-module multiple de $\m_1$ et de $\m_2$ tel que $H_3:=H_1\cap I_K(\m_3)=H_2\cap I_K(\m_3)$,
c'est possible (cf. DŽfinition \argdf). On pose
$$H=H_3\cdot P_\m.$$
On va voir que $H$ a toutes les bonnes propriŽtŽs. Tout d'abord, $P_\m\subset H\subset I_K(\m)$, car
$H_3\subset I_K(\m_3)\subset I(_K(\m)$. On va montrer que 
$$H\cap I_K(\m_3)=H_3.\eqno{(*)}$$
Cela suffit pour montrer le lemme. En effet, considŽrons les groupes $H_1$ et $H\cap I_K(\m_1)$. Il sont les
deux des groupes de congruence pour $\m_1$ (pour $H_1$, c'est l'hypothse et pour $H\cap I_K(\m_1)$, on a
$H=H_3\cdot P_\m\supset P_{\m_1}$). Par $(*)$ et puisque $I_K(\m_3)\subset I_K(\m_1)$, on a $I_K(\m_3)\cap H_1=I_K(\m_3)\cap(H\cap I_K(\m_1))=H_3$. En vertu du Lemme \argdg, on a donc
$H_1=H_3\cdot P_{\m_1}=H\cap I_K(\m_1)$. Et on montre de la mme manire que $H_2=H\cap I_K(\m_2)$.

Reste ˆ prouver la relation $(*)$. Le fait que $H_3\subset H\cap I_K(\m_3)$ est trivial. RŽciproquement, soit
$\aa\cdot (\alpha)\in H\cap I_K(\m_3)$, avec $\aa\in H_3$ et $\alpha\in K^*_\m$. Puisque $\aa$ et $\aa\cdot
(\alpha)\in I_K(\m_3)$, $\aa^{-1}$ aussi et donc $(\alpha)$ aussi. Nous cherchons $\beta\in K^*$ tel que 
$$\beta\in K^*_{\m_1}\hbox{ et } \alpha\beta^{-1}\in K^*_{\m_2}\hbox{ et }(\beta)\in I_K(\m_3).\eqno{(**)}$$
Supposons l'existence d'un tel $\beta$. D'une part, $\aa\cdot(\beta)\in H_3\cdot P_{\m_1}\buildrel{\rm
Lemme\ \argdg}\over =H_1$. D'autre part, $\aa\cdot (\beta)\in I_K(\m_3)$. On en dŽduit donc que $\aa\cdot
(\beta)\in H_1\cap I_K(\m_3)=H_3$. Puis, $\aa\cdot (\alpha)=\aa\cdot (\beta)\cdot (\alpha\beta^{-1})\in H_3
P_{\m_2}\buildrel {\rm Lemme\ \argdg}\over =H_2$. Et enfin,
$$\aa\cdot (\alpha)\in H_2\cap I_K(\m_3)=H_3.$$

Montrons l'existence de $\beta$ satisfaisant $(**)$. Supposons que la place $\P|\m_3$ mais ne divise pas
$\m_1$ ni $\m_2$, alors on demande que $\beta\equiv 1\pmodast\P$. Si $\P|\m_1$ ou $\m_2$, posons $\P^{a_i}$
l'exacte puissance de $\P$ qui divise $\m_i$, $i=1,2$. Si $a_1>a_2$, on demande que $\beta\equiv
1\pmodast{\P^{a_1}}$. Dans ce cas-lˆ, $\P^{a_2}$ est l'exacte puissance de $\P$ qui divise $\m={\rm
pgcd}(\m_1,\m_2)$. Donc on a $\beta\equiv 1\pmodast{\P^{a_2}}$. Mais puisque $\alpha\equiv 1\pmodast
{\P^{a_2}}$, on a $\alpha\beta^{-1}\equiv 1\pmodast{\P^{a_2}}$ (si $a_2=0$, on n'a pas besoin de cette Žquivalence). Enfin, si $a_1\leq a_2$, on demande que $\beta\equiv\alpha\pmodast {\P^{a_2}}$. Dans ce cas,
$\P^{a_1}$ est l'exacte puissance de $\P$ qui divise $\m$, donc $\alpha\equiv 1\pmodast{\P^{a_1}}$, et comme
{\it a fortiori} $\beta\equiv \alpha\pmodast{\P^{a_1}}$, on a $\beta\equiv 1\pmodast{\P^{a_1}}$. Un tel $\beta$
existe en vertu du thŽorme d'approximation dŽbile (ThŽorme \argc) et il satisfait clairement les
conditions $(**)$ ce qui montre notre lemme.\qed

\bigskip

{\soustitre Corollaire-DŽfinitions ({\the\chapnomb}.{\the\nomb})}\global\advance\nomb by 1
\medskip
{\sl Soit $\bbH$ une classe d'Žquivalence de sous-groupes de congruence. Alors il existe un $K$-module $\ff$
appelŽ le {\it conducteur de } $\bbH$ qui a la propriŽtŽ suivante~: pour chaque multiple $\m$ de $\ff$, $\bbH$
contient un unique sous-groupe de congruence, notŽ $H(\m)$, dŽfini pour $\m$, et

$$\bbH=\{H(\m)\mid \ff|\m\}\hbox{ et en particulier }H(\m)=H(\ff)\cap I_K(\m) \hbox{ pour tout $\m$
multiple de $\ff$}.$$

D'autre part, si $\bbH$ et $\bbH'$ sont des classes d'Žquivalence de sous-groupes de congruence, on dira que
$\bbH\subset \bbH'$ s'il existe un $K$-module $\m$ tel que $H(\m)\subset H'(\m)$. On a
 alors
$$\bbH\subset \bbH'\iff \ff' |\ff$$
et dans ce cas, $H(\m)\subset H'(\m)$ pour tout $\ff |\m$. }
\goodbreak

\newcount\gaagae\gaagae=\gaga
{\bf Preuve}

C'est un corollaire facile des Lemmes \argdg\ et \argdh~: si $H_1$ et $H_2\in \bbH$ sont deux groupes dŽfinis pour $\m$,
le Lemme \argdh, appliquŽ ˆ $\m=\m_1=\m_2$, nous montre que $H_1=H_2$. D'autre part, posons $\ff$ le pgcd de tous
les $K$-modules pour lesquels un sous-groupe de congruence est dans $\bbH$. Par le Lemme \argdh, il existe un
unique groupe de congruence $H(\ff)$ dŽfinit pour $\ff$ tel que pour tout $H(\m)\in \bbH$ dŽfinit pour $\m$, on
ait $H(\m)=H(\ff)\cap I_K(\m)$. Enfin, pour tout multiple $\m$ de $\ff$, on a en vertu de la remarque qui
prŽcde la DŽfinition \argdf, que  $H(\ff)\cap I_K(\m)$ est un sous-groupe de congruence pour $\m$ Žquivalent ˆ
$H(\ff)$.

Supposons maintenant que $\bbH\subset \bbH'$, donc l'existence d'un $\m$ tel que
$H(\m)\subset H'(\m)$. On sait gr‰ce ˆ la premire partie que $H(\m)=H(\ff)\cap
I_K(\m)$ et gr‰ce au Lemme \argdg\ que $I_K(\ff)/H(\ff)\simeq I_K(\m)/H(\m)$ et cet
isomorphisme est donc induit par les inclusions. On est donc dans la situation~:
$$\matrix{P_\ff&\subset &H(\ff)&\subset &?&\subset&I_K(\ff)\cr 
\cup&&\cup&&&&\cup\cr
P_\m&\subset &H(\m)&\subset&
H'(\m)&\subset& I_K(\m)\cr}$$
Il existe donc un unique sous-groupe $V$  tel que $H(\ff)\subset V\subset I_K(\ff)$
tel que $H'(\m)=V\cap I_K(\m)$. Donc, $V$ est un sous-groupe de congruence pour $\ff$ Žquivalent ˆ $H'(\m)$, ainsi $V=H'(\ff)\in \bbH'$ ce qui prouve que
$\ff'|\ff$ et que $H(\ff)\subset H(\ff')$. Finalement,  $H(\m)=H(\ff)\cap I_K(\m)\subset H(\ff')\cap I_K(\m)= H'(\m)$ pour tout $\ff |\m$.
\qed

\bigskip

{\soustitre Application-DŽfinition ({\the\chapnomb}.{\the\nomb})}\global\advance\nomb by 1

Soit $L/K$ une extension abŽlienne de corps de nombres de groupe de Galois $G$. Soit $\m$ un $K$-module
admissible (cf. DŽfinition \argcn). Posons $H(\m)$ le noyau de l'application d'Artin $\Phi_{L/K}\;
:\; I_K(\m)\to G$. C'est un sous-groupe de congruence pour $\m$. Il est clair que pour tout autre $K$-module
admissible $\m'$, les sous-groupes $H(\m)$ et $H(\m')$ sont Žquivalents; donc tous ces groupes dŽfinissent  une
unique classe d'Žquivalence de sous-groupes de congruence qu'on notera $${\bbH}(L/K).$$
Le $K$-module qui est le conducteur de cette classe d'Žquivalence s'appellera {\it le conducteur de} $L/K$, qu'on notera $\ff(L/K)$.
\bigskip
Remarquons que si $\m$ est un $K$-module avec $\ff(L/K)|\m$, cela n'implique pas forcŽment que $\m$ est
admissible, mme si le groupe $H(\m)$ existe, mais n'est pas forcŽment le noyau de l'application d'Artin
restreinte ˆ $I_K(\m)$. NŽanmoins on a le rŽsultat~:
\newcount\gaagaf\gaagaf=\gaga
\bigskip

\lem
\medskip

{\sl Soit $L/K$ une extension abŽlienne de corps de nombres. ConsidŽrons ${\bbH}(L/K)=\{H(\m)\mid \ff
|\m\}$, o $\ff=\ff(L/K)$. Supposons que  $\m$ soit un $K$-module divisible par toutes les places (finies
ou infinies) qui ramifient dans $L$ et que $\ff|\m$. Alors $\m$ est admissible pour $L/K$.

}
{\bf Preuve}

Le thŽorme de rŽciprocitŽ d'Artin (cf. ThŽorme \argda) nous assure l'existence d'un $K$-module $\n$ 
admissible, qui est un multiple de $\m$ et divisible par les mme place que $\m$, ainsi,
$I_K(\m)=I_K(\n)$, et donc, en vertu du Corollaire-DŽfinitions \argdi, $H(\m)=I_K(\m)\cap H(\ff)=I_K(\n)\cap H(\ff)=H(\n)$. Pour la mme raison que $I_K(\m)=I_K(\n)$, on a aussi $\ker(\Phi_{L/K}|I_K(\n))=\ker(\Phi_{L/K}|I_K(\m))$.
Puisque $\n$ est admissible, alors $H(\n)=\ker(\Phi_{L/K}|I_K(\n))$, et donc,
$H(\m)=\ker(\Phi_{L/K}|I_K(\m))$. Or, puisque $H(\m)$ est un sous-groupe de congruence pour $\m$, on a
$P_\m\subset H(\m)=\ker(\Phi_{L/K}|I_K(\m))$, ce qui veut dire que $\m$ est admissible.\qed

\bigskip

Voici un petit lemme facile qui sera utile dans bien longtemps (au ThŽorme \argek ) mais qui peut sans autre tre ŽnoncŽ ici
\medskip
\lem
\medskip

{\sl Soit $L/K$ une extension abŽlienne de corps de nombres. Soit  $\ff=\ff(L/K)$ le conducteur de $L/K$. Soit $\P\in\gfP(K)$ une place. Supposons que $\P | \ff$, alors $\P$ ramifie dans $L$

}

{\bf Preuve}

C'est aussi un corollaire du thŽorme de rŽciprocitŽ d'Artin (cf. ThŽorme \argda) : supposons que $\P |\ff$ et que $\P$ ne ramifie pas dans $L$. Par le dit thŽorme de rŽciprocitŽ, il existe un $K$-module admissible $\m$ tel que $\P\notdiv \m$. Mais puisqu'il est admissible, $\ff |\m$ et donc $\P|\m$, ce qui est une contradiction. On a donc montrŽ que $\P$ doit ramifier dans $L$.\qed  
\medskip

Voici un rŽsultat important~: 
\medskip\goodbreak
\prop
\medskip

{\sl Soit $L/K$ une extension abŽlienne de corps de nombres. Alors la correspondance 
$$L/K\mapsto \bbH(L/K)$$
est injective. Mieux~: si $L/K$ et $L'/K$ sont des extensions abŽliennes
contenues dans une mme cl™ture algŽbrique de $K$, on a
$$L\subset L'\iff \bbH(L/K)\supset \bbH(L'/K).$$ 

}

{\bf Preuve}

Si $L\subset L'$, et si $\m$ est un $K$-module admissible pour $L'/K$, alors 
$\m$ est aussi admissible pour $L/K$ (Lemme \argcq ) et on a clairement que
$H'(\m)\subset H(\m)$, car $\Phi_{L/K}=\Phi_{L'/K}|_{L}$, ce qui montre que
$\bbH(L'/K)\subset \bbH(L/K)$. RŽciproquement, supposons que $\bbH(L'/K)\subset
\bbH(L/K)$. Puisque les idŽaux premiers qui se dŽcomposent totalement dans $L'$ sont
dans le noyau de $\Phi_{L'/K}$, alors (sauf  Žventuellement pour un nombre fini~:
ceux qui diviseraient le conducteur de $L'/K$) ces idŽaux dŽcomposent compltement
dans $L$. Cela prouve, gr‰ce au ThŽorme \argar, que $L\subset
L'$.\qed
\bigskip

Nous pouvons alors Žnoncer un des thŽormes centraux de la thŽorie du corps de
classe qui dit essentiellement que la correspondance qui prŽcde est surjective~:
\bigskip
\th {\soustitre  (ThŽorme d'existence du corps de classe)}

\medskip

{\sl Soit $K$ un corps de nombres. Alors pour toute
classe $\bbH$ d'Žquivalence de sous-groupes de congruence de $K$, il existe une
extension abŽlienne $L/K$ telle que $\bbH=\bbH(L/K)$. On dit alors que $L/K$ est {\it
le corps de classe de $\bbH$} (on devrait plut™t dire le corps de la classe $\bbH$).

}
\bigskip
On ne va pas prouver ce thŽorme en un clin d'oeil. On va d'abord faire des
rŽductions pour la preuve finale.
\bigskip
\goodbreak
\defi

Si $\bbH=\{H(\m)\mid \ff |\m\}$ est une classe d'Žquivalence de sous-groupes de
congruence. On a vu au Lemme \argdh\ que les groupes $I_K(\m)/H(\m)$ sont isomorphes
canoniquement. On notera ce groupe $I_K/\bbH$.

\newcount\gaagag\gaagag=\gaga
\bigskip
\goodbreak
\th
\medskip

{\sl Soit $K$ un corps de nombres. Soit $\bbH_0$ et
$\bbH_1$ des classes de sous-groupes de congruence de $K$. Supposons que
$\bbH_0\subset \bbH_1$ et que $\bbH_0$ possde un corps de classe $L$. Alors $\bbH_1$
en possde aussi un (c'est une sous-extension de $L$).

}

{\bf Preuve}

Supposons donc que $\bbH_0=\bbH(L/K)$ et posons, pour $i=1,2$, $\ff_i$ le conducteur
de $\bbH_i$. Si $\m$ est un $K$-module admissible pour $L/K$. Alors, en vertu du
ÒCorollaire-DŽfinitions \argdi", on a $\ff_1|\ff_0 |\m$ et $P_\m\subset H_0(\m)\subset
H_1(\m)\subset I_K(\m)$ o $H_i(\m)$ est le sous-groupe de congruence de $\bbH_i$
dŽfini pour $\m$. Soit $G={\rm Gal}(L/K)$ et posons $G_1=\Phi_{L/K}(H_1(\m))$ et
$E$ le corps fixe pour $G_1$. L'application d'Artin $\Phi_{E/K}$ est dŽfinie sur
$I_K(\m)$ (car $\m$ contient tous les idŽaux qui ramifient dans $L$, donc dans $E$)
et $\Phi_{E/K}=R\circ\Phi_{L/K}$ o $R$ est la restriction $G\to{\rm
Gal}(E/K)=G/G_1$. Si $\aa\in H_1(\m)$, on a, par dŽfinition de $G_1$,
$\Phi_{E/K}(\aa)=1\in G/G_1$, donc $H_1(\m)\subset \ker(\Phi_{E/K})$. Mais on a
$$\eqalign{[I_K(\m):\ker(\Phi_{E/K})]&=|G/G_1|=[G:G_1]\cr &\buildrel \Phi_{L/K}\rm \,
surjective\over =[\Phi_{L/K}(I_K(\m)):\Phi_{L/K}(H_1(\m))]\cr &\buildrel \rm Lemme\ \argcj\over ={[I_K(\m):H_1(\m)]\over [H_0(\m):\underbrace{H_0(\m)\cap
H_1(\m)]}_{=H_0(\m)}}=[I_K(\m):H_1(\m)].\cr}$$
Donc, $H_1(\m)=\ker(\Phi_{E/K})$ (restreint ˆ $I_K(\m)$) et donc
$\bbH_1=\bbH(E/K)$.\qed

\bigskip

\defi
\medskip
Soit $E/K$ une extension de corps de nombres et $\bbH=\{H(\m)\mid\ff|\m\}$ une classe
de sous-groupes de congruence de $K$. Si $H(\m)\in \bbH$, on pose 
$$H_E(\widetilde{\m})=\{\aa\in I_E(\widetilde{\m})\mid N_{E/K}(\aa)\in H(\m)\}.$$
Le lemme suivant montrera que tous les $H_E(\widetilde{\m})$ sont des sous-groupes de
congruence pour $\widetilde{\m}$ et sont tous Žquivalents. On pose
$\bbH_E$ la classe de sous-groupes de congruence engendrŽs par les
$H_E(\widetilde{\m})$. Le conducteur de cette classe est donc un diviseur de
$\widetilde{\ff}$ o $\ff$ est le conducteur de $\bbH$. On vŽrifie aussi que si
$K\subset E\subset F$ est une Òtour" de corps de nombres, alors on a
$(\bbH_E)_F=\bbH_F$. 
\bigskip
\newcount\gaagah\gaagah=\gaga
\lem
\medskip

{\sl Sous les mmes hypothses que pour la dŽfinition prŽcŽdente, alors tous les $H_E(\widetilde{\m})$ sont des sous-groupes de
congruence pour $\widetilde{\m}$ et sont tous Žquivalents.

}

{\bf Preuve}

Puisque $N_{E/K}(E^*_{\widetilde{\m}})\subset K^*_\m$ (partie b) du Lemme \aargbu), on en dŽduit alors que $N_{E/K}(P_{\widetilde{\m}})\subset P_\m$
et donc, de $P_\m\subset H(\m)\subset I_K(\m)$, suit $P_{\widetilde{\m}}\subset
H_E(\widetilde{\m})\subset I_E(\widetilde{\m})$. Donc $H_E(\widetilde{\m})$ est bien
un sous-groupe de congruence pour $\widetilde{\m}$.

Supposons que $\m|\m'$ (par exemple on peut prendre $\m=\ff$). On a donc
$H(\m')=H(\m)\cap I_K(\m')$. Alors on a   
$$H_E(\widetilde{\m'})=\{\aa\in I_E(\widetilde{\m'})\mid N_{E/K}(\aa)\in
H(\m')=H(\m)\cap I_K(\m')\}\subset I_E(\widetilde{\m'})=I_E(\widetilde{\m'})\cap
H_E(\widetilde{\m}).$$ 
En effet~: puisque $\m |\m'$, on a $I_K(\m)\supset I_K(\m')$
et donc $I_E(\widetilde{\m})\supset I_E(\widetilde{\m'})$, il est clair que
$H_E(\widetilde{\m'})\subset I_E(\widetilde{\m'})\cap H_E(\widetilde{\m})$.
Inversement, si
$\aa\in I_E(\widetilde{\m'})\cap H_E(\widetilde{\m})$, alors $N_{E/K}(\aa)\in
H(\m)\cap I_K(\m')$, car
$N_{E/K}(I_E(\widetilde{\m'}))\subset I_K(\m')$ par dŽfinition et $N_{E/K}(\aa)\in
H(\m)$ par hypothse. Cela prouve notre lemme.\qed
\bigskip
\goodbreak
\th
\medskip
{\sl Soit $E/K$ une extension cyclique de corps de nombres et $\bbH$ une classe de
sous-groupes de congruence pour $K$. Supposons que $\bbH_E$ ait un corps de classe
$L/E$, alors $\bbH$ a aussi un corps de classe.

}

{\bf Preuve}

Si $\ff$ est le conducteur de $\bbH$, choisissons un $K$-module $\m$ tel que
$\ff|\m$ et tel que $\widetilde{\m}$ soit admissible pour $L/E$ (c'est possible en
vertu du thŽorme de rŽciprocitŽ d'Artin (ThŽorme \argda)). Alors 
$$ H_E(\widetilde{\m})=\ker(\Phi_{L/E}\, : I_E(\widetilde{\m})\longrightarrow {\rm
Gal}(L/E)).$$

Nous allons montrer que l'extension $L/K$ est galoisienne (mme abŽlienne) et que
$\bbH(L/K)\subset \bbH$ ce qui achevera la preuve en vertu du ThŽorme \argdo.

\art{a)}$L/K$ est galoisienne. En effet, soit $\sigma\, : L\to \C$ un
$K$-morphisme (on peut voir $L$ comme un sous-corps de $\C$). Il suffit de montrer
que $\sigma(L)=L$. D'abord, $\sigma(E)=E$, puisque $E/K$ est galoisienne. On a  que $\sigma(L)/E$
est le corps de classe de $\sigma(\bbH_E)$ ( on applique $\sigma$ aux idŽaux des
groupes de $\bbH_E$). Alors, puisque $\m$ est un $K$-module, $\widetilde{\m}$ et
invariant par $\sigma$, donc $I_E(\widetilde{\m})$ l'est aussi, donc~:
$$\eqalign{\sigma(H_E(\widetilde{\m}))&=\{ \sigma(\aa)\mid \aa\in
I_E(\widetilde{\m})\hbox{ et }N_{E/K}(\aa)\in H(\m)\}\cr &=\{\aa\in
I_E(\widetilde{\m})\mid N_{E/K}(\sigma|_E^{-1}(\aa))\in H(\m)\}\cr&=\{\aa\in
I_E(\widetilde{\m})\mid N_{E/K}(\aa)\in H(\m)\})=H_E(\widetilde{\m}),\cr}$$
L'avant dernire ŽgalitŽ vient du fait que $\sigma|_E$ est un $K$-automorphisme de
$E$, donc un ŽlŽment du groupe de galois de $E/K$ qui ne change pas la norme d'un
idŽal. Cela prouve que $\sigma(\bbH_E)=\bbH_E$, donc que $\sigma(L)=L$ par
l'injectivitŽ de la correspondance vue ˆ la Proposition \1argdm.

\art{b)}$L/K$ est abŽlienne. En effet, choisissons $\sigma\in {\rm Gal}(L/K)$ tel
que $\sigma|_E$ engendre ${\rm Gal}(E/K)$ (qui est cyclique par hypothse). Ainsi,
tout ŽlŽment $\mu\in {\rm Gal}(L/K)$ est tel que $\mu|_E=\sigma^a|_E$ pour un
certain $a\in \Z$. Ainsi, $\mu\cdot\sigma^{-a}=\tau\in {\rm  Gal}(L/E)$. Ce qui
veut dire que $\mu=\sigma^a\cdot\tau$. Donc, pour montrer que ${\rm Gal}(L/K)$ est
abŽlien, il suffit de montrer que $\sigma\cdot\tau=\tau\cdot\sigma$ pour tout
$\tau\in {\rm Gal}(L/E)$. Soit donc un tel $\tau$. Puisque l'application d'Artin est
surjective (cf. ThŽorme \argaq), il existe $\aa\in I_E(\widetilde{\m})$ tel que $\Phi_{L/E}(\aa)=\tau$.
On a donc
$\sigma\tau\sigma^{-1}=\sigma\Phi_{L/E}(\aa)\sigma^{-1}=\Phi_{L/E}(\sigma(\aa))$.
Or, on a $N_{E/K}({\aa\over\sigma(\aa)})=O_K=1_{H(\m)}\in H(\m)$. Cela prouve que
${\aa\over\sigma(\aa)}\in H_E(\m)=\ker(\Phi_{L/E})$ par hypothse. Et donc, on a
$\Phi_{L/E}(\sigma(\aa))=\Phi_{L/E}(\aa)=\tau$. On a alors prouvŽ que
$\sigma\cdot\tau\cdot\sigma^{-1}=\tau$ donc, $\sigma\cdot\tau=\tau\cdot\sigma$. Ce
qui montre que $L/K$ est abŽlienne. Elle est jolie cette preuve, vous ne trouvez
pas ?

Terminons la preuve. Le raisonnement que nous venons de faire est valable pour tout
$K$-module multiple de $\m$ (a ne change pas la classe). Donc on peut supposer en
plus que $\m$ soit admissible pour $L/K$. C'est-ˆ-dire que $P_\m\subset
\ker\left (\Phi_{L/K}\, :\, I_K(\m)\to {\rm Gal}(L/K)\right)\subset I_K(\m)$.
D'autre part, puisque $\bbH$ est une classe de sous-groupe de congruence, on a
aussi $P_\m\subset H(\m)\subset I_K(\m)$. Comme $\widetilde{\m}$ est admissible
pour $L/E$ et que $H_E(\widetilde{\m})=\ker(\Phi_{L/E}\, :
I_E(\widetilde{\m})\to {\rm Gal}(L/E))$, on a en vertu du corollaire du ThŽorme
\argg\ que $N_{L/E}(I_L(\widetilde{\widetilde{\m}}))\subset
H_E(\widetilde{\m})$, o $\widetilde{\widetilde{\m}}$ est le $L$-module au-dessus
de $\m$ et de $\widetilde{\m}$. Ce qui veut dire par dŽfinition et par propriŽtŽ
de la norme que $N_{L/K}(I_L(\widetilde{\widetilde{\m}}))\subset H(\m)$. Ainsi, on
a~:
$$P_\m\subset P_\m\cdot N_{L/K}(I_L(\widetilde{\widetilde{\m}}))\subset
H(\m)\subset I_K(\m).$$
Puisque $\m$ est admissible pour $L/K$, le Lemme \argco\ montre que 
$P_\m\cdot N_{L/K}(I_L(\widetilde{\widetilde{\m}}))=\ker \left(\Phi_{L/K}|_{
I_K(\m)}\right ).$ On a donc prouvŽ que $\bbH(L/K)\subset \bbH$ et
donc, en vertu du ThŽorme \argdo, que $\bbH$ a aussi un corps de classe.\qed

\bigskip

\th{\soustitre (rŽduction)}

\medskip
{\sl Le thŽorme d'existence (notre ThŽorme \argdm) est vrai si on montre le rŽsultat
suivant~:

Si $\bbH$ est une classe de sous-groupes de congruence de $K$ et s'il existe un
entier $n$ tel que

\art{a)}$K$ contient une racine primitive $n$-ime de l'unitŽ 

\art{b)}$n$ est un exposant de $I_K/\bbH$ (rappel~: l'exposant d'un groupe est le
ppcm des ordres de ses ŽlŽments, par extension, un exposant est un entier $n$ tel
que $x^n$ est l'ŽlŽment unitŽ pour tout $x$ dans le groupe),

alors $\bbH$ admet un corps de classe.

}

{\bf Preuve}

Supposons donc le rŽsultat dŽmontrŽ et soit $\bbH$ une classe de sous-groupes de
congruence d'un corps de classe $K$. Posons $n$ l'exposant de $I_K/\bbH$, soit
$\beta$ un racine primitive $n$-ime de l'unitŽ et $E:=K(\beta)$. Montrons que $n$
est un exposant pour $I_E/\bbH_E$. On choisit des corps $K_i$ tels que
$$K=K_1\subset K_2\subset\cdots\subset K_t=K(\beta)=E,$$
et tels que $K_i/K_{i-1}$ soit cyclique, pour tout $i=2,\ldots ,t$ (partant de $E$,
on prend le corps fixe par le sous-groupe engendrŽ par un ŽlŽment, qui donne
$K_{t-1}$, etc..., jusqu'ˆ obtenir $K$). On dŽfinit alors les classes suivantes~:
$\bbH_1=\bbH, \bbH_2,\ldots ,\bbH_t=\bbH_E$, $\bbH_i$ Žtant la classe de sous-groupes
de congruence de $K_i$ dŽfinies pour tout $i>1$ par
$$\bbH_i=(\bbH_{i-1})_{K_i}=\bbH_{K_i}.$$
La dernire ŽgalitŽ vient de la remarque vue dans la dŽfinition de $\bbH_E$.
Montrons que $I_{K_i}/\bbH_i$ est d'exposant $n$ pour chaque $i$. C'est vrai par
hypothse si $i=1$. Supposons donc par rŽcurrence que c'est vrai pour $i\geq 1$. 
Soit $H_i(\m)\in \bbH_i$ le sous-groupe de congruence correspondant au $K_i$-module
$\m$. Supposons que le conducteur de $\bbH_{i+1}$ divise $\widetilde{\m}$. Alors on
a
$H_{i+1}(\widetilde{\m})=\{\aa\in I_{K_{i+1}}(\widetilde{\m})\mid
N_{K_{i+1}/K_i}(\aa)\in H_i(\m)\}\in \bbH_{i+1}$. Soit $\aa\in
I_{K_{i+1}}(\widetilde{\m})$. Alors,
$N_{K_{i+1}/K_i}(\aa^n)=(N_{K_{i+1}/K_i}(\aa))^n\in I_{K_i}(\m)^n\buildrel\rm h.r.\over \subset H_i(\m)$.
Donc, $\aa^n\in H_{i+1}(\widetilde{\m})$. Cela prouve que $n$ est un exposant de
$I_{K_{i+1}}/\bbH_{i+1}\simeq I_{K_{i+1}}(\widetilde{\m})/H_{i+1}(\widetilde{\m})$.  Ainsi, $n$ est un exposant de $I_E/\bbH_E$. Nous sommes donc bien dans l'hypothse du rŽsultat acceptŽ. Donc, en cascade, chaque $\bbH_i$ (et en particulier $\bbH$) admet un corps de classe, en vertu du ThŽorme \argdr.\qed

\bigskip

\vfill\eject
\global\advance\chapnomb by 1
\nomb=1

\centerline{\para Chapitre 9 :}
\medskip
\centerline{\para Quelques rŽsultats sur la thŽorie des $\taille{16}n$-extensions de Kummer, }
\medskip
\centerline{\para et calcul d'un nouvel indice}

\bigskip
Interrompons  un instant le propos pour donner quelques rŽsultats que nous utiliserons par la suite comme lemmes.
\bigskip
\defi
\medskip

Une extension $L/K$ de corps est dite {\it $n$-extension de Kummer} si $K$  contient une racine primitive
$n$-ime de l'unitŽ, si $L/K$ est abŽlienne et si $n$ est un exposant de ${\rm Gal}(L/K)$. Il est Žvident que si $K$ 
contient une racine primitive $n$-ime de l'unitŽ, elle contient toutes les autres racines $n$-ime de l'unitŽ (d'une mme
cl™ture algŽbrique de $K$).

\bigskip

\newcount\gaagai\gaagai=\gaga

\th
\medskip
{\sl Si $L/K$ est une $n$-extension de Kummer de degrŽ fini et $M=M(L/K)=\{\alpha\in L^*\mid \alpha^n\in K^*\}$, alors $M^n$ est un
sous-groupe de $K^*$ contenant ${K^*}^n$, $M^n/{K^*}^n$ est fini et $L=K(\{\root n\of {x}\mid x\in M^n\})$. De plus, l'application $L/K\mapsto M(L/K)^n$ est une bijection de l'ensemble des $n$-extensions de Kummer finies de $K$ et l'ensemble des
sous-groupes $W$ de $K^*$ tels que $W\supset {K^*}^n$ et $W/{K^*}^n$ est fini. En outre, on a 

$$W/{K^*}^n\buildrel\rm canonique\over\simeq \widehat{{\rm Gal}(L/K)}\buildrel\rm non\ canonique\over\simeq {\rm
Gal}(L/K),$$
o $\widehat{{\rm Gal}(L/K)}$ est le groupe de caractre de ${\rm Gal}(L/K)$.

}

\headline={\hfill \phantom{ouuh}\hfill}
{\bf Preuve}

Soit donc $L/K$ une $n$-extension de Kummer et $M=\{\alpha\in L^*\mid \alpha^n\in K^*\}$. Il est clair que $M$ est un
sous-groupe de $L^*$ contenant ${K^*}^n$. Posons $G={\rm Gal}(L/K)$ et $\mu_n$ le groupe de racines $n$-ime de l'unitŽ.
Pour chaque $\alpha\in M$, on dŽfinit
$$\eqalign{\psi(\alpha)\, :\, G&\longrightarrow \mu_n\cr \sigma&\longmapsto
\psi(\alpha)(\sigma):={\sigma(\alpha)\over\alpha}.\cr}$$ 
C'est une application bien dŽfinie, puisque $\alpha\in M$~: $\left
({\sigma(\alpha)\over\alpha}\right )^n={\sigma(\alpha^n)\over\alpha^n}={\alpha^n\over\alpha^n}=1$. On vŽrifie facilement
que $\psi(\alpha)$ est un homomorphisme (un caractre de $G$) pour tout $\alpha$. En effet,
$$\psi(\alpha)(\sigma\tau)={\sigma(\tau(\alpha))\over\alpha}=
{\sigma(\tau(\alpha))\over\tau(\alpha)}\cdot{\tau(\alpha)\over\alpha}\buildrel (*)\over={\sigma(\alpha)\over\alpha}\cdot
{\tau(\alpha)\over\alpha}=\psi(\alpha)(\sigma)\cdot \psi(\alpha)(\tau).$$
L'ŽgalitŽ $(*)$ vient du fait que ${\sigma(\alpha)\over\alpha}$ est une racine $n$-ime de l'unitŽ donc appartient ˆ $K$
et donc, pour tout $\tau\in G$,
${\sigma(\alpha)\over\alpha}=\tau({\sigma(\alpha)\over\alpha})={\tau(\sigma(\alpha))\over\tau(\alpha)}
={\sigma(\tau(\alpha))\over\tau(\alpha)}$. Et on voit donc facilement que $\psi\, :\, M\to\widehat{G}$ est un homomorphisme de
groupe. Voyons son noyau~: $\alpha\in\ker(\psi)\iff{\sigma(\alpha)\over\alpha}=1$ pour tout $\sigma\in G\iff\alpha\in
K^*$. Ainsi, $\psi$ induit un homomorphisme injectif
$$\overline{\psi}\, :\, M/K^*\longrightarrow \widehat{G}$$
qui est parfois appelŽ ÒKummer pairing". L'ŽlŽvation ˆ la puissance $n$ et le passage au quotient induisent une
application surjective $M \to\hskip-7pt\to M^n\to\hskip-7pt\to M^n/{K^*}^n$. L'ŽlŽment $x\in M$ est dans le noyau si et
seulement s'il existe  $y\in K^*$ tel que $x^n=y^n$. Cela veut dire que ${x\over y}:=z$ est une racine $n$-ime de l'unitŽ,
donc $z\in K^*$, ce qui veut dire que $x=z\cdot y\in K^*$ et donc que $M/K^*\simeq M^n/{K^*}^n$. Montrons que $\psi$
(donc $\overline{\psi}$) est surjectif~: soit $\chi$ un caractre de $\widehat{G}/{\rm im}(\psi)$. On peut voir $\chi$
comme un caractre de $\widehat{G}$ trivial sur ${\rm im}(\psi)$. Dans tout ouvrage traitant des caractres [cf. Fr.-Tay., (A9). p. 331 ],
on peut voir que les caractres de $\widehat{G}$ sont les Žvaluations de $G$, i.e. il existe $\sigma\in G$ tel que
$\chi(\omega)=\omega(\sigma)$, pour tout $\omega\in \widehat{G}$. Dans notre cas, ce $\sigma$ a donc la propriŽtŽ que
$\omega(\sigma)=1$ pour tout  $\omega\in {\rm im}(\psi)$. C'est-ˆ-dire $\psi(\alpha)(\sigma)=1$ pour tout $\alpha\in M$;
ou encore $\sigma(\alpha)=\alpha$ pour tout $\alpha\in M$. Pour conclure que $\sigma=1$ (donc que $\chi=1$ et donc que
$\widehat{G}={\rm im}(\psi)$), il suffit de montrer que $M$ engendre $L$ sur $K$. Le groupe $G$ Žtant abŽlien, l'extension
$L/K$ est une composition de sous-extensions cycliques $L_1/K,\ldots ,L_k/K$ (voir le dŽbut de la dŽmonstration de la
rŽciprocitŽ d'Artin pour plus de dŽtail (ThŽorme \argda)) et il suffit donc de montrer que pour tout $i=1,\ldots , k$, $M$ contient un
gŽnŽrateur de $L_i/K$. On va le faire pour $L_1$ (les autres cas seront identiques). Soit $d=[L_1:K]$. On a $d|n$; soit
$\zeta$ une racine primitive $d$-ime de l'unitŽ (par hypothse, $\zeta\in K$). Soit $C={\rm Gal}(L_1/K)$ et $\sigma_1$
un gŽnŽrateur de $C$. Puisque, $\zeta\in K$, $N_{L_1/K}(\zeta)=\zeta^d=1$, par le ThŽorme de Hilbert 90 (cf. [La1, Thm.
6.6.1, p. 298]), il existe $\alpha\in L_1$ tel que $\sigma_1(\alpha)=\zeta\cdot\alpha$. En Žlevant ˆ la puissance $d$, on
obtient $\sigma_1(\alpha^d)=\alpha^d$, c'est-ˆ-dire que $a:=\alpha^d\in K$. A fortiori, $\alpha^n\in K$, ce qui veut dire
que $\alpha\in M$. ConsidŽrons le polyn™me $X^d-a$. On veut montrer que ce polyn™me est irrŽductible dans $K[X]$. Dans
$L_1[X]$, on a $X^d-a=\prod_{i=0}^{d-1}(X-\zeta^i\alpha)$. Supposons que $X^d-a$ soit divisible par un polyn™me $f\in
K[X]$, unitaire de degrŽ $l$, avec $0<l<d$. Alors le coefficient constant de ce polyn™me doit tre de la forme
$\pm \varepsilon\cdot\alpha^l$, o $\varepsilon$ est une racine $d$-ime de l'unitŽ; ce qui implique (puisque les racines
$d$-ime de l'unitŽ sont dans $K$) que $\alpha^l\in K$, cela est impossible puisque
$\sigma_1(\alpha^l)=\zeta^l\cdot\alpha^l$ et $\zeta^l\ne 1$. On a donc montrŽ que $X^d-a$ est irrŽductible dans $K[X]$,
ce qui montre que $\alpha$ engendre $L_1$ sur $K$. En dŽfinitive, on a montrŽ que 
$$M^n/{K^*}^n\simeq M/K^*\buildrel\overline{\psi}\over\simeq\widehat{G}\simeq G\eqno{(*)} $$
donc que $M^n/{K^*}^n$ est fini. On a aussi montrŽ que
$$L=K(\{x\mid x\in M\})=K(\{\root n\of x\mid x\in M^n\}),$$
la dernire ŽgalitŽ venant Žvidemment du fait que les racines de l'unitŽ sont dans $K$.

Enfin, soit $W$ un sous-groupe de $K^*$ contenant ${K^*}^n$ tel que $W/{K^*}^n$ soit fini. Posons $L=K(\{\root n\of x\mid
x\in W\})$. La finitude de $W/{K^*}^n$ fait que $L=K(\root n\of{\alpha_1},\ldots ,\root n\of{\alpha_m})$, o
$\alpha_1,\ldots ,\alpha_m$ engendrent $W$ modulo ${K^*}^n$ (on utilise toujours que $\mu_n$, l'ensemble des racines
$n$-ime de l'unitŽ, est dans $K$). Donc, $L/K$ est une extension finie. Clairement, $L/K$ est galoisienne, car les
conjuguŽs des racines $n$-ime des $\alpha_i$ sont dans $L$. Si $\sigma\in {\rm Gal}(L/K)$, on a, pour tout $i$,
$\sigma(\root n\of{\alpha_i})=\omega_{i,\sigma}\cdot \root n\of{\alpha_i}$, avec $\omega_{i,\sigma}\in \mu_n\subset K$.
Gr‰ce ˆ cela, on montre que trs facilement que ${\rm Gal}(L/K)$ est abŽlien et admet $n$ comme exposant. En bref, $L/K$
est une $n$-extension de Kummer.

Pour achever la preuve, il suffit de montrer que $W=M^n$, o $M=M(L/K)$. Trivialement, $W\subset M^n$. On sait dŽjˆ, gr‰ce
ˆ $(*)$ que $[M^n:{K^*}^n]=[L:K]$, il suffit donc de montrer que $[L:K]\leq [W : {K^*}^n]$. Supposons que l'ordre de
$\alpha_i$ modulo ${K^*}^n$ est $d_i$ et que les $\alpha_i$ correspondent ˆ une prŽsentation de $W/{K^*}^n$ comme produit
de groupes cycliques~:
$$W/{K^*}^n\simeq <\alpha_1\bmod{{K^*}^n}>\times\cdots\times <\alpha_m\bmod{{K^*}^n}>.$$
Alors, $[W : {K^*}^n]=d_1\cdots d_m$. Par dŽfinition de $d_i$, on a $\alpha_i^{d_i}=\beta_i^n$, avec $\beta_i\in K^*$,
donc on peut supposer (quitte ˆ remplacer $\beta_i$ par $\beta_i$ fois une racine $n$-ime de l'unitŽ) que $(\root
n\of{\alpha_i})^{d_i}=\beta_i$. Ainsi $\root n\of{\alpha_i}$ est une racine de $X^{d_i}-\beta_i$ et donc le degrŽ de $\root
n\of{\alpha_i}$ sur $K$ est $\leq d_i$. Enfin, on a
$$\eqalign{[L:K]=&[K(\root n\of{\alpha_1},\ldots ,\root n\of{\alpha_m}):K(\root n\of{\alpha_1},\ldots ,\root
n\of{\alpha_{m-1}})]\cdot[K(\root n\of{\alpha_1},\ldots ,\root n\of{\alpha_{m-1}}):K(\root n\of{\alpha_1},\ldots ,\root
n\of{\alpha_{m-2}})]\cr &\cdots [K(\root n\of {\alpha_1}):K]\leq [K(\root n\of {\alpha_m}):K]\cdots [K(\root n\of
{\alpha_1}):K]\leq d_m\cdots d_1=[W : {K^*}^n]\cr},$$
ce qui montre le thŽorme.\qed

\headline={\hfill \smcap $n$-extensions de Kummer  et calcul d'un nouvel indice \hfill}
\bigskip
Passons maintenant ˆ un tout autre rŽsultat qui gŽnŽralise quelque peu le thŽorme des unitŽs de Dirichlet.
\bigskip

\defi
\medskip

Soit $S$ un ensemble de places d'un corps de nombre $K$. On dŽfinit
$${K^*}^S=\{a\in K^*\mid v_\P(aO_K)\ne 0\Rightarrow \P\in S\}.$$
Il va sans dire que si on ajoute ou enlve une place infinie ˆ $S$, on de change pas l'ensemble ${K^*}^S$. On supposera donc que $S$ contienne toutes le place infinie (ainsi, l'ŽnoncŽ du thŽorme suivant sera plus ŽlŽgant). On pose 
$$I_K[S]\ \hbox{ le sous-groupe de $I_K$ engendrŽ par les idŽaux premiers $\P\in S$.}$$

\newcount\gaagaj\gaagaj=\gaga
\bigskip

\th {\bf (ThŽorme de Dirichlet-Chevalley-Hasse)}
\medskip
{\sl Soit $K$ un corps de nombres et $S$ un ensemble fini de places de $K$. On suppose que $S$ contienne toutes les places infinies. Alors ${K^*}^S$ est le produit direct d'un groupe cyclique fini (le groupe des racines de l'unitŽ de $K$) et d'un groupe abŽlien libre de rang $|S|-1$. On remarque que si $S$ est l'ensemble des places infinies, alors ${K^*}^S=U_K$ et qu'on retrouve le thŽorme classique des unitŽs de Dirichlet }

{\bf Preuve}

Posons $S_0$ l'ensemble des places finies de $S$. On a une suite exacte

$$1\lra  U_K\lra {K^*}^S\buildrel\iota\over \lra I_K[S_0]$$
o l'application $\iota$ est dŽfinie par $\iota(a)=a\cdot O_K$. Soit $h=h_K=|I_K/P_K|$. Il est clair que $I_K[S_0]^h\subset \iota({K^*}^S)\subset I_K[S_0]$ (cf. [Sam, \S 4.2, ThŽorme 2, p. 71 ). Ainsi, $\iota({K^*}^S)$ est un groupe abŽlien libre de mme rang que $I_K[S_0]$ qui est $|S_0|$. Puisque la suite exacte ci-dessus est formŽe de groupes abŽliens, elle est scindŽe, donc ${K^*}^S\simeq \iota({K^*}^S)\times U_K$. Or, le thŽorme classique de Dirichlet sur les unitŽs (cf. [Sam, ThŽorme 1, \S 4.4, p. 72]) nous apprend que $U_K$ est le produit direct du groupe des racines de l'unitŽs de $K$ et d'un groupe abŽlien libre de rang $|S\setminus S_0|-1$. Cela prouve le thŽorme.\qed

\bigskip\goodbreak

\coro
\medskip

{\sl Sous les mmes hypothses que le thŽorme prŽcŽdent, supposons de plus que $K$ contienne une racine primitive $n$-ime de l'unitŽ. Alors on a

$$[{K^*}^S:({K^*}^S)^n]=n^{|S|}.$$

}

{\bf Preuve}

Le thŽorme prŽcŽdent nous montre que ${K^*}^S\simeq<g>\times \Z^{|S|-1}$, o $g$ est une racine de l'unitŽ. Par hypothse, $g$ est d'ordre un multiple de $n$. Alors $({K^*}^S)^n\simeq <g^n>\times (n\Z)^{|S|-1}$. D'o, par propriŽtŽ de $g$,
$${K^*}^S/({K^*}^S)^n\simeq \Z/n\Z\times (\Z/n\Z)^{|S|-1}=(\Z/n\Z)^{|S|}.$$\qed

\bigskip\goodbreak

\centerline {\soustitre Un nouveau calcul d'indice }

\bigskip

\defi
\medskip

Supposons que $K$ contienne une racine primitive $n$-ime de l'unitŽ et soit $\m$  un $K$-module. On pose

$$c(\m)=[K^*:{K^*}^nK^*_\m].$$

\bigskip\goodbreak

\newcount\gaagak\gaagak=\gaga
\th
\medskip

{\sl Soit $K$ un corps de nombres contenant une racine primitive $n$-ime de l'unitŽ

\art{a)} Si $\m_1$ et $\m_2$ sont des $K$-modules premiers entre eux, alors,
$$c(\m_1\m_2)=c(\m_1)\cdot c(\m_2).$$

Il suffit donc de calculer $c(\m)$ lorsque $m=\P^t$ o $\P$ est une place finie ou infinie de $K$.

\art{b)}Si $\m=\P$ est une place infinie rŽelle, alors 
$$c(\m)=2.$$

\art{c)}Si $\m=\P^t$ o $t\in\N$ et $\P$ est un idŽal premier de $K$, alors, si $t$ est assez grand, on a~: 
$$c(\m)=n^2\cdot [O_{\P}:n\cdot O_{\P}]=[O_{(\P)}:n\cdot O_{(\P)}],$$

o $\okp$ est le complŽtŽ localisŽ de $O_K$ en $\P$ et $\okpp$ est le localisŽ de $O_K$ en $\P$.

}

{\bf Preuve}

Le thŽorme d'approximation dŽbile (ThŽorme \argc) nous dit que $K^*/K^*_{\m_1\m_2}\simeq K^*/K^*_{\m_1}\times K^*/K^*_{\m_2}$. Or, les puissances $n$-imes se correspondent, donc
$${K^*}^nK^*_{\m_1\m_2}/K^*_{\m_1\m_2}\simeq {K^*}^nK^*_{\m_1}/K^*_{\m_1}\times {K^*}^nK^*_{\m_2}/K^*_{\m_2}.$$
En quotientant, on trouve
$$K^*/({K^*}^nK^*_{\m_1\m_2})\simeq K^*/({K^*}^nK^*_{\m_1})\times K^*/({K^*}^nK^*_{\m_2}),$$
ce qui montre la partie a).

\art{b)}Si $\m=\P$ est une place infinie rŽelle, alors $n=2$ (puisque $K$ contient une racine $n$-ime de l'unitŽ et possde un plongement rŽel). Si $\sigma$ est le plongement associŽ ˆ $\P$, l'homomorphisme surjectif $K^*\to\{ \pm 1\}$, $x\mapsto {\rm sgn}(\sigma(x))$ a pour noyau $\{x\in K^*\mid\sigma(x)>0\}=K^*_\m={K^*}^2 K^*_\m$, ce qui montre la partie b).

\art{c)}Supposons que $\m=\P^t$ o $t\in\N$ et $\P$ est un idŽal premier de $K$. Coupons notre indice en deux~:
$$c(\m)=[K^*:{K^*}^nK^*(\m)]\cdot [{K^*}^nK^*(\m):{K^*}^nK^*_\m],$$
o $K^*(\m)$ est le groupe des unitŽs du localisŽ de $O_K$ en $\P$ (autrement dit : $K^*(\m)=U(\okpp)$ ). Calculons le premier terme~:  nous savons que $\okpp$ est un anneau de Dedekind ˆ quotient fini n'ayant qu'un seul idŽal premier, c'est donc un anneau de valuation discrte. Il existe donc une {\it uniformisante} $\pi$ c'est -ˆ-dire que  $\pi\okpp=\P\okpp$ et  ainsi, tout ŽlŽment de $K^*$ s'Žcrit de manire unique $\pi^k\cdot u$, avec $k\in\Z$ et $u\in K^*(\m)$. On a donc un isomorphisme de groupe $K^*\simeq \Z\times K^*(\m)$. Via cet isomorphisme, on a que ${K^*}^n\simeq n\Z\times K^*(\m)^n$ et $K^*(\m)\simeq \{0\}\times K^*(\m)$. Ce qui donne ${K^*}^nK^*(\m)\simeq n\Z\times K^*(\m)$ et donc $K^*/({K^*}^nK^*(\m))\simeq \Z/n\Z$, ce qui prouve que $$[K^*:{K^*}^nK^*(\m)]=n.$$ Calculons le second terme. La composition d'applications naturelles $K^*(\m)\to {K^*}^n K^*(\m)\to ({K^*}^nK^*(\m))/({K^*}^nK^*_\m)$ est surjective. Son noyau est $K^*(\m)\cap {K^*}^nK^*_\m= K^*(\m)^nK^*_\m$ (car $K^*_\m\subset K^*(\m)$). Donc,
$$[{K^*}^nK^*(\m):{K^*}^nK^*_\m]=[K^*(\m):K^*(\m)^nK^*_\m].$$

Notons $\widetilde{\P}=\P\okpp$. On se souvient que l'homomorphisme
surjectif $\okpp\to \okpp/\widetilde{\P}^t$ induit un homomorphisme $O^*_{K_{(\P)}}=K^*(\m)\to
(O_{(\P)}/\widetilde{\P}^t)^*$ (cf. preuve du Lemme \argbz). Puisque $\okpp$ est local, alors $1+\widetilde{\P}^t\subset K^*(\m)$,
donc cet homomorphisme est surjectif et son noyau est $1+\widetilde{\P}^t=K^*_\m$ . On a aussi $(\okpp/\widetilde{\P}^t)^*\simeq
(\okp/\widehat{\P}^t)^*=U_\P/(1+\widehat{\P}^t)=U_\P/U_\P^{(t)}$, o  $U_\P=O^*_{\P}$, $U_\P^{(t)}=1+\widehat{\P}^t$ et l'avant dernire ŽgalitŽ vient aussi du fait que $\okp$ est local. Ainsi $K^*(\m)/K^*_\m\simeq U_\P/U_\P^{(t)}$. On a montrŽ au passage que $K^*(\m)/K^*_\m$ Žtait fini. Comme avant, passant au quotient par les puissances $n$-imes, on obtient~:
$$K^*(\m)/(K^*(\m)^n\cdot K^*_\m)\simeq U_\P/(U^n_\P\cdot U^{(t)}_\P).$$
Or, si $t$ est assez grand, on a prouvŽ ˆ la Proposition \argcc\ que $U_\P^{(t)}\subset U_\P^n$. Ainsi, si $t$ est assez grand, on a
$$c(\m)=n\cdot [U_\P:U_\P^n].\eqno{(*)}$$
Il nous reste donc ˆ calculer l'indice $[U_\P:U_\P^n]$. Pour cela, nous allons utiliser la $q$-machine  ($ q$ pour Òquotient de Herbrand")~: considŽrons $G$ un groupe cyclique d'ordre $n$. Tous les $G$-modules $M$ considŽrŽs seront considŽrŽs comme triviaux ($\sigma\cdot x=x$ pour tout $x\in M$ et $\sigma\in G$).  Soit donc $M$ un tel module, notŽ multiplicativement. Dans ce cas, $H^0(M)=\ker\Delta/{\rm Im}N=M/M^n$ et $H^1(M)=\{x\in M\mid x^n=1\}/\{1\}$. Ainsi, $|H^0(U_\P)|=[U_\P:U_P^n]$ et $|H^1(U_P)=n$, car par hypothse $K$ contient les racines $n$-ime de 1. Cela implique alors gr‰ce ˆ $(*)$ que 
$$c(\m)={n^2\over q(U_\P)}.\eqno{(**)}$$
Puisque, pour tout $s$ entier, on a $U_\P/U_\P^{(s)}\simeq U(\okp/\widehat{\P}^s)$, on a $[U_\P: U_\P^{(s)}]<\infty$ ainsi, en vertu du Corollaire~\argbl , on a $q(U_\P)=q(U_\P^{(s)})$. Le Corollaire \argcb\ nous montre que si $s$ est assez grand, alors l'application $\log$ est un isomorphisme de $U_\P^{(s)}$ sur $\widehat {\P}^s$. Et puisque $[\okp: \widehat {\P}^s]<\infty$, le  Corollaire \argbl nous montre que 
$$q(U_\P)=q( \widehat {\P}^s)=q(\okp),$$
o la structure de $ \widehat {\P}^s$ et de $\okp$ est maintenant additive, mais toujours triviale. On peut calculer $q(\okp)$ directement~: $H^0(\okp)=\okp/n\okp$ et $H^1(\okp)=\{0\}$. Ainsi, $q(\okp)={1\over [\okp:n\okp]}$. Finalement, la relation $(**)$ nous donne
$$c(\m)=n^2\cdot [\okp:n\okp]=n^2\cdot [\okpp:n\okpp]$$
la dernire ŽgalitŽ venant du fait que $\okp/n\okp\simeq \okpp/n\okpp$ [cf. Fr-Tay Th. 11 + Cor, p.77].\qed
\bigskip
\centerline{\soustitre Interlude}
\bigskip

Nous allons prouver en passant un rŽsultat bien connu. Il  s'agit de la rŽciprocitŽ quadratique. Ceci explique pourquoi on appelle le ThŽorme \argda\ le ThŽorme de {\it rŽciprocitŽ} d'Artin. Evidemment, ce n'est pas la preuve la plus directe de ce rŽsultat...
\bigskip

\defi

Soit $n\in \N$ et $K$ un corps de nombres contenant une racine primitive $n$-ime de l'unitŽ $\zeta_n$. Soit $\P\in \gfP_0(K)$ premier ˆ $n$,  et
$\alpha\in O_K\setminus \P$. Il est clair que $\alpha^{\N(\P)-1}\equiv 1\pmod \P$. Donc
l'image de $\alpha$ dans $\F:=O_K/\P$ notŽe $\overline{\alpha}$ est telle que
$\overline{\alpha}^{\N(\P)-1}=1$. Ainsi, $\overline{\alpha}^{\N(\P)-1\over n}$ est une racine $n$-ime de
l'unitŽ dans $\F$. Par le Lemme~\argo, il existe une unique racine $n$-ime de l'unitŽ dans $K$,
notŽe $\left ({\alpha\over \P}\right )_n$ telle que 

$$\left ({\alpha\over \P}\right )_n\equiv \alpha^{\N(\P)-1\over m}\pmod\P.$$

\newcount\gaagal\gaagal=\gaga

On Žtend  cette application (appelŽe {\it symbole de puissance $n$-ime rŽsiduelle}) multiplicativement au niveau des dŽnominateurs. Cela donne un homomorphisme de groupe
$$\eqalign{\left ({\alpha\over\cdot}\right )\ :\ I_K(\m)&\longrightarrow \mu_n\cr \aa&\longmapsto \left ({\alpha\over \aa}\right )_n\cr}$$
en posant  $\mu_n$ l'ensemble des racines $n$-ime de l'unitŽ de $\C$ et en supposant que $\m$ est n'importe quel $K$-module divisible pas des idŽaux premiers qui ne contiennent pas $n\cdot\alpha$.
\bigskip

\lem
\medskip

{\sl Sous les mmes hypothses que la dŽfinition prŽcŽdente, on a 

\art{a)}$\left ({\alpha\over \P}\right )_n=\left ({\beta\over \P}\right )_n$ si
$\alpha\equiv\beta\pmod\P$.

\art{b)}$\left ({\alpha\over \P}\right )_n\equiv \alpha^{\N(\P)-1\over n}\pmod\P$ pour tout $\alpha\in O_K$.

\art{c)}$\left ({\alpha\beta\over \P}\right )_n=\left ({\alpha\over \P}\right )_n\cdot\left
({\beta\over \P}\right )_n$.

\art{d)}$\left ({\alpha\over \P}\right )_n=1$ si et seulement s'il existe $\beta\in O_K\setminus\P$
tel que $\alpha\equiv \beta^{n}\pmod \P$.

\art{e)}Si $n=2$, $\zeta_2=-1$ et $K=\Q$, on retrouve le symbole de Legendre : soit $q\in\gfP(\Q)$ et $a\in \Z$ avec $(a,q)=1$, $\left ({a\over q}\right )=1$ si l'Žquation $x^2\equiv a\pmod q$ est rŽsoluble et $-1$ sinon.

\art{f)}l'application $a\mapsto \left ({a\over q}\right )$ est un homomorphisme surjectif de $(\Z/2p\Z)^*=\F_p^*$ sur $\{\pm 1\}$. 

}

{\bf Preuve}

Les parties a), b) et c) dŽcoulent de la dŽfinition.  Pour la partie d), posons $q=\N(\P)$. S'il existe $\beta\in O_k\setminus\P$ tel que $\alpha\equiv \beta^{n}\pmod \P$, alors $\alpha^{q-1\over
n}\equiv(\beta^n)^{q-1\over n}=\beta^{q-1}\equiv 1\pmod{\P}$. RŽciproquement, si $\left ({\alpha\over
\P}\right )_n=1$, alors $\alpha^{q-1\over n}\equiv 1\pmod\P$.  On se souvient que $\F^*$ est un
groupe cyclique engendrŽ par un ŽlŽment disons $\overline{\gamma}$. Donc
$\alpha\equiv\gamma^s\pmod\P$ pour un certain $1\leq s\leq q-1$. Ainsi $\alpha^{q-1\over n}\equiv
\gamma^{s\cdot {q-1\over n}}\equiv 1\pmod\P$. Ainsi, l'ordre de $\gamma$ qui est $q-1$ divise
$s\cdot {q-1\over n}$, c'est ˆ dire que
$n$ divise $s$, disons, $s=kn$. Finalement, $\alpha\equiv (\gamma^k)^n=\beta^n\pmod\P$ en posant
$\beta=\gamma^k$. La partie e) est un corollaire immŽdiat de la partie d). La partie f) suit de la partie c) et du fait que l'application de $x\mapsto x^2$ est un endomorphisme de $\F_p^*$ de noyau $\pm 1$; donc cet endomorphisme n'est pas injectif donc pas surjectif.\qed

\bigskip

\th
\medskip
{\sl Soit $K$ un corps de nombres contenant une racine primitive $n$-ime de l'unitŽ, $\alpha\in O_{K}$, $\m$ un $K$-module divisible par $n\cdot\alpha O_{K}$. Notons $L=K(\root n\of \alpha)$. Puisque $K$ contient les racines $n$-ime de l'unitŽ, il est Žvident que $L/K$ est une extension galoisienne, c'est mme une $n$-extension de Kummer, donc on se souvient de l'homomorphisme injectif $\psi(\root n\of \alpha)\ :\ {\rm Gal}(L/K)\to \mu_{n}$ vu lors de la preuve du ThŽorme \argdu. Alors, le diagramme suivant commute~:
\vglue 2.5cm 
\psset{xunit=1cm,yunit=1cm}

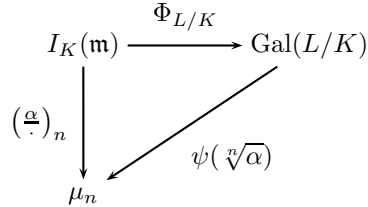
\rput(6,0){\rnode{b1}{$\mu_{n}$}}
\rput(6,2){\rnode{a1}{$I_{K}(\m)$}}
\rput(9,2){\rnode{a2}{${\rm Gal}(L/K)$}}
\ncline[nodesep=3pt]{->}{a1}{b1}
\Bput{$\left ({\alpha\over \cdot}\right )_n$} 
\ncline[nodesep=3pt]{->}{a1}{a2}
\Aput{$\Phi_{L/K}$} 
\ncline[nodesep=3pt]{->}{a2}{b1}
\Aput{$\psi(\root n\of \alpha)$}
\medskip

}
{\bf Preuve}

Remarquons tout d'abord que les idŽaux premiers de $K$ qui ramifient dans $L$ divisent $n\cdot\alpha$. En effet,  supposons que  $\P\in \gfP_{0}(K)$ ramifie dans $L$. Cela veut dire que $\P$  divise le discriminant de $L/K$ (cf. [Sam,Thm.1, Chap. 5, p.88]). Or, puisque $\alpha$ est entier, le discriminant de
$L/K$ divise le discriminant de $O_K[\root n\of{\alpha}]/O_K$ qui lui-mme divise le
discriminant de $f:=X^n-\alpha$ (cf. [Sam, prop. 1, Chap. 2,p.46] et la dŽfinition du discriminant d'un polyn™me [Fr-Tay,  rel. 1.5 et 1.9, pp.10-11]). Et finalement, si $f=\prod_{i=1}^n(X-\alpha_i)$ avec
$\alpha_1=\root n\of\alpha$, on a
$$\pm {\rm Disc}(f)=\prod_{i\ne
j}(\alpha_i-\alpha_j)=\prod_{i=1}^nf'(\alpha_i)=n^n\cdot\prod_{i=1}^n\alpha_i^{n-1}=
n^n\cdot\alpha^{n-1}\cdot\hbox{(rac. de l'unitŽ)}.$$
Donc, l'application $\Phi_{L/K}$ est bien dŽfinie. Soit $\P\in I_{K}(\m)$. Pour prouver le thŽorme, il suffit de prouver que 
$${\rm Frob}_{L/K}(\P)(\root n\of\alpha)=\left ({\alpha\over \P}\right )_n\cdot\root n\of \alpha.\eqno{(*)}$$
D'une part, on a ${\rm Frob}_{L/K}(\P)(\root n\of\alpha)=\omega\cdot\root n\of \alpha$, avec $\omega\in \mu_{n}$. D'autre part, par dŽfinition de l'automorphisme de Frobenius, ${\rm Frob}_{L/K}(\P)\equiv \left (\root n\of \alpha\right )^{\N(\P)}\pmod{\P\cdot O_{L}}$. Mais $ \left(\root n\of \alpha\right )^{\N(\P)}=\alpha^{\N(\P)-1\over n}\cdot \root n\of \alpha\equiv \left ({\alpha\over \P}\right )_n\cdot\root n\of \alpha\pmod{\P\cdot O_{L}}$. Donc $\omega\cdot \root n\of \alpha\equiv \left ({\alpha\over \P}\right )_n\cdot\root n\of \alpha\pmod{\P\cdot O_{L}}$. Or, $\root n\of\alpha$ est premier ˆ $\P\cdot O_{L}$ par hypothse, donc $\omega\equiv \left ({\alpha\over \P}\right )_n\pmod{\P\cdot O_{L}}$ et finalement, $\omega=\left ({\alpha\over \P}\right )_n$ en vertu du Lemme \argo.\qed
\bigskip
\coro

{\sl Sous les mmes hypothses, on suppose de plus que $\m$ est admissible pour $L/K$ (c'est possible, en vertu du thŽorme de rŽciprocitŽ d'Artin (ThŽorme \argda)) et on note $\overline{G}$ l'image de ${\rm Gal}(L/K)$ par $\psi(\root n\of \alpha)$. Alors $ \left ({\alpha\over \cdot}\right )_n$ induit un homomorphisme surjectif $I_{K}(\m)/P_{\m}\to\overline{G}$, qu'on notera encore $ \left ({\alpha\over \cdot}\right )_n$.

}
{\bf Preuve}

Puisque $\m$ est admissible pour $L/K$, $P_\m\subset \ker (\Phi_{L/K})$, donc en vertu du thŽorme prŽcŽdent, $P_\m\subset \ker ( \left ({\alpha\over \cdot}\right )_n)$, et on conclut en se souvenant que $\Phi_{L/K}$ est surjective (cf. ThŽorme \argaq).\qed
\bigskip
\coro {\soustitre (rŽciprocitŽ quadratique)}

{\sl Si $p$ et $q$ sont des nombres premiers impairs, alors on a
$$\left ({p\over q}\right )\cdot \left ({q\over p}\right )=(-1)^{{p-1\over 2}\cdot {q-1\over 2}}.$$

}

{\bf preuve}

Le ThŽorme \argm\ nous apprend que $\Q(\zeta_p)$ est le corps de la classe d'Žquivalence du groupe de congruence $P_\n$, o $\n$ est le $\Q$-module $(p)\cdot\infty$, o $\infty$ est l'unique place infinie de $\Q$. Or, on sait que $\Q(\zeta_p)/\Q$ est une extension cyclique d'ordre $p-1$ pair. Par la thŽorie de Galois, $\Q(\zeta_p)$ contient une unique extension quadratique de $\Q$. Dans cette extension, $p$ est le seul nombre premier qui ramifie (car c'est le seul qui ramifie dans $\Q(\zeta_p)$). Or, on sait que le discriminant de $\Q(\sqrt{m})=\cases{4m &si $m\equiv 2,3\pmod 4$\cr m& si $m\equiv 1\pmod 4$\cr}$ (cf.  [Sam, exemple, p. 89]) et qu'un nombre premier ramifie si et seulement s'il divise ce discriminant (cf. [Sam, Thm. 1, p. 88]). Ainsi, puisque 2 ne ramifie pas, ce discriminant vaut $\pm p\equiv 1\pmod 4$. Donc ce sous-corps est $\Q(\sqrt{p^*})$, o $p^*=(-1)^{p-1\over 2}\cdot p$.

On a vu lors de la remarque qui suit le Lemme \argcp\ que que $(p)\cdot\infty$ Žtait admissible pour l'extension $\Q(\zeta_{p})/\Q$. Donc, par le Lemme \argcq, $(p)\cdot\infty$ est aussi admissible pour $\Q(\sqrt{p^*})/\Q$, et donc, par le Lemme~\argcp, $(2p)\cdot\infty=:\m$ est admissible pour $\Q(\sqrt{p^*})/\Q$. On peut donc appliquer le corollaire prŽcŽdent avec $K=\Q$, $n=2$ et $\alpha=p^*$. Dans ce cas, l'application $\psi:=\psi(\sqrt{p^*})\ :{\rm Gal}(\Q(\sqrt{p^*})/\Q)\to\mu_{2}=\{\pm1\}$ est un isomorphisme, car si $\sigma\ne{\rm Id}$, on a $\sigma(\sqrt{p^*})=\psi(\sigma)\cdot \sqrt{p^*}=(-1)\cdot \sqrt{p^*}$. Donc, $\overline{G}=\mu_{2}$. On a donc par le corollaire prŽcŽdent un homomorphisme surjectif $\left ({p^*\over \cdot }\right )_{2}\ : I_{\Q}(\m)/P_{\m}\to\mu_{2}$. Or, on a montrŽ que $ I_{\Q}(\m)/P_{\m}\simeq(\Z/2p\Z)^*\simeq (\Z/p\Z)^*$ (cf. ThŽorme \argm). D'o un homomorphisme surjectif $\varphi\, :\, (\Z/p\Z)^*\to\mu_{2}$. Soit $q$ un premier impair diffŽrent de $p$. Suivons ˆ la trace l'image de $q\bmod p$ par $\varphi$. On a  $q\bmod p\mapsto  q\bmod {2p}\mapsto q\Z\bmod {P_{\m}}\mapsto \left ({p^*\over q\Z}\right)_{2}\buildrel {\rm Lemme}\ \argeb\over =\left ({p^*\over q}\right )$. En bref, $\varphi(q\bmod p)=\left ({p^*\over q}\right )$. D'autre part, la partie f) du Lemme \argeb\ montre que $\left ({\cdot\over p}\right )\, :\, (\Z/p\Z)^*\to \mu_{2}$ est aussi un homomorphisme surjectif. Le noyau de cette application est un sous-groupe d'indice 2 de $(\Z/ p\Z)^*$ qui est cyclique (cf. [Jac1, Theorem 2.8, p.132]). Or, il n' y a qu'un sous-groupe d'indice donnŽ dans un groupe cyclique. Cela montre que $\left ({\cdot\over p}\right )$ et $\varphi$ ont le mme noyau. Puisque leur image est $\mu_{2}$, c'est forcŽment les mmes applications. Cela montre que $\left ({p^*\over q}\right )=\left ({q\over p}\right )$. De cela et des parties b) et c) du Lemme \argeb, on tire finalement
$$\left ({p\over q}\right )\cdot \left ({q\over p}\right )=\left ({p\over q}\right )\cdot \left ({p^*\over q}\right )=\left ({p\over q}\right )\cdot \left ({p\over q}\right )\cdot \left ({(-1)^{p-1\over 2}\over q}\right )= \left ({-1\over q}\right )^{p-1\over 2}=(-1)^{{p-1\over 2}\cdot {q-1\over 2}}.$$\qed
\bigskip
Les esprits chagrins rŽtorquerons que c'est la preuve la plus compliquŽe de ce fameux thŽorme et qu'on ne montre mme pas que $\left ({2\over p}\right )=(-1)^{p^2-1\over 8}$. Et ils auraient raisons ! Mais au moins, on comprend le lien entre rŽciprocitŽ d'Artin et rŽciprocitŽ quadratique. Nous ne faisons donc ici qu'expliquer un mot.

\vfill\eject
\global\advance\chapnomb by 1
\nomb=1

\centerline{\para Chapitre 10 }
\medskip
\centerline{\para Le thŽorme principal du corps de classe }
\bigskip
Rappelons l'ŽnoncŽ de ce thŽorme, que nous avons dŽjˆ ŽnoncŽ au Chapitre 8 (ThŽorme \argdm)~: 
\bigskip

{{{\soustitre ThŽorme ({\the\chapnomb}.{\the\nomb})}\global\advance\nomb by 1} {\soustitre (ThŽorme d'existence du corps de classe)}}
\medskip
{\sl Soit $K$ un corps de nombres. Alors pour toute
classe $\bbH$ d'Žquivalence de sous-groupes de congruence de $K$, il existe une
extension abŽlienne $L/K$ telle que $\bbH=\bbH(L/K)$. On dit alors que $L/K$ est {\it
le corps de classe de $\bbH$} (on devrait plut™t dire le corps de la classe $\bbH$).

}
\bigskip
Rappelons encore exactement de quoi il s'agit~: un sous-groupe $H$ de $I_K$ est dit Òde congruence" s'il existe un $K$-module $\m$ tel que $P_m \subset H\subset I_K(\m)$. Deux groupes de congruences $H_1$ et $H_2$ sont dit Žquivalents s'il existe un $K$-modules $\m''$ tel que $H_1\cap I_K(\m'')=H_2\cap I_K(\m'')$. On a montrŽ au  Corollaire-DŽfinitions \argdi\ que si $\bbH$ est une classe d'Žquivalence  alors il existe un $K$-module $\ff$ tel que $\bbH=\{H(\m)\mid \ff | \m\}$, avec $H(\m)=H(\ff)\cap I_K(\m)$, et que tous les groupes $I_K(\m)/H(\m)$ avec $\ff |\m$ sont isomorphe, on note ce groupe $I_K/\bbH$ (cf. DŽfinition \argdn) . D'autre part, si $L/K$ est une extension abŽlienne de $K$ de degrŽ fini, alors les noyaux des applications d'Artin $\Phi_{L/K}\ :\ I_K(\m)\to {\rm Gal}(L/K)$ forment sont Žquivalents et la classe de telle groupe ce note $\bbH(L/K)$. Puisque l'application d'Artin est surjective (cf. ThŽorme \argaq), dans ce cas, on a donc $I_K/\bbH(L/K)\simeq {\rm Gal}(L/K)$. Au Chapitre 8 (Proposition \1argdm) on a montrŽ que l'application
$$L/K\longmapsto \bbH(L/K)$$
Žtait une application injective de l'ensemble des extensions abŽliennes finies de $K$ (incluses dans $\C$) dans l'ensemble des classes d'Žquivalences de sous-groupes de congruence pour $K$. Le ThŽorme \argef\ dit simplement que cette application est surjective.

\headline={\hfill \phantom{ouuh}\hfill}
\medskip
Rappelons aussi qu'en vertu du ThŽorme \argds\ qu'il suffit de prouver ce rŽsultat dans le
cas o $K$ possde une racine primitive $n$-ime de l'unitŽ et que $n$ est un exposant du groupe
$I_K/\bbH$.

Mais avant cela, il va falloir dŽmontrer un gros thŽorme (dž vraisemblablement ˆ Herbrand)
qui aura comme corollaire immŽdiat ce qu'on recherche. Donnons donc le cadre de ce thŽorme qui est, il faut bien l'avouer, un peu surprenant~:
\medskip
On suppose que $K$ possde une racine primitive $n$-ime de l'unitŽ. On pose $S_1$ et $S_2$ deux
ensembles finis de places telles que $S_1\cap S_2=\emptyset$ et $S=S_1\cup S_2$ remplit les trois
conditions suivantes~:
\art{i)}$S$ contient toutes les places infinies.

\art{ii)}$S$ contient toutes les places finies $\P$ telles que $\P|n\cdot O_K$

\art{iii)}$S$ contient toutes les places finies $\P$ telles que  $\P|\aa_1\cdots \aa_k$, o $\aa_1,\ldots
,\aa_k$ est un systme fixŽ de reprŽsentants de toutes les classes de $I_K$ modulo $P_K$
\medskip
Soit $\m_1$ (resp. $\m_2$) un $K$-module tel que l'ensemble de places qui divisent $\m_1$ (resp. $\m_2$)
soit exactement les places non complexes de $S_1$ (resp. $S_2$).

On dŽfinit deux sous-groupes de congruence

$$\eqalign{H_1=P_{\m_1}\cdot (I_K(\m_1))^n\cdot I_K[S_2]\subset I_K(\m_1)\cr 
H_2=P_{\m_2}\cdot (I_K(\m_2))^n\cdot I_K[S_1]\subset I_K(\m_2)\cr}$$

Il est clair que $H_1$ est dŽfinit modulo $\m_1$ et que $H_2$ est dŽfinit modulo $\m_1$. On dŽfinit aussi
les sous-groupes de $K^*$ suivants

$$\eqalign{W_1={K}^S\cdot {K^*}^n\cap K^*_{\m_2}\cr
W_2={K}^S\cdot {K^*}^n\cap K^*_{\m_1}.\cr}$$
Le lecteur attentif aura remarquŽ qu'on a ÒcroisŽ" les indices ! Et finalement, pour $i=1,2$, on pose
$$L_i=K(\root n\of{W_i}).$$
On remarque dŽjˆ que les extensions $L_i/K$ sont de degrŽ fini. En effet~: on voit facilement que l'image
de $W_i$ par l'application $K^*\to K^*/{K^*}^n$ est $W_i{K^*}^n/{K^*}^n$. Et on a
$W_i{K^*}^n/{K^*}^n\subset {K}^S{K^*}^n/{K^*}^n\buildrel\rm thm.d'isom\over\simeq {K}^S/{K}^S\cap
{K^*}^n={K}^S/({K}^S)^n$, qui est un groupe d'ordre $n^{|S|}$ par le thŽorme Dirichlet-Chevalley-Hasse (cf. Corollaire \argdx). Ainsi, on voit aisŽment (en utilisant un mme raisonnement qu'au ThŽorme \argdu) que les extensions $L_i/K$ sont des $n$-extensions de Kummer. Par le mme ThŽorme \argdu, on a~:
$${\rm Gal}(L_i/K)\simeq W_i {K^*}^n/{K^*}^n \buildrel\rm thm.d'isom\over\simeq W_i/(W_i\cap
{K^*}^n).\eqno{(*)}$$
\bigskip\goodbreak
\headline={\hfill \smcap Le thŽorme principal du corps de classe\hfill}
\th
\medskip
{\sl Sous les mmes hypothses et en supposant que les idŽaux premiers finis qui divisent $\m_1$ et
$\m_2$ apparaissent dans $\m_1$ et $\m_2$ ˆ des puissances suffisamment grandes, alors, pour $i=1,2$, on
a~:

$$H_i\in \bbH(L_i/K)\quad \hbox{et $\m_i$ est admissible pour $L_i/K$}.$$

}

{\bf Preuve}

\art{(I)}On suppose que pour chaque idŽal premier $\P$ qui divise $\m_1$ (resp. $\m_2$), la puissance
$\P^t|\m_1$ (resp. $\m_2$) soit assez grande pour que $U_\P^{(t)}\subset U_\P^n$. (cf. Proposition \argcc). Alors on affirme que toute place de $S_1$ (resp. $S_2$) se dŽcompose compltement dans
$L_2$ (resp. $L_1$). En effet, prouvons-le pour $\P\in S_1$ (la preuve pour un ŽlŽment de $S_2$ est
identique). Il suffit de montrer que $[{L_2}_\gP:\bbK_\P]=1$ pour tout place $\gP$ au-dessus de $\P$, car $[{L_2}_\gP:\bbK_\P]=e(\gP/\P)\cdot f(\gP/\P)$ (cf. [Fr-Tay, 1.14, p. 111] pour les places finies et aux DŽfinitions \argbn\ pour les places infinies). Supposons que
$\P$ est infinie complexe, alors ${L_2}_\gP=\bbK_\P=\C$, donc c'est en ordre. Si $\P$ est infinie rŽelle,
alors $\bbK_\P=\R$  et $n=2$ (car $K$ a un plongement rŽel et une racine $n$-ime de l'unitŽ, forcŽment
$-1$). Comme $W_2\subset K^*_{\m_1}$ et que $\P|\m_1$, on a $W_2\subset \R_+^*$. Donc, puisque les
racines carrŽes de nombres positifs existent dans $\R$, on a
${L_2}_\gP=\bbK_\P(\sqrt{W_2})=\R(\sqrt{W_2})=\R$, donc c'est aussi bon pour ce cas-lˆ. Enfin, si $\P$ est
une place finie et si $\P^t|\m_1$, alors $W_2\subset K^*_{\m_1}\subset U_\P^{(t)}\subset U_\P^n$ et donc
$\root n\of{W_2}\subset \bbK_\P$ et donc ${L_2}_\gP=\bbK_\P$, ce qui rgle la partie (I). 

\art{(II)}Sous les mmes hypothses que (I), alors les places qui ramifient dans $L_i$ sont dans $S_i$
(pour $i=1,2$). Montrons-le pour $i=1$. Soit $\P$ une place qui ramifie dans $L_1$. Il suffit de montrer
que $\P\in S$ (par la partie (I)). Si $\P$ est infini, c'est vrai par les propriŽtŽs postulŽes sur $S$
(qui possde toutes les places infinies). On suppose donc $\P$ finie. Il suffit de montrer que si
$\alpha\in {K^*}^S$ et si $\P\not\in S$, alors $\P$ ne ramifie pas dans $K(\root n\of\alpha)$, car
$W_1\subset {K^*}^S\cdot {K^*}^n$ et $L_1$ est un produit de tels extensions et on sait que si un idŽal ne
ramifie pas dans deux corps, alors il ne ramifie pas dans le produit de ces deux corps (cf. [Mar,
Thm. 31, p.107]). De plus, on peut supposer que $\alpha\in O_K\cap {K^*}^S$. En effet, l'idŽal $\alpha
O_K={\aa\over\bb}$, avec $\aa$ et $\bb$ des idŽaux entiers uniquement divisibles par des premiers de
$S$. Soit $h=|I_K/P_K|$. Alors on a $\alpha O_K={\aa\cdot \bb^{h-1}\over \bb^h}={\aa \cdot
\bb^{h-1}\over\beta O_K}$, avec $\beta\in O_K\cap {K^*}^S$. Ainsi, $\alpha={\gamma\over\beta}$, avec
$\gamma\in O_K\cap {K^*}^S$. Ainsi, $K(\root n\of\alpha)=K(\root
n\of{\gamma\over\beta})=K(\root\n\of{\beta^{n-1}\gamma})$, car $\root n\of{\beta^n}\in K$ puisque $K$
possde toutes les racines $n$-imes de~1. Pour montrer que $\P$ ne ramifie pas dans $K(\root
n\of{\alpha})$, il suffit de montrer que $\P$ ne divise pas le discriminant du polyn™me $f:=X^n-\alpha$.
En effet (on a dŽjˆ donnŽ ce raisonnement, mais il Žtait dans un interlude, donc pas obligatoire; maintenant il le devient !), il est bien connu que si $\P$ ne divise pas le discriminant de $K(\root n\of{\alpha})/K$ alors
il ne ramifie pas (cf. [Sam, Thm.1, Chap. 5, p.88]). Or, puisque $\alpha$ est entier, le discriminant de
$K(\root n\of{\alpha})/K$ divise le discriminant de $O_K[\root n\of{\alpha}]/O_K$ qui lui-mme divise le
discriminant de $f$ (cf. [Sam, prop. 1, Chap. 2, p.46] et la dŽfinition du discriminant d'un polyn™me
[Fr-Tay, rel. 1.5 et 1.9, pp.10-11]). Et finalement, si $f=\prod_{i=1}^n(X-\alpha_i)$ avec
$\alpha_1=\root n\of\alpha$, on a
$$\pm {\rm Disc}(f)=\prod_{i\ne
j}(\alpha_i-\alpha_j)=\prod_{i=1}^nf'(\alpha_i)=n^n\cdot\prod_{i=1}^n\alpha_i^{n-1}=
n^n\cdot\alpha^{n-1}\cdot\hbox{(rac. de l'unitŽ)}.$$
Il est Žvident que $\P$ ne divise pas $\alpha\cdot O_K$, car $\alpha\in {K^*}^S$ et $\P$ ne divise pas
$n\cdot O_K$, par la relation (ii) dŽfinissant $S$. Cela montre $(II)$.

Si on augmente les exposants des diviseurs premiers de $\m_1$ et de $\m_2$, cela a pour effet  de
remplacer les $L_1$ et $L_2$ originaux par des sous-extensions de ces  $L_i/K$. Et comme il n'y a qu'un nombre fini de sous-extensions de ces $L_1/K$ et $L_2/K$, il arrive une situation o augmenter encore ces exposants ne change plus les corps $L_1$ et $L_2$ correspondants (quitte ˆ arriver ˆ $K$, a ne fait rien). On supposera dans la suite que les exposants sont assez grands pour remplir cette condition de stabilitŽ, et alors, en augmentant encore les exposants si nŽcessaire, gr‰ce au thŽorme de rŽciprocitŽ d'Artin, on peut supposer que $\m_1$ (resp. $\m_2$) est admissible pour $L_1/K$ (resp. pour $L_2/K$). Posons alors
$$H_i^*:=\ker(\Phi_{L_i/K}|_{I_K(\m_i)})=P_{\m_i}\cdot N_{L_i/K}(I_{L_i}(\widetilde{\m_i}))\in\bbH(L_i/K).$$

Il reste ˆ montrer que $H_i^*=H_i$ ($i=1,2$) pour terminer la preuve. Regardons dŽjˆ pour $i=1$. On se souvient que $H_1=P_{\m_1}\cdot (I_K(\m_1))^n\cdot I_K[S_2]$. Soit $\P\in S_2$. On a vu en (I) que $\P$ se dŽcompose totalement dans $L_1$, donc est la norme d'un idŽal de $L_1$ qui est premier et dans $I_{L_1}(\widetilde{\m_1})$. Donc, $I_K[S_2]\subset N_{L_1/K}(I_{L_1}(\widetilde{\m_1}))$. D'autre part, $I_K(\m_1)/H_1^*\buildrel \rm appl. d'Artin\over\simeq {\rm Gal}(L_1/K)$  qui est par hypothse d'exposant $n$. Donc $I_K(\m_1)^n\subset H_1^*$. Ce qui montre que $H_1\subset H_1^*$. On montre de mme que $H_2\subset H_2^*$. Il reste ˆ voir, pour terminer le thŽorme, que les indices de $H_i$ et $H_i^*$ dans $I_K(\m_i)$ sont Žgaux (pour $i=1,2$). Autrement dit, on sait que 
$$[I_K(\m_i):H_i]\geq [I_K(\m_i):H_i^*]=[L:K]\buildrel(*)\over=[W_i:W_i\cap {K^*}^n],$$
et il faut montrer qu'il y a ŽgalitŽ. Tout se rŽsume donc ˆ montrer que 
$${[I_K(\m_1):H_1]\cdot [I_K(\m_2):H_2]\over [W_1:W_1\cap {K^*}^n]\cdot [W_2:W_2\cap {K^*}^n]}=1.\eqno{(**)}$$
En vertu de la Remarque b), du Chapitre 8 et du Lemme \argdg\ b), on a pour $i=1,2$,
$$I_K(\m_i)/H_i\simeq I_K^S/(I_K^S\cap H_i),$$
o, bien sžr,  $I_K^S=I_K(\m_1\m_2)$.
\bigskip
\goodbreak
{\soustitre Lemme A}

{\sl On a, pour $i=1,2$

$$K^*/({K^*}^n\cdot {K^*}^S\cdot K^*_{\m_i})\simeq I_K^S/(I_K^S\cap H_i).$$

}

{\bf Preuve du Lemme A}

Regardons l'application $f$ composŽe

$$K^*\buildrel \iota\over\longrightarrow I_K\buildrel j\over\longrightarrow I_K^S$$
o, comme toujours, l'application $\iota$ est dŽfinie par $\iota(a)=a\cdot O_K$, qu'on note parfois $(a)$ quand il n'y a pas d'ambigu•tŽ et $j(\P)=\cases{\P&si  $\P\not\in S$\cr O_K& si $\P\in S$\cr}$.

\art{a)}$f$ est surjective~: soit $\aa\in I_K^S$. Par la propriŽtŽ (iii) de $S$, il existe $\alpha\in K^*$ et $\bb\in I_K[S]$ tels que $\aa=\bb\cdot (\alpha)$. Ainsi, $f(\alpha)=j(\aa\cdot\bb^{-1})=\aa$. 

\art{b)}Il reste ˆ voir que $\{\alpha\in K^*\mid f(\alpha)\in H_i\}={K^*}^n\cdot {K^*}^S\cdot K^*_{\m_i}$ pour $i=1,2$. Montrons dŽjˆ $Ò\subset"$. Soit $\alpha\in K^*$ tel que  $f(\alpha)\in H_i$. Ecrivons $\iota(\alpha)=\aa_0\cdot\aa_1$, avec $\aa_0\in I_K[S]$ et $\aa_1\in I_K^S$; on a donc $f(\alpha)=\aa_1\in H_i$. Ainsi, par dŽfinition de $H_i$,  $\aa_1=\bb^n\cdot (\beta)\cdot \cc$, avec $\bb\in I_K(\m_i)$, $\beta\in K^*_{\m_i}$ et $\cc\in I_K[S_{3-i}]$, et, comme en a), on peut Žcrire $\bb=\bb_0\cdot (\theta)$ avec $\bb_0\in I_K[S]$ et $\theta\in K^*$. Alors, $(\alpha\cdot\theta^{-n}\cdot\beta^{-1})=\aa_0\cdot\bb_0^n\cdot\cc\in I_K[S]$, i.e. $\alpha\cdot\theta^{-n}\cdot\beta^{-1}\in {K^*}^S$, ce qui montre que $\alpha\in {K^*}^n\cdot {K^*}^S\cdot K^*_{\m_i}$.  Montrons maintenant $Ò\supset" $. 
\artart{i)}Si $\alpha\in K^*_{\m_i}$, alors Žcrivons $(\alpha)=\aa\cdot\bb$, avec $\aa\in I_K^{S_{3-i}}$ et $\bb\in I_K[S_{3-i}]$. Comme $\alpha$, donc aussi $\aa$ est premier ˆ $\m_i$, on a $\aa\in I_K^S$. Et alors $f(\alpha)=\aa=(\alpha)^{-1}\cdot \bb\in P_{\m_i}\cdot I_K[S_{3-i}]\subset H_i$.

\artart{ii)}$f({K^*}^S)=\{O_K\}\subset H_i$.

\artart{(iii)}$f({K^*}^n)=f(K^*)^n\subset (I_K^S)^n\subset I_K(\m_i)^n\subset H_i$.

Ce qui achve la preuve du Lemme A.\hfill$\qed$ (Lemme A)
\bigskip
\goodbreak
{\soustitre Lemme B}
\medskip

{\sl Soient $A,B,C$ des sous-groupes d'un groupe abŽlien $T$ (notŽ multiplicativement). On suppose que $B\subset A$ et 
$[A:B]<\infty$. Alors
$$[A:B]=[AC:BC]\cdot [A\cap C :B\cap C].$$

}

{\bf Preuve du Lemme B}

 Prenant $\beta : T\to T/C$ restreint ˆ $A$, le Lemme \argcj\ nous apprend que 
$$[A:B]=[\beta(A):\beta(B)]\cdot [\ker(\beta):\ker(\beta)\cap B].$$
Dans notre cas, $\beta(A)=AC/C$ et $\beta(B)=BC/C$. Donc, $[\beta(A):\beta(B)]=[AC:BC]$ et 
$\ker(\beta)=A\cap C$ et $B\cap\ker(\beta)=B\cap A\cap C=B\cap C$. \hfill$\qed$ (Lemme B)
\bigskip

 Reprenons le fil de notre calcul~:  il est Žvident (transitivitŽ des indices) que pour $i=1,2$~:
 
 $$[K^*:{K^*}^n\cdot {K^*}^S\cdot K^*_{\m_i}]={[K^*:{K^*}^n\cdot K^*_{\m_i}]\over 
[{K^*}^n\cdot {K^*}^S\cdot K^*_{\m_i}:{K^*}^n\cdot K^*_{\m_i}]}.$$
 On applique le Lemme B ˆ $A={K^*}^n\cdot {K^*}^S, B={K^*}^n$ et  $C=K^*_{\m_{i}}$, ce qui donne
 
 $$[{K^*}^n\cdot {K^*}^S:{K^*}^n]=[{K^*}^n\cdot {K^*}^S\cdot K^*_{\m_i}:{K^*}^n\cdot K^*_{\m_i}]\cdot [\underbrace{({K^*}^n\cdot {K^*}^S)\cap K^*_{\m_{i}}}_{=W_{3-i}}:{K^*}^n\cap K^*_{\m_{i}}].$$
 Rappelons que 
 $${K^*}^n\cdot {K^*}^S/{K^*}^n\buildrel\rm thm.\ d'isom\over\simeq {K^*}^S/({K^*}^S\cap {K^*}^n)=
{K^*}^S/({K^*}^S)^n$$ et le dernier de ces groupes ˆ $n^{|S|}$ ŽlŽments en vertu du Corollaire \argdx. Et on a $W_{3-i}={K^*}^n\cdot{K^*}^S\cap K^*_{\m_i}$ et $W_{3-i}\cap {K^*}^n={K^*}^n\cap
K^*_{\m_i}$. Enfin, on rŽsume tout ce que nous venons de voir depuis la relation~$(**)$~:
$$\eqalign{[I_K(\m_i):H_i]&=[I_K^S:I_K^S\cap H_i]\buildrel\rm Lemme\ A\over=[K^*:{K^*}^n\cdot{K^*}^S\cdot K^*_{\m_i}]\cr
&={[K^*:{K^*}^n\cdot K^*_{\m_i}]\over [{K^*}^n\cdot{K^*}^S:{K^*}^n]}\cdot [({K^*}^n\cdot
{K^*}^S)\cap{K^*}_{\m_i}:{K^*}^n\cap K^*_{\m_i}]\cr
&={[K^*:{K^*}^n\cdot K^*_{\m_i}]=c(\m_i)\over n^{|S|}}\cdot [W_{3-i}:W_{3-i}\cap {K^*}^n].\cr }$$

Finalement,  pour prouver la relation $(**)$, il reste ˆ montrer la relation~:

$$c(\m_1\cdot\m_2):=c(\m)=n^{2\cdot |S|}.\eqno{(***)}$$
Posons alors $\m=\prod_{\P\in S'}\P^{n_\P}$, o $S'$ est l'ensemble des places non complexes de $S$.   Alors on a $c(\m)=\prod_{\P\in S'}c(\P^{n_\P})$ (cf. ThŽorme \argdz\ a)). Soit $S_0\subset S'$ le sous-ensemble des places finies de $S'$. Soit $\P\in S_0$. Dans la partie c) du mme thŽorme, on a vu que $c(\P^{n_\P})=n^2\cdot [\okp:n\cdot\okp]$ (si $n_{\P}$ est assez grand). Ainsi,

$$\eqalign{\prod_{\P\in S_0}c(\P^{n_\P})&=n^{2|S_0|}\cdot \prod_{\P\in S_0}[\okp :n\cdot\okp]\cr &
=n^{2|S_0|}\cdot \prod_{\P\in S_0}[\okp : \P^{v_\P(n)}\okp]\cr
&=n^{2|S_0|}\cdot \prod_{\P\in S_0}[O_K:\P^{v_\P(n)}].\cr}$$
Par le thŽorme chinois, on a
$$\prod_{\P\in S_0}\left |O_K/\P^{v_\P(n)}\right |=\left | O_K/\prod_{\P\in S_0}\P^{n_\P(n)}\right |\buildrel\rm cond.\ ii)\over =\left |O_K/n\cdot O_K \right | =n^{[K:\Q]}.$$
Ce qui montre que 
$$\prod_{\P\in S_0}c(\P^{n_\P})=n^{2|S_0|+[K:\Q]}.$$
Soit $r$ (resp. $s$) le nombre de places rŽelles (resp. complexes) de $K$. On sait que $r+2s=[K:\Q]$.

Si $r=0$, alors $[K:\Q]=2s$ et $S'=S_0$. Ainsi, $c(\m)=\prod_{\P\in S_0}c(\P^{n_\P})=n^{2|S_0|+2s}=n^{2(|S_0|+s)}=n^{2\cdot |S|}$, car on se souvient que $S$ contient toutes les places infinies (condition i)). 

Si $r>0$, alors $n=2$ (se souvenir pourquoi...) et en vertu du ThŽorme \argdz\ b), on a $c(\m)=\prod_{\P\in S_0} c(\P^{n_\P})\cdot 2^r=(2^{2|S_0|+r+2s})\cdot 2^r=2^{2\cdot|S|}$. Ce qui achve la dŽmonstration du thŽorme.\qed
\bigskip
Enfin, nous pouvons achever la preuve de ce grand et beau thŽorme
\bigskip
{\soustitre Fin de la preuve du thŽorme d'existence du corps de classe}
\medskip
Soit donc $\bbH$ une classe de sous-groupes de congruence. On rappelle qu'on peut supposer que $K$ contient les racines $n$-imes de l'unitŽ et que $n$ est un exposant de $I_K/\bbH$.  Soit $\m$ un $K$-module multiple de $\ff(\bbH)$ et $H(\m)\in\bbH$ le sous-groupe de congruence dŽfini modulo $\m$.  

Appliquons le ThŽorme \argeg\ au cas o $S_2=\emptyset$ et $\m_2={\euf 1}=O_K\cdot\emptyset$. Suppose que $S_1$ soit construit de telle manire que $S=S_1$ satisfasse les conditions i), ii), iii) du ThŽorme \argeg\ et contiennent en plus toutes les places qui divisent $\m$. Et on dŽfinit $\m_1$ toujours comme dans le ThŽorme \argeg\ avec en plus $\m|\m_1$, donc $H(\m_{1})$ existe. Le groupe $H_1$ du ThŽorme \argeg\ est ici
$$H_1=P_{\m_1}\cdot (I_K(\m_1))^n\cdot( 1 )\subset H(\m_1),\eqno{ (+)}$$
 car, par hypothse $I_K(\m_1)/H(\m_1)$ est d'exposant $n$. Par le ThŽorme \argeg, il existe $L_1/K$ une extension abŽlienne telle que $H_1\in \bbH(L_1/K)$. Par la relation $(+)$ et par le Corollaire-DŽfinition \argdi, cela implique que  $\bbH(L_1/K)\subset \bbH$. Donc, en vertu du ThŽorme \argdo, il existe $L/K$ une extension abŽlienne telle que $\bbH=\bbH(L/K)$. Et c'est ce qu'il fallait dŽmontrer.\qed
 \bigskip\goodbreak
 \centerline{\soustitre Application ˆ des extensions non abŽliennes}
 \medskip
 Soit $E/K$ une extension galoisienne de corps de nombres de groupe de Galois $G$. Notons $G'=D(G:G)$ le sous-groupe engendrŽ par les commutateurs $[a,b]=aba^{-1}b^{-1}$, $a,b\in G$. Dans la littŽrature, les gens notent plut™t $[G:G]$ plut™t que $D(G,G)$. Mais $[\cdot,\cdot]$ dŽnote dŽjˆ deux notions : l'indice d'un groupe dans un autre et la dimension d'un corps dans un autre, nous n'allons pas encore Òcharger" la notation ! On voit facilement que $G'$ est normal dans $G$ et $G^{\rm ab}=G/G'$ est  l'abŽlianisŽ de $G$ (qui est le plus grand quotient abŽlien de $G$). On dŽfinit une application d'Artin 

$$\Phi_{E/K}\, :\, I_K(\m)\longrightarrow G^{\rm ab},$$
 \newcount\gaagam\gaagam=\gaga
pourvu que le $K$-module $\m$ soit divisible par toutes les places de $K$ qui ramifient dans $E$ ainsi : si $\P\in I_K(\m)$, on choisit un idŽal $\gP$ de $E$ au-dessus de $\P$ et on pose $\Phi_{E/K}(\P)={\rm Frob}(\gP/\P)\cdot G'\in G^{\rm ab}$. L'application est bien dŽfinie pour $\P$, car si on avait pris un autre premier $\gP'$ au-dessus de $\P$, l'ŽlŽment ${\rm Frob}(\gP'/\P)$ est un conjuguŽ de ${\rm Frob}(\gP/\P)$ et on voit facilement qu'ils reprŽsentent la mme classe dans $G^{\rm ab}$; et on prolonge comme toujours par multiplicativitŽ. Si $L$ est la sous-extension de $E/K$ fixe par $G'$, alors $L/K$ est la plus grand sous-extension abŽlienne de $E/K$. Soit $\m$, comme avant,  un $K$-module divisible par les places de $K$ qui ramifient dans $E$ et supposons que les exposants des diviseurs finis de $\m$ soient assez grand pour que $\m$ soit admissible pour toutes les sous-extensions abŽliennes de $E/K$ (c'est possible en  vertu du thŽorme de rŽciprocitŽ d'Artin (ThŽorme \argda) et puisqu'il n'y a qu'un nombre fini de telles sous-extensions. ConsidŽrons la classe de sous-groupes de congruence qui contient le sous-groupe $P_{\m}\cdot N_{E/K}(I_{E}(\widetilde{\m}))=H(\m)$ (qui est dŽfini modulo $\m$). Soit $F/K$ l'extension abŽlienne correspondant ˆ cette classe de sous-groupes de congruence, i.e. $P_{\m}\cdot N_{E/K}(I_{E}(\widetilde{\m}))\in \bbH(F/K)$ (cette extension existe, bien sžr, en vertu du thŽorme d'existence du corps de classe).

\bigskip

\th
\medskip

{\sl Sous les mmes hypothses que prŽcŽdemment, alors on a

\art{1)} $\ker(\Phi_{L/K})=\ker(\Phi_{E/K})$ 

\art{2)}$F=L$.

\art{3)}$$[I_K(\m):P_\m\cdot N_{E/K}(I_E(\widetilde{\m}))]=[G:G']\leq G=[E:K],$$ o $\widetilde{\m}$ est l'extension de $\m$ ˆ $E$; et l'inŽgalitŽ est Žvidemment stricte si $G$ n'est pas abŽlien.

On voit donc que la premire inŽgalitŽ du corps de classe reste vraie dans les extension non abŽlienne, mais en aucun cas l'ŽgalitŽ, ce qui rend impossible le thŽorme de rŽciprocitŽ d'Artin et donc, dans une certaine mesure la possibilitŽ par cette mŽthode de bien cerner les extension non abŽliennes.

}

{\bf Preuve}

Prouvons 1). L'application de restriction $G\to{\rm Gal}(L/K)$ induit un isomorphisme  $T\, :\, G/G'\to {\rm Gal}(L/K)$ et $\Phi_{L/K}=T\circ\Phi_{E/K}$. En effet, si $\P\in I_{K}(\m)$, $\gP_{0}$ dans  $L$ au-dessus de $\P$ et $\gP$ dans $E$ au-dessus de $\gP_{0}$. On sait que ${\rm Frob}_{L/K}(\gP_{0})={\rm Frob}_{E/K}(\gP)|_{L}$, donc $\Phi_{L/K}(\P)=T\circ\Phi_{E/K}(\P)$ et on conclut par multiplicativitŽ. Cela prouve en particulier que $\ker(\Phi_{L/K})=\ker(\Phi_{E/K})$ .

Prouvons  2). Soit $\n$ un $K$-module admissible pour $F/K$ et multiple de $\m$. Alors on a $N_{E/K}(I_{E}(\widetilde{\n}))\subset N_{E/K}(I_{E}(\widetilde{\m}))\subset N_{E/K}(I_{E}(\widetilde{\m}))\cdot P_{\m}=H(\m)$, o $\widetilde{\m}$ est le $E$-module qui prolonge $\m$. Donc
$$N_{E/K}(I_{E}(\widetilde{\n}))\subset H(\m)\cap I_{K}(\n)\buildrel {\rm \n\ multiple\ de\ \m} \over=H(\n)\buildrel {\rm \n\ admissible\ pour\ }F/K\over =P_{\n}\cdot N_{F/K}(I_{F}(\widetilde{\n}'))$$
o $\widetilde{\n}'$ est le $F$-module qui prolonge $\n$. Cela montre, en vertu du Corollaire \argdc\ que $F\subset E$. Ainsi, par hypothse,  $\m$ est admissible pour $F/K$, ce qui veut dire que 
$$P_{\m}\cdot N_{E/K}(I_{E}(\widetilde{\m}))=H(\m)=P_{\m}\cdot N_{F/K}(I_{F}(\widetilde{\m}'))=\ker(\Phi_{F/K}| I_{K}(\m)) =\ker(\Phi_{E/K}| I_{K}(\m)),\eqno{(*)}$$
la dernire ŽgalitŽ se montrant comme lors de la partie 1), puisque $F/K$ est une sous-extension (abŽlienne) de $E/K$.  Soit maintenant $K\subset F\subset L_{1}\subset E$, avec $L_{1}/K$ abŽlienne. Alors $N_{E/K}(I_{E}(\widetilde{\m}))\subset N_{L_{1}/K}(I_{L_{1}}(\widetilde{\m}''))\subset  N_{F/K}(I_{F}(\widetilde{\m}'))$, o $\widetilde{\m}''$ est le $L_{1}$-module qui prolonge $\m$. Ainsi,
$$P_{\m}\cdot N_{E/K}(I_{E}(\widetilde{\m}))\subset P_{\m}\cdot N_{L_{1}/K}(I_{L_{1}}(\widetilde{\m}''))\subset P_{\m}\cdot  N_{F/K}(I_{F}(\widetilde{\m}')).$$

Or, la relation $(*)$ nous montre que le premier et le troisime groupe sont Žgaux. Donc il y a ŽgalitŽ partout. En particulier, 

$$N_{F/K}(I_{F}(\widetilde{\m}'))\subset P_{\m}\cdot N_{L_{1}/K}(I_{L_{1}}(\widetilde{\m}'')),$$
ce qui prouve, en vertu du Corollaire \argdc\ et puisque par hypothse $\m$ est admissible pour $L_{1}/K$ que $L_{1}\subset F$. Ainsi, $L_{1}=F$, ce qui veut dire que $F$ est la sous-extension abŽlienne maximale de $E/K$, et donc que $L=F$.

Prouver 3) est alors un jeu d'enfant~:  la partie 2) et l'Žquation $(*)$ nous montrent que $\ker(\Phi_{L/K}| I_{K}(\m))=P_{\m}\cdot N_{E/K}(I_{E}(\widetilde{\m}))$. Ainsi,

$$I_{K}(\m)/P_{\m}\cdot N_{E/K}(I_{E}(\widetilde{\m}))\simeq G/G'=G^{\rm ab}$$
ou encore,
$$[I_{K}(\m):P_{\m}\cdot N_{E/K}(I_{E}(\widetilde{\m}))]=[G:G']\leq |G|,$$
 et c'est ce qu'il fallait dŽmontrer.\qed

\vfill\eject
\global\advance\chapnomb by 1
\nomb=1

\centerline{\para Chapitre 11}
\medskip
\centerline{\para Symbole de restes normiques, conducteur }
\medskip
\centerline{\para et corps de classe de Hilbert }
\bigskip

Dans ce chapitre nous allons montrer que le conducteur d'une extension abŽlienne est admissible, donc il correspond ˆ un sous-groupe de congruence qui est le noyau de l'application d'Artin d'une extension. Cela impliquera l'existence du corps de Hilbert (que nous construisons ˆ la fin de ce chapitre).

Pour parvenir ˆ la preuve de tout cela, nous allons dŽfinir et donner quelques propriŽtŽs d'une application importante~: l'application $\theta_\P$ qu'on nomme {\it symbole de restes normiques} qui interviendra ˆ la fin du chapitre 13 (DŽfinition \arggq) et qui sera crucial dans les rŽsultats du corps de classe local (Proposition-DŽfinition \arggx).

Fixons pour ce chapitre $L/K$ une extension abŽlienne de corps de nombres de groupe $G$, et posons $\ff=\ff(L/K)$, le conducteur de cette extension.

\headline={\hfill \phantom{ouuh}\hfill}
Commenons par donner un rŽsultat qui a dŽjˆ ŽtŽ partiellement prouvŽ (Lemme \argcr )
\bigskip

{{{\soustitre ThŽorme ({\the\chapnomb}.{\the\nomb})}\global\advance\nomb by 1}\soustitre (Lemme de translation)}
\medskip
{\sl Sous les mmes hypothses, considŽrons $E/K$ une extension quelconque finie de $K$. Alors, $\bbH(EL/E)$ est l'extension ˆ $E$ par les normes de $\bbH(L/K)$. Plus prŽcisŽment, si $\m$ est un $K$-module multiple  de $\ff(L/K)$, alors l'extension $\widetilde{\m}$ de $\m$ ˆ $L$ est un multiple de $\ff(EL/E)$, et on a la relation suivante qui lie les groupes de congruences~:
$$H(\widetilde{\m})=\{\aa\in I_E(\widetilde{\m})\mid N_{E/K}(\aa)\in H(\m)\},$$
o $H(\widetilde{\m})=H(\widetilde{\m},EL/E)$ est le sous-groupe de congruence pour $\widetilde{\m}$ de la classe $\bbH(EL/E)$ et  $H(\m)=H(\m,L/K)$ est le sous-groupe de congruence pour ${\m}$ de la classe $\bbH(L/K)$.

}

{\bf Preuve}

Soit $\m$ un $K$-module admissible pour $L/K$ et $\widetilde{\m}$ l'extension de $\m$ ˆ $L$. On sait (cf. ThŽorme \argf ) que $\Phi_{EL/E}$ est dŽfinie sur $I_E(\widetilde{\m})$ et que 
$R\circ \Phi_{EL/E}=\Phi_{L/K}\circ N_{E/K}$, o $R\, :\, {\rm Gal}(EL/E)\to G$ est l'injection donnŽe par la restriction ˆ $L$ des $E$-automorphismes de $EL$. Le Lemme \argcr\ nous dit que $\widetilde{\m}$ est admissible pour $EL/E$. Ainsi, 
$$\eqalign{H(\widetilde{\m})&=\ker(\Phi_{EL/E}\mid I_E(\widetilde{\m}))=\{\aa\in I_E(\widetilde{\m})\mid N_{E/K}(\aa)\in\ker\left (\Phi_{L/K}\mid I_K(\m)\right )\}=\cr
&=\{\aa\in I_E(\widetilde{\m})\mid N_{E/K}(\aa)\in H(\m)\}.\cr}$$
Et donc le rŽsultat est prouvŽ pour les $\m$ admissibles. Supposons maintenant que $\m$ soit un $K$-module tel que $\ff |\m$ et posons $A_{\widetilde{\m}}=\{\aa\in I_E(\widetilde{\m})\mid N_{E/K}(\aa)\in H(\m).\}$. On vient de voir que $A_{\widetilde{\m}}=H({\widetilde{\m}})$ si $\m$ est admissible. Or, on montre trs facilement que si  $\ff |\m$, $A_{\widetilde{\m}}=A_{\widetilde{\ff}}\cap I_E(\widetilde{\m})$. Ainsi $\{ A_{\widetilde{\m}}\mid \ff |\m\}$ est dans une mme classe d'Žquivalence de sous-groupes de congruence. Donc, $\{ A_{\widetilde{\m}}\mid \ff |\m\}\subset \bbH(EL/E)$, puisque c'est le cas pour les $\m$ admissibles. On en dŽduit que $\ff(EL/E)|\widetilde{\m}$ et, par unicitŽ du groupe de congruence pour $\m$,  $A_{\widetilde{\m}}=H(\widetilde{\m})$ pour tout $\m |\ff$.\qed

\bigskip\goodbreak

\defi
\medskip

Un $K$-module $\n$ est dit {\it $\P$-admissible} s'il est admissible (pour $L/K$ que nous avons supposŽe abŽlienne de groupe $G$) et si $\P |\n$. On peut alors Žcrire $\n=\P^a\cdot \m$, avec $a\geq 1$ et $\P\notdiv \m$. Supposons donc $\n$ $\P$-admissible. On appelle $\Theta$ $(=\Theta(L/K,\P^a,\m))$, l'homomorphisme obtenu par composition des applications 
$$\thboxed 25 {\Theta\, :\, K^*_\m\buildrel \iota\over\to I_K\buildrel j\over \to I_K(\n)\buildrel\Phi_{L/K}\mid I_K(\n)\over\longrightarrow G},$$
 \newcount\gaagan\gaagan=\gaga
o $\iota(x)=x\cdot O_K$ et $j=j_\n$ est l'homomorphisme dŽfini sur les idŽaux premiers $\euf q$ de la manire suivante~: $j({\euf q})=\cases{{\euf q}&si ${\euf q}\notdiv \n$\cr O_K&si $\euf q |\n$\cr}$.
Rappelons que $Z(\P)=Z(L/\P)=\{\sigma\in G\mid \sigma(\gP)=\gP\}$ o $\gP$ est n'importe quel idŽal premier de $L$ au-dessus de $\P$.

\bigskip\goodbreak

\th
\medskip

{\sl Sous les mmes hypothses que pour la dŽfinition prŽcŽdente, on a

$${\rm Im} (\Theta)=Z(\P).$$

}

\headline={\hfill \smcap Symbole de restes normiques, conducteur  et corps de classe de Hilbert \hfill}
{\bf Preuve}

Montrons $Ò\subset "$~:  Posons $F$ le sous-corps de $L$ fixe par $Z(\P)$. On se souvient que $F$ est la plus grande sous-extension de $L$ dans laquelle $\P$ se dŽcompose compltement (cf. [Mar, Thm. 29 (1), p.104]) (c'est aussi valable pour les places infinies, mais dans ce cas dŽcomposer compltement veut simplement dire ne pas ramifier). En particulier, $\P$ ne ramifie pas dans $F$. Cela implique, en vertu du Lemme \argdl\ que $\P\notdiv \ff(F/K)$. D'autre part, puisque $F\subset L$, on a $\ff(F/K) |\ff(L/K)$ (cf. Proposition \1argdm). Ainsi,  $\ff(F/K)|\m$. De plus, les places de $K$ qui ramifient dans $F$ divisent $\m$ (puisque  c'est une partie de celles qui ramifient dans $L$ et que $\P$ n'en fait pas partie). Donc, en vertu du Lemme \argdk, cela veut dire que $\m$ est admissible pour $F/K$ et donc (Lemme \argco )~:
$$\ker(\Phi_{F/K}|I_K(\m))=P_\m\cdot N_{F/K}(I_F(\widetilde{\m})),\eqno{(*)}$$
o $\widetilde{\m}$ est l'extension de $\m$ ˆ $F$. 

Supposons $\P$ finie. Puisque $\P$ dŽcompose compltement dans $F$ et $\P\notdiv\m$, alors $\P\in \ker (\Phi_{F/K}\mid I_K(\m))$. Soit maintenant $\alpha\in K^*_\m$. Par dŽfinition, on a $\iota(\alpha)=\P^s\cdot\aa$, o $s\in\Z$ et $\aa\in I_K(\n)$. Alors $j(\iota(\alpha))=\aa=\P^{-s}\cdot \iota(\alpha)\in \ker (\Phi_{F/K}\mid I_K(\m))$, car $\P$ est dans ce noyau et  $\iota(\alpha)\in P_\m\buildrel(*)\over \subset \ker(\Phi_{F/K}|I_K(\m))$. 

Supposons $\P$ infinie.  Dans ce cas,  $I_K(\m)=I_K(\n)$ et pour tout $\alpha\in K^*_\m$, on a $j(\iota(\alpha))=\iota(\alpha)\in P_\m\buildrel(*)\over \subset \ker(\Phi_{F/K}|I_K(\m))$. Donc dans les deux cas (fini et infini), on a
$$1_G=\Phi_{F/K}(j(\iota(\alpha)))=\Phi_{L/K}(j(\iota(\alpha)))|_F=\Theta(\alpha)|_F.$$
Cela montre que $\Theta(\alpha)\in Z(\P)$, par la thŽorie de Galois.  Ainsi, on a prouvŽ que ${\rm Im} (\Theta)\subset Z(\P)$.

La preuve de l'inclusion inverse $Ò\supset"$ se fait par l'absurde. Supposons donc que ${\rm Im} (\Theta)\subsetneqq  Z(\P)$. Un raisonnement facile sur les groupes abŽliens finis nous assure l'existence d'un sous-groupe $G_0$ tel que  ${\rm Im}(\Theta)\subset G_0\subset Z(\P)$, avec $[Z(\P):G_0]=q\in \bbP(\Q)$. Posons $E$ le sous-corps de $L$ fixe par $G_0$ et comme pour la partie $``\subset"$, $F$ le sous-corps de $L$ fixe par $Z(\P)$. On a donc $K\subset F\subset E\subset L$, avec $[E:F]=q$. Soit $\zeta$ une racine primitive $q$-ime de l'unitŽ. On pose $F'=F(\zeta)$ et $E'=E(\zeta)$. On est donc dans la situation~:
\vglue 2cm

\rput(6,0){\rnode{K}{$K$}}
\rput(7.8,0){\rnode{F}{$F$}}
\rput(9.6,0){\rnode{E}{$E$}}
\psset{nodesep=3pt} 
\nccurve[ncurv=1.2,angleA=-100,angleB=-85]{K}{E}
\Bput{$G/G_{0}$}
\nccurve[ncurv=0.4,angleA=-90,angleB=-90]{K}{F}
\Bput{$G/Z(\P)$}
\nccurve[ncurv=0.4,angleA=-90,angleB=-90]{F}{E}
\Bput{$Z(\P)/G_{0}$}
\rput(6.9,0){\rnode{S1}{$\subset$}}
\rput(8.7,0){\rnode{S2}{$\subset$}}
\rput(7.5,1.4){\rnode{F'}{$F'=F(\zeta)$}}
\rput(7.8,0.7){\rnode{S3}{$\cup$}}
\rput(10.2,1.4){\rnode{E'}{$E'=EF'$}}
\rput(9.6,0.7){\rnode{S4}{$\cup$}}
\rput(8.7,1.4){\rnode{S5}{$\subset$}}
\rput(10.5,0){\rnode{S6}{$\subset$}}
\rput(11.4,0){\rnode{L}{$L$}}
\nccurve[ncurv=0.4,angleA=-90,angleB=-90]{E}{L}
\Bput{$G_{0}$}
\vglue 3cm

On veut appliquer le ThŽorme \argeg , les r™les de $K$ et de $n$ Žtant tenus ici par $F'$ et $q$ respectivement. De plus, posons $S_1$, l'ensemble des places de $F'$ au-dessus de $\P$ et $S_2$ un ensemble fini de places de $F'$ disjoint de $S_{1}$ tel que $S=S_1\cup S_2$ vŽrifie les conditions $i)-iii)$ du ThŽorme \argeg. Prenons encore $\m_1$, l'extension ˆ $F'$ de $\P^b$, avec $b\geq 1$ assez grand si $\P$ est finie; et $\m_2$ le produit de toutes les places non complexes de $S_2$, avec des exposants assez grands pour les places finies (toujours pour appliquer le mme ThŽorme \argeg).  On supposera de plus que $\widetilde{\m}|\m_2$ et que $\widetilde{\n}|\m_1\cdot\m_2$, o $\widetilde{\m}$ et $\widetilde{\n}$ sont les extensions ˆ $F'$ de $\m$ et de $\n$. 

ConsidŽrons alors le $H_2$ et le $L_2$ du mme ThŽorme \argeg. Rappelons que
$$H_2=I_{F'}(\m_2)^q\cdot P_{F',\m_2}\cdot I_{F'}[S_1],\hbox{ et que $L_2$ est une $q$-extension de Kummer de $F'$.}$$ 
Le mme ThŽorme \argeg\ nous apprend que $H_2=H(\m_{2},L_{2}/F')$ est le sous-groupe de congruence dans $\bbH(L_2/F')$ dŽfini modulo $\m_2$ et que $\m_2$ est admissible pour l'extension $L_2/F'$. On va appliquer le lemme de translation (ThŽorme \argei) au diagramme~:

\input pstricks
\input pst-node

 \vglue 2.5cm

\psset{xunit=0.8cm,yunit=0.8cm}

\rput(7,0){\rnode{K}{$K$}}
\rput(7,2){\rnode{E}{$E$}}
\rput(9,1){\rnode{ff}{$F'$}}
\rput(9,3){\rnode{ee}{$E'$}}
\ncline[nodesep=3pt]{K}{E}
\ncline[nodesep=3pt]{K}{ff}
\ncline[nodesep=3pt]{E}{ee}
\rput(10,3){\rnode{eee}{$=EF'$}}
\ncline[nodesep=3pt]{ff}{ee}

\medskip

Puisque $\n$ est admissible pour $L/K$ alors, en vertu du Lemme \argcq , $\n$ est admissible pour $E/K$. Par le lemme de translation et sa preuve, l'extension $\widetilde{\n}$ de $\n$ ˆ $F'$ est aussi admissible pour $E'/F'$ et (rappelons-le)
$$H(\widetilde{\n},E'/F')=\{\aa\in I_{F'}(\widetilde{\n})\mid N_{F'/K}(\aa)\in H(\n,E/K)\},\eqno{(**)}$$
o $H(\widetilde{\n},E'/F')$ est le sous-groupe de congruence pour $\widetilde{\n}$ de la classe $\bbH(E'/F')$ et  $H(\n,E/K)$ est le sous-groupe de congruence pour ${\n}$ de la classe $\bbH(E/K)$.

D'autre part, par hypothse (absurde), $\Theta(K^*_{\m})\subset G_{0}$, c'est-ˆ-dire que $\Phi_{L/K}| I_{K}(\n)$ envoie $j(\iota(K^*_{\m}))$ dans $G_{0}$, qui est le sous-groupe de $G$ formŽ des ŽlŽments qui sont l'identitŽ sur $E$. Ainsi,
$$j(\iota(K^*_{\m}))\subset \ker (\Phi_{E/K}|I_{K}(\n))=H(\n,E/K).\eqno{(***)}.$$

On va dŽmontrer que 

$$H_{2}\cap I_{F'}(\widetilde{\n})\subset H(\widetilde{\n},E'/F').\eqno{(\dag)}$$

ConsidŽrons donc, ${\euf X}\in H_{2}$, un idŽal premier ˆ $\widetilde{\n}$. Par dŽfinition de $H_{2}$, on a ${\euf X}={\cal B}^q\cdot (\alpha)\cdot {\cal C}$, avec ${\cal B}\in I_{F'}(\m_{2})$, $\alpha\in {F'}^*_{\m_{2}}$ et ${\cal C}\in I_{F'}[S_{1}]$ (cf. Chapitres 9 et 10 pour la dŽfinition de ces objets) . Il faut donc montrer (en vertu de $(**)$ ) que $N_{F'/K}({\euf X})\in H(\n,E/K)$.

Supposons $\P$ fini. On peut Žcrire ${\cal B}={\cal B}_{0}\cdot{\cal B}_{1}$ avec ${\cal B}_{0}$ premier ˆ $\widetilde{\n}$ et $N_{F'/K}({\cal B}_{1})=\P^r$, $\iota(\alpha)=(\alpha)=\aa_{0}\cdot\aa_{1}$, avec $\aa_0$ premier ˆ $\widetilde{\n}$ et $N_{F'/K}(\aa_{1})=\P^s$ et, par dŽfinition de $S_{1}$, $N_{F'/K}({\cal C})=\P^t$. On a $r\cdot q+s+t=0$, car ${\euf X}\in I_{F'}(\widetilde{\n})$, et donc $N_{F'/K}({\euf X})$ est premier ˆ $\P$. Cela implique que $N_{F'/K}({\euf X})=N_{F'/K}({\cal B}_{0})^q\cdot N_{F'/K}(\aa_{0})$. 

On sait que $H(\n,E/K)\subset H(\n,F/K)$, o $H(\n,F/K)$ est le sous-groupe de congruence pour ${\n}$ de la classe $\bbH(F/K)$ (cf. Proposition \1argdm ). L'indice de $H(\n,E/K)$ dans $H(\n,F/K)$ vaut $q$. En effet, 
puisque $\n$ est admissible pour $F/K$ et $E/K$, 
$$\eqalign{H(\n,F/K)&=\ker(\Phi_{F/K}|I_{K}(\n))=\{\aa\in I_{K}(\n)\mid \Phi_{F/K}(\aa)=\Phi_{E/K}(\aa)|_{F}={\rm Id}_{F}\}\cr &=\{\aa\in I_{K}(\n)\mid \Phi_{E/K}(\aa)\in {\rm Gal}(E/F)\simeq Z(\P)/G_{0}\}.\cr}$$
De mme, $H(\n,E/K)=\ker(\Phi_{E/K}|I_{K}(\n))=\{\aa\in I_{K}(\n)\mid \Phi_{E/K}(\aa)={\rm Id}_{E}\}$. On en dŽduit que $[H(\n,F/K):H(\n,E/K)]=[Z(\P):G_{0}]=q$, car l'application d'Artin $\Phi_{E/K}$ est surjective (cf. ThŽorme \argaq).

Par dŽfinition de ${\cal B}_{0}$, on a $N_{F'/F}({\cal B}_{0})\in I_{F}(\widetilde{\n}')$, o $\widetilde{\n}'$ est l'extension de $\n$ ˆ $F$. Puisque $\n$ est admissible pour $L/K$, il l'est aussi pour $F/K$ (toujours en vertu du Lemme \argcq ), et alors on a~:
$$N_{F'/F}({\cal B}_{0})=N_{F/K}(N_{F'/F}({\cal B}_{0})\in N_{F/K}(I_{F}(\widetilde{\n}'))\subset P_{\n}\cdot N_{F/K}(I_{F}(\widetilde{\n}'))=H(\n,F/K).$$
Cela implique, puisque $[H(\n,F/K):H(\n,E/K)]=q$, que $N_{F'/F}({\cal B}_{0})^q\in H(\n,E/K)$. D'autre part,  puisque $\widetilde{\m}|\m_{2}$, 
$$N_{F'/K}(\alpha)\in N_{F'/K}({F'}^*_{\m_{2}})\subset N_{F'/K}({F'}^*_{\widetilde{\m}})\subset K^*_{\m}.$$
La dernire inclusion provient de la partie b) du Lemme \aargbu. D'autre part, on a~:
$$j(\iota(N_{F'/K}(\alpha))=j(N_{F'/K}(\aa_{0}))\cdot \overbrace{j(\underbrace{N_{F'/K}(\aa_{1})}_{=\P^s})}^{=O_{K}}=N_{F'/K}(\aa_{0}).$$
Ainsi,
$$N_{F'/K}(\aa_{0})\in j(\iota(K^*_{\m})\buildrel (***)\over\subset H(\n,E/K).$$
Ce qui montre que $N_{F'/K}({\euf X})\in H(\n,E/K)$, prouvant la relation $(\dag)$ dans le cas o $\P$ est fini. Si $\P$ est infini, on refait exactement le mme calcul, avec quelques simplifications~: dans ce cas, ${\cal C}=\aa_{1}={\cal B}_{1}=O_{K}$.

Soient $H(\m_{1}\m_{2},L_{2}/F')$ le sous-groupe de congruence pour $\m_{1}\m_{2}$ de la classe $\bbH( L_{2}/F')$, $\ff(L_{2}/F')$, le conducteur de l'extension $L_{2}/F'$ et $H(\ff(L_{2}/F'))$ le sous-groupe de congruence associŽ ˆ $\ff(L_{2}/F')$. Calculons~:
$$\eqalign{H(\m_{1}\m_{2},L_{2}/F')&=H(\ff(L_{2}/K))\cap I_{F'}(\m_{1}\m_{2})\cr
&=H(\ff(L_{2}/K))\cap  I_{F'}(\m_{2})\cap  I_{F'}(\widetilde{\n})\cap I_{F'}(\m_{1}\m_{2})\cr
&=H_{2}\cap I_{F'}(\widetilde{\n})\cap I_{F'}(\m_{1}\m_{2})\cr
&\buildrel (\dag)\over\subset H(\widetilde{\n},E'/F')\cap  I_{F'}(\m_{1}\m_{2})=H(\m_{1}\m_{2},E'/F'). \cr}$$

Cela implique, en vertu de la Proposition \1argdm, que $\thboxed15{E'\subset L_{2}}$.

D'aprs la partie $(I)$ du ThŽorme \argeg, les places de $S_{1}$ sont compltement dŽcomposŽes dans $L_{2}$, donc aussi dans $E'$. Soit $\P_{F}$, une place de $F$ au-dessus de $\P$, $\gP$ une place de $F'$ au-dessus de $\P_{F}$, $\gP'$, une place $E'$ au-dessus de $\gP$ et $\gP_{E}=\gP'\cap E$. On a donc la situation~:
\vglue 1cm

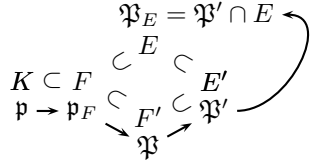
\rput(5,0){\rnode{K}{$K$}}
\rput(5.5,0){\rnode{S1}{$\subset$}}
\rput(6,0){\rnode{F}{$F$}}
\rput[b]{30}(6.7,0.2){$\subset$}
\rput[b]{-30}(6.5,-0.4){$\subset$}
\rput(7.1,0.6){\rnode{E}{$E$}}
\rput(7.1,-0.6){\rnode{F'}{$F'$}}
\rput[b]{30}(7.7,-0.5){$\subset$}
\rput(8.2,0){\rnode{E'}{$E'$}}
\rput[b]{-30}(7.6,0.2){$\subset$}
\rput(8.2,0){\rnode{E'}{$E'$}}
\rput(5,0){\rnode{K}{$K$}}
\rput(5,-0.5){\rnode{P}{$\P$}}
\rput(6,-0.5){\rnode{PF}{$\P_{F}$}}
\rput(7.1,-1.1){\rnode{gP}{$\gP$}}
\rput(8.2,-0.5){\rnode{gP'}{$\gP'$}}
\rput(7.9,1.1){\rnode{PE}{$\gP_{E}=\gP'\cap E$}}
\ncline[nodesep=3pt]{->}{P}{PF}
\ncline[nodesep=3pt]{->}{PF}{gP}
\ncline[nodesep=3pt]{->}{gP}{gP'}
\nccurve[ncurv=1,angleA=-0,angleB=0]{->}{gP'}{PE}

\vglue 1.5cm

Comme $\gP\in S_{1}$, on a l'indice de ramification $e(\gP'/\gP)=1$ et le degrŽ de $\gP'$ sur $F'$ $f(\gP'/\gP)=1$. Ainsi, 
$$e(\gP'/\gP_{E})\cdot e(\gP_{E}/\P_{F})=e(\gP'/\P_{F})=e(\gP'/\gP)\cdot e(\gP/\P_{F})=e(\gP/\P_{F}).$$
Mais, $[F':F]\mid[\Q(\zeta):\Q]=q-1$ (thŽorie de Galois), donc $e(\gP'/\gP_{E})\cdot e(\gP_{E}/\P_{F})\mid q-1$. D'autre part, $[E:F]=q$, donc $e(\gP_{E}/\P_{F})\mid q$. Cela implique que $e(\gP_{E}/\P_{F})=1$. Le mme calcul s'applique aux $f$, donc $f(\gP_{E}/\P_{F})=1$. Mais alors cela veut dire que $\P$ est compltement dŽcomposŽ dans $E$, ce qui contredit le fait que $F$ est le plus grand sous-corps de $L$ dans lequel $\P$ se dŽcompose compltement. Cette contradiction implique que ${\rm Im}(\Theta)\supset Z(\P)$ et achve la (longue) preuve de ce thŽorme.\qed

\bigskip\goodbreak

\defi
\medskip

Soit $K$ un corps de nombres. On rappelle que $\bbK_{\P}$ est le localisŽ complŽtŽ de $K$ en la place $\P$ et $\bbK_{\P}^*=\bbK_{\P}\setminus\{0\}$. On note aussi $U_{\P}=\okp^*$, le groupe des unitŽs de l'anneau des entiers de $\bbK_{\P}$ et $\widehat{\P}=\P\okp$, l'idŽal maximal de $\okp$.

Soit $b\geq 0$, on dŽfinit 
$$U_{\P}^{(b)}=\cases{U_{\P}&si $b=0$ et $\P$ fini \cr 1+\widehat{\P}^b&si $b>0$ et $\P$ fini\cr
U_{\P}=\bbK^*_{\P}=\R^*&si $b=0$ et $\P$ infini rŽel\cr \bbK_{\P,+}^*=\R_{+}^*& si $b>0$ et $\P$ infini rŽel\cr
U_{\P}=\bbK^*_{\P}=\C^*&si $\P$ est infinie complexe.\cr}$$

On avait dŽjˆ dŽfini cet objet au Chapitre 5 (DŽfinition \aargby), mais ici la dŽfinition est Žtendue dans le cas des places infinies et si $b=0$. Il est Žvident que si $b'\geq b$, alors $U_{\P}^{(b')}\subset  U_{\P}^{(b')}$. Il s'agit aussi d'Žtendre la dŽfinition de $\pmodast$. On l'avait dŽjˆ fait partiellement ˆ la DŽfinition \aargby, mais nous allons ici encore un peu plus loin~: si $x,y\in\bbK_{\P}^*$, on Žcrit  $x\equiv y\pmodast {\widehat{\P}^b}$ si ${x-y\over y}\in \widehat{\P}^b$ pour les places finies et $b>0$. Si $b=0$ ou $\P$ complexe, cela veut dire que ${x\over y}\in U_{\P}$ et enfin si $b>0$ dans le cas des places infinie rŽelle, cela veut dire que $x$ et $y$ ont le mme signe. Autrement dit, et pour faire court, si $x,y\in\bbK_{\P}^*$, $\P\in \gfP(K)$ et $b\in\N$, alors on a~:
$$x\equiv y\pmodast{\widehat{\P}^b}\iff {x\over y}\in U_{\P}^{(b)}.$$

\bigskip\goodbreak

\prop
\medskip
{\sl Soit $K$ un corps de nombres, $\m$ un $K$-module, $\P\in\gfP(K)\setminus \gfP_{\C}(K)$ avec $\P\notdiv \m$ et $b\in\N$. Alors la composition d'homomorphisme
$$K_{\m}^*\hra \bbK_{\P}^*\longrightarrow \hskip-0.5cm\longrightarrow \bbK_{\P}^*/U_{\P}^{(b)} ,$$
est surjective (o la premire application est l'inclusion et la seconde, la projection canonique). Elle induit un isomorphisme
$$K_{\m}^*/(K_{\m}^*\cap U_{\P}^{(b)})\lra \bbK_{\P}^*/U_{\P}^{(b)}.$$

}

{\bf Preuve}

Soit $x\in \bbK^*_{\P}$. Par densitŽ, il existe $y\in K^*$ tel que $x\equiv y\pmodast {\widehat{\P}^b}$. Le thŽorme d'approximation dŽbile (ThŽorme \argc) nous assure l'existence d'un ŽlŽment $z\in K^*$ tel que $z\equiv 1\pmodast \m$ et $z\equiv y\pmodast {\P^b}$. La deuxime Žquivalence veut dire que $z\in K_{\m}^*$ et le deux autres veulent dire que la classe de $z$ modulo $U_{\P}^{(b)}$ est la mme que la classe de $x$ modulo $U_{\P}^{(b)}$, ce qui montre la surjectivitŽ. Enfin, le noyau de cet homomorphisme est Žvidemment  $K_{\m}^*\cap U_{\P}^{(b)}$.\qed

\bigskip
\defi
\medskip
Soit $K$ un corps de nombres, $\P$ une place non-complexe et $\n=\P^{a}\cdot\m$ un $K$-module $\P$-admissible ($\P\notdiv\m$). Supposons que $b\geq 0$ est tel que $K_{\m}^*\cap U_{\P}^{(b)}\subset \ker(\Theta(L/K,\P^{a},\m))$, alors on dŽfinit $\theta_{\P}(L/K,\P^{a},\m,b)\, :\, \bbK_{\P}^*\to Z(\P)$, l'homomorphisme composŽ
$$\thboxed 25 {\theta_{\P}\, :\,\bbK_{\P}^*\lra \bbK_{\P}^*/U_{\P}^{(b)}\buildrel\simeq\over\lra K_{\m}^*/(K_{\m}^*\cap U_{\P}^{(b)})\buildrel\overline{\Theta}\over \lra Z(\P)}$$
 \newcount\gaagao\gaagao=\gaga
o $\overline{\Theta}=\overline{\Theta}(L/K,\P^{a},\m)$ est dŽfini par $\Theta(L/K,\P^{a},\m)$.

Si $\P$ est une place complexe, on dŽfinit $\theta_{\P}\, :\, \bbK_{\P}\to G$ l'homomorphisme trivial.
\bigskip\goodbreak
{\soustitre Remarques}
\medskip
\art{a)}On vŽrifie facilement que 
$$K_{\m}^*\cap U_{\P}^{(b)}=\cases{K^*_{\P^b\cdot\m}&si $b>0$\cr
K^*_{\m}(\P)&si $b=0$ et $\P$ finie\cr
K^*_{\m}&si $b=0$ et $\P$ infinie\cr}$$
o $K^*_{\m}(\P)=\{x\in K^*_{\m}\mid v_{\P}(x)=0\}$.

\art{b)}Dans la dŽfinition prŽcŽdente, si $\n=\P^{a}\cdot\m$, on peut prendre $b=a$, car $K_{\m}^*\cap U_{\P}^{(a)}=K^*_{\P^{a}\cdot\m}=K_{\n}^*$ est dans le noyau de $\Theta$ (cf. ThŽorme \argco\ et $j(\iota(K_{\n}^*))=P_{\n}$).

\art{c)}Les anglophones appellent l'applications $\theta_\P$ le ``norm-residue symbol".

\bigskip
\prop
\medskip
{\sl Sous les mmes hypothses que pour la dŽfinition prŽcŽdente, alors $\theta_{\P}(L/K,\P^{a},\m,b)$ est indŽpendant du choix du $K$-module $\P$-admissible $\P^{a}\cdot\m$ et de l'entier $b\geq 0$ tel que  $K_{\m}^*\cap U_{\P}^{(b)}\subset \ker(\Theta(L/K,\P^{a},\m))$.

}

{\bf Preuve}

D'abord, si $\n$ est fixŽ, $\Theta$ aussi et $\theta_{\P}$ est indŽpendant du choix de $b$. En effet, supposons que $U_{\P}^{(b')}\cap K_{\m}^*\subset\ker(\Theta)$. Sans limiter la gŽnŽralitŽ, on peut supposer que $b<b'$ et on a clairement le diagramme commutatif~:
\vglue 3cm
\psset{xunit=0.9cm,yunit=0.9cm}
\rput(4,2){\rnode{kp}{$\bbK_{\P}^* $}} 
\rput(7,3){\rnode{kpb}{$\bbK_{\P}^*/U_{\P}^{(b)}$}}
\rput(10,3){\rnode{kmb}{$K^*_{\m}/(K_{\m}^*\cap U_{\P}^{(b)})$}}
\rput(13,2){\rnode{zp}{$Z(\P)$}}
\rput(10,1){\rnode{kmbpr}{$K^*_{\m}/(K_{\m}^*\cap U_{\P}^{(b')})$}}
\rput(7,1){\rnode{kpbpr}{$\bbK_{\P}^*/U_{\P}^{(b')}$}}
\ncline[nodesep=3pt]{->}{kp}{kpb}
\Aput{$\rm nat. $} 
\ncline[nodesep=3pt]{->}{kp}{kpbpr}
\Bput{$\rm nat.$} 
\ncline[nodesep=3pt]{->}{kpb}{kmb}
\Aput{$\simeq $} 
\ncline[border=2pt, nodesep=3pt]{->}{kmb}{zp}
\Aput{$\overline{\Theta}_{b}$} 
\ncline[nodesep=3pt]{->}{kpbpr}{kmbpr}
\Aput{$\simeq $} 
\ncline[border=2pt , nodesep=3pt]{->}{kmbpr}{zp}
\Bput{$\overline{\Theta}_{b'}$} 
\ncline[border=2pt , nodesep=3pt]{->}{kpbpr}{kpb}
 \Aput{$\rm nat. $} 
 \ncline[border=2pt , nodesep=3pt]{->}{kmbpr}{kmb}
 \Aput{$\rm nat. $}

Soit maintenant $\n=\P^{a}\cdot\m$ et $\n'=\P^{a'}\cdot \m'$ deux $K$-modules $\P$-admissibles tels que $\n |\n'$. Donc $\m |\m'$ et $a\leq a'$, et on peut prendre $b=a$ et $b'=a'$. On a un diagramme commutatif~:

\vglue 3cm
\psset{xunit=0.9cm,yunit=0.9cm}
\rput(4,3){\rnode{kn}{$K^*_{\n}$}} 
\rput(4,1){\rnode{knpr}{$K^*_{\n'}$}} 
\rput(7,3){\rnode{km}{$K^*_{\m}$}}
\rput(10,3){\rnode{ikm}{$I_{K}(\m)$}}
\rput(13,3){\rnode{ikn}{$I_{K}(\n)$}}
\rput(10,1){\rnode{ikmpr}{$I_{K}(\m')$}}
\rput(7,1){\rnode{kmpr}{$K_{\m'}$}}
\rput(13,1){\rnode{iknpr}{$I_{K}(\n')$}}
\rput(15,2){\rnode{g}{$G$}}
\ncline[nodesep=3pt]{->}{knpr}{kn}
\Aput{$\rm incl. $} 
\ncline[nodesep=3pt]{->}{kn}{km}
\Bput{$\rm incl.$} 
\ncline[nodesep=3pt]{->}{knpr}{kmpr}
\Aput{$\rm incl. $} 
\ncline[nodesep=3pt]{->}{kmpr}{km}
\Aput{$\rm incl. $} 
\ncline[border=2pt, nodesep=3pt]{->}{km}{ikm}
\Aput{$\iota$}
\ncline[border=2pt, nodesep=3pt]{->}{kmpr}{ikmpr}
\Aput{$\iota$}  
\ncline[nodesep=3pt]{->}{ikmpr}{ikm}
\Aput{$\rm incl. $} 
\ncline[border=2pt , nodesep=3pt]{->}{ikm}{ikn}
\Aput{$j_{\n}$} 
\ncline[border=2pt , nodesep=3pt]{->}{ikmpr}{iknpr}
 \Aput{$j_{\n'} $} 
 \ncline[border=2pt , nodesep=3pt]{->}{iknpr}{ikn}
 \Aput{$\rm incl $} 
 \ncline[border=2pt , nodesep=3pt]{->}{ikn}{g}
 \Aput{$ \Phi_{L/K}$} 
 \ncline[border=2pt , nodesep=3pt]{->}{iknpr}{g}
 \Bput{$ \Phi_{L/K} $}

 Donc, en passant aux quotients, un diagramme commutatif~:
 
 \vglue 3cm
\psset{xunit=0.9cm,yunit=0.9cm}
\rput(5,3){\rnode{kmn}{$K^*_{\m}/K^*_{\n}$}} 
\rput(5,1){\rnode{kmnpr}{$K^*_{\m'}/K^*_{\n'}$}} 
\rput(8,2){\rnode{zp}{$Z(\P)$}}
\ncline[nodesep=3pt]{->}{kmnpr}{kmn}
\Aput{$\rm nat. $} 
\ncline[nodesep=3pt]{->}{kmn}{zp}
\Aput{$\rm \overline{\Theta}(\m,\P^a) $} 
\ncline[nodesep=3pt]{->}{kmnpr}{zp}
\Bput{$\rm \overline{\Theta}(\m',\P^{a'}) $} 

Que l'on complte en un diagramme toujours commutatif~:

\vglue 3cm
\psset{xunit=0.9cm,yunit=0.9cm}
\rput(4,2){\rnode{kp}{$\bbK_{\P}^* $}} 
\rput(7,3){\rnode{kpa}{$\bbK_{\P}^*/U_{\P}^{(a)}$}}
\rput(10,3){\rnode{kmn}{$K^*_{\m}/K^*_\n$}}
\rput(13,2){\rnode{zp}{$Z(\P)$}}
\rput(10,1){\rnode{kmnpr}{$K^*_{\m'}/K^*_{\n'}$}}
\rput(7,1){\rnode{kpapr}{$\bbK_{\P}^*/U_{\P}^{(a')}$}}
\ncline[nodesep=3pt]{->}{kp}{kpa}
\Aput{$\rm nat. $} 
\ncline[nodesep=3pt]{->}{kp}{kpapr}
\Bput{$\rm nat.$} 
\ncline[nodesep=3pt]{->}{kpa}{kmn}
\Aput{$\simeq $} 
\ncline[border=2pt, nodesep=3pt]{->}{kmn}{zp}
\Aput{$\overline{\Theta}(\m,\P^a)$} 
\ncline[nodesep=3pt]{->}{kpapr}{kmnpr}
\Aput{$\simeq $} 
\ncline[border=2pt , nodesep=3pt]{->}{kmnpr}{zp}
\Bput{$\overline{\Theta}(\m',\P^{a'})$} 
\ncline[border=2pt , nodesep=3pt]{->}{kpapr}{kpa}
 \Aput{$\rm nat. $} 

D'o l'indŽpendance en le $K$-module $\n$ (le cas de deux $K$-modules gŽnŽraux $\n$ et $\n'$ se rŽduit au cas prŽcŽdent en passant par l'intermŽdiaire du $\rm ppcm$ de $\n$ et $\n'$).\qed

\bigskip\goodbreak
{\soustitre Remarque}
\medskip
Pour calculer $\theta_\P$ lorsque $\P$ est un idŽal premier, on peut donc procŽder comme suit~: on choisit d'abord un $K$-module $\P$-admissible $\n=\P^a\cdot\m$ ($\P\notdiv\m$) et un $b\in \N$ tel que  $K_{\m}^*\cap U_{\P}^{(b)}\subset \ker(\Theta(L/K,\P^{a},\m))$ (par exemple $b=a$) et si $x\in \bbK_\P^*$, on choisit (par densitŽ et gr‰ce au thŽorme d'approximation dŽbile, ThŽorme \argc),  $y\in K^*$ tel que $y\equiv 1\pmodast{\m}$ et $y\equiv x\pmodast{\widehat{\P}^b}$, et alors, 
$$\theta_\P(x)=\Phi_{L/K}(j_\n((y)))=\Phi_{L/K}\left((y)\cdot\P^{-v_\P(y)}\right ),$$
o, $\Phi_{L/K}=\Phi_{L/K}|_{I_K(\n)}$. En particulier,

\bigskip
\prop
\medskip
{\sl Soit $L/K$ une extension abŽlienne de corps de nombres. Si $\P$ est un idŽal premier non-ramifiŽ dans $L$, alors 
$$\theta_\P(x)={\rm Frob}_{L/K}(\P)^{-v_\P(x)}.$$

} 
{\bf Preuve}

Puisque la Proposition \aargel\ est vraie, on peut prendre (en vertu du thŽorme de rŽciprocitŽ d'Artin, ThŽorme \argda) pour $\m$ un $K$-module admissible tel que $\P\notdiv \m$, et choisir $\n=\P\cdot \m$ et $a=b=1$. Soit $x\in \bbK_\P$. Comme dans la remarque prŽcŽdente, choisissons $y\in K_\m^*$ tel que $y\equiv x\pmodast{\widehat{\P}}$. Alors, l'idŽal $(y)\in P_\m$, et puisque $\m$ est admissible, $\Phi_{L/K}((y))=1$. Ainsi, en vertu de la remarque prŽcŽdente, 
$$\theta_\P(x)=\Phi_{L/K}\left ((y)\cdot\P^{-v_\P(y)}\right )=\Phi_{L/K}(\P)^{-v_\P(x)}.$$
\qed
\bigskip
\prop
\medskip
{\sl Soit $L/K$ une extension abŽlienne de corps de nombres. Soit $\P\in\gfP(K)$, $\n=\P^a\cdot\m$ un $K$-module $\P$-admissible et $c\geq 0$ un nombre entier. ConsidŽrons $\Theta=\Theta(L/K,\P^a,\m)$ et $\theta_\P$. ConsidŽrons enfin $\ff(L/K)$, le conducteur de $L/K$ (cf. Application-DŽfinition \argdj\ pour la dŽfinition). Alors les trois conditions sont Žquivalentes~:

\art{i)}$U_\P^{(c)}\subset\ker(\theta_\P)$,

\art{ii)}$U_\P^{(c)}\cap K^*_\m\subset\ker(\Theta)$,

\art{iii)}$v_\P(\ff(L/K))\leq c$.

}

{\bf preuve}

La preuve de  $\sl ii)\Rightarrow i) $ est Žvident, car on peut alors dŽfinir $\theta_\P$ avec $b=c$ dans ce cas.

Pour prouver $\sl i) \Rightarrow  ii)$, on peut supposer $c<a$, car si $c\geq a$, i) et ii) sont tous les deux vrais. Par $\theta_\P$, $U_\P^{(c)}$ est envoyŽ successivement sur (on utilise $b=a$ pour dŽfinir $\theta_\P$) $U_\P^{(c)}/U_\P^{(a)}\subset \bbK_\P/U_\P^{(a)}$ puis sur $(U_\P^{(c)}\cap K^*_\m)/(U_\P^{(a)}\cap K^*_\m)\subset K_\m^*/(U_\P^{(a)}\cap K^*_\m)$ et enfin sur $\Theta(U_\P^{(c)}\cap K^*_\m)$ qui doit tre 1 par hypothse. On en dŽduit bien que $U_\P^{(c)}\cap K^*_\m\subset\ker(\Theta)$.

Il reste donc ˆ voir que   $\sl ii)\iff iii) $ . 

\art{a)}Supposons $c\geq 1$. Si $\P$ est une place finie, alors $U_\P^{(c)}\cap K_\m^*=K_{\P^c\cdot\m}$ (cf. Remarque suivant la DŽfinition \argel) et $\iota(K_{\P^c\cdot\m})=P_{\P^c\cdot\m}\subset I_K(\n)$. Donc
$$\eqalign{{\sl ii)}&\iff P_{\P^c\cdot\m}\subset\ker(\Phi_{L/K} | I_{K}(\n))\cr
&\iff P_{\P^c\cdot\m}\subset\ker(\Phi_{L/K} | I_{K}(\P^{c}\cdot\m))\ \hbox{ car $I_{K}(\n)=I_{K}(\P^{c}\cdot\m)$}\cr
&\iff \P^c\cdot\m\ \hbox{ est admissible (car $\P^{a}\cdot\m$ l'est et donc est div. par les places qui ram.)}\cr
&\iff c\geq v_{\P}(\ff(L/K)).\cr}$$
Pour la dernire Žquivalence, $\Rightarrow$ est clair et $\Leftarrow$ vient du fait que $\P^{a}\cdot\m$ est admissible, donc est un multiple de $\ff(L/K)$. Si $\P$ est infini rŽel, alors dans notre cas, $a=1=c$ et {\sl ii)} et {\sl iii)} sont vraies.

\art{b)}Supposons $c=0$. Si $\P$ est finie, alors ici $U_\P^{(c)}\cap K_\m^*= U_\P\cap K_\m^*=K^*_\m(\P)=\{x\in K^*_\m\mid\hbox{$x$ est premier ˆ $\P$}\}$. Et on a clairement $\iota(U_\P\cap K^*_\m)=P_\m\cap I_K(\n)$. Donc la condition $\sl ii)$ est Žquivalente ˆ 
$$P_\n\subset P_\m\cap I_K(\n)\subset \ker(\Phi_{L/K}| I_{K}(\n))\subset I_{K}(\n) \eqno{(*)}$$ 
ConsidŽrons les deux classes de groupes de congruences~: $\bbH'=\{H'(\overline{\m})\mid \ff' |\overline{\m}\}$, la classe d'Žquivalence de $H'(\n):=P_{\m}\cap I_{K}(\n)$ et de $H'(\m)=P_{\m}$, et la classe $\bbH(L/K)$. La condition $\sl ii)$ est alors Žquivalente ˆ dire (par $(*)$)   que $\bbH'\subset \bbH(L/K)$ qui est Žquivalente ˆ $\ff(L/K) |\ff'$ (cf. Corollaire-DŽfinition \argdi). Or, $\ff' |\m$, donc (petit raisonnement facile) $\sl ii)$ est Žquivalent au fait que $\ff(L/K) |\m$ et donc que $v_\P(\ff(L/K))=0$. Enfin, supposons que $\P\in\gfP_{\R}(K)$. Dans ce cas, $\sl i)$ veut dire que $\ker(\theta_{\P})=U_{\P}=\R^*=\bbK_{\P}^*$, donc $\theta_{\P}$ est l'application triviale. En suivant ˆ la trace les applications successives qui dŽfinissent $\theta_{\P}$, cela veut dire que $P_{\m}=j(\iota(K^*_\m))\subset \ker(\Phi_{L/K}|I_{K}(\n))= \ker(\Phi_{L/K}|I_{K}(\m))$, ce qui veut dire que $\m$ est admissible, donc que $\ff(L/K) |\m$ et finalement que  $v_\P(\ff(L/K))=0$.\qed

\bigskip
\prop {\soustitre (Premier lemme de naturalitŽ)}
\medskip
{\sl Soit $E/K$ une extension abŽlienne de corps de nombres, $L$ un sous-corps intermŽdiaire et $\P\in\gfP(K)$. Alors on a~:
$$\theta_{\P}(L/K)=R\circ\theta_\P(E/K),$$
o $R\, :\, {\rm Gal}(E/K)\to {\rm Gal}(L/K)$ est la restriction habituelle ˆ $L$.

}
{\bf Preuve}

Si $\P$ est complexe, c'est Žvident, puisque $\theta_\P$ est l'homomorphisme trivial.

Supposons $\P$ non complexe et soit $\n=\P^a\cdot\m$ un $K$-module $\P$-admissible pour $E/K$ et $L/K$. Alors le diagramme suivant~:

\vglue 3cm
\psset{xunit=0.9cm,yunit=0.9cm}
\rput(4,2){\rnode{km}{$K_{\m}^* $}} 
\rput(6,2){\rnode{ikm}{$I_K(\m)$}}
\rput(8,2){\rnode{ikn}{$I_K(\n)$}}
\rput(11,3){\rnode{galek}{${\rm Gal}(E/K)$}}
\rput(11,1){\rnode{gallk}{${\rm Gal}(L/K)$}}
\ncline[nodesep=3pt]{->}{km}{ikm}
\Aput{$\iota $} 
\ncline[nodesep=3pt]{->}{ikm}{ikn}
\Aput{$j_\n $} 
\ncline[nodesep=3pt]{->}{ikn}{galek}
\Aput{$\Phi_{E/K} $} 
\ncline[nodesep=3pt]{->}{ikn}{gallk}
\Bput{$\Phi_{L/K} $} 
\ncline[nodesep=3pt]{->}{galek}{gallk}
\Aput{$R$} 

commute, donnant la relation $\Theta(L/K,\P^{a},\m)=R\circ\Theta(E/K,\P^{a},\m)$. Puis, en quotientant par $K_{\m}^*\cap U_{\P}^{(a)}$, $\overline{\Theta}(L/K,\P^{a},\m)=R\circ\overline{\Theta}(E/K,\P^{a},\m)$. Et, composant avec $\bbK_{\P}^*\lra \bbK_{\P}^*/U_{\P}^{(a)}\buildrel\simeq\over\lra K_{\m}^*/(K_{\m}^*\cap U_{\P}^{(a)})$, on trouve le rŽsultat.\qed
\bigskip\goodbreak
\prop {\soustitre (Deuxime lemme de naturalitŽ)}
\medskip
{\sl Soit $L/K$ une extension abŽlienne de corps de nombres (de groupe $G$) et $E/K$ une extension quelconque de corps de nombres. La thŽorie de Galois montre que l'extension $EL/E$ est aussi abŽlienne et que $H:={\rm Gal}(EL/E)$ est identifiable ˆ un sous-groupe de $G$ via la restriction ˆ $L$. Notons $R\, :\, H\to G$ cette restriction. Soit $\P\in\gfP(K)$ et $\gP\in\gfP(E)$ telle que $\gP |\P$. Alors on a~:
$$\theta_\P(L/K)\circ N_{\bbE_\gP/\bbK_\P}=R\circ\theta_\gP(LE/E).$$

}

{\bf Preuve}

Si $\P$ est complexe, alors $\gP$ aussi et donc notre ŽgalitŽ est vraie puisque $\theta_{\P}$ et $\theta_{\gP}$ sont triviaux.

Si $\P$ est rŽelle et $\gP$ complexe ($\P$ ramifie), le membre de droite de l'ŽgalitŽ est trivial. Pour l'autre membre, l'image de $N_{\bbE_{\gP}/\bbK_{\P}}\, (=N_{\C/\R})$ est l'ensemble ${\bbK_{\P}}_{+}^{*}=\R_{+}^{*}$. D'autre part, ici, $\theta_{\P}$ est un homomorphisme qui part de $\bbK^{*}=\R^*$ pour arriver dans $Z(\gP'/\P)$ (o $\gP'$ est n'importe quel idŽal premier de $L$ au-dessus de $\P$) qui est ici d'ordre 1 ou 2. Or, dans tout homomorphisme de ce type, $\R^*_{+}$ est dans le noyau (car tout ŽlŽment est un carrŽ). Donc l'ŽgalitŽ ˆ prouver est vraie dans ce cas.

Si $\P$ et $\gP$ sont rŽelles, alors $\P$ et $\gP$ ramifient ou non simultanŽment dans $L$ respectivement $LE$. En effet, si $\P$ ramifie dans $L$ alors $\sigma_{\gP'}(L)\not\subset \R$ pour tout $\gP'\in\gfP_{\infty}(L)$, $\gP' |\P$ et donc si $\gP''\in\gfP_{\infty}(EL)$, $\gP''|\gP$, on a $\sigma_{\gP''}(EL)\not\subset \R$ car dŽjˆ $R\circ  \sigma_{\gP''}(L)\not\subset \R$. RŽciproquement, si $\gP$ ramifie dans $LE$, soit $\gP''\in\gfP_{\infty}(LE)$, $\gP''|\gP$ et posons $\gP'=\gP''\cap L$. On a donc $\sigma_{\gP''}(LE)\not\subset\R$. Supposons par l'absurde que  $\sigma_{\gP'}(L)\subset \R$, alors puisque $\sigma_{\gP}(E)\subset \R$, on a aussi $\sigma_{\gP''}(LE)\subset\R$, car $LE=L(\alpha_{1},\ldots ,\alpha_{k})$, avec des $\alpha_{i}\in E$ pour tout $i$; c'est une contradiction. Donc $\sigma_{\gP'}(L)\not\subset \R$ ce qui veut dire que $\P$ ramifie.

Supposons donc que  $\P$ et $\gP$ soient rŽelles et ne ramifient pas. Alors l'image de $\theta_{\P}$ et de $\theta_{\gP}$ est triviale (car dans ce cas, $Z(\P)=Z(\gP)=\{1\}$). Donc, l'ŽgalitŽ cherchŽe est vraie. Maintenant, s'ils ramifient les deux, en prenant $b=a=1$, $\theta_{\P}$ et $\theta_{\gP}$ induisent des isomorphismes
$$\R/\R_{+}^*\longrightarrow\cases{Z(L/\P)\simeq\{\pm 1\}&\cr Z(LE/\gP)\simeq\{\pm 1\}\cr}$$
Or, dans notre cas, $N_{\bbE_\gP/\bbK_\P}=N_{\R/\R}={\rm Id}_{\R}$ et l'homomorphisme de restriction $R$ est injectif et envoie en toute gŽnŽralitŽ $Z(LE/\gP)$ dans $Z(L/\P)$. Cela montre donc l'ŽgalitŽ dans ce cas.

Supposons $\P$ et $\gP$ finies. Pour calculer $\theta_\P$, on choisit (comme on l'a vu ˆ la remarque prŽcŽdent la Proposition \3argel ) un $K$-module $\P$-admissible $\n=\P^a\cdot\m$ ($a>0$, $\P\notdiv \m$ ) et $b=a$. Il est Žvident (Lemme \argcr) que $\widetilde{\n}$, le $E$-module extension ˆ $E$ de $\n$, est $\gP$-admissible. En fait, $\widetilde{\n}=\gP^{ae}\cdot\m'$ o $e=e(\gP/\P)$, $\gP\notdiv\m'$ et, clairement, $\m'=\widetilde{\m}\cdot\gP_2^{ae_2}\cdots\gP_r^{ae_r}$, o $\gP_1=\gP$, $\gP_2,\ldots ,\gP_r$ sont les idŽaux premiers de $E$ au-dessus de $\P$, et $e_i=e(\gP_i/\P)$. Souvenons-nous que $N_{E/K}(E^*_{\tilde\m})\subset K^*_\m$ (Lemme \aargbu\ b)), donc {\it a fortiori}, $N_{E/K}(E^*_{\m'})\subset K^*_\m$. On a aussi la relation $N_{E/K}(x\cdot O_E)=N_{E/K}(x)\cdot O_K$ (voir en page \the\gaagc ). Enfin, ce qu'on vient de voir, une vŽrification facile et le ThŽorme \argf\ montre que le diagramme suivant est commutatif~:

\vglue 3cm
\psset{xunit=0.9cm,yunit=0.9cm}
\rput(4,3){\rnode{empr}{$E^*_{\m'}$}} 
\rput(4,1){\rnode{km}{$K^*_{\m}$}} 
\rput(7,3){\rnode{iempr}{$I_E(\m')$}}
\rput(10,3){\rnode{ientl}{$I_{E}(\tilde\n)$}}
\rput(13,3){\rnode{h}{$H$}}
\rput(10,1){\rnode{ikn}{$I_K(\n)$}}
\rput(7,1){\rnode{ikm}{$I_K(\n)$}}
\rput(13,1){\rnode{g}{$G$}}
\ncline[nodesep=3pt]{->}{empr}{km}
\Aput{$N_{E/K}$} 
\ncline[nodesep=3pt]{->}{iempr}{ikm}
\Aput{$N_{E/K}$} 
\ncline[nodesep=3pt]{->}{ientl}{ikn}
\Aput{$N_{E/K}$} 
\ncline[nodesep=3pt]{->}{empr}{iempr}
\Aput{$\iota $} 
\ncline[border=2pt, nodesep=3pt]{->}{iempr}{ientl}
\Aput{$j_{\tilde\n}$}
\ncline[border=2pt, nodesep=3pt]{->}{ientl}{h}
\Aput{$\Phi_{LE/E}$}  
\ncline[nodesep=3pt]{->}{km}{ikm}
\Aput{$\iota$} 
\ncline[border=2pt , nodesep=3pt]{->}{ikm}{ikn}
\Aput{$j_{\n}$} 
\ncline[border=2pt , nodesep=3pt]{->}{ikn}{g}
 \Aput{$\Phi_{L/K}$} 
\ncline[border=2pt , nodesep=3pt]{->}{h}{g}
 \Aput{$R$} 

Soit $x\in\bbE_{\gP}^*$. On choisit, pour calculer $\theta_{\gP}$,  $y\in E_{\m'}^*$ tel que $x\equiv y\pmodast{\widehat{\gP}^{ae}}$. On a donc
$$R\circ\theta_{\gP}(x)=R\circ\Phi_{LE/E}((y))\cdot\gP^{-v_{\gP}(y)})\buildrel \rm diag.\ prec.\over =\Phi_{L/K}((N_{E/K}(y))\cdot\P^{-v_{\P}(N_{E/K}(y))}).$$
Donc, en vertu de la remarque prŽcŽdent la Proposition \3argel, si on montre que
$$ N_{\bbE_\gP/\bbK_\P}(x)\equiv N_{E/K}(y)\pmodast{\widehat{\P}^{a}},$$
on montre la proposition. Par une preuve similaire ˆ la preuve du Lemme \aargbu, on voit que $N_{\bbE_\gP/\bbK_\P}(1+\widehat{\gP}^{ae})\subset 1+\widehat{\P}^{a}$. Cela montre que $N_{\bbE_\gP/\bbK_\P}(x)\equiv N_{\bbE_\gP/\bbK_\P}(y)\pmodast{\widehat{\P}^{a}}$. De plus, par dŽfinition de $y$ et de $\m'$, si $i=2,\ldots r$, on a $y\equiv 1\pmodast{\gP_{i}^{ae_{i}}}$. Ainsi, pour la mme raison que tout ˆ l'heure, on a $N_{\bbE_{\gP_{i}}/\bbK_\P}(y)\equiv N_{\bbE_{\gP_{i}}/\bbK_\P}(1)=1\pmodast{\widehat{\P}^{a}}$. Enfin, il est bien connu (cf. [Fr-Tay, III, 1.10, p. 110]) que $N_{E/K}(y)=\prod_{i=1}^rN_{\bbE_{\gP_{i}}/\bbK_\P}(y)$. En rŽsumŽ, on a~:
$$N_{E/K}(y)=\prod_{i=1}^rN_{\bbE_{\gP_{i}}/\bbK_\P}(y)\equiv N_{\bbE_\gP/\bbK_\P}(y)\equiv N_{\bbE_\gP/\bbK_\P}(x)\pmodast{\widehat{\P}^{a}}.$$\qed
\goodbreak

\bigskip
\coro
\medskip
{\sl Soit $L/K$ une extension abŽlienne de corps de nombres. Soit $\P\in\gfP(K)$ et $\gP\in\gfP(L)$ telle que $\gP |\P$. Alors, $N_{\bbL_{\gP}/\bbK_{\P}}(\bbL_{\gP}^*)\subset\ker(\theta_{\P})$.

}
{\bf Preuve}

En choisissant $E=L$ dans le thŽorme prŽcŽdent, on voit que $\theta_\P(L/K)\circ N_{\bbL_\gP/\bbK_{\P}}=R\circ\theta_{\gP}(L/L)={\rm Id}_{\bbL_\gP}$. Cela prouve le corollaire.\qed
\bigskip
\prop  {\soustitre (Troisime lemme de naturalitŽ)}
\medskip
{\sl Soit $L/K$ une extension abŽlienne de corps de nombres. Soit $\P\in\gfP(K)$ et $\sigma\, : L\to\C$ un plongement. Alors, il est clair que $\sigma$ se prolonge en un homomorphisme qu'on note encore $\sigma\, :\, L_{\gP}\to\C$ pour tout $\gP\in\gfP(L)$. Notons  $K^\sigma, L^\sigma, \P^\sigma $ pour $\sigma(K),\sigma(L)$ respectivement $\sigma(\P)$, Alors $L^\sigma/K^\sigma$ est clairement aussi une extension abŽlienne. Alors, pour tout $y=\sigma(x)\in \bbK^\sigma_{\P^\sigma}$, on a  
$$\sigma\circ\theta_\P(L/K)(x)\circ\sigma^{-1}=\theta_{\P^{\sigma}}(L^\sigma/K^\sigma)(y).$$

}
{\bf Preuve}

C'est une vŽrification Žvidente sachant que $Z(\P^\sigma)=\sigma Z(\P)\sigma^{-1}$ et que $\sigma$ transporte tout, faisant commuter tout diagramme utile pour la preuve.\qed

\bigskip
\th
\medskip

{\sl Sous les mmes hypothses : si $L/K$ est une extension abŽlienne de corps de nombres, $\P$ une place non complexe de $K$ et  $\gP$ une place de $L$ au-dessus de $\P$, alors on a 
$$\ker(\theta_{\P})=N_{\bbL_\gP/\bbK_{\P}}(\bbL_{\gP}^*).$$

}

{\bf Preuve}

On a montrŽ que $\ker(\theta_{\P})\supset N_{\bbL_\gP/\bbK_{\P}}(\bbL_{\gP}^*)$ au Corollaire \7argel.
Montrons donc l'inclusion inverse~: la surjectivitŽ de $\theta_{\P}$ (ThŽorme \argek), l'inclusion que nous venons de prouver et  [Fr-Tay, 1.14+4.2,pp. 111 et 143]) montrent que 

$$[\bbL_{\gP}:\bbK_{\P}]=|Z(\P)|=[\bbK^*_{\P}:\ker(\theta_{\P})]\leq [\bbK^*_{\P}:N_{\bbL_\gP/\bbK_{\P}}(\bbL^*_{\P})].$$

Ainsi, pour montrer l'inclusion inverse, il suffit de prouver que 

$$[\bbK^*_{\P}:N_{\bbL_\gP/\bbK_{\P}}(\bbL^*_{\gP})]\leq [\bbL_{\gP}:\bbK_{\P}].\eqno{(***)}$$

Si $\P$ est infinie, on se souvient que $\bbK_{\P}=\bbL_{\gP}=\R$ si $\P$ ne ramifie pas dans $L$, et dans ce cas, $N_{\bbL_\gP/\bbK_{\P}}$ est l'identitŽ et donc les indices cherchŽs valent 1; et si $\P$ ramifie dans $L$, alors dans ce cas, $\bbK_{\P}=\R$, $L_{\gP}=\C$ et $N_{\bbL_\gP/\bbK_{\P}}(\C^*)=\R_{+}^*$ et donc les indices cherchŽs valent 2, car $N_{\bbL_\gP/\bbK_{\P}}$ est la norme complexe.

Supposons $\P$ finie. Prouvons ce rŽsultat par rŽcurrence sur le nombre de facteurs premiers de $[L:K]$. Si ce nombre est 1, alors l'extension $L/K$ est cyclique (un groupe abŽlien d'ordre premier est cyclique...). Et dans ce cas lˆ, on a dŽjˆ prouvŽ que  $[\bbL_{\gP}:\bbK_{\P}]=[\bbK^*_{\P}:N_{\bbL_\gP/\bbK_{\P}}(\bbL^*_{\gP})]$ (Proposition \argci). Supposons maintenant que le nombre de facteurs premiers de $[L:K]$ soit strictement supŽrieur ˆ 1 et que le thŽorme est prouvŽ pour toute extension d'indice plus petit. Choisissons une extension intermŽdiaire $K\subsetneq E\subsetneq L$, et on pose $\gP_{0}=\gP\cap E$. On a $\bbK_{\P}\subset \bbE_{\gP_{0}}\subset \bbL_{\gP}$. Si on pose $N_{1}=N_{\bbL_{\gP}/ \bbE_{\gP_{0}}}$ et $N_{2}=N_{\bbE_{\gP_{0}}/\bbK_{\P}}$, alors bien sžr, $N_{\bbL_\gP/\bbK_{\P}}=N_{2}\circ N_{1}$. D'autre part, le lemme technique (Lemme \argcj) montre facilement que $[N_{2}(\bbE^*_{\gP_{0}}):N_{2}(N_{1}(\bbL^*_{\gP}))]\leq  [\bbE^*_{\gP_{0}}:N_{1}(\bbL^*_{\gP})]$. Ainsi, on a~:
$$\eqalign{[\bbK^*_{\P}:N_{\bbL_\gP/\bbK_{\P}}(\bbL^*_{\P})]&=[\bbK^*_{\P}:N_{2}(\bbE^*_{\gP_{0}})]\cdot [N_{2}(\bbE^*_{\gP_{0}}):N_{2}(N_{1}(\bbL^*_{\P}))]\cr &\buildrel\rm hyp.\ de\ rec+rem.\over\leq [\bbE_{\gP_{0}}:\bbK_{\P}]\cdot [\bbE^*_{\gP_{0}}:N_{1}(\bbL^*_{\gP})]\buildrel \rm hyp.\ de\ rec\over \leq  [\bbE_{\gP_{0}}:\bbK_{\P}]\cdot [\bbL_{\gP}:\bbE_{\gP_{0}}].\cr}$$
Cela montre la relation $(***)$ et donc le thŽorme.\qed

\bigskip\goodbreak

\coro
\medskip

{\sl Sous les mmes hypothses, i.e. si $L/K$ est est une extension abŽlienne de corps de nombres, $\P$ une place non complexe de $K$  et $\gP|\P$ est une place de $L$ au-dessus de $\P$, alors on a~:

$$[\bbK^*_{\P}:N_{\bbL_\gP/\bbK_{\P}}(\bbL^*_{\gP})]= [\bbL_{\gP}:\bbK_{\P}].$$

De plus, pour $c\geq 0$, on a 
$$U_{\P}^{(c)}\subset N_{\bbL_\gP/\bbK_{\P}}(\bbL^*_\gP)\iff c\geq v_{\P}(\ff(L/K)).$$

}

{\bf Preuve}

L'ŽgalitŽ $[\bbK^*_{\P}:N_{\bbL_\gP/\bbK_{\P}}(\bbL^*_{\gP})]= [\bbL_{\gP}:\bbK_{\P}]$ ˆ ŽtŽ montrŽ dans la preuve du thŽorme prŽcŽdent. La second affirmation est un corollaire immŽdiat du thŽorme prŽcŽdent combinŽ avec la Proposition \4argel.

\qed

\bigskip
\lem
\medskip

{\sl Sous les mmes hypothses : si $L/K$ est une extension abŽlienne de corps de nombres, $\P$ une place non complexe de $K$ et  $\gP$ une place de $L$ au-dessus de $\P$, alors on a

$$U_{\P}=N_{\bbL_\gP/\bbK_{\P}}(U_{\gP})\iff U_{\P}\subset N_{\bbL_\gP/\bbK_{\P}}(\bbL^*_\gP)\iff  \P\hbox{ est non ramifiŽe dans }L$$

}

{\bf Preuve}

Si $\P$ est une place infinie, c'est une vŽrification~: si $\P$ est ramifiŽ, $N_{\bbL_\gP/\bbK_{\P}}(U_\gP^{(0)})=N_{\C/\R}(\C^*)=\R^*_{+}\subsetneq \R^*=U_{\P}^{(0)}$; et si $\P$ est non ramifiŽ, $N_{\P}(U_\gP^{(0)})=N_{\R/\R}(\R^*)=\R^*= U_{\P}^{(0)}$.

Supposons $\P$ finie. La partie ``$\Rightarrow$'' de la premire Žquivalence est Žvidente. La partie  ``$\Leftarrow$'' aussi~:  soit $\pi$ (resp. $\pi'$) une uniformisante de $\bbL_\gP$ (resp. $\bbK_\P$). On sait que $N_\P(\pi)=u\cdot \pi'^f$, o $f=f(\gP/\P)$ et $u\in U_\P$. Donc dire que $U_{\P}\subset N_{\bbL_\gP/\bbK_{\P}}(L^*_\gP)$, implique puisque $L^*_\gP=<\pi>\times\, U_\gP$ et $K^*_\P=<\pi'>\times\, U_\P$, que $U_\P\subset N_{\bbL_\gP/\bbK_{\P}}(U_\gP)$, mais on a toujours  $N_{\bbL_\gP/\bbK_{\P}}(U_\gP)\subset U_{\P}$ (car la norme d'un entier est un entier) et on trouve bien $N_{\bbL_\gP/\bbK_{\P}}(U_\gP)= U_{\P}$.

Montrons la seconde Žquivalence~:  Supposons que $U_{\P}\subset N_{\bbL_\gP/\bbK_{\P}}(\bbL^*_\gP)$, on a vu que c'est Žquivalent au fait que $U_{\P}=N_{\bbL_\gP/\bbK_{\P}}(U_{\gP})$. Ainsi, toujours puisque $N_\P(\pi)=u\cdot \pi'^f$, on a que $\pi'^f$ est une norme (puisque $u$ en est une). Cela prouve que $N_{\bbL_\gP/\bbK_{\P}}(\bbL^*_\gP)=<\pi'^f>\times\, U_\P$. Ainsi $[\bbK^*_\P:N_{\bbL_\gP/\bbK_{\P}}(\bbL^*_\gP)]=f$. Or, on a montrŽ au Corollaire \argen\ que $[\bbK^*_\P:N_{\bbL_\gP/\bbK_{\P}}(\bbL^*_\gP)]=[\bbL_\gP:\bbK_\P]=f\cdot e$, o $e=e(\gP/\P)$. Donc $e=1$, ce qui montre que $\P$ n'est pas ramifiŽ dans $L$. Inversement, si $\P$ est non ramifiŽe, le Lemme \argdl\ nous montre que $\P\notdiv \ff(L/K)$, donc $v_{\P}(\ff(L/K))=0$ ce qui implique en vertu du Corollaire \argen\ que $U_{\P}^{(0)}=U_{\P}\subset N_{\bbL_\gP/\bbK_{\P}}(\bbL^*_\gP)$. Le lemme est alors dŽmontrŽ.\qed

 \bigskip\goodbreak

Voilˆ enfin un des thŽormes que nous visions depuis un moment~:
\bigskip
\th
\medskip

{\sl Soit $L/K$ une extension abŽlienne de corps de nombres. Alors les places ramifiŽes divisent le conducteur. Plus prŽcisŽment, $\P$ ramifie dans $L\iff\P|\ff(L/K)$. Cela implique en vertu du Lemme \argdk\ que $\ff(L/K)$ est admissible.
}

{\bf Preuve}

Soit $\P$ une place de $K$ non complexe. Posons $c$ l'exposant de $\P$ dans la dŽcomposition de $\ff$. Alors on a la sŽrie d'Žquivalence~:

$$\eqalign{\P\notdiv\ff&\iff c=0\buildrel \rm Prop.\ \4argel\over \iff U_{\P}\subset \ker(\theta_{\P})\cr &\buildrel \rm Thm. \argem\over\iff U_{\P}\subset N_{\bbL_\gP/\bbK_{\P}}(\bbL^*_\gP)\cr &\buildrel \rm Lemme \argeo\over\iff  \P\hbox{ est non ramifiŽe dans }L.\cr}$$

\qed

\bigskip

{{{\soustitre Corollaire ({\the\chapnomb}.{\the\nomb})}\global\advance\nomb by 1}\soustitre (construction et existence du corps de Hilbert)}
\medskip
{\sl Soit $K$ un corps de nombres. Alors il existe un (unique) corps de nombres $H\supset K$  (appelŽ {\it corps de Hilbert de $K$}) possŽdant les propriŽtŽs suivantes~:

\art{a)}L'extension $H/K$ est abŽlienne (de groupe disons $G$).

\art{b)}Aucune place de $K$ ne ramifie dans $H$.

\art{c)}L'application d'Artin $\Phi_{H/K}\, : I_{K}\to G$ est de noyau $P$, l'ensemble des idŽaux fractionnaires principaux de $K$; ainsi, puisque $\Phi_{H/K}$ est surjective (ThŽorme \argaq), cela implique que $G$ est isomorphe au groupe des classes ${\cal CL}_{K}=I_{K}/P$.

En outre, L'extension $H/K$ contient toute extension de $K$ abŽlienne non ramifiŽe (i.e. satisfaisant les propriŽtŽs a) et b) de la mme cl™ture algŽbrique. 

}

{\bf Preuve}

ConsidŽrons le $K$-module $O_{K}$ $(=(O_{K},\emptyset )$) et la classe d'Žquivalence de sous-groupes de congruences $$\bbH=\{ P({\m})=P\cap I_{K}(\m)\mid \m\hbox{ est un $K$-module}\}. $$ Il est Žvident que $\bbH$ est bel et bien une classe d'Žquivalence et que le conducteur de cette classe est $\ff=O_{K}$ et le groupe de congruence associŽ ˆ $\ff$ est $P$. Par le thŽorme d'existence du corps de classe (cf. ThŽorme \argef), il existe une extension abŽlienne $H/K$ telle que $\bbH=\bbH(H/K)$. Le thŽorme prŽcŽdent nous montre que $O_{K}$ est admissible, donc ${\rm Gal}(H/K)\simeq I_{K}/P$, d'o la partie c). Puisque $\ff=O_{K}$, le thŽorme prŽcŽdent nous dit que l'extension $H/K$ est non ramifiŽe, d'o la partie b). Pour la dernire affirmation, soit $L/K$ une extension abŽlienne non ramifiŽe. A nouveau, gr‰ce au thŽorme prŽcŽdent, le conducteur $\ff=\ff(L/K)$ de cette extension vaut $O_{K}=\ff(H/K)$. Soit $H(\ff,L/K)\in \bbH(L/K)$, le sous-groupe de congruence pour $\ff(L/K)$. On a par dŽfinition $P_{\ff}\subset H(\ff,L/K)$. Or, puisque $\ff=O_{K}$, on a $P_{\ff}=P(\ff)=H(\ff,H/K)$. Cela montre que $H(\ff,H/K)\subset H(\ff,L/K)$ et donc, $\bbH(H/K)\subset \bbH(L/K)$ et donc $L\subset H$ en vertu de la Proposition \1argdm. Ce qui achve la preuve de ce corollaire.\qed

\bigskip

{\bf Remarque}

Le corps de Hilbert a encore une propriŽtŽ remarquable : si $K$ est un corps de nombres, alors tout idŽal fractionnaire de $K$ devient principal dans le corps de Hilbert de $K$. La preuve de ce rŽsultat est assez longue. Le chapitre suivant est consacrŽ ˆ la preuve de ce rŽsultat
\vfill\eject


\global\advance\chapnomb by 1
\nomb=1

\centerline{\para Chapitre 12 : }
\smallskip
\centerline{\para Capitulation des idŽaux d'un corps nombres}
\smallskip
 \centerline{\para dans son corps de Hilbert }
\bigskip

Ici, nous allons montrer que tout idŽal d'un corps de nombres devient principal vu dans le corps de Hilbert, on dit qu'il {\it capitule}. Evidemment, dans la vraie vie, il vaut mieux que les idŽaux ne capitulent pas, mais en mathŽmatique, cela simplifie bien les choses. Avant de s'attaquer de front au problme, nous devons faire une petite incursion dans la thŽorie des groupe pour dŽfinir un homomorphisme dit de Òtransfert''. Nous montrerons ensuite que cet homomorphisme est trivial sous certaines hypothses, puis nous utiliserons ce rŽsultat pour rŽsoudre notre problme. Notons que ce rŽsultat a ŽtŽ prouvŽ par FurtwŠngler en 1930.

\bigskip

\defis
\medskip
Soit $G$ un groupe, $H\subset G$ un sous-groupe de $G$ et $J\subset H$ un sous-groupe normal dans $H$. On suppose que $[G:H]<\infty$ et que $H/J$ soit abŽlien. Une {\it transversale (ˆ droite) de $H$ dans $G$} est une partie $T\subset G$ telle que 
$$G=\bigsqcup_{t\in T}H\cdot t.\eqno{(1)}$$
Cette rŽunion est finie par hypothse. Soit donc $T=\{t_1,\ldots ,t_n\}$ ($n=[G:H]$) une transversale de $H$ dans $G$. Soit $g\in G$. Pour chaque $i=1,\ldots ,n$, il existe $h_i\in H$ et $j(i)\in \{1,\ldots ,n\}$ tel que $t_i\cdot g=h_i\cdot t_{j(i)}$. L'application $i\mapsto j(i)$ est une permutation de $\{1,\ldots ,n\}$, car $Tg=\{t_1\cdot g,\ldots ,t_n\cdot g\}$ est aussi une transversale (il suffit de multiplier (1) par $g$). Alors on pose
$$\eqalign{{\rm Ver}\, :\, G&\lra H/J\cr g&\longmapsto \prod_{i=1}^n J\cdot h_i=J\cdot \prod_{i=1}^n h_i.\cr}$$
Cette application est Òpresque" l'homomorphisme de transfert, mais pas tout-ˆ-fait. NŽanmoins, nous allons voir qu'elle est indŽpendante de la transversale $T$ et que c'est un homomorphisme de groupe.
\bigskip
\headline={\hfill \phantom{ouuh}\hfill}
\lem
\medskip
{\sl Sous les mmes hypothses,  l'application ${\rm Ver}\, :\,  G\to H/J$ est indŽpendante de $T$ et c'est un homomorphisme de groupe.

}
{\bf Preuve}

Soit $T'=\{t_1',\ldots ,t_s'\}$ une autre transversale de $H$ dans $G$. Choisissons une numŽrotations de $t_i'$ telle que $H\cdot t_i=H\cdot t_i'$ pour tout $i=1,\ldots ,n$. Posons, pour ces mmes $i$, $h_i''\in H$ tel que $t_i'=h_i''\cdot t_i$. Soit $g\in G$. De $t_i\cdot  g_i=h_i\cdot t_{j(i)}$ suit (en multipliant par $h_i''$) $h_i''\cdot t_i\cdot  g_i=h_i''\cdot h_i\cdot h_{j(i)}''^{-1}\cdot h_{j(i)}''\cdot t_{j(i)}$, c'est-ˆ-dire 
$$t_i'\cdot g=h_i'\cdot t_{j(i)}'\quad \hbox{avec }h_i'=h_i''\cdot h_i\cdot h_{j(i)}''^{-1}$$
On en dŽduit (en utilisant plusieurs fois que $H/J$ est commutatif et que $i\mapsto j(i)$ est une permutation de  $\{1,\ldots ,n\}$)~:
$$\eqalign{\prod_{i=1}^nJ\cdot h_i'&=\prod_{i=1}^nJ\cdot h_i''\cdot h_i\cdot h_{j(i)}''^{-1}=
\prod_{i=1}^nJ\cdot h_i''\cdot\prod_{i=1}^n J\cdot h_i\cdot \prod_{i=1}^nJ\cdot h_{j(i)}''^{-1}\cr
&=\prod_{i=1}^nJ\cdot h_i''\cdot\prod_{i=1}^n J\cdot h_i\cdot \prod_{i=1}^nJ\cdot h_i''^{-1}=
\prod_{i=1}^nJ\cdot h_i''\cdot\prod_{i=1}^n J\cdot h_i\cdot \left (\prod_{i=1}^nJ\cdot h_i''\right )^{-1}\cr
&=\prod_{i=1}^n J\cdot h_i.\cr}$$
Cela montre que ${\rm Ver}(g)$ ne dŽpend pas de la transversale $T$.

Montrons ˆ prŽsent que c'est un homomorphisme de groupe~: soit donc $T=\{t_1,\ldots ,t_n\}$ une transversale et $g,h\in G$. On a comme avant $t_i\cdot g=h_i\cdot t_{j(i)}$ et $t_i h=h_i'\cdot t_{k(i)}$, pour $i=1,\ldots n$, $h_i, h_i'\in H$, $i\mapsto j(i)$, $i\mapsto k(i)$ sont des permutations de $\{1,\ldots ,n\}$. On a alors $t_i\cdot g\cdot h=h_i\cdot t_{j(i)}\cdot h=h_i\cdot h'_{j(i)}t_{k(j(i))}$. Ainsi~(toujours avec les mmes propriŽtŽs de $J/H$, $i\mapsto j(i)$ )~:
$$\eqalign{{\rm Ver}(g\cdot h)&=\prod_{i=1}^nJ\cdot h_i\cdot h_{j(i)}'=\prod_{i=1}^n J\cdot h_i\cdot \prod_{i=1}^n J\cdot h_{j(i)}'\cr &= \prod_{i=1}^nJ\cdot h_i\cdot \prod_{i=1}^n J\cdot h_i'={\rm Ver}(g)\cdot{\rm Ver}(h)\cr}$$
\qed
\bigskip

\defi

Soit $G$ un groupe. On rappelle que $G'=D(G:G)$ est le sous-groupe engendrŽ par les commutateurs $[a,b]=aba^{-1}b^{-1}$, $a,b\in G$ et $G/G'$, l'abŽlianisŽ de $G$, est le plus grand quotient abŽlien de $G$. Sous les mmes hypothses que le lemme prŽcŽdent, mais en supposant que $J=H'$. On a que l'application ${\rm Ver}\, :\, G\to H/H'$ est bien dŽfinie. Puisque $H/H'$ est abŽlien et par dŽfinition de $G'$,  il est clair que $G'\subset \ker ({\rm Ver})$. Ainsi, on peut dŽfinir une application 
$$V_{G\to H}\, :\, G/G'\longrightarrow H/H'$$
 \newcount\gaagap\gaagap=\gaga
qu'on appelle l'{\it homomorphisme de transfert de $G$ sur $H$}. 

Le transfert est trs utile pour montrer de jolis thŽormes sur les groupes. Nous ne pouvons pas rŽsister ˆ la tentation d'en citer quelques-uns, mme s'il ne sont pas utiles pour la suite~:
\bigskip
\headline={\hfill \smcap Capitulation des idŽaux d'un corps nombres dans son corps de Hilbert \hfill}
\th \ {\soustitre (ThŽorme de Schur)}
\medskip
{\sl Si $G$ est groupe tel que le centre $Z(G)=\{x\in G\mid yx=xy\ \forall y\in G\}$ est tel que $[G:Z(G)]<\infty$, alors $|G'|<\infty$.

}

\bigskip

\th 
\medskip

{\sl Soit $G$ un groupe fini et $p$ le plus petit diviseur de $|G|$. Si un des $p$-sous-groupe de Sylow $P$ de $G$ est cyclique (donc tous...), alors il existe un sous-groupe normal $N$ de $G$ tel que $G/N\simeq P$.

}
\bigskip
\coro
\medskip
{\sl Si $G$ est un sous-groupe simple fini non abŽlien, alors les 2-sous-groupe de Sylow de $G$ ne sont pas cycliques.

}
\bigskip
\coro
\medskip
{\sl Soit $G$ un groupe fini. Si tous les sous-groupe de Sylow de $G$ sont cycliques, alors $G$ est rŽsoluble.

}
\bigskip
En revanche le thŽorme suivant sera le rŽsultat crucial de ce chapitre. Il s'appelle le ÒthŽorme de l'idŽal principal de la thŽorie des groupe'' (visiblement, car l'unique corollaire connu de ce rŽsultat est prŽcisŽment le sujet de ce chapitre) ~:
\bigskip\goodbreak
\th 
\medskip
{\sl Soit $G$ un groupe. Supposons que $[G:G']<\infty$ et que $G'/G''$ soit de gŽnŽration finie. Alors l'homomorphisme de transfert $V_{G\to G'}$ est trivial (c'est-ˆ-dire tout est envoyŽ sur 1).

}

\medskip

La preuve de ce thŽorme est assez longue et nous devrons tout d'abord faire ce qu'on peut appeler une ``version additive du transfert''. D'abord une
\bigskip
\defi
\medskip
Soit $G$ un groupe et $H\subset G$ tel que $[G:H]<\infty$. On introduit $\Z[G]$ (l'anneau de groupe de $G$, qui est l'ensemble des sŽries formelles $\sum_{g\in G}m_{g}\cdot g$, o $m_{g}\in \Z$ pour tout $g\in G$). Le noyau de l'homomorphisme~:
$$\eqalign{\Z[G]&\lra \Z\cr \sum_{g\in G}m_{g}\cdot g&\longmapsto \sum_{g\in G}m_{g}\cr}$$
 \newcount\gaagaq\gaagaq=\gaga
est un idŽal, appelŽ {\it l'idŽal d'augmentation}, qu'on note $I_{G}$. Remarquons que $\Z[H]\subset \Z[G]$ et que $I_{H}\subset I_{G}$. On observe aussi que $(g-1)_{g\in G\setminus\{1\}}$ forme une $\Z$-base de $I_{G}$. En effet , si $0=\sum_{g\ne 1}m_{g}(g-1)$, on a $(\sum_{g\ne1}m_{g})\cdot 1=\sum_{g\ne 1}m_{g}\cdot g$, ce qui implique que $m_{g}=0$ pour tout $g\ne 1$, donc la famille est libre. De plus, si $\sum_{g\in G}m_{g}\in I_{G}$, on a $\sum_{g\in G}m_{g}=\sum_{g\ne 1}m_{g}(g-1)+\underbrace{(\sum_{g\in G}m_{g})}_{=0}\cdot 1$, donc on a la gŽnŽration. On montre de mme que si $g_{0}\in G$, les ŽlŽments $(g-1)\cdot g_{0}$ forment aussi une base de $I_{G}$, cela vient de la relation $g-1=(g\cdot g_{0}^{-1}-1)\cdot g_{0}-(g_{0}^{-1}-1)\cdot g_{0}$. 

Notons $d$ l'application $\Z[G]\to\Z[G]$ telle que $d(g)=g-1$. Dans $\Z[G]$, on peut former les idŽaux $I_G^2$ et $I_H\cdot I_G$ (ce dernier n'Žtant qu'un idŽal ˆ droite, et donc un sous-groupe additif) ainsi que le sous-groupe $I_H+I_HI_G$.
\bigskip\goodbreak
\lem {\soustitre (version additive du transfert)}

{\sl Soit $G$ un groupe et $H\subset G$ un sous-groupe d'indice fini, ainsi que $T$ une transversale de $H$ dans $G$. Alors il existe des isomorphismes, notŽs $\log$ et une application $S$ tels que le diagramme suivant soit commutatif~:

\vglue 2.5cm 
\psset{xunit=1cm,yunit=1cm}

\rput(6,0){\rnode{IGI}{$I_G/I_G^2$}}
\rput(6,2){\rnode{GG}{$G/G'$}}
\rput(10,0){\rnode{IHI}{$(I_H+I_H\cdot I_G)/(I_H\cdot I_G)$}}
\rput(10,2){\rnode{HH}{$H/H'$}}
\ncline[nodesep=3pt]{->}{GG}{HH}
\Aput{$V_{G\to H} $} 
\ncline[nodesep=3pt]{->}{IGI}{IHI}
\Aput{$S $} 
\ncline[nodesep=3pt]{->}{GG}{IGI}
\Aput{$\log$} 
\ncline[nodesep=3pt]{->}{HH}{IHI}
\Aput{$\log$} 
\pcline(6,0.5)(6,1.5) \aput{:U}{$\simeq$} 
\pcline(10,0.5)(10,1.5) \aput{:U}{$\simeq$} 
\medskip

o $S(x\bmod I_G^2)=(\sum_{t\in T} t)\cdot x\bmod I_H\cdot I_G$.

}

{\bf Preuve}

Tout d'abord, on vŽrifie que $d(x\cdot y)=d(x)+d(y)+d(x)\cdot d(y)$ pour tout $x,y\in G$. Ainsi, l'application $h\mapsto d(h) \bmod I_{H}\cdot I_{G}$ est un homomorphisme de groupe de $H\to (I_{H}+I_{H}\cdot I_{G})/(I_H\cdot I_G)$, le groupe de droite est bien entendu additif (donc abŽlien). Cela donne donc un homomorphisme 
$$\log\, :\, H/H'\to (I_H+I_H\cdot I_G)/(I_H\cdot I_G).$$
 Par les thŽormes d'isomorphismes $(I_H+I_H\cdot I_G)/(I_H\cdot I_G)\simeq I_{H}/(I_{H}\cap I_{H}\cdot I_{G})$, et donc $\log$ est surjective, car les $d(h)$ engendrent $I_{H}$. Nous allons montrer que $\log$ est en fait un isomorphisme. Soit $T$ une transversale de $H$ dans $G$. On notera $\tau$, l'unique ŽlŽment de $T\cap H$. On affirme que les ŽlŽments $d(h)\cdot t$ avec $t\in T$ et $h\in H\setminus\{1\}$ forment une base de $I_H+I_H\cdot I_G$. En effet, on vient de voir que les $d(h)\cdot \tau$ engendrent  $I_H$. D'autre part, $I_H\cdot I_G$ est engendrŽ par les $(h-1)\cdot (h'\cdot t-1)=(h\cdot h'-1)\cdot t-(h'-1)\cdot t-(h-1)=d(h\cdot h')\cdot t-d(h')\cdot t-d(h\cdot\tau^{-1})\cdot\tau+d(\tau^{-1})\cdot\tau$. Pour l'indŽpendance, si 
$$0=\sum_{\eqalign{\noalign{\vskip-3pt}\scriptstyle t\in T\phantom{88}\cr\noalign{\vskip-9pt} \scriptstyle 1\ne h\in H\cr}}n_{t,h}\cdot d(h)\cdot t=\sum_{\eqalign{\noalign{\vskip-3pt}\scriptstyle t\in T\phantom{88}\cr\noalign{\vskip-9pt} \scriptstyle 1\ne h\in H\cr}}n_{t,h}\cdot h\cdot t-\sum_{t\in T}\left (\sum_{1\ne h\in H} n_{t,h}\right )\cdot t.$$
Les $h\cdot t$, $t\in T$, $1\ne h\in H$ et les $t=1\cdot t$ forment exactement tous les ŽlŽments de $G$, donc les $n_{t,h}=0$ pour tout $t,h$. C'est donc une famille libre.
\smallskip
On dŽfinit un homomorphisme (par l'image de la base qu'on vient de trouver)
$$\eqalign{I_{H}+I_{H}\cdot I_{G}&\lra H/H'\cr d(h)\cdot t&\longmapsto h\bmod H'.\cr}$$
On vient de voir que $I_{H}\cdot I_{G}$ Žtait engendrŽ par les ŽlŽments $d(h)\cdot d(h'\cdot t)=d(h\cdot h')\cdot t-d(h')\cdot t-d(h\cdot\tau^{-1})\cdot\tau+d(\tau^{-1})\cdot\tau$ qui est envoyŽ sur $h\cdot\tau\cdot h^{-1}\cdot \tau^{-1}\in H'$. Donc $I_{H}\cdot I_{G}$ est dans le noyau, ce qui veut dire qu'on a l'homomorphisme
$$\exp\,:\, I_{H}+I_{H}\cdot I_{G}/(I_{H}\cdot I_{G})\lra H/H'$$
caractŽrisŽ par $d(h)\cdot t \mapsto h\bmod H'$ . Il est clair que $\rm exp$ est l'inverse de $\log$~: $\exp(\log(h\bmod H'))=\exp(d(h))=\exp(d(h\cdot \tau^{-1})\cdot\tau-d(\tau^{-1})\cdot\tau)=h\cdot\tau^{-1}\cdot\tau=h$. Et inversement  $\log(\exp(d(h)\cdot t))=\log(h)=d(h)$, et on observe que $d(h)\cdot t-d(h)=d(h)\cdot d(t)\in I_{H}\cdot I_{G}$, donc $d(h)\bmod I_{H}\cdot I_{G}=d(h)\cdot t\bmod I_{H}\cdot I_{G}$. On a ainsi montrŽ que $\log$ (et $\exp$) Žtait un isomorphisme. En particulier, en appliquant cela ˆ $H=G$, on obtient aussi l'isomorphisme~:
$$\log\, :\, G/G'\longmapsto I_{G}/I_{G}^2.$$

Traduisons maintenant le transfert $V_{G\to H}$ en un homomorphisme
$$S\, :\, I_{G}/I_{G}^2\lra  I_{H}+I_{H}\cdot I_{G}/(I_{H}\cdot I_{G}).$$
Rappelons que $V_{G\to H}(g\bmod G')=\prod_{t\in T}h_{t}\bmod H'$, o $\forall t$, on a $t\cdot g=h_{t}\cdot t'$, avec $t'\in T$ et $h_{t}\in H$. Cela veut dire que $S(d(g)\bmod I_{G}^2)=\sum_{t\in T}d(h_{t})\bmod I_{H}\cdot I_{G}$. Finalement, la relation $d(t)+t\cdot d(g)=d(h_{t})+d(t')+d(h_{t})\cdot d(t')$, donne en sommant sur tous les $t$,
$$\underbrace{\sum_{t\in T}d(t)}_{=(*)}+\left (\sum_{t\in T}t\right )\cdot d(g)=\sum_{t\in T}d(h_{t})+\underbrace{\sum_{t'\in T}d(t')}_{=(*)}+\underbrace{\sum_{t\in T}d(h_{t})\cdot d(t')}_{\in I_{H}\cdot I_{G}}.$$
On a donc montrŽ que $S(d(g)\bmod I_{G}^2)=(\sum_{t\in T}t )\cdot d(g)\bmod I_{H}\cdot I_{G}$, ce qui montre notre lemme (les $d(g)$ engendrent $I_{G}$).\qed
\bigskip
{\soustitre Preuve du ThŽorme \argfb}
\medskip
En rempla\c cant $G$ par $G/G''$, on se ramne ˆ prouver le thŽorme sous l'hypothse que $G''=\{1\}$, en effet,
il est clair que $(G/G'')'=G'/G''$, que $(G/G'')''=G''/G''=\{1\}$, que $G/G'\buildrel\rm thm.d'isom.\over \simeq (G/G'')/(G'/G'')$ et que le carrŽ suivant commute~:

\vglue 2.5cm 
\psset{xunit=1cm,yunit=1cm}

\rput(6,0){\rnode{ggprr}{$(G/G'')/(G'/G'')$}}
\rput(6,2){\rnode{GG}{$G/G'$}}
\rput(10,0){\rnode{ggsecc}{$(G'/G'')/(G''/G'')$}}
\rput(10,2){\rnode{HH}{$G'/G''$}}
\ncline[nodesep=3pt]{->}{GG}{HH}
\Aput{$V_{G\to G'} $} 
\ncline[nodesep=3pt]{->}{ggprr}{ggsecc}
\Aput{$V_{G/G''\to G'/G'' }$} 
\ncline[nodesep=3pt]{->}{GG}{ggprr}
\ncline[nodesep=3pt]{->}{HH}{ggsecc}
\pcline(6,0.5)(6,1.5) \aput{:U}{$\simeq$} 
\pcline(10,0.5)(10,1.5) \aput{:U}{$\simeq$} 
\medskip
\goodbreak
car les transversales ne posent pas non plus de problme. Cela permet de supposer que $G'$ et donc $G$ est de gŽnŽration finie. Pour la preuve, on fixe $g_1,\ldots ,g_r$ un systme de gŽnŽrateurs de $G$ et $T$ une transversale ˆ droite de $G'$ dans $G$. En vertu du Lemme \argfd , il s'agit donc de dŽmontrer que pour tout $g\in G$, on a
$$\left (\sum_{t\in T} t\right )\cdot d(g)\equiv 0\pmod{I_{G'}\cdot I_G},\eqno{(*)}$$
car les $d(g)$ engendrent $I_G$.

Remarquons d'abord la chose suivante~:  si $H\subset G$ est un sous-groupe normal, alors  on a un isomorphisme
$$\Z[G/H]\simeq\Z[G]/(I_H\cdot \Z[G]).\eqno{(**)}$$
En effet, d'abord $I_H\cdot \Z[G]$ est un idŽal bilatre~: si $h\in H$ et $g\in G$, alors si $h'\in H$ est l'ŽlŽment tel que $g\cdot h=h'\cdot g$, alors $g\cdot (h-1)=(h'-1)\cdot g$. Ensuite, si $T$ est une transversale ˆ droite de $H$ dans $G$, alors $\Z[G]=\bigoplus_{t\in T}\Z[H]\cdot t$ et $I_H\cdot \Z[G]=\bigoplus_{t\in H}I_H\cdot t$, les dŽcompositions Žtant compatibles (i.e $I_H\cdot t\subset \Z[H]\cdot t\ \forall\, t$). De sorte que
$$\Z[G]/(I_H\cdot \Z[G])\simeq \bigoplus_{t\in T}(\Z[H]\cdot t)/(I_H\cdot t)\simeq \bigoplus_{t\in T} \left (\Z[H]/I_H\right )\cdot\overline{t}\simeq \bigoplus_{t\in T}\Z\cdot\overline{t},\eqno{(***)}$$
o $\overline{g}$ est l'image de $g$ dans le quotient $\Z[G|/(I_H\cdot\Z[G])$, pour tout $g\in G$; et le dernier isomorphisme vient du fait que $h-1\in I_H$ et donc $\overline{h}=1$, pour tout $h\in H$. Cette dernire ŽgalitŽ montre que l'homomorphisme multiplicatif  $G\to \Z[G]/(I_H\cdot \Z[G])$ envoie $h\cdot t\mapsto \overline{t}$; ce qui permet d'identifier $G/H$ avec $\{\overline{t}\mid t\in T\}$ et cela fourni l'isomorphisme (d'anneau) entre $\Z[G/H]$ et $\bigoplus_{t\in T}\Z\cdot\overline{t}$ et donc entre $\Z[G/H]$ et $\Z[G]/(I_H\cdot \Z[G])$ en vertu de $(***)$. Donc $(**)$ est prouvŽ. Nous appliquerons ce rŽsultat ˆ $H=G'$. Dans ce cas-lˆ, $\Z[G/G']$ (et donc $\Z[G]/(I_{G'}\cdot\Z[G])$) est un anneau commutatif (puisque $G/G'$ est un groupe commutatif). Donc nous pourrons faire un peu d'algbre linŽaire dans ce cas-lˆ.

Voici encore un rŽsultat~: soit $g\in G$. Si on Žcrit $g=x_{1}\cdots x_{N}$, o $x_{i}\in \{g_{1},\ldots,g_{n}, g_{1}^{-1}, \ldots , g_{n}^{-1}\}$, pour $i=1,\ldots , N$, alors d(g) peut s'Žcrire
$$d(g)=\sum_{i=1}^ny_{i}\cdot d(g_{i})\eqno{(+)}$$
pour des $y_{i}\in \Z[G]$ tels que $y_{i}\equiv n_{i}^+-n_{i}^-\pmod{I_{G}}$, o $n_{i}^+=\#$ de $j$ tels que $x_{j}=g_{i}$ et $n_{i}^-=\#$ de $j$ tels que $x_{j}=g_{i}^{-1}$. De mme, $d(g)=\sum_{i=1}^n d(g_{i})\cdot z_{i}$ pour certains $z_{i}\in \Z[G]$, mais pour nous n'aurons pas besoin de prŽcisions supplŽmentaires sur les $z_{i}$. Nous prouvons $(+)$ par rŽcurrence sur $N$. Supposons $N=1$. Si $g=g_{k}$ pour un $k=1,\ldots , n$, il n'y a rien ˆ prouver. Si $g=g_{k}^{-1}$, on observe que $d(g_{k}^{-1})=(-1-d(g_{k}^{-1}))\cdot d(g_{k})$ et $(-1-d(g_{k}^{-1}))\equiv -1\pmod{I_{G}}$. Supposons le rŽsultat vrai pour $N-1$ et supposons $g=x\cdot g_{k}$ de longueur $N$. Alors on a
$$\eqalign{d(x\cdot g_{k})&=d(x)+d(g_{k})+d(x)\cdot d(g_{k})\buildrel{\rm H.R}\over =\sum_{i=1}^ny_{i}'\cdot d(g_{i})+d(g_{k})+d(x)\cdot d(g_{k})\cr
&=\sum_{i\ne k}y_{i}'\cdot d(g_{i})+((y_{k}'+1)+d(x))\cdot d(g_k).\cr}$$
En posant $y_{i}=y_{i}'$ si $i\ne k$ et $y_{k}=y_{k}'+1+d(x)\equiv y'_k+1\pmod{I_G}$, le rŽsultat est prouvŽ dans ce cas-lˆ. Enfin, et de mme, si $g=x\cdot g_{k}^{-1}$, on a
$$d(x\cdot  g_{k}^{-1})\buildrel{\rm H.R}\over =\sum_{i\ne k} y_{i}'\cdot d(g_{i})+\underbrace{((y_{k}'-1)-d(g_{k}^{-1})+d(x)\cdot (-1-d(g_{k}^{-1})))}_{\equiv (y_{k}'-1)\pmod I_{G}}\cdot d(g_{k}),$$
ce qui montre aussi le rŽsultat dans ce cas-lˆ. Pour le fait que $d(g)=\sum_{i=1}^n d(g_{i})\cdot z_{i}$, la preuve est similaire.
\smallskip
ConsidŽrons maintenant l'application $\Z^n\buildrel\alpha\over\lra G/G'$ dŽfinie par $\alpha\pmatrix{\l_{1}\cr\vdots\cr l_{n}}=\prod_{i=1}^n g_{i}^{l_{i}}\pmod{G'}$. Cette application est un homomorphisme surjectif de groupes, car $G/G'$ est abŽlien. Son noyau est un sous-groupe de $\Z^n$ d'indice $[G:G']<\infty$; c'est donc aussi un groupe abŽlien libre de rang $n$. Donc c'est l'image d'une application $\Z$-linŽaire injective $\Z^n\buildrel \beta\over \lra\Z^n$ (on dit que la suite exacte $0\to\Z^n\buildrel \beta\over \lra\Z^n\buildrel \alpha\over \lra G/G'\to 1$ est une {\it prŽsentation de $G/G'$}). l'application $\beta$ est dŽcrite par une matrice $(m_{ij})\in M_n(\Z)$ de dŽterminant $[G:G']$. On a donc, pour chaque $j=1,\ldots ,n$ une relation 
$$\tau_j\cdot\prod_{i=1}^n g_i^{m_{ij}}=1\quad\hbox{pour un }\tau_j\in G'.$$
Appliquant $d$ ˆ ces relations, en se souvenant que $d(1)=0$, et en utilisant la relation $(+)$, on trouve pour chaque $j=1,\ldots ,n$ une relation 
$$\sum_{k=1}^n\tau_{jk}\cdot d(g_k)=0\hbox{ avec des }\tau_{jk}\in \Z[G]\hbox{ tels que }\tau_{jk}\equiv m_{kj}\pmod{I_G},\eqno {(+*)}$$  
(les $\tau_j$ de la relation prŽcŽdente ne contribuent pas aux $n_i^{+}-n_i^{-}$, car Žtant dans $G'$, ce sont des produits de commutateurs $ghg^{-1}h^{-1}$ qui ont autant de $+$ que de $-$). Posons $\tau=\det(\tau_{jk})$ dŽfini par la formule $\sum_{\sigma\in S_n}{\rm sgn}(\sigma)\tau_{1\sigma(1)}\cdots\tau_{n\sigma(n)}$, et on dŽfinit $(\widetilde{\tau}_{ik})$ la co-matrice de $(\tau_{ik})$ telle que $\widetilde{\tau}_{ik}=(-1)^{i+k}\det(\tau(i,k))$, o $\tau(i,k)$ est la matrice obtenue en biffant dans $(\tau_{ik})$ la $i^e$ colonne et la $j^e$ ligne. Si on rappelle ces choses-lˆ, c'est que $\Z[G]$ n'est pas forcŽment commutatif. Pour appliquer les thŽormes d'algbre linŽaire classiques, il faut passer ˆ l'anneau commutatif (on vient de le voir) $\Z[G]/(I_{G'}\cdot\Z[G])\simeq \Z[G/G']$. On a donc, pour tout $i,k$~:
$$\sum_{j=1}^n\widetilde{\tau}_{ij}\cdot \tau_{jk}=\delta_{ik}\cdot\tau+d_{ik}\hbox{ avec }d_{ik}\in I_{G'}\cdot\Z[G]\hbox{ et $\delta_{ij}$ est le symbole de Kronecker.}$$
En multipliant cette relation par $d(g_k)$ et en sommant sur $k$, on obtient 
$$\tau\cdot d(g_i)=\sum_{j=1}^n \widetilde{\tau}_{ij}\underbrace{\left (\sum_{k=1}^n\tau_{jk}\cdot d(g_k)\right )}_{=0\ \rm rel. (+*)}-\sum_{k=1}^n d_{ik}d(g_k),$$
et donc $\tau\cdot d(g_i)\equiv 0\pmod{\underbrace{I_{G'}\cdot\Z[G]\cdot I_G}_{=I_{G'}\cdot I_G}}$ pour tout $i$. Soit $g\in G$, on se souvient que juste aprs la relation $(+)$, on a vu qu'il existait de $z_i\in \Z[G]$ tels que $d(g)=\sum_{i=1}^n d(g_i)\cdot z_i$, donc $\tau\cdot d(g)=\sum_{i=1}^n\underbrace{\tau\cdot d(g_i)}_{\equiv 0}\cdot z_i$. On trouve alors la relation~:
$$\tau \cdot d(g)\equiv 0\pmod {I_{G'}\cdot I_G}\quad\hbox{pour tout }g\in G.\eqno{(++)}$$
Notant toujours $\overline{x}$ la classe de $x\in \Z[G]$ modulo $I_{G'}\cdot\Z[G]$, on a {\it a priori}  (cf. preuve de $(**)$) $\overline{\tau}=\sum_{t\in T}n_{\overline{t}}\cdot\overline{t}$, avec $n_{\overline{t}}\in \Z$, uniques. Soit $g\in G$. La relation $\tau\cdot d(g)=\tau \cdot (g-1)\equiv 0\pmod {I_{G'}\cdot I_G}$ donne $\tau\cdot g\equiv \tau\pmod{I_{G'}\cdot I_G}$ donc {\it a fortiori} $\pmod {I_{G'}\cdot\Z[G]}$; et donc $\overline{\tau}\cdot\overline{g}=\overline{\tau}$ ou encore 
$$\sum_{t\in T}n_{\overline{t}}\cdot\overline{t}=\sum_{t\in T}n_{\overline{t}}\cdot\overline{t}\cdot\overline{g}=\sum_{t\in T}n_{\overline{t}\cdot\overline{g}^{-1}}\cdot\overline{t},$$
la dernire ŽgalitŽ venant du fait que $T\cdot g$ est une transversale et que $G$ agit transitivement sur les classes ˆ droite de $G$ modulo $G'$; ainsi, $n_{\overline{t}\cdot\overline{g}^{-1}}=n_{\overline{t}}$ pour tout $t\in T$ et $g\in G$. Cela prouve que les $n_{\overline{t}}$ sont tous Žgaux ˆ un mme entier, disons $m$. On a donc 
$$\tau\equiv m\cdot\sum_{t\in T}t\pmod{I_{G'}\cdot\Z[G]},\eqno{(+++)}$$
donc {\it a fortiori} modulo $I_G$. Mais comme $\tau_{ij}\equiv m_{ji}\pmod {I_G}$, pour tout $i,j$, on a $\tau\equiv \det(m_{ij})=[G:G']\pmod{I_G}$. D'autre part $t\equiv 1\pmod{I_G}$ pour tout $t\in T$ et que $|T|=[G:G']$, on a donc 
$$[G:G']=\det(m_{ij})\equiv\tau\buildrel(+++)\over\equiv m\cdot\sum_{t\in T} t\equiv m\cdot [G:G']\pmod {I_G}.$$
Ce qui prouve que $m\equiv 1\pmod {I_G}$ ($\Z[G]/I_G$ est isomorphe ˆ $\Z$ donc intgre), d'o $m=1$, puisque la congruence modulo $I_G$ est l'ŽgalitŽ pour les entiers. L'Žquivalence $(+++)$ devient alors 
$\tau\equiv \sum_{t\in T}t\pmod{I_{G'}\cdot\Z[G]}$, et donc $\tau\cdot d(g)\equiv \left (\sum_{t\in T} t\right )\cdot d(g)\pmod{\underbrace{I_{G'}\cdot\Z[G]\cdot I_G}_{=I_{G'}\cdot I_G}}$. Finalement, l'Žquivalence $(++)$ devient 
$$ \left (\sum_{t\in T} t\right )\cdot d(g)\equiv 0\pmod{I_{G'}\cdot I_G}$$
qui est l'Žquivalence $(*)$ cherchŽe, ce qui prouve le thŽorme.\qed
\bigskip

Maintenant que ce long et joli thŽorme de thŽorie de groupe est prouvŽ, nous pouvons revenir attaquer le problme de front

\bigskip

\th
\medskip

{\sl Soit $K\subset E\subset L$ des corps de nombres. On suppose $L/K$ galoisien. Posons $G={\rm Gal}(L/K)$ et $H={\rm Gal}(L/E)\subset G$. Soit encore $S$, l'ensemble des idŽaux de $K$ ramifiŽs dans $L$ et $S'$ l'ensemble des idŽaux de $E$ ramifiŽs dans $L$. Alors le diagramme suivant commute~:

\vglue 2.5cm 
\psset{xunit=1cm,yunit=1cm}

\rput(6,0){\rnode{b1}{$I_E^{S'}$}}
\rput(6,2){\rnode{a1}{$I_K^S$}}
\rput(9,0){\rnode{b2}{$H/H'=H^{\rm ab}$}}
\rput(9,2){\rnode{a2}{$G/G'=G^{\rm ab}$}}
\ncline[nodesep=3pt]{->}{a1}{a2}
\Aput{$\Phi_{L/K} $} 
\ncline[nodesep=3pt]{->}{b1}{b2}
\Aput{$\Phi_{L/E}$} 
\ncline[nodesep=3pt]{->}{a1}{b1}
\Aput{$i$}
\ncline[nodesep=3pt]{->}{a2}{b2}
\Aput{$V_{G\to H}$}
\medskip

o l'application $i(\aa)=\aa\cdot O_E$ pour tout $\aa\in I_K^S$, et les $\Phi$ on ŽtŽ dŽfini ˆ la fin du Chapitre 10, dans la section Òapplication ˆ des extensions non abŽliennes".

}
{\bf Preuve}

Soient $\P$ un idŽal premier de $K$, $\P\not\in S$, $\gP$ un idŽal premier de $L$ au-dessus de $\P$ et $\sigma={\rm Frob}(\gP/\P)$. Le ThŽorme \aargp\ nous apprend que l'ensemble des classes ˆ droites de $G$ modulo $H$ se dŽcompose sous l'action de $Z(\gP/\P)=<\sigma>$ en $r$ orbites~:
$$C_i=\{H\cdot\tau_i,H\cdot\tau_i\cdot\sigma,\ldots, H\cdot\tau_i\cdot\sigma^{f_i-1}\}\quad i=1,\ldots ,r,$$
o, en posant $\P_i=\tau_i(\gP)\cap E$, $ i=1,\ldots ,r,$ on a $i(\P)=\P\cdot O_E=\P_1\cdots\P_r$, et pour chaque $i$, $f_i=f(\P_i/\P)$. On a vu ˆ la fin du Chapitre 10 que $\Phi_{L/K}(\P)=\sigma\bmod {G'}$.  En s'inspirant de la preuve du ThŽorme~\argf, on a, pour $i=1,\ldots ,r$
$${\rm Frob}_{L/E}(\tau_{i}(\gP)/\P_{i})={\rm Frob}_{L/K}((\tau_{i}(\gP)/\P)^{f_{i}}=(\tau_{i}\cdot\sigma\cdot\tau_{i}^{-1})^{f_{i}}=\tau_{i}\cdot\sigma^{f_{i}}\cdot\tau_{i}^{-1}.$$
Ainsi,
$$\Phi_{L/E}(i(\P))=\prod_{i=1}^r\tau_{i}\cdot\sigma^{f_{i}}\cdot\tau_{i}^{-1}\bmod H'.\eqno{(*)}$$
D'autre part, on se souvient que $V_{G\to H}(\sigma\bmod G')=\prod h_{t}\bmod H'$, o pour chaque $t\in T$, $t\cdot\sigma=h_{t}\cdot t'$ ($h_{t}\in H, t'\in T$), o $T$ est n'importe quelle transversale ˆ droite de $H$ dans $G$. Alors prenons par exemple
$$T=\{\tau_{1},\tau_{1}\sigma,\ldots ,\tau_{1}\sigma^{f_{1}-1},\ldots ,\tau_{r},\tau_{r}\sigma,\ldots ,\tau_{r}\sigma^{f_{r}-1}\}.$$
Remarquons maintenant que, pour $0\leq i<f_{j}-1$, $j=1,\ldots , r$, on a $(\tau_{j}\cdot\sigma^{i})\cdot\sigma=\tau_{j}\cdot\sigma^{i+1}$, et donc $h_{\tau_{j}\cdot\sigma^{i}}=1$. Sinon, pour $j=1,\ldots ,r$, on a $(\tau_{j}\cdot\sigma^{f_{j}-1})\cdot\sigma=\tau_{j}\cdot\sigma^{f_{j}}=h_{\tau_{j}\cdot \sigma^{f_{j}-1}}\cdot \tau_{j}$. Ce qui veut dire que $h_{\tau_{j}\cdot \sigma^{f_{j}-1}}=\tau_{j}\cdot\sigma^{f_{i}}\cdot \tau_{j}^{-1}$. Et enfin, 
$$V_{G\to H}(\Phi_{L/K}(\P))=V_{G\to H}(\sigma\bmod G')=\prod_{i=1}^r\tau_{i}\cdot\sigma^{f_{i}}\cdot\tau_{i}^{-1}\bmod H'.$$
Ce qui donne, combinŽ avec $(*)$,  $\Phi_{L/E}(i(\P))=V_{G\to H}(\Phi_{L/K}(\P))$.\qed

\bigskip
\th {\bf  (ThŽorme des idŽaux principaux de la thŽorie du corps de classe)}
\medskip
{\sl Soit $K$ un corps de nombres. Notons $E$ le corps de Hilbert de $K$. Alors tout idŽal de $K$ {\it capitule} dans $E$, i.e. devient principal dans $E$.

}

{\bf Preuve}

Posons $L$ le corps de Hilbert de $E$. L'extension $L/K$ est abŽlienne. En effet, soit $L^{\rm alg}$ la cl™ture algŽbrique de $L$ et $\varphi\, :\, L\to L^{\rm alg}$, un $K$-morphisme. Il faut voir que $\varphi(L)\subset L$. L'extension $E/K$ Žtant galoisienne, on a donc $\varphi(E)=E$. Donc $\varphi(L)$ est une extension de $E$. En faisant des allers et venues avec $\varphi$ pour des groupes d'automorphismes (et d'inerties), on voit que $\varphi (L)/E$ est une extension abŽlienne non ramifiŽe (${\rm Gal}(\varphi(L)/E)=\varphi\cdot{\rm Gal}(L/E)\cdot\varphi^{-1}$). Donc $\varphi(L)\subset$ le corps de Hilbert de $E=L$ (en vertu du Corollaire \argeq), ce qu'il fallait voir.

Notons alors $G={\rm Gal}(L/K)$ et $H={\rm Gal}(L/E)$. La thŽorie de Galois nous dit alors que ${\rm Gal}(E/K)\simeq G/H$. L'extension $E/K$ est la plus grande sous-extension abŽlienne de $L/K$;  c'est Žvident, car $L/K$ est une extension non ramifiŽe (il en est donc de mme de toute-sous extension) et nous savons que $L/K$ contient toute extension abŽlienne non-ramifiŽe de $K$ (cf. Corollaire \argeq). Donc, $G/H$ est le plus grand quotient abŽlien de $G$, et alors, par dŽfinition, $H=G'$. De plus $H=G'$ est abŽlien, donc $G''=\{1\}$. Puisque $L/E$ est abŽlienne $\Phi_{L/E}$ est l'application d'Artin usuelle; et on a aussi $\Phi_{L/K}=\Phi_{E/K}'$ (ˆ nouveau, car $E$ est la plus grande sous-extension abŽlienne de $L/K$ et vu la fin du Chapitre 10). On veut utiliser le ThŽorme \argfe. Dans notre cas, $S=S'=\emptyset$, ainsi $I_K^S=I_K$ et $I_E^{S'}=I_E$. Ce mme thŽorme nous dit qu'on a le diagramme commutatif~:

\vglue 2.5cm 
\psset{xunit=1cm,yunit=1cm}

\rput(6,0){\rnode{b1}{$I_E$}}
\rput(6,2){\rnode{a1}{$I_K$}}
\rput(9,0){\rnode{b2}{$G'/G''\simeq G'={\rm Gal}(L/E)$}}
\rput(9,2){\rnode{a2}{$G/G'\simeq {\rm Gal}(E/K)$}}
\ncline[nodesep=3pt]{->}{a1}{a2}
\Aput{$\Phi_{E/K} $} 
\ncline[nodesep=3pt]{->}{b1}{b2}
\Aput{$\Phi_{L/E}$} 
\ncline[nodesep=3pt]{->}{a1}{b1}
\Aput{$i$}
\ncline[nodesep=3pt]{->}{a2}{b2}
\Aput{$V_{G\to G'}$}
\medskip

Le ThŽorme \argaq\ nous dit que $\Phi_{E/K}$ et $\Phi_{L/E}$ sont surjectives. Le Corollaire \argeq\  nous dit que $\ker (\Phi_{E/K})=P_K$ et $\ker (\Phi_{L/E})=P_E$, de plus $i(P_K)\subset P_E$. Donc le diagramme prŽcŽdent Òpasse au quotients" donnant un diagramme commutatif~:

\vglue 2.5cm 
\psset{xunit=1cm,yunit=1cm}

\rput(6,0){\rnode{b1}{$I_E/P_E$}}
\rput(6,2){\rnode{a1}{$I_K/P_K$}}
\rput(9,0){\rnode{b2}{$G'/G''$}}
\rput(9,2){\rnode{a2}{$G/G'$}}
\ncline[nodesep=3pt]{->}{a1}{a2}
\Aput{$\overline{\Phi_{E/K} }$} 
\ncline[nodesep=3pt]{->}{b1}{b2}
\Aput{$\overline{\Phi_{L/E}}$} 
\ncline[nodesep=3pt]{->}{a1}{b1}
\Aput{$\overline{i}$}
\ncline[nodesep=3pt]{->}{a2}{b2}
\Aput{$V_{G\to G'}$}
\rput(7.5,-0.3){$\simeq$}
\rput(7.5,1.7){$\simeq$}
\medskip

 L'hypothse du ThŽorme \argfb\ est satisfaite, donc l'homomorphisme $V_{G\to G'}$ est l'homomorphisme trivial, donc $\overline{i}$ est aussi trivial, ce qui veut dire que $i(I_K)\subset P_E$, donc le thŽorme.\qed
 \bigskip
 Maintenant nous pouvons dire que nous avons construit le corps de Hilbert compltement. C'est une certaine satisfaction et on espre que le lecteur aura apprŽciŽ la lecture de ce texte jusqu'ici et est d'accord que cette construction n'est pas totalement ÒpŽdestre''.
 
\vskip4cm 
 


\vfill\eject

\global\advance\chapnomb by 1
\nomb=1

\centerline{\para Chapitre 13 : }
\bigskip

\centerline{\para InterprŽtation idŽlique}
\bigskip
La manire de prŽsenter habituellement  la thŽorie du corps de classe passe par un nouvel objet qu'on appelle les idles. Nous verrons que le ThŽorme \argef\ se traduit dans ce cadre de manire trs simple (cf. Corollaire \arggj\ et ThŽorme \arggp). Seulement, nous trouvons (mais c'est un avis personnel) que cette interprŽtation du corps de classe est plus obscure. C'est-ˆ-dire que l'ŽnoncŽ est trs court (c'Žtait probablement la volontŽ de faire ainsi), mais toutes les difficultŽs ont ŽtŽ pour ainsi dire ``mises sous le tapis''. NŽanmoins, il nous a paru important de faire le lien avec la vision classique de cette magnifique thŽorie.\bigskip
\defi
\medskip
Soit $K$ un corps de nombres. Rappelons qu'on note $\gfP(K), \gfP_{0}(K),\gfP_{\infty}(K)$ l'ensemble des places, des places finies respectivement infinies de $K$. Pour chaque $\P\in\gfP_{0}(K)$, on note $\bbK_{\P}$ le complŽtŽ de $K$ en $\P$, $O_{\P}$ son anneau de valuation, $v_{\P}$ la valuation $\P$-adique et $|x|_{\P}=\N(\P)^{-v_{\P}(x)}$. Si $\P\in\gfP_{\infty}(K)$, on note aussi $\bbK_{\P}$ le complŽtŽ et si $\sigma_{\P}$ est un plongement qui dŽfinit $\P$, on pose
$$|x|_{\P}=\cases{|\sigma_{\P}(x) |&si $\P$ est rŽelle\cr |\sigma_{\P}(x) |^2& si $\P$ est complexe.\cr}$$
 \newcount\gaagar\gaagar=\gaga
Remarquons que dans le cas complexe, on avait dŽfinit $|x|_{\P}$ autrement (cf. page \the\gaagf ) et que maintenant ce n'est plus vraiment une valeur absolue (on n'a pas $|x+y|_{\P}\leq |x|_{\P}+|y|_{\P}$), mais  sans cette dŽfinition (un peu malheureuse, il est vrai), la proposition suivante ne serait pas vraie...

On note aussi
$$O_{\P}^*=\{x\in\bbK_{\P}\mid |x|_{\P}=1\}\quad\forall\P\in\gfP(K).$$
Si $\P\in\gfP_{0}(K)$, on a $O_{\P}^*=U_{\P}$, le groupe des unitŽs de $O_{\P}$; et si $\P\in\gfP_{\infty}(K)$, $O_{\P}^*$ est le cercle unitŽ si $\P$ est complexe et $\{\pm 1\}$ si $\P$ est rŽelle.
\bigskip
\prop
\medskip
{\sl 
Pour tout $x\in K^*$, on a 
$$\prod_{\P\in \gfP(K)}|x|_{\P}=1.$$

}

{\bf Preuve}

Soit donc $x\in K^*$. Il existe $k_{1},\ldots ,k_{r}\in \Z$ et $\P_{1},\ldots,\P_{r}\in\gfP_{0}(K)$  tels que $x\cdot O_{K}=\P_{1}^{k_{1}}\cdots \P_{r}^{k_{r}}$. D'autre part, $|x|_{\P}=1$ si $\P\ne\P_{i} 
$, $i=1,\ldots ,r$. Et donc, $\prod_{\P\in \gfP(K)}|x|_{\P}=\prod_{i=1}^{r}|x|_{\P_{i}}\cdot\prod_{\P\in \gfP_{\infty}(K)}|x|_{\P}$. Or, d'une part, $\prod_{i=1}^{r}|x|_{\P_{i}}=\prod_{i=1}^{r}(\N(\P_{i}))^{-k_{i}}=\N(\P_{1}^{k_{1}}\cdots\P_{r}^{k_{r}})^{-1}=|N_{K/\Q}(x)|^{-1}$. D'autre part, si $S$ est l'ensemble de tous les plongements de $K$ dans $\C$, on a $\prod_{\P\in \gfP_{\infty}(K)}|x|_{\P}=\prod_{\sigma\in S}|\sigma(x)|=\left |\prod_{\sigma\in S}\sigma(x)\right |=\left | N_{K/\Q}(x)\right |$. Finalement,
$$\prod_{\P\in \gfP(K)}|x|_{\P}=|N_{K/\Q}(x)|^{-1}\cdot \left | N_{K/\Q}(x)\right |=1$$
\qed
\bigskip
\headline={\hfill \phantom{ouuh}\hfill}
\defi
\medskip
Soit $K$ un corps de nombres. Un {\it adle} de $K$ est une famille $(\alpha_{\P})_{\P\in \gfP(K)}$ telle que $\alpha_{\P}\in \bbK_{\P}$ pour tout $\P\in \gfP(K)$ et $\alpha_{\P}\in O_{\P}$ pour presque tout $\P\in\gfP_{0}(K)$. Le syntagme {\it pour presque tout} signifie Òpour tout ŽlŽment de l'ensemble considŽrŽ sauf, Žventuellement, un nombre fini''. Nous rŽsumerons tout cela par ``ppt''. L'ensemble des adles de $K$ se note $\bbA_{K}$. En rŽsumŽ, on a 
$$\bbA_{K}=\left \{ (\alpha_{\P})_{\P}\in \prod_{\P\in \gfP(K)}\bbK_{\P}\mid \alpha_{\P}\in O_{\P}\ {\rm ppt}\ \P\right \}.$$
De mme, on dŽfinit l'ensemble des {\it idles} de $K$
$$\bbI_{K}=\left \{ (\alpha_{\P})_{\P}\in \prod_{\P\in \gfP(K)}\bbK_{\P}^*\mid \alpha_{\P}\in O_{\P}^*\ {\rm ppt}\ \P\in\gfP_{0}(K)\right \}.$$
 \newcount\gaagas\gaagas=\gaga
Remarquons qu'historiquement, le mot ``idle'' est antŽrieur au mot  ``adle''. Idle vient de $\underline{\rm id}$eal  $\underline{\rm ele}$ments et adle vient de $\underline{\rm ad}$ditive id$\underline{\hbox{le}}$. L'accent grave vient de la ressemblance avec le prŽnom fŽminin et que les inventeurs du concept Žtaient certainement francophones. Mais attention, on dit  ``{\bf un} adle'' et ``{\bf un} idle''.

Clairement, $\bbA_{K}$ est un sous-anneau de l'anneau produit $ \prod_{\P\in \gfP(K)}\bbK_{\P}$ et $\bbI_{K}$ est un sous-groupe du groupe produit $ \prod_{\P\in \gfP(K)}\bbK_{\P}^*$.

Si $S$ est une partie finie de $\gfP(K)$ qui contient $\gfP_{\infty}(K)$, on note 
$$\bbA_{K}(S)=\prod_{\P\in S}\bbK_{\P}\times \prod_{\P\not\in S}O_{\P}\quad\hbox{et}\quad\bbI_{K}(S)=\prod_{\P\in S}\bbK_{\P}^*\times \prod_{\P\not\in S}O_{\P}^*.$$
Dans le cas o $S=\gfP_{\infty}(K)$, on Žcrira $\bbA_{K}(\infty)$ et $\bbI_{K}(\infty)$.

Les $\bbA_{K}(S)$ sont des sous-anneaux des $\bbA_{K}$ et les $\bbI_{K}(S)$ sont des sous-groupes des $\bbI_{K}$. Ils sont filtrants supŽrieurement (toute rŽunion finie possde un ŽlŽment qui contient cette rŽunion) et leur rŽunion donnent $\bbA_{K}$ respectivement $\bbI_{K}$.

\headline={\hfill \smcap InterprŽtation idŽlique\hfill}
On munit $\bbA_{K}$ et $\bbI_{K}$ d'une topologie comme suit~:

Pour $\bbA_{K}$ une base d'ouverts est l'ensemble des parties de la forme  $\prod_{\P}V_{\P}$, o $V_{\P}$ est un ouvert de $\bbK_{\P}$ pour tout $ \P\in\gfP(K)$, et $V_{\P}=O_{\P}$ ppt $\P$. Remarquons  que si $\P\in\gfP_{0}(K)$,  $O_{\P}=\{x\in \bbK_{\P}\mid |x|_{\P}\leq 1\}$ est un ouvert malgrŽ les apparences :  pour la topologie $\P$-adique de $\bbK_{\P}$, ``les boules fermŽes sont ouvertes et rŽciproquement'' on dit que c'est un espace topologique {\it totalement discontinu}. De plus, dans $\bbK_{\P}$, $O_{\P}$ est compact (cf. [Fr-Tay, II.4, rel. (3.29), p.86]).

Pour $\bbI_{K}$, de mme, une base d'ouverts est l'ensemble des parties de la forme  $\prod_{\P}V_{\P}$, o $V_{\P}$ est un ouvert de $\bbK_{\P}^*$ pour tout $ \P\in\gfP(K)$, et $V_{\P}=O_{\P}^*$ ppt $\P$.

Pour cette topologie, la topologie induite sur $\bbA_{K}(S)$ (resp. $\bbI_{K}(S)$) est identique ˆ la topologie produit, et $\bbA_{K}(S)$ (resp. $\bbI_{K}(S)$) est ouvert dans $\bbA_{K}$ (resp. dans $\bbI_{K}$). Comme les $\bbA_{K}(S)$  et les $\bbI_{K}(S)$) sont localement compacts (ils sont les produits finis d'espaces localement compacts avec un produit d'espaces compacts), alors $\bbA_{K}$ et $\bbI_{K}$ sont localement compacts (pour tout $x\in \bbA_{K}$, il existe $S$ tel que $x\in\bbA_{K}(S)$). Remarquons que dans notre acceptation de termes ``compact'' et ``localement compact'', nous incluons la propriŽtŽ d'tre sŽparŽ. La  rŽunion des $\bbA_{K}(S)$ (resp. des $\bbI_{K}(S)$) est filtrante et on peut voir $\bbA_{K}$, (resp. $\bbI_{K}$) comme limite inductive des $\bbA_{K}(S)$, (resp. des $\bbI_{K}(S)$). En outre, cette topologie induit sur $\bbA_{K}$ (resp. sur $\bbI_{K}$) une structure de d'anneau (resp. de groupe) topologique, car chaque $\bbA_{K}(S)$ (resp. $\bbA_{K}(S)$) en est un.
\bigskip
{\soustitre Remarque}
\medskip
La topologie de $\bbI_{K}$ n'est pas celle induite par celle de $\bbA_{K}$. En effet, choisissons pour chaque $\P\in \gfP_{0}(K)$ une uniformisante $\pi_{\P}$ de $\bbK_{\P}$ et posons $x_{\P}$ l'idle $(1,\ldots ,1,\pi_{\P},1,\ldots ,1)$, alors, pour la topologie des adles,  n'importe quel voisinage de $(1,\ldots ,1)$ contient presque tous les $x_{\P}$, donc $(1,\ldots ,1)$ est un point d'accumulation de la famille des $x_{\P}$. En revanche, pour la topologie des idles, $\bbI_{K}(\infty)$ est un voisinage de $(1,\ldots ,1)$ et il ne contient aucun des $x_{\P}$.

Cette remarque est importante pour nous pousser ˆ la prudence quand nous raisonnerons sur ces objets. En revanche, nous avons un rŽsultat qui donne un lien entre les deux topologies~:
\bigskip
\lem
\medskip
{\sl
L'application injective
$$\eqalign{j\, :\;\bbI_{K}&\longrightarrow\bbA_{K}\times \bbA_{K}\cr x&\longmapsto (x,x^{-1})\cr}$$
induit la topologie de $\bbI_{K}$ (identifiŽ ˆ un sous-espace de $\bbA_{K}\times \bbA_{K}$).

}

{\bf Preuve}

Il suffit de montrer que $j$ est continue et que tout ouvert de $\bbI_{K}$ est l'image rŽciproque d'un ouvert de $\bbA_{K}\times \bbA_{K}$. On remarque que 
$$j^{-1}\left (\left (\prod_{\P\in S} V_{\P}\times \prod_{\P\not\in S} O_{\P}\right )\times \left (\prod_{\P\in S} V'_{\P}\times \prod_{\P\not\in S} O_{\P}\right )\right )=\prod_{\P\in S}(V_{\P}\setminus\{0\})\cap (V'_{\P}\setminus\{0\})^{-1}\times\prod_{\P\not\in S}O^*_{\P}$$
qui est un ouvert de base de $\bbI_{K}$, car $(V_{\P}\setminus\{0\})$ et $(V'_{\P}\setminus\{0\})^{-1}$ sont des ouverts de $\bbK_{\P}^*$ si $V_{\P}$ et $V'_{\P}$ sont des ouverts de $\bbK_{\P}$ (le passage de $x$ ˆ $x^{-1}$ est une application bi-continue). Donc l'application $j$ est continue. D'autre part,
$$\prod_{\P\in S}V_{\P}\times\prod_{\P\not\in S}O_{\P}^*=j^{-1}\left (\left (\prod_{\P\in S} V_{\P}\times \prod_{\P\not\in S} O_{\P}\right )\times \left (\prod_{\P\in S} V_{\P}^{-1}\times \prod_{\P\not\in S} O_{\P}\right )\right ),$$
ce qui montre le lemme.\qed
\bigskip
\lem
\medskip
{\sl Soit $G$ un groupe topologique sŽparŽ et $H$ un sous-groupe de $G$. Si $H$ est discret, alors il est fermŽ

}
{\bf Preuve}

Mettons que $G$ soit notŽ multiplicativement. Soit $a\in G\setminus H$ et $V$ un ouvert tel que $V\cap H=\{1\}$ (c'est possible puisque $H$ est discret). On choisit un ouvert $U$, voisinage de 1 tel que $a\cdot U\cdot U^{-1}\cdot a\subset V$. C'est toujours possible, car l'application $(x,y)\mapsto a\cdot x\cdot y^{-1}\cdot a^{-1}$ est une application continue de $G\times G$ dans $G$. Supposons que $x,y\in H\cap a\cdot U$. Alors $xy^{-1}\in H$ et $xy^{-1}\in a\cdot U\cdot U^{-1}\cdot a\subset V$. Donc $xy^{-1}=1$, i.e. $x=y$. Donc le voisinage $a\cdot U$ de $a$ contient au plus un ŽlŽment de $H$. Puisque $G$ est sŽparŽ, on peut, en restreignant $a\cdot U$ si nŽcessaire, trouver un voisinage de $a$ qui ne rencontre pas $H$. Donc $G\setminus H$ est ouvert et donc $H$ est fermŽ.\qed
\bigskip
\defi
\medskip
Pour chaque $\P\in\gfP(K)$, on suppose choisi une identification de $K$ avec un sous-corps de $\bbK_{\P}$. Si $x\in K$, on lui associe diagonalement l'ŽlŽment $(x,x,\ldots, x)\in\prod_{\P}\bbK_{\P}$. C'est un adle, car $|x|_{\P}\leq 1$ ppt $\P$. On obtient un homomorphisme d'anneau $K\to \bbA_{K}$. De mme, si $x\in K^*$, $|x|_{\P}=1$ ppt $\P$. Donc $x\mapsto (x,\ldots ,x)$ dŽfinit un homomorphisme de groupe de $K^*$ dans $\bbI_{K}$. Nous associerons $K$ (resp. $K^*$) avec son image dans $\bbA_{K}$ (resp. dans $\bbI_{K}$) que nous appelleront les adles (resp. les idles) {\it principaux}.
 \newcount\gaagat\gaagat=\gaga
\bigskip\goodbreak
\lem
\medskip
{\sl Les adles principaux forment un sous-anneau discret (donc fermŽ en vertu du lemme prŽcŽdent) de $\bbA_{K}$; et les idles principaux forment un sous-groupe discret (et donc fermŽ) de $\bbI_{K}$.

}
{\bf Preuve}

L'inclusion $\bbI_{K}\to \bbA_{K}$ est continue (on peut la voir comme la composition des applications $x\mapsto (x,x^{-1})\mapsto x$ qui est continue en vertu du Lemme \argfr ), donc il suffit de montrer que $K$ est discret dans $\bbA_{K}$. En effet,  l'image rŽciproque de $K$ de cette inclusion qui est $K^*$ serait alors discret dans $\bbI_{K}$ (l'image rŽciproque d'un ouvert ne rencontrant pas $(1,\ldots, 1)$ est un ouvert ne rencontrant pas $(1,\ldots, 1)$). 

ConsidŽrons alors l'ensemble 
$$N=\prod_{\P\in \gfP_{\infty}(K)}\{\alpha\in \bbK_{\P}\mid |\alpha|_{\P}<1\}\times \prod_{\P\in \gfP_{0}(K)}O_{\P}.$$
Clairement $N$ est un voisinage de $(0,0,\ldots, 0)$ dans $\bbA_{K}$. Si $x\in K\cap N$, alors $\prod_{\P\in\gfP(K)}|x|_{\P}<1$ (car les ŽlŽments de $O_{\P}$ sont tels que $|\alpha |_{\P}\leq 1$). Mais en vertu  de la formule du produit (Proposition \argfp) le seul $x\in K$ possible  est $x=(0,\ldots, 0)$. Cela montre que  $(0,\ldots ,0)$ est isolŽ, donc par translation que $K$ est discret dans $\bbA_{K}$.\qed
\bigskip
\lem
\medskip
{\sl Soit $G$ un groupe topologique notŽ multiplicativement. Supposons que $\{1\}$ soit fermŽ. Alors $G$ est sŽparŽ.

}
{\bf Preuve}

Puisque $\{1\}$ est fermŽ, alors $\{x\}$ est fermŽ pour tout $x\in G$ par continuitŽ.  On montre que  $\{1\}$ et  $\{x\}$ peuvent tre sŽparŽs.  Posons $V=G\setminus\{x\}$ qui est ouvert. Par continuitŽ de l'application $(x,y)\mapsto x\cdot y$, il existe $U$ un voisinage de 1 tel que $U\cdot U\subset V$. Quitte ˆ remplacer $U$ par $U\cap U^{-1}$, on peut supposer que $U=U^{-1}$. Alors $U$ et $x\cdot U$ sont des voisinages de 1 et $x$ sont disjoints~: si $U\ni u=x\cdot v\in x\cdot U$, alors $x=u\cdot v^{-1}\in U\cdot U\subset V$, ce qui est une contradiction.\qed
\bigskip\goodbreak

{{\soustitre Proposition-DŽfinition ({\uppercase\expandafter{\the\chapnomb}}.{\the\nomb})}\global\advance\nomb by 1}
\medskip
{\sl On note $C_{K}$ le groupe (multiplicatif) quotient $\bbI_{K}/K^*$ muni de la topologie quotient. On appelle $C_{K}$ le {\it groupe des classes d'idles}. On considre de mme le groupe (additif) quotient $\bbA_{K}/K$, qu'on appellera le {\it groupe des classes d'adles}. On affirme alors que ces deux groupes sont localement compacts.

}

 \newcount\gaagau\gaagau=\gaga
{\bf Preuve}

Ce fait vient du fait que $K$ (resp. $K^*$) est fermŽ (et normal) dans $\bbA_{K}$ (resp. dans $\bbI_{K})$ en vertu du lemme prŽcŽdent.  Supposons en toute gŽnŽralitŽ que $G$ soit un groupe topologique localement compact et que $H$ soit  un sous-groupe normal fermŽ. Rappelons que la topologie quotient sur $G/H$ est la plus fine telle que l'application $\pi\, :\, G\to G/H$ soit continue et donc $U\subset G/H$ est ouvert si et seulement si $\pi^{-1}(U)$ est ouvert. Il faut dŽjˆ montrer que $G/H$ est sŽparŽ. Pour cela, il suffit de montrer comme vu au lemme prŽcŽdent (Lemme \aargfu) que $\{\overline{1}\}$ est fermŽ dans $G/H$. C'est Žvident car $\pi^{-1}(G/H\setminus \{\overline{1}\})=G\setminus H$ qui est ouvert puisque $H$ est fermŽ. Donc $G/H$ est sŽparŽ. Maintenant, l'application $\pi$ est ouverte, car $\pi^{-1}(\pi(U))=H\cdot U=\bigcup_{x\in H}x\cdot U$ qui est ouvert si $U$ est ouvert. Enfin, pour montrer que $G/H$ est localement compact, il suffit par translation de trouver un voisinage compact de $\overline{1}$. Puisque, par hypothse, $G$ est localement compact, on considre $U$ un voisinage compact de $1$. Puisque $\pi$ est une application ouverte, $\pi(U)$ est un voisinage de $\overline{1}$, il est en outre compact en vertu de la propriŽtŽ bien connue que l'image directe d'un compact par une application continue dans un espace sŽparŽ est compact. Donc $G/H$ est localement compact.\qed
\bigskip
Maintenant nous allons faire quelques investigations en vue de montrer que $\bbA_{K}/K$ est en fait compact. Nous montrerons aussi que $\bbI_{K}/K^*$ ne l'est en revanche pas. Tout d'abord voici un forme particulire du thŽorme chinois~:
\bigskip
\lem
\medskip
{\sl Soit $\P_1,\ldots ,\P_n\in \gfP_0(K)$, $\varepsilon_1,\ldots ,\varepsilon_n$ des nombres rŽels positifs et, pour chaque $i=1,\ldots ,n$, $\alpha_i\in \bbK_{\P_i}$. Alors il existe $\beta\in K$ tel que 
$$|\beta-\alpha_i|_{\P_i}\leq\varepsilon_i,\quad i=1,\ldots , n\hbox{ et }|\beta|_\qq\leq 1\hbox{ pour tout  idŽal premier $\qq\ne\P_i$, $ i=1,\ldots , n$.} $$

}
{\bf Preuve}

Puisque $K$ est dense comme dans chaque $K_\P$, on peut supposer que $\alpha_i\in K$ pour tout $ i=1,\ldots , n$. Il existe $m\in\Z$ et $\beta_1,\ldots ,\beta_n\in O_K$ tel que $\alpha_i={\beta_i\over m}$, pour $i=1,\ldots , n$. En effet, montrons-le pour $\alpha_1$, on prend ensuite un dŽnominateur commun~: puisque $\alpha_1$ est algŽbrique, il existe $a_k,a_{k-1},\ldots ,a_0\in\Z$ tels que $a_k\alpha_1^k+a_{k-1}\alpha^{k-1}+\cdots +a_0=0$. En multipliant cette dernire ŽgalitŽ par $a_k^{k-1}$, on montre que $a_k\cdot \alpha_1\in O_K$ et le tour est jouŽ. ConsidŽrons $\qq_1,\ldots ,\qq_s$ les idŽaux premiers distincts des $\P_i$ qui divisent $m$. Par le thŽorme chinois, il existe $\gamma\in O_K$ tel que $| \gamma-\beta_i|_{\P_i}\leq |m|_{\P_i}\cdot\varepsilon_i$, $\forall\, i=1,\ldots ,n$ et $|\gamma|_{\qq_j}\leq |m|_{\qq_j}$, pour $\forall\, j=1,\ldots , s$. Alors, on voit facilement que $\beta={\gamma\over m}$ rŽpond aux exigences du lemme.\qed
\bigskip\goodbreak
\lem
\medskip
{\sl On a les deux ŽgalitŽs~:
\art{a)}$\bbA_{K}(\infty)+K=\bbA_{K}$

\art{b)}$\bbA_{K}(\infty)\cap K=O_{K}$ (o $O_{K}$ est Žvidemment vu comme le plongement diagonal de $O_{K}$ dans $\bbA_{K}$).

}
{\bf Preuve}

Montrons b). L'inclusion $O_{K}\subset \bbA_{K}(\infty)\cap K$ est claire. Inversement, si $x\in \bbA_{K}(\infty)\cap K$, alors $x\in O_\P$ pour tout $\P\in\gfP_0(K)$, c'est-ˆ-dire $|x|_\P\leq 1$ pour tout $\P\in\gfP_0(K)$. Cela veut dire que  $x\in O_K$.

Montrons a). Il faut donc montrer que $\forall\, (a_\P)_\P\in\bbA_K$, $\exists x\in K$ tel que $|a_\P-x|_\P\leq 1$, $\forall \P\in \gfP_{0}(K)$. L'ensemble des $\P\in\gfP_0(K)$ tel que $| a_\P|_\P>1$ est fini. Notons 
$$T:=\{ p\in\gfP_0(\Q)\mid\exists\, \P\in\gfP_0(K), \P|p,  |a_\P|_\P>1\}\hbox{ et }S:=\{\P\in\gfP_0(K)\mid \P|p\hbox{ pour un } p\in T\}.$$
Alors $T$ et $S$ sont finis. Posons $m=\left (\prod_{p\in T} p\right )^k$, avec $k\in\N$ assez grand pour que $|m\cdot a_\P|_\P\leq 1$ pour tout $\P\in \gfP_0(K)$. Par le lemme prŽcŽdent, il existe $\beta\in K$ tel que $|m\cdot a_\P-\beta|_\P\leq |m|_\P$ pour tout $\P\in S$ et $|\beta|_\qq\leq 1$ pour tout $\qq\in\gfP_0(K)\setminus S$. Alors $x:={\beta\over m}$ rŽpond ˆ la question~: si $\P\in S$, $|m\cdot a_\P-\beta|_\P\leq |m|_\P$ implique bien sžr  $|a_\P-x|_\P\leq 1$; et si $\qq \in \gfP_0(K)\setminus S$, $|a_\qq-{\beta\over m}|_\qq={1\over| m|_\qq}\cdot|m\cdot a_\qq-\beta|_\qq=|m\cdot a_\qq-\beta|_\qq\leq \max(|m\cdot a_\qq|_\qq,|\beta|_\qq)\leq 1$. Ce qui achve la preuve du lemme.\qed
\bigskip
\th

{\sl Soit $K$ un corps de nombres. Alors $\bbA_K/K$ est compact (on dit alors que $K$ est {\it co-compact} dans $\bbA_K$).

}

{\bf Preuve}

Rappelons le fait suivant~: supposons que $[K:\Q]=r+2s$. L'application
$$\eqalign{v\ : K&\to \R^r\times \C^s\simeq \prod_{\P\in \gfP_{\infty}(K)}\bbK_{\P}\cr
x&\mapsto (\sigma_{1},\ldots ,\sigma_{r},\sigma_{r+1},\ldots ,\sigma_{r+s})\cr}$$
est une application telle que $v(O_{K})$ est un $\Z$-rŽseau plein (les $\sigma_{i}$ sont les plongements de $K$ dans $\C$.  Et si $\omega_{1},\ldots ,\omega_{n}$ est une $\Z$-base de $O_{K}$, alors l'ensemble
$$F_{\infty}=\{x\in\prod_{\P\in \gfP_{\infty}(K)}\bbK_{\P}\mid x=\sum_{i=1}^nt_{i} \cdot v(\omega_{i})\ 0\leq t_{i}<1\}$$
est un parallŽlotope fondamental (voir Lemme \argt). Soit $F=F_{\infty}\times\prod_{\P\in \gfP_{0}(K)}O_{\P}$. Il est Žvident que l'adhŽrence $\overline{F}$ de $F$ est compacte dans $\bbA_{K}$ par le lemme de Tychonov. 

Soit $x\in \bbA_{K}$. Notons $\overline{x}$ la classe de $x$ dans $\bbA_{K}/K$. La partie a) du  Lemme \argfx, nous assure l'existence de $y\in \bbA_{K}(\infty)$ tel que $\overline{x}=\overline{y}$. Notons $y=(y_{\infty},y_{0})$, avec $y_{\infty}\in \prod_{\P\in \gfP_{\infty}(K)}\bbK_{\P}$. Par ce qui prŽcde, il existe $l\in O_{K}$ et $z_{\infty}\in F_{\infty}$ tel que $y_{\infty}=z_{\infty}+l_{\infty}$, o $l_{\infty}=v(l)$. On a donc $y=(y-l)+l$ et donc $\overline{x}=\overline{y}=\overline{y-l}$ et $y-l\in F\subset\overline{F}$. Ainsi, la restriction de la projection $\bbA_{K}\to \bbA_{K}/K$ ˆ $\overline{F}$ est surjective et bien sžr continue. Puisque $\bbA_{K}/K$ est localement compact (Proposition-DŽfinition \argfv) et l'image d'une application continue d'un compact d'un espace sŽparŽ dans un autre est compact. Cela montre que $\bbA_{K}/K$ est compact.\qed

\bigskip

Nous allons maintenant prouver que $C_{K}=\bbI_{K}/K^*$ n'est pas compact.
\bigskip
\defi

Soit $a=(a_{\P})_{\P}\in \bbI_{K}$. Alors le produit $\prod_{\P\in\gfP_{K}}|a_{\P}|_{\P}:=|a|$ est bien dŽfini, car c'est un produit fini, puisque $|a_{\P}|_{\P}=1$ ppt $\P$. On appellera $|a|$ le {\it volume} de $a$. L'application $\bbI_{K}\to \R_{+}^*$, $a\mapsto |a|$ est clairement un homomorphisme (puisque chacune des normes est un homomorphisme), surjective ($K$ possde au moins une place infinie $\P$ et on a $\R_{+}^*\subset \bbK_{\P}^*$ dans lequel on choisit un ŽlŽment et on met 1 aux autres places), continue (il suffit de montrer que la restriction aux ouverts fondamentaux, qui sont du type $\prod_{\P\not\in S}O_{\P}^*\times\prod_{\P\in S}N_{\P}$, o les $N_{\P}$ sont des ouverts de $\bbK_{\P}^*$ est continue, et c'est clairement le cas, puisque chacune des normes est continue). Le noyau de cette application est un sous-groupe fermŽ, car dans dans des groupes topologiques sŽparŽs, la prŽ-image d'un fermŽ est un fermŽ. On note ce noyau $\bbI_{K}^0$ et on appelle ce sous-groupe les {\it idles spŽciaux}. Par la formule du produit (Proposition \argfp), on a $K^*\subset \bbI_{K}^0$. On note alors $C_{K}^0=\bbI_{K}^0/K^*$ qui est un sous-groupe fermŽ de $C_{K}=\bbI_{K}/K^*$.

 \newcount\gaagav\gaagav=\gaga
\bigskip

\lem

{\sl Il existe un isomorphisme de groupe topologique $C_{K}\simeq C_{K}^0\times \R_{+}^*$. En particulier, $C_{K}$ n'est pas compact.

}

{\bf Preuve}

L'application $\bbI_{K}\to \R_{+}^*$, $a\mapsto |a|$ vue ˆ la dŽfinition prŽcŽdente admet une section continue~: si $n=[K:\Q]$, 
$$\eqalign{\R_{+}^*&\to \bbI_{K}\cr
t&\mapsto (\underbrace{t^{1\over n},\ldots , t^{1\over n}}_{\rm places\ inf.},\underbrace{1,\ldots ,1}_{\rm places\ finies}).\cr}$$
Donc, $\bbI_{K}\simeq \bbI_{K}^0\times\R_{+}^*$. Puisque $K^*$ est dans le noyau, l'application $|\cdot |$ induit aussi un homomorphisme surjectif continu $C_{K}\to \R_{+}^*$, et la section prŽcŽdente donne aussi une section ici, ce qui montre notre lemme.\qed

\bigskip

Nous voyons donc (mme si cela a dŽjˆ ŽtŽ vu) que les topologies idŽliques et adŽliques sont bien diffŽrentes. Nous allons maintenant prouver que $C_{K}^0$ est en revanche compact. Nous allons voir que la compacitŽ de cet espace est Žquivalente au fait que le groupe des classe $I_{K}/P_{K}$ est fini et au thŽorme de Dirichlet sur les unitŽs de $K$, deux rŽsultats que nous connaissons bien ! mais tout d'abord deux petits lemmes de topologie des groupes~:
\bigskip

\lem

{\sl Soit $G$ un groupe topologique, $K$ une partie compacte de $G$ et $O\supset K$ un voisinage de $K$. Alors il existe $U$ un voisinage ouvert de $1$ tel que $UK\subset O$.
}

{\bf Preuve}

Soit $x\in K$. Alors il existe $V_{x}$ voisinage de $1$ tel que $V_{x}x\subset O$ (par exemple $V_{x}= Ox^{-1}$ ). Puisque la multiplication $(x,y)\mapsto x\cdot y$ est par dŽfinition continue, il existe $U_{x}$ voisinage ouvert de $1$ tel que $U_{x}\cdot U_{x}\subset V_{x}$. Il est clair que $\{ U_{x} x\}_{x\in K}$ est un recouvrement de $K$. Par compacitŽ, il existe $x_{1},\ldots ,x_{n}$ tels que $K\subset \bigcup_{i=1}^n U_{x_{i}}x_{i}$. Posons $U=\bigcap_{i=1}^n U_{x_{i}}$. Alors $UK\subset O$. En effet, soit $t=u\cdot k\in UK$. Puisque $k\in K\subset  \bigcup_{i=1}^n U_{x_{i}}x_{i}$, il existe $i$ et $u_{i}\in U_{x_{i}}$ tel que $k=u_{i}\cdot x_{i}$. Donc $t=u\cdot u_{i}\cdot x_{i} \in U_{x_{i}}U_{x_{i}}x_{i}\subset V_{x_{i}}x_{i}\subset O$. Cela prouve le lemme.\qed
\bigskip\goodbreak
\lem

{\sl Soit
$$1\to H\buildrel f\over \lra G \buildrel g\over \lra L\to 1$$
une suite exacte de groupes topologiques. On suppose $f,g$ continue  et $g$ ouverte. Supposons $H$ et $L$ compacts et $G$ sŽparŽ. Alors $G$ est aussi compact.

}

{\bf Preuve}

Comme $H$ est compact et $G$ est sŽparŽ, $f$ est un homŽomorphisme sur $f(H)$. On peut donc identifier $H$ ˆ $f(H)$ et $f$ comme l'inclusion. Soit $(U_{i})_{i\in I}$ un recouvrement ouvert de $G$. Pour chaque $x\in G$, $Hx$ est compact. Donc, il existe $I_{x}\subset I$, fini, tel que $Hx\subset \cup_{i\in I_{x}} U_{i}$, ce qui veut dire que $H\subset\bigcup_{i\in I_{x}}U_{i} x^{-1}$. En vertu du lemme prŽcŽdent, il existe $U_{x}$ un voisinage de 1 tel que $U_{x}\cdot H\subset \bigcup_{i\in I_{x}}U_{i} x^{-1}$. Alors, $U_{x}Hx$ est un voisinage de $Hx$ contenu dans $\bigcup_{i\in I_{x}}U_{i}$. Or, puisque $H$ est un sous-groupe normal, on a $U_{x} Hx=U_{x}(xHx^{-1})x=U_{x} xH=\bigcup_{y\in U_{x}x}yH$ et comme $H$ est le noyau de $g$, on a $g^{-1}(g(U_{x}Hx))=U_{x}Hx$. Maintenant, puisque $g$ est ouverte et surjective, l'ensemble $\{ g(U_{x}Hx)\}_{x\in G}$ est un recouvrement ouvert de $L$. Puisque $L$ est compact, il existe $x_{1},\ldots ,x_{n}\in G$ tels que $L=\bigcup_{i=1}^n g(U_{x_{i}}Hx_{i})$. En prenant le $g^{-1}$ et en utilisant l'ŽgalitŽ vue avant, on a $G=\bigcup_{i=1}^n U_{x_{i}}Hx_{i}\subset \bigcup_{i=1}^n\bigcup_{j\in I_{x_{i}}}U_{i}\subset G$. Ce qui montre que $G$ est compact.\qed
\bigskip

\th

{\sl Le sous-groupe des classes d'idles spŽciaux $C^{0}_{K}$ (cf. DŽfinition \argfz) est compact.

}

{\bf Preuve}

ConsidŽrons $I_{K}$ vu comme groupe topologique (avec la topologie discrte) et l'homomorphisme~:
$$\eqalign{\psi\ :\ \bbI_{K}&\lra I_{K}\cr a=(a_{\P})_{\P}&\longmapsto\prod_{\P\in \gfP_{0}(K)}\P^{v_{\P}(a_{\P})}.\cr}$$
Il est bien dŽfini (les idles n'ont qu'un nombre finis de $a_{\P}$ de valuation non nulle) et surjectif. Son noyau est $\bbI_{K}(\infty)=\prod_{\P\in \gfP_{\infty}(K)}\bbK_{\P}^*\times \prod_{\P\in\gfP_{0}(K)}O_{\P}^*$ qui est un ouvert de $\bbI_{K}$, donc cet homomorphisme est aussi continu. De plus, l'image de $K^*$ est clairement $P_{K}$. D'o un isomorphisme topologique~:
$\bbI_{K}/(\bbI_K(\infty)\cdot K^*)\simeq I_K/P_K.$
Si on restreint l'homomorphisme $\bbI_K\to I_K$ ˆ $\bbI_K^0$, il est encore continu et surjectif (on choisit judicieusement les places finies (comme pour celui de dŽpart) puis on s'arrange avec les places infinies pour que le produit des normes donne 1). Le noyau est $\bbI_K^0(\infty):=\bbI_K^0\cap \bbI_K(\infty)$, et, comme avant,
$$\bbI_{K}^0/(\bbI^0_K(\infty)\cdot K^*)\simeq I_K/P_K,\eqno{(*)}$$
qui est Žvidemment ouverte. L'application $\bbI_K^0/K^*\to \bbI_K^0/(\bbI_K^0(\infty)\cdot K^*)$ est continue par dŽfinition et ouverte (cf. raisonnement dans la preuve la Proposition-DŽfinition \argfv) son noyau est Žvidemment $(\bbI_K^0(\infty)\cdot K^*)/K^*$ qui est ouvert dans $\bbI_K^0/K^*$ (car c'est l'image rŽciproque de $\{1\}$ qui est ouvert, car on vient de voir que $\bbI_K^0/(\bbI_K^0(\infty)\cdot K^*)$ Žtait muni de la topologie discrte). Ce qui nous donne une suite exacte~:
$$1\to (\bbI_K^0(\infty)\cdot K^*)/K^*\buildrel\rm continue\over\lra\bbI_K^0/K^*\buildrel\rm continue\ et\ ouvert\over\lra \bbI_K^0/(\bbI_K^0(\infty)\cdot K^*)\to 1.\eqno{(**)}$$
D'autre part, l'inclusion $\bbI_K^0(\infty)\hookrightarrow \bbI_K^0(\infty)\cdot K^*$ est continue par dŽfinition, mais elle est aussi ouverte, car $\bbI_K^0(\infty)\cdot K^*$ et $\bbI_K^0(\infty)$ sont des ouverts de $\bbI_K^0$. D'autre part, la projection $\bbI_K^0(\infty)\cdot K^*\to (\bbI_K^0(\infty)\cdot K^*)/K^*$ est aussi ouverte et continue (cf. preuve de la Proposition-DŽfinition \argfv). Ainsi, la composŽe de ces deux applications est aussi ouverte et continue. D'o, par passage au quotient, un isomorphisme continu et ouvert~:
$$\bbI_{K}^0(\infty)/(\bbI_{K}^0(\infty)\cap K^*)\simeq (\bbI_K^0(\infty)\cdot K^*)/K^*.\eqno{(***)}$$
On voit facilement que $\bbI_{K}^0(\infty)\cap K^*=U_{K}$ et rappelons que $C_{K}^0=\bbI_{K}^0/K^*$. En combinant les relations $(*)$, $(**)$ et $(***)$ on a la suite exacte~:
$$1\to \bbI_K^0(\infty)/U_{K}\buildrel\rm continue\over\lra C_{K}^0\buildrel\rm continue\ et\ ouvert\over\lra I_{K}/P_{K}\to 1.$$
En vertu du lemme prŽcŽdent, il suffit, pour achever la dŽmonstration, de montrer que 

\art{i)}$I_{K}/P_{K}$ est compact, mais cela nous le savons car c'est un groupe fini (cf. [Sam, Thm. 2, p.71]) muni de la topologie discrte, donc Žvidemment compact.

\art{ii)}$\bbI_K^0(\infty)/U_{K}$ est compact. Pour cela, nous allons travailler encore un petit peu et utiliser le mme lemme avec d'autres homomorphismes. RedŽfinissons une vielle connaissance vue au chapitre 1~:

\font\tet=cmr12 at 7pt
$$\eqalign{l\ :\ \bbI_{K}^0(\infty)&\lra \R^{r+s}\cr (a_{\P})_{\P}&\longmapsto (\underbrace{\log|a_{\P_{1}}|_{\P_{1}},\ldots, \log|a_{\P_{r}}|_{\P_{r}}}_{\rm places\ r\hbox{\tet Ž}elles},\underbrace{\log|a_{\P_{r+1}}|_{\P_{r+1}},\ldots, \log|a_{\P_{r+s}}|_{\P_{r+s}}}_{\rm places\ complexes}),\cr}$$
o $r$ et $s$ sont comme toujours le nombre de plongements rŽels respectivement complexes de $K$. L'image de $l$ est $H:=\{(x_{i})\in \R^{r+s}\mid\sum_{i=1}^{r+s} x_{i}=0\}$ et le noyau de $l$ est $\prod_{\P\in\gfP(K)}O_{\P}^*$ qui est compact (en vertu du Lemme de Tychonov). De plus, l'application $l$ est continue et ouverte, car les applications $a+bi\mapsto \sqrt{a^2+b^2}$, $a\mapsto \log(a)$ sont continues et ouvertes (lˆ o elles sont dŽfinies). En passant aux quotients, on obtient une suite exacte~:
$$1\to \left (\prod_{\P\in\gfP(K)}O_{\P}^*\right )\big /\left (U_{K}\cap \prod_{\P\in\gfP(K)}O_{\P}^*\right )\lra \bbI_{K}^0(\infty)/U_{K}\buildrel \overline{l}\over\lra H/l(U_{K})\to 1.$$
comme avant, l'homomorphisme injectif est continu $\overline{l}$ est continu et ouvert. Il s'agit donc de montrer que les groupes ``extŽrieurs'' sont compacts. Celui de gauche l'est facilement~: $U_{K}$ est inclu dans $K^*$ qui est discret dans $\bbI_{K}$ (cf. Lemme \argfu), donc $U_{K}$ est en particulier fermŽ et alors $\left (\prod_{\P\in\gfP(K)}O_{\P}^*\right )\big /\left (U_{K}\cap \prod_{\P\in\gfP(K)}O_{\P}^*\right )$ est compact, puisque c'est un espace compact quotientŽ par un fermŽ. Enfin, le thŽorme de Dirichlet sur les unitŽs nous dit que $l(U_{K})$ est un $\Z$-rŽseau de rang $n+r-1={\rm dim}_{\R}(H)$. Cela implique que $H/l(U_{K})$ est compact. Cela prouve le thŽorme.\qed
\bigskip

\defi
\medskip
Soit $K$ le corps de nombres que nous tra"nons depuis le dŽbut de ce chapitre et $\m=\prod_{\P\in\gfP(K)}\P^{m_{\P}}=\m_{0}\cdot\m_{\infty}$ un $K$-module. Soit $b\in\N$ et $\P\in\gfP(K)$. Souvenons-nous des $U_{\P}^{(b)}$ vus lors de la DŽfinition \aargek. On notera 
$$\bbI_{\m}=\prod_{\P\in\gfP(K)}U_{\P}^{(m_{\P})}\subset \bbI_{K}$$
 \newcount\gaagaw\gaagaw=\gaga
Puisque dans les corps non-archimŽdiens toute boule fermŽe est ouverte, $\bbI_{\m}$ est un sous-groupe ouvert de $\bbI_{K}$. Il est clair que $\bigcap_{\m}\bbI_{\m}=\prod_{\P\in\gfP_{\R}(K)}\R_{+}^*\times\prod_{\P\in\gfP_{\C}(K)}\C^*\times \prod_{\P\in\gfP_{0}(K)}\{ 1\}$. Remarquons que l'ensemble $\{\bbI_{\m}\mid \m$ est un $K$-module $\}$ est un systme fondamental de sous-groupes ouverts~: si $U$ est un sous-groupe ouvert de $\bbI_{K}$, alors il existe $S$ un ensemble fini de places contenant les places infinies tel que $U\supset \prod_{\P\not\in S}U_{\P}\times \prod_{\P\in S}V_{\P}$, o $V_{\P}$ est un voisinage ouvert de $1$ dans $\bbK_{\P}^*$. Soit $\P\in S$. Si $\P$ est fini, on choisit $m_{\P}\in\N$ assez grand pour que $U_{\P}^{(m_{\P})}\subset V_{\P}$. Si $\P$ est infini complexe, alors on peut prendre $V_{\P}=\C^*$, car tout sous-groupe ouvert est fermŽ (vŽrification facile) et donc, puisque $\C^*$ est connexe, c'est forcŽment $\C^*$ lui-mme; et pour la mme raison, si $\P$ est infini rŽel, on peut prendre $V_{\P}=\R_{+}^*$, ainsi en prenant $\m$ le $K$-module contenant toutes les places infinies rŽelles et dont les $m_{\P}$ sont ceux donnŽs plus haut pour les places finies est bien tel que $\bbI_{\m}\subset U$. Enfin, on voit facilement que si $\m_{1}$ et $\m_{2}$ sont des $K$-modules, $\bbI_{{\rm pgcd}(\m_{1},\m_{2})}=\bbI_{\m_{1}}\cdot \bbI_{\m_{2}}$ et que si $\m_{1}|\m_{2}$, alors $\bbI_{\m_{2}}\subset\bbI_{\m_{1}}$.

On dŽfinit aussi 
$$\bbI'_{\m}=\{(a_{\P})_{\P}\in\bbI_{K}\mid a_{\P}\in U_{\P}^{(m_{\P})}\ \forall \P |\m\}.$$
Il est clair que $\bbI_{\m}\subset \bbI_{\m}'$ et que $\bbI_{\m}'$ est un ouvert de $\bbI_{K}$. 

Enfin, on pose 
$$C_{\m}=(\bbI_{\m}\cdot K^*)/K^*.$$
\bigskip\goodbreak
\th
\medskip
{\sl 
Soit $\m$ un $K$-module. Alors on a les isomorphismes de groupes topologiques~:
$$\bbI_{K}/(\bbI_{\m}\cdot K^*)\simeq \bbI'_{\m}/(\bbI_{\m}\cdot K^*_{\m})\simeq I_{K}(\m)/P_{\m},$$
o $K^*_{\m}, I_{K}(\m)$ et $P_{\m}$ sont les groupes connus de longue date, dŽfinis au Chapitre 0, $K^*_{\m}$ Žtant bien entendu associŽ ˆ l'idle principal correspondant. Le premier de ces isomorphismes est donnŽ par l'inclusion $\bbI'_{\m}\subset \bbI_{K}$ et le second est donnŽ par l'application $\psi\ :\ \bbI_{K}\to I_{K}$ vue ˆ la preuve du ThŽorme \arggd\ restreinte ˆ $\bbI'_{\m}$ qu'on notera dŽsormais $\psi_{\m}$.

}

{\bf Preuve}

L'application 
$$\eqalign{\psi_{\m}\ :\ \bbI'_{\m}&\lra I_{K}(\m)\cr (a_{\P})_{\P}&\longmapsto \prod_{\P\in\gfP_{0}(K)}\P^{v_{\P}(a_{\P})}\cr}$$
est un homomorphisme surjectif. Son noyau est $\bbI_{\m}$ qui est ouvert, donc cet homomorphisme est continu. Il est clair que $P_{\m}$ est l'image de $K^*_{\m}$ vu comme idle principal. Donc, en composant avec la projection $I_{K}(\m)\to I_{K}(\m)/P_{\m}$, on obtient un homomorphisme surjectif $\bbI'_{\m}\to I_{K}(\m)/P_{\m}$ dont le noyau est $\bbI_{\m}\cdot K^*_{\m}$ qui est ouvert. On obtient donc le second isomorphisme
$$\bbI'_{\m}/(\bbI_{\m}\cdot K^*_{\m})\simeq I_{K}(\m)/P_{\m}.$$
Or, c'est une vŽrification de voir que $K^*_{\m}=\bbI'_{\m}\cap K^*$. Donc $\bbI_{\m}\cdot K^*_{\m}=\bbI_{\m}\cdot (\bbI'_{\m}\cap K^*)=\bbI'_{\m}\cap (\bbI_{\m}\cdot K^*)$. D'o,
$$\bbI'_{\m}/(\bbI_{m}\cdot K^*_{\m})=\bbI'_{\m}/(\bbI'_{\m}\cap (\bbI_{\m}\cdot K^*))\buildrel (*)\over\simeq \bbI'_{\m}\cdot(\bbI_{\m}\cdot K^*) /(\bbI_{\m}\cdot K^*)= (\bbI'_{\m}\cdot K^*)/(\bbI_{\m}\cdot K^*),$$
l'isomorphisme $(*)$ venant du troisime thŽorme d'isomorphisme en observant de plus qu'il est continu et ouvert, car $\bbI'_{\m}$ est un ouvert dans $\bbI_{K}$, donc dans $\bbI'_{\m}\cdot(\bbI_{\m}\cdot K^*)$. Enfin, on montre que $\bbI'_{\m}\cdot K^*=\bbI_{K}$. En effet, soit $(a_{\P})_{\P}\in \bbI_{K}$. Par densitŽ et gr‰ce au thŽorme d'approximation dŽbile, il existe $\alpha\in K^*$ tel que ${a_\P \over \alpha}\equiv 1\pmodast {\widehat{\P}^{m_\P}}$ pour tout $\P|\m$, ce qui veut dire que ${a_\P \over \alpha}\in U_\P^{(m_\P)}$ pour tout $\P |\m$ et donc que $({a_\P\over \alpha})_\P\in \bbI'_\m$. Et cela prouve le premier isomorphisme~:
$$\bbI'_{\m}/(\bbI_{\m}\cdot K^*_{\m})\simeq\bbI_{K}/(\bbI_{\m}\cdot K^*).$$
\qed
\bigskip
\coro
\medskip
{\sl Tout sous-groupe ouvert de $\bbI_K$ contenant $K^*$ doit contenir un $\bbI_\m$ et est nŽcessairement d'indice fini. 
De manire similaire, les sous-groupes ouverts de $C_K$ sont ceux qui contiennent un sous-groupe $C_\m$. Il sont tous d'indice fini et pour tout $K$-module $\m$, on a un isomorphisme
$$C_K/C_\m\simeq I_K(\m)/P_\m.$$
Et rŽciproquement, si $H$ est un sous-groupe de $C_K$ tel que $H\supset C_\m$, alors il est forcŽment ouvert.

}
{\bf Preuve}

Si $H$ est un sous-groupe ouvert de $\bbI_K$, on a vu ˆ la DŽfinition \argge\ qu'il existe $\m$ un $K$-module tel que $\bbI_\m\subset H$. Donc, par hypothse, on a $\bbI_K\supset H\supset \bbI_\m\cdot K^*$. En utilisant la finitude de $I_K(\m)/P_\m$ (cf. ThŽorme \argl) et le ThŽorme \arggf, on conclut que $H$ est d'indice fini dans $\bbI_K$. La seconde assertion est Žvidente au vue du ThŽorme \arggf\ et du deuxime thŽorme d'isomorphisme (qui prŽserve la continuitŽ). La dernire assertion est aussi Žvidente, car si $H\supset C_\m$, il est isomorphe ˆ un sous-groupe de $I_K(\m)/P_\m$ qui est fini avec la topologie discrte, donc forcŽment ouvert.\qed
\bigskip\goodbreak

\defi
\medskip
Si $H$ est un sous-groupe ouvert de $\bbI_{K}$ contenant $K^*$, on dit que le $K$-module $\m$ est {\it admissible pour $H$}, si $\bbI_{\m}\subset H$. On vu ˆ la DŽfinition \argge\ qu'un tel $\m$ existait toujours. On vŽrifie facilement que si $\m$ et $\n$ sont admissibles pour $H$, alors  $\rm pgcd(\m,\n)$ est aussi admissible pour $H$ (car on a vu ˆ la DŽfinition \argge\ que $\bbI_{\rm pgcd(\m,\n)} =\bbI_{\m}\cdot\bbI_{\n}$). Il existe donc un $K$-module $\ff$ (appelŽ le {\it conducteur de $H$}) tel que 

$$\m \hbox{ est admissible pour $H\iff\m$ divise $\ff$.}$$
On rappelle les applications $$\eqalign{\psi\ :\ \bbI_{K}&\lra I_{K}\cr a=(a_{\P})_{\P}&\longmapsto\prod_{\P\in \gfP_{0}(K)}\P^{v_{\P}(a_{\P})}\cr}\quad\ \hbox{et}\quad \eqalign{\psi_{\m}=\psi |_{\bbI_{m}'}\ :\ \bbI'_{\m}&\lra I_{K}(\m)\cr (a_{\P})_{\P}&\longmapsto \prod_{\P\in\gfP_{0}(K)}\P^{v_{\P}(a_{\P})}.\cr}$$
Soit $H$ comme ci-dessus et $\m$ un $K$-module admissible pour $H$. On pose $$H(\m)=\psi(H\cap \bbI_{\m}').$$
Cette notation est la mme que celle dŽfinissant le sous-groupe de congruence dŽfini pour $\m$. Ce n'est pas un hasard~:

\newcount\gaagax\gaagax=\gaga

\bigskip
\th
\medskip

{\sl Soit $K$ un corps de nombres. Alors l'application 
$$H\longmapsto \{H(\m)\mid \m\hbox{ est admissible pour } H\}$$
est une bijection de l'ensemble des sous-groupes ouverts de $\bbI_{K}$ contenant $K^*$ sur l'ensemble des classes d'Žquivalence de sous-groupes de congruences (voir Chapitre 8 pour les dŽfinitions). En outre, le conducteur de $H$ est Žgal au conducteur (au sens du Corollaire-DŽfinitions \argdi) de la classe de sous-groupes de congruence correspondante. Enfin, on a encore l'ŽgalitŽ $\bbI_{K}/H\simeq I_{K}(\m)/H(\m)$.

}

{\bf Preuve}
 
Soit $H$ un sous-groupe ouvert de $\bbI_{K}$ contenant $K^*$ et $\m$ un $K$-module admissible pour $H$. Il y a bijection entre les sous-groupes ouverts de $\bbI_{K}$ qui contiennent $\bbI_{\m}\cdot K^*$ et les sous-groupes ouverts de $\bbI_{K}/(\bbI_{\m}\cdot K^*)$ qui sont en bijection, par le ThŽorme \arggf , avec les sous-groupe ouverts de $I_{K}(\m)/P_{\m}$, qui correspondent eux-mme aux ouverts (pour la topologie discrte) de $I_{K}(\m)$ qui contiennent $P_{\m}$. Ces bijections Žtant donnŽes par $\psi$ et l'inclusion, on a donc que $P_{\m}\subset H(\m)\subset I_{K}(\m)$. Plus prŽcisŽment, la rŽciproque de l'application $H\mapsto H(\m)$ est l'application $H'\mapsto \psi_{\m}^{-1}(H')\cdot K^*$. En effet, $\psi_{\m}^{-1}(H(\m))\cdot K^*=\psi_{\m}^{-1}(\psi_{\m}(H\cap \bbI_{\m}'))\cdot K^*=(H\cap \bbI_{m}')\cdot \bbI_{\m}\cdot K^*=H$. La dernire ŽgalitŽ se montre en utilisant que $\bbI_{K}=\bbI_{\m}'\cdot K^*$ et que $H\supset \bbI_{\m}\cdot K^*$. RŽciproquement, on voit que $\psi_{\m}((\psi_{\m}^{-1}(H')\cdot K^*)\cap \bbI_{\m}')=H'$, car $K^*\cap \bbI_{\m}'=K_{\m}^*$ et que $\psi_{\m}(K_{\m}^*)=P_{\m}\subset H'$. Cela montre que $H(\m)$ est  un sous-groupe de congruence pour $\m$ (au sens de la DŽfinition \argde). On a donc un diagramme

\vglue 4cm 

\psset{xunit=0.8cm,yunit=0.8cm}

\rput(5,0){\rnode{a1}{$\bbI_{\m}\cdot K^*$}}
\rput(8,0){\rnode{a2}{$\bbI_{\m}\cdot K^*_{\m}$}}
\rput(11,0){\rnode{a3}{$P_{\m}$}}
\rput(5,2){\rnode{b1}{$H$}}
\rput(8,2){\rnode{b2}{$H\cap\bbI_{\m}'$}}
\rput(11,2){\rnode{b3}{$H(\m)$}}
\rput(5,4){\rnode{c1}{$\bbI_{K}$}}
\rput(8,4){\rnode{c2}{$\bbI_{\m}'$}}
\rput(11,4){\rnode{c3}{$I_{K}(\m)$}}
\ncline[nodesep=3pt]{->}{a2}{a1}
\ncline[nodesep=5pt]{->}{a1}{b1}
\rput{90}(5.013,0.35){\rnode{u1}{$\taille{20}\lhook$}}
\ncline[nodesep=5pt]{->}{a2}{b2}
\rput{90}(8.013,0.36){\rnode{u1}{$\taille{20}\lhook$}}
\ncline[nodesep=5pt]{->}{a3}{b3}
\rput{90}(11.013,0.36){\rnode{u1}{$\taille{20}\lhook$}}
\ncline[nodesep=3pt]{->}{b2}{b1}
\ncline[nodesep=7pt]{->}{c2}{c1}
\rput(7.5,3.99){\rnode{u1}{$\taille{20}\rhook$}}
\ncline[nodesep=3pt]{->}{c2}{c3}
\Aput{$\psi_{\m}$}
\ncline[nodesep=5pt]{->}{b1}{c1}
\rput{90}(5.013,2.36){\rnode{u1}{$\taille{20}\lhook$}}
\ncline[nodesep=5pt]{->}{b2}{c2}
\rput{90}(8.013,2.38){\rnode{u1}{$\taille{20}\lhook$}}
\ncline[nodesep=5pt]{->}{b3}{c3}
\rput{90}(11.013,2.38){\rnode{u1}{$\taille{20}\lhook$}}
\ncline[nodesep=5pt]{->}{a2}{a3}
\ncline[nodesep=5pt]{->}{b2}{b3}
\medskip

En observant ce diagramme, on remarque aisŽment (gr‰ce au premier thŽorme d'isomorphisme) que $\bbI_{K}/H\simeq I_{K}(\m)/H(\m)$ et que tout sous-groupe de congruence (pour $K$) est de la forme $H(\m)$ pour un $H$ et un $\m$ adŽquat. De plus, pour les mme raisons, si $H_{1}$ et $H_{2}$ sont des sous-groupes ouverts de $\bbI_{K}$ contenant $\bbI_{\m}\cdot K^*$, alors $H_{1}(\m)=H_{2}(\m)$ implique que $H_{1}=H_{2}$ $(*)$.

Supposons que $\m$ et $\m'$ soient admissibles pour $H$. Alors $H(\m)$ et $H(\m')$ sont Žquivalents (au sens de la DŽfinition \argdf). En effet, on peut supposer sans limiter la gŽnŽralitŽ en passant par le pgcd, que $\m | \m'$ (dans ce cas, il est clair que $\bbI_{\m'}'\subset \bbI_{\m}'$). Alors on a~:
$$H(\m')=\psi( \bbI_{\m'}'\cap H)=\psi((\bbI_{\m}'\cap H)\cap \bbI_{\m'}')\buildrel (**)\over =\psi(\bbI_{\m}'\cap H)\cap\psi ( \bbI_{\m'}')=H(\m)\cap I_{K}(\m'),\eqno{(***)}$$
ce qui montre que $H(\m)$ et $H(\m')$ sont Žquivalents. L'ŽgalitŽ $(**)$ vient du rŽsultat suivant~: si $A,B\subset G$, $G$ est un groupe, $A$ un sous-groupe de $G$, $f$ un homomorphisme dŽfini sur $G$ tel que $\ker (f | G)\subset A$, alors  $f(A\cap B)=f(A)\cap f(B)$. L'inclusion $\subset$ est toujours vraie et triviale. RŽciproquement, soit $z\in f(A)\cap f(B)$. Alors il existe $x\in A$ et $y\in B$ tels que $f(x)=f(y)=z$, ce qui veut dire que $x^{-1}\cdot y\in \ker(f)\subset A$, et donc $y=x\cdot (x^{-1} y)\in A$.

Enfin, pour achever la preuve du thŽorme, il faut encore voir la chose suivante~: soit $H_{1},H_{2}$ des sous-groupes ouverts de $\bbI_{K}$ contenant $K^*$, et $\m_{1},\m_{2}$, des $K$-modules admissibles pour $H_{1}$ et $H_{2}$ respectivement, alors on a 
$$H_{1}=H_{2}\iff H_{1}(\m_{1})\sim H_{2}(\m_{2}),$$
o $\sim$ est l'Žquivalence des sous-groupes de congruences. En effet, soit $\m$ un multiple commun de $\m_{1}$ et de $\m_{2}$. La relation $(***)$ montre que $H_{1}(\m_{1})\sim H_{1}(\m)$ et $H_{2}(\m_{2})\sim H_{2}(\m)$. Ainsi 
$$H_{1}(\m_{1})\sim H_{2}(\m_{2})\iff H_{1}(\m)\sim H_{2}(\m)\buildrel \rm cor-def\argdi\over \sim H_{1}(\m)=H_{2}(\m)\buildrel (*)\over \iff H_{1}=H_{2}.$$
\qed
\bigskip
Nous pouvons maintenant Žnoncer la premire version de la version idŽlique du corps de classe~:
\bigskip\goodbreak
\coro
\medskip
{\sl Soit $K$ un corps de nombres. Alors il existe une bijection entre l'ensemble des sous-groupes ouverts $H$ de $\bbI_{K}$ contenant $K^*$ et l'ensemble des extensions abŽliennes $L$ de $K$ (contenus dans une mme cl™ture algŽbrique). On a en outre une bijection entre le groupe de Galois  ${\rm Gal} ( L/K)$ et $\bbI_{K}/H$.

}
{\bf Preuve}

C'est un corollaire immŽdiat du thŽorme prŽcŽdent et du ThŽorme \argef.\qed
\bigskip
Nous allons maintenant prŽciser encore un peu de quelle nature est cette bijection.
\bigskip
\defi
\medskip
Soit $L/K$ une extension de corps de nombres. Si $\P\in \gfP(K)$, on note $\bbL_{\P}=\prod_{\gP|\P}\bbL_{\gP}$. Alors on peut voir $\bbA_{L}$ comme le produit rŽduit des $\bbL_{\P}$ par rapport aux $\prod_{\gP|\P}O_{\P}$ ($\P$ fini). De mme, $\bbI_{L}$ peut tre vu comme le produit des $\bbL_{\P}^*=\prod_{\gP |\P}\bbL_{\gP}^*$ par rapport aux $\prod_{\gP |\P}U_{\gP}$ ($\P$ finis). En fait, on regroupe par ``paquets''. On dŽfinit entre $\bbI_{L}$ et $\bbI_{K}$ une norme qu'on note encore $N_{L/K}$~:
$$\eqalign{ N_{L/K}\ :\ \bbI_{L}&\lra \bbI_{K}\cr x=(x_{\gP})_{\gP\in\gfP(L)}&\longmapsto N_{L/K}(x)=y=(y_{\P})_{\P\in\gfP(K)}\cr}$$
tel que pour tout $\P\in \gfP(K)$ on ait $y_{\P}=\prod_{\gP|\P}N_{\bbL_{\gP}/\bbK_{\P}}(x_{\gP})$.
\newcount\gaagay\gaagay=\gaga
\goodbreak

Si $x\in L^*$, il est bien connu que pour tout $\gP\in \gfP(K)$, on a $N_{L/K}(x)=\prod_{\gP|\P}N_{\bbL_{\gP}/\bbK_{\P}}(x)$ (cf. [Fr-Tay, Ch. III, 1.10,p.110]). Cela montre que le premier carrŽ est commutatif (l'autre l'est plus trivialement)~:
 \vglue 2cm
 
\psset{xunit=0.8cm,yunit=0.8cm}

\rput(4,0){\rnode{a1}{$ K^*$}}
\rput(7,0){\rnode{a2}{$\bbI_{K}$}}
\rput(11,0){\rnode{a3}{$K^*$}}
\rput(14,0){\rnode{a4}{$\bbI_{K}$}}
\rput(4,2){\rnode{b1}{$L^*$}}
\rput(7,2){\rnode{b2}{$\bbI_{L}$}}
\rput(11,2){\rnode{b3}{$L^*$}}
\rput(14,2){\rnode{b4}{$\bbI_{L}$}}
\ncline[nodesep=3pt]{->}{b1}{a1}
\Bput{$N_{L/K}$}
\ncline[nodesep=3pt]{->}{b2}{a2}
\Aput{$N_{L/K}$}
\ncline[nodesep=6pt]{->}{b1}{b2}
\rput(4.5,1.99){\rnode{u1}{$\taille{20}\lhook$}}
\ncline[nodesep=6pt]{->}{a1}{a2}
\rput(4.5,-0.01){\rnode{u1}{$\taille{20}\lhook$}}
\ncline[nodesep=6pt]{->}{b3}{b4}
\rput(11.5,1.99){\rnode{u1}{$\taille{20}\lhook$}}
\ncline[nodesep=6pt]{->}{a3}{a4}
\rput(11.5,-0.01){\rnode{u1}{$\taille{20}\lhook$}}
\ncline[nodesep=5pt]{->}{a3}{b3}
\Aput{$\rm incl.$}
\rput{90}(11.013,0.36){\rnode{u1}{$\taille{20}\lhook$}}
\ncline[nodesep=5pt]{->}{a4}{b4}
\Bput{$\rm incl.$}
\rput{90}(14.013,0.36){\rnode{u1}{$\taille{20}\lhook$}}

Enfin, puisque $N_{L/K}(L^*)\subset K^*$, $N_{L/K}$ induit un homomorphisme de $C_{L}\to C_{K}$ qu'on note encore $N_{L/K}$.
\bigskip
\lem
\medskip
{\sl Soit $L/K$ une extension de corps de nombres. Alors le diagramme suivant est commutatif

 \vglue 2cm
 
\psset{xunit=0.8cm,yunit=0.8cm}

\rput(7,0){\rnode{a1}{$ \bbI_{K}$}}
\rput(9,0){\rnode{a2}{$I_{K}$}}
\rput(7,2){\rnode{b1}{$ \bbI_{L}$}}
\rput(9,2){\rnode{b2}{$I_{L}$}}
\ncline[nodesep=3pt]{->}{b1}{a1}
\Bput{$N_{L/K}$}
\ncline[nodesep=3pt]{->}{b2}{a2}
\Aput{$N_{L/K}$}
\ncline[nodesep=3pt]{->}{b1}{b2}
\Aput{$\psi $}
\ncline[nodesep=3pt]{->}{a1}{a2}
\Aput{$\psi $}

Il commute ainsi~:

\vglue 3cm

\rput(5,0){\rnode{a1}{$\left (\prod_{\gP|\P}N_{\bbL_{\gP}/\bbK_{\P}}(x_{\gP})\right )_{\P}$}}
\rput(11,0){\rnode{a2}{$\dst\prod_{\P\in\gfP_{0}(K)}\P^{\sum_{\gP|\P}v_{\P}(N_{\bbL_{\P}/\bbK_{\P}}(x_{\gP}))}$}}
\rput(5,3){\rnode{b1}{$(x_{\gP})_{\gP}$}}
\rput(11,3){\rnode{b2}{$\dst\prod_{\gP\in\gfP_{0}(L)}\gP^{v_{\gP}(x_{\gP})}$}}
\ncline[nodesep=3pt]{|->}{b1}{a1}
\Bput{$N_{L/K}$}
\ncline[linestyle=dashed,nodesep=3pt]{|->}{b2}{a2}
\Aput{$N_{L/K}$}
\ncline[nodesep=3pt]{|->}{b1}{b2}
\Aput{$\psi $}
\ncline[nodesep=3pt]{|->}{a1}{a2}
\Aput{$\psi $}
\medskip
Et la flche en traitillŽ est bien dŽfinie

}
{\bf Preuve}

A priori la flche en traitillŽ est $\prod_{\gP\in\gfP_{0}(L)}\gP^{v_{\gP}(x_{\gP})}\longmapsto \prod_{\P\in\gfP_{0}(K)}\P^{\sum_{\gP|\P}f(\gP/\P) v_{\gP}(x_{\gP})}$.

Donc, la seule chose qu'il faut voir est que pour tout $\gP|\P$, on a $$v_{\P}(N_{\bbL_{\gP}/\bbK_{\P}}(x))=f(\gP/\P)\cdot v_{\gP}(x)$$ pour tout $x\in L_{\gP}^*$. On rappelle que $x=\pi^k\cdot u$, o $\pi$ est une uniformisante, $k\in \Z$ et $u\in U_{\gP}$. Par multiplicativitŽ de la norme, il suffit de voir le rŽsultat pour $u$ et $\pi^k$ sŽparŽment. Pour $u$, il  est clair que $N_{\bbL_{\gP}/\bbK_{\P}}(U_{\gP})\subset U_{\P}$ et que $v_{\gP}(u)=0$, donc les deux membres de l'ŽgalitŽ valent 0. Enfin, $v_{\P}(N_{\bbL_{\gP}/\bbK_{\P}}(\pi^k))=v_{\P}(N_{\bbL_{\gP}/\bbK_{\P}}(\pi^k)\cdot O_{\P})=v_{\P}(N_{\bbL_{\gP}/\bbK_{\P}}(\pi\cdot O_{\gP})^k)=v_{\P}(\P^{f(\gP/\P)\cdot k})=f(\gP/\P)\cdot k=f(\gP/\P)\cdot v_{\gP}(\pi^k)$.\qed
\bigskip

\prop
\medskip
{\sl Si $L/K$ est une extension abŽlienne de corps de nombres, il existe un $K$-module $\m$ tel que $$N_{L/K}(\bbI_{L})\supset \bbI_{\m}.$$ Il est Žvident que tout multiple de $\m$ fait aussi l'affaire. Cela implique que $N_{L/K}(\bbI_{L})$ est un ouvert de $\bbI_{K}$ (car alors $N_{L/K}(\bbI_{L})=\bigcup_{x\in N_{L/K}(\bbI_{L})}x\cdot \bbI_{\m}$).

}

{\bf Preuve}

Il suffit de trouver un $K$-module $\m$ tel que pour chaque $\P\in\gfP(K)$, il existe $\gP|\P$ tel que $N_{\bbL_{\gP}/\bbK_{\P}}(\bbL_{\gP})\supset U_{\P}^{(v_{\P}(\m))}$ (pour les autres $\gP|\P$, on pose 1) et $N_{L_{\gP}/K_{\P}}(U_{\gP})\supset U_{\P}^{(v_{\P}(\m))}$ ppt les $\P$ et les $\gP|\P$. On sait que si $\P$ est non ramifiŽ, alors $N_{\bbL_{\gP}/\bbK_{\P}}(U_{\gP})=U_{\P}$ (Lemme \argeo). Pour les places ramifiŽes (qui sont en nombre fini, notons $R$ l'ensemble de ces places) la Proposition \argcc\ montre que si $m$ est assez grand, $N_{\bbL_{\gP}/\bbK_{\P}}(\bbL_{\gP})\supset U_{\P}^{(m)}$ (pour les places infinies, il suffit de prendre $m=1$). On choisit donc pour $\m$ un $K$-module $\prod_{\P\in R}\P^{m_{\P}}$ avec $m_{\P}>0$ et suffisamment grand si $\P$ est fini.\qed
\bigskip
\defi
\medskip
Soit $L/K$ une extension de corps de nombres et $\m$ un $K$-module. On note
$$\bbI''_{L}(\m)=\{(a_{\gP})_{\gP}\in \bbI_{L}\mid a_{\gP}=1\ \forall \gP\ \hbox{ tel que } \gP |\widetilde{\m}\}$$

\newcount\gaagaz\gaagaz=\gaga
Rappelons que $\widetilde{\m}$ a ŽtŽ dŽfini ˆ la page \the\gaagi.
\bigskip
\lem
\medskip
{\sl 
Soit $L/K$ une extension de corps de nombres et $\m$ un $K$-module. Alors on a les ŽgalitŽs~:
$$\displaylines{\psi^{-1}(P_{\m}\cdot N_{L/K}(I_{L}(\widetilde{\m})))=K^*_{\m}\cdot \bbI_{\m}\cdot N_{L/K}(\bbI''_{L}(\m))\cr
K^*\cdot\bbI_{\m}\cdot N_{L/K}(\bbI''_{L}(\m))= K^*\cdot\bbI_{\m}\cdot N_{L/K}(\bbI_{L}).\cr}$$

}
{\bf Preuve}

Pour la premire ŽgalitŽ.

Montrons ``$\supset$Ó. On sait que $\psi(K^*_{\m})=P_{\m}$ et $\psi(\bbI_{\m})=1$. D'autre part, $\psi (N_{L/K}(\bbI''_{L}(\m)))\buildrel\rm Lemme\ \arggl\over=N_{L/K}(\psi(\bbI''_{L}(\m)))\subset N_{L/K}(I_{L}(\widetilde{\m}))$ (la dernire inclusion est Žvidente). 

Montrons ``$\subset$''.  Soit $a\in \psi^{-1}(P_{\m}\cdot N_{L/K}(I_{L}(\widetilde{\m})))$. Alors $\psi (a)=(\alpha)\cdot N_{L/K}(\aa)$, avec $\alpha\in K^*_{\m}$ et $\aa\in I_{L}(\widetilde{\m})$. Soit $A=(A_{\gP})_{\gP}\in \bbI''_{L}(\m)$ tel que $A_{\gP}=1$ si $\gP$ est infini ou si $\gP\notdiv\aa$. Si $\gP |\aa$, on choisit $A_{\gP}$ tel que $v_{\gP}(A_{\gP})=v_{\gP}(\aa)$. Ainsi,  $\psi(A)=\aa$ et $\psi(N_{L/K}(A))=N_{L/K}(\psi(A))=N_{L/K}(\aa)$. Donc, $\psi(\alpha\cdot N_{L/K}(A))=(\alpha)\cdot N_{L/K}(\aa)=\psi(a)$. Donc, $a^{-1}\cdot\alpha\cdot N_{L/K}(A)\in\ker(\psi_{\m})=\bbI_{\m}$.

Pour la seconde ŽgalitŽ, ``$\subset$'' est triviale. Pour ``$\supset$'', il suffit de montrer que $N_{L/K}(\bbI_{L})\subset K^*\cdot\bbI_{\m}\cdot N_{L/K}(\bbI''_{L}(\m))$. Soit $b=(b_{\gP})_{\gP}\in \bbI_{L}$. Pour tout $x\in L^*$ (plongŽ diagonalement dans $\bbI_{L}$), on a $N_{L/K}(b)=N_{L/K}(x)\cdot N_{L/K}(b\cdot x^{-1})$. Ecrivons $b\cdot x^{-1}=b'\cdot b''$, o $b'=(b'_{\gP})_{\gP}\in \bbI''_{L}(\m)$ et $b''=(b''_{\gP})_{\gP}$, avec $b'_{\gP}=\cases{b_{\gP}\cdot x^{-1}&si $\gP\notdiv \widetilde{\m}$\cr 1&sinon\cr}$ et $b''_{\gP}=\cases{1&si $\gP\notdiv \widetilde{\m}$\cr b_{\gP}\cdot x^{-1}&sinon\cr}$. Et donc 
$$N_{L/K}(b)=\underbrace{N_{L/K}(x)}_{\in K^*}\cdot \underbrace{N_{L/K}(b')}_{\in N_{L/K}(\bbI''_{L}(\m))}\cdot N_{L/K}(b'').$$
Maintenant, il s'agit de choisir $x$ de sorte que $N_{L/K}(b'')\in \bbI_{\m}$. Or, puisque $N_{L/K}(b'')_{\P}=1$ si $\P\notdiv\m$, il suffit de de demander que $N_{L/K}(b'')_{\P}\in U_{\P}^{(v_{\P}(\m))}$ si $\P|\m$, c'est-ˆ-dire  $N_{L/K}(b''_{\gP})_{\P}\in U_{\P}^{(v_{\P}(\m))}$ si $\gP|\P$ et $\P|\m$. Ceci est rŽalisŽ si $x\equiv b_{\gP}\pmodast {\widehat{\gP}^m}$ pour $m$ assez grand, si $\gP|\P$ et $\P|\m$. Mais cela est vrai par densitŽ et en vertu du thŽorme d'approximation dŽbile (ThŽorme \argc).\qed
\bigskip
On peut maintenant Žnoncer le thŽorme principal du corps de classe un peu plus affinŽe que le Corollaire \arggj.
\bigskip\goodbreak
\th
\medskip
{\sl Soit $K$ un corps de nombres. On a une correspondance bijective entre les extensions abŽliennes de $K$ et les sous-groupes ouverts de $\bbI_{K}$ contenant $K^*$ . Cette correspondance est donnŽe par 
$$L/K\longmapsto K^*\cdot N_{L/K}(\bbI_{L}).$$
En passant aux classes, on peut dire la mme chose en disant que les extensions abŽliennes de $K$ sont en correspondances bijective avec les sous-groupes ouverts de $C_{K}$ via $L/K\mapsto N_{L/K}(C_{L})$. En outre on a des isomorphismes 
$$C_{K}/N_{L/K}(C_{L})\simeq \bbI_{K}/(K^*\cdot N_{L/K}(\bbI_{L}))\simeq{\rm Gal}(L/K),$$
obtenus en composant les isomorphismes
$$\bbI_{K}/(K^*\cdot N_{L/K}(\bbI_{L}))\simeq I_{K}(\m)/(P_{\m}\cdot N_{L/K}(I_{L}(\widetilde{\m})))\buildrel \rm Artin\over \simeq {\rm Gal}(L/K).$$

}

{\bf Preuve}

Si $L/K$ est une extension abŽlienne, on considre $\bbH=\bbH(L/K)$ la classe d'Žquivalence de cette extension. Choisissons $\m$, un $K$-module admissible. Le noyau de l'application d'Artin est le sous-groupe de congruence pour $\m$ et il est Žgal ˆ $P_{\m}\cdot N_{L/K}(I_{L}(\widetilde{\m}))$ (cf. Lemme \argco). On a vu lors de la preuve du ThŽorme \arggi\ que le sous-groupe ouvert de $\bbI_{K}$ associŽ ˆ ce sous-groupe de congruence Žtait $\psi^{-1}(P_{\m}\cdot N_{L/K}(I_{L}(\widetilde{\m})))\cdot K^*$. En vertu du lemme prŽcŽdent \arggo, on a alors~:
$$\eqalign{\psi^{-1}(P_{\m}\cdot N_{L/K}(I_{L}(\widetilde{\m})))\cdot K^*&=\underbrace{K^*\cdot K^*_{\m}}_{=K^*}\cdot \bbI_{\m}\cdot N_{L/K}(\bbI''_{L}(\m))\cr
&=K^*\cdot\bbI_{\m}\cdot N_{L/K}(\bbI_{L})\cr
&=K^*\cdot N_{L/K}(\bbI_{L}).\cr}$$
La dernire ŽgalitŽ vient du fait qu'en prenant un multiple adŽquat de $\m$, on peut supposer que $\bbI_{\m}\subset N_{L/K}(\bbI_{L})$ (cf. Proposition \arggm), et il reste Žvidemment admissible.

RŽciproquement, si $H$ est un sous-groupe ouvert de $\bbI_{K}$ contenant $K^*$, il lui correspond (en vertu du ThŽorme \arggi), une classe d'Žquivalence de sous-groupes de congruences, qui lui fait correspondre un corps de classe $L$ (ThŽorme \argef), et en ``revenant'' comme ˆ la premire partie, on voit que $H=K^*\cdot N_{L/K}(\bbI_{L})$. La dernire partie est Žvidente au vu de ce qui prŽcde.\qed
\bigskip
Pour terminer la partie idŽlique de cette thŽorie, nous allons chercher une description ``directe'' de l'homomorphisme $\Upsilon\ : \bbI_{K}\to {\rm Gal}(L/K)$. qui engendre l'isomorphisme $ \bbI_{K}/(K^*\cdot N_{L/K}(\bbI_{L}))\simeq{\rm Gal}(L/K)$.
\bigskip
\defi
\medskip
Soit $L/K$ une extension abŽlienne de corps de nombres et $\P$ une place non complexe de $K$. Souvenons-nous de l'application $\theta_{\P}\ :\ \bbK_{\P}^*\to Z(\P)\subset {\rm Gal}(L/K)$ vue lors de la DŽfinition \argel.  On pose 
$$\eqalign{\Upsilon\ :\ \bbI_{K}&\lra {\rm Gal}(L/K)\cr (a_{\P})_{\P}&\longmapsto \prod_{\P\in \gfP(K)}\theta_{\P}(a_{\P}).\cr}$$
\newcount\gaagba\gaagba=\gaga
Pour presque tout $\P$, $\P$ est non ramifiŽ et $a_{\P}\in U_{\P}$. Donc (Lemme \argeo) $a_{\P}\in N_{L_{\gP}/K_{\P}}(\bbL_{\gP})\buildrel \rm Cor.\ \7argel\over\subset \ker(\theta_{\P})$, donc $\theta_{\P}(a_{\P})=1$, donc $\Upsilon$ est bien dŽfinie. C'est un homomorphisme, car chaque $\theta_{\P}$ l'est.
\bigskip\goodbreak
\lem {\soustitre (rŽciprocitŽ pour le symbole des restes normiques)}
\medskip
{\sl Les idles principaux $K^*\subset \ker(\Upsilon)$, i.e. si $x\in K^*$,
$$\prod_{\P\in\gfP(K)}\theta_{\P}(x)=1.$$

}

{\bf Preuve}

Posons $S$ l'ensemble des places qui ramifient ou qui divisent $x$. Alors $\prod_{\P\in\gfP(K)}\theta_{\P}(x)=\prod_{\P\in S}\theta_{\P}(x)$. Supposons $S=\{\P_{1},\ldots ,\P_{t},\P_{t+1},\ldots ,\P_{s}\}$, o $\P_{1},\ldots ,\P_{t}$ sont des places finies et $\P_{t+1},\ldots ,\P_{s}$ sont des places infinies rŽelles. Supposons que $(x)=\prod_{i=1}^t\P_{i}^{a_{i}}$. Soit $\m$ un $K$-module admissible pour $L/K$ avec $\m=\prod_{i=1}^s\P_{i}^{t_{i}}$ tel que $a_{i}\leq t_{i}$ pour tout $i=1,\ldots , t$ (c'est toujours possible en vertu du thŽorme de rŽciprocitŽ d'Artin \argda). On a $\prod_{\P\in\gfP(K)}\theta_{\P}(x)=\prod_{i=1}^s\theta_{\P_{i}}(x)$, et pour chaque $i$, on a $\m=\P_{i}^{t_{i}}\cdot \m_{i}$ ($\P_{i}\notdiv \m_{i}$ et $\m$ est $\P_{i}$-admissible. Pour chaque $i$, on peut choisir en vertu du ThŽorme d'approximation dŽbile \argc\ $y_{i}\in K^*$ $y_{i}\equiv x\pmodast{\P_{i}^{t_{i}}}$ et $y_{i}\equiv 1\pmodast {\m_{i}}$. On a donc (cf. Remarque prŽcŽdant la Proposition \3argel)  $\theta_{\P_{i}}(x)=\Phi_{L/K}(j_{\m}((y_{i})))$ (avec l'abus  habituel $\Phi_{L/K}=\Phi_{L/K}|_{I_{K}(\m)}$). Remarquons que $x^{-1}y_{1}\cdots y_{s}\in K^*_{\m}$ et que $(y_{i})=\cases{\P_{i}^{a_{i}}\cdot \aa_{i}&avec $\aa_{i}$ premier ˆ $\m$ pour $i=1,\ldots ,t$\cr \aa_{i}&avec $\aa_{i}$ premier ˆ $\m$ pour $i=t+1,\ldots ,s$\cr}$. Ainsi, $j_{\m}((y_{i}))=\aa_{i}$, pour tout $i=1,\ldots , s$. Enfin,
$$\eqalign{\prod_{\P\in\gfP(K)}\theta_{\P}(x)&=\prod_{i=1}^s\Phi_{L/K}(\aa_{i})=\Phi_{L/K}\left (\prod_{i=1}^s \aa_{i}\right )\cr
&=\Phi_{L/K}\left ({(y_{1}\ldots y_{s})\over \prod_{i=1}^t\P_{i}^{a_{i}}=(x)}\right )=\Phi_{L/K}(\underbrace{(x^{-1}y_{1}\cdots y_{s})}_{\in P_{m}})=1.\cr}$$
\qed
\bigskip
\th
\medskip
{\sl Soit $L/K$ une extension abŽlienne de corps de nombres et $\m$ un $K$-module admissible pour $L/K$. Alors l'application $\Upsilon\ :\ \bbI_{K}\to {\rm Gal}(L/K)$ $(a_{\P})_{\P}\mapsto \prod_{\P\in\gfP(K)}\theta_{\P}(a_{\P})$ induit un isomorphisme de $\bbI_{K}/(K^*\cdot N_{L/K}(\bbI_{L}))\to {\rm Gal}(L/K)$, obtenu en composant les homomorphismes~:
$$\bbI_{K}\to \bbI_{K}/(K^*\cdot N_{L/K}(\bbI_{L}))\simeq \bbI'_{\m}/(K^*_{\m}\cdot \bbI_{\m}\cdot N_{L/K}(\bbI''_{L}(\m)))\simeq I_{K}(\m)/(P_{\m}\cdot N_{L/K}(I_{L}(\widetilde{\m})))\buildrel{1\over \Phi_{L/K}}\over\lra {\rm Gal}(L/K).$$

Vous aurez remarquŽ que pour la dernire flche, on a pris ${1\over \Phi_{L/K}}$ au lieu de $\Phi_{L/K}$.

}
{\bf Preuve}

On a $\Upsilon((a_{\P})_{\P})=1$ si $(a_{\P})_{\P}\in N_{L/K}(\bbI_{L})$, car alors $a_{\P}\in N_{L_{\gP}/K_{\P}}(\bbL_{\gP})=\ker(\theta_{\P})$ pour tout $\P$ et $\gP|\P$ (ThŽorme \argem). De mme si $(a_{\P})_{\P}=x\in K^*$ (lemme prŽcŽdent \arggr) et puisque le premier isomorphisme est induit par l'injection, il suffit de vŽrifier l'assertion pour un $a=(a_{\P})_{\P}\in\bbI'_{\m}$.  Alors $\Upsilon((a_{\P})_{\P})=\prod_{\P\in\gfP(K)}\theta_{\P}(a_{\P})=\prod_{\P\in S}\theta_{\P}(a_{\P})$, o $S$ est l'ensemble fini des $\P$ tels que $\theta_{\P}(a_{\P})\ne 1$. Ces $\P$-lˆ sont non ramifiŽs. En effet, si $\P$ est ramifiŽ, alors $\P|\m$. Or, $a\in\bbI'_{\m}$, donc $a\in U_{\P}^{(v_{\P}(\m))}\subset \ker(\theta_{\P})$ au vu de la dŽfinition des $\theta_{\P}$ (DŽfinition \argel) et de la Remarque b) qui suit ($b=a=v_{\P}(\m)$). Par multiplicativitŽ, il suffit de prouver le rŽsultat lorsque $a=(a_{\P})_{\P}$ est tel que $a_{\P}=1$ sauf si $\P=\P_{0}$ non ramifiŽ. Dans ce cas,  $\Upsilon(a)=\theta_{\P_{0}}(a_{\P_{0}})\buildrel \rm Prop. \3argel\over ={\rm Frob}_{L/K}(\P_{0})^{-v_{\P_{0}}(a_{\P_{0}})}$. Et d'autre part, en suivant la suite d'applications de l'ŽnoncŽ du thŽorme, on a $(a_{\P})_{\P}\mapsto \Phi_{L/K}(\psi((a_{\P})_{\P}))^{-1}=\Phi_{L/K}(\P_{0}^{v_{\P_{0}}(a_{\P_{0}})})^{-1}={\rm Frob}_{L/K}(\P_{0})^{-v_{\P_{0}}(a_{\P_{0}})}$. Cela montre le thŽorme.\qed

\vfill\eject


\global\advance\chapnomb by 1
\nomb=1

\centerline{\para Chapitre 14 : }
\bigskip

\centerline{\para Corps de classe local}
\bigskip

La thŽorie du corps de classe se dŽmontre maintenant assez facilement, en utilisant notamment les rŽsultats sur $\theta_{\P}$, le symbole de norme rŽsiduel, vus au Chapitre 11. Pour fixer les idŽes, nous allons donner une dŽfinition pour nous mettre d'accord sur la notion de corps locaux (pas comme ce coquin de Serre qui fait un livre dessus et qui ne prend mme pas la peine de dire ce que c'est...)

\bigskip

\defi
\medskip
On appellera {\it corps local } l'un des corps suivant~:
$$\R,\C,\hbox{toute extension finie d'un corps $\Q_{p}$}$$
(en fait il s'agit de la liste complte des corps topologiques  localement compacts (non discrets) de caractŽristique 0). Pour simplifier, on notera $\|\cdot\|$ la valeur absolue du corps local considŽrŽ. Il est Žvident que si $K$ est un corps de nombres et $\P\in\gfP(K)$, alors $\bbK_{\P}$ est un corps local. On a mieux~:
\bigskip
\headline={\hfill \phantom{ouuh}\hfill}
\th
\medskip
{\sl Si $\bbK$ est un corps local, alors il existe un corps de nombres $K$ et $\P$ une place de $K$ tels que $\bbK=\bbK_{\P}$.

}
{\bf Preuve}

Si $\bbK=\C$ ou $\R$, c'est Žvident. On aura besoin de 3 lemmes et d'une dŽfinitions pour prouver ce rŽsultat dans le cas non archimŽdien. 
\bigskip
\lem {\soustitre (Lemme de Krasner)}
\medskip
{\sl Soit $\bbK$ un corps local non archimŽdien, $\bbK^{\rm alg}$ une cl™ture algŽbrique de $\bbK$ et $\alpha,\beta\in \bbK^{\rm alg}$. Notons  $\alpha=\alpha_1,\alpha_2,\ldots ,\alpha_n$ les conjuguŽes de $\alpha$ sur $K$. Supposons que $\|\alpha-\beta\|<\|\alpha_i-\beta\|\ \forall i=2,\ldots ,n$. Alors $\bbK(\alpha)\subset \bbK(\beta)$.

}

{\bf Preuve}

Supposons d'abord connu le rŽsultat suivant~: si $\bbL/\bbK$ est une extension algŽbrique de corps locaux, alors la valeur absolue de $\bbK$ se prolonge de manire unique sur $\bbL$ par la formule 
$$\forall\gamma\in \bbL \ \| \gamma \| =\|  N_{\bbK(\gamma)/\bbK}(\gamma)\|^{1\over [\bbK(\gamma):\bbK]}\eqno{(*)}$$
(cf. [Fr-Tay, 1.15, p.113]) ce qui implique que si $\gamma_1$ et $\gamma_2\in \bbL$ sont conjuguŽs (i.e. $m_{\gamma_1/\bbK}=m_{\gamma_2/\bbK}$), alors $\|\gamma_1\|=\|\gamma_2\|$.
\bigskip
Prouvons maintenant le lemme. Supposons par l'absurde que  $\alpha\not\in \bbK(\beta)$. Alors le polyn™me minimal de $\alpha$ sur $\bbK(\beta)$ est de degrŽ $>1$ et divise le polyn™me minimal de $\alpha$ sur $\bbK$. Il admet donc parmi ses racines un $\alpha_{i_0}$ pour un $i>1$. Ainsi, $\alpha$ et $\alpha_{i_0}$ sont conjuguŽs sur $\bbK(\beta)$ et donc $\beta-\alpha$ et $\beta-\alpha_{i_0}$ sont aussi conjuguŽs sur $\bbK(\beta)$ (car le changement de variable $x\mapsto -x+\beta$ est unimodulaire). Par l'ŽgalitŽ $(*)$, on a donc que $\|\beta-\alpha\|=\|\beta-\alpha_{i_0}\|$, contredisant l'hypothse.\qed
\bigskip

\lem

{\sl Soit $\bbK$ un corps local non archimŽdien, $f=x^n+a_{n-1} x^{n-1}+\cdots + a_0\in \bbK[x]$ et $\alpha$ une racine de $f$ dans une extension de $\bbK$. Alors on a~:
$$\|\alpha\|\leq \max (1,\|a_0\|,\ldots ,\|a_{n-1}\|).$$

}

{\bf Preuve}

Supposons {\it ab absurdo} que $\|\alpha\|> \max (1,\|a_0\|,\ldots ,\|a_{n-1}\|).$ Alors
$$\left\|\sum_{i=0}^{n-1} a_i\cdot \alpha^i\right \|\leq  \max_{0\leq i\leq n-1}(\|a_i\|\cdot \|\alpha\|^i)<\|\alpha\|\cdot \|\alpha\|^{n-1}=\|\alpha\|^n.$$
Cela montre que $0\ne \alpha^n+\sum_{i=0}^{n-1} a_i\cdot \alpha^i=f(\alpha)=0$. C'est une contradiction .\qed
\bigskip
\defi

Soit $f=\sum_{i=0}^{n} a_i\cdot x^i\in \bbK[x]$,  $n\geq 1$. On dŽfinit $\|f\|=\max_{i=0}^n(\|a_i\|)$. On vŽrifie facilement que c'est une norme sur $\bbK[x]$. \bigskip

\headline={\hfill\smcap Corps de classe local\hfill}

\lem{\soustitre (Lemme de continuitŽ des racines)}
\medskip
{\sl  

Soit $\bbK$ un corps local non archimŽdien, $\bbK^{\rm alg}$ une cl™ture algŽbrique de $\bbK$, $f=x^n+a_{n-1} x^{n-1}+\cdots + a_0\in \bbK[x]$ un polyn™me irrŽductible dans $\bbK[x]$. Soit encore $\alpha=\alpha_1,\alpha_2,\ldots ,\alpha_n$ les racines de $f$ dans $\bbK^{\rm alg}$ qui sont bien sžr toutes distinctes, car on est en caractŽristique 0. Alors pour tout $\varepsilon>0$, il existe $\delta>0$ tel que pour tout $g=x^n+b_{n-1} x^{n-1}+\cdots + b_0\in\bbK[x]$ satisfaisant $\|f-g\|<\delta$, alors $g$ admet une racine $\beta$ telle que $\|\alpha-\beta\|<\varepsilon$.

Remarquons qu'en prenant $\delta$ encore plus petit, on peut avoir  le rŽsultat pour tout $\alpha_i$.
}

{\bf Preuve}

On peut supposer que $2\cdot \varepsilon<\rho:= \min_{i\ne j}\|\alpha_i-\alpha_j\|$. On va chercher un $\delta<1$ satisfaisant la conclusion. Pour tout $g=x^n+b_{n-1} x^{n-1}+\cdots + b_0\in \bbK[x]$ satisfaisant $\|f-g\|<\delta<1$ et pour tout $\mu$ racine de $g$, on aura
$$\eqalign{\|\mu\|&\buildrel \rm Lem.\ \aaarggv\over\leq \max(1,\|b_0\|,\ldots,\|b_{n-1}\|)\cr
&\leq \max(1,\|a_0\|+1,\ldots,\|a_{n-1}\|+1)\cr
&\leq \|f\|+1\cr}\eqno{(*)}$$
qui est indŽpendant de $g,\beta$ et $\delta$ $(<1)$. Avec ces hypothses sur $g,\beta$ et $\delta$, on a  encore~:
$$\eqalign{\prod_{i=1}^n\|\mu-\alpha_i\|&=\left \|\prod_{i=1}^n(\mu-\alpha_i)\right \|=\|f(\mu)\|\cr &=\|f(\mu)-g(\mu)\|=\left\|\sum_{i=0}^{n-1}(a_i-b_i)\cdot \mu^i\right\|\cr
&\leq \|f-g\|\cdot \max_{0\leq i\leq n-1}\|\mu\|^i\cr
&\leq \|f-g\|\cdot \max(1,\|\mu\|^{n-1})\leq\delta\cdot \max(1,\|\mu\|)^{n-1}\cr
&\buildrel (*)\over \leq \delta\cdot(\|f\|+1)^{n-1}\cr}$$
Si $\delta$ est suffisamment petit, alors $\delta\cdot(\|f\|+1)^{n-1}<\varepsilon^n$. Donc, il existe $i$ tel que $$\|\mu-\alpha_i\|<\varepsilon,\eqno{(**)}$$ et cet $i$ est unique. En effet, s'il existe un autre $j$, on aura $\|\alpha_i-\alpha_j\|\leq \|\mu-\alpha_i\|+\|\mu-\alpha_i\|<2\cdot\varepsilon<\rho=\min_{i\ne j}\|\alpha_i-\alpha_j\|$, ce qui est impossible. Mais hŽlas, notre $\alpha_{i}$, n'est peut-tre pas $\alpha$.

Soit $\bbL/\bbK$ une extension galoisienne finie de $\bbK$ qui contient  les $\alpha_{i}$ et $\mu $. Soit $\sigma\in {\rm Gal}(\bbL/\bbK)$ tel que $\sigma(\alpha_{i})=\alpha$ ($\sigma$ existe, car $f$ est irrŽductible). On trouve enfin 
$$\|\sigma(\mu)-\sigma(\alpha_{i})\|=\|N_{\bbL/\bbK}(\sigma(\mu-\alpha_i))\|^{1\over[\bbL:\bbK]}=\|N_{\bbL/\bbK}(\mu-\alpha_i)\|^{1\over[\bbL:\bbK]}=\|\mu-\alpha_i\|<\varepsilon$$
 et enfin, $\beta:=\sigma(\mu)$ est une racine de $g$ et satisfait la conclusion de notre lemme.\qed
\bigskip\goodbreak
\rem
\medskip
Ce rŽsultat est vrai aussi en caractŽristique positive, de mme si $f$ n'est pas irrŽductible, mais la preuve est un poil plus longue.
\bigskip
{\soustitre Fin de la preuve du ThŽorme \arggu}

Soient $\alpha$ tel que $\bbK=\Q_{p}(\alpha)$, $f=m_{\alpha/\Q_{p}}$ le polyn™me minimal de $\alpha$, $n=[\bbK:\Q_{p}]$ le degrŽ de $f$, $\alpha=\alpha_{1},\alpha_{2},\ldots ,\alpha_{n}$ les racines de $f$ dans une cl™ture algŽbrique de $\bbK$. Soit $\varepsilon>0$ tel que $\varepsilon <\min_{{i\ne j}}\|\alpha_{i}-\alpha_{j}\|$ et soit $\delta >0$ comme dans l'ŽnoncŽ du lemme prŽcŽdent (Lemme \aaarggx). Puisque $\Q$ est dense dans $\Q_{p}$, il existe un polyn™me $g\in \Q[x]$ unitaire de degrŽ $n$ tel que $\|f-g\|<\delta$. Par le lemme prŽcŽdent (Lemme \aaarggx), $g$ a une racine $\beta$ telle que $\|\alpha-\beta\|<\varepsilon$. On prŽtend que $\|\beta-\alpha_{i}\|>\|\beta-\alpha\|$ si $i=2,\ldots ,n$. En effet, sinon on aurait $\|\alpha-\alpha_{i}\|\leq \max (\|\alpha-\beta\|,\|\beta-\alpha_{i}\|)=\|\alpha-\beta\|<\varepsilon$, ce qui est contraire au choix de $\varepsilon$. Par le Lemme de Krasner (Lemme \aaarggu), on a $\Q_{p}(\beta) \supset\Q_{p}(\alpha)=\bbK$. Mais, on a $[\Q_{p}(\beta):\Q_{p}]\leq {\rm deg}(g)=n=[\Q_{p}(\alpha) :\Q_{p}]$. Cela montre alors que $\bbK=\Q_{p}(\beta)$. Il est alors clair que $K:=\Q(\beta)$ est dense dans $\bbK$, donc $\bbK$ est le complŽtŽ de $K$ et donc $\bbK=\bbK_{\P}$ pour un idŽal $\P |p$.\qed

\bigskip
Un corollaire un peu plus fort :
\bigskip
\th 
\medskip
{\sl Soit $\bbL/\bbK$ une extension galoisienne de corps locaux (donc forcŽment finie). Alors il existe une extension $L/K$ galoisienne de corps de nombres telle que $K\subset \bbK$, $L\subset \bbL$ et des places $\P\in\gfP(K)$ et $\gP\in\gfP(L)$ telles que $\bbL=\bbL_{\gP}$ et $\bbK=\bbK_{\P}$. On peut mme supposer que ${\rm Gal}(L/K)={\rm Gal}(\bbL/\bbK)$.

}

{\bf Preuve}

Soit $G={\rm Gal}(\bbL/\bbK)$. On choisit $L$ un corps de nombres et $\gP$ une place de $L$ telle $\bbL=\bbL_{\gP}$ (c'est possible en vertu du ThŽorme\arggu). Alors $\bbL=\bbK\cdot L$, car $\bbK\cdot L$ est bien un sous-corps local de $\bbL$ dans lequel $L$ dense. Rempla\c cant $L$  par $\prod_{\sigma\in G}\sigma(L)$, on peut supposer que $L$ est stable par les ŽlŽments de $G$. Alors $\sigma\mapsto\sigma|_{L}$ est un homomorphisme injectif de $G$ dans le groupe des automorphismes de $L$. Soit $K$ le sous-corps de $L$ fixe par les $\sigma|_{L}$. Alors $L/K$ est une extension galoisienne de groupe de galois $\simeq G$, et $K\subset \bbK$. Alors posons $\P=\gP|_{K}$ et $\bbK_{\P}$ l'adhŽrence de $K$ dans $\bbK$, on a Žvidemment $\bbK_{\P}\subset\bbK$. Mais, $\bbL/\bbK_{\P}$ est une extension galoisienne de groupe de Galois $\simeq Z(\gP/\P)$. Donc, $[\bbL:\bbK_{\P}]=|Z(\gP/\P)|\leq |G|=[\bbL :\bbK]$, ce qui prouve que $\bbK=\bbK_{\P}$.\qed
\bigskip
\lem
\medskip

{\sl Soit $\bbL/\bbK$ une extension abŽlienne de corps locaux, et pour $i=1,2$ $L_{i}/K_{i}$ deux extensions abŽliennes de corps de nombres avec des places $\gP_{i}$ de $L_{i}$ et $\P_{i}$ de $K_{i}$ telles que $\gP_{i}|\P_{i}$ et $\bbL_{1\gP_{1}}=\bbL_{2\gP_{2}}=\bbL$ et $\bbK_{1\P_{1}}=\bbK_{2\P_{2}}=\bbK$. Alors on a 
$$\theta_{\P_{1}}(L_{1}/K_{1})=\theta_{\P_{2}}(L_{2}/K_{2}),$$
moyennant bien sžr les identifications de $Z(\gP_{i}/\P_{i})$ avec ${\rm Gal}(\bbL/\bbK)$.

}
{\bf Preuve}
\vglue 2.5cm
On a le diagramme~: 

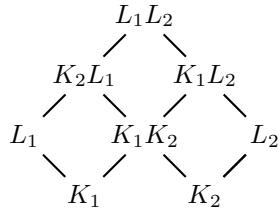
\rput(7,0){\rnode{a1}{$ K_{1}$}}
\rput(9,0){\rnode{a2}{$K_{2}$}}
\rput(6,1){\rnode{b1}{$L_{1}$}}
\rput(8,1){\rnode{b2}{$K_{1}K_{2}$}}
\rput(10,1){\rnode{b3}{$L_{2}$}}
\rput(7,2){\rnode{c1}{$K_{2}L_{1}$}}
\rput(9,2){\rnode{c2}{$K_{1}L_{2}$}}
\rput(8,3){\rnode{d1}{$L_{1}L_{2}$}}
\ncline[nodesep=3pt]{a1}{b1}
\ncline[nodesep=3pt]{a1}{b2}
\ncline[nodesep=3pt]{a2}{b2}
\ncline[nodesep=3pt]{a2}{b3}
\ncline[nodesep=3pt]{b1}{c1}
\ncline[nodesep=3pt]{b2}{c1}
\ncline[nodesep=3pt]{b2}{c2}
\ncline[nodesep=3pt]{b3}{c2}
\ncline[nodesep=3pt]{c1}{d1}
\ncline[nodesep=3pt]{c2}{d1}
\ncline[nodesep=3pt]{a2}{b3}

Les places correspondantes sont 

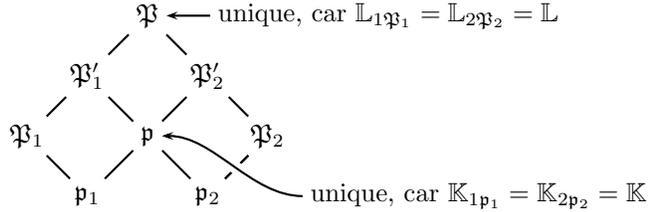
\vglue 2.5cm

\rput(7,0){\rnode{a1}{$ \P_{1}$}}
\rput(9,0){\rnode{a2}{$\P_{2}$}}
\rput(6,1){\rnode{b1}{$\gP_{1}$}}
\rput(8,1){\rnode{b2}{$\P$}}
\rput(10,1){\rnode{b3}{$\gP_{2}$}}
\rput(7,2){\rnode{c1}{$\gP'_{1}$}}
\rput(9,2){\rnode{c2}{$\gP'_{2}$}}
\rput(8,3){\rnode{d1}{$\gP$}}
\rput(12,3){\rnode{e}{unique, car $\bbL_{1\gP_{1}}=\bbL_{2\gP_{2}}=\bbL$}}
\rput(13.5,0){\rnode{f}{unique, car $\bbK_{1\P_{1}}=\bbK_{2\P_{2}}=\bbK$}}
\ncline[nodesep=3pt]{a1}{b1}
\ncline[nodesep=3pt]{a1}{b2}
\ncline[nodesep=3pt]{a2}{b2}
\ncline[nodesep=3pt]{a2}{b3}
\ncline[nodesep=3pt]{b1}{c1}
\ncline[nodesep=3pt]{b2}{c1}
\ncline[nodesep=3pt]{b2}{c2}
\ncline[nodesep=3pt]{b3}{c2}
\ncline[nodesep=3pt]{c1}{d1}
\ncline[nodesep=3pt]{c2}{d1}
\ncline[nodesep=3pt]{a2}{b3}
\ncline[nodesep=3pt]{->}{e}{d1}
\nccurve[border=2pt,nodesep=3pt,angleA=180]{->}{f}{b2}
\bigskip

En appliquant le deuxime lemme de naturalitŽ (Proposition \6argel) avec $K=K_2$, $L=L_2$ et $E=K_1K_2$, on trouve $\theta_{\P_2}(L_2/K_2)\circ N_{\bbK/\bbK}=R\circ \theta_\P(K_1L_2/K_1K_2)$, o $R\ : {\rm Gal}(K_1L_2/K_2K_2)\to {\rm Gal}(L_2/K_2)$ est la restriction ˆ $L_2$. Or, l'image de $\theta_\P(K_1L_2/K_1K_2)$ est $Z(\gP'_2/\P)\simeq {\rm Gal}(\bbL/\bbK) $ et l'image de $\theta_{\P_2}(L_2/K_2)$ est $Z(\gP_2/\P_2)\simeq {\rm Gal}(\bbL/\bbK)$ (cf. ThŽorme \argek). Dans ces conditions, on peut supposer (avec un petit abus) que $R={\rm id}$ et donc que $\theta_{\P_2}(L_2/K_2)=\theta_\P(K_1L_2/K_1K_2)$. Par un mme raisonnement, on a $\theta_{\P_1}(L_1/K_1)=\theta_{\P}(K_2L_1/K_1K_2)$. En appliquant les premier lemme de naturalitŽ (Proposition \5argel) avec $K=K_1K_2$, $L=K_1L_2$ et $E=L_1L_2$, on a $\theta_\P(K_1L_2/K_1K_2)=R\circ \theta_\P(L_1L_2/K_1K_2)$, o (comme avant) $R\ : Z(\gP/\P)\to Z(\gP'_2/\P)$ qui peut aussi tre vu comme l'identitŽ car ˆ nouveau chacun de ces groupes est isomorphe ˆ ${\rm Gal}(\bbL/\bbK)$. Puis, on montre de manire identique que $\theta_\P(K_2L_1/K_1K_2)=\theta_\P(L_1L_2/K_1K_2)$. Enfin,
$$\theta_{\P_1}(L_1/K_1)=\theta_{\P}(K_2L_1/K_1K_2)=\theta_\P(L_1L_2/K_1K_2)=\theta_\P(K_1L_2/K_1K_2)=\theta_{\P_2}(L_2/K_2).$$
\qed
\bigskip\goodbreak
{{\soustitre Proposition-DŽfinition ({\uppercase\expandafter{\the\chapnomb}}.{\the\nomb})}\global\advance\nomb by 1}
\medskip
{\sl Soit $\bbL/\bbK$ une extension abŽlienne de corps locaux. Alors il existe un homomorphisme 
$$\theta(\bbL/\bbK)\ :\ \bbK^*\lra {\rm Gal}(\bbL/\bbK)$$
surjectif de noyau $N_{\bbL/\bbK}(\bbL^*)$ ayant les propriŽtŽs suivantes~:

}
\newcount\gaagbb\gaagbb=\gaga
{\sl
\art{a)}Si $L/K$ est une extension abŽlienne de corps de nombres et $\gP|\P$ des places de $L$ respectivement de $K$ tels que $\bbL_\gP=\bbL$ et $\bbK_\P=\bbK$. Alors 
$$\theta(\bbL/\bbK)=\theta_\P(L/K).$$

\art{b)}Si en plus $\bbL/\bbK$ est une extension non ramifiŽe et si $\pi$ est une uniformisante de $\bbK$, alors, pour tout $k\in \Z$ et $u\in U_\P$ (les unitŽs de $\bbK$)~:
$$\theta(\bbL/\bbK)(\pi^k\cdot u)={\rm Frob}(\bbL/\bbK)^{-k}.$$

\art{c)}Si $\bbE/\bbK$ et $\bbL/\bbK$ sont des extensions de corps locaux tels que $\bbL/\bbK$ est abŽlienne, alors 
$$\theta(\bbL/\bbK)\circ N_{\bbE/\bbK}=R\circ \theta(\bbE\bbL/\bbE),$$
o $R\ :\ {\rm Gal}(\bbE\bbL/\bbE)\to{\rm Gal}(\bbL/\bbK)$ est la restriction ˆ $\bbL$.

\art{d)}Soit $\bbE/\bbL$ et $\bbL/\bbK$ des extensions de corps locaux telles que $\bbE/\bbK$ est abŽlienne. Alors 
$$\theta(\bbL/\bbK)=R\circ \theta(\bbE/\bbK),$$
o $R\ :\ {\rm Gal}(\bbE/\bbK)\to{\rm Gal}(\bbL/\bbK)$ est la restriction ˆ $\bbL$.

}
{\bf Preuve}

L'existence de cet homomorphisme dŽcoule des ThŽormes \argek\ et \argem.  La partie a) vient du ThŽorme \arggv\ et du lemme prŽcŽdent \arggw. La partie b) provient de la Proposition \3argel. La partie c) est le deuxime lemme de naturalitŽ \6argel. Et la partie d) est le premier lemme de naturalitŽ \5argel. \phantom{ihig}\hfill\qed
\bigskip
Maintenant, nous devons faire un petit dŽtour par les extensions de Kummer de corps locaux (voir DŽfinition \argdt\ pour la dŽfinition d'extension de Kummer).
\bigskip
\prop
\medskip
{\sl Soit $\bbK$ un corps local et $n\in\N$ un nombre entier. Alors $\bbK^*/(\bbK^*)^n$ est fini. 

Supposons que $\bbK$ contienne une racine primitive $n$-ime de l'unitŽ. Posons $\bbL=\bbK(\root n\of{\bbK^*})$. Alors c'est une extension de Kummer de degrŽ fini (c'est l'extension de Kummer maximale de $\bbK$) et on a $N_{\bbL/\bbK}(\bbL^*)=(\bbK^*)^n$.

}
{\bf Preuve}

Soit $\pi$ une uniformisante de $\bbK$. On sait que $\bbK^*\simeq <\pi>\hskip-4pt\times U_\P\simeq \Z\times U_\P$ o $\P$ est l'idŽal maximal de $\bbK$ et $U_{\P}$ est le sous-groupe des unitŽs de $\bbK$. Et $(\bbK^*)^n\simeq n\Z\times U_\P^n$. Montrons que $U_\P/U_\P^n$ est fini. La Proposition \argcc\ nous assure l'existence d'un entier $m$ tel que $U_\P^{(m)}\subset U_\P^n\subset U_\P$. Or, $U_\P/U_\P^{(m)}$ est fini, nous l'avons montrŽ au dŽbut de la preuve du Lemme \argbz. En vertu du ThŽorme \argdu\ et de ce qui prŽcde, $\bbL/\bbK$ est l'extension de Kummer maximale de $\bbK$. Le mme ThŽorme \argdu\ nous montre que $\bbK^*/(\bbK^*)^n\simeq {\rm Gal}(\bbL/\bbK)$. La proposition-dŽfinition prŽcŽdente \arggx\ montre aussi que $\bbK^*/N_{\bbL/\bbK}(\bbL^*)\simeq {\rm Gal}(\bbL/\bbK)$. Cela montre que $N_{\bbL/\bbK}(\bbL^*)=(\bbK^*)^n$, en effet, on a d'autre part  $(\bbK^*)^n\subset N_{\bbL/\bbK}(\bbL^*)$, car puisque $n$ est un exposant de ${\rm Gal}(\bbL/\bbK)$ (cf. ThŽorme \argdu), on a $(\bbK^*)^n\subset \ker(\theta(\bbL/\bbK))=N_{\bbL/\bbK}(\bbL^*)$.\qed
\bigskip
Maintenant un petit lemme de thŽorie des groupes~: 
\bigskip
{{\soustitre Proposition-DŽfinition ({\uppercase\expandafter{\the\chapnomb}}.{\the\nomb})}\global\advance\nomb by 1}
\medskip
{\sl Soit $G$ un groupe fini et $H$ un sous-groupe. On note $D(G,H)$ le sous-groupe engendrŽ par les commutateurs $ghg^{-1}h^{-1}$, ($g\in G, h\in H$). Alors $D(G,H)$ est un sous-groupe normal de $G$. De plus, l'image de $H$ via l'homomorphisme canonique $G\to \overline{G}=G/D(G,H)$ est dans le centre de $\overline{G}$. Enfin, si $H$ est normal dans $G$ et que $G/H$ est cyclique, alors $G/D(G,H)$ est abŽlien.

}
{\bf Preuve}

Tout d'abord, on montre que $x D(G,H)x^{-1}=D(G,H)$ pour tout $x\in G$. Soit $ghg^{-1}h^{-1}\in D(G,H)$. On a
$$xghg^{-1}h^{-1}x^{-1}=\underbrace{\left [(xg)h(xg)^{-1}h^{-1}\right ]}_{\in D(G,H)}\cdot \underbrace{\left (xhx^{-1}h^{-1}\right )^{-1}}_{\in D(G,H)}.$$
Le fait que l'image de $H$ dans $\overline{G}$ est dans le centre est Žvident. Montrons la dernire partie~: on note $\overline{g}$ l'image de $g\in G$ dans $\overline{G}$. On choisit $b\in G$ tel que sa classe modulo $H$ engendre $G/H$. Alors tout ŽlŽment de $G$ s'Žcrit $h\cdot b^m$ ($h\in H, m\in \Z$). Donc tout ŽlŽment de $\overline{G}$ s'Žcrit $\overline{h}\cdot\overline{b}^m$. Soit $\overline{h'}\cdot \overline{b}^k$ un autre ŽlŽment de $\overline{G}$. Alors 
$$(\overline{h}\cdot\overline{b}^m)\cdot (\overline{h'}\cdot \overline{b}^k)=(\overline{h'}\cdot \overline{b}^k)\cdot (\overline{h}\cdot\overline{b}^m)$$
car $\overline{h}$ et $\overline{h'}$ sont dans le centre de $\overline{G}$ et que $\overline{b}^m\cdot \overline{b}^k=\overline{b}^{m+k}=\overline{b}^k\cdot \overline{b}^m$.\qed
\bigskip\goodbreak
\th
\medskip
{\sl
Soit $\bbL/\bbK$ une extension galoisienne de degrŽ fini de corps locaux. Supposons que ${\rm Gal}(\bbL/\bbK)$ soit rŽsoluble. Notons $\bbL^{\rm ab}$ la sous-extension abŽlienne maximale de $\bbL/\bbK$. Alors 
$$N_{\bbL/\bbK}(\bbL^*)=N_{\bbL^{\rm ab}/\bbK}({\bbL^{\rm ab}}^*).$$
(En fait le rŽsultat est vrai mme dans le cas non rŽsoluble, mais nous n'avons pas besoin de le montrer,  ce serait plus long et la situation dans laquelle nous utiliserons ce rŽsultat sera clairement rŽsoluble).
}
{\bf Preuve}

On a clairement $N_{\bbL^{\rm ab}/\bbK}({\bbL^{\rm ab}}^*)\supset N_{\bbL/\bbK}(\bbL^*)$.

Pour montrer l'inclusion inverse, on procde par rŽcurrence sur $n=[\bbL:\bbK]$. Si $n=1$, c'est Žvident. Supposons donc $n>1$ et le thŽorme vrai pour tout $m<n$. Puisque ${\rm Gal}(\bbL/\bbK)$ est rŽsoluble, on voit facilement qu'il admet un quotient cyclique non trivial (si vous n'tes pas convaincu, voyez [La1, Prop. 1.3.1, p.20]). Autrement dit, il existe une sous-extension $\bbE/\bbK$ de $\bbL/\bbK$ cyclique de degrŽ $>1$. Soit $\bbM/\bbE$ la sous-extension abŽlienne maximale de $\bbL/\bbE$. Il est Žvident que ${\rm Gal}(\bbL/\bbM)$ est le sous-groupe des commutateurs de ${\rm Gal}(\bbL/\bbE)$ (on avait dŽjˆ fait ce raisonnement au ThŽorme \argeh) et   ${\rm Gal}(\bbL/\bbE)$ est normal dans ${\rm Gal}(\bbL/\bbK)$. Cela implique que ${\rm Gal}(\bbL/\bbM)$ est normal dans ${\rm Gal}(\bbL/\bbK)$ (vŽrification facile, sinon, voir [Jac1, rel. 26, p. 246]). Donc l'extension $\bbM/\bbK$ est galoisienne. Notons $G={\rm Gal}(\bbM/\bbK)$ et $H={\rm Gal}(\bbM/\bbE)$. Ainsi, $H$ est un sous-groupe normal, abŽlien de $G$ et $G/H={\rm Gal}(\bbE/\bbK)$ est cyclique. Soit $\bbF$ le corps fixe par $D(G,H)$. Par le lemme prŽcŽdent \arggz, l'extension $\bbF/\bbK$ est abŽlienne. On est donc dans la situation suivante~:

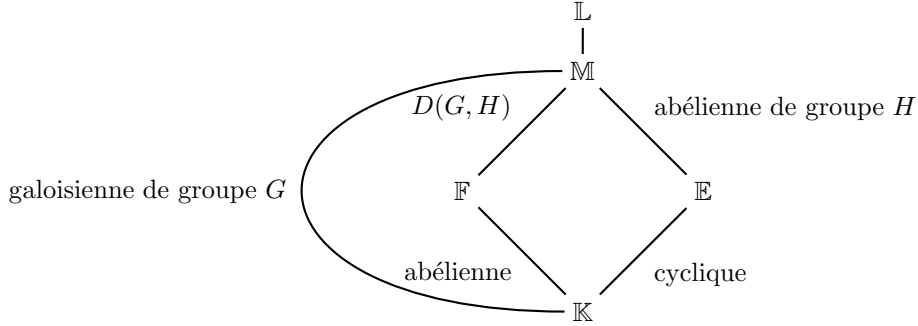
\vglue4cm

\rput(9,0){\rnode{a1}{$ \bbK$}}
\rput(7,2){\rnode{b1}{$\bbF$}}
\rput(11,2){\rnode{b2}{$\bbE$}}
\rput(9,4){\rnode{c1}{$\bbM$}}
\rput(9,5){\rnode{d1}{$\bbL$}}
\ncline[nodesep=3pt]{a1}{b1}
\Aput{abŽlienne}
\ncline[nodesep=3pt]{a1}{b2}
\Bput{cyclique}
\ncline[nodesep=3pt]{b2}{c1}
\Bput{abŽlienne de groupe $H$}
\ncline[nodesep=3pt]{b1}{c1}
\Aput{$D(G,H)$}
\ncline[nodesep=3pt]{c1}{d1}
\nccurve[ncurv=2.9,border=2pt,nodesep=3pt,angleA=180,angleB=180]{a1}{c1}
\Aput{galoisienne de groupe $G$}
\bigskip
De la partie c) de la Proposition-DŽfinition \arggx\ appliquŽe ˆ $\bbL=\bbF$, $\bbE=\bbE$ et $\bbK=\bbK$, on trouve $R_{\bbF\bbE\to\bbF}\circ\theta(\bbF\bbE/\bbE)=\theta(\bbF/\bbK)\circ N_{\bbE/\bbK}$. La partie d) de cette mme proposition \arggx\ appliquŽ ˆ $\bbE=\bbM$, $\bbL=\bbF\bbE$ et $\bbK=\bbE$, montre que $\theta(\bbF\bbE/\bbE)=R_{\bbM\to\bbF\bbE}\circ\theta(\bbM/\bbE)$. En combinant tout cela, on trouve
$$R_{\bbM\to\bbF}\circ\theta(\bbM/\bbE)=\theta(\bbF/\bbK)\circ N_{\bbE/\bbK}.\eqno{(*)}$$

Rappelons que nous voulons montrer que $N_{\bbL^{\rm ab}/\bbK}({\bbL^{\rm ab}}^*)\subset N_{\bbL/\bbK}(\bbL^*)$. Soit donc $x\in N_{\bbL^{\rm ab}/\bbK}({\bbL^{\rm ab}}^*)$. Puisque $\bbL^{\rm ab}\supset \bbF,\bbE$, alors $x\in N_{\bbF/\bbK}(\bbF^*)\cap N_{\bbE/\bbK}(\bbE^*)$. Choisissons $y\in\bbE^*$ tel que $N_{\bbE/\bbK}(y)=x$. Appliquant l'ŽgalitŽ $(*)$ ˆ $y$, on trouve
$$\theta(\bbM/\bbE)(y)|_\bbF=\theta(\bbF/\bbK)(x)={\rm id}_{\bbF},$$
car $x\in N_{\bbF/\bbK}(\bbF^*)=\ker (\theta(\bbF/\bbK))$. Cela montre que $\theta(\bbM/\bbE)(y)\in D(G,H)$, ce qui veut dire que $\theta(\bbM/\bbE)(y)$ est un produit d'ŽlŽments de la forme $\rho \theta(\bbM/\bbE)(z)\rho^{-1}\theta(\bbM/\bbE)^{-1}(z)$, avec $\rho\in G$ et $z\in \bbE^*$, ceci par dŽfinition de $D(G,H)$ et parce que $\theta(\bbM/\bbE)$ est une surjection sur $H$. Par le troisime lemme de naturalitŽ (Proposition \8argel) on a $\rho \theta(\bbM/\bbE)(z)\rho^{-1}=\theta(\bbM/\bbE)(\rho(z))$ (car $\rho(\bbM)=\bbM$ et $\rho(\bbE)=\bbE$, puisque $\bbE/\bbK$ est galoisien). Donc  $\theta(\bbM/\bbE)(y)$ est un produit d'ŽlŽments de la forme $\theta(\bbM/\bbE)({\rho(z)\over z})$, disons,
$$\theta(\bbM/\bbE)(y)=\theta(\bbM/\bbE)\left (\prod_{i=1}^r{\rho_{i}(z_{i})\over z_{i}}\right ),$$
avec $z_{1},\ldots ,z_{r}\in \bbE^*$ et $\rho_{1},\ldots ,\rho_{r}\in G$. Ainsi, 
$$y':=y\cdot \left (\prod_{i=1}^r{\rho_{i}(z_{i})\over z_{i}}\right )^{-1}\in \ker(\theta(\bbM/\bbE))=N_{\bbM/\bbE}(\bbM^*)\buildrel\rm hyp\ de\ rec.\over =N_{\bbL/\bbE}(\bbL^*).$$
Il existe donc $l\in\bbL^*$ tel que $N_{\bbL/\bbE}(l)=y'$. D'autre part, puisque $\rho_{i}|_{\bbE}\in {\rm Gal}(\bbE/\bbK)$, on a $N_{\bbE/\bbK}({\rho_{i}(z_{i})\over z_{i}})=1$ pour tout $i$. Donc, $N_{\bbE/\bbK}(y')=N_{\bbE/\bbK}(y)\cdot \underbrace{N_{\bbE/\bbK} \left (\prod_{i=1}^r{\rho_{i}(z_{i})\over z_{i}}\right )^{-1}}_{=1}=x$. Finalement,
$$x=N_{\bbE/\bbK}(y)=N_{\bbE/\bbK}(y')=N_{\bbE/\bbK}(N_{\bbL/\bbE}(l))=N_{\bbL/\bbK}(l),$$
ce qui achve la (jolie) preuve du thŽorme.\qed
\bigskip
\th {\soustitre \ (Corps de classe local)}
\medskip
{Soit $\bbK$ un corps local. La correspondance
$$\bbL\lra N_{\bbL/\bbK}(\bbL^*)$$
est une bijection (qui renverse l'inclusion) de l'ensemble des extensions abŽliennes de degrŽ fini $\bbL$ de $\bbK$ (contenue dans une mme cl™ture algŽbrique de $\bbK$), et les sous-groupes ouverts d'indices finis de $\bbK^*$.

}

{\bf Preuve}

Remarquons au passage que dans les cas archimŽdiens, il n'y a pas besoin de toute cette thŽorie : $\C^*$ est le seul sous-groupe ouvert de $\C^*$ et $\R^*$ et $\R_{+}^*$ sont les seuls sous-groupes ouverts de $\R^*$.

D'autre part, il serait judicieux de s'assurer que si $\bbL/\bbK$ est une extension abŽlienne, alors $N_{\bbL/\bbK}(\bbL^*)$ est un ouverts de $\bbK^*$. En effet, si $\bbK$ est archimŽdien, c'est Žvident. Sinon, on sait que, via une uniformisante,  $\bbK^*\simeq \Z\times U_{\P}$ (topologie produit, $\Z$ est muni de la topologie discrte et $U_{\P}$ de la topologie hŽritŽe de $\bbK$) qui est un homŽomorphisme de groupes topologiques. Or, en vertu de la Proposition \argcc, il existe un entier $b$ tel que 
$$N_{\bbL\bbK}(\bbL^*)\supset (\bbK^*)^n\simeq n\Z\times  U_{\P}^n\supset n\Z\times U_{\P}^{(b)},$$
et l'ensemble de droite est clairement un sous-groupe ouvert, donc $(\bbK^*)^n\simeq \bigcup_{x\in U_{\P}^n}(n\Z\times x U_{\P}^{(b)})$ et $N_{\bbL\bbK}(\bbL^*)=\bigcup_{x\in N_{\bbL\bbK}(\bbL^*)} x (\bbK^*)^n$ sont aussi ouverts.

Ensuite, montrons que la correspondance $\bbL\to N_{\bbL/\bbK}(\bbL^*)$ est injective et inverse l'inclusion~:  supposons que $N_{\bbL_{2}/\bbK}(\bbL_{2}^*)\subset N_{\bbL_{1}/\bbK}(\bbL_{1}^*)$. La partie c) de la Proposition-DŽfinition \arggx\ nous montre que $\underbrace{\theta (\bbL_{1}/\bbK)\circ N_{\bbL_{2}/\bbK}}_{=1}=R\circ \theta(\bbL_{1}\bbL_{2}/\bbL_{2})$. Cela montre (puisque $R$ est injective) que l'application $\theta(\bbL_{1}\bbL_{2}/\bbL_{2})$ qui est surjective sur ${\rm Gal}(\bbL_{1}\bbL_{2}/\bbL_{2})$ est l'application triviale, cela montre que ${\rm Gal}(\bbL_{1}\bbL_{2}/\bbL_{2})$ est le groupe trivial, donc que $\bbL_{1}\bbL_{2}=\bbL_{2}$ et alors $\bbL_{1}\subset \bbL_{2}$.

Enfin, montrons la surjection de la correspondance~: soit $N$ un sous-groupe ouvert d'indice fini de $\bbK^*$. Soit $n\geq 1$ tel que $(\bbK^*)^n\subset N$ (par exemple $n=[\bbK^*:N]$). Soit $\mu_n$ l'ensemble des racines $n$-imes de l'unitŽs dans la cl™ture algŽbrique considŽrŽe de $\bbK$. Et soit enfin, $\bbL$ la $n$-extension de Kummer maximale de $K(\mu_n)$. L'extension $\bbL/\bbK$ est galoisienne. En effet, soit $\sigma\, :\, \bbL\to\overline{\bbL}$, un $\bbK$-plongement de $\bbL$ dans une cl™ture algŽbrique. Soit $\alpha\in\bbL$. En vertu de la Proposition \arggy, il existe $a\in \bbK(\mu_{n})$ tel que $\alpha^n-a=0$. Donc, $\sigma(\alpha)^n-\sigma(a)=0$. Puisque $\bbK(\mu_{n})/\bbK$ est galoisienne, $\sigma(a)\in\bbK(\mu_{n})$. Et puisque $\bbL$ contient toutes les racines $n$-imes d'ŽlŽment de $\bbK(\mu_{n})$ (toujours gr‰ce ˆ la Proposition \arggy), on en dŽduit que $\sigma(\alpha)\in\bbL$, ce qui prouve que $\bbL/\bbK$ est galoisienne. Et puisque $\bbK(\mu_n)/\bbK$ est cyclique, ${\rm Gal}(\bbL/\bbK)$ est rŽsoluble; on est donc dans les hypothses du thŽorme prŽcŽdent \argha. Ainsi, posons $N_0:=N_{\bbL^{\rm ab}/\bbK}({\bbL^{\rm ab}}^*)=N_{\bbL/\bbK}(\bbL^*)$. D'autre part, par la Proposition   \arggy, $N_{\bbL/\bbK(\mu_n)}(\bbL^*)=(\bbK(\mu_n)^*)^n$. Donc,
$$\eqalign{N_0&=N_{\bbL/\bbK}(\bbL^*)=N_{\bbK(\mu_n)/\bbK}(N_{\bbL/\bbK(\mu_n)}(\bbL^*))=N_{\bbK(\mu_n)/\bbK}((\bbK(\mu_n)^*)^n)\cr &=N_{\bbK(\mu_n)/\bbK}((\bbK(\mu_n)^*))^n\subset (\bbK^*)^n\subset N.\cr}$$
D'aprs ce qu'on a prouvŽ pour l'injection, l'application $\bbE\to N_{\bbE/\bbK}(\bbE^*)$ est une injection des sous-extensions $\bbE/\bbK$ de $\bbL^{\rm ab}/\bbK$ dans l'ensemble des sous-groupes de $\bbK^*$ qui contiennent $N_0$. Or, $\theta(\bbL^{\rm ab}/\bbK)$ induit un isomorphisme $\bbK^*/N_0\simeq {\rm Gal}(\bbL^{\rm ab}/\bbK)$. Donc, l'ensemble de ces sous-groupes est en bijection avec les sous-groupes de ${\rm Gal}(\bbL^{\rm ab}/\bbL)$, donc, par la thŽorie de Galois avec l'ensemble des sous-extensions de $\bbL^{\rm ab}/\bbK$. En particulier, il existe une sous-extension $\bbE/\bbK$ de $\bbL^{\rm ab}/\bbK$ telle que $N_{\bbE/\bbK}(\bbK^*)=N$. Cela dŽmontre le thŽorme.\qed

\centerline{\soustitre Voilˆ !!!}


\vfill\eject
\def\propa {{\soustitre Proposition ({\uppercase\expandafter{\romannumeral \the\chapnomb}}.{\romannumeral\the\nomb})}\global\advance\nomb by 1}
\def\tha {{\soustitre ThŽorme ({\uppercase\expandafter{\romannumeral \the\chapnomb}}.{\romannumeral\the\nomb})}\global\advance\nomb by 1}
\def\defia {{\soustitre DŽfinition ({\uppercase\expandafter{\romannumeral \the\chapnomb}}.{\romannumeral\the\nomb})}\global\advance\nomb by 1}
\def\defisa {{\soustitre DŽfinitions ({\uppercase\expandafter{\romannumeral \the\chapnomb}}.{\romannumeral\the\nomb})}\global\advance\nomb by 1}
\def\coroa {{\soustitre Corollaire ({\uppercase\expandafter{\romannumeral \the\chapnomb}}.{\romannumeral\the\nomb})}\global\advance\nomb by 1}
\def\lema {{\soustitre Lemme ({\uppercase\expandafter{\romannumeral \the\chapnomb}}.{\romannumeral\the\nomb})}\global\advance\nomb by 1}

\newcount\chapnomb \chapnomb=1
\nomb=1

\centerline{\para Appendice \uppercase\expandafter{\romannumeral 1} }
\bigskip
\centerline{\para Deux mots sur les corps quadratiques}
\smallskip
 \centerline{\para et sur des reprŽsentations de nombres premiers}
\bigskip

\defisa
\medskip
Soit $K$ un corps de nombres. Un {\it ordre sur $K$} est un sous-anneau $\O$ de $K$
qui est un $\Z$-module de gŽnŽration finie contenant une $\Q$-base de $K$. Un
raisonnement classique sur l'intŽgralitŽ ($\alpha\cdot \O\subset \O\Rightarrow \alpha $ est entier) montre que  tout ordre  est inclu dans $O_K$. C'est pourquoi $O_K$ est souvent appelŽ l'ordre maximal de $K$. Evidemment, si $\O$ est inclu strictement dans $O_K$, il n'est pas intŽgralement clos, donc
l'ensemble de ses idŽaux fractionnaires ne forme pas un groupe comme c'est le cas
pour $O_K$. En revanche, $\O$ est tout de mme ˆ quotients finis, il est donc noethŽrien, et tout
idŽal premier non nul est maximal. Clairement, si $\alpha\in O_K$, alors $\Z[\alpha]$ est un
ordre. 

Concentrons-nous sur les corps quadratiques. Rappelons les rŽsultats classiques sur
ces corps~: soit $m\in\Z\setminus\{1\}$ sans facteurs carrŽs et $K=\Q(\sqrt{m})$. Notons $d_K$ le discriminant de $K/\Q$. Alors on a~:

$$d_K=\cases{4m&si $m\equiv 2,3\pmod 4$\cr m&si $m\equiv 1\pmod
4$\cr}\qquad\hbox{et}\qquad
O_K=\Z\left[{d_K+\sqrt{d_K}\over 2}\right]=\cases{\Z[\sqrt{m}]&si $m\equiv 2,3\pmod
4$\cr
\Z\left[{1+\sqrt{m}\over 2}\right]&si
$m\equiv 1\pmod 4.$\cr}$$

Posons $w_K={d_k+\sqrt{d_K}\over 2}$. Si $\alpha, \beta\in K$ sont linŽairement
indŽpendants, le $\Z$-module $\Z\alpha\oplus \Z\beta$ se note $[\alpha,\beta]$. Il est
donc clair que $O_K=\Z[w_K]=[1,w_K]$.
\bigskip

\headline={\hfill \phantom{ouuh}\hfill}
\bigskip
{{\soustitre DŽfinition-Lemme ({\uppercase\expandafter{\romannumeral \the\chapnomb}}.{\romannumeral\the\nomb})}\global\advance\nomb by 1}

{\sl Pour tout entier $f\geq 1$, il existe un unique ordre $\O_f$ de $\Q(\sqrt{m})$
tel que $[O_K:\O_f]=f$ et de fait, $\O_f=\Z[fw_K]=[1,fw_K]$. On dit que $f$ est le
{\it conducteur} de l'ordre $\O_f$. Et tous les ordres sont donc de ce type

}
\newcount\gaagbc\gaagbc=\gaga

{\bf Preuve}

Puisque $(\Z\oplus \Z w_K)/(\Z\oplus \Z fw_K)\simeq \Z/ f\Z$, $\Z[fw_K]$ est un ordre
d'indice $f$ dans $O_K$.

Inversement, soit $\O$ un ordre d'indice $f$ dans $O_K$. Alors $\Z\subset \O$ (car
$\O$ est un sous-anneau de $K$) et
$f\cdot O_K\subset \O$. Ainsi, $\Z+f\cdot O_K=\Z\oplus \Z fw_k\subset \O$. Donc il y a ŽgalitŽ,
puisqu'ils ont le mme indice.\qed

\bigskip
{{\soustitre DŽfinition-ThŽorme ({\uppercase\expandafter{\romannumeral \the\chapnomb}}.{\romannumeral\the\nomb})}\global\advance\nomb by 1}
\medskip

{\sl Soit $m\in\Z\setminus\{1\}$ sans facteurs carrŽs et $K=\Q(\sqrt{m})$. Soit un entier $f>0$ et $\O=[1,fw_K]$ l'unique ordre d'indice $f$ de $K$. On pose $D:= f^2\cdot d_K$ le {\it discriminant }de $\O$. Il est Žvident que $D\equiv 0$ ou $1\pmod 4$ et que $\O=[1,\sqrt{D}]$ si $D\equiv 0\pmod 4$ et $\O=[1,{1+\sqrt{D}\over 2}]$ si $D\equiv 1\pmod 4$. RŽciproquement, si $D\equiv 0,1\pmod 4$, $D\ne 0,1$ n'est pas un carrŽ, (on dit alors que $D$ est {\it un discriminant}), alors il existe un unique entier $f>0$ et un unique corps quadratique $K$, tel que $D=f^2\cdot d_K$ (et donc, en vertu du lemme prŽcŽdent, un unique ordre $\O$ tel que $[O_K:\O]=f$. 

}
{\bf Preuve}

C'est ˆ peu prs Žvident~: soit $D$ un discriminant. Il existe un unique entier positif $g$ tel que $D=g^2\cdot m$, avec $m$ sans facteur carrŽ. Si $m\equiv 1\pmod 4$, alors $K=\Q(\sqrt{m})$ et si $m\equiv 2,3\pmod 4$, alors forcŽment $g$ est pair (par hypothse sur $D$), et dans ce cas, $K=\Q(\sqrt{m})$ et $f={g\over 2}$. L'unicitŽ vient du fait que $\Q(\sqrt{h^2\cdot m})=\Q(\sqrt{m})$, pour tout $h$ entier non nul.\qed
\bigskip\goodbreak
\defisa
\medskip
Soit $m\in\Z\setminus\{1\}$ sans facteurs carrŽs et $K=\Q(\sqrt{m})$. Soit un entier $f>0$ et $\O=[1,fw_K]$ l'unique ordre d'indice $f$ de $K$. On pose $H_{\O}=\{{\alpha\over \beta}\cdot O_K \mid\alpha,\beta\in \O,$ premiers ˆ $f\}$,  le sous-groupe de $I_K(f\cdot O_K)$ engendrŽ par les idŽaux principaux de la forme $\alpha\cdot O_K$, avec $\alpha\in \O$, premier ˆ $f$. Alors, il est Žvident que $P_{f\cdot O_K}\subset H_{\O}\subset I_K(f\cdot O_K)$. Si $K$ est rŽel, on pose $H_{\O}^+$ le mme sous-groupe avec la contrainte supplŽmentaire que $\alpha$ soit totalement positif (i.e $\sigma(\alpha)>0$ pour tout plongement de $K$ dans $\R$). Alors on a $P_\m\subset  H_{\O}^+\subset I_K(\m)$, o $\m=f O_K\cdot \{\sigma_1,\sigma_2\}$, avec $\sigma_1,\sigma_2$ les deux places infinies de $K$ (l'identitŽ et la conjugaison).

On appelle le {\it corps de classe de $\O$} le  corps de la classe de $H_{\O}$. Et on appelle  le {\it corps de classe Žtendu de $\O$} le  corps de la classe de $H_{\O}^+$. On peut construire ces corps de classe en vertu du thŽorme d'existence du corps de classe (ThŽorme \argdm).

Si $D$ est un discriminant, on dŽsigne par $Q_D$ la forme quadratique~:
$$Q_D(x,y)=\cases{x^2-{D\over 4} y^2&si $D\equiv 0\pmod 4$\cr x^2+xy+{1-D\over 4}y^2& si $D\equiv 1\pmod 4$.\cr}$$
\newcount\gaagbd\gaagbd=\gaga
Clairement, c'est une forme quadratique entire (i.e. une application du type $(x,y)\mapsto ax^2+bxy+cy^2,\ a,b,c\in\Z$), primitive (${\rm pgcd}(a,b,c)=1$), de discriminant $b^2-4ac=D$, dŽfinie positive si $D<0$.
\bigskip
\headline={\hfill \smcap Deux mots sur les corps quadratiques et sur des reprŽsentations de nombres premiers\hfill}
\tha
\medskip
{\sl Soit $D$ un discriminant, $K=\Q(\sqrt{D})$, et $\O$, l'ordre de $K$ de discriminant $D\ (=f^2\cdot d_K)$. Soit $L$ le corps de classe de $\O$ (Žtendu dans la cas rŽel ($D>0$)). Soit $p$ un nombre premier (en particulier positif), $p\notdiv D$. Alors on ˆ l'Žquivalence~:
$$\exists\, x,y\in \Z\hbox{ tels que }p=Q_D(x,y)\Longleftrightarrow \hbox{$p$ est compltement dŽcomposŽ dans $L$}.$$

}
{\bf Preuve}

On remarque dŽjˆ immŽdiatement que $Q_D(x,y)=N_{K/\Q}(x+y \cdot w_D)$, o $w_D=\cases{\sqrt{D}&si $D\equiv 0\pmod 4$\cr {1+\sqrt{D}\over 2}&si $D\equiv 1\pmod 4$,\cr}$ en vertu de la dŽfinition de la norme et de la DŽfinition-ThŽorme \argfi. Donc, $\exists\, x,y\in \Z\hbox{ tels que }p=Q_D(x,y)\iff \exists\, \alpha\in O$ tel que $p=N_{K/\Q}(\alpha)$. Ceci est Žquivalent ˆ $p=\alpha\cdot\alpha'$  o $\alpha'$ est le conjuguŽ d'$\alpha$; donc ils ont le mme signe si $K$ est rŽel. Donc, quitte ˆ remplacer $\alpha$ par $-\alpha$, on peut supposer que $\alpha$ est totalement positif si $K$ est rŽel. Ainsi, $\exists\, \alpha\in \O$ tel que $p=N_{K/\Q}(\alpha)\iff p\cdot O_K=\P\cdot\P'$, avec $\P=\alpha O_K\in \gfP_0(K)$ (car $p$ est un nombre premier), et $\P'\in \gfP_0(K)$ est le conjuguŽ de $\P$. Or $\P$ et $\P'$ sont premiers ˆ $f$, puisque $p$ l'est. Donc $\P$ et $\P'$ sont dans $H_{\O}$ (resp. dans $H_{\O}^+$ si $K$ est rŽel). Mais par dŽfinition, $H_{\O}$ (resp.  $H_{\O}^+$ si $K$ est rŽel) est le noyau de l'application d'Artin pour l'extension $L/K$, donc ils sont totalement dŽcomposŽs dans $L$. Finalement,
$$\eqalign{\exists\, x,y\in \Z\hbox{ tels que }p=Q_D(x,y)&\iff p\cdot O_K=\P\cdot\P'\hbox{ avec }\P,\P'\in H_{\O}\hbox{ (resp.  $H_{\O}^+$ si $K$ est rŽel)}\cr
&\iff p\cdot O_K=\P\cdot\P'\hbox{ avec $\P,\P'$ totalement dŽcomposŽs dans $L$}\cr
&\iff \hbox{$p$ est compltement dŽcomposŽ dans $L$.}\cr}$$
\qed
\bigskip
\coroa
\medskip
{\sl Si $D$ est un discriminant, alors l'ensemble des nombre premiers $p$ reprŽsentŽ par la forme $Q_D$ a une densitŽ de Dirichlet strictement positive (en particulier, il y en a une infinitŽ).

}
{\bf Preuve}

C'est Žvident, en vertu du thŽorme prŽcŽdent et du fait que l'ensemble des premiers qui dŽcompose compltement ˆ une densitŽ strictement positive (cf. Lemme \argao).\qed
\goodbreak
\bigskip
\coroa
\medskip
{\sl Si $D$ est un discriminant, alors il existe $f_D\in \Z[X]$, unitaire et irrŽductible tel que pour tout $p$ nombre premier tel que $p\notdiv D$ et $p\notdiv {\rm disc}(f_{D})$, on ait l'Žquivalence~:
$$\eqalign{\exists x,y\in \Z\hbox{ tels que } p=Q_D(x,y)&\iff f_D\hbox{ est totalement scindŽ }\bmod p\cr &\iff f_{D}(x)\equiv 0\pmod p\hbox{ a une solution}.\cr}$$

}

{\bf Preuve}

ConsidŽrons $K=\Q(\sqrt{D})$ et $L$ le corps de classe construit au thŽorme prŽcŽdent. Supposons que $L=\Q(\alpha)$, avec $\alpha\in O_L$ (c'est possible en vertu de [Sam,Corollaire, p. 41], pour un $\alpha\in K$, mais on peut le supposer dans $O_{K}$  gr‰ce au raisonnement du dŽbut de la preuve du Lemme \argfw ). On prend $f_D$ le polyn™me minimal de $\alpha$ sur $\Q$. Remarquons tout d'abord que l'extension $L/\Q$ est galoisienne. En effet, soit $\sigma\, :\, L\to\C$ un plongement. On regarde comme toujours $L\subset \C$. Il s'agit de prouver que $\sigma(L)=L$. Si $\sigma|_K={\rm id}$ alors $\sigma(L)=L$, puisque $L/K$ est galoisienne. Si $\sigma|_K\ne{\rm id}$, alors $\sigma(x)=x'$ o $x'$ est le conjuguŽ de $x$ dans $K$, pour tout $x\in K$ (car $\sigma(K)=K$, puisque $K/\Q$ est galoisienne). Le corps $\sigma(L)$ est une extension de $\sigma(K)=K$. Donc $\sigma(L)/K$ est une extension abŽlienne, celle qui correspond ˆ $\sigma(H_{\O_D})=\{\sigma(\alpha O_K)\mid\alpha O_K\in H_{\O_D}\}$. Mais $\sigma(\alpha O_K)=\alpha' O_K$. La dŽfinition de $H_{\O_D}$ montre que $\sigma(H_{\O_D})=H_{\O_D}$, car  $\alpha$ est premier ˆ $f$ si et seulement si $\alpha'$ l'est et $\alpha'$ est totalement positif si et seulement si $\alpha'$ l'est. Et finalement $\sigma(L)=L$ et donc $L/\Q$ est galoisienne. 

Terminons la preuve~: le thŽorme prŽcŽdent nous apprend que $p=Q_D(x,y)\iff p$ est totalement dŽcomposŽ dans $L$ c'est ˆ dire $f(\P/p)=1$ pour tout $\P\in\gfP_0(L)$ tel que $\P |p$. Cela est Žquivalent ˆ dire que $O_L/\P=\Z/p\Z=\bbF_p$ et ainsi toutes les racines de $f_D$ (qui sont dans $O_L$, car $L/\Q$ galoisienne)  modulo $\P$  sont dans $\bbF_p$ c'est-ˆ dire que $f_D$ est totalement scindŽ modulo $p$. 

Reste ˆ voir que si $f_{D}(x)\equiv 0\pmod p$ a une solution $\bmod\ p$, alors il est totalement scindŽ $\bmod\ p$. Deux moyens de voir ceci~: 

\art{1)} Si $p\notdiv [O_L:\Z[\alpha]]$ (ce qui est le cas ici, car $p\notdiv |{\rm disc}(f_{D})|= [O_L:\Z[\alpha]]\cdot{\rm disc}(O_L)$), on sait que les $f(\P/p)$ sont les degrŽs des facteurs irrŽductibles de la dŽcomposition de $f_D$ dans $\bbF_p$ (cf. [Mar, Thm. 27, p. 79]). Donc dire que  $f_{D}(x)\equiv 0\pmod p$ a une solution $\bmod\ p$, veut dire qu'il existe $\P|p$ tel que $f(\P/p)=1$, et alors $f(\P/p)=1$ pour tout premier $\P|p$, car dans une extension galoisienne, tous les $f(\P/p)$ sont Žgaux pour $p$ fixŽ. En traduction, $f_D$ est totalement dŽcomposŽ $\bmod\ p$.

\art{2)}Un autre manire de voir serait de se souvenir que $L\otimes \Q_p=\Q_p[x]/(f_D)=\prod_{\P|p}\bbL_\P$ (cf. [Fr-Tay, 1.6.a, p. 109]) et que si $f_D=f_1\cdots f_r$ est la dŽcomposition de $f_D$ dans $\Q_p$, on a $\Q_p[x]/(f_i)\simeq \bbL_{\P_i}$ pour un certain $\P_i |p$. Or, $[\bbL_\P:\Q_p]=f(\P/p)$ (puisque $p$ ne ramifie pas dans $L$ (cf. [Fr-Tay, 1.14, p. 111])). Donc, dire que $f$ a une racine $\bmod\ p$ veut dire que $f_i$ a une racine $\bmod\ \widehat{p}$ ($\widehat{p}=p\cdot\Z_{p}$). Donc par le Lemme de Hensel (qui s'applique ici, puisque $p\notdiv {\rm disc}(f_D)$(cf. [Fr-Tay, 3.25+1.9, pp. 84+10])) $f_i$ a une racine dans $\Q_p$, ce qui veut dire que $\bbL_{\P_i}=\Q_p$, et donc que $f(\P_i/p)=1$ et par suite, puisque $L/\Q$ est galoisienne, $f(\P/p)=1$ pour tout $\P|p$ ce qui veut dire que $p$ est totalement dŽcomposŽ sur $\bbF_p$ et donc $f_D$ est  totalement dŽcomposŽ dans $\F_p$, comme on vient de le voir.\qed
\bigskip
On a encore un petit raffinement si le discriminant est nŽgatif~:
\bigskip
\coroa
\medskip
{\sl Si $D$ est un discriminant strictement nŽgatif, alors il existe $g_D\in \Z[X]$, unitaire et irrŽductible de degrŽ moitiŽ de celui du corollaire prŽcŽdent tel que pour tout $p$ nombre premier tel que $p\notdiv D$ et $p\notdiv {\rm disc}(g_{D})$, on ait l'Žquivalence~:
$$\exists x,y\in \Z\hbox{ tels que } p=Q_D(x,y)\iff \left ({D\over p}\right )=1\hbox{ et }
g_D(x)\equiv 0\pmod p\hbox{ a une solution,}$$
o $\left ({D\over p}\right )$ est bien sžr le symbole de Legendre.

}
{\bf Preuve}

Soit $K=\Q(\sqrt{D})$ et $L$ comme pour les rŽsultats prŽcŽdents. Remarquons d'abord que l'extension $L/\Q$ est totalement complexe, car non rŽelle, ˆ cause de $K$, et on a vu au corollaire prŽcŽdent que $L/\Q$ Žtait galoisienne. Donc, la conjugaison complexe est un $\Q$ automorphisme de $L$. Soit $E$ le sous-corps de $L$ fixe par la conjugaison complexe. Alors $E\cap K=\Q$ et $[E:\Q]={1\over 2}\cdot [L:\Q]=[L:K]$. Donc $L=E\cdot K$, et si $E=\Q(\alpha)$, on a $L=K(\alpha)$ et ${\rm min}(\alpha/K)={\rm min}(\alpha/\Q)\in \Z[X]$, unitaire et irrŽductible. Prenons donc $g_D={\rm min}(\alpha/\Q)\in \Z[X]$ et $\alpha\in O_E\subset O_L$. On a vu au ThŽorme \argfk\ que $\exists\, x,y\in \Z\hbox{ tels que } p=Q_D(x,y)$ si et seulement si $p$ se dŽcompose totalement dans $L$. Or $p$ se dŽcompose totalement dans $K$ si et seulement si $\left ({D\over p}\right )=1$ (car dans ce cas, $x^2-D$ possde une solution modulo $p$ et on raisonne comme pour le corollaire prŽcŽdent). Et de mme, si $pO_K=\P\cdot\overline{\P}$, $\P$ se dŽcompose compltement dans $L$ si et seulement si $g_D(x)\equiv 0\pmod \P$  a une solution. Mais, comme $g_D\in\Z[X]$ et $O_K/\P\simeq \Z/p\Z$, cela revient ˆ dire que $g_D(x)\equiv 0\pmod p$  a une solution. On fait le mme raisonnement pour $\overline{\P}$ et cela montre le corollaire.\qed
\bigskip
\rem
\medskip
La ThŽorie de la multiplication complexe permet d'exhiber un $\alpha$ comme dans le corollaire prŽcŽdent (c'est-ˆ-dire $\alpha\in O_L$ rŽel tel que $L=K(\alpha)$) comme suit~: soit $\Lambda$ un rŽseau dans $\C$ (c'est-ˆ-dire un sous-$\Z$-module de $\C$ engendrŽ par des ŽlŽment de $\R$). On dŽfinit $g_2(\Lambda)=60\cdot\dst\sum_{0\ne\lambda\in\Lambda}{1\over \lambda^4}$, $g_3(\Lambda)=140\cdot\dst\sum_{0\ne\lambda\in\Lambda}{1\over \lambda^6}$, $\Delta(\Lambda)=g_2(\Lambda)^3-27 g_3(\Lambda)^2$ qu'on prouve tre le discriminant de $\Lambda$, donc $\ne 0$ et $j(\Lambda)=1728\cdot{g_2(\Lambda)^3\over\Delta(\Lambda)}$. Si $\O$ est un ordre dans le corps quadratique imaginaire $K$, on peut voir $\O$ comme un rŽseau dans $\C$ et on peut montrer que $j(\O)$ est un entier algŽbrique tel que $L=K(j(\O))$. 

Par exemple, si $D=-56$, et donc $\O=\Z\oplus\Z[\sqrt{-14}]=O_{\Q(\sqrt{-14})}$, alors on peut calculer que 
$$j(\O)=2^3\left (323+228\sqrt{2}+(231+161\sqrt{2})\sqrt{2\sqrt{2}-1}\right )^3,$$
Òce qui implique que" $L=K(\sqrt{2\sqrt{2}-1} )$. On montre facilement que le polyn™me minimal de $\sqrt{2\sqrt{2}-1}$ est $x^4+2x^2-7=(x^2+1)^2- 8$. Ce qui montre que la densitŽ des nombres premiers impairs  $p\ne 7$ tels que $p=x^2+14y^2$ est strictement positive et 
$$p=x^2+14y^2\iff \left({-14\over p}\right)=1\ \hbox{ et  }(x^2+1)^2\equiv 8\pmod p\ \hbox{ a une solution.}$$
Pour plus de dŽtails, voir le livre [Cox].

\vfill\eject

\newcount\chapnomb \chapnomb=2
\nomb=1

\centerline{\para Appendice \uppercase\expandafter{\romannumeral 2} }
\bigskip
\centerline{\para Deux mots sur le symbole de Hilbert}
\bigskip
\defia
\medskip
Soit $n\in\N$, $K$ un corps de nombres contenant une racine primitive $n$-ime de l'unitŽ. Soit $\P$ une place de $K$ et $\bbK_{\P}$ le complŽtŽ localisŽ de $K$ en $\P$. On considre encore $\bbL=\bbK_{\P}(\root n\of{\bbK_{\P}^*})$ qui est la $n$-extension maximale de Kummer de $\bbK_{\P}$. On sait de plus que $\ker(\theta(\bbL/\bbK_{\P}))=(\bbK_{\P}^*)^n$ (cf. Propositions \arggx \ et \arggy). Notons $G={\rm Gal}(\bbL/\bbK_{\P})$. On a donc un isomorphisme $\overline{\theta}(\bbL/\bbK_{\P})\, :\,\bbK_{\P}^*/(\bbK_{\P}^*)^n\simeq G$. Notons $\mu_{n}$ l'ensemble des racines $n$-imes de l'unitŽ de $\bbK_{\P}^*$. Notons comme toujours $\widehat{G}$ l'ensemble des homomorphismes de $G$ sur $\mu_{n}$. Puisque l'extension $\bbL/\bbK_{\P}$ est une $n$-extension de Kummer, nous savons (cf. \argdu, $M=K^*$) que l'application $\chi\, :\, \bbK_{\P}^*\to \widehat{G}$, $b\mapsto (\sigma\mapsto {\sigma(\root n\of b)\over \root n\of b})$ est un homomorphisme de noyau $(\bbK^*)^n$. Ainsi on a un isomorphisme $\overline{\chi}\, :\,\bbK_{\P}^*/(\bbK_{\P}^*)^n\simeq \widehat{G}$. Puis par composition avec le couplage $G\times \widehat{G}\to \mu_{n}$, $(\sigma,\chi)\mapsto \chi(\sigma)$, on peut (enfin) dŽfinir le couplage non dŽgŽnŽrŽ suivant~:
$$\eqalign{\bbK_{\P}^*/(\bbK_{\P}^*)^n\times \bbK_{\P}^*/(\bbK_{\P}^*)^n&\lra \mu_{n}\cr (\overline{a},\overline{b})&\longmapsto \left ({a,b\over\P}\right )_{n}:=\overline{\chi}(\overline{b})(\overline{\theta}(\bbL/\bbK_{\P})(\overline{a}))\cr}$$
qu'on appelle {\it $n$-ime symbole de Hilbert}. Par restriction ˆ $\bbK_{\P}(\root n\of b)$ (Proposition-DŽfinition \arggx, d)), le symbole $\left ({a,b\over\P}\right )_{n}$, peut tre dŽfini par la relation
$$\theta(\bbK_{\P}(\root n\of b)/\bbK_{\P})(a)(\root n\of b)=\left ({a,b\over\P}\right )_{n}\cdot (\root n\of b).\eqno{(*)}$$
Les esprits observateurs remarquerons que cette relation ressemble dr™lement ˆ la relation $(*)$ de la preuve du ThŽorme \argec. Ce n'est pas un hasard~:

\headline={\hfill \phantom{ouuh}\hfill}
\newcount\gaagbe\gaagbe=\gaga

\bigskip
\propa
\medskip
{\sl Soit $n$, $K$, $\P$ comme pour la dŽfinition prŽcŽdente. Posons $\pi$ une uniformisante de $\bbK_{\P}$ et soit $\alpha\in \bbK_{\P}$. Supposons de plus que $\P\notdiv n\cdot\alpha$. Alors
$$\left ({\pi,\alpha\over\P}\right )_{n}=\left ({\alpha\over\P}\right )^{-1}_{n},$$
o $\left ({\alpha\over\P}\right )_{n}$ est le symbole de puissance $n$-ime rŽsiduelle dŽfini lors de la DŽfinition \argea, et Žtendu de manire naturelle sur $\okp$ ($=O_{K_{\P}}$) puisque $\okp/\widehat{\P}\simeq O_{K}/\P$ (cf. [Fr-Tay, Theorem 11, p. 77]). 

}
{\bf Preuve}

Puisque $\widehat{\P}\notdiv n\cdot\alpha$, alors $\widehat{\P}$ ne ramifie pas dans $\bbK_{\P}(\root n\of b)$ (cf. [Nar, Thm. 4.17, \S 3,  chap IV, p. 184]). Dans ce cas, la Proposition \3argel\ s'applique, et donc  (en adaptant un petit peu la preuve) on obtient~: $\theta(\bbK_{\P}(\root n\of b)/\bbK_{\P})(\pi)={\rm Frob}_{\bbK_{\P}(\root n\of b)/\bbK_\P}(\P)^{-1}$. Et on conclut en vertu de la relation $(*)$ de la dŽfinition prŽcŽdente et de la relation $(*)$ du ThŽorme \argec.\qed
\bigskip\goodbreak

\propa {\soustitre (RŽciprocitŽ de Hilbert)}
\medskip
{\sl Soit $K$ un corps de nombres contenant une racine $n$-ime primitive de l'unitŽ, et $a,b\in K^*$. Alors
$$\prod_{\P\in\gfP(K)}\left ({a,b\over\P}\right )_{n}=1.$$

}
{\bf Preuve}

La pertinence de ce produit et la preuve de la proposition est une application directe du Lemme de rŽciprocitŽ pour les symboles des restes normiques \arggr.\qed 
\bigskip

\tha
\medskip
{\sl Soit $n\in\N$, $K$ un corps de nombres contenant une racine primitive $n$-ime de l'unitŽ, $\P$ une place de $K$ et $\bbK_{\P}$ le complŽtŽ localisŽ de $K$ en $\P$. Alors l'application $\left ({\cdot,\cdot\over\P}\right )_{n}\, :\, \bbK_{\P}^*\times  \bbK_{\P}^*\to\mu_{n}$ a les propriŽtŽs suivantes~($a,b,a',b'\in\bbK_{\P}^*$)~:
\art{a)}$\left ({aa',b\over\P}\right )_{n}=\left ({a,b\over\P}\right )_{n}\cdot \left ({a',b\over\P}\right )_{n}$

\art{b)}$\left ({a,bb'\over\P}\right )_{n}=\left ({a,b\over\P}\right )_{n}\cdot \left ({a,b'\over\P}\right )_{n}$

\art{c)}$\left ({a,b\over\P}\right )_{n}=1\iff a\in N_{\bbK_\P(\root n\of b)/\bbK_\P}(\bbK_\P(\root n\of b)^*)$.

\art{d)}$\left ({a,b\over\P}\right )_{n}=1 \ \forall\,  b \iff a\in (\bbK_{\P}^*)^n$

\art{e)}$\left ({a,b\over\P}\right )_{n}\cdot \left ({b,a\over\P}\right )_{n}=1$

\art{f)}$\left ({a,1-a\over\P}\right )_{n}=1$ si $a\ne 1$

\art{g)}$\left ({a,-a\over\P}\right )_{n}=1$

Enfin, si $n=2$, $\left ({a,b\over\P}\right )_{2}\in\{\pm 1\}$. Donc $\left ({a,b\over\P}\right )_{2}=\left ({b,a\over\P}\right )_{2}$ en vertu de e).  Et donc $\left ({\cdot,\cdot\over\P}\right )_{2}\, :\, \bbK_{\P}^*\times  \bbK_{\P}^*\to\{\pm 1\}$ est une $\bbF_{2}$-forme bilinŽaire symŽtrique non dŽgŽnŽrŽe. Mais on a en plus,
$$\left ({a,b\over\P}\right )_{2}=\cases{1&si $ax^2+by^2-z^2=0$ a une solution non triviale en $(x,y,z)\in {\bbK_{\P}^*}^3$\cr -1& sinon.\cr}$$

}
\medskip
{\bf Preuve}

Les parties  a)-d)  sont immŽdiates par dŽfinition du symbole de Hilbert. On dŽmontre les parties f) et g) ensemble~: soit $\xi\in \bbK_{\P}^*$ tel que $\xi^n-b\ne 0$. Alors on a $\xi^n-b=\prod_{i=0}^{n-1}(\xi-\zeta^{i}\root n\of b)$, pour $\root n\of b$ une racine $n$-ime de $b$ fixŽe et $\zeta$ une racine primitive $n$-ime de l'unitŽ fixŽe aussi. Soit $d$ le plus grand diviseur de $n$ pour lequel $\bbK_\P^*$ possde une racine $d$-ime de $b$. Alors $\bbK_{\P}(\root n\of b)/\bbK_\P$ est une extension de degrŽ ${n\over d}:=m$. Les conjuguŽs de $\xi-\zeta^{i}\root n\of b$ sont les $\xi-\zeta^{j}\root n\of b$ avec $j\equiv i\pmod d$. En effet, $\root n\of b=\root m\of{\root d\of b}:=\root m\of c$ et $x^m-c\in \bbK_\P[x]$ irrŽductible. Et si $\sigma\in {\rm Gal}(\bbK_\P(\root n\of b)/\bbK_\P)$,  $\sigma(\zeta^{i}\cdot\root n\of b)=\zeta^i\cdot\zeta_m^l\cdot\root m\of c=\zeta^{i}\cdot \zeta^{dl}\cdot\root m\of c=\zeta^{i+dl}\cdot \root n\of b$. Ainsi, 
$$\xi^n-b=\prod_{i=0}^{d-1}N_{\bbK_\P(\root n\of b)/\bbK_\P}(\xi-\zeta^{i}\cdot\root n\of b)\in  N_{\bbK_\P(\root n\of b)/\bbK_\P}(\bbK_\P(\root n\of b)^*),$$
Donc, en vertu de la partie c), on a $\left ({\xi^n-b,b\over\P}\right )_{n}=1$ Prenant $\xi=1, b=1-a$, puis $\xi=0$ et $b=-a$, on obtient f) et g). La partie e) rŽsulte du calcul~:
$$\left ({a,b\over\P}\right )_{n}\cdot \left ({b,a\over\P}\right )_{n}=\left ({a,-a\over\P}\right )_{n}\cdot \left ({a,b\over\P}\right )_{n}\cdot\left ({b,a\over\P}\right )_{n}\cdot \left ({b,-b\over\P}\right )_{n}=\left ({a,-ab\over\P}\right )_{n}\cdot \left ({b,-ab\over\P}\right )_{n}=\left ({ab,-ab\over\P}\right )_{n}=1$$
Enfin, si $n=2$, et que $\left ({a,b\over\P}\right )_{2}=1$, alors, la partie c) nous apprend que $a$ est une norme de $\bbK_{\P}(\root n\of b)/\bbK_{\P}$, donc $a=z^2-by^2$ pour certains $x,y\in\bbK_{\P}$. RŽciproquement, si, $ax^2+by^2-z^2=0$ a une solution non triviale $(x_{1},y_{1},z_{1})$, alors, si $x_{1}\ne 0$, alors en divisant par $x_{1}$, on voit que $a$ est une norme de $\bbK_{\P}(\root n\of b)/\bbK_{\P}$ et si $y_{1}\ne 0$, alors en divisant par $y_{1}$, on voit que $b$ est une norme de $\bbK_{\P}(\root n\of b)/\bbK_{\P}$, donc ˆ chaque fois $\left ({a,b\over\P}\right )_{2}=1$. Cela montre le thŽorme.
\qed
\bigskip
Si $n=2$ et $K=\Q$, on peut faire des calculs plus explicites pour obtenir les mmes rŽsultats (cf. [Ser, Chapitre III]).

\headline={\hfill \smcap Deux mots sur le symbole de Hilbert\hfill}


\vfill\eject

\centerline{\para Glossaire et symboles (dans l'ordre d'apparition )}
\bigskip
\headline={\hfill \smcap \phantom{ouuh}\hfill}
\settabs 3 \columns 
\+$O_{K}$ &l'anneau des entiers d'un corps de nombres $K$ & \hskip3cm  page\ \the\gaaga    \cr 
\+$U_{K}$&unitŽs de $O_{K}$& \hskip3cm page\ \the\gaaga    \cr 
\+ $f(\gP/\P)$ &degrŽ rŽsiduel de $\gP|\P$\  & \hskip3cm page\ \the\gaaga   \cr
\+  $e(\gP/\P)$ &indice de ramification de $\gP|\P$\  & \hskip3cm page\ \the\gaaga   \cr
\+${\rm Gal}(L/K)$&groupe de Galois de $L/K$&\hskip3cm page\ \the\gaagb\cr
\+ $\N(\aa), N_{L/K}$ &normes absolues et relatives&\hskip3cm page \the\gaagc \cr 
\+Z(\gP/\P)&groupe de dŽcomposition&\hskip3cm page \the\gaagd \cr 
\+T(\gP/\P)&groupe d'inertie&\hskip3cm page \the\gaage \cr
\+${\rm Frob}(\gP/\P)$&automorphisme de Frobenius&\hskip3cm page \the\gaage \cr
\+${\rm Fr}_{L/K}(\P)$&classe de conjugaison des ${\rm Frob}(\gP/\P)$&\hskip3cm page \the\gaage \cr
\+$v_{\P}(x), |x|_{\P}$&valuation et valeur absolue $\P$-adique &\hskip3cm page \the\gaagf \cr
\+$\gfP_{0}(K)$, $\gfP_{\infty}(K)$, $\gfP_{\R}(K)$,&& \cr 
\+$\gfP_{\C}(K)$, $\gfP(K)$&ensembles de places&\hskip3cm page \the\gaagf \cr
\+$\bbL_{\gP}$, $\bbK_{\P}$&corps $\gP$, $\P$-adiques&\hskip3cm page \the\gaagf \cr
\+ $O_{\P}=O_{\bbK_{\P}}$ et $O_{(\P)}=O_{\P}\cap K$&complŽtŽ localisŽ et localisŽ de $O_K$ en $\P$&\hskip3cm page \the\gaagf \cr
\+ $\widehat{\P}$ et $\widetilde{\P}$&idŽaux maximaux de $O_{\P}$ et $O_{(\P)}$ &\hskip3cm page \the\gaagf \cr
\+\m&$K$-module&\hskip3cm page \the\gaagg \cr
\+${\euf 1}$&le $K$-module unitŽ &\hskip3cm page \the\gaagg \cr
\+$S(\m),S_{0}(\m),S_{\infty}(\m)$&places divisant $\m$ &\hskip3cm page \the\gaagg \cr
\+$K^*_{\m}$&$K$ Žtoile $\m$&\hskip3cm page \the\gaagg \cr
\+$x\equiv y\pmodast\m$&congruence Žtoile&\hskip3cm page \the\gaagg \cr
\+$I_{K},I_{K}^S,I_{K}(\m)$&groupes d'idŽaux&\hskip3cm page \the\gaagh \cr
\+$\Phi_{L/K}$&l'application d'Artin&\hskip3cm page \the\gaagh \cr
\+$\widetilde{S}, \widetilde{\m}, I_{L}(\widetilde{\m})$&extensions ˆ $L$ de $S$, $\m$ et $I_{K}(\m)$&\hskip3cm page \the\gaagi \cr
\+$P_{K},P(\m)$&groupes d'idŽaux principaux&\hskip3cm page \the\gaagj \cr
\+$K^*(\m)$&$\cap_{\P\in S_0(\m)}O^*_{(\P)}$&\hskip3cm page \the\gaagj \cr
\+$P_{\m}$&groupe d'idŽaux principaux&\hskip3cm page \the\gaagj \cr
\+$h_{\m}$&cardinal $\left |I(\m)/P_\m\right |$ &\hskip3cm page \the\gaagj \cr
\+$U_{\m}$&$U_{K}\cap K^*_{\m}$&\hskip3cm page \the\gaagj \cr
\+$j(x,{\euf K})$&idŽaux de $\euf K$ de norme inf. ˆ $x$&\hskip3cm page \the\gaagk \cr
\+$v$ et $l$ & plongements canonique et log.& \hskip3cm page \the\gaagk \cr
\+$l_{0}$ et $N_{0}$ & applications&\hskip3cm page \the\gaagk \cr
\+$\partial A, \buildrel\hskip 3pt \circ\over A,\overline{A}$&bord, intŽrieur, adhŽrence de $A$& \hskip3cm page \the\gaagl \cr
\+$\Re(z),\Im(z)$&partie rŽelle et imaginaire de $z$ &\hskip3cm page \the\gaagl \cr
\+$\zeta(s)$&la fonction zeta de Riemann&\hskip3cm page \the\gaagm \cr
\+$\zeta_\m(s,{\euf K})$& fonction zeta de $\m$ et $\euf K$&\hskip3cm page \the\gaagn \cr
\+$\chi(g)\in\widehat{G}$& caractre d'un groupe abŽlien $G$&\hskip3cm page \the\gaagn \cr
\+${\bf 1}$&caractre unitŽ&\hskip3cm page \the\gaagn \cr
\+$L_\m(s,\chi)$& fonction $L$ pour $\m$ et $\chi$&\hskip3cm page \the\gaago \cr
\+$\zeta_\m(s)$& $= L_\m(s,{\bf 1})$&\hskip3cm page \the\gaago \cr
\+${\rm Log}(z)$& branche principale du logarithme&\hskip3cm page \the\gaago \cr
\+$\llog  L_\m (s,\chi)$& fonction logarithme de $ L_\m (s,\chi)$&\hskip3cm page \the\gaagp \cr
\+$f \sim g$& Žquivalence de $f$ et $g$ &\hskip3cm page \the\gaagp \cr
\+$\delta(S)$& densitŽ de Dirichlet de $S$ &\hskip3cm page \the\gaagq \cr
\+$T(m,n)$& nombre d'ŽlŽments de $C_m$ d'ordre multiple de $n$  &\hskip3cm page \the\gaagr \cr
\+$S(L/K),\widetilde{S}(L/K)$& idŽaux de $K$ qui se dŽcomposent &\hskip3cm page \the\gaags \cr
\+$\Delta,N,\Delta|A,N|A$& $1-\sigma$, norme-trace  &\hskip3cm page \the\gaagt \cr
\+$H^0(A),H^1(A)$& groupes de cohomologie  &\hskip3cm page \the\gaagt \cr
\+$q(A)$& quotient de Herbrand  &\hskip3cm page \the\gaagu \cr
\+$f_\P,\, e_\P$& degrŽs rŽsiduels, indices de ramification&\hskip3cm page \the\gaagv \cr
\+$\iota,\, j_\m,\, f_\m$& homomorphismes de $G$-modules &\hskip3cm page \the\gaagw \cr
\+${L^*}^S$& les $S$-unitŽs de $L$&\hskip3cm page \the\gaagw \cr
\+$a(\m)$& l'indice $[K^*:N(L^*)K_\m^*]$ &\hskip3cm page \the\gaagx \cr 
\+$U_\P, U_\P^{(k)}$& $U_\P=$ les inversibles de $O_\P$, $U_\P^{(k)}=1+\widehat{\P}^k\subset U_\P$ &\hskip3cm page \the\gaagy \cr 
\+$\exp(x),\ \log(x+1)$& dŽfinition formelle des fonctions  $\exp$ et $\log$ &\hskip3cm page \the\gaagz
 \cr 
\+ $n(\m)$ &$[K^*_\m\cap\iota^{-1}(N(I_L(\widetilde{\m}))): K^*_\m\cap N(L^*)]$&\hskip3cm page \the\gaagaa\cr
\+ $\m$ admissible&$\m$ est admissible si $P_\m\subset \ker(\Phi_{L/K}|I_K(\m))$ &\hskip3cm page \the\gaagac\cr
\+ $H$ sous-groupe de congruence & $H$ est ainsi si $P_\m\subset H\subset I_K(\m)$ pour un $\m$ &\hskip3cm page \the\gaagad\cr
\+$\bbH$&  classe d'Žquivalence de sous-groupe de congruence&\hskip3cm page \the\gaagae\cr
\+$\ff$& conducteur d'une classe $\bbH$&\hskip3cm page \the\gaagae\cr
\+$H(\m)\subset\bbH$ &groupe de congruence dŽfini pour $\m$ dans  $\bbH$&\hskip3cm page \the\gaagae\cr
\headline={\hfill \smcap Glossaire et symboles\hfill}
\+$\bbH(L/K)$&classe d'Žquivalence dŽf. par les noyaux de $\Phi_{L/K}$&\hskip3cm page \the\gaagaf\cr
\+$\ff(L/K)$&conducteur de la classe $\bbH(L/K)$ (ou de $L/K$)&\hskip3cm page \the\gaagaf\cr
\+$I_K/\bbH$&$I_K(\m)/H(\m)$ pour n'importe quel $H(\m)$ de $\bbH$&\hskip3cm page \the\gaagag\cr
\+$H_E(\widetilde{\m})\in \bbH_E$& classe de $E$ dŽfinie par une classe $\bbH$ de $K$&\hskip3cm page \the\gaagah\cr
\+$L/K $&$n$-extension de Kummer &\hskip3cm page \the\gaagai\cr
\+${K^*}^S$&$\{a\in K^*\mid v_\P(aO_K)\ne 0\Rightarrow \P\in S\}$, les $S$-unitŽs de $K$ &\hskip3cm page \the\gaagaj\cr
\+$I_K[S]$&le sous-groupe de $I_K$ engendrŽ par les  $\P\in S$&\hskip3cm page \the\gaagaj\cr
\+$c(\m)$&$[K^*:{K^*}^nK^*_\m]$&\hskip3cm page \the\gaagak\cr
\+$\left ({\alpha\over \P}\right )_n$& symbole de puissance $n$-ime rŽsiduelle&\hskip3cm page \the\gaagal\cr
\+$\Phi_{E/K}$& application d'Artin pour $E/K$ non abŽlienne&\hskip3cm page \the\gaagam\cr
\+$\Theta$& application de $K_\m\to {\rm Gal}(L/K)$&\hskip3cm page \the\gaagan\cr
\+$U_\P^{(b)}$& extension de la dŽf. de $U_\P^{(b)}$ vu p. \the\gaagy&\hskip3cm page \the\gaagao\cr
\+$\theta_\P$& symbole de norme rŽsiduelle&\hskip3cm page \the\gaagao\cr
\+$V_{G\to H}$& hom. de transfert de $G$ sur $H$&\hskip3cm page \the\gaagap\cr
\+$I_G$& idŽal d'augmentation&\hskip3cm page \the\gaagaq\cr
\+$d(g)$& $d(g):=g-1$&\hskip3cm page \the\gaagaq\cr
\+$|x|_{\P}$& Ònorme $\P$-adique"&\hskip3cm page \the\gaagar\cr
\+ppt& Òpour presque tout"&\hskip3cm page \the\gaagas\cr
\+$\bbA_{K},\bbI_K$& adles et idles de $K$&\hskip3cm page \the\gaagas\cr
\+$\bbA_{K}(S)$& $\prod_{\P\in S}\bbK_{\P}\times \prod_{\P\not\in S}O_{\P}$&\hskip3cm page \the\gaagas\cr
\+$\bbI_{K}(S)$& $\prod_{\P\in S}\bbK_{\P}^*\times \prod_{\P\not\in S}O_{\P}^*$&\hskip3cm page \the\gaagas\cr
\+$K,K^*$& adles et idles principaux&\hskip3cm page \the\gaagat\cr
\+$C_K$& $\bbI_{K}/K^*$&\hskip3cm page \the\gaagau\cr
\+$|a|$&$\prod_{\P\in\gfP(K)}|a|_\P$ volume d'un idle &\hskip3cm page \the\gaagau\cr
\+$\bbI_K^0$& noyau de l'application volume, idles spŽciaux &\hskip3cm page \the\gaagav\cr
\+$C_K^0$& $\bbI_K^0/K^*$ &\hskip3cm page \the\gaagav\cr
\+$\bbI_{\m}$& $\prod_{\P\in\gfP(K)}U_{\P}^{(m_{\P})}$, o $\m=\prod_{\P\in\gfP(K)}\P^{m_\P}$ &\hskip3cm page \the\gaagaw\cr
\+$\bbI'_{\m}$& $\{(a_{\P})_{\P}\in\bbI_{K}\mid a_{\P}\in U_{\P}^{(m_{\P})}\ \forall \P |\m\}$, o $\m=\prod_{\P\in\gfP(K)}\P^{m_\P}$ &\hskip3cm page \the\gaagaw\cr
\+$C_\m$& $(\bbI_\m\cdot K^*)/K^*$ &\hskip3cm page \the\gaagaw\cr
\+$\psi $& $(a_{\P})_{\P}\mapsto\prod_{\P\in \gfP_{0}(K)}\P^{v_{\P}(a_{\P})}$ &\hskip3cm page \the\gaagax\cr
\+$\psi_\m$& $\psi |_{\bbI'_\m} $ &\hskip3cm page \the\gaagax\cr
 \+$H(\m)$& $\psi(H\cap \bbI'_\m)$ si $I_\m\subset H\subset \bbI_K$ est un sous-groupe ouvert &\hskip3cm page \the\gaagax\cr
\+$N_{L/K}$& Extension de la dŽfinition de norme de $\bbI_L$ sur $\bbI_K$ &\hskip3cm page \the\gaagay\cr
\+ $\bbI''_{L}(\m)$& $\{(a_{\gP})_{\gP}\in \bbI_{L}\mid a_{\gP}=1\ \forall \gP\ \hbox{ tel que } \gP |\widetilde{\m}\}$ &\hskip3cm page \the\gaagaz\cr 
\+$\Upsilon$& $\bbI_{K}\to {\rm Gal}(L/K)$, $ (a_{\P})_{\P}\mapsto \prod_{\P\in \gfP(K)}\theta_{\P}(a_{\P})$ &\hskip3cm page \the\gaagba\cr 
\+$\theta(\bbL/\bbK) $& Symbole de norme rŽsiduel pour $\bbL/\bbK$ &\hskip3cm page \the\gaagbb\cr
\+$\O_f$& $\Z[fw_K]=[1,fw_K] $, ordre de d'indice $f$ &\hskip3cm page \the\gaagbc\cr
\+$Q_D(x,y)$& forme quadratique entire liŽe au discriminant $D$ &\hskip3cm page \the\gaagbd\cr
\+$ \left ({a,b\over\P}\right )_{n}$& symbole de Hilbert &\hskip3cm page \the\gaagbe\cr
 
\vfill\eject

\centerline{\para Index}
\bigskip
\headline={\hfill \smcap \phantom{ouuh}\hfill}

adles	\qquad	136	\par
anneau de Dedekind	\qquad	2	\par
anneau des entiers d'un corps de nombres	\qquad	2	\par
application d'Artin	\qquad	11	\par
application d'Artin pour des extensions non abŽliennes	\qquad	108	\par
automorphisme de Frobenius	\qquad	3	\par
branche principale du logarithme	\qquad	31	\par
caractres	\qquad	30	\par
classes d'Žquivalence de sous-groupes de congruences	\qquad	85	\par
conducteur de $L/K$	\qquad	87	\par
conducteur d'un ordre	\qquad	161	\par
conducteur d'un sous-groupe ouvert	\qquad	145	\par
conducteur d'une classes d'Žquivalence de sous-groupes de congruences	\qquad	87	\par
corps cyclotomiques	\qquad	2	\par
corps de classe	\qquad	89	\par
corps de classe d'un ordre	\qquad	162	\par
corps de classe Žtendu d'un ordre	\qquad	162	\par
corps de Hilbert	\qquad	123	\par
corps de nombres	\qquad	1	\par
corps local	\qquad	152	\par
dŽcomposition d'idŽaux premiers	\qquad	48	\par
degrŽ rŽsiduel	\qquad	2	\par
densitŽ de Dirichlet	\qquad	33	\par
discriminant d'un corps	\qquad	21	\par
discriminant d'un corps quadratique	\qquad	101+161	\par
discriminant d'un ordre	\qquad	161	\par
discriminant d'un polyn™me	\qquad	46	\par
ŽgalitŽ fondamentale du corps de classe (ext. cycliques)	\qquad	74	\par
ensemble $J$-mesurable	\qquad	22	\par
ensemble rŽgulier	\qquad	33	\par
entiers indŽpendants modulo $m$	\qquad	79	\par
exponentiel formel	\qquad	64	\par
extension cyclotomique de corps de nombres	\qquad	16	\par
extension de Kummer	\qquad	94	\par
fonction lipschitzienne	\qquad	22	\par
fonctions zeta	\qquad	27+29	\par
$G$-modules	\qquad	50	\par
groupe de classes radiales	\qquad	13	\par
groupe de dŽcomposition	\qquad	4	\par
groupe des classes d'adles	\qquad	139	\par
groupe des classes d'idŽaux	\qquad	13	\par
groupe des classes d'idles	\qquad	139	\par
groupe d'inertie	\qquad	4	\par
homomorphisme de transfert	\qquad	126	\par
idŽal d'augmentation	\qquad	127	\par
idles	\qquad	136	\par
idles spŽciaux	\qquad	141	\par
indice de ramification	\qquad	2	\par
$K$-modules	\qquad	7	\par
$K$-modules $\P$-admissibles	\qquad	112	\par
$K$-modules admissibles	\qquad	77	\par
$K$-modules admissibles pour un sous-groupe ouvert	\qquad	145	\par
lemme d'Artin	\qquad	80	\par
lemme de continuitŽ des racines	\qquad	153	\par
lemme de Krasner	\qquad	152	\par
lemme de l'hexagone	\qquad	51	\par
lemme de naturalitŽ (deuxime)	\qquad	119	\par
lemme de naturalitŽ (premier)	\qquad	119	\par
lemme de naturalitŽ (troisime)	\qquad	121	\par
lemme de translation	\qquad	111	\par
logarithme formel	\qquad	64	\par
module de permutation	\qquad	54	\par
normes, absolues, relatives, d'un idŽal	\qquad	3	\par
ordre sur $K$	\qquad	161	\par
partie $(n-1)$-Lipschitz paramŽtrisable	\qquad	22	\par
pavŽ	\qquad	22	\par
place (finie, infinie rŽelle ou complexe)	\qquad	5	\par
pour presque tout	\qquad	136	\par
quotient  de Herbrand	\qquad	53 et suiv.	\par
ramification des places infinies	\qquad	55+57	\par
rŽduction d'un polyn™me modulo $p$	\qquad	45	\par
rŽseau plein	\qquad	21	\par
sŽrie $L$ de Dirichlet	\qquad	30	\par
sŽrie de Dirichlet	\qquad	26	\par
sous-groupe de congruence	\qquad	85	\par
sous-groupe de congruence Žquivalents	\qquad	85	\par
sous-groupe des classes d'idles spŽciaux	\qquad	141	\par
sous-groupes des commutateurs	\qquad	108+1157	\par
S-unitŽs	\qquad	56	\par
symbole de Hilbert	\qquad	166	\par
symbole de Legendre	\qquad	100+164	\par
symbole de puissance $n$-ime rŽsiduelle	\qquad	99	\par
symbole de reste normique	\qquad	116	\par
symbole de reste normique pour les corps locaux	\qquad	156	\par
thŽorme $[\bbK^*_\P:N_{\bbL_\gP/\bbK_{\P}}(\bbL^*_\gP)]=[\bbL_\gP:\bbK_\P]$	\qquad	69+122	\par
thŽorme $L_\m(1,\chi)\ne 0$	\qquad	36	\par
thŽorme $U_\P=N_\P(U_{\gP})$	\qquad	69+122	\par
thŽorme 90 de Hilbert	\qquad	7	\par
thŽorme chinois	\qquad	6	\par
thŽorme d'approximation dŽbile	\qquad	9	\par
thŽorme de $\check {\rm C}$ebotarev	\qquad	39 et suiv.	\par
thŽorme de Bauer	\qquad	49	\par
thŽorme de Dirichlet sur les progressions arithmŽtiques	\qquad	40	\par
thŽorme de Dirichlet-Chevalley-Hasse	\qquad	96	\par
thŽorme de Kronecker-Weber	\qquad	84	\par
thŽorme de la norme de Hasse	\qquad	75	\par
thŽorme de la premire inŽgalitŽ du corps de classe	\qquad	38	\par
thŽorme de l'admissibilitŽ du conducteur	\qquad	123	\par
thŽorme de l'existence du corps de Hilbert	\qquad	123	\par
thŽorme de l'idŽal principal de la thŽorie des groupes	\qquad	127	\par
thŽorme de prŽsentation des corps locaux	\qquad	152	\par
thŽorme de prŽsentation d'un nombre premier par une forme quadratique	\qquad	163+164	\par
thŽorme de rŽciprocitŽ d'Artin	\qquad	83	\par
thŽorme de rŽciprocitŽ de Hilbert	\qquad	167	\par
thŽorme de rŽciprocitŽ pour le symbole des restes normiques	\qquad	151	\par
thŽorme de rŽciprocitŽ quadratique	\qquad	101	\par
thŽorme de surjectivitŽ de l'application d'Artin	\qquad	35	\par
thŽorme des idŽaux principaux de la thŽorie du corps de classe	\qquad	133	\par
thŽorme des normes d'une extension de Kummer locale	\qquad	156	\par
thŽorme des unitŽs de Dirichlet	\qquad	7	\par
thŽorme d'existence du corps de classe	\qquad	89+103	\par
thŽorme d'existence du corps de classe (version idŽlique)	\qquad	147+150	\par
thŽorme d'existence du corps de classe local	\qquad	159	\par
thŽormes de Galois	\qquad	3	\par
thŽormes d'isomorphismes, 	\qquad	1	\par
transversale de $H$ dans $G$	\qquad	125	\par
uniformisante	\qquad	63	\par
valeur absolue (archimŽdienne, non archimŽdienne)	\qquad	5	\par
valuation $\P$-adique	\qquad	63	\par
volume d'un idle	\qquad	141	\par
\vfill\eject

\centerline{\para Bibliographie}
\vskip2cm

\medskip
[Apo] : T. APOSTOL, {\it Mathematical Analysis}, Addison-Wesley, 1974
\medskip
[Con] : J. B. CONWAY, {\it Functions of One Complex Variable}, Springer, 1973
\medskip
[Cox] :  D. A. COX, {\it Primes of the form $x^2+ny^2$}, Wiley interscience, 1997
\medskip
[Fr-Tay] : A. FR¬\"OHLICH \& TAYLOR, {Algebraic number theory}, Cambridge studies in advanced mathematics 27, 1991
\medskip
[Has] : H. HASSE, {\it Beweis eines Satzes und Widerlegung einer Vermutung Ÿber des allgemeine Normenrestsymbol}, Nachrichten der Geselschaft den Wissenschaften zu Gšttingen, Math.-Phys. KI.HI (1931) 64-69
\medskip
[Jac1] : N. JACOBSON, {\it Basic Algebra 1}, Second Edition. New York, W.H. Freeman, 1989.
\medskip
[Jac2] : N. JACOBSON, {\it Basic Algebra 2}, Second Edition. New York, W.H. Freeman, 1989.
\medskip
[Jan] : G.J. JANUSZ, {Algebraic Number Fields}, Second Edition, AMS, 1996.
\medskip
[La1] : S. LANG, {\it Algbre}, troisime Ždition, Dunod, 2004.
\medskip
[La2] : S. LANG, {\it Algebraic Number Theory}, Second Edition, 1994
\medskip
[Mar] : D. MARCUS, {\it Number Fields}, Springer, 1977.
\medskip
[Nar] : W. NARKIEWICZ, {\it Elementary and Analytic Theory of Algebraic Numbers}, Springer, 1990.
\medskip
[Neu] : J. NEUKIRCH, {Algebraic Number Theory}, Springer, 1999
\medskip
[Ru] : W. RUDIN,    {\it Principles of Mathematical Analysis}, International Series
in Pure \& Applied Mathematics, McGraw-Hill, 1964.
\medskip
[Sam] : P. SAMUEL, {\it ThŽorie algŽbrique des nombres}, Hermann, 1971.
\medskip
[Ser] : J.-P. SERRE,  {\it Cours d'arithmŽtique}, Collection SUP No.2, Presses Universitaires de France,
Paris 1970.

\end

\vfill\eject
[Apo]= T. Apostol : Mathematical Analysis, Addison-Wesley, 1974
[Ru]=W.Rudin : Real and Complex analysis, Mc Graw Hil 1966
[Has]Beweis eines Satzes und Widerlegung einer Vermutung Ÿber des allgemeine Normenrestsymbol, Nachrichten der Geselschaft den Wissenschaften zu Gšttigen, Math.-Phys. KI.HI (1931) 64-69
[Cox] = Primes of the form $x^2+ny^2$ Wiley interscience
\uppercase\expandafter{\romannumeral\year}